%% file: greenfunctions.tex
\documentclass[11pt, a4paper]{amsart}
\usepackage{etex}
\usepackage{amscd,amsmath,amsthm,amssymb}
\usepackage{mathtools}
\usepackage{mathrsfs}
\usepackage{comment}
\usepackage[all,cmtip]{xy}
\usepackage{xcolor}
\usepackage{enumitem} 
\usepackage{array}
\usepackage{caption}
\usepackage{mathabx}

\usepackage{graphicx}

\usepackage{aurical}
\usepackage{amsbsy}
\usepackage{bm}

\usepackage{circuitikz}
\ctikzset{bipoles/resistor/height=0.15}
\ctikzset{bipoles/resistor/width=0.4}

\definecolor{cadmiumgreen}{rgb}{0.0, 0.42, 0.24}
\usepackage[
colorlinks, citecolor=cadmiumgreen,
pdfauthor={Omid Amini, Noema Nicolussi}, 
pdftitle={Hybrid moduli space II: tropical and hybrid Laplacians},
pdfstartview ={FitV},
]{hyperref}

\usepackage{manfnt}

\usepackage[
alphabetic,
msc-links,
nobysame,
lite,
]{amsrefs} 

\usepackage{tikz, float} 
\usetikzlibrary {positioning}

\usepackage{pgfplots}
\usetikzlibrary{patterns}
\pgfplotsset{compat=1.10}
\usepgfplotslibrary{fillbetween}

\usetikzlibrary{matrix, arrows}
\usetikzlibrary{patterns}

\usetikzlibrary{calc,decorations.markings}
\usetikzlibrary{shapes,snakes}

\usepackage{a4wide}
\usepackage[T1]{fontenc}
\usepackage[OT1]{fontenc}
\usepackage[utf8]{inputenc}
\usepackage[english]{babel}
\usepackage{epsfig}
\usepackage[leqno]{amsmath}
\usepackage{tabularx}
\usepackage{subfigure}
\usepackage{url}
\usepackage{hyperref}
\usepackage{tikz-cd}
\usepackage{csquotes}

\usepackage{mathtools}
\usepackage{manfnt}
\usepackage{mathrsfs}
\usepackage{scalerel}

\DeclareRobustCommand{\SkipTocEntry}[5]{}

\usepackage{xcolor}
\definecolor{darkgreen}{HTML}{00AA00}

\usepackage[disable]{todonotes}

\newtheorem{thm}{Theorem}[section]

\newtheorem{lem}[thm]{Lemma}

\newtheorem{prop}[thm]{Proposition}

\newtheorem{cor}[thm]{Corollary}

\newtheorem{defn}[thm]{Definition}

\newtheorem{question}[thm]{Question}
\theoremstyle{definition}
\numberwithin{equation}{section}

\theoremstyle{definition}
\newtheorem{defii}[thm]{Definition}
\newenvironment{defi}
  {\pushQED{\qed}\defii}
  {\popQED\enddefii}
\newtheorem{remm}[thm]{Remark}
\newenvironment{remark}
  {\pushQED{\qed}\remm}
  {\popQED\endremm}
\newtheorem{exx}[thm]{Example}
\newenvironment{example}
  {\pushQED{\qed}\exx}
  {\popQED\endexx}

\numberwithin{equation}{section}

\newcommand{\cf}{cf.}
\newcommand{\ie}{i.e.}

\renewcommand{\~}{\widetilde}

\newcommand{\R}{\mathbb{R}}
\newcommand{\Z}{{\mathbb Z}}
\newcommand{\N}{{\mathbb N}}
\newcommand{\eR}{\mathbb T}   
\newcommand{\eRm}{\widehat{\R}}   

\newcommand{\abs}[1]{\lvert #1\rvert}
\renewcommand{\hom}{\Hom}     
\newcommand{\rquot}[2]{#1\big/#2}       

\newcommand{\Div}{{\rm Div}}
\newcommand{\SiegAbel}{{\mathcal A}}    
\newcommand{\localsystem}{{\mathcal H}}   
\newcommand{\sH}{\mathcal H}               

\newcommand{\exact}{{\mathcal E}}  

\DeclareMathOperator{\sedbis}{sed}  
\newcommand{\sedind}[1]{{\scaleto{\sedbis=#1}{5.4pt}}}     
\newcommand{\combind}[1]{{\scaleto{#1}{5.4pt}}}
\newcommand{\grind}[1]{{\scaleto{#1}{4pt}}}
\newcommand{\TP}{\mathbb{TP}}   

\newcommand{\mass}{\scaleto{\mathscr M\hspace{-.1cm}\it{ass}}{6pt}}  
\newcommand{\Zh}{{_{\scaleto{\mbox{Zh}}{4.5pt}}}} 
\newcommand{\Ar}{{_{\scaleto{\mbox{Ar}}{4.5pt}}}} 
\newcommand{\Can}{\mathit{can}}
\newcommand{\can}{{{\scaleto{\Can}{2.8pt}}}}

\DeclareMathOperator{\Res}{Res}

\DeclareMathOperator{\Hom}{Hom}

\DeclareMathOperator{\id}{Id}

\renewcommand{\Im}{\operatorname{Im}}

\newcommand{\C}{{\mathbb C}}

\newcommand{\E}{{\mathbb E}}

\renewcommand{\H}{{\mathbb H}}

\renewcommand{\P}{{\mathbb P}}

\newcommand{\gr}{\mathit{gr}}

\newcommand{\grup}{{^{\mathit{gr}}}}
\newcommand{\surfup}{{^{\mathit{surf}}}}

\newcommand{\mA}{ {\mathcal{A}} } 
\newcommand{\mB}{ {\mathcal{B}} } 

\newcommand{\mg}{\mathscr M} 
\newcommand{\mgbar}{\comp{\mathscr M}} 
\newcommand{\mgg}[1]{\mg_{{\hspace{-.04cm}#1}}^{ }} 
\newcommand{\mgbarg}[1]{\mgbar_{{\hspace{-.09cm}#1}}^{ }} 
\newcommand{\unicurve}[1]{\hcurve_{{{\hspace{-.05cm}#1}}}^{}}

\newcommand{\mgtrop}[1]{\mg_{{\hspace{-.04cm}#1}}^\trop}   
\newcommand{\mgtropr}[2]{\mg_{{\hspace{-.04cm}#1}}^{^{\scaleto{\Tropit(#2)}{4.8pt}}}}  
\newcommand{\mgbartropr}[2]{\mgbar_{{\hspace{-.09cm}#1}}^{^{\scaleto{\Tropit(#2)}{4.8pt}}}} 
\newcommand{\mgtropcombin}[1]{\tilde{\mg}_{\hspace{-.02cm}{#1}}^{^{\scaleto{\Tropit}{3.8pt}}}} 
\newcommand{\unicurvetrop}[1]{\hcurve_{{\hspace{-.04cm}#1}}^\trop}
\newcommand{\mggraph}[1]{\mg_{{\hspace{-.04cm}#1}}^\metgra} 
\newcommand{\umggraph}[1]{\mathscr U\hspace{-.12cm}\mg_{{\hspace{-.04cm}#1}}^\metgra}
\newcommand{\umggraphcombin}[1]{\mathscr U\hspace{-.12cm}\tilde{\mg}_{\hspace{-.02cm}{#1}}^\metgra}
\newcommand{\mggraphcombin}[1]{\tilde{\mg}_{\hspace{-.02cm}{#1}}^\metgra}

\newcommand{\mghyb}[1]{\mg_{{\hspace{-.04cm}#1}}^\hyb}   
 
\newcommand{\mgbarhybr}[2]{\mgbar_{{\hspace{-.09cm}#1}}^{^{\scaleto{\Hyb(#2)}{4.8pt}}}}

\makeatletter
\newsavebox\myboxA
\newsavebox\myboxB
\newlength\mylenA

\newcommand*\overbar[2][0.75]{%
    \sbox{\myboxA}{$\m@th#2$}%
    \setbox\myboxB\null
    \ht\myboxB=\ht\myboxA%
    \dp\myboxB=\dp\myboxA%
    \wd\myboxB=#1\wd\myboxA
    \sbox\myboxB{$\m@th\overline{\copy\myboxB}$}
    \setlength\mylenA{\the\wd\myboxA}
    \addtolength\mylenA{-\the\wd\myboxB}%
    \ifdim\wd\myboxB<\wd\myboxA%
       \rlap{\hskip 1\mylenA\usebox\myboxB}{\usebox\myboxA}%
    \else
        \hskip -0.5\mylenA\rlap{\usebox\myboxA}{\hskip 0.5\mylenA\usebox\myboxB}%
    \fi}

\renewcommand\part{%
   \vspace*{2.5 ex}  
   \if@noskipsec \leavevmode \fi
   \@afterindentfalse
 \secdef \@part\@spart}

\def\@part[#1]#2{%
    \ifnum \c@secnumdepth >\m@ne
      \refstepcounter{part}%
      \addcontentsline{toc}{part}{\partname \nobreakspace \thepart:\hspace{1em}#1} 
    \else
      \addcontentsline{toc}{part}{#1}%
    \fi 
    {\parindent \z@ \centering 
     \interlinepenalty \@M
     \normalfont
     \ifnum \c@secnumdepth >\m@ne
	\LARGE \partname  \nobreakspace\thepart:
       \fi 
       \LARGE {\Fontskrivan{#2}}%
     \par}%
    \nobreak
    \vskip 5ex
    \@afterheading}
\def\@spart#1{%
    {\parindent \z@ \centering 
     \interlinepenalty \@M
     \normalfont
     \huge \bfseries #1\par}%
     \nobreak
     \vskip 5ex
     \@afterheading}

\makeatother

\renewcommand{\thepart}{\Roman{part}} 

\newcommand{\comp}[1]{\overbar[.5]{#1}} 

\newcommand{\st}{\bigm|} 

\newcommand{\dual}{\star}
\newcommand{\RMod}{{\R_+}} 

\newcommand{\zerocone}{{\underline0}}
\newcommand{\cancomp}[1]{\overline{#1}} 

\newcommand{\Tropcomp}[1]{\overline{#1}^{^{\mathit{trop}}}}

\newcommand{\varC}{\scaleto{\mathrm{C}}{6.4pt}}
\newcommand{\keg}{\varC}            
\newcommand{\canchart}{\~}                    
\newcommand{\inn}{\mathring}          
\newcommand{\face}{\scaleto{\mathrm{F}}{6.4pt}}
\newcommand{\sface}{\scaleto{\mathrm{f}}{6.4pt}} 

\newcommand{\rsf}{\mathcal S} 
\newcommand{\base}{B} 
\newcommand{\bp}{\flat} 
\newcommand{\Plit}{\mathit{pl}}                
\newcommand{\pl}{{^{\scaleto{\Plit}{4.5pt}}}}   
\newcommand{\plumb}{\rsf^{\pl}}  

\newcommand{\Tropit}{\mathit{trop}}                
\newcommand{\trop}{{^{\scaleto{\Tropit}{4pt}}}}   
\newcommand{\tropr}[1]{{^{\scaleto{\Tropit(#1)}{4.5pt}}}}  
\newcommand{\metgra}{{^{\scaleto{\mathit{gr}}{3.5pt}}}}
\newcommand{\tropdiag}{{{\scaleto{\Tropit}{4pt}}}} 
\newcommand{\phom}{\mathbb{P}\mathrm{Hom}}

\newcommand{\Hyb}{{\mathit{hyb}}}         
\newcommand{\hyb}{{^{\scaleto{\Hyb}{4.4pt}}}} 
\newcommand{\hybr}[1]{{^{\scaleto{\Hyb(#1)}{4.8pt}}}}  
\newcommand{\hybdiag}{{{\scaleto{\Hyb}{4pt}}}} 

\newcommand{\Log}{\mathrm{Log}} 
\newcommand{\LOG}{\mathrm{Log}} 
\newcommand{\Lognoind}{{\text{\Fontauri\slshape L}\scriptstyle{\textrm{og}}}} 
\newcommand{\Logsp}[1]{{\text{\Fontauri\slshape L}\scriptstyle{\textrm{og}}\strut^{\hspace{-2ex}{#1}}}} 
\newcommand{\Logspdiag}[1]{{\text{\Fontauri\slshape L}\scriptstyle{\scaleto{\textrm{og}}{3.5pt}}\strut^{\hspace{-1.2ex}{#1}}}} 

\newcommand{\logtrop}{\Logsp{\trop}} 
\newcommand{\loghyb}{\Logsp{\hyb}} 
\newcommand{\logtropind}[1]{\logtrop_{\hspace{-.4ex}#1}}

\newcommand{\loghybdiag}{\Logspdiag{\hybdiag}} 
\newcommand{\loghybdiagpi}[1]{\Logspdiag{\hybdiag}_{^{\hspace{-0.5ex}\scaleto{#1}{2pt}}}} 

\newcommand{\logtropdiag}{\Logspdiag{\tropdiag}} 

\newcommand{\varS}{\scaleto{\mathrm{S}}{6.8pt}}
\newcommand{\realtor}{\varS}  

\newcommand{\cT}{\mathcal{T}}   
\newcommand{\genusfunction}{\mathfrak g} 
\newcommand{\graphgenus}{{h}} 
\newcommand{\umgr}{\mathscr G} 
\newcommand{\mgr}{\mathcal G}  
\newcommand{\grm}[2]{\mathit{gr}_{_{\hspace{-.08cm}#1}}^{#2}} 
\newcommand{\proj}{{\kappa}} 
\newcommand{\cont}[1]{{{\proj_{\hspace{-.04cm}{\scaleto{#1}{5.4pt}}}}}}

\newcommand{\pr}{\mathfrak p} 

\newcommand{\forget}{\mathfrak q}

\newcommand{\marking}{{\mathfrak m}}
\newcommand{\countmarking}{{\mathfrak n}}

\newcommand{\filter}{\mathscr F} 
\newcommand{\dfilter}{\mathscr E} 

\newcommand{\mc}{{\mathcal{M}\curve}} 
\newcommand{\hcurve}{\mathscr C}  
\newcommand{\hcurveg}[1]{\hcurve_{_{\hspace{-.04cm}#1}}}
\newcommand{\tropcurve}{\mathscr C^{\trop}} 

\newcommand{\mccan}{\Sigma}
\renewcommand{\curve}{\mathcal  C}
\newcommand{\lf}{\bm{f}} 
\newcommand{\layg}{\bm{g}} 
\newcommand{\layh}{\bm{h}} 
\newcommand{\lz}{\bm{\zeta}} 

\newcommand{\lmu}{{\bm{\mu}}}
\newcommand{\lnu}{{\bm{\nu}}}

\newcommand{\lalpha}{\bm{\alpha}}

\newcommand{\thy}{\mathbf{t}} \newcommand{\shy}{\mathbf{s}}  

\newcommand{\subface}{\preceq}

\renewcommand{\~}{\widetilde}

\newcommand{\s}{\mathbf{s}} 
\newcommand{\contract}[2]{#1\big/#2}
\newcommand{\aut}{\mathrm{Aut}}
\newcommand{\Ical}{\mathcal I}

\newcommand{\rest}[1]{\raisebox{-1pt}{$\vert$}_{#1}}

\newcommand{\hp}[2]{\langle#2\rangle_{_{_{\hspace{-.05cm}#1}}}}
\newcommand{\lhp}[2]{{\pmb\langle}#2{\pmb\rangle}_{_{_{\hspace{-.05cm}#1}}}}

\newcommand{\innone}[2]{(#2)_{_{\hspace{-.02cm}#1}}}
\newcommand{\innoneind}[3]{(#3)_{_{_{{\hspace{-.05cm}#1}_{\hspace{-.02cm}{\scaleto{#2}{3.2pt}}}}}}}

\newcommand{\intprod}[2]{\langle#2\rangle_{_{_{\hspace{-.05cm}#1}}}}

\newcommand{\grs}{\mathrm{g}}  
\newcommand{\grg}{\mathrm{g}}  
\newcommand{\gri}[1]{{\mathrm{g}_{_{#1}}}} 
\newcommand{\grprimei}[1]{{\mathrm{g}'_{_{#1}}}} 
\newcommand{\grihat}[1]{{\hat{\mathrm{g}}_{_{#1}}}}
\newcommand{\gritilde}[1]{{\~{\mathrm{g}}_{_{#1}}}}
\newcommand{\grireg}[1]{{\mathrm{g}^\reg_{_{#1}}}} 
\newcommand{\gritildereg}[1]{{\~{\mathrm{g}}}^\reg_{_{#1}}}
\newcommand{\lgr}{\mathbf{g}}
\newcommand{\lgri}[1]{\lgr_{_{#1}}}
\newcommand{\jfunc}[1]{j_{_{{\hspace{-.05cm}\scaleto{#1}{4.1pt}}}}}  
\newcommand{\ljfunc}[1]{{\bm{j}}_{_{{\hspace{-.05cm}\scaleto{#1}{7 pt}}}}} 
\newcommand{\ljfuncbis}[1]{{\bm{j}}_{_{{\hspace{-.05cm}\scaleto{#1}{5.6pt}}}}} 
\newcommand{\jfunci}[1]{j_{_{{\hspace{-.04cm}\scaleto{#1}{5pt}}}}}
\newcommand{\jvide}{j}   

\newcommand{\supp}[1]{{|{#1}|}}
\newcommand{\tiret}{\text{-}}
\newcommand{\LS}{\mathcal H} 
\newcommand{\nr}{\mathrm{nr}} 
\newcommand{\weight}{{\scaleto{W}{4pt}}}   
\newcommand{\bipm}{{\widehat \Omega}} 
\newcommand{\group}{\textsc{G}} 
\newcommand{\gl}{\textrm{GL}} 
\newcommand{\col}{\mathrm{Col}} 
\newcommand{\row}{\mathrm{Row}}
\newcommand{\pdomain}{\mathscr{D}} 
\newcommand{\dimms}{N}  
\newcommand{\constant}{C}
\newcommand{\ExtM}{\widehat \rmM} 
\newcommand{\ELambda}{\widehat \Lambda}

\newcommand{\ind}[1]{_{_{\hspace{-.02cm}{\scaleto{#1}{4.1pt}}}}} 
\newcommand{\indbis}[1]{_{_{\hspace{-.02cm}{\scaleto{#1}{4.9pt}}}}} 
\newcommand{\slp}{\mathrm{sl}}

\newcommand\bH{{\mathbf{H}}}

\newcommand{\dom}{\operatorname{dom}}

\newcommand{\Deltatrop}{\bm{\Delta}}
\newcommand{\Deltahyb}{\bm{\Delta}}
\newcommand{\Deltaind}[1]{\Delta_{_{{\hspace{-.02cm}#1}}}}
\newcommand{\Deltahybind}[1]{\bm{\Delta}_{_{{\hspace{-.02cm}#1}}}}
\newcommand{\consts}[2]{\scaleto{K}{6pt}_{{\scaleto{#1}{4pt}}{_{,\scaleto{#2}{3pt}}}}} 
\newcommand{\const}[2]{\scaleto{K}{6pt}_{{\scaleto{#1}{3.2pt}}}^{\scaleto{#2}{5.5pt}}} 
\newcommand{\constwi}{\scaleto{K}{6pt}}
\newcommand{\smallc}{{_{\hspace{-0.04cm}\scaleto{\mathbb C}{4pt}}}}
\newcommand{\smallcc}{{{\scaleto{\mathbb C}{4pt}}}}

\newcommand{\ld}{\bm{d}}

\makeatletter
\let\@oldinfty\infty
\newcommand{\@sminfty}{{\hspace{-.02cm}\scaleto{\@oldinfty}{2.8pt}}} 
\newcommand{\@smsminfty}{_{\hspace{-.04cm}\scaleto{\@oldinfty}{2.8pt}}} 
\renewcommand{\infty}{{\mathchoice%
  {\displaystyle{\@oldinfty}}%
  {\textstyle{\@oldinfty}}%
  {\scriptstyle{\@sminfty}}%
  {\scriptscriptstyle{\@smsminfty}}}
}
\makeatother

\makeatletter
\let\@oldpi\pi
\newcommand{\@smsmpi}{_{\hspace{-.04cm}\scaleto{\@oldpi}{3.2pt}}} 
\renewcommand{\pi}{{\mathchoice%
  {\displaystyle{\@oldpi}}
  {\textstyle{\@oldpi}}%
  {\scriptstyle{\@oldpi}}%
  {\scriptscriptstyle{\@smsmpi}}}
}
\makeatother

\newcommand{\fin}{{\scaleto{\mathfrak f}{8pt}}}

\makeatletter
\let\@oldfin\fin
\newcommand{\@smfin}{{\hspace{-.01cm}\scaleto{\@oldfin}{5.4pt}}} 
\newcommand{\@smsmfin}{{\hspace{-.01cm}\scaleto{\@oldfin}{5pt}}} 
\renewcommand{\fin}{{\mathchoice%
  {\displaystyle{\@oldfin}}
  {\textstyle{\@oldfin}}%
  {\scriptstyle{\@smfin}}%
  {\scriptscriptstyle{\@smsmfin}}}
}
\makeatother

\newcommand{\sed}{{\scaleto{\mathrm{sed}}{7pt}}}

\makeatletter
\let\@oldsed\sed
\newcommand{\@smsed}{{\hspace{-.01cm}\scaleto{\@oldsed}{4.4pt}}} 
\newcommand{\@smsmsed}{{\hspace{-.01cm}\scaleto{\@oldsed}{4pt}}} 
\renewcommand{\sed}{{\mathchoice%
  {\displaystyle{\@oldsed}}
  {\textstyle{\@oldsed}}%
  {\scriptstyle{\@smsed}}%
  {\scriptscriptstyle{\@smsmsed}}}
}
\makeatother

\newcommand{\EE}{\mathbb E}
\newcommand{\oo}{\mathfrak o}
\newcommand{\head}{\mathrm{h}}
\newcommand{\tail}{\mathrm{t}}
\renewcommand{\div}{\mathrm{div}}
\newcommand{\divind}[2]{\mathrm{div}_{_{\hspace{-.05cm}\scaleto{#1\prec #2}{5pt}}}}

\tikzstyle{Cwhite}=[scale = .8,circle, fill = white, minimum size=3mm] 
\tikzstyle{Cgray}=[scale = .4,circle, fill = gray, minimum size=3mm] 
\tikzstyle{Cblack2}=[scale = .4,circle, fill = black, minimum size=5mm] 
\tikzstyle{Cblack}=[scale = .7,circle, fill = black, minimum size=3mm]
\tikzstyle{C0}=[scale = .9,circle, fill = black!0, inner sep = 0pt, minimum size=3mm]
\tikzstyle{C1}=[scale = .7,circle, fill = black!0, inner sep = 0pt, minimum size=3mm]
\tikzstyle{Cred}=[scale = .4,circle, fill = red, minimum size=3mm]

\newcommand{\cube}{\mbox{\,\mancube}}
\newcommand{\one}{{\mathbf 1}}

\newcommand{\rmw}{\mathrm{w}}
\newcommand{\rmW}{\mathrm{W}}

\newcommand{\rmZ}{\mathrm{Z}}
\newcommand{\rmM}{\mathrm{M}}

\newcommand{\reg}{{^{\scaleto{\mathrm{reg}}{4pt}}}} 

\newcommand{\fsed}{{\scaleto{\mathrm{fulsed}}{4pt}}}
\newcommand{\Pihat}{{\widehat{\Pi}}}
\newcommand{\Piall}{{\widehat{\Pi}}}
\newcommand{\Pifs}{\Pi}

\newcommand{\tameclass}{\mathfrak C}

\newcommand{\Ap}{{\mathcal{A}}}
\newcommand{\whAp}{\widehat{\Ap}}
\newcommand{\Apc}{{\mathcal{L}}}
\newcommand{\whApc}{\widehat{\Apc}}

\newcommand{\fell}{{\mathfrak{l}}}
\newcommand{\fEll}{{\mathfrak{L}}}

\newcommand{\suppdiv}{\Lambda}

\newcommand{\nmark}{m}

\usepackage{stmaryrd}
\usepackage{trimclip}

\makeatletter
\DeclareRobustCommand{\shortto}{%
  \mathrel{\mathpalette\short@to\relax}%
}

\newcommand{\short@to}[2]{%
  \mkern2mu
  \clipbox{{.5\width} 0 0 0}{$\m@th#1\vphantom{+}{\shortrightarrow}$}%
  }
\makeatother

\newcommand{\cmc}{{\smallcc, \scaleto{\mathrm{mc}}{2.5pt}}}

\newcommand{\projhar}{\mathrm{P}_{_{\hspace{-.14cm}\perp}}^{^{\hspace{-.03cm}\diamond}}}
\newcommand{\projexact}{\mathrm{P}_{_{\hspace{-.14cm}\perp}}^{^{\hspace{-.03cm}\mathrm{ex}}}}
\newcommand{\transpose}{{\scaleto{\mathrm{T}}{4.5pt}}}

\begin{document}
\title[Moduli of hybrid curves II: Tropical and hybrid Laplacians]{Moduli of hybrid curves II: \\ Tropical and hybrid Laplacians}

\author{Omid Amini}
\address{CNRS - Centre de math\'ematiques Laurent Schwartz, \'Ecole Polytechnique}
\email{\href{omid.amini@polytechnique.edu}{omid.amini@polytechnique.edu}}

\author{Noema Nicolussi}
\address{Faculty of Mathematics, University of Vienna}
\email{\href{noema.nicolussi@univie.ac.at}{noema.nicolussi@univie.ac.at}}

\date{\today}

\begin{abstract} 
The present paper is a sequel to our work on hybrid geometry of curves and their moduli spaces. We introduce a notion of \emph{hybrid Laplacian}, formulate a \emph{hybrid Poisson equation}, and give a mathematical meaning to the convergence both  of the Laplace operator and the solutions to the Poisson equation on Riemann surfaces.  As the main theorem of this paper, we then obtain a \emph{layered  description} of the asymptotics of Arakelov Green functions on Riemann surfaces close to the boundary of their moduli spaces. This is done in terms of a suitable notion of \emph{hybrid Green functions}.

\smallskip

As a byproduct of our approach, we obtain other results of independent interest. In particular, we introduce  \emph{higher rank canonical compactifications of fans and polyhedral spaces} and use them to define \emph{the moduli space  of higher rank tropical curves}. Moreover, we develop the first steps of a function theory in \emph{higher rank non-Archimedean, hybrid, and tame analysis}. Furthermore, we establish the convergence of the Laplace operator on metric graphs toward the \emph{tropical Laplace operator} on limit tropical curves in the corresponding moduli spaces, leading to new perspectives in operator theory on metric graphs.  

\smallskip

 Our result on the Arakelov Green function is inspired by the works of several authors, in particular those of Faltings, de Jong, Wentworth and Wolpert, and solves a  long-standing open problem arising from the Arakelov geometry of Riemann surfaces.  The hybrid layered behavior close to the boundary of moduli spaces is expected to be a broad phenomenon and will be explored in our forthcoming work.

  \end{abstract}

\maketitle

\setcounter{tocdepth}{1}

\tableofcontents


\section{Introduction}

The main motivation behind the study undertaken in this paper is to understand the behavior of the Arakelov Green functions on degenerating families of Riemann surfaces. Recall that any smooth Riemann surface $S$ can be equipped with its  Arakelov-Bergman measure $\mu_{\Ar}$ defined in terms of the space of holomorphic one-forms on $S$,
\begin{align*}
\mu_\Ar \coloneqq \frac i{2\pi} \sum_{j=1}^g \omega_j\wedge \bar \omega_j 
\end{align*}
where $\omega_1, \dots, \omega_g$ form an orthonormal basis of the space of holomorphic one-forms on $S$ for the hermitian inner product
\[
\innone{S}{\eta_1, \eta_2} \coloneqq \frac i2 \int_S \eta_1 \wedge \bar \eta_2.
\]
The measure $\mu_{\Ar}$ has total mass $g$, and we define the canonical measure of $S$ by normalization
\[\mu^\can \coloneqq \frac 1g \mu_\Ar.\]

To the canonically measured Riemann surface $S$ we can associate its Green function
\[
\begin{array}{cccc}
\gri{S} \colon &S \times S - \mathrm{diag}_S  &\longrightarrow & \R
\end{array}
\]
 defined as the unique solution of the equations
\begin{align*}
\frac{1}{\pi i}\partial_z \partial_{\bar{z}}\, \grs(p, \cdot) &= \delta_p - \mu^{\can}    , \qquad \int_S \grs(p, y) \, d\mu^{\can}(y) = 0\end{align*}
valid for all $p \in S$. Here, $\mathrm{diag}_S$ is the diagonal embedding of $S\hookrightarrow S\times S$.

\smallskip
Arakelov Green functions naturally arise in Arakelov geometry where they allow to appropriately include the contributions of the places at infinity when doing intersection theory on arithmetic surfaces~\cite{Ara74, Fal84, Bost99, Lang, Zhang} (the picture moreover extends  to higher dimensional arithmetic varieties~\cite{GS90, BGS94}). Due to their connections to fundamental questions in arithmetic geometry and number theory, and in physics (in Polyakov formulation of string theory and in perturbative quantum field theory), they have been extensively studied from different perspectives in the past, see e.g.~\cite{JK06, JK09, Zhang2, Cinkir, Went91, Faltings21, deJong, Smit88, Went08, JW92, HGB19, Wil17}.

\medskip

 Let $\mgg{\grind{g}}$ be the moduli space of curves of genus $g$ and denote by $\mgbarg{\grind{g}}$ the Deligne--Mumford compactification 
consisting of the stable curves of genus $g$. The question underlying this paper can be innocently stated as follows:

\begin{question}\label{question:main-surfaces}
Consider a sequence of smooth compact Riemann surfaces $S_1, S_2, \dots$ of genus $g$ such that the corresponding points $s_1, s_2, \dots$ converge to a point in $\mgbarg{\grind{g}}$.

\smallskip

What is the limit of the canonical Green functions $\gri{S_1}, \gri{S_2}, \dots$?
\end{question}

\smallskip

To formalize our question, we will work with the moduli space of curves with two marked points $\mgg{\combind{g,2}}$. The Arakelov Green function can be naturally defined on $\mgg{\combind{g,2}}$. Namely, each point in $\mgg{\combind{g,2}}$ corresponds to a Riemann surface $S$ together with two distinct points $q, y \in S$. Hence we can naturally define
\[
\begin{array}{cccc}
\grs \colon &\mgg{\combind{g,2}}  &\longrightarrow & \R \\
& (S, q,y) &\mapsto & \gri{S}(q,y),
\end{array}
\]
with $\gri{S}(\cdot, \cdot)$ the Green function of $S$. 

Consider the forgetful map $\pi$ from $\mgg{\combind{g,2}}$ to $\mgg{\grind{g}}$, which simply forgets the markings,
\[
\begin{array}{cccc}
\pi \colon &\mgg{\combind{g,2}}  &\longrightarrow & \mgg{\grind g} \\
& (S, q,y) &\mapsto & S.
\end{array}
\]
Using this map, we can naturally ask for the shape of the Arakelov Green function $\gri{s}\coloneqq \grs\rest{\pi^{-1}(s)}$, $s \in \mgg{\grind{g}}$, when $s$ is close to a fixed boundary point $s_0$ on the boundary $\partial \mgbarg{\grind{g}} = \mgbarg{\grind{g}} \setminus \mgg{\grind{g}}$. 

\smallskip

However, in our previous work~\cite{AN} it became clear that the Deligne--Mumford compactification is not the right compactification to study the variations of the Arakelov--Bergman measures $\mu_{\Ar}$. Since the Arakelov Green function of a Riemann surface is defined relative to the measure, essentially the same problem arises in the study of the Green functions.

\smallskip

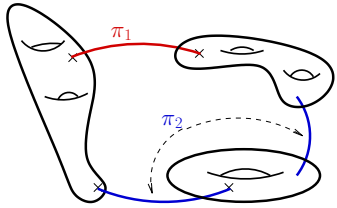
\begin{figure}[!t]
\centering
    \scalebox{.35}{\input{example3-2.tikz}}
\caption{An example of a hybrid curve. The graph of the underlying stable Riemann surface is a $3$-cycle. Its edges are partitioned in two sets $\pi_1$ and $\pi_2$.}
\label{fig:HybridCurveIntro}
\end{figure}

On the other hand, we resolved these difficulties in~\cite{AN} by working with a larger compactification of $\mgg{\grind{g}}$: \emph{the moduli space $\mghyb{\grind{g}}$ of hybrid curves of genus $g$}. Each point $\thy \in \mghyb{\grind{g}}$ represents a hybrid curve, that is, a geometric object associated to a stable Riemann surface $S$ together with an ordered partition of the edge set $E$ of its dual graph and a (suitably normalized) edge length function on different parts of the partition. One can think of constructing a hybrid curve as taking the components of the stable Riemann surface $S$, and connecting them with intervals attached at the nodes of $S$. The ordered partition sorts these intervals into groups and the edge length function assigns them lengths,  leading to an associated (layered) metric graph with normalized edge lengths in each layer (see Figure~\ref{fig:HybridCurveIntro}). The number of parts in the ordered partition is called the \emph{rank} of the hybrid curve.

\smallskip

Hybrid curves arise naturally in the study of multiparameter families of algebraic curves. The situation to have in mind is a generically smooth family of complex projective curves $\mathcal X$ over a multidimensional base $B$, with a discriminant locus $D$, which one can suppose to be a simple normal crossing divisor after a base change to a log-resolution. In this situation, hybrid curves replace singular fibers of the family and give rise to a hybrid family defined on the hybrid replacement $B^\hyb$ of the pair $(B, D)$.

\smallskip

The main result of~\cite{AN} states that the Arakelov--Bergman measure on families of degenerating Riemann surfaces converges to a certain \emph{hybrid measure} on a hybrid curve, that is, the limit measure is a combination of canonical measures of Riemann surfaces and canonical measures of (layered) metric graphs. This measure can be thought of as a hybrid object, in the sense that it mixes notions from analysis on Riemann surfaces and graphs, or more broadly speaking, from complex geometry and tropical geometry, reflecting a hybrid of both continuous and discrete features.

\smallskip

Continuing with this approach, we rephrase Question~\ref{question:main-surfaces} as follows:
\begin{question} \label{question:main-surfaces2} 
What is the shape of $\gri{s}$, $s \in \mgg{\grind{g}}$, for $s$ close to a fixed point $\thy$ in the hybrid moduli space $\mghyb{\grind{g}}$?
\end{question}

\smallskip

The main goal of the present paper is to answer essentially fully this question, by demonstrating that the Arakelov Green function exhibits a \emph{hybrid structure} as well.

\medskip

Our central contribution consists in introducing a {\em hybrid Laplacian} and formulating the {\em hybrid Poisson equation}
\begin{equation} \label{eq:IntroHybridPoissonEquation}
\Deltahyb \lf = \lmu 
\end{equation}
on hybrid curves. The latter, in particular, leads to a suitable notion of \emph{hybrid Green functions}.

\medskip

The hybrid Laplacian is a \emph{layered measure valued} operator, defined on functions on hybrid curves. Functions on hybrid curves are vector valued, with values taken in the space $\R^{r+1}$, for $r$ the rank of the hybrid curve. The hybrid Laplacian contains both the Archimedean and non-Archimedean notions of Laplacians, on the corresponding  Riemann surfaces and metric graph layers of the hybrid curve, and additionally measures the impact of \emph{higher order infinity layers} on \emph{lower order layers}. Altogether, the hybrid Poisson equation can be thought of as a system of $r+1$ Poisson equations, posed on the Riemann surface and metric graph layers of the hybrid curve, which are coupled by contributions of higher order to lower order layers.

\smallskip

Getting a well-posed Poisson equation in this setting is a delicate problem, as the space of solutions has dimension higher than expected, due to the possibility of disconnected geometric objects appearing in layers. We solve this by introducing the space of \emph{harmonically arranged functions on hybrid curves}, which ensures existence and uniqueness of solutions, modulo global constants.   

\smallskip
 
Once this has been achieved, we pull-back the hybrid Green functions to Riemann surfaces, by using an appropriate notion of \emph{logarithm  map} on moduli spaces and their universal families, and show that the Arakelov Green function admits an approximate \emph{layered development} in terms of the components of the hybrid Green function.

\medskip

Although our main motivation is to understand the limit of the Arakelov Green function, one should view the \emph{formulation of the hybrid Poisson equation} as the principal result of the present work. In fact, in our approach we prove that an approximative layered development, using solutions of hybrid Poisson equations, holds for \emph{solutions of general Poisson equations} on degenerating Riemann surfaces. Specializing to the equation defining the Arakelov Green function, the main theorem can be viewed as a consequence of this principle,  holding more generally for Green functions of measures with suitable continuity properties on the moduli space of hybrid curves. In the case of Arakelov Green functions, this continuity comes from our previous work~\cite{AN}.

\medskip 

\begin{figure}[!t]
\centering \scalebox{.3}{\input{example2-2.tikz}}
\caption{The map $\Psi$ sends an almost singular Riemann surface to a metric graph.}
\label{fig:IntroTropicalizationMap}
\end{figure}
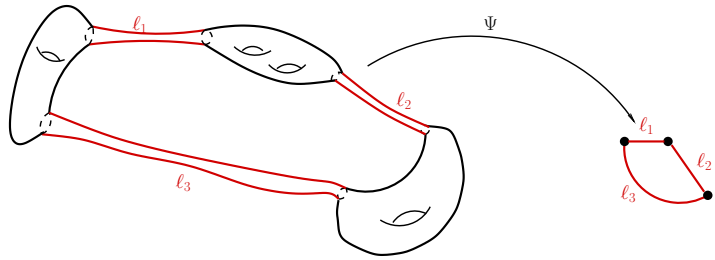

Our approach relies crucially on tools from \emph{analysis on metric graphs}. Close to a stable Riemann surface $S$ of genus $g$, there is a \emph{log map} from $\mgg{\grind{g}}$ to the moduli space of genus $g$ metric graphs $\mggraph{\grind g}$ (defined in local étale charts), which sends each smooth Riemann surface to a metric graph over the dual graph of $S$ (see Figure~\ref{fig:IntroTropicalizationMap}). If the smooth Riemann surface degenerates, the associated metric graph degenerates in the sense that the edges become infinitely long. It turns out that, in order to understand the behavior of Poisson solutions on degenerating Riemann surfaces, one  has to address first the analogous problem on degenerating metric graphs.

\smallskip

We solve this problem by introducing a new \emph{compactification of the metric graph moduli spaces} which adds new geometric objects, the so-called \emph{tropical curves}, as boundary points. These tropical curves contain the ones introduced in our previous work~\cite{AN}, and provide a higher rank refinement of the current notion in the literature. Defining as well a \emph{tropical Laplace operator} and \emph{tropical Poisson equation} on tropical curves, we obtain a {layered approximative development} of solutions to Poisson equations and Green functions on degenerating metric graphs, in terms of respective solutions and Green functions on tropical curves. In short, this framework is essentially parallel to the one for degenerating Riemann surfaces, with \emph{hybrid curves replaced by tropical curves}. The considerations on metric graphs form a self-contained part of our paper, which can be read independently of the rest. Since these results are of independent interest, we describe them in the introduction as well.

\medskip

In the rest of the introduction, we introduce the appropriate framework and language for stating the main theorems, and explain the precise meaning of the above. Comparison to the other works, such as the (recent) works of Faltings~\cite{Faltings21}, de Jong~\cite{deJong}, and the (pioneering) work of Wentworth~\cite{Went91}, as well as discussions of future developments will appear at relevant sections of this introduction.

\subsection{Metric graphs} In order to address Question~\ref{question:main-surfaces2}, we consider first the related problem of describing the canonical Green functions on \emph{degenerating families of metric graphs}. This leads to the introduction of several important concepts.

\smallskip

Let $\mgr$ be a metric graph corresponding to a pair $(G, l)$ of a finite graph $G=(V, E)$ and an edge length function $l\colon E \to (0, + \infty)$. Recall that $\mgr$ is the metric space obtained by identifying every edge $e\in E$ with an interval of length $l(e)$ in $\R$, and then gluing these intervals at their extremities following the incidence relations between edges and vertices of the graph.

\smallskip

In context of the present work, metric graphs arise as approximations of Riemann surfaces, and as so, they are \emph{augmented} meaning that they come with integral weights associated to the vertices.  This is encoded in the \emph{genus function} $\genusfunction \colon V \to \N \cup\{0\}$ which associates to each vertex $v$ its \emph{genus} $\genusfunction(v)$. The genus $g$ of the metric graph $\mgr$ can be then defined as the genus of the underlying graph $G$, added up with the sum of the genera of vertices. When dealing with (global) moduli spaces, we moreover have to impose the so-called stability condition, that is, $2\genusfunction(v)+\deg(v) -3 \geq0$ at each vertex $v$ of $G$. (This is not necessary if we are only interested in metric graphs of given combinatorial type.)

\smallskip

The {\em Zhang measure} on an (augmented) metric graph $\mgr$ of genus $g$ is the measure
\begin{equation} \label{eq:FosterCoefficient}
\mu_\Zh :=   \sum_{e \in E} \frac{\mu (e)}{l(e)} d\lambda_e +\sum_{v\in V}\genusfunction(v)\delta_v,
\end{equation}
where $d\lambda_e$ denotes the uniform Lebesgue measure on the edge $e \in E$, and  $\mu(e)$ is the {\em Foster coefficient at $e$}, $e \in E$, given in terms of the spanning trees of $G$ by the following equation
 \[\mu(e) =  \frac{\sum_{T: e\notin T} w(T)}{\sum_{T} w(T)}. \]
 In the above sum, $T$ runs over the spanning trees of $G$ and the weight of $T$ is the product $w(T) = \prod_{e\in T} l(e)$. Altogether, $\mu(e)$ is the probability that a random spanning tree, chosen according to a hardcore distribution given by the edge lengths, does not contain the edge $e$.
 
 The measure $\mu_\Zh$ has total mass $g$, and we define the {\em canonical  measure} of $\mgr$ by normalization
 \[\mu^\can:=\frac 1g \mu_\Zh.\]

\smallskip

The {\em canonical Green function} on  $\mgr$ is the unique function $\gri{\mgr} \colon \mgr \times \mgr \to \R$ such that for each point $p \in \mgr$, the function $\gri{\mgr}(p, \cdot)$ satisfies the Poisson equation with normalization
\begin{equation} \label{eq:DefGFGraph}
\Delta \gri{\mgr}(p, \cdot)=  \delta_p - \mu^{\can}, \qquad \int_\mgr \gri{\mgr}(p, y) \, d\mu^{\can}(y) = 0.
\end{equation}
The operator $\Delta=\Deltaind{\mgr}$ in this equation is the metric graph Laplacian, which is the measure-valued operator on $\mgr$ mapping a function $f\colon \mgr \to \R$ to the measure
\[
	\Delta f = - f'' d \lambda - \sum_{x \in \mgr} \Big( \sum_{\nu \in T_x} \slp_\nu f(x) \Big) \delta_x.
\]
In the above definition, $T_x$ is the set of (outgoing) unit tangential directions at $x$ and $\slp_\nu f(x)$ is the slope of $f$ at $x$ in the direction of $\nu$. The canonical Green function $\gri{\mgr}(p, \cdot)$ is a piecewise quadratic function on $\mgr$ for each $p \in \mgr$. It moreover encodes important arithmetic and geometric information when the metric graph is associated to arithmetic surfaces~\cite{Zhang, Zhang2, Cinkir, DKY20}.

\medskip

For a given graph $G$, the \emph{space of metric graphs of combinatorial type $G$} can be identified with the open cone $\mggraphcombin{\grind{G}}:= \R_{>0}^E$. Gluing these cones appropriately, when $G$ varies over stable graphs of given genus $g$, one obtains the moduli space of metric graphs $\mggraph{\grind{g}}$. 

The moduli space $\mggraph{\grind g}$ of genus $g$ stable metric graphs is not compact. Allowing the edge lengths to take values infinity leads to the compactified moduli space $\cancomp{\mggraph{\grind g}}$ studied in~\cite{ACP}.

\smallskip

In context with questions addressed in this paper, one can think of $\cancomp{\mggraph{\grind g}}$ as the metric graph analogue of the Deligne--Mumford compactification $\mgbarg{\grind{g}}$. Close to a stable Riemann surface $S$ of genus $g$, there is a natural log map from $\mgbarg{\grind{g}}$ to $\cancomp{\mggraph{\grind g}}$ (defined in local étale charts). The non-smooth Riemann surfaces are exactly the ones which are sent to metric graphs having some edges of infinite length.

\smallskip

Hence, the metric graph analogue of Question~\ref{question:main-surfaces} reads as follows.
 
\begin{question}\label{question:main-graphs}
Consider a sequence of stable metric graphs $\mgr_1, \mgr_2, \dots$ of genus $g$ such that the corresponding points in $\mggraph{\grind g}$ converge to a point in $\cancomp{\mggraph{\grind g}}$.

\smallskip

What is the limit of the canonical Green functions $\gri{\mgr_1}, \gri{\mgr_2}, \dots$?
\end{question}

\begin{figure}[!t]
\centering \centering
    \scalebox{.5}{\input{example8.tikz}}
\caption{The metric graph $\mgr$ discussed in Example~\ref{ex:CounterExampleCanonicalMeasure}. Depending on the way the edge lengths $\ell_1, \ell_2, \ell_3$ converge to infinity, one gets different limit measures.}
\label{fig:BB}
\end{figure}
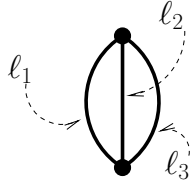

Just as in the case of Riemann surfaces, it turns out that the space  $\cancomp{\mggraph{\grind g}}$ is not the right compactification for studying variations of the canonical  measures $\mu^{\can}$. Simple examples, such as the one depicted in Figure~\ref{fig:BB}, already show that different limit measures may arise if the edge lengths go to infinity, depending on the relative speed of convergence to infinity. We refer to Example~\ref{ex:CounterExampleCanonicalMeasure} discussed in Section~\ref{sec:canonical_green_functions_tropical}.

Moreover, since the canonical Green function is again defined relative to $\mu^{\can}$, we have to confront essentially the same problem as in the study of Riemann surfaces.

\smallskip

This raises the question of constructing an \emph{appropriate compactification} of the moduli space $\mggraph{\grind g}$, suitable to study the variation of the canonical measure and Green function.  The answer leads us to the \emph{botany of canonical compactifications}, a large classe of compactifications of fans and cone complexes with origin in higher rank valuation theory.

\smallskip

The finest among all these compactifications is the absolute tropical compactification of $\mggraph{\grind g}$, which turns out to be the \emph{moduli space} for a new geometric object: \emph{tropical curves  of genus $g$} (\emph{of arbitrary rank}). So we begin by presenting these objects.

\subsection{Higher rank tropical curves} The key concept in introducing hybrid compactifications in our work~\cite{AN} are ordered partitions. In order to define a relevant notion of higher rank tropical curves and tropical compactifications, we need to extend the set-up of~\cite{AN} by allowing a (possibly empty) \emph{finitary part}.  

\smallskip

 Let $E$ be a non-empty finite set. An \emph{ordered partition} of $E$ is an ordered pair $\pi = (\pi_\infty, \pi_\fin)$ consisting of 
 \begin{itemize}
 \item a subset $\pi_\fin \subseteq E$ called the \emph{finitary part} of $\pi$, and 
 \item an ordered sequence $\pi_\infty = (\pi_1, \pi_2, \dots, \pi_r)$,  $r\in \N \cup\{0\}$, consisting of non-empty subsets $\pi_j \subseteq E$ which form a partition of $E\setminus \pi_\fin$. The ordered sequence $\pi_\infty$ is called the \emph{sedentarity part} of $\pi$.
  \end{itemize}
 Equivalently, an ordered partition is an ordered sequence $\pi = (\pi_1, \pi_2, \dots, \pi_r, \pi_\fin)$, $r\in \N \cup\{0\}$, of  pairwise disjoint subsets of $E$ such that the following holds:
 \begin{itemize}
 \item $\pi_1, \dots, \pi_r$ are non-empty, and 
 \item we have 	$E = \bigsqcup_{i=1, \dots, r, \fin} \pi_i.$
 \end{itemize}
Note that we allow $\pi_\fin$ to be empty. Such a partition is then called of \emph{full sedentarity}.  Moreover, the integer $r$, which is called the \emph{rank} of $\pi$, is allowed to be zero. In this case we obtain $\pi_\infty =\varnothing$ and $\pi$ consists of $\pi_\fin =E$. 

\smallskip

 This notion of sedentarity can be viewed as a higher rank analog of sedenarity currently in use in tropical geometry. It leads to higher rank canonical compactifications of fans and polyhedral and cone complexes, and extends to higher rank the constructions which have recently played a central role in the developments of tropical geometry in relation to Hodge theory~\cite{AP-homology,AP,AP-CS}.  We will elaborate on this later in Section~\ref{sec:tropical_moduli}. 

\smallskip

Conceptually, one can imagine the sedentarity part $\pi_{\infty}=(\pi_1, \dots, \pi_r)$ as combining parts of the ordered partition which lie at  $r$ distinguished \emph{layers of infinities}, separated \emph{infinitely far away} from each other, yielding a vertical tower which is ordered top to bottom by the integers $1, \dots, r$. 

\smallskip

In the degeneration picture for metric graphs, when the edge lengths on a given graph $G=(V,E)$ vary with some of them going to infinity, the part $\pi_1$ in the limit will contain the edges with the fastest convergence to infinity, $\pi_2$ those among the remaining ones with the fastest convergence to infinity, and so on. On the other hand, the finitary part $\pi_\fin$ will contain all the edges whose lengths remained bounded.

\smallskip

We call the limit geometric object emerging from the above picture a \emph{tropical curve},  defined formally as follows. The data underlying the definition  consists of
 \begin{itemize}
 \item a stable graph $G=(V, E, \genusfunction)$,
 \item an ordered partition $\pi=(\pi_\infty, \pi_\fin)$, $\pi_\infty = (\pi_1,\dots, \pi_r)$, and
 \item an edge length function $l\colon E \to (0, + \infty)$ which is well-defined only up to positive rescaling of the layers, that is, modulo multiplication of $l_j:=l\rest{\pi_j}$ by a positive real number, for $j=1, \dots, r$.
 (Imposing $\sum_{e\in\pi_j} l(e)=1$, $j=1, \dots, r$, leads to a unique choice of $l$.)
  \end{itemize}
In the above degeneration picture for metric graphs, the edge length function $l$ encodes the \emph{relative speed of convergence to infinity} of the edge lengths in each set $\pi_j$, $j= 1, \dots, r$, of the sedentarity part.
  
 With this data, we define the tropical curve $\curve$ as the metric realization of the pair $(G, l)$ enriched with the data of the ordered partition $\pi$. By definition, the \emph{rank of a tropical curve} is the rank of the underlying ordered partition. 
   
 Alternatively, a tropical curve can be defined via \emph{tropical uniformization} as a \emph{conformal equivalence class} of layered metric graphs, with layering induced by the ordered partition $(\pi_\infty, \pi_\fin)$ and the conformal equivalence defined by multiplication of $l_j:=l\rest{\pi_j}$ on infinitary layers by a positive real number, for $j=1, \dots, r$. In the first definition, we associated at the same time a \emph{canonical metric graph} $\Gamma$ to $\curve$, the metric realization of the normalized edge length function $l$.
 
 As in~\cite{AN}, we define the \emph{graded minors} $\Gamma^1, \dots, \Gamma^r, \Gamma^\fin$ of $\curve$ as follows. The $j$th graded minor is the metric graph $\Gamma^j$ obtained by removing all the edges which appear in lower indexed layers $\pi_1, \dots, \pi_{j-1}$ and then contracting all the edges which appear in upper indexed layers $\pi_{j+1}, \dots, \pi_r, \pi_\fin$.
 
\subsection{Moduli of higher rank tropical curves} We define $\mgtrop{\grind{g}}$ as the moduli space of tropical curves of genus $g$  (of arbitrary rank). That is, the points in $\mgtrop{\grind{g}}$ correspond to the isomorphism classes of tropical curves of genus $g$. The space $\mgtrop{\grind{g}}$ coincides with the \emph{absolute  canonical compactication of the cone complex $\mggraph{\grind g}$}, and belongs to the class of higher rank canonical compactifications of $\mggraph{\grind g}$. Let us briefly explain what we mean by this and then refer to the relevant section in the text for more information.

\smallskip

For each pair $(p,q)$ with $p+q -1 \leq N:=3g-3$, we define a canonical compactification $\mgbartropr{\grind{g}}{p,q}$ of $\mggraph{\grind g}$, coming therefore with an embedding $\mggraph{\grind g} \hookrightarrow \mgbartropr{\grind{g}}{p,q}$, 
 with the following properties:

\begin{itemize}
\item for each pair $(p,q)$ with $p+q=N+1$, we have $\mgbartropr{\grind{g}}{p,q} = \mgtrop{\grind{g}}$.
\smallskip

\item for two pairs $(p,q)$ and $(p',q')$ with  $p' \le p$ and $q' \le q$, we have a surjective forgetful map from $\mgbartropr{\grind{g}}{p,q}$ to $\mgbartropr{\grind{g}}{p',q'}$.
\end{itemize}
Altogether,  we obtain  a commutative diagram of maps depicted in Figure~\ref{fig:CompactificationTower-intro}. The original compactification  $\cancomp{\mggraph{\grind g}}$ appears as the element  $\cancomp{\mggraph{\grind g}} = \mgbartropr{\grind{g}}{0,1}$

\begin{center}
\begin{figure}[h!] 
\begin{tikzpicture}[scale=0.8, every node/.style={scale=0.9}]

\foreach \x in {1,..., 2}{
	\node at (2*\x, 0) {$\mgbartropr{\grind{g}}{\x, 0}$};
	\draw[->>]  (2*\x, 1.25) -- (2*\x, 0.75) ;
	\draw[->>]  (2*\x + 1.25, 0) -- (2*\x + 0.75, 0) ;
}

\foreach \x in {0,..., 2}{
	\node at (2*\x, 2) {$\mgbartropr{\grind{g}}{\x, 1}$};
	\draw[->>]  (2*\x, 3.25) -- (2*\x, 2.75) ;
	\draw[->>]  (2*\x + 1.25, 2) -- (2*\x + 0.75, 2) ;
}

\foreach \x in {0,..., 1}{
	\node at (2*\x, 4) {$\mgbartropr{\grind{g}}{\x, 2}$};
	\draw[->>]  (2*\x + 1.25, 4) -- (2*\x + 0.75, 4) ;
		\draw[->>]  (2*\x, 5.25) -- (2*\x, 4.75) ;
}

	\node at (0, 9) {$\mgbartropr{\grind{g}}{0, N+1}$};
	\node[right] at (0.75, 9) {$= \mgtrop{\grind{g}}$};

	\node at (0, 7) {$ \mgbartropr{\grind{g}}{0, N}$};
	\draw[->>]  (0, 9 - 0.75) -- (0, 9-1.25) ;
	
	\node at (2, 7) {$\mgbartropr{\grind{g}}{1, d} $};
	\node[right] at (2.75, 7) {$= \mgtrop{\grind{g}}$};

	\draw[->>]  (1.25, 7) -- (0.75, 7) ;

	\draw[->>, dotted]  (2, 7 - 0.75) -- (2, 7-2.25) ;
	\draw[->>, dotted]  (0, 7 - 0.75) -- (0, 7-2.25) ;
	\node at (9, 0) {$\mgbartropr{\grind{g}}{N+1, 0}$};
	\node[right] at (9.75, 0) {$= \mgtrop{\grind{g}}$};

	\node at (7, 0) {$ \mgbartropr{\grind{g}}{N, 0}$};
	\draw[->>]  (9 - 0.75, 0) -- (9-1.25, 0) ;
	
	\node at (7, 2) {$\mgbartropr{\grind{g}}{N, 1} $};
	\node[right] at (7.75, 2) {$= \mgtrop{\grind{g}}$};

	\draw[->>]  (7, 1.25) -- (7, 0.75 ) ;

	\draw[->>, dotted]  (7 - 0.75, 2) -- (7-2.25, 2) ;
	\draw[->>, dotted]  (7 - 0.75, 0) -- (7-2.25, 0) ;

\end{tikzpicture}
\caption{The spaces $\mgbartropr{\grind{g}}{p, q}$ and the corresponding forgetful maps.}\label{fig:CompactificationTower-intro}
\end{figure}
\end{center}

The double tower of tropical spaces given in Figure~\ref{fig:CompactificationTower-intro} interpolates between $\mgtrop{\grind{g}}$ and $\cancomp{\mggraph{\grind g}}$. Moreover, it should be compared with the tower of hybrid spaces
\begin{align}\label{eq:tower-moduli-intro}
\mgbarg{\grind g} \longleftarrow \mgbarhybr{\grind g}{1} \longleftarrow \mgbarhybr{\grind g}{2} \longleftarrow \dots \longleftarrow \mgbarhybr{\grind g}{N}= \mghyb{\grind{g}},
\end{align}
for $N = 3g-3$ which we introduced in~\cite{AN}. The second last line in Figure~\ref{fig:CompactificationTower-intro} is precisely the tropical analogue of the tower \eqref{eq:tower-moduli-intro}.

\medskip

Similarly to $\mgtrop{\grind{g}}$, we can define the tropical moduli space $\mgtrop{\grind{g,n}}$, that is, the \emph{moduli space of tropical curves of genus $g$ with $n$ marked points}. The latter coincides with the absolute canonical compacticiation of  $\mggraph{\grind{g,n}}$ in the family of higher rank canonical compactifications $\mgbartropr{\grind{g,n}}{p,q}$, $p+q \leq N =3g-3+n =\dim(\mggraph{\grind{g,n}})$.

\medskip

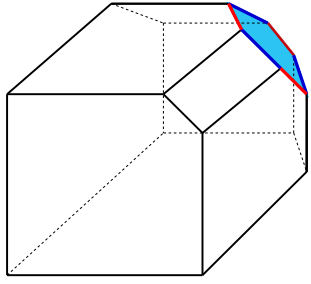
\begin{figure}[!t]
\centering
    \scalebox{.27}{\input{example5.tikz}}
\caption{Canonical compactification of a  three-dimensional  simplicial cone. The part in color is the locus of points of full sedentarity.}
\label{fig:compactified_cone}
\end{figure}

The moduli space of tropical curves of combinatorial type $G$ defined in our previous work~\cite{AN} appears naturally inside $\mgtrop{\grind{g}}$, namely, as the part consisting of all the tropical curves  of combinatorial type $G$ with full sedentarity, that is, with underlying ordered partitions of full sedentarity (see Figures~\ref{fig:compactified_cone} and \ref{fig:full_sedentarity} for an example).

\medskip

The moduli spaces of tropical curves $\mgtrop{\grind{g}}$ and $\mgtrop{\grind{g,n}}$ come with the corresponding universal tropical curves $\tropcurve_{\grind{g}}$ and $\tropcurve_{\grind{g,n}}$, respectively (defined over suitable \'etale charts $\mgtropcombin{\grind{G}}$).

\subsection{Canonical measures on tropical curves} Coming back to Question~\ref{question:main-graphs},  the second task consists in providing a description of all the possible limits of canonical measures on metric graphs. 

\smallskip

Consider thus a tropical curve $\curve$ defined as the metric realization of the pair $(G, l)$, of a graph $G$ and an edge length function $l$, equipped with an ordered partition $\pi = (\pi_1, \dots, \pi_r, \pi_\fin)$ on the edges of the graph $G$.

We define the \emph{canonical measure} $\mu^\can$ on the tropical curve $\curve$ of genus $g$ as the sum 
\[\mu^\can := \frac{1}{g} \Big( \mu^1_\Zh + \dots + \mu^r_\Zh + \mu^\fin_\Zh  \Big)\]
where the measure $\mu^j_\Zh$,  \emph{the $j^{\it th}$ graded  part of the canonical measure}, $j=1, \dots, r, \fin$, has support in the union of the intervals $\Ical_e$ in $\curve$ for edges $e\in \pi_j$, and on each such interval $\Ical_e$, it coincides with the canonical Zhang measure of the $j$-th graded minor $\Gamma^j$. In case that $\curve$ comes with a genus function $\genusfunction \colon V \to \N \cup\{0\}$, we additionally add point masses at vertices.

\medskip

The relevance of these constructions is described in the following theorem.

\begin{thm}[Continuity of the canonical measure on the moduli space of tropical curves] \label{thm:tropical-canonical-measures-intro} The universal family of canonically measured tropical curves $(\tropcurve_{\grind{g}}, \mu^\can)$ forms a continuous family of measured spaces over the tropical moduli space $\mgtrop{\grind{g}}$.
\end{thm}
As in~\cite{AN}, continuity in the above statement is defined in a weak sense. That is, for every continuous function $f :\tropcurve_{\grind{g}} \to \R$, the function $F : \mgtrop{\grind{g}}\to \R$ defined by integration along fibers
\begin{align*}
&F(\thy) := \int_{\tropcurve_{\thy}} f_{|_{\tropcurve_{\thy}}} \, d\mu^\can_\thy, \qquad  \thy \in \mgtrop{\grind{g}},
\end{align*}
is continuous on $\mgtrop{\grind{g}}$. Here $\mu_\thy^\can$ refers to the canonical measure on the tropical curve $\tropcurve_\thy$ which is the fiber in the universal family $\tropcurve_{\grind{g}}$ over the point $\thy\in \mgtrop{\grind{g}}$.

\medskip

The above theorem is thus the tropical analogue of the main result in~\cite{AN}. It follows from this theorem that, if a sequence $s_1, s_2, \dots$ in $\mggraph{\grind g}$, representing metric graphs $\mgr_1, \mgr_2, \dots$, converges  to a point $\thy$ in $\mgtrop{\grind{g}}$, representing the tropical curve $\tropcurve_\thy$, then the canonical measure $\mu_j^\can$ of $\mgr_j$ converges  to the canonical measure of $\tropcurve_\thy$.
 The converse is also true, and follows from the fact that the tropical moduli space $\mgtrop{\grind{g}}$ is a compact Hausdorff space.

\subsection{Formulation of the question on variations of Green functions}
Taking into account the above, it is natural to study the variation of the canonical Green function  in the tropical moduli space $\mgtrop{\grind{g}}$. We replace Question~\ref{question:main-graphs} by the following problem:
\begin{question}\label{question:main-graphs-precise}
What is the shape of $\gri{s}$, $s \in \mggraph{\grind g}$, for $s$ close to a fixed point $\thy$ in the tropical moduli space $\mgtrop{\grind{g}}$?
\end{question}

Before answering Question~\ref{question:main-graphs-precise}, we first need to introduce several important concepts. In what follows, we develop a \emph{function theory on tropical curves}, define a \emph{tropical Laplace operator}, and introduce \emph{tropical Poisson equations and Green functions}. The latter leads to a precise description of the canonical Green function close to the boundary of the tropical moduli space, giving a complete answer to Question~\ref{question:main-graphs-precise}. Moreover, the framework allows to consider Green functions associated to other natural measures on metric graphs.

\subsection{Function theory on tropical curves and the tropical Laplacian} \label{ss:IntroFunctionTheoryTropicalCurves}
Let $\curve$ be a tropical curve of rank $r$ with underlying model $(G, l)$ and with ordered partition $\pi=(\pi_\infty, \pi_\fin)$, $\pi_\infty = (\pi_1, \dots, \pi_r)$ of rank $r$.  For $j=1, \dots, r, \fin$, denote by $\Gamma^j$ the $j$-th graded minor of $\curve$.

\subsubsection{Functions} A \emph{real valued function} $\lf$ on the tropical curve $\curve$ of rank $r$ is an $r+1$-tuple $\lf = (f_1, \dots, f_r, f_\fin)$ consisting of functions $f_j \colon \Gamma^j \to \R$  for each $j=1, \dots, r, \fin$; $\lf$ is called continuous, piecewise smooth, etc. if the corresponding components $f_j$ on $\Gamma^j$ verify the same property, i.e., are continuous, piecewise smooth, etc., respectively.

\subsubsection{Measures} A \emph{layered measure} $\lmu = \lmu_\pi$ on the tropical curve $\curve$ is an $r+1$-tuple 
\[\lmu=(\mu_1, \dots, \mu_r, \mu_\fin)\] consisting of measures $\mu_j$ on the graded minors $\Gamma^j$, for $j=1, \dots, r, \fin$.  The set of all layered measures on $X$ is denoted by $\mathcal{M}(\curve)$. A layered measure comes with the corresponding \emph{mass function}, with domain of definition the connected components of the graded minors of $\curve$. A layered measure $\lmu$ on $\curve$ has mass zero if the function $\mass$ is the zero function. The set of layered measures of mass zero on $X$ is denoted by $\mathcal{M}^0(\curve)$.

\subsubsection{Tropical Laplacian} We define the \emph{tropical Laplacian} $\Deltatrop=\Deltatrop_\curve$ on $\curve$ to be a measure-valued operator which maps piecewise smooth functions $\lf =(f_1, \dots, f_r, f_\fin)$ on $\curve$ to layered measures of mass zero on $\curve$. Informally speaking, $\Deltatrop$ takes the Laplacian on each graded minor and additionally transfers the function slopes from \emph{larger} infinities to \emph{smaller} infinitary and finitary layers , that is, from an infinity layer of lower index to all layers of higher index. Formally, it is defined as the sum
\begin{equation} \label{eq:TropicalLaplacian-intro}
\Deltatrop (\lf) = \Deltatrop_1(f_1) + \Deltatrop_2(f_2) + \dots \Deltatrop_r(f_r) + \Deltatrop_\fin(f_\fin)
\end{equation}
where for each $j=1, \dots, r, \fin$, and for any function $f_j$ on the graded minor $\Gamma^j$, the layered measure $\Deltatrop_j(f_j)$ on the tropical curve $\curve$ is given by
\[
\Deltatrop_j(f_j)  := \Bigl(0, \dots, 0, \Delta_j (f_j),  \divind{j}{j+1}(f_j), \dots,  \divind{j}{r}(f_j),\divind{j}{\fin} (f_j) \Bigr).
\]
Here, $\Delta_j (f_j) \in \mathcal{M}^0(\Gamma^j)$ denotes the Laplacian on the $j$-th graded minor $\Gamma^j$. Moreover, for each $i > j$, the point measure $\divind{j}{i} (f_j)$ is defined by first calculating the sum of incoming slopes of $f_j$ around vertices in the graph $G$, and then pushing out the resulting point measure (with support in $V$) to the $i$-th graded minor (see Section~\ref{sec:tropical_laplacian} for details).

\subsubsection{Tropical Laplacian as a weak limit of metric graph Laplacians}

We will establish the following theorem (see Section~\ref{sec:tropical_laplacian}).
\begin{thm}\label{thm:GraphLaplacianConvergence-intro}
The tropical Laplacian $\Deltatrop$ is the weak limit of metric graph Laplacians. 
\end{thm}
Here, we briefly explain the meaning of this convergence and refer to Section~\ref{sec:tropical_laplacian} for a more precise statement.

 We fix a graph $G=(V, E)$ and consider the space $\mggraphcombin{\grind{G}}$ of metric graphs of combinatorial type $G$, that is the cone $\mggraphcombin{\grind{G}} = \R_{>0}^E$ of possible metrics on the edges. The canonical compactification $\mgtropcombin{\grind{G}}$ coincides with the space of tropical curves of combinatorial type $G$. Each tropical curve $\curve_\thy$ for $\thy \in \mgtropcombin{\grind{G}}$ gives rise to a \emph{pull-back operator} which transforms functions $\lf$ on $\curve_\thy$ into functions $f_t$, $t \in \mggraphcombin{\grind{G}}$, on each metric graph $\mgr_t$ for $t\in \mggraphcombin{\grind{G}}$. This is described in Section~\ref{sss:PullbackExplanation}.
 
 Broadly speaking, we can identify the first piece $f_1$ of $\lf = (f_1, \dots, f_r, f_\fin)$ with a function $f_1$ on the space $\Gamma$, which is constant on all edges of the higher index layers $\pi_i$, $i > 1$. More generally, the $j$th piece $f_j$ of $\lf$ can be identified with a function on $\Gamma$, which is constant on edges of higher index layers, equal to the $f_j$ on the $j$th layer $f_j$, and linear on all edges of lower index layers. Using homotecies between edges, we can map the metric graph $\mgr_t$ to $\Gamma$ and pull-back these functions to $r+1$ functions $f_1^\ast, \dots, f_r^\ast, f_\fin^\ast$  on $\mgr_t$. The total \emph{pull-back of the function} $\lf_\thy$ is given by
  \[
\lf_t^\ast := \sum_{j=1}^r L_j(t) f_j^\ast + f_\fin^\ast
 \]
 for $L_j(t) := \sum_{e \in \pi_j} t(e)$. Note in particular that $\lf_t^\ast$ is unbounded when $t$ converges to $\thy$, since some edge lengths are going to infinity. For any function $\lf$ on $\curve_\thy$, the following convergence of measures holds in the weak sense
\[
\Delta(\lf_t^\ast)  \to \Deltatrop(\lf),
\]
as the point $t$ converges to $\thy$ in $\mggraphcombin{\grind{G}}$. This shows $\Deltatrop$ can be interpreted as a weak measure-theoretic limit of metric graph Laplacians.

\subsection{Tropical Poisson equation} With a notion of Laplacian at hand, we can formulate the Poisson equation on a tropical curve. 
 Suppose that $\curve $ is a tropical curve with underlying graph $G$, edge length function $l$, and the ordered partition $\pi$ of rank $r$.   We suppose the edge length function $l \colon E(G) \to (0, + \infty)$ is normalized on each layer $\pi_j$, for $j\in [r]$, and associate the metric graph $\Gamma$ in the equivalence class of the tropical curve $\curve$.

   \smallskip

 Let further $\lmu$ and $\lnu$ be two layered Borel measures on $\curve$ such that
 \begin{itemize}
 \item the measure $\lmu$ is of mass zero, that is $\lmu \in \mathcal{M}^0(\curve)$.
 
 \smallskip
 
 \item the mass function $\mass_{\lnu}$ of the measure $\lnu$ has mass one on the first graded minor $\Gamma^1$, and has mass zero on all components of all higher indexed graded minors. 
  \end{itemize}
On such a  {\em bimeasured tropical curve} $(\curve, \lmu, \lnu)$, consider the following {\em tropical Poisson equation}
 \begin{equation}  \label{eq:PoissonTropicalCurve-intro}
 \begin{cases}
 \Deltatrop \lf = \lmu \\
\int_{\curve} \lf \, d \lnu = 0,
 \end{cases}
\end{equation}
 where the integration condition means that
\[
\int_\Gamma f_j \, d\nu = 0 \qquad \text{for all }j=1, \dots, r, \fin.
\]
Here $\nu$ is the measure on the metric graph $\Gamma$ obtained from the layered measure $\lnu$ (which has total mass $\nu(\Gamma) = 1$). In the integration condition, the function $f_j$ is extended from $\Gamma^j$ to the full metric graph $\Gamma$ by linear interpolation. 
\begin{thm}[Existence and uniqueness of solutions of the tropical Poisson equations] Notations as above, the tropical Poisson equation \eqref{eq:PoissonTropicalCurve-intro} has a solution $\lf$. The solution is moreover unique if we require $\lf$ to be harmonically arranged. 
\end{thm}

In the above theorem, harmonically arranged functions $\lf$ are defined as follows. A function $f_j \colon \Gamma^j \to \C$ on the $j$-th graded minor $\Gamma^j$ is called \emph{lower harmonic}, if, after extending $f_j$ linearly to the full metric graph  $\Gamma$, the differential $d f_j $ on $\Gamma$ is harmonic on all lower index layers. That is, for all $i < j$, the slopes of $f_j$ around each vertex in $\Gamma^i$ sum up to zero. A function $\lf =(f_j)_j$ on the tropical curve $\curve$ is called \emph{harmonically arranged} if all $f_j$'s are lower harmonic.

\smallskip

The property of being harmonically arranged is crucial in our considerations and one of the main discoveries in this work. It allows to capture the asymptotic behavior of solutions to the Poisson equation on degenerating metric graphs (and later on Riemann surfaces), and can be motivated from the perspective of higher rank valuation theory.

\smallskip

An important result about this concept is the \emph{Rearrangement Theorem}~\ref{thm:HarmonicExtension}, which says every function  on $\curve^\trop$ has a \emph{harmonic rearrangement} unique up to a global additive constant (vector) on the tropical curve. We refer to Section~\ref{sec:harmonic_rearrangement} for more details.

\subsection{Tropical log maps and tameness of solutions to Poisson equations} We are interested in describing the behavior of the solution of tropical Poisson equation close to the tropical limit, that is, under degenerations of metric graphs. In order to give a mathematical meaning to this question, we introduce two notions of tameness for functions on families of tropical curves: the {\em stratumwise tameness} and {\em weak tameness}. 

\smallskip

The idea behind these notions is to give a mathematical meaning to the convergence of functions $f$ on metric graphs $\mgr$ to functions $\lf$ on higher rank tropical curves $\curve$. Since in contrast to a usual function $f \colon \mgr \to \C$, a tropical function $\lf =(f_1, \dots, f_r, f_\fin)$ has several components, classical continuity notions obviously do not apply. Tameness of a family of functions means that, when a metric graph $\mgr$ degenerates to a tropical curve $\curve$, the respective function $f_\mgr$ admits an \emph{asymptotic expansion in terms of the components of tropical curve functions}. We next make this precise.

\subsubsection{Tropical log maps} The definitions are both based on the use of appropriate log maps on the corresponding moduli spaces.  
 
 The  problem of constructing log maps on families plays a central role in the present paper both in the tropical and in the hybrid settings. 
In the former case, these are maps which allow to compare the geometry of tropical curves, appearing on the boundary of moduli spaces, to the geometry of metric graphs. The \emph{tropical log map} is a retraction map 
\[\logtrop\colon \mgtrop{\grind{g,n}} \setminus \umggraph{\grind{g,n}} \to \partial_\infty \mgtrop{\grind{g,n}}, \]
defined on a neighborhood of the boundary of the moduli space of tropical curves (of genus $g$ with $n$ marked points), where $\umggraph{g,n}$ is the part of the moduli space of metric graphs with edge lengths all bounded by one. It is defined by using a foliation of the metric graph moduli space in terms of \emph{generalized permutohedra},  a special class of polytopes related to ordered partitions. We refer to  Section~\ref{sec:tropical_moduli} for the details of the construction. 

\smallskip

Let $G = (V,E)$ be a fixed (augmented) graph (with marking) and consider the space $\mgtropcombin{\grind{G}}$ of tropical curves of combinatorial type $G$. Let $\unicurvetrop{\grind{G}}$ be the corresponding universal family of tropical curves. Denote by $\mggraphcombin{\grind{G}} \subseteq \mgtropcombin{\grind{G}}$ the space of metric graphs of combinatorial type $G$. The above framework leads to a tropical log map
\[
\logtrop \colon \mggraphcombin{\grind{G}} \to \partial_\infty \mgtropcombin{\grind{G}} = \mgtropcombin{\grind{G}}  \setminus \mggraphcombin{\grind{G}},
\]
which sends any point $t \in  \mggraphcombin{\grind{G}} $ representing a metric graph to a point $\thy \in \partial_\infty \mgtropcombin{\grind{G}}$ representing a tropical curve. Using  simply homotecies of edges, we can also map the corresponding geometric objects into each other. Altogether, we get a map from the family $\unicurvetrop{\grind{G}} /\mgtropcombin{\grind{G}} $ to $\unicurvetrop{\grind{G}} / \partial_\infty \mgtropcombin{\grind{G}} $, arriving at a commutative diagram
\[
\begin{tikzcd}
\unicurvetrop{\grind{G}} / {\mggraphcombin{\grind{G}}} \arrow[d]\arrow[r, ""] & \arrow[d] \unicurvetrop{\grind{G}} / \partial_\infty \mgtropcombin{\grind{G}} \\
\mggraphcombin{\grind{G}} \arrow[r, "\logtropdiag"] &\partial_\infty \mgtropcombin{\grind{G}}
\end{tikzcd}
\]
Note that the boundary $\partial_\infty \mgtropcombin{\grind{G}}$ has a natural stratification
 \[
\partial_\infty \mgtropcombin{\grind{G}} = \bigsqcup_{\pi} \mgtropcombin{\combind{(G, \pi)}},
 \]
 where $\mgtropcombin{\combind{(G, \pi)}}$ is the space of tropical curves of combinatorial type $(G, \pi)$ for an ordered partition $\pi =(\pi_\infty, \pi_\fin)$ of $E$ with $\pi_\infty \neq \varnothing$. For each boundary stratum $\mgtropcombin{\combind{(G, \pi)}}$, there is a natural projection map $\pr_\pi \colon  \mggraphcombin{\grind{G}} \to \mgtropcombin{\combind{(G, \pi)}}$. The latter sends each point $t$ representing a metric graph $\mgr_t$ of combinatorial type $G$ to the point $\pr_\pi(t)$ representing the conformal equivalence class of the layered metric graph $(\mgr_t, \pi)$.  Using again homotecies of edges,  this leads to a map from the family $\unicurvetrop{\grind{G}} /{\mggraphcombin{\grind{G}}}$ to $\unicurvetrop{\grind{G}} / \mgtropcombin{\combind{(G, \pi)}}$ and a commutative diagram
 
\[
\begin{tikzcd}
\unicurvetrop{\grind{G}} / {\mggraphcombin{\grind{G}}} \arrow[d]\arrow[r, ""] & \arrow[d] \unicurvetrop{\grind{G}} / {\mgtropcombin{\combind{(G,\pi)}}  } \\
\mggraphcombin{\grind{G}} \arrow[r, "\pr_\pi"] &\mgtropcombin{\combind{(G, \pi)}}
\end{tikzcd}
\]

\subsubsection{Stratumwise tameness} 

Let $\lf = (\lf_\thy)_\thy$ be a {\em family of tropical functions} on $\unicurvetrop{\grind{G}}$. That is, for each $\thy \in \mgtropcombin{\grind{G}}$, we have a function $\lf_\thy$ on the tropical curve represented by $\thy$.

\medskip

We call the family $\lf = (\lf_\thy)_\thy$ {\em tame at the stratum} $\mgtropcombin{\combind{(G, \pi)}}$ if for any point $\thy$ belonging to the \emph{closure of $\mgtropcombin{\combind{(G, \pi)}}$} in $\mgtropcombin{\grind{G}}$, we have
\begin{equation} \label{eq:StratumwiseTameness-intro}
	\lim_{t \to \thy} \sup_{x \in \unicurvetrop{t}} \big |  \lf_t  (x)  - \lf_{\pr_\pi(t)}^\ast  (x) \big | = 0
\end{equation}
as $t \in  \mggraphcombin{\grind{G}}$ converges to $\thy$. In the above statement, $ \lf_{\pr_\pi(t)}^\ast  $ denotes the pull-back of the function $\lf_{\pr_\pi(t)} $ on the hybrid curve $\unicurvetrop{\pr_\pi(t)}$ to the metric graph $\unicurvetrop{t}$,
\[
 \lf_{\pr_\pi(t)}^\ast = L_1(t) f_{\pr_\pi(t), 1}^\ast + \dots + L_r(t) f_{\pr_\pi(t), r}^\ast + f_{\pr_\pi(t), \fin}^\ast
\]
for $L_i(t) = \sum_{e \in \pi_i} t(e)$, $i \in [r]$. Note that the pull-back becomes unbounded when $t$ converges to $\thy$, since some edge lengths are going to infinity and hence $L_j(t) \to \infty$ for $j \in [r]$.

\smallskip

A family of functions $\lf_\thy$, $\thy \in \mgtropcombin{\grind{G}}$, is called {\em stratumwise tame} if it is tame at all boundary strata $\mgtropcombin{\combind{(G, \pi)}}$.

\subsubsection{Weak tameness} \label{ss:TameFunctions-intro}  Let $\lf = (\lf_\thy)_\thy$ be a family of tropical functions on $\unicurvetrop{\grind{G}}$, in the sense of the previous section.

\smallskip

 The boundary $\partial_\infty \mgtropcombin{\grind{G}} = \mggraphcombin{\grind{G}}  \setminus \mggraphcombin{\grind{G}}$ is stratified as $\partial_\infty \mgtropcombin{\grind{G}}  = \bigsqcup_\pi  \mgtropcombin{\combind{(G, \pi)}}$, which induces the following decomposition
\[
\mgtropcombin{\grind{G}} = \bigsqcup_\pi \inn R_\pi, \qquad \qquad \inn R_\pi := \logtrop^{-1}\Big ( \mgtropcombin{\combind{(G, \pi)}} \Big ).
\]
Here the union is taken over all ordered partitions $\pi$ with non-empty sedentarity of $E$.

\smallskip

Assume further that $\layh_\thy$, $\thy \in  \mgtropcombin{\grind{G}}$, are tropical functions satisfying the following properties:
\begin{itemize}
\item[(i)] for each $\thy \in \mgtropcombin{\grind{G}}$, the function $\layh_\thy$ is a tropical function \emph{on the tropical curve} $\unicurvetrop{\shy}$ for $\shy = \logtrop(\thy)$.
\item[(ii)] The equality  $\layh_{\thy} = \lf_\thy$ holds for all $\thy \in \partial_\infty \mgtropcombin{\grind{G}}  $.
\item[(iii)] Consider the boundary stratum $\mgtropcombin{\combind{(G, \pi)}}$ for an ordered partition $\pi =(\pi_1, \dots, \pi_r, \pi_\fin)$ on $E$. By property (i), the tropical function $\layh_\thy$ is of the form $\layh_\thy = (h_{\thy, 1}, \dots, h_{\thy, r}, h_{\thy, \fin})$  for each base point $\thy$ in the region $\inn R_\pi$.  Then all components $h_{\thy, j}$ , $j=1, \dots, r, \fin$,  depend continuously on the parameter $ \thy \in \inn R_\pi$.
\end{itemize}

For any point $t \in \mggraphcombin{\grind{G}}$, denote by $\layh^\ast_t$ the pull-back
\[
\layh^\ast_t = L_1(t) h_{\thy, 1}^\ast + \dots + L_{r_\thy}(t) h_{\shy, r_\thy}^\ast + h_{\thy, \fin}^\ast
\]
of $\layh_t$ from the tropical curve $\unicurvetrop{\thy}$, $\thy = \logtrop(t)$, to the metric graph $\unicurvetrop{t}$.

\begin{defi}[Weak tameness] \label{def:WeakGlobalTameTropical}
A family of tropical functions $\lf_\thy$, $\thy \in \mgtropcombin{\grind{G}}$, is {\em weakly tame} if there are tropical functions $\layh_\thy$, $\thy \in \mgtropcombin{\grind{G}}$, as above such that, in addition, for any point $\thy \in \partial_\infty\mgtropcombin{\grind{G}}$, we have
\begin{equation} \label{eq:Differencev}
 \lim_{t \to \thy} \sup_{x \in \unicurvetrop{t}} |\lf_t (x) -  \layh^\ast_t (x)| = 0
 \end{equation}
if the point $t \in \mggraphcombin{\grind{G}}$ converges to $\thy$.
\end{defi}

For a more precise definition, we refer to Section~\ref{ss:TameFunctions}.

\smallskip

Using the properties of the tropical logarithm map, it can be shown that every continuous, stratumwise tame family of functions is also weakly tame.

\subsubsection{Tameness of solutions of the tropical Poisson equation}  Our main theorems on the tropical Poisson equation can be now stated as follows.

\begin{thm}[Stratumwise tameness of solutions] \label{thm:main_graph_strong-intro}
Let $(\lmu_\thy)_\thy$ and $(\lnu_\thy)_\thy$ be continuous families of measures on $\unicurvetrop{\grind{G}} / \mgtropcombin{\grind{G}}$, making each fiber a bimeasured tropical curve $(\unicurvetrop{\thy}, \lmu_\thy, \lnu_\thy)$. Suppose that, for every boundary stratum $\mgtrop{\combind{(G, \pi)}}$, we have that
\begin{equation} \label{eq:condition_strong_tame-intro}
L(t) \big (\lmu_t - \lmu_{\pr_\pi(t)}^\ast\big) \qquad \text{and} \qquad L(t) \big (\lnu_t - \lnu_{\pr_\pi(t)}^\ast\big)
\end{equation}
converge weakly to zero provided that the point $t \in \mggraphcombin{\grind{G}}$ tends to a point $\thy$ in the closure of $\mgtropcombin{\combind{(G, \pi)}}$.

\smallskip

Then the solutions $\lf_\thy$, ${\thy \in \mgtropcombin{\grind{G}}}$, to the tropical Poisson equation form a stratumwise tame family of functions.
\end{thm}

\medskip

In the above, $L(t)$ denotes the total length of the fiber at $t$ in the universal curve $\unicurvetrop{\grind{G}}$, that is the associated metric graph  $\mgr_t$, and $\lmu_{\pr_\pi(t)}^\ast$ and $\lnu_{\pr_\pi(t)}^\ast$ denote the natural pull-backs of the measures $\lmu_{\pr_\pi(t)}$ and $\lnu_{\pr_\pi(t)}$ to this metric graph. Moreover,  the weak convergence to zero in \eqref{eq:condition_strong_tame-intro} is understood in the distributional sense, that is we require
\[
\lim_{t \to \thy} L(t) \int_{\mgr_t} f(q) \, d(\lmu_t - \lmu_{\pr_\pi (t)}^\ast)(q)  = 0
\]
for all continuous functions $f \colon \unicurvetrop{\grind{G}} \to \R$, and similar for the measures $\nu_\thy$.

\smallskip

We emphasize that the stratumwise tameness of the solutions implies that, if a sequence of metric graphs $(\mgr_n)_n$ converges to a tropical curve $\curve= (G, \pi, l)$, then
\[
	\lim_{n \to \infty} \sup_{x \in \mgr_n} \big |  \lf_{t_n}  (x)  - \lf_{\pr_\pi(t_n)}^\ast  (x) \big | = 0,
\]
for the pullback $ \lf_{\pr_\pi (t_n)}^\ast$ obtained by using the projection map $\pr_\pi \colon \mggraphcombin{\grind{G}} \to \mgtropcombin{\combind{(G, \pi)}}$.
\smallskip

Hence, the above theorem provides a full description of the asymptotic of the solutions.

\medskip

Applied to discrete measures, the above result reads as follows.

\begin{thm}[Stratumwise tameness for discrete measures] \label{thm:MainGraphDiscrete-intro}
Let $D \in \Div^0(G)$ be a degree zero divisor and $v \in V$ a fixed vertex. Denote by $v_j$ the image of $v$ in the graded minor $\Gamma^j$ under contraction. For any point $\thy \in \mgtropcombin{\grind{G}}$, let $\lf_\thy$ be unique solution to the equations
\[
\begin{cases}
 \Deltatrop (\lf) = D \\
 \lf \text{ is harmonically arranged} \\
 f_j (v_j) = 0, \qquad j=1, \dots, r, \fin,
 \end{cases}
\]
on its associated tropical curve  $\unicurvetrop{\thy}$. 

Then the solutions $\lf_\thy$, ${\thy \in \mgtropcombin{\grind{G}}}$, form a stratumwise tame family of functions.
\end{thm}

We will complement this statement in Theorem~\ref{thm:layered_expansion_tropical_laplacian}, where we give an expansion of the solutions around each point $\curve$ of the tropical moduli space $\mgtropcombin{\grind{G}}$. 
\medskip

In the above results, we assume that the measures are converging rather fast, formalized in \eqref{eq:condition_strong_tame-intro}. This condition is always satisfied in the situation of Theorem~\ref{thm:MainGraphDiscrete-intro}. Dropping the condition \eqref{eq:condition_strong_tame-intro}, we get the following general result.

\begin{thm}[Global weak tameness of solutions] \label{thm:MainGraphWeak-intro}
Let $(\mu_\thy)_\thy$ and $(\nu_\thy)_\thy$ be continuous families of measures on $\unicurvetrop{\grind{G}}/\mgtropcombin{\grind{G}}$, making each $\unicurvetrop{\thy}$ a bimeasured tropical curve. Then the solutions $\lf_\thy$ to the tropical Poisson equation form a weakly tame family of functions on $\mgtropcombin{\grind{G}}$.
\end{thm}

We note that stratumwise tameness might fail if the convergence of the metric graph measures $\mu_t$ and $\nu_t$ to tropical curve measures $\lmu_\thy$ and $\lnu_\thy$ is not fast enough. Since this convergence can be arbitrarily slow, the weak tameness result appears to be the most general statement one can hope for.  Moreover, by choosing $\mu_t$ and $\nu_t$ as discrete measures with masses within the edges of the first layer $\pi_1$, one can show that the condition for stratumwise tameness in \eqref{eq:condition_strong_tame-intro} is in a certain sense sharp. 
\smallskip

 On the other hand, even in the situation where the convergence assumption \eqref{eq:condition_strong_tame-intro} fails, it is sometimes possible to find the full asymptotic of solutions by using higher order terms to approximate the measures. We will demonstrate this approach by the examples of Green functions associated to the Lebesgue and the canonical measures (see Theorem~\ref{thm:Lebesgue_Green_functions}  and Proposition~\ref{prop:simple_layering_Green_functions}).

\subsection{Tropical Green functions}
Let $\curve$ be a tropical curve of rank $r$. Let $\lmu$ be a measure of total mass one on $\curve$. That is the mass function $\mass_{\lmu}$ takes value one on the first minor, and it takes value zero on all the components of any other graded minor. 

\smallskip

Denote by $\Gamma$ the layered metric graph corresponding to $\curve$ with normalized edge length functions in infinitary layers.  Any point $p \in \Gamma$ induces a layered Dirac measure $\bm{\delta}_p$ on $\curve$. In particular, the following Poisson equation
\begin{equation} \label{eq:TropicalGreenFunction-intro}
\begin{cases}
 \Deltatrop \lf = \bm{\delta}_p - \lmu, \\[1 mm]
 \lf \text{ is harmonically arranged}, \\[1mm]
 \int_\curve \lf \, d \lmu = 0
 \end{cases}
\end{equation}
has a unique solution $\lf = (f_1, \dots, f_r, f_\fin)$ on $\curve$ for every $p \in \Gamma$.

\smallskip

We denote the unique solution of the equation \eqref{eq:TropicalGreenFunction-intro} by
\[
\lgri{\lmu} (p, \cdot) = \big ( \gri{\lmu, 1}(p, \cdot), \dots  ,\gri{\lmu, r}(p, \cdot), \gri{\lmu, \fin}(p, \cdot) \big).
\]
By linear interpolation, each component $\gri{\lmu, j}(p, \cdot)$ can be viewed as a function on the full metric graph $\Gamma$. Therefore we get a {\em tropical Green function} 
\begin{equation} \label{eq:IntroTropicalGreenFunction}
\begin{array}{cccc}
\gri{\curve, \mu} \colon &\Gamma \times \Gamma  &\longrightarrow & \R^{r+1} = \R^{[r]}\times \R^\fin
\end{array}
\end{equation}
where $\mu$ is the measure on $\Gamma$ obtained from $\lmu$. The tropical Green function of the canonical measure $\lmu = \lmu^\can$ is simply denoted by $\gri{\curve}$.

\subsection{Main theorem on canonical Green functions in the tropical setting}
Our main theorem on variations of canonical Green functions on metric graphs reads as follows. It provides an asymptotic expansion for the Green function on degenerating metric graphs.

\begin{thm}[Expansion of the canonical Green function] Fix a combinatorial graph $G= (V,E)$ and consider the space $\mggraphcombin{\grind{G}}$ of metric graphs $\mgr$ over $G$. Let further $\mgtropcombin{\combind{(G, \pi)}}$ be a fixed boundary stratum in the space $\mgtropcombin{\grind{G}}$ of tropical curves over $G$.

Then there exists a family of functions $\grihat{\mgr,j} \colon \mgr \times \mgr \to \R$, for $j=1, \dots, r, \fin$ (where $r$ is the rank of the partition $\pi$) and $\mgr \in \mggraphcombin{\grind{G}}$, with the following properties.

\smallskip

\begin{itemize}
\item[$(1)$] Denoting $L_i(\mgr) = \sum_{e \in \pi} l_\mgr(e)$, $i \in [r]$, the difference 
\[
\gri{\mgr} - L_1(\mgr)\grihat{\mgr,1} - L_2(\mgr)\grihat{\mgr,2} - \dots - L_r(\mgr)\grihat{\mgr,r}- \grihat{\mgr, \smallcc}\]
 extends uniformly by zero over the boundary stratum $\mgtropcombin{\combind{(G, \pi)}}$.
 
 That is, for a sequence of metric graphs $\mgr$ which converges to a tropical curve $\curve \in \mgtropcombin{\combind{(G, \pi)}}$ in $\mgtropcombin{\grind{G}}$, the above difference tends to zero uniformly on $\mgr \times \mgr$.

\smallskip

\item[$(2)$]
Moreover, the $j$-th function $\grihat{\mgr,j}$ converges uniformly to the $j$-th part of the canonical Green function on $\curve$,
 \[\lim_{\mgr \to \curve} \grihat{\mgr,j} = \gri{\curve ,j},\] 
in the following sense. In the above situation, the uniform difference
 \[\sup_{x,y \in \mgr} \|\grihat{\mgr,j} (x,y)  - \mathrm{g}^*_{_{\curve,j}}  (x,y)\| \] 
to the pull-back $ \mathrm{g}^*_{_{\curve,j}} \colon \mgr \times \mgr \to \R$ of $\gri{\curve ,j}$ from $\curve$ to $\mgr$ goes to zero.

\smallskip

\item[$(3)$] Moreover, the tuple of functions $(\grihat{t,1},\dots, \grihat{t,r}, \grihat{t,\fin})$ is precisely the pullback to $\mgr$ of the tropical Green function $\gri{\curve(\mgr), \mu(\mgr)}$ on the tropical curve $\curve(\mgr) = \pr_\pi(\mgr)$, associated to the measure $\mu(\mgr)$ on $\curve(\mgr)$ obtained by the pushing out the canonical measure on $\mgr$ to $\curve(\mgr)$.

\end{itemize}
\end{thm}

\subsection{Hybrid curves and their canonical measures}
Having summarized our considerations related to the degenerations of metric graphs and tropical curves, we proceed to present our results concerning Riemann surfaces and hybrid curves. For convenience, we recall first some of the most important notions introduced in our previous work \cite{AN}.

\medskip
Recall that to each stable Riemann surface $S$, we can associate its dual graph $G = (V,E, \genusfunction)$. The vertices of $G$ are the components $C_v$, $v \in V$, of $S$, and the edges between $u,v \in V$ correspond to the intersection points $p^e$ in $ C_u \cap C_v$. The genus function $\genusfunction$ maps every $v \in V$ to the genus of the component $C_v$.

A (stable) \emph{hybrid curve} $\curve$ is defined in terms of the following data:
 \begin{itemize}
 \item a stable Riemann surface $S$ with dual graph $G=(V, E, \genusfunction)$,
 \item an ordered partition $\pi= (\pi_1, \dots, \pi_r)$ of full sedentarity on the edge set $E$,
 \item an edge length function $l\colon E \to (0, + \infty)$ which is well-defined only up to positive rescaling of the layers, that is, modulo multiplication of $l^j:=l\rest{\pi_j}$ by a positive real number, for $j=1, \dots, r$.
 (Imposing $\sum_{e\in\pi_j} l(e)=1$, $j=1, \dots, r$, leads to a unique choice of $l$.)
  \end{itemize}
  With this data, we define the hybrid curve $\curve$ as the \emph{metrized complex} associated to the pair $(S, l)$ enriched with the data of the ordered partition $\pi$. That is, $\curve$ is obtained by taking the components $C_v$, $v \in V$, of $S$, and then replacing for each edge $e \in E$ the intersection point $p^e$ by two new points $p^e_u \in C_u$ and $p^e_v \in C_v$, connected by an interval $\Ical_e$ of length $l(e)$ (see Figure~\ref{fig:HybridCurveIntro}).  By definition, the \emph{rank of a hybrid curve} is the rank of the underlying ordered partition. 
   
\smallskip
   
 Alternatively, a hybrid curve can be defined via uniformization as a conformal equivalence class of layered metrized complexes, with layering induced by the ordered partition $\pi = (\pi_1, \dots, \pi_r)$. In the first definition, we associated at the same time a \emph{canonical metrized complex} $\Sigma$ to $\curve$, the metric realization of the normalized edge length function $l$.
 
\smallskip

Each hybrid curve $\curve$ naturally comes with an \emph{underlying tropical curve} $\curve^\trop$ and the corresponding layered metric $\Gamma$, defined by the graph $G = (V,E)$, the ordered partition $\pi$ and the normalized  edge length function $l$.  The \emph{graded minors }  $\Gamma^1, \dots, \Gamma^r$,  of $\curve$ are defined as the graded minors of the underlying tropical curve $\curve^\trop$. The role of the finite minor $\Gamma^\fin$ in the tropical context exposed in the previous sections is taken over by the disjoint union of components $\pi_\smallcc = \bigsqcup_v C_v$.

\medskip

We define the \emph{canonical measure} $\mu^\can$ on a hybrid $\curve$ of genus $g$ as the sum 
\[
\mu^\can := \frac{1}{g} \Big ( \mu^1_\Zh + \dots + \mu^r_\Zh + \sum_{v \in V} \mu_{C_v}  \Big ).
\]
The measure $\mu_{C_v}$ has support in the component $C_v$, $v \in V$ of $\curve$ and its restriction coincides with the Arakelov--Bergman measure of the smooth Riemann surface $C_v$.  The measure $\mu^j_\Zh$, $j=1, \dots, r$, has support in the union of the intervals $\Ical_e$ in $\curve$ for edges $e\in \pi_j$, and on each such interval $\Ical_e$, it coincides with the Zhang measure of the $j$-th graded minor $\Gamma^j$. 

\medskip

The \emph{moduli space of hybrid curves} of genus $g$, denoted by $\mghyb{\grind{g}}$, and its corresponding universal family ${\mathscr C}^{\hyb}_{\grind{g}}$ of hybrid curves were introduced in our previous work  \cite{AN}. The moduli space $\mgg{\grind{g}}$ of smooth Riemann surfaces of genus $g$ is naturally embedded in $\mghyb{\grind{g}}$. Moreover, the natural map from $\mghyb{\grind{g}}$ to the Deligne--Mumford compactification $\mgbarg{\grind{g}}$, sending each (equivalence class of a) hybrid curve $\curve$ to the (equivalence class of the) underlying stable Riemann surface $S$, is continuous.

\medskip

The main result of  \cite{AN} reads as follows.

\begin{thm}[Continuity of the canonical measure on the moduli space of hybrid curves] \label{thm:hybrid-canonical-measures-intro} The universal family of canonically measured hybrid curves $({\mathscr C}^{\hyb}_{\grind{g}}, \mu^\can)$ forms a continuous family of measured spaces over the hybrid moduli space $\mghyb{\grind{g}}$.
\end{thm}

\subsection{Function theory on hybrid curves and the hybrid Laplacian} 
\label{sec:hybrid_function_theory_intro}

In the third part of our paper, we develop a framework for \emph{function theory in the hybrid setting} and introduce a \emph{hybrid Laplacian}, by mixing notions of analysis on Riemann surfaces and tropical curves (as defined in the previous sections). 

\subsubsection{Functions} 

Let $\curve=\curve^\hyb$ be a hybrid curve of rank $r$ with underlying graph $G=(V, E)$, ordered partition $\pi = (\pi_1, \dots, \pi_r)$ of $E$, and edge length function $l\colon E \to (0, + \infty)$. By construction, $\curve$ has an underlying tropical curve $\curve^\trop$, which is of full sedentarity and has normalized edge lengths on each infinity layer $\pi_j$, $j \in [r]$.  For each vertex $v \in V$, we moreover have a corresponding smooth Riemann surface $C_v$. 

\smallskip

A  function on a hybrid curve $\curve$ of rank $r$ is an $(r+1)$-tuple $\lf = (f_1, \dots, f_r, f_\smallc)$ consisting  of functions $f_j\colon \Gamma^j \to \C$ on the graded minors $\Gamma^j$, $j\in[r]$, and a function $f_\smallc \colon \pi_\smallc \to \C$ on $\pi_\smallc = \bigsqcup_{v\in V} C_v$. 

\smallskip

The tropical (or non-Archimedean) part of $\lf$ is the function $\lf^\trop := (f_1, \dots, f_r)$ on the underlying tropical curve $\curve^\trop$, and  $f_\smallc$ is called the complex (or Archimedean) part of $\lf$. A hybrid function $\lf$ on $\curve^\hyb$ is called continuous, smooth, etc., if the corresponding functions on $\Gamma^j$, $j\in[r]$, and on the $C_v$'s are continuous, smooth, etc., respectively. 

\smallskip

In context with the hybrid Laplacian, we will consider hybrid functions $\lf=(f_j)_j$ on $\curve$ such that each $f_j$, $j \in [r]$, belongs to the Zhang space $D(\Delta_j)$ (see Section~\ref{sec:potential_theory_metric_graphs}). Furthermore, we usually suppose that the complex part $f_\smallc$ is $\mathcal C^2$.

\smallskip When treating the hybrid Poisson equation on hybrid curves, it is necessary to slightly relax these assumptions. In this case, the complex part $f_\smallc$ is allowed to have logarithmic poles at finitely many points of $\bigsqcup_{v\in V} C_v$.  This means there are finitely many points $p_1, \dots, p_n$ on the disjoint union of the $C_v$'s, such that in a local holomorphic coordinate $z_j$ around each point $p_j$, there exists a real number $r_j$ such that the function $f_\smallc - r_j \log|z_j|$ is $\mathcal C^2$.   We define the \emph{extended hybrid Arakelov--Zhang space} $D_{\log}(\Deltahybind{\curve})$ as the set of all hybrid functions $\lf = (f_1, \dots, f_r, f_\smallc)$ with these properties.

\subsubsection{Hybrid Laplacian}  Consider a hybrid curve $\curve$ as above. We define the hybrid Laplacian on $\curve$ by extending the definition of the tropical Laplacian. 

\medskip

The \emph{hybrid Laplacian} $\Deltahyb=\Deltahybind{\curve}$ on $\curve$ is a measure-valued operator which maps hybrid functions $\lf =(f_1, \dots, f_r, f_\smallc) \in D(\Deltahybind{\curve})$ to layered measures  on $\curve$.  It is defined as the sum
\begin{equation} \label{eq:HybridLaplacian}
\Deltahyb(\lf) = \Deltahyb_1(f_1) + \Deltahyb_2(f_2) + \dots \Deltahyb_r(f_r) + \Deltahyb_\smallcc(f_\smallc)
\end{equation}
where the components in the above sum are defined as follows. For $j\in [r]$, and for a function $f_j$ on the graded minor $\Gamma^j$, the layered measure $\Deltatrop_j(f_j)$ on the hybrid curve $\curve$ is given by
\[
\Deltatrop_j(f_j)  := \Bigl(0, \dots, 0, \Delta_j (f_j),  \divind{j}{j+1}(f_j), \dots,  \divind{j}{r}(f_j),\divind{j}{\smallcc} (f_j) \Bigr).
\]
Here, the first $r$ components of $\Deltatrop_j$ coincide with the definition of the tropical Laplacian. The last component is given by setting for $j\in [r]$,
\[
\divind{j}{\smallcc}(f_j)  := -  \sum_{\substack{e \in \pi_j} } \sum_{\substack{v \in e}} \slp_e f_j (v) \, \delta_{p^e_v},
\]
where $p^e_v$ is the point of $C_v$ corresponding to the edge $e$ in the graph. Finally, 
\[\Deltahyb_{\smallcc}(f_\smallc) := \Bigl(0, \dots, 0, \Delta_\smallcc (f_\smallc)\Bigr)\] 
where $\Delta_\smallcc =\frac 1{\pi i}\partial_z \bar\partial_z$ is the Laplacian on the disjoint union $\bigsqcup_v C_v$.

\smallskip

It is easy to see that the measure $\Deltahyb(\lf)$ has mass zero.

\subsubsection{Hybrid Laplacian as a weak limit of Laplacians on Riemann surfaces} Just as in the case of tropical curves, we show that the hybrid Laplacian can be, in a precise sense, obtained as a limit of Laplacians on degenerating Riemann surfaces. Let  $\thy \in \mghyb{g}$ be a hybrid point and denote by $\curve$ the corresponding hybrid curve. For a point $t\in \mgg{g}$ denote by $\unicurve{t}$ the Riemann surface associated to $t$.

Let $\lf$ be a function on $\curve$. Using an extension of the tropical log maps to the hybrid setting, we pullback $\lf$ to nearby Riemann surfaces $\unicurve{t}$ for $t$ living in an appropriate part of the moduli space $\mgg{g}$; see Section~\ref{sec:pullback_rs-intro} below for a more detailed discussion. For any such $t$, denote by $f_t$ the corresponding function on $\unicurve{t}$.   
\begin{thm} \label{thm:RiemannSurfaceLaplacianConvergence-intro}
As $\unicurve{t}$ degenerates to $\curve$, that is as $t$ tends to $\thy$ in $\mghyb{g}$, 
\[
\Delta(f_t)  \to \Deltatrop(\lf)
\]
in the weak sense.
\end{thm}

Convergence in the weak sense is defined as before. Namely, we require that
\[
 \lim_{t \to \thy} \int_{\unicurve{t}} h \,  \Delta f_t = \int_{\curve} h  \, \Deltahyb \lf
\]
for any continuous function $h$ defined on the universal hybrid curve.

\subsection{Potential theory in the hybrid setting}\label{sec:hybrid_green-intro}
As in the tropical setting, we study the solutions of the Poisson equation on hybrid curves and relate them to limits of solutions of the Poisson equation on degenerating families of Riemann surfaces.  

\subsubsection{Poisson equation on a hybrid curve}
 With a notion of Laplacian at hand, we can formulate the Poisson equation on a hybrid curve. 
 
 \smallskip
 
 Suppose that $\curve $ is a hybrid curve with underlying tropical curve $\curve^\trop$ and graph $G = (V,E)$, Riemann surface parts $C_v$, $v \in V$, and ordered partition $\pi$ of rank $r$. We associate to $\curve$ the metric graph $\Gamma$ in the equivalence class of the tropical curve $\curve^\trop$ which has normalized edge lengths in each layer. Let $\mccan$ be the metrized complex with underlying metric graph $\Gamma$ and Riemann surface parts $C_v$, $v\in V$.

 \smallskip

We consider bimeasured hybrid curves $(\curve, \lmu, \lnu)$ for $\lmu$ and $\lnu$ two layered Borel measures on $\curve$ such that $\lmu$ is of mass zero, and $\lnu$ has mass one on the first graded minor $\Gamma^1$ and mass zero on all the components of any higher indexed graded minor, including on each Riemann surface $C_v$, $v\in V$. Analogous to the case of tropical curves (see Section~\ref{ss:IntroFunctionTheoryTropicalCurves}), a layered measure $\lmu$ on $\curve$ is an $(r+1)$-tuple $\lmu = (\mu_1, \dots,  \mu_r, \mu_\smallcc)$ consisting of measures $\mu_j$, $j \in [r]$ on the graded minors $\Gamma^j$, and a measure $\mu_\smallcc$ on $\pi_\smallcc = \bigsqcup_{v \in V} C_v$.

\smallskip

Using the hybrid Laplacian, we introduce the following {\em hybrid Poisson equation}
 \begin{equation}  \label{eq:PoissonHybridCurve-intro}
 \begin{cases}
 \Deltatrop \lf = \lmu \\
\int_{\curve} \lf \, d \lnu = 0
 \end{cases}
\end{equation}
for a function $\lf = (f_1, \dots, f_r, f_\smallcc)$ on $\curve$, where the integration condition means that
\[
\int_{\mccan} f_j \, d\nu = 0 \qquad \text{for all }j=1, \dots, r, \smallcc.
\]
Here $\nu$ is the measure of total mass one on the metrized complex $\mccan$ obtained from the layered measure $\lnu$, and the functions $f_j$ are extended to the full metrized complex $\mccan$ by linear interpolation.

\smallskip

In order to have a well-posed equation in the extended Arakelov--Zhang space $D_{\log}(\Deltahybind{\curve})$, we impose the following natural assumptions on the measures $\lmu$ and $\lnu$. We suppose that the complex part,  that is the part supported on $\pi_\smallcc = \bigsqcup_v C_v$, of any of the two measures can be written as a sum $\alpha + D$ for some continuous $(1,1)$-forms and (real) divisors $D$ on the $C_v$'s. Moreover, the two divisors appearing have disjoint supports.

\begin{thm}[Existence and uniqueness of solutions to the hybrid Poisson equation] The hybrid Poisson equation \eqref{eq:PoissonHybridCurve-intro} has a solution $\lf$ in the extended hybrid Arakelov--Zhang space $D_{\log}(\Deltahybind{\curve})$. Moreover, the solution $\lf$ is unique if we require it to be harmonically arranged. 
\end{thm}

Note that to make the statement of the theorem precise, we need to introduce a relevant notion of being \emph{harmonically arranged} for functions $\lf$  living in the extended hybrid Arakelov--Zhang space. 

In order to do this, we shall choose a local parameter $z^e_v$ around each point $p^e_v$ of $C_v$ for $v\in V$ and incident edge $e$. 

(If the hybrid curve $\curve$ arises as the limit of a family of Riemann surfaces, these choices will be naturally given by fixing an \emph{adapted system of coordinates}, discussed below.)

\smallskip

The complex part $f_\smallcc$ of a function $\lf = (\lf^\trop, f_\smallcc)$ in $D_{\log}(\Deltahybind{\curve})$ can be written as
\[f_{\smallc{_{, v}}}(z^e_v) =  c^e_v\log|z_v^e| + f^\reg_{\smallc{_{, v}}}(z^e_v)\]
locally around the point $p^e_v$ on a component $C_v$, $v \in V$, where $f^\reg_{\smallc{_{, v}}}(z^e_v)$ is a harmonic function defined in the neighborhood of $p^e_v$. 
In particular, \emph{the regularized valued} $f^\reg_{\smallc{_{, v}}}(p^e_v)$ is well-defined.

A hybrid function $\lf = (\lf^\trop, f_\smallc)$ in the extended hybrid Arakelov--Zhang space $D_{\log}(\Deltahybind{\curve})$ is called \emph{harmonically arranged} (with respect to the local parameters $z^e_v$) if
\begin{itemize}
\item the tropical part $\lf^\trop$ is harmonically arranged on $\curve^\trop$, and
\item the one-form $\alpha_\smallc$ on $G$ defined on each oriented edge $e =wv$ by
\[\alpha_\smallc (e) := \frac{f^\reg_{\smallc{_{, v}}}(p^e_v) -f^\reg_{\smallc{_{, w}}}(p^e_w)}{l(e)}\]
belongs to the space of harmonic one-forms on any graded minor of $G$. That is, for each vertex $v$ of the $j$-th graded minor, $j\in[r]$, we have
\[\sum_{e= vw \in \pi_j} \alpha_{\smallc}(e) =0.\]
\end{itemize}

\subsection{Hybrid log maps and pullback of hybrid functions to Riemann surfaces} Let $(S_0, q_1, \dots, q_n)$ be a stable marked Riemann surface of genus $g$ with $n$ markings. Consider the versal deformation space $B$ of $S_0$, which is a polydisc of dimension $3g-3+n$. Denote by $B^*$ the locus of all points corresponding to deformations of $S_0$ which are smooth. $B$ provides an \'etale chart around the point $s_0$ of the moduli space $\mgbar_{\grind{g,n}}$.  Given the local nature of the questions we consider, we can basically work with $B$ and the versal family of stable marked Riemann surfaces $\rsf \to B$ over it. 

\smallskip

In~\cite{AN} we introduced the hybrid replacement $B^\hyb$ of $B$ and over it, the family of hybrid curves
\[\rsf^\hyb \to B^\hyb.
\]
The hybrid boundary $\partial_\infty B^\hyb  = B^\hyb \setminus B^\ast$ is naturally stratified into a disjoint union $\partial_\infty B^\hyb = \bigsqcup_\pi D_\pi^\hyb$, where $\pi$ runs over all ordered partitions of full sedentarity on subsets of the node set $E$ of $S_0$. The points of the hybrid stratum $D_\pi^\hyb$ represent hybrid curves with underlying ordered partition $\pi$.

For a given system of coordinates on $B$, there is a natural projection map $\loghyb_\pi \colon B^\ast \to D_\pi^\hyb$ to each hybrid stratum $D_\pi^\hyb$. As we discuss in Section~\ref{ss:LogMapHybridTropicalCompactification}, there is another continuous map
\[
	B^\hyb \longrightarrow \mgtropcombin{\grind{G}}
\]
where $\mgtropcombin{\grind{G}}$ is the space of tropical curves whose underlying combinatorial graph $G = (V,E)$ coincides with the dual graph of $S_0$. The above map sends the hybrid boundary $\partial_\infty B^\hyb = B^\hyb \setminus B^\ast$ to the tropical boundary $\partial_\infty \mgtropcombin{\grind{G}} = \mgtropcombin{\grind{G}} \setminus \mggraphcombin{\grind{G}}$. On the open part $B^\ast$, it is given by
\begin{align*}
t \in B^\ast \qquad  \mapsto \qquad  \big (-\log\abs{z_e(t)} \big)_{e\in E} \in \mggraphcombin{\grind{G}},
\end{align*}
where $z_e$, $e \in E$, is a parameter on $B$ corresponding to the smoothing of the node $p_e$ of $S_0$. 

 Using the tropical log map $\logtrop$ on $\mgtropcombin{\grind{G}}$, we construct a retraction map $\loghyb \colon B^* \to \partial_\infty B^\hyb$ to the full hybrid boundary $\partial_\infty B^\hyb$.  The two \emph{hybrid log maps} $\loghyb$ and $\loghyb_\pi$ are precisely the hybrid analogues of the tropical log maps $\logtrop$ and $\pr_\pi$ in the setting of tropical curves.

\smallskip

\begin{figure}[!t]
\centering
   \scalebox{.28}{\input{example-adcoord.tikz}}
\caption{The decomposition given by the rectification theorem.}
\label{fig:example-adcoord}
\end{figure}
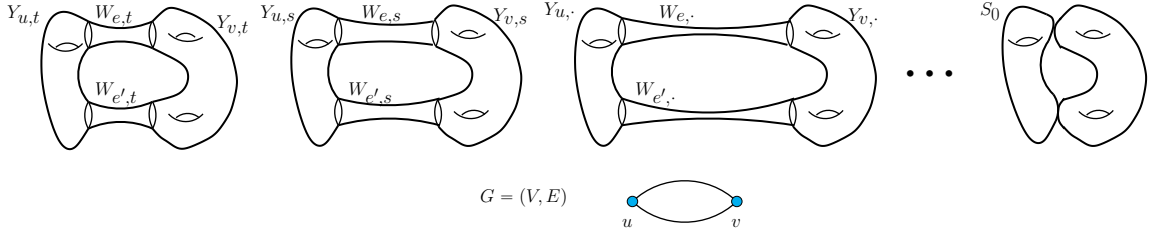

We then construct log maps from the family $\rsf$ over $B$ to the restriction of the hybrid family $\rsf^\hyb$ to the boundary  $\partial_\infty B^\hyb$. Roughly speaking, this amounts in defining a log map which goes from each smooth Riemann surface corresponding to a point $s$ in $B$ to the hybrid curve corresponding to the point $\shy =\loghyb(s)$ (or $\shy =\loghyb_\pi(s)$). Our construction is based on a rectification theorem, which provides an \emph{adapted system of coordinates}.  The latter allows to decompose the smooth Riemann surfaces into cylinders indexed by graph edges of the graph and smaller genus Riemann surfaces with boundaries, indexed by vertices of the graph, Figure~\ref{fig:example-adcoord}. Using this surgery, we define the hybrid log maps by taking logarithmic polar coordinates on the cylinder. The existence of adapted coordinate systems is based on the work of Hubbard-Koch~\cite{HK14} which, using the plumbing constructions~\cite{Mas76,Wol90}, gives an analytic description of the Deligne--Mumford compactification.

\smallskip

Using these constructions, we extend the log maps on the base to the hybrid versal family, and arrive at the \emph{stratumwise} and \emph{global}  hybrid log maps. More precisely, we define two maps
\begin{align*}
\loghyb_{\pi} \colon \rsf^\ast \to \rsf_\pi^\hyb &&\loghyb \colon \rsf^\hyb\rest{\mathscr U}  \to  \rsf^\hyb\rest{\partial_\infty B^\hyb}
\end{align*}
where $\rsf^\ast=\rsf^\hyb\rest{B^\ast}$ and $\rsf^\hyb_\pi =\rsf^\hyb\rest{D^\hyb_\pi}$, such that the following two diagrams commute:
 \begin{equation} \label{eq:CDHybridLogs}
\begin{tikzcd}
 \rsf^\ast \arrow[d]\arrow[r, "\loghybdiagpi{\pi}"] & \arrow[d] \rsf_\pi^\hyb\\
B^\ast\arrow[r, "\loghybdiagpi{\pi}"] &D_\pi^\hyb
\end{tikzcd}
\qquad  \qquad \qquad
\begin{tikzcd}
 \rsf^\hyb\rest{\mathscr U} \arrow[d]\arrow[r, "\loghybdiag"] & \arrow[d] \rsf^\hyb\rest{\partial_\infty B^\hyb}\\
\mathscr U\arrow[r, "\loghybdiag"] &\partial_\infty B^\hyb 
\end{tikzcd}
\end{equation}
with $\mathscr U$ a neighborhood of the hybrid boundary $\partial_\infty B^\hyb$.

\smallskip

Using the log maps, we  pull-back functions on hybrid curves to nearby Riemann surfaces. This procedure is parallel to pulling back functions from tropical curves to metric graphs. However, the presence of logarithmic poles in the complex part $f_\smallcc$ of a hybrid function $\lf = (f_1, \dots, f_r, f_\smallcc)$ requires some modifications.

\subsubsection{Extending hybrid functions to the metrized complex} We first explain how the components $f_j$ of a hybrid function $\lf=(f_j)_j$ on a hybrid curve can be viewed as functions on the associated metrized complex.

\smallskip

Let $\curve$ be a hybrid curve and $\mccan$ the associated metrized complex. Denote by $l:E \to (0, + \infty)$ the edge length function in $\mccan$ and let $\Gamma$ be the corresponding metric graph.  Let $\lf = (\lf^\trop, f_\smallc)$ be a hybrid function on $\curve^\hyb$. Its tropical part $\lf^\trop = (f_1, \dots, f_r)$ is a function on the tropical curve $\curve^\trop$ underlying $\curve$. We allow the complex part $f_\smallc$ to have logarithmic poles at a finite number of points. 

For each $j\in[r]$, we define the pullback $f_j^\ast \colon \mccan \to \C$ by taking first the pullback of $f_j$ to $\Gamma$, and then by composing it with the projection map $\mccan \to \Gamma$ (which contracts each Riemann surface $C_v$ to the vertex $v$). 

The pullback of $f_\smallc$ to $\mccan$ is the function 
\[f_{\smallc}^{\ast}:\mccan \to \C\]
obtained by linear interpolation of the regularized values $f^\reg_\smallc(p^e_u)$ on the edges. More precisely, one each edge $e=uv\in E$, it is given by
\[
f_{\smallc}^\ast\rest e (x) = f_\smallc^\reg (p^e_u) + \frac{f^\reg_\smallc(p^e_v) - f^\reg_\smallc(p^e_u)}{l(e)}x.
\]
Here $x$ is the parametrization of the oriented edge $e=uv$ in  $\Gamma$ (identified with the interval $[0,l(e)]$) and $p^e_u$ and $p^e_v$ are the corresponding marked points on $C_u$ and $C_v$, respectively.

\subsubsection{Pullback of hybrid functions to Riemann surfaces} \label{sec:pullback_rs-intro} Consider the previous setting, with $\rsf/B$ the analytic versal deformation space of a stable marked Riemann surface $(S_0, q_1, \dots, q_n)$, and let  $\rsf^\hyb/B^\hyb$ be the corresponding hybrid family. We assume that an adapted coordinate system $\underline z = (z_e)_{e\in E}$, is fixed on $B$, and consider the associated hybrid log maps.

\smallskip

Fix a point boundary point $\thy \in \partial_\infty B^\hyb$, lying in the hybrid stratum $D_\pi^\hyb$ of an ordered partition $\pi=(\pi_1, \dots, \pi_r)$, and consider the preimage $\mathscr U_\thy := \bigr({\loghyb}\bigr)^{-1}(\thy)$. A hybrid function $\lf = (f_1, \dots, f_r, f_\smallc)$ on $\rsf^\hyb_\thy$ can be pulled back to the Riemann surface fibers $\rsf_t$, $t \in \mathscr U_\thy\cap B^\ast$ as follows.

\smallskip
For each $j=1, \dots, r$, introduce the logarithmic coefficient
\[
L_j(t) = - \sum_{e\in \pi_j} \log|z_e(t)|
\]
and consider the extended function $f_j^* \colon \rsf^\hyb_{\thy} \to \C$. Note 
that here and in the definition of the log map, $\rsf^\hyb_\thy$ is viewed as a metrized complex with normalized edge lengths on each layer $\pi_j$, $j\in[r]$.  

\smallskip

For a point $t \in \mathscr U_\thy \cap B^\ast$, we consider the hybrid log map between fibers
\[\loghyb_{t} \colon \rsf_{t} \to \rsf^\hyb_{\thy}.\]

We define the pullback $\lf^*_{t}\colon \rsf_{t} \to \C$ by
\[\lf^*_{t} := \sum_{j=1}^r L_j(t) f_j^*\circ \loghyb_{t}  + f_\smallc^*\circ \loghyb_{t}.\]

Note that if the hybrid function $\lf$ is real-valued, then the pullback $\lf^*_t$ is also real-valued.

\medskip

 Similarly as above, we can define an alternate pullback of hybrid functions by considering stratumwise log maps. That is, consider the logarithm map $\loghyb_{\pi} \colon B^\ast \to D_\pi^\hyb$ associated to the stratum $D_\pi^\hyb$ of the fixed point $\thy$. Then the expression
\[\lf^*_{t} := \sum_{j=1}^r L_j(t) f_j^*\circ \loghyb_{\pi,t}  + f_\smallc^*\circ \loghyb_{\pi,t}\]
gives a pullback to all Riemann surfaces $\rsf_t$, $t \in B^\ast$, satisfying $\loghyb_{\pi}(t) =\thy$, where again
\[\loghyb_{\pi,t} \colon \rsf_{t} \to \rsf^\hyb_{\thy}\]
is the map between the fibers.  The two notions of pullbacks are obviously distinct. In the sequel, we will use both the definitions, and each time, it will be clear from the context which of the two pullbacks we are referring to.

\subsection{Tameness of solutions of the hybrid Poisson equation} Using the hybrid log maps and pullbacks, we can define the notions of \emph{stratumwise tameness and weak tameness} for families of hybrid functions on $\rsf^\hyb/B^\hyb$. The definitions are analogous to the case of tropical curves and details can be found in Section~\ref{ss:TameHybridFunctions}.

\smallskip

We state our main theorems concerning the variation of solutions of the Poisson equations on Riemann surfaces.

\subsubsection{The case of discrete measures} Consider a stable marked Riemann surface  with $n+1$ markings $(S_0, p, q_1, \dots, q_n)$, and let $D=d_1q_1+\dots + d_nq_n$ be a divisor of degree zero on $S_0$ with real coefficients.   Let $B$ be the versal deformation space of the marked curve $S_0$.
 
 For each point $\thy\in B^\hyb,$ consider the degree zero divisor $D_\thy = \sum_{j=1}^n d_j q_j(\thy)$ on the hybrid curve $\rsf^\hyb_\thy$. Let $\lf_\thy$, $\thy \in B^\hyb$, be the unique harmonically arranged solution to the Poisson equation
  \begin{equation} \label{eq:JDFcthyb-discrete-intro}
\begin{cases}
 \Deltatrop (\lf) = D_\thy \\
 f_j (\proj_j(p_\thy)) = 0, \qquad j=1, \dots, r, \smallc,
 \end{cases}
\end{equation}
on the hybrid curve $\rsf_\thy$. Here $\kappa_j$, $j=1, \dots, r, \smallcc$, are the contraction maps from $\pi_\smallcc = \bigsqcup_v C_v$ to the graded minors $\Gamma^1, \dots, \Gamma^r$ and $\pi_\smallcc$. 

\begin{thm}[Stratumwise tameness of solutions of the Poisson equation: discrete measures] \label{thm:TameContinuityHybridCurves-intro} Notations as above, the family of hybrid functions $\lf_\thy$, $\thy\in B^\hyb$, form a stratumwise tame family of hybrid functions on $\rsf^\hyb$.
\end{thm} 

By definition of stratumwise tameness, the theorem thus leads to the following local expression of the solutions.

\begin{thm}[Local expression of solutions of the Poisson equation: discrete measures]  
Notations as above, let $\thy\in B^\hyb$ be a hybrid point with underlying ordered partition $\pi=(\pi_1, \dots, \pi_r)$. Let $\lf = (f_1, \dots, f_r, f_\smallc)$ the solution of the Poisson equation on the hybrid curve $\rsf^\hyb_\thy$. Then the difference 
\[f_t - \sum_{j=1}^r L_j(t) f_j^*\circ \loghyb_{\pi,t}  - f_\smallc^*\circ \loghyb_{\pi,t}, \qquad t\in (\loghyb_{\pi})^{-1}(\thy),\]
extends continuously to zero over $\thy$. That is, for $t \in (\loghyb_{\pi})^{-1}(\thy)$ converging to $\thy$, the difference tends to zero uniformly on $\rsf_t$.  Here $L_j(t) = - \sum_{e\in \pi_j} \log\abs{z_e(t)}$.
\end{thm}

The condition $t \in (\loghyb_{\pi})^{-1}(\thy)$ can be dropped by considering the pullbacks for the hybrid function $\lf_\thy$, $\thy = \loghyb_\pi(t)$, globally, see the definition of stratumwise tameness in Section~\ref{ss:TameHybridFunctions}.

\subsubsection{Weak tameness for general measures} 
We consider now the more general case of bimeasured hybrid curves. Let thus $(S_0,q_1, \dots, q_n)$ be a stable curve with $n$ markings, and let $\rsf^\hyb/B^\hyb$ be the associated hybrid versal family. We consider two continuous families of measures $\lmu_\thy$ and $\lnu_\thy$ on fibers $\rsf^\hyb_\thy$, $\thy \in B^\hyb$,  satisfying fiberwise the conditions in Section~\ref{sec:hybrid_green-intro}. Moreover, we require that the discrete parts of the complex pieces $\mu_{\thy, \smallcc}$ and $\nu_{\smallcc, \thy}$ are supported on two disjoint subsets $Q_\mu$ and $Q_\nu$ of $\{q_1(\thy), \dots, q_n(\thy) \}$ and the continuous parts satisfy certain natural bounds (see Section~\ref{ss:SurfacesPoissonGeneralLimit}).

\smallskip

 We are interested in the behavior of the solutions $f_\thy$, $\thy \in B^\hyb$, to the hybrid Poisson equation
 \begin{equation} \label{eq:poisson_general_measure-intro}
\begin{cases}
 \Deltatrop (\lf_\thy) = \lmu_\thy  \\[1 mm]
 \lf_\thy \text{ is harmonically arranged} \\[1mm]
\int_{\rsf^\hyb_\thy} \lf_\thy \, d \lnu_\thy = 0
 \end{cases}
\end{equation}
close to the hybrid boundary $\partial_\infty B^\hyb$.

\begin{thm}[Weak tameness of solutions of the Poisson equation] Notations as above, the family of functions $\lf_\thy$, ${\thy\in B^\hyb}$, is weakly tame on $B^\hyb$.
\end{thm}  

In addition we will establish the following concrete asymptotics. Fix a boundary stratum $D_\pi^\hyb$ in $B^\hyb$. For a point $t\in B^\ast$, let $\thy = \loghyb_\pi(t)$, and consider the hybrid log map
\[\loghyb_{\pi,t}\colon \rsf_t \to \rsf^\hyb_\thy.\]

Denote by $\loghyb_{t_*}(\mu_t)$ and $\loghyb_{t_*}(\nu_t)$ the push-out measures on $\rsf^\hyb_\thy$, and let $\hat\lf_t$ be the unique solution to the Poisson equation
 \begin{equation} \label{eq:JDFcthyb}
\begin{cases}
 \Deltatrop (\hat \lf_t) = \loghyb_{t_*}\mu_t \\[1 mm]
 \hat \lf_t \text{ is harmonically arranged} \\[1mm]
\int_{\rsf^\hyb_\thy} \hat \lf_t \, \,  d  \loghyb_{t_*}\nu_t  = 0
 \end{cases}
\end{equation}
on the hybrid curve $\rsf^\hyb_\thy$, $\thy = \loghyb_\pi(t)$. Let
\[
\hat\lf^\ast_t = L_1(t) \hat f_{t,1}^\ast + L_2(t) \hat f_{t,1}^\ast + \dots + L_r(t) \hat f_{t,r}^\ast + \hat f_{t, \smallcc}^\ast
\]
be the pullback of $\hat\lf_t$ to the smooth Riemann surface $\rsf_t$.

\begin{thm} Notations as above, suppose that $t \in B^\ast$ approaches $D_\pi^\hyb$, that is, $t$ converges to a point $\shy$ in $D_\pi^\hyb$.

Then, the difference $f_t -\hat\lf^\ast_t$ goes to zero uniformly on $\rsf_t$.
\end{thm}

\subsection{Hybrid Green functions}

Let $\curve$ be a hybrid curve of rank $r$ and $\lmu$ a (layered) measure of total mass one on $\curve$, which on the Riemann surface part of $\curve$ is given by a continuous $(1,1)$-form. Denote by $\mccan$ the metrized complex associated to $\curve$.

\smallskip

Any point $p \in \mccan$ induces a layered Dirac measure $\bm{\delta}_p$ on $\curve$. The pair $(\lmu, \bm{\delta}_p -\lmu)$ gives a bimeasured hybrid curve, and the hybrid Poisson equation
\begin{equation} \label{eq:HybridGreenFunction-intro}
\begin{cases}
 \Deltatrop \lf = \bm{\delta}_p - \lmu \\[1 mm]
 \lf \text{ is harmonically arranged} \\[1mm]
 \int_\curve \lf \, d \mu = 0
 \end{cases}
\end{equation}
has a unique solution $\lf = (f_1, \dots, f_r, f_\smallc)$ that we denote by
\[
\lgri{\lmu}(p, \cdot) = \big ( \gri{\lmu, 1}(p, \cdot), \dots  ,\gri{\lmu, r}(p, \cdot), \gri{\lmu, \smallcc}(p, \cdot) \big).
\]
Note that the component $\gri{\lmu, \smallcc}(p, \cdot)$ usually has logarithmic singularities at a finite number of points, and thus, the Green function belongs to the extended hybrid Arakelov--Zhang space. 
  
By the previous considerations, each component $\gri{\lmu, j}(p, \cdot)$, $j\in [r]\cup\{\smallcc\}$, can be viewed as a function on the full space $\mccan$. 
Therefore, we get a {\em hybrid Green function} 
\begin{equation} \label{eq:IntroHybridGreenFunction}
\begin{array}{cccc}
\gri{\curve, \mu} \colon &\mccan \times \mccan  &\longrightarrow & \R^{r+1} = \R^{[r]}\times \eRm, \qquad \eRm :=\R \cup \{-\infty\},
\end{array}
\end{equation}
where $\mu$  is the measure on $\mccan$ obtained from the layered measure $\lmu$. The hybrid Green function of the canonical measure $\lmu = \lmu^\can$ is simply denoted by $\gri{\curve}$.

\subsection{Asymptotics of the Arakelov Green function} \label{sec:general_asymptotics_intro}
We now give a presentation of our results concerning the Arakelov Green function.  Denote by  $\mgg{\combind{g,2}}$ the moduli space of Riemann surfaces of genus $g$ with two marked points. Consider the global Arakelov Green function defined naturally on $\mgg{\combind{g,2}}$,
\[
\begin{array}{cccc}
\grs \colon &\mgg{\combind{g,2}}  &\longrightarrow & \R \\
& (S, p,y) &\mapsto & \gri{s}(p,y),
\end{array}
\]
where $s$ denotes the point of $\mgg{\grind{g}}$ which represents the Riemann surface $S$ and $\gri{s}(\cdot, \cdot)$ denotes the Green function on $S$.

\medskip

We first describe the asymptotics of the Arakelov Green function in $\mghyb{\combind{g,1}}$, the moduli space of hybrid curves of genus $g$ with one marked point. In this setting, we obtain a full description in terms of Green functions on hybrid curves with one marked point.

For each point $\thy \in \mghyb{\combind{g,1}}$, we denote by $\gri{\thy}(p_\thy, \cdot)$ the canonical Green function on the corresponding hybrid curve $\curve_\thy$ with marking $p_\thy$.

Our theorem can be stated as follows.
\begin{thm} \label{thm:main_Arakelov_intro} The family of functions $\gri{\thy}(p_\thy, \cdot)$ on $\mghyb{\combind{g,1}}$ is weakly tame.
\end{thm}

More precisely, we establish the following \emph{asymptotic expansion of the Green function}.

Consider a point $\thy  \in \mg^\hyb_{\combind{g,1}}$ and let $\curve^\hyb$ be the corresponding hybrid curve with one marking. Denote by $(S_0, p_0)$ the underlying stable marked Riemann surface, $G=(V,E)$ the dual graph of $S_0$, and let $\rsf/B$ be the family over the versal deformation space $B$ of $(S_0, p_0)$. We fix an adapted system of parameters on $\rsf/B$. 

Let $\pi = (\pi_1, \dots, \pi_r)$ be the ordered partition of $E$ underlying the hybrid curve $\curve$. For each point $t\in B^\ast$, define the \emph{logarithmic coefficients}
\begin{equation} \label{eq:IntroLj}
L_j(t) = -\sum_{e\in \pi_j}\log\abs{z_e(t)}, \qquad j=1, \dots, r,
\end{equation}
where $z_e$, $e \in E$, is the parameter on $B$ describing the smoothings of the node $p_e$ in $S_0$.

\begin{thm}[Expansion of the Arakelov Green function] \label{thm:MainArakelovIntro} There exists a family of functions $\grihat{t,j}(p_t,\cdot)$ on $\rsf_t$, $j=1, \dots, r, \smallcc$, and $t\in B^*$, such that the following asymptotics holds. 

\begin{itemize}
\item[$(1)$] The difference 
\[\gri{t}(p_t, \cdot) - \Big ( L_1(t)\grihat{t,1}(p_t,\cdot) + L_2(t)\grihat{t,2}(p_t,\cdot) + \dots + L_r(t)\grihat{t,r}(p_t,\cdot) + \grihat{t, \smallcc}(p_t, \cdot) \Big )
\]
 extends uniformly by zero to $\thy$. That is, for a sequence of points $t \in B^\ast$ converging to $\thy$ in $B^\hyb$, the above difference tends to zero uniformly on $\rsf_t$.

\smallskip

\item[$(2)$] The convergence
 \[\lim_{t\to \thy} \grihat{t,j}(p_t,\cdot) = \gri{\thy,j}(p_\thy, \cdot)
 \] 
holds in the following sense. If $t \in (\loghyb_\pi)^{-1}(\thy)$ converges to $\thy$, then the difference
 \[
 \grihat{t,j}(p_t, y) - \mathrm{g}^*_{_{\thy, j}}(p_t, y), \qquad j=1, \dots,r,
 \] 
goes to zero uniformly for $y \in \rsf_t$. Moreover,
\[
 \grihat{t, \smallcc}(p_t, y) - \mathrm{g}^*_{_{\thy, \smallcc}}(p_t, y)
\]
goes to zero uniformly for $y \in \rsf_t$ remaining uniformly separated from the appearing nodes.
\smallskip

\item[$(3)$] Moreover, the tuple of functions $(\grihat{t,1}(p_t,\cdot),\dots, \grihat{t,r}(p_t,\cdot), \grihat{t,\smallcc}(p_t,\cdot))$ is precisely the pullback to $\rsf_t$ of the hybrid Green function on $\rsf^\hyb_{\shy}$ for $\shy=\loghyb_\pi(t)$, associated to the measure on $\rsf^\hyb_\shy$ which is obtained by pushing out the canonical measure of $\rsf_t$ via the log map $\loghyb_{\pi,t} \colon \rsf_t \to \rsf^\hyb_{\shy}$. 
\end{itemize}
\end{thm}
We refer to Section~\ref{sec:MainArakelovDetails} for a more detailed statement of the theorem.

\smallskip

We emphasize that the quantities $L_1(t), \dots, L_r(t)$ go to infinity when $t$ goes to $\thy$ in the situation of the above theorem. Moreover, by the definition of the topology on  $\mg^\hyb_{\combind{g}}$, we have that $L_1(t) \gg L_2(t) \gg \dots \gg L_r(t)$ under this degeneration. 
The theorem can be thus viewed as providing an analytic description of the Arakelov Green function close to the hybrid limit.

\smallskip

Note that from the above theorem, we can deduce similar asymptotic results for other algebraic families of Riemann surfaces. Consider a generically smooth family of stable Riemann surfaces $\pi: \rsf \to X$. Taking a log-resolution if necessary, we can assume  the discriminant locus of the family, that is, the locus of points $t$ in $X$ with a non-smooth fiber in the family,  forms a simple normal crossing divisor  $D$. To the pair $(X, D)$, we can associate the corresponding hybrid space $X^\hyb$, the hybrid compactification of $X^*:= X \setminus D$ and the corresponding family of hybrid curves $\rsf^\hyb$. It follows essentially from our theorem above that

\begin{thm} \label{thm:families} Fix a section $(p_\thy)_{\thy\in X^\hyb}$ of $\rsf^\hyb /X^\hyb$ such that $p_\thy$ belongs to the smooth part of the hybrid curve $\rsf_\thy^\hyb$ for all $\thy$. Then the family of functions $\gri{\thy}(p_\thy, \cdot)$ is weakly tame. \end{thm}

A special case of the theorem concerns families which vary over a one-dimensional base, i.e., when $\dim(X)=1$. The divisor $D$ in this situation corresponds to a finite set of points, and our hybrid construction replaces the stable curves $\rsf_{t}$, for $t\in D$, with the corresponding metrized complex $\mccan_{t}$. We can suppose without loss of generality that $X =\Delta$ is a unite disc, and $t =0$.
Denote by $\mgr$ the underlying stable metric graph of $\mccan$, the metrized complex at $0$. We have a projection map $\mccan \to \mgr$, and the push-out of the canonical measure on $\mccan$ coincides with the canonical measure on the metric graph $\mgr$. Denote by $\mu^\can$ the canonical measure on $\mccan$.

We get the following particularly nice statement on the asymptotics of the Arakelov Green functions in this case. 
\begin{thm}[Expansion of the Arakelov Green function: one-parameter families] \label{thm:main_one_parameter} Notations as above, let $p\colon \Delta \to \rsf$ be a section of the family with $p_0$ in the smooth part of some component $C_v$ of the stable Riemann surface $\rsf_0$. Then, the difference 
\[\gri{t}(p_t, \cdot) + \log\abs{t} \gri{\mgr}(v,\cdot) + h_{t}(v,\cdot)- \gri{\smallcc}(p_t,\cdot), \qquad t\in \Delta^*,\]
extends by zero over $0$. Here $(\gri{\mgr}(v,\cdot), \gri{\smallcc}(p_t,\cdot))$ is the solution to the Poisson equation on the bimeasured hybrid curve $\mccan$ with pairs of measures $(\delta_{p_t} - \mu^\can, \mu^\can)$, and $h_{t}(v,\cdot)$ is a correction term originating from the higher order terms in the convergence of the Riemann surface canonical measures to the canonical measure on the metric graph.
\end{thm}
The term $h_{t}(v,\cdot)$ is constant in the interior of the Riemann surface part of the hybrid curve, sufficiently away from the nodes. So in particular, we get the full asymptotic up to a constant in the regime where points converge to a point on the smooth part of the hybrid curve. It might be possible to get an explicit form of these constants  via the matricial analysis carried out in the second part of the paper. 

\smallskip

Finally, we address the asymptotics of the Arakelov Green function in the moduli space of hybrid curves $\mghyb{\grind{g}}$. The subsequent considerations provide an answer to Question~\ref{question:main-graphs-precise}, posed at the beginning of the introduction.

\smallskip

Consider a hybrid curve $\curve$ having an underlying stable Riemann surface $S_0$ with dual graph $G=(V,E)$, an ordered partition $\pi = (\pi_1, \dots, \pi_r)$ on the edge set $E$, and a normalized edge length function $l$. Let $\rsf^\hyb \to B^\hyb$ be the family of hybrid curves associated to a versal family $\rsf \to B$ of $S_0$, equipped with adapted coordinates. The hybrid curve $\curve$ appears as the fiber $\rsf^\hyb_\shy$ over some hybrid point $\shy$ in the hybrid stratum $D_\pi^\hyb$.

We obtain the following asymptotic expansion of the Arakelov Green function, in terms of hybrid Green functions and an additional correction term $\varepsilon_\pi(p,y )$.

\begin{thm}[Full expansion of the Arakelov Green function] \label{thm:MainArakelovIntro2}  There exists a family of functions $\grihat{t,j}(\cdot,\cdot)$ on $\rsf_t \times \rsf_t$, $j=1, \dots, r, \smallcc$ and $t\in B^*$, such that the following asymptotics holds. 

\begin{itemize}
\item[$(1)$] The difference 
\[\gri{t} -  \Big ( L_1(t)\grihat{t,1} + \dots + L_r(t)\grihat{t,r} + \grihat{t, \smallcc}- \varepsilon_\pi \Big )
\]
 extends uniformly by zero to $\thy$. That is, for a sequence of points $t \in B^\ast$ converging to $\thy$ in $B^\hyb$, the above difference tends to zero uniformly on $\rsf_t \times \rsf_t$.
 
 Here $L_j(t)$ is the logarithmic coefficient given by \eqref{eq:IntroLj} and $\varepsilon_\pi(\cdot,\cdot)$ is an explicit correction term defined in Section~\ref{ss:MainArakelovDetails2}.

\smallskip

\item[$(2)$] The convergence
 \[\lim_{t\to \thy} \grihat{t,j}= \gri{\thy,j}
 \] 
 holds in the following sense. If $t \in (\loghyb_\pi)^{-1}(\thy)$ converges to $\thy$, then the difference
 \[
 \grihat{t,j}(p, y) - \mathrm{g}^*_{_{\thy, j}}(p, y), \qquad j=1, \dots,r,
 \] 
goes to zero uniformly for $p,y \in \rsf_t$. Moreover,
\[
 \grihat{t, \smallcc}(p, y) - \mathrm{g}^*_{_{\thy, \smallcc}}(p, y)
\]
goes to zero uniformly for $p,y \in \rsf_t$ remaining uniformly separated from the appearing nodes.
\smallskip

\item[$(3)$] Moreover, the tuple of functions $(\grihat{t,1},\dots, \grihat{t,r}, \grihat{t,\smallcc})$ is precisely the pullback to $\rsf_t \times \rsf_t$ of the hybrid Green function on $\rsf^\hyb_{\shy} \times \rsf^\hyb_{\shy} $ for $\shy=\loghyb_\pi(t)$, associated to the measure on $\rsf^\hyb_\shy$ which is obtained by pushing out the canonical measure of $\rsf_t$ via the log map $\loghyb_\pi\colon \rsf_t \to \rsf^\hyb_{\shy}$. 
\end{itemize}
\end{thm}

We refer to Section~\ref{ss:MainArakelovDetails2} for a more careful statement of Theorem~\ref{thm:MainArakelovIntro2} and the explicit description of the correction term $\varepsilon_\pi(p,y)$. The latter is a sum of contributions which are reminiscent of the height pairing between divisors on a Riemann sphere.

\subsection{Comparison to the works of Faltings, de Jong, and Wentworth}

Our theorem on Arakelov Green functions, stated  above, solves essentially completely a long-standing open problem in Arakelov geometry by providing higher rank refinements and generalizations of the works of Faltings~\cite{Faltings21}, de Jong~\cite{deJong}, and Wentworth~\cite{Went91}.

The work of de Jong~\cite{deJong} concerns degenerations of Riemann surfaces on the punctured disc, that is one parameter algebraic families $\rsf/\Delta^*$ of Riemann surfaces. The main theorem of~\cite{deJong}, on what concerns the asymptotics of Green functions, is a slightly weaker version of Theorem~\ref{thm:main_one_parameter}, namely, the statement that for distinct sections $p_t$ and $q_t$ of $\rsf_t$ which converge to smooth points $p_0$ and $q_0$ of $\rsf_0$, the Arakelov Green function has the asymptotic behavior 
\[\gri{\rsf_t}(p_t, q_t) \sim - \gri{\mgr}(v_p, v_q)\log\abs{t}.\]
In this equation, $\mgr$ is the dual stable metric graph associated to $\rsf_0$ and $v_p, v_q$ are the vertices of $\mgr$ which accommodate the points $p_0$ and $q_0$. 
A similar statement can be phrased in the case where one considers the regularization of the Arakelov Green function on the diagonal. Note that Theorem~\ref{thm:main_one_parameter} treats as well the case where the points $p_t$ and $q_t$ go to the singular point of the special fiber. 

\smallskip

The statement of Theorem~\ref{thm:main_one_parameter} for one-parameter families of Riemann surfaces with only one node in the special fiber $\rsf_0$ was proved in the pioneering work of Wentworth~\cite{Went91}. Note that he obtains a full asymptotic description in this case (that is, with a description of the correction constants in Theorem \ref{thm:main_one_parameter}) by taking into account higher order vanishing terms in the asymptotic of the canonical measure.

\smallskip

More recently, Faltings studies the degenerations of Arakelov Green functions over the moduli space of Riemann surfaces, and explains that the asymptotics could be described in terms of the corresponding metric graph limits in the situation where the parameters converge with the same speed to zero, that is, formulated in the language of ours papers, if the corresponding hybrid limit has only one infinity layer. The precise asymptotics is not worked out there due to the absence of the asymptotic description of the canonical measures. The methods of~\cite{Faltings21} are based on the study of degenerations of Abelian varieties.

\subsection{Future developments} This paper and its preceding companion~\cite{AN} are part of a series of works dedicated to an analytic study of the geometry of moduli spaces and their corresponding geometric objects close to their boundaries.  

In our upcoming work~\cite{AN-AG-hybrid}, we develop an algebraic geometry for hybrid curves, and show that the geometry of Riemann surfaces survives in the hybrid limit, that is, hybrid curves enjoy fundamental theorems governing the geometry of smooth compact Riemann surfaces. 
 
We will furthermore study the spectral theory of hybrid curves and use them to control the behavior of eigenvalues of Riemann surfaces close to the boundary of their moduli spaces.  
 
In another direction, we will consider the asymptotics of modular graph functions arising in string theory~\cite{HGB19}. The results of~\cite{AN} on convergence of canonical measures and our Theorem~\ref{thm:MainArakelovIntro2} allow to obtain a refined description of the asymptotic behavior of the modular graph functions close to the hybrid limit.

\smallskip

From a larger perspective, the results in our two available papers and their forthcoming companions suggest that the layered behavior close to the boundary of moduli spaces is a broad phenomenon. This will be further developed in our future work.

\smallskip

 The study undertaken in this paper and its sequels are furthermore related to a higher rank version of non-Archimedean geometry. In this regard, a geometric study of higher rank valuations is undertaken in~\cite{AI}, and the underlying higher rank polyhedral geometry has been further developed in~\cite{Iri22}. The precise connection to our results will appear elsewhere.

\subsection{Basic notations}  By curve in this paper we mean a complex projective algebraic curve. The analytification of a curve is a  compact Riemann surface.  We use the terminology curve and Riemann surface deliberately to switch back and forth between the algebraic and analytic setting. 

For a non-negative integer $r$, the symbol $[r]$ means the set $\{1,\dots, r\}$. For a pair of non-negative integers $k\leq r$, we denote by ${[r] \choose k}$ the family of subsets  $S\subset [r]$ of size $k$.

By \emph{almost all} in this paper we mean for all but a finite number of exceptions.

Given an ambient set $[N]$ for an integer $N$, and a subset $E$ of $[N]$, we denote by $E^c$ the complement of $E$ in $[N]$.

The letter $G$ is used for graphs.  We use the symbols $\mgr$,  $\curve^\trop$, $\mc$ and $\curve^\hyb$ for metric graphs, tropical curves, metrized complexes and hybrid curves, respectively. Moreover, when working with the family of curves over a versal deformation space, we use the notation $\rsf^\hyb$ for the corresponding family of hybrid curves over the hybrid deformation space.

In context with asymptotics, we will also make use of the standard Landau symbols. Let $X$ be a topological space, $x$ a point in $X$ and $g \colon U \setminus \{x\} \to \R_+$ a function on a punctured neighborhood of $x$. Recall that the expressions $O(g)$ and $o(g)$ are used to abbreviate functions $f$ such that the quotient $|f(y)| / g(y)$ remains bounded and goes to zero, respectively, as $y$ goes to $x$ in $X$ (here, $f$ is defined on a punctured neighborhood of $x$ as well).

 For a $g\times g$ matrix $P$, with rows and columns indexed by $[g]$, and for $I, J \subset [g]$, we denote by $P[I, J]$ the matrix with rows in $I$ and columns in $J$.

In this manuscript, $\log(x) := \log_e(x)$ denotes the logarithm of $x>0$ with respect to the standard base $e$. That is, $e^{\log(x)} = x$ for all $x >0$.

\newpage
\section{Glossary of notation}
This section summarizes the most important notations used in this paper.

\subsection*{Basic notation}
\quad\\[1mm]
$\N = \Z \cap [1, + \infty) = \{1, 2, \dots \}$ is the set of natural numbers. \\
$\Z_{\ge 0} = \Z \cap [0, + \infty) = \N \cup \{0 \}$ is the set of non-negative integers. \\
$[n] = \{1, 2, \dots, n \}$ for $n \in \N$. \\
$\R_+ = [0, + \infty)$ and $\R_{>0} =( 0, +\infty)$. \\
$\eR = \R \cup \{+ \infty\}$. \\
$\eRm = \R \cup \{- \infty\}$. \\
For a finite set $A$, $\keg_A = \R_+^A$ and $\inn \keg_A = \R_{>0}^A$ are closed and open cones in $\R^A$. \\
$\sigma_A$ and $\inn \sigma_A$ are the open and closed  standard simplex in $\R^A$. \\
For a finite set $A$, we sometimes denote the number of elements in $A$ by $|A|$.

\subsection*{Compactifications of fans}
\quad\\[1mm]
$\eta$ denotes a cone in $N_\R$; $\eta^\vee$ and $\eta^\perp$ are the dual cone and orthogonal complement of $\eta$. \\
$\Sigma$ denotes a fan in $N_\R$. \\
$\cancomp \Sigma$ is the canonical compactification of $\Sigma$. \\ \smallskip
$\inn \keg_\eta^\delta \cong \eta^\delta$ is the open stratum in $\cancomp \Sigma$ associated with two faces $\delta \subseteq \eta$. \\[2mm]
$\cancomp \Sigma^\trop$ is the higher rank compactification of $\Sigma$. \\
$\pi $ denotes a flag of faces in $\Sigma$ of the form $\pi\colon \zerocone = \delta_0 \subsetneq \delta_1 \subsetneq \dots \delta_{r} \subseteq \delta_{\infty + \fin}$; $r$ is the depth of $\pi$. \\
$\mathscr{F}(\Sigma)$ denotes the set of all such flags.\\
"$\preceq$" is the partial order on $\mathscr{F}(\Sigma)$ given by refinement of flags. \\
$\inn \keg_{\pi}^\trop$ is the open stratum in $\cancomp \Sigma$ associated with a flag $\pi \in \mathscr{F}(\Sigma)$.

\subsection*{Metric graphs and tropical curves}
\quad\\[1mm]
$G= (V,E)$ is a finite graph with vertex set $V$ and edge set $E$. \\
$\deg(v)$ is the degree of a vertex $v \in V$. \\
$H_1(G, \Z)$ is the first homology of $G$; $\graphgenus$ usually denotes the genus of $G$. \\[2mm]
$\marking: [n] \to V$ is marking function, placing $n$ markings at the vertices of $G$. \\
$\genusfunction: V \to \mathbb N \cup\{0\}$ is genus function. \\
$G = (V, E, \marking, \genusfunction)$ is an augmented graph with markings; $g$ is the genus of an augmented graph. \\[2mm]
$l \colon E \to \R_{>0}$ denotes an edge length function. \\
$\mgr$ is a metric graph, a metric realization of $G = (V,E)$ and an edge length function $l$. \\[2mm]
$\pi = (\pi_\infty, \pi_\fin)$ denotes an ordered partition of a finite set $F$. \\
$\pi_\infty = (\pi_1, \dots, \pi_r)$ is the sedentarity part of $\pi$; the integer $r$ is called the depth of $\pi$. \\
$\pi_\fin$ is called the finitary part of $\pi$. \\
$\Piall(F)$ is the set of all ordered partitions of a finite set $F$.\\
$\Pifs(F)$ is the subset of ordered partitions of full sedentarity. \\
For a fixed, finite set $E$, $\Pihat :=\bigcup_{F\subseteq E} \Piall(F)$ denotes the set of all ordered partitions of subsets $F \subseteq E$; "$\preceq$" is the partial order on $\Pihat$ given by refinement of ordered partitions. \\[2mm]
$(G, \pi)$ denotes a layered graph. \\
$G_\pi^j$ is the $j$-th graph in the decreasing  sequence of spanning subgraphs of $G$ coming from $\pi$. \\
$\grm{\pi}{j}(G)$ is the $j$-th graded minor of $(G, \pi)$; the vertex set of $\grm{\pi}{j}(G)$ is denoted by $V^j$ or $V_\pi^j$. \\
$\cont{j}:  V \to V^j_\pi$ is the contraction map. \\
$G^\fin$ is the finitary minor of $(G, \pi)$; it has vertex set $V$. \\
$\graphgenus^j_\pi$ is the genus of the minor $\grm{\pi}{j}(G)$.\\
and $J_\pi^j$ is the $j$-th subset in the decomposition of $[\graphgenus]$ induced by $\graphgenus = \graphgenus_\pi^1 + \dots \graphgenus_\pi^\fin$.
\\[2mm]
$\curve^\trop$ (or sometimes $\curve$) denotes a tropical curve.\\
$\Gamma$ is the canonical layered metric graph in the class of $\curve$, which arises from the normalized edge length function $l$ of $\curve$; \\
$\Gamma^i$ is the $i$-th graded minor of $\Gamma$. \\[2mm]
$L(\mgr)$ is the total length of a metric graph $\mgr$. \\
$L_i(\mgr)$ is the total length of the $i$-th minor of a layered metric graph $\mgr$.

\subsection*{Moduli spaces of graphs and tropical curves}
\quad\\[1mm]
$\mggraph{\grind{g,n}}$ is the moduli space of stable metric graphs of genus $g$ with $n$ marked points; \\
$\mgtrop{\grind{g,n}}$ is the moduli space of stable tropical curves of genus $g$ with $n$ marked points;  \\
$\mggraph{\grind{g}}$ and $\mgtrop{\grind{g}}$ are the corresponding moduli spaces for $n = 0$ marked points. \\
$\mgtrop{\grind{G}}$ and $\mgtrop{\combind{(G, \pi)}}$ are the moduli spaces of tropical curves with underlying graph $G$ and combinatorial type $(G, \pi)$, respectively. \\
$\umggraph{\grind{g}}$ is the moduli space of metric graphs of genus $g$ with lengths bounded by one. \\[2mm]
$\mggraph{\grind{G}}$, $\mgtrop{\grind{G}}$ are the spaces of metric graphs and tropical curves with underlying graph $G$. \\
$\partial_\infty \mgtrop{\grind{G}} = \mgtrop{\grind{G}} \setminus \mggraph{\grind{G}}$ is the subspace corresponding to curves of non-empty sedentarity. \\
$\mgtropcombin{\combind{(G, \pi)}}$ and $\mgtropcombin{\combind{(G, \subface \pi)}}$ are the spaces of tropical curves with combinatorial type equal to and coarser than $(G, \pi)$, respectively.\\[2mm]
$\hcurveg{G}^\trop$, $\unicurvetrop{\combind{(G,\subface \pi)}}$, and $\unicurvetrop{g,n}$ are the corresponding universal curves over $\mggraphcombin{G}$, $\mgtropcombin{\combind{(G, \subface \pi)}}$, and $\mgtrop{\grind{g,n}}$. \\
$\unicurvetrop{\thy}$ is the fibre over a base point $\thy$ in such a universal family. \\
$\umgr$ is the universal curve over $\mggraph{\grind{G}}$; its fibres are denoted by $\umgr_t$.  \\

\subsection*{Analysis on metric graphs and tropical curves}
\quad\\[1mm]
The star $\ast$ indicates a pulled back object (e.g., from a tropical curve to a metric graph). \\
$\lf = (f_1, \dots, f_r, f_\fin)$ denotes a function on a tropical curve $\curve$ of rank $r$; its $j$-th part $f_j$ is a function on the minor $\Gamma^j$ and can be viewed as a function on $\Gamma$. \\
$\lf^\ast$, $f_j^\ast$ are the pullbacks of $\lf$, $f_j$ from $\curve$ to a metric graph $\mgr$ over the same graph $G = (V,E)$. \\
$\mathcal{P}_{jk}$ for two natural numbers $j \le k$, denotes the set of strictly increasing paths from $j$ to $k$. \\
$A_p(\curve)$, $A_p(\mgr)$ are the matrices associated to $p \in \mathcal{P}_{jk}$ and a tropical curve $\curve$ or a layered metric graph $\mgr$ (see \eqref{eq:AsmyptoticsApTropical} and \eqref{eq:MatrixApGraph}). \\[2mm]
$\Lognoind = \logtrop{}$ is the log map from $\mggraph{\grind{G}}$ to $\partial_\infty \mgtrop{\grind{G}}$. \\
$\pr_\pi \colon  \mggraphcombin{\grind{G}} \to \mgtrop{\grind{(G, \pi)}}$ is the projection map from $\mggraph{\grind{G}}$ to a fixed boundary stratum $\mgtropcombin{\combind{(G, \pi)}}$.\\[2mm]
$\mathcal{M}(\mgr)$ is the space of complex-valued Borel measures on a metric graph $\mgr$. \\
$\mathcal{M}^0(\mgr) \subset \mathcal{M}(\mgr)$ is the subspace consisting of mass zero measures. \\
$\widetilde{\mathcal{M}}^0(\mgr) \subset \mathcal{M}(\mgr)$ is the subspace consisting of mass zero measures of the form \eqref{eq:SpecialMus}. \\
$\lmu = (\mu_1, \dots, \mu_r, \mu_\fin)$ denotes a layered measure on a tropical curve $\curve$ of rank $r$. \\
$\mass$ is the mass function of such a measure $\lmu$. \\
$\mathcal{M}_\pi(\curve)$ is the space of layered measures on $\curve$. \\
$\widetilde{\mathcal{M}}_\pi^0(\curve) \subset \mathcal{M}_\pi(\curve)$ is the subspace of mass zero measures whose pieces $\mu_j$ are of the form $ \eqref{eq:SpecialMus}$.\\[2mm]
$\innone{}{\cdot\,,\cdot}$ is the pairing between $1$-chains in the (metric) graph.\\ 
$\rmM(\mgr) = \rmM_l = \rmM(G, l)$ is the period matrix \eqref{eq:M_l} for a metric graph $\mgr$ (arising from $(G,l)$).\\
$\mu_\Zh$ denotes the Zhang measure on a metric graph; it has mass equal to the genus. \\
$\lmu^\can$ is the canonical measure on a tropical curve. \\[2mm]
$(\mgr, \mu, \nu)$, $(\curve, \lmu, \lnu)$ denote bimeasured metric graphs and bimeasured tropical curves. \\
$D(\Delta)$ or $D(\Deltaind{\mgr})$ is the Zhang space of a metric graph $\mgr$. \\
$\operatorname{BDV}(\mgr)$ is the space of functions with bounded differential variation on $\mgr$. \\
$\jfunc{p \tiret q, x} (y)$ denotes the $\jvide$-function of a metric graph $\mgr$. \\
$\slp_e f(v)$ is the slope of a function $f \colon \mgr \to \C$ at a vertex $v$ along an incident (half-)edge $e$. \\
$\Deltaind{\mgr}$ (or simply $\Delta$) is the Laplacian on a metric graph $\mgr$. \\
$\Deltatrop_{\grind{\curve}}$ (or simply $\Deltatrop$) is the Laplacian on a tropical curve $\curve$. \\[2mm]
$\grg_{\mu} \colon \mgr \times \mgr \to \R$ is the Green function of a measure $\mu$ on a metric graph $\mgr$. \\
$\gri{\curve, \mu}  \colon \Gamma \times \Gamma \to \R^{r+1}$ is the Green function of a measure $\lmu$ on a tropical curve $\curve$ of rank $r$. \\[2mm]
$C^1(G, \R)$ is the space of one-forms on a graph $G = (V,E)$. \\
$\partial \colon C^1(G, \R) \to \Div(G)$ is the boundary operator. \\
$\Omega^1(G) = \ker(\partial)$ is the space of harmonic one-forms on $G$. \\
$df$ is the differential of a function $f \colon V \to \R$ on a graph $G$ (w.r.t. to some $l \colon E \to \R_{>0}$). \\
$\exact(G, l)$ is the space of exact one-forms on $(G,l)$. \\
$ \projhar$ is the orthogonal projection onto $\Omega^1(G)$ (w.r.t. an edge length function $l \colon E \to \R_{>0}$). \\
$ \projexact$ is the the orthogonal projection onto $\exact(G,l)$.\\[2mm]
$C^1(\curve, \R)$ is the space of one-forms on a tropical curve $G = (V,E)$. \\
$\ld(\lf)$ is the differential of a  function $\lf$ on $\curve$. \\
$\exact(\curve)$ is the space of exact one-forms on $\curve$. \\[2mm]
$\Div(G)$, $\Div^0(G)$ denote the spaces of (real) divisors and degree zero divisors on $G = (V,E)$.\\
$\hp{\mgr}{\cdot\,, \cdot}$, $\lhp{\curve}{\cdot\,, \cdot}$ denote the height pairings on metric graphs $\mgr$ and tropical curves $\curve$. \\

\subsection*{Riemann surfaces, hybrid curves and their moduli spaces}
\quad\\[1mm]
$S$ is a stable (marked) Riemann surface with stable dual graph $G = (V,E, \genusfunction, \marking)$. \\
$C_v$, $v \in V$ are the components of the normalization of $S$. \\
$\mc$ is a metrized complex. \\
$\curve$ or $\curve^\hyb$ is a hybrid curve with underlying stable Riemann surface $S$, ordered partition $\pi = (\pi, \dots, \pi_r)$ in $\Pi(E)$ and (normalized) edge length function $l \colon E \to \R_{>0}$. \\
$\Gamma^1, \Gamma^2, \dots, \Gamma^r$ are the graded minors of $\curve$. \\
$\pi_\smallcc = \bigsqcup_{v} C_v$ is the complex (or finitary) part of the hybrid curve $\curve$. \\
$\mccan$ is the canonical layered metrized complex in the class of $\curve$, which arises from the normalized edge length function $l$ of $\curve$  \\[2mm]
$\mgg{\grind{g,n}}$ is the moduli space of Riemann surfaces with $n$ marked points. \\
$\mgbarg{\grind{g,n}}$ is the moduli space of stable Riemann surfaces with $n$ marked points, the Deligne--Mumford compactification of $\mgg{\grind{g,n}}$. \\
$\mgg{\grind g}$ and $\mgbarg{\grind{g}}$ are the corresponding moduli spaces for $n = 0$ marked points. \\
$\mghyb{\grind{g,n}}$ is the moduli space of hybrid curves with $n$ marked points. \\
$\mghyb{\grind{g}} = \mghyb{\combind{g,0}}$ is the moduli space of hybrid curves. \\[2mm]
$B$ and $D$ denotes a pair of a complex manifold $B$ and simple normal crossing divisor $D$; usually, $B = \Delta^N$ is an $N$-dimensional polydisc and $D$ is a union of coordinate hyperplanes $D_e = \{ z_e(t) = 0 \}$, $e \in E$ for a subset $E \subseteq [N]$.\\
This situation appears when $B$ is the deformation space of  a stable Riemann surface $S$. \\
$B^\ast = B \setminus D$ is called the finite part. \\
$\rsf \to B$ is the versal deformation family associated to a stable Riemann surface $S$. \\
$\rsf^\ast \to B^\ast$ is the restriction of $\rsf$ to $B$. \\
$B^\hyb$ is the hybrid space associated to a pair $(B,D)$. \\
$\partial_\infty B^\hyb$ is the boundary at infinity of $B^\hyb$, the locus of points of non-empty sedentarity.\\
$D_\pi^\hyb$ is a hybrid stratum in the hybrid space $B^\hyb$. \\
$\rsf^\hyb \to B^\hyb$ is the hybrid versal deformation family associated to a stable Riemann surface $S$. \\
$\rsf^\hyb_\pi \to D^\hyb_\pi$ is the restriction of $\rsf$ to the hybrid stratum $D_\pi^\hyb$. \\

\subsection*{Analysis on Riemann surfaces and hybrid curves} \quad \\[1mm]
$\innone{S}{\cdot\,,\cdot}$ is the Hermitian inner product between holomorphic one forms on $S$.\\
$\intprod{}{\cdot\,,\cdot}$ is the intersection product between cycles in $S$.\\
The star $\ast$ indicates a pulled back object (e.g., from a hybrid curve to a smooth Riemann surface). \\
$\lf = (f_1, \dots, f_r, f_\smallcc)$ denotes a function on a hybrid curve $\curve$ of rank $r$. \\
$\loghyb$ is the log map from $B^\ast$ to $\partial_\infty B^\hyb$. The corresponding map from $\rsf^\ast$ to $\rsf^\hyb$ is denoted by $\loghyb$ as well.  \\[1mm]
$\loghyb_{\pi} \colon  B^\ast \to D_\pi^\hyb$ is the projection map from $B^\ast$ to a fixed boundary stratum $D_\pi^\hyb$. The corresponding map from $\rsf^\ast \to B^\ast$ to $\rsf^\hyb_\pi \to D_\pi^\hyb$ is denoted by $\loghyb_\pi$ as well.\\
$\lf^\ast$ denotes a pullback of a function $\lf$ on a hybrid curve to a smooth Riemann surface, using one of the maps $\loghyb$ and $\loghyb_{\pi}$. \\[2mm]
$\mu_{\Ar}$ is the  Arakelov--Bergman measure on a smooth Riemann surfaces $S$, it has mass the genus. $\mu^\can =\frac 1g\mu_{\Ar}$ is the canonical measure. \\
$\lmu^\can$ is the canonical measure of a hybrid curve. \\[2mm]
$\Deltaind{S} = (\pi i)^{-1} \partial_z \partial_{\overline{z}}$ (or simply $\Delta$) is the Laplacian on a smooth Riemann surface $S$. \\
$\Deltahyb_{\grind{\curve}}$ (or simply $\Deltahyb$) is the Laplacian on a hybrid curve $\curve$. \\[2mm]
$\gri{S} \colon S \times S \to  \R$ is the Arakelov Green function on a smooth Riemann surface $S$. \\
$\lgri{\curve}\colon \curve \times \curve \to \R^{r}\times \eRm$ is the canonical Green function of a hybrid curve $\curve$ of rank $r$. \\[2mm]
$\Div$, $\Div^0$ denote spaces of divisors and degree zero divisors with real coefficients, respectively.\\
$\hp{S}{\cdot\,, \cdot}$, $\hp{\mc}{\cdot\,, \cdot}$, $\lhp{\curve}{\cdot\,, \cdot}$ are the height pairings on smooth Riemann surfaces $S$, metrized complexes $\mc$, and hybrid curves $\curve$. \\[2mm]
$\varepsilon_\pi(p,x)$ and $\varepsilon_\pi(p,q,x,y)$ are correction terms appearing in the description of height pairings and Green functions (see~\eqref{eq:CorrectionTwoPoints} and~\eqref{eq:CorrectionFourPoints}).

\newpage


\section{Preliminaries} \label{sec:preliminaries} The aim of this section is to introduce the main objects of study in this paper. Building on our previous work~\cite{AN}, we will generalize the definition of tropical curves studied in~\cite{AN} by associating a notion of \emph{sedentarity}, which encodes in which direction and in which depth at \emph{infinity} a tropical curve lives while at the same time allowing to have finitary edges. Those tropical curves considered in~\cite{AN} correspond to the ones of \emph{full sedentarity} in this paper.

\subsection{Graphs}  As far as graphs are concerned, the terminology and notation in this paper is the same as in~\cite{AN}. We briefly recall the relevant notions in the present section. For further information, we refer to~\cite{AN} and classical text books either in graph theory~\cite{Bol, BM, Diestel} or in algebraic geometry~\cite{ACGH}.

\smallskip

All graphs which appear in this paper are finite. We allow parallel edges and loops. We use the notation $V(G)$ and $E(G)$ for the vertex set and the edge set of a graph $G$, respectively, and if $G$ is understood from the context, we simply use $V$ and $E$. We also write $G=(V, E)$ to indicate the vertex and edge sets. 
 
 Two vertices $u$ and $v$ in $G$ are called \emph{adjacent} if there is an edge in $E$ with extremities $u$ and $v$, in such a case we write $u \sim v$. A vertex $v$ is \emph{incident} to an edge $e$ and we write $e\sim v$ if $v$ is an extremity of $e$.  The \emph{degree} of a vertex $v$ denoted by $\deg(v)$ counts the number of half-edges incident to $v$ and is formally defined as 
  \[\deg(v) := \Bigl|\bigl\{e\in E \,\bigl|\, e \sim v\bigr\}\Bigr| + \#\textrm{ loops based at $v$}.\]

 A \emph{subgraph} in $G=(V,E)$ is a graph $H=(U, F)$ with $U \subseteq V$ and $F \subseteq E$. A subgraph $H=(U, F)$ is called \emph{spanning} if it has the same vertex set as $G$, that is, if $U=V$. A spanning subgraph is thus uniquely determined by its set of edges.  For spanning subgraphs $H_1=(V, E_1), \dots, H_r = (V, E_r)$ of $G$, the union $H_1 \cup \dots \cup H_r$ is defined as the spanning subgraph of $G$ with edge set $E_1 \cup \dots \cup E_r$.

\smallskip

For a subset of edges $F \subseteq E$ and a subset $U \subseteq V$, the notation $G[F]:=(V, F)$ refers to the spanning subgraph of $G$ with edge set $F$, while $G[U]$ is the \emph{induced subgraph of $G$ on $U$} which by definition has vertex set $U$ and the edge set consists of all the edges in $E$ with both extremities belonging to $U$.

\smallskip

A \emph{spanning tree} in $G=(V,E)$ is by definition a spanning subgraph $T$ of $G$ which has no cycle and which has the same number of connected components as $G$. This requires $T$ to contain the maximum possible number of edges without containing any cycle. The set of all spanning trees of $G$ is denoted by $\cT(G)$. 
  
\smallskip

We now recall the definition of a \emph{minor} of a given graph $G$. The terminology refers to the Robertson-Seymour theory of graph minors, which has been central in the development of graph theory and its applications. As in the previous work~\cite{AN},  special kinds of minors associated to the layerings of edges that we call graded minors, and we recall below, play an important role in this paper. All the minors of $G$ will appear (possibly several times) as a graded minor. 

 The definition of minors in a graph is based on two simple operations on graphs called \emph{deletion} and \emph{contraction}. Given $G=(V, E)$ and an edge $e$ in $E$, we denote by $G-e$ the spanning subgraph $H = (V, E\setminus \{e\})$ of $G$ obtained by deleting the edge $e$. Moreover, by $\contract{G}{e}$ we mean the graph obtained by contracting the edge $e$: the vertex set in $\contract{G}{e}$ is the set $U$ obtained by identifying the two extremities of  $e$, and the edge set is $E \setminus \{e\}$ that we view in $U$ via the canonical projection map $V \to U$. More generally, for a subset $F \subseteq E$, we denote by $G-F$ and $\contract{G}{F}$ the graphs obtained by deletion and contraction of $F$ in $G$, respectively: $G-F$ is the spanning subgraph of $G$ with edge set $E\setminus F$, and $\contract{G}{F}$ is the graph obtained by contracting one by one all the edges in $F$. (One shows that the order does not matter.)  For the purpose of clarification, we mention that later on we will introduce a variant of the construction of contractions for \emph{augmented graphs}, which are graphs endowed with an  integer valued \emph{genus function} on vertices. In this case, the contraction remembers the \emph{genera of connected components of the contracted part}, thus keeping the genus of the augmented graph constant before and after contraction.

A \emph{minor} of a graph $G$  is a graph $H$ which can be obtained by contracting a subset of edges in a subgraph $H'=(W, F)$ of $G$.

\medskip

For a graph $G$, we denote by $H_1(G)$ the first homology group of $G$ with $\Z$ coefficients. The rank of $H_1(G)$, which coincides with the first Betti number of $G$ viewed as a topological space, is called \emph{genus} of the graph. In this paper we usually use the letter $\graphgenus$ when referring to the genus of graphs, and reserve the letter $g$ and $\genusfunction$ for the genus of augmented graphs, tropical curves, algebraic curves, and their analytic and hybrid variants introduced later. 

\noindent The genus $\graphgenus$ of a graph $G=(V, E)$ is equal to $\abs{E} - \abs{V}+ c(G)$ where $c(G)$ denotes the number of connected components of $G$.

\subsection{Ordered partitions} \label{subsec:OrdPart} We recall the definition of the ordered partitions, central combinatorial objects in dealing with tropical and hybrid compactifications, by extending the set-up of~\cite{AN} and allowing a (possibly empty) \emph{finitary part}.  

\smallskip

 Let $E$ be a finite set and let $F \subseteq E$ be a non-empty subset. An \emph{ordered partition} of $F$ is an ordered pair $\pi = (\pi_\infty, \pi_\fin)$ consisting of 
 \begin{itemize}
 \item a subset $\pi_\fin \subseteq F$ called the \emph{finitary part} of $\pi$, and 
 \item an ordered sequence $\pi_\infty = (\pi_1, \pi_2, \dots, \pi_r)$,  $r\in \N \cup\{0\}$, consisting of non-empty subsets of $F$ which form a partition of $F\setminus \pi_\fin$. We call $\pi_\infty$ the \emph{sedentarity part} of $\pi$.
  \end{itemize}
 Equivalently, an ordered partition is an ordered sequence $\pi = (\pi_1, \pi_2, \dots, \pi_r, \pi_\fin)$, $r\in \N \cup\{0\}$, of  pairwise disjoint subsets of $F$ such that the following holds
 \begin{itemize}
 \item $\pi_1, \dots, \pi_r$ are non-empty, and 
 \item we have 	$F = \bigsqcup_{i=1}^r \pi_i \sqcup \pi_\fin.$
 \end{itemize}
 
 Note that we allow $\pi_\fin$ to be empty.  The integer $r$ is called the \emph{rank} of $\pi$. It is allowed to be zero, in which case $\pi_\infty =\varnothing$ and $\pi$ consists of $\pi_\fin =F$. 

Note that $F = \varnothing$ has a unique ordered partition with sedentarity and finitary parts $\varnothing$.

For an ordered partition $\pi = (\pi_\infty, \pi_\fin)$ as above, the sedentarity part $\pi_{\infty}=(\pi_1, \dots, \pi_r)$ is viewed as parts of the ordered partition lying at infinity, consisting of different layers of infinity distinguished according to the integers $1, \dots, r$.

An ordered partition is called of \emph{full sedentarity} if $\pi_\fin=\varnothing$. 

The set of ordered partitions of a subset $F \subseteq E$ is denoted by $\Piall(F)$. Those of full sedentarity form a subset of $\Piall(F)$ that we denote by $\Pifs(F)$.

The ordered partitions which appear in our previous work~\cite{AN} are precisely obtained from those of full sedentarity in this paper by omitting the mention of $\pi_\fin=\varnothing$.  In practice, we can thus view an ordered partition $\pi =(\pi_1, \dots, \pi_r)$ of~\cite{AN} as an ordered partition of full sedentarity by adding $\pi_\fin = \varnothing$.   

\smallskip

Each $\Piall(F)$ has a distinguished element denoted by $\pi_\sedind{\varnothing}=(F)$ which is the ordered partition of \emph{empty sedentarity}, that is of rank zero with finitary part $\pi_\fin=F$. This is sometimes denoted by $\pi_\varnothing(F)$ or simply $\pi_\varnothing$ if $F$ is understood. 

Note in particular that $\Piall(\varnothing) = \Pifs(\varnothing) = \{\pi_\varnothing=(\varnothing)\}$.

Given an ordered partition $\pi=(\pi_\infty=(\pi_1, \dots, \pi_r), \pi_\fin)$, we set $E_\pi := \pi_1 \cup \dots \cup \pi_r\cup \pi_\fin$, the subset of $E$ underlying the ordered partition. We denote by $E_{\pi_\infty}$ or simply $E_{\infty}$ if $\pi$ is understood, the union $E_{\pi_\infty} := \pi_1\cup \dots \cup \pi_r$ and call it the \emph{sedentarity set} or the \emph{infintary set} of the ordered partition. 

\smallskip

We have the \emph{sedentarity map} $\sed$ which sends $\pi=(\pi_\infty, \pi_\fin)$ to $\sed(\pi) = \pi_\infty$,
\begin{align*}
\sed\colon \Piall:=\bigsqcup_{F\subseteq E}\Piall(F) &\to \Pifs:=\bigsqcup_{F\subseteq E} \Pifs(F)\\
\pi=(\pi_\infty, \pi_\fin) &\mapsto \pi_\infty.
\end{align*}
\smallskip

There is a bijective correspondence between the ordered partitions and specific types of \emph{filtrations}.  By a \emph{filtration} $\filter$ on a (non-empty) subset $F \subseteq E$ of \emph{rank} $r$ we mean an increasing sequence of non-empty subsets $F_1 \subsetneq F_2  \subsetneq F_3 \subsetneq \dots \subsetneq F_r\subseteq F_{\infty+\fin}=F$ for $r \in \mathbb N\cup\{0\}$. We set $F_\infty :=F_r$.  Note that we allow $F_\infty = F_{\infty+\fin}$.

Each filtration $\filter$ as above defines an ordered partition $\pi = (\pi_1, \dots, \pi_r, \pi_\fin)$ of $F$ by setting
\[
	\pi_i = F_i \setminus F_{i-1}, \qquad i=1,...,r, \,\textrm{ and } \pi_\fin =F_{\infty+\fin}\setminus F_{\infty},
\]
 with $F_0:=\varnothing$ by convention. In order to emphasize that the filtration is the one associated to the ordered partition $\pi$, we sometime use the notation $\filter^\pi$. For $F = \varnothing$, again, we allow $\filter^\varnothing := \varnothing$ as its only filtration. 
The maximum element of the filtration $\filter$ is denoted by $\max(\filter)$ or $F_{\infty+\fin}$. We thus have $E_\pi = \max(\filter^\pi)$.

\medskip

We now define a natural partial order  on the set of ordered partitions (equivalently, filtrations) of subsets of a given set $E$.

\smallskip

 A filtration $\filter$ is a \emph{refinement} of the filtration $\filter'$ if we have 
 \begin{itemize}
 \item $F_{\infty+\fin} \subseteq F'_{\infty+\fin}$, and
 \item as a set we have $\filter' \setminus\{F'_{\infty+\fin}\} \subseteq \filter \setminus\{F_{\infty+\fin}\}$.
 \end{itemize} 
 The filtration $\filter$ is called a \emph{tame refinement} of $\filter'$ if the above properties hold and moreover, we have $F_j = F'_j$ for $j=1, \dots, r'$ with $r'$ denoting the rank of $\filter'$. 
 
 \smallskip
An ordered partition $\pi$ of a subset of $E$ is called a \emph{refinement} of an ordered partition $\pi'$ of a (possibly different) subset of $E$, and we write $\pi' \preceq \pi$, if this holds true for the corresponding filtrations  $\filter$ and  $\filter'$. That is, denoting $\pi' = (\pi_1',\dots, \pi_r', \pi_\fin)$, then $\pi' \preceq \pi$ exactly when $\pi$ is obtained from $\pi'$ performing the following two steps. 
\begin{itemize}
\item[(i)] Replacing each set $\pi_i'$ in the ordered sequence $\pi'_\infty = (\pi_i')_{i=1}^r$ by an ordered partition ${\varrho}^i = ({\varrho}_k^i)_{k=1}^{s_i}$ of $\pi_i'$.
\item[(ii)] Possibly removing a subset $B$ of $\pi'_\fin$ and adding an ordered partition $\varrho=(\varrho_1, \dots, \varrho_s, \varrho_\fin)$ of  $\pi'_\fin \setminus B$ at the end of the ordered sequence $\pi'_\infty$. We thus have $\pi_\fin = \varrho_\fin$.
\end{itemize}
An ordered partition $\pi$ is called a {\em tame refinement} of an ordered partition $\pi'$ if the corresponding filtration $\filter^\pi$ to $\pi$ is a tame refinement of the filteration $\filter^{\pi'}$ corresponding to $\pi'$. This means $\pi$ is a refinement of $\pi'$, and moreover, for any $j=1, \dots, r'$ with $r'$ the rank of $\pi'$, we have $\pi_j = \pi'_j$.

\begin{remark} We stress once more that in the definition of a refinement, the ordered partitions  are not necessarily defined over the same subset $F \subseteq E$. However, if $\pi$ is a refinement of $\pi'$, then by the properties above, we get $E_{\pi'_\infty} \subseteq E_{\pi_\infty} \subseteq E_\pi \subseteq E_{\pi'}$. 
 This means, refinements should be regarded as referring to the refinements of the sedentarity part. The set $E_{\pi'}$ provides the ambient set containing the underlying sets of all the refinements. In cannot be increased.
\end{remark}

One easily verifies that "$\preceq$" defines a partial order on $ \Piall =\bigcup_{F\subseteq E} \Piall(F)$, the set of all the ordered partitions of subsets $F \subseteq E$.
\begin{example} \label{ex:op}
Assume that $E = \{e,f \}$ for two elements $e \neq f$. Then the ordered partitions of subsets $F \subseteq E$ are given by
\begin{align*}
	\Piall ( \varnothing ) &= \big \{ \pi_\varnothing \big  \} \qquad \Pi ( \{e\}) = \Big \{ \big (\{e \} \big ), \big (\{e \}, \varnothing \big ) \Big  \} \qquad
	\Piall( \{f\}) =\Big  \{\big  (\{f \} \big ),  \big(\{f \}, \varnothing\big ) \Big \} \\
	\Piall(E) &= \Big \{ \big (\{e\}, \{ f\} \big ), \big (\{e\}, \{ f\}, \varnothing \big ),\big  (\{f\}, \{ e\} \big ), \big  (\{f\}, \{ e\} , \varnothing \big ), \big  (\{e, f\}\big  ), \big  (\{e, f\}, \varnothing\big  )  \Big  \}.
\end{align*}
The partial order "$\preceq$" on $\Piall = \bigcup_{F\subseteq E} \Piall(F)$ consists of the following relations

\begin{minipage}[t]{0.5\textwidth}
\begin{align*}
\big ( \{e,f\} \big ) & \preceq \big  (\{e\}, \{ f\} \big ) \preceq \big  (\{e\}, \{ f\}, \varnothing \big ), \\
	\big ( \{e,f\}\big ) &\preceq \big (\{f\}, \{ e\} \big ) \preceq \big (\{f\}, \{ e\}, \varnothing \big ),\\
		\big ( \{e,f\} \big ) & \preceq \big  (\{e, f\}, \varnothing \big ) \preceq \big  (\{e\}, \{ f\}, \varnothing \big ), \\
				\big ( \{e,f\} \big ) & \preceq \big  (\{e, f\}, \varnothing \big ) \preceq \big  (\{f\}, \{ e\}, \varnothing \big ), \\
	\big ( \{e,f\}\big ) &\preceq \big (\{f\}\big ) \preceq \big (\{f\}, \varnothing \big ),\\
	\big ( \{e,f\}\big ) &\preceq \big (\{e\}\big ) \preceq \big (\{e\}, \varnothing \big ),
\end{align*}
\end{minipage}
\begin{minipage}[t]{0.5\textwidth}
\begin{align*}
\big (\{f\},\{e\}\big ) &\preceq \big (\{f\}, \varnothing \big )\\
	\big (\{e\},\{f\}\big ) &\preceq \big (\{e\}, \varnothing \big )\\
	\big ( \{f\}\big ) & \preceq \big  (\{f\}, \varnothing \big ) \\
		\big ( \{e\}\big ) & \preceq \big  (\{e\}, \varnothing \big ) \\
	\big ( \{e,f\}\big ) &\preceq \varnothing \\
	\big ( \{e\}\big ) &\preceq \varnothing, \,\, \big ( \{f\}\big ) \preceq \varnothing
\end{align*}
\end{minipage}
and the trivial relations that  $\pi \preceq \pi$ for all $\pi \in \Piall$. Among the above relations, some are of the form $\pi_{\sed=\varnothing}(F) \preceq \pi$ for an ordered partition $\pi$ of a subset $A$ of $F$.
\end{example}

\subsection{Layered graphs and their graded minors} \label{sec:PreliminariesLayeredGraphs} A \emph{layered graph}  is the data of a graph $G=(V, E)$ and an ordered partition $\pi=(\pi_1 \dots, \pi_r, \pi_{\fin}) \in \Piall(E)$ on the set of edges $E$.

The elements $\pi_j$ are called \emph{layers}; those $\pi_j$ with $j\neq \fin$ are called \emph{unbounded}; the part $\pi_\fin$ is called the \emph{bounded} or the \emph{finitarty part} of $\pi$. 

The number $r$ is called the \emph{rank} of the layered graph $G$ and equals the rank of the ordered partition underlying unbounded edges. By an abuse of the notation, we use the same letter $G$ to denote the layered graph $(V, E, \pi)$. The genus of a layered graph is by definition equal to that of its underlying graph. 

The terminology is justified as follows. In the degeneration picture for metric graphs, where we also give an edge length to each edge, the parts $\pi_1, \dots, \pi_r$ correspond to those edges whose lengths go to infinity, with the same \emph{speed} in each layer, while the part $\pi_\fin$ contains those edges whose lengths remain bounded all the time. Note in particular that in this way, there is a distinction between the two partitions $(E, \varnothing)$ and $(E)$: in the former, the edges in $E$ have lengths which go to infinity with the same speed, while in the latter, all the edge lengths remain bounded. 

\smallskip

Let now $G=(V, E, \pi)$ be a layered graph with the ordered partition $\pi = (\pi_1, \dots, \pi_r, \pi_\fin)$ of $E$. We define the \emph{graded minors} of $G$ as follows. 

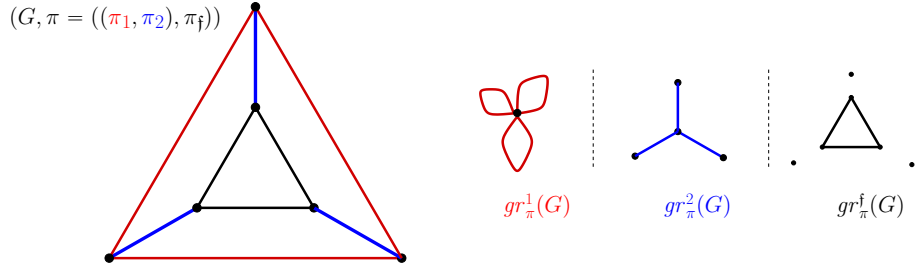
\begin{figure}[!t]
\centering
    \scalebox{.35}{\input{example1.tikz}}
\caption{A layered graph with an ordered partition $\pi=(\pi_\infty, \pi_\fin)$, $\pi_\infty=(\pi_1, \pi_2)$, and its corresponding graded minors.}
\label{fig:layered_graph}
\end{figure}

Consider the filtration $\filter^\pi_\bullet$ associated to $\pi$
\[ F_{0} :=\varnothing \subsetneq F_1 \subsetneq F_2 \subsetneq \dots \subsetneq F_r \subseteq F_{\infty+\fin} := E\] 
 defined by 
\[F_j := \bigsqcup_{1\leq i\leq j} \pi_i, \qquad j=1, \dots, r.\]
 For each $j \ge 1$, we let 
\[E^j_\pi: = E\setminus F_{j-1} = \pi_j \cup \pi_{j+1} \cup \dots \cup \pi_r \cup \pi_\fin \]
and naturally arrive at a decreasing filtration $\dfilter_\pi^\bullet$: 
\[ E_\pi^1 := E \supsetneq E_\pi^2 \supsetneq E_\pi^3 \subsetneq \dots E_\pi^r \supsetneq E_\pi^{r+1} := \pi_\fin.\]
This leads to a decreasing sequence of spanning subgraphs of $G$ 
\[G =: G^1_\pi \supseteq G^2_\pi \supsetneq G^3_\pi \supsetneq \dots \supsetneq G_\pi^{r} \supsetneq G_\pi^{r+1}= (V, \pi_\fin)\]
where for each $j=1, \dots, r$,  $G_\pi^j:=(V, E_\pi^j)$ has edge set $E_\pi^j$.

\smallskip

For each integer $j\in [r]$, the $j^{\mathrm{th}}$ \emph{graded minor} of $G$, denoted by $\grm{\pi}{j}(G)$, is obtained by contracting all the edges of $E_{\pi}^{j+1}$ 
in $G_\pi^j$, that is., 
\[\grm{\pi}{j}(G) := \contract{G^{j}_\pi}{E^{j+1}_\pi} = \contract{G[F_{j-1}^c]}{F_j^c} = \contract{G[\pi_j\cup \dots \cup\pi_r \cup \pi_\fin]}{\pi_{j+1} \cup \dots \cup \pi_r\cup \pi_\fin}.\]

It follows that $\grm{\pi}{j}(G)$ has edge set equal to $\pi_j$. (More precisely, we should view the elements of $\pi_j$ in the contracted graph, but we simplify and identify the two sets in the writing.) We denote by $V_\pi^j$ the vertex set of $\grm{\pi}{j}(G)$. Sometimes we will drop $\pi$ or $G$ and simply write $\grm{}{j}(G)$, $\grm{\pi}{j}$, $V^j$, \emph{etc.} if $\pi$ or $G$ is understood from the context. We denote by $\cont{j}:  V \to V^j_\pi$ the contraction map. 

We extend the definition to $\pi_\fin$ by setting $G_\fin$ to be the \emph{finitary graded minor} defined as the spanning subgraph of $G$ with vertex set $V$ and edge set $\pi_\fin$. In order to be more specific, we refer sometimes to the graded minors $\grm{\pi}{1}(G), \grm{\pi}{2}(G), \dots, \grm{\pi}{r}(G)$ as the \emph{infinitary graded minors} of $G=(V, E, \pi)$.

\smallskip

By the genus formula we proved in~\cite{AN}, we have 
\[\graphgenus = h^1_\pi + \dots + h^r_\pi + h^\fin_\pi\]
where $h^j_\pi$ is the genus of the graded minor $\grm{\pi}{j}(G)$, that is the sum of the genera  of the connected components of  $\grm{\pi}{j}(G)$.

\subsection{Admissible homology bases for layered graphs} \label{sec:admissible_basis_layered_graphs}
We recall in this section the important notion of admissible basis for the homology $H_1(G, \R)$ of a graph $G=(V, E)$ endowed with the ordered partition $\pi = (\pi_\infty, \pi_\fin)$, $\pi_\infty=(\pi_1, \dots, \pi_r)$, introduced in~\cite{AN}.

Let $G=(V, E, \pi)$ be a layered graph of genus $\graphgenus$, as above. Let $\filter^\pi_\bullet$ and $\dfilter_\pi^\bullet$ be the corresponding increasing and decreasing filtrations of $E$, 
\[\filter^\pi_{\bullet}: \quad F_0 = \varnothing \subsetneq F_1 \subsetneq F_2 \subsetneq \dots \subsetneq F_r\subseteq F_{\infty+\fin} = E\]
\[\dfilter_\pi^\bullet:\quad E^1_\pi =E \supsetneq E^2_\pi \supsetneq \dots \supsetneq E^r_\pi \supsetneq E^{r+1}_\pi =\pi_\fin,\]
with the associated decreasing sequence of spanning subgraphs,
\[G =: G^1_\pi \supset G^2_\pi \supset  \dots \supset G_\pi^{r} \supset G_\pi^{r+1}= (V,\pi_\fin)\]
and $\grm{\pi}{j}(G)$, $j\in [r] \cup\{\fin\}$ the corresponding graded minors of $G$.

From the contraction map $\cont{j}: G^{j}_\pi \to \grm{\pi}{j}(G)$, we obtain the map on homology $\proj_j : H_1\bigl(G^{j}_\pi, \Z\bigr) \to H_1\bigl(\grm{\pi}{j}(G), \Z\bigr)$, which is easily shown to be surjective.

We consider the partition 
\begin{equation} \label{eq:ordered_partition_edges}
[\graphgenus] = J^1_\pi \,\sqcup\, J^2_\pi \sqcup \dots \sqcup J^r_\pi \,\sqcup\,J^\fin_\pi
\end{equation}
into intervals $J^j_\pi$ of size $h^j_\pi$ given by 
\[J^j_\pi := \Bigl\{1+\sum_{i=1}^{j-1}h^i_\pi,\,\, 2+\sum_{i=1}^{j-1}h^i_\pi,\,\, \dots,\,\,\sum_{i=1}^{j}h^i_\pi\Bigr\}.\]

An \emph{admissible basis} for $G$ is a basis $\gamma_1, \dots, \gamma_h$ of $H_1(G, \Z)$ which verifies the following condition.

\begin{itemize}
\item For each $j=1, \dots, r, \fin$,  the following properties hold:
\medskip

\begin{itemize}
\item[$(i)$] For $k \in J^j_\pi$, the cycle $\gamma_k$ lies in the spanning subgraph $G^{j}_\pi$ of $G$. This means all the edges of $\gamma_k$ belong to $E^{j}_\pi \subseteq E$.

\smallskip

\item[$(ii)$] The collection of cycles $\proj_j(\gamma_k)$ with $k\in J^j_\pi$ forms a basis of the homology $H_1(\grm{\pi}{j}(G), \Z)$.
\end{itemize}
\end{itemize}

In~\cite{AN} we showed that every layered graph has an admissible basis.


\subsection{Marked graphs and counting function} We recall the definition of marked graphs \cite{AN}. Let $n$  be a non-negative integer. A \emph{graph with $n$} (\emph{labelled}) \emph{markings} is by definition a graph $G=(V, E)$ endowed with the \emph{marking function} $\marking: [n] \to V$ which associates the vertex $\marking(j) \in V$ to each label $j\in [n]$. 

To a given graph with $n$ marked points, one associates the \emph{counting function} $\countmarking: V \to \mathbb N \cup \{0\}$ 
which  takes value $\countmarking(v)$ at vertex $v$ equal to the number of labels placed at $v$, i.e.,  $\countmarking(v)$ is the number of elements $j \in [n]$ with $\marking(j)=v$.


\subsection{Stable graphs}  An \emph{augmented graph} is a graph $G=(V, E)$ endowed with a function $\genusfunction: V \to \mathbb N \cup\{0\}$ called the  \emph{genus function} defined on the set of vertices. The value $\genusfunction(v)$ is called the \emph{genus} of the vertex $v \in V$. The \emph{genus} of $G=(V, E, \genusfunction)$ is  
\[g:= h + \sum_{v\in V} \genusfunction(v)\]
where $h$ denotes the genus of the underlying graph $(V, E)$. 

An \emph{augmented graph with $n$ markings} is an augmented graph $(V, E, \genusfunction)$ endowed with a marking $\marking: [n] \to V$. 

\medskip

 In this paper, we sometime consider augmented graphs or augmented marked graphs which verify the \emph{stability condition} below. In this case, they are called \emph{stable graphs}. So by definition, a stable graph comes with a genus function on vertices.

A \emph{stable graph with $n$ markings} is a quadruple $G=(V, E, \genusfunction, \marking)$ with marking and genus functions  $\marking$ and $\genusfunction$, respectively,  which verifies the following \emph{stability condition}: 
 
 \begin{itemize}
\item for a vertex $v$ of genus zero,  $\deg(v) + \countmarking(v) \geq 3$, and
\item for a vertex $v$ of genus one,  $\deg(v) + \countmarking(v)\geq 1$.
\end{itemize}

 A \emph{stable graph}  is a triple $(V, E, \genusfunction)$ which verifies the stability condition above for the null counting function $\countmarking \equiv0$.
 
By a slight abuse of the notation, the notation $G$ is used both for the stable graph (with markings)  and the underlying graph $(V, E)$, thus understanding the genus function (and the marking) from the context.

 \medskip

The definition of deletion and contraction can be extended in a natural way to stable marked graphs. We refer to~\cite{AN} for a precise formulation. For a stable (marked) graph $G$ of genus $g$, and an edge $e \in E$, the stable graph $\contract{G}{e}$ obtained by contracting $e$ is again of genus $g$.

 We define layered stable graphs (with markings) as stable graphs $G=(V, E, \genusfunction, \marking)$ (with markings)  with the data of an ordered partition $\pi = (\pi_\infty=(\pi_1, \dots, \pi_r), \pi_\fin)$ on the edge set $E$.
 
 \subsection{Partial order on stable graphs of given genus}
Let $n$ be a non-negative integer. We define a partial order on the set of stable marked graphs of genus $g$ with $n$ markings as follows. We say $G \subface H$ for two stable graphs $G$ and $H$ (with markings) if $G$ can be obtained  from $H$ by a sequence of edge-contractions.

 \subsection{Partial order on layered stable (marked) graphs of given genus} 
 Let $H$ be a stable graph with markings and with layering $\pi_H \in \Piall(E(H))$. Let $G$ be a stable graph with markings and with layering $\pi_G \in \Piall(E(G))$.  We say $(G, \pi_G) \subface (H, \pi_H)$ if $G \subface H$ and $\pi_H \preceq \pi_G$ as ordered partitions, in the sense of Section~\ref{subsec:OrdPart}.

\subsection{Layered metric graphs} \label{ss:MetricGraphs} We briefly recall the definition of a metric graph. Suppose $G=(V,E)$ is a finite graph and let $l: E \rightarrow \mathbb R_{>0}$ be a length function, which assigns a positive real number $l(e)$ to every edge $e \in E$. To the pair $(G, l)$, we can associate a metric space $\mgr$ by assigning each edge a direction, identifying every edge $e\in E$ with a copy of the interval $\Ical_e= [0,l(e)]$ (the left and right endpoints corresponding to the initial and terminal vertices of the directed edge), by defining $\mgr$ as the result of identifying the ends of intervals corresponding to the same vertex $v$ (in the sense of a topological quotient). The topology on $\mgr$ is metrizable by the so-called {\em path metric}: the distance between two points $x,y \in\mgr$ is defined as the arc length of the shortest path connecting them.

A \emph{metric graph} is a compact metric space which can arise from the above construction for some pair $(G, l)$ of a graph $G$ and length function $l$. In this case,  $\mgr$ is called the \emph{metric realization} of the pair $(G, l)$. On the other hand, a pair $(G, l)$ whose metric realization is isometric to $\mgr$  is called a \emph{finite graph model} of $\mgr$.

 One can naturally define augmented metric graphs (with markings) and stable metric graphs (with markings)~\cite{AN}. Any stable metric graph has a minimal finite graph model.

 \medskip

Let $\mgr$ be a metric graph obtained as the metric realization of the pair $(G, l)$ consisting of a finite graph $G=(V, E)$ and an edge length function $l\colon E \to (0,+\infty)$. We define the {\em total length} $L(\mgr)$ of $\mgr$ as
\begin{equation} \label{eq:DefLength}
	L(\mgr)  := \sum_{e \in E} l(e).
\end{equation}
If the edge set $E$ is equipped with an ordered partition $\pi = (\pi_1, \dots, \pi_r, \pi_\fin)$, then we denote the {\em total length of the $i$-th graded minor} denoted by $L_i$, $i=1,\dots, r, \fin$, by
\begin{equation} \label{eq:DefGrLength}
 L_i(\mgr) :=  L(\mgr_\pi^i) = \sum_{e \in \pi_i} l(e)
\end{equation}
By definition, the equality $L(\mgr) = \sum_i L_i(\mgr)$ holds true.

\subsection{Tropical curves} We define  on layered metric graphs the notion of being \emph{conformally equivalent at infinity} that we simply call  \emph{the conformal equivalence relation} (\emph{at infinity} being assumed at all time) as follows. Two layered metric graphs $\mgr$ and $\mgr'$ are called conformally equivalent if the following three properties hold:
\begin{itemize}
\item $\mgr$ and $\mgr'$ have the same combinatorial type $(G, \pi)$ for some ordered partition $\pi=(\pi_1, \dots, \pi_r, \pi_\fin)$ with $r\in \N \cup\{0\}$.
\item for each layer $\pi_j$ of $\pi$, for $j\in[r]$, there is a positive number $\lambda_j >0$ such that
\[l_{|_{\pi_j}} = \lambda_j \, l_{|_{\pi_j}}'
\]
for the corresponding edge length functions $l, l' \colon  E(G) \to (0, + \infty)$. 

\item the two length functions $l$ and $l'$ coincide on the finitary part $\pi_\fin$, i.e., 
\[l_{|_{\pi_\fin}} = l_{|_{\pi_\fin}}'
\]
\end{itemize}
Saying it differently, the layered metric graphs $(\mgr', \pi')$ equivalent to some fixed $(\mgr, \pi)$ are simply obtained by multiplying the edge lengths for every layer $\pi_j$ for $j\in [r]$ by some constant $\lambda_j >0$.

\smallskip

A \emph{tropical curve} $\curve$ is a conformal equivalence class of layered metric graphs. Equivalently, we can view a tropical curve  as a pair $\curve= (\mgr, \pi)$ consisting  of an augmented metric graph $\mgr$  with an ordered partition $\pi=(\pi_1, \dots, \pi_r, \pi_\fin)$ on the edge set $E$ of a finite graph model $(G, l)$ of $\mgr$, $r\in \N\cup
\{0\}$, such that in addition, the normalization equality 
\[\sum_{e\in \pi_j} l_e =1, \quad j \in[r]\]
holds true. 

We call the set of edges $E \setminus \pi_\fin$ the part of the tropical curve which \emph{lives at infinity} and sometimes denote it by $E_\infty$. The \emph{sedentarity} of the tropical curve $\curve$ that we denote by $\sed(\curve)$ is the ordered partition $(\pi_1, \dots, \pi_r)$ of the part at infinity $E_\infty$ so that there are $r$ different layers of infinity, $\pi_1$ being above $\pi_2$, $\pi_2$ above $\pi_3$, etc.

\smallskip

We say a tropical curve $\curve$ \emph{ has full sedentarity} if all its edges  live at infinity, that is, if $E_\infty = E$.  The tropical curves introduced in~\cite{AN} all had full sedentarity. For the purpose of what we wish to do in this paper, all the possible sedentarities are needed, and considering them all at the same time will allow to define a \emph{compactified moduli space of tropical curves}.  

The notion of sedentarity introduced above can be viewed as a higher rank version of the current-in-use notion of sedenarity in tropical geometry, and leads to higher rank canonical compactifications of fans and extended cone complexes, extending to higher rank the construction which appeared for example in~\cites{ACP, AP-homology, AP, OR11}.
 We will elaborate on this later in Section~\ref{sec:tropical_moduli}.

\subsection{Layered metrized complexes and hybrid curves} \label{sec:hybrid_curves}  We recall the definition of \emph{metrized complexes} and \emph{hybrid curves} from~\cite{AB15} and~\cite{AN}.

A metrized curve complex $\mc$ consists of the data of 
\begin{itemize}
\item a finite graph $G=(V, E)$.
\item a metric graph $\mgr$ with a model $(G, l)$ for a length function $l \colon E \to (0, + \infty)$.  
\item for each vertex $v\in V$, a smooth projective complex curve $C_v$. 
\item for each vertex $v \in V$, a bijection $e \mapsto p^e_v$ between the edges of $G$ incident to $v$ (with loop edges counted twice) and a subset 
$\mathcal A_v = \{ p^e_v \}_{e \ni v}$ of $C_v(\C)$.
\end{itemize}
We use the same notation $\mc$ for the metrized complex and its \emph{geometric realization} that we define as follows.  As for metric graph geometric realization, for each edge $e\in E$, let $\Ical_e$ be an interval of length $\ell_e$.  For each vertex, by an abuse of the notation, let $C_v$ be the compact Riemann surface obtained as the analytification of the complex curve $C_v$. For each vertex $v$ of an edge $e$, we identify the corresponding extremity of $\Ical_e$ with the marked point $p^e_v$. This allows to identify $\mc$ as the disjoint union of the Riemann surfaces $C_v$ and intervals $\Ical_e$, for $v\in V$ and $e\in E$, quotiented by these identifications. We endow $\mc$ with the quotient topology.

\smallskip

To any metrized complex, we naturally associate the underlying graph, that we can endow with the genus function $\genusfunction: V \to \mathbb N\cup\{0\}$ which to the vertex $v$ associates $\genusfunction(v)$ defined as the genus of $C_v$. The metric realization of this augmented graph with the edge length function $l$ is called the underlying metric graph of $\mc$. Similarly, we associate to any metrized complex, the underlying semistable curve by gluing the curves $C_v$ along the marked points associated to the same edges.

\smallskip

A \emph{layered metrized complex} is a pair $(\mc, \pi)$ consisting of a metrized complex and an ordered partition $\pi=(\pi_1, \dots, \pi_r)$ of the edge set $E$ in the underlying metric graph $G=(V, E)$, and $r\in \N$ provided that $E \neq \varnothing$, so that all the parts $\pi_1, \dots, \pi_r$ are non-empty. In analogy with the definition of layered metric graphs above, it is possible to enrich $\pi$ to an ordered partition $(\pi_1, \dots, \pi_r, \pi_\smallc)$ of the full metrized complex by setting $\pi_\smallc:= \bigsqcup_{v\in V} C_v$. In this way, we can view the layers $\pi_1, \dots, \pi_r$ as the \emph{metric graph unbounded} or \emph{infinitary} layers of the layered metrized complex and see $\pi_\smallc$ as the \emph{complex finitary} or \emph{bounded} part of the layered metrized complex. 

To any layered metrized complex $(\mc, \pi)$ with $\pi = (\pi_1, \dots, \pi_r)$, as above, we associate the \emph{underlying layered metric graph} $(\mgr, \pi)$ where $\mgr$ is the metric graph underlying $\mc$ and we view $\pi$ as an ordered partition of full sedentarity by putting $\pi_\fin=\varnothing$, that is, $\pi = (\pi_1, \dots, \pi_r, \pi_\fin = \varnothing)$. 

We define the \emph{conformal equivalence relation (at infinity)} on layered metrized complexes as follows. Two layered metrized complexes $\mc$ and $\mc'$ are called conformally equivalent if there exists an isomorphism of the semistable curves underlying $\mc$ and $\mc'$ so that under this isomorphism, they have the same underlying graph $G=(V, E)$ (this will be automatic), the same ordered partition $\pi =(\pi_1, \dots, \pi_r)$ on $E$, and for each layer $\pi_j$ of $\pi$, for $j \in [r]$, there is a positive number $\lambda_j >0$ such that
\[l_{|_{\pi_j}} = \lambda_j \, l_{|_{\pi_j}}'
\]
for the corresponding edge length functions $l, l' \colon  E(G) \to (0, + \infty)$ in the two metrized complexes. In other words, the two underlying layered metric graphs $(\mgr, \pi)$ and $(\mgr', \pi')$ are conformally equivalent.

\begin{figure}[!t]
\centering
    \scalebox{.35}{\input{example7.tikz}}
\caption{Example of a hybrid curve with the corresponding infinitary and finitary worldsheets.}
\label{fig:worldsheets}
\end{figure}
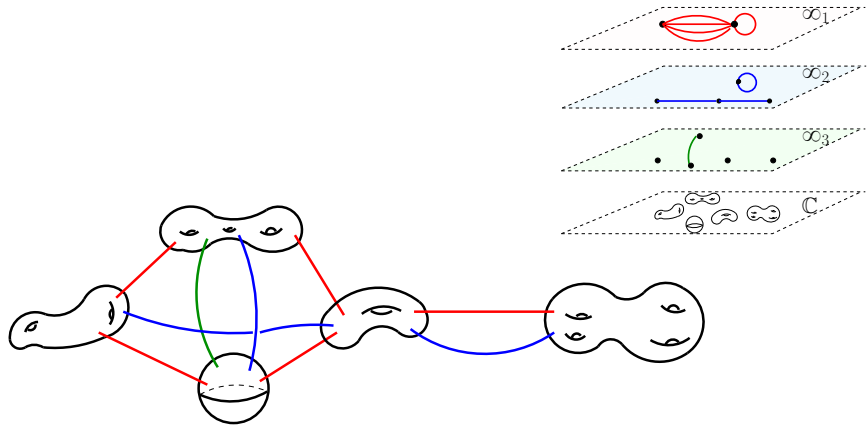

A \emph{hybrid curve}  $\curve$ is a conformal equivalence class of layered metrized complexes.  We define the \emph{sedentarity} of the hybrid curve $\curve$ to be the ordered partition $\pi=(\pi_1, \dots, \pi_r)$ of its edge set.

Any hybrid curve $\curve$ naturally gives rise to the corresponding tropical curve $\curve$ of full sedentarity $\sedbis(\curve^\trop)=\pi$ obtained by taking the underlying augmented metric graph $\mgr$ of $\mc$ endowed with the layering given by the ordered partition $\pi=(\pi_1, \dots, \pi_r, \pi_\fin=\varnothing)$  of the edge set $E$ of the graph $G$.

\subsection{Continuous families of measured spaces} \label{ss:FamiliesMeasuredSpaces} 
In what follows, we will frequently consider families of geometric objects $(X_t)_{t \in B}$, where each geometric object $X_t$ is equipped with a measure $\mu_t$. In particular, we will need a notion of continuity of such families $(X_t, \mu_t)_{t \in B}$ of measured spaces. The following definition turns out to be useful for our purposes.

\smallskip

Let $X$ and $B$ be topological spaces together with a surjective, continuous map $\pr \colon X \to B$. Suppose further that each fiber $X_t := \pr^{-1}(\{t\})$, $t \in B$, is equipped with a Radon measure $\mu_t$. If for any continuous function $f :X \to \R$, the function $F : Y \to \R$ defined by \emph{integration along fibers},
\begin{align*}
&F(t) := \int_{X_t} f\rest{X_t}(x) \, d \mu_t(x), \qquad  t \in B,
\end{align*}
is continuous on $B$, then we call $(X_t, \mu_t)_{t \in B}$ a {\em continuous family of measured spaces} over the base space $B$, and call the measures $\mu_t$, $t \in B$, a {\em continuous family of measures} on $X/B$.

\newpage


\part{Tropical compactifications and log maps}\label{part:compactifications}

The aim of this part is to introduce the polyhedral geometry materials which arise naturally in higher rank valuation theory, in the study of degenerations of multiparameter families of algebraic varieties. In the first section we introduce higher rank compactifications of polyhedral fans. This will allow us in the second part to define the moduli space of higher rank tropical curves. We then introduce global and stratumwise log maps. These maps are crucial in our study as they allow to compare the geometry of objects living at infinities of the corresponding moduli spaces to the geometry of nearby objects, living in the finitary part, both in the tropical and in the hybrid settings.

\section{Higher rank canonical compactifications of fans} \label{sec:higher_rank_compactifications}

Canonical compactifications of fans and polyhedral complexes have been playing a central role in recent developments of tropical geometry, in particular, in connection with Hodge theory~\cite{AP-homology,AP,AP-CS}.

Our aim in this section is to introduce higher rank analogues of these canonical compactifications. Later on, in Parts~\ref{part:TropicalCurves} and~\ref{part:HybridCurves}, we will use these compactifications to define the \emph{moduli space of higher rank tropical curves}, as well as \emph{tropical} and \emph{hybrid log maps}.

 As we will see, the existing notion of canonical compactification in the literature  can be viewed within the general framework of this section.

\subsection*{Basic terminology} In what follows, we denote by $N$ a free $\Z$-module of finite rank $d$ and by $M=N^\dual = \hom(N, \Z)$ the dual of $N$. The corresponding vectors spaces are denoted by $N_\R$ and $M_\R$, respectively, so we have $M_\R = N_\R^\dual$. For a given rational polyhedral cone $\sigma$ in $N_\R$, we use the notation $N_{\sigma, \R}$ for the real vector subspace of $N_\R$ generated by elements of $\sigma$, and  use $N^\sigma_\R := \rquot{N_\R}{N_{\sigma, \R}}$ for the quotient. Since $\sigma$ is rational, we get natural lattices $N^\sigma$ and $N_\sigma$ in $N^\sigma_\R$ and $N_{\sigma, \R}$, respectively, which are both of full rank.

\smallskip

 For the ease of reading, we adopt the following convention  from~\cite{AP-homology, AP}. We use $\sigma$ (or any other face of $\Sigma$) as a superscript where referring to the quotient of some space by $N_{\sigma, \R}$ or to the elements related to this quotient. In contrast, we use $\sigma$ as a subscript for subspaces of $N_{\sigma,\R}$ or for elements associated to these subspaces.

\smallskip

We define $\eR := \R \cup \{\infty\}$ to be the extended real line endowed with the topology induced by that of $\R$ and a basis of open neighborhoods of infinity given by intervals $(a, \infty]$ for $a\in \R$. The addition of $\R$ extends naturally to $\eR$ in a natural way, by setting $\infty + a = \infty$ for all $a \in \eR$, and gives $\eR$ the structure of a topological monoid called the monoid of \emph{tropical numbers}. The submonoid of non-negative tropical numbers is $\eR_+ := \R_+ \cup\{\infty\}$ and is endowed with the induced topology. Both monoids  carry a natural $\R_+$-semimodule structure.

\smallskip

\subsection{Canonical compactification of fans and their conical stratification} We first recall the definition of the canonical compactification $\cancomp \Sigma$ of a fan $\Sigma$. A more detailed presentation is given in~\cites{AP-homology, AP, OR11}.

\smallskip

Let $\Sigma$ be a fan in $N_\R$. For any cone $\eta$, its \emph{dual cone} $\eta^\vee$ is defined as
\[\eta^\vee := \Bigl\{m \in M_\R \:\bigm|\: \langle m, a \rangle \geq 0 \:\textrm{ for all } a \in \eta\Bigr\}. \]
Clearly, the orthogonal complement 
\[\eta^\perp := \Bigl \{ m \in M_\R \: | \: \langle m, a \rangle = 0 \:\textrm{ for all } a \in \eta \Bigr\}\]
 is contained in the dual cone $\eta^\vee$.

We begin by compactifying each cone separately. The \emph{canonical compactification} or the {\em extended cone} $\cancomp\eta$ of a cone $\eta$ in $\Sigma$ is given by 
\[\cancomp\eta :=\hom_{\RMod}(\eta^\vee, \eR_+), \]
where the $\hom$ is taken in the category of $\R_+$-semimodules. The original cone $\eta$ can be naturally identified with $\hom_{\RMod}(\eta^\vee, \R_+)$, which is the subset of $\cancomp\eta$ consisting of the morphisms with image contained in $\R_+$.

Equipped with the topology of pointwise convergence of morphisms (w.r.t.\ the topology on $\eR_+$), $\cancomp\eta$ is a compact Hausdorff space. The induced topology on $\eta $ coincides with the original Euclidean topology on $\eta$ and hence $\cancomp \eta$ is a compactification of $\eta$.

\smallskip

Each extended cone $\cancomp\eta$ contains a special point denoted by $\infty_\eta$, which, when viewed as an element of $\hom_{\RMod}(\eta^\vee, \eR_+)$, is represented by the map from $\eta^\vee$ to $ \eR_+$
\[
\infty_\eta(m) := \begin{cases} 0 & \text{if } m \in \eta^\perp \\ + \infty    & \text{if } m \in \eta^\vee \setminus \eta^\perp. \end{cases}
\]
Note that for the cone $\zerocone :=\{0\}$, we have $\infty_{\zerocone} = 0$.  

\smallskip

For an inclusion of cones $\delta \subseteq \eta$ in $\Sigma$, we have the dual inclusion $\eta^\vee \subseteq \delta^\vee$, and this yields an inclusion of compactifications $\cancomp \delta \subseteq \cancomp \eta$. Under this identification, $\cancomp \delta$ coincides with the topological closure of $\delta$ in $\cancomp \eta$.

\smallskip

Finally, we define the canonical compactification $\cancomp\Sigma$ as the union of the extended cones $\cancomp\eta$, $\eta\in \Sigma$, where $\cancomp \delta$ is identified with a subset of $\cancomp \eta$ whenever $\delta \subseteq \eta$ in $\Sigma$. Formally, $\cancomp\Sigma$ is defined as the quotient
\[
\cancomp\Sigma :=  \bigsqcup_{\eta \in \Sigma} \rquot{\cancomp\eta}{\sim},
\]
where $"\sim"$ denotes the above identifications, and is endowed with the quotient topology. Note that each extended cone $\cancomp\eta$ naturally embeds as a subspace of $\cancomp\Sigma$.

\smallskip

We recall that the space $\cancomp\Sigma$ lives in a partial compactification $\TP_\Sigma$ of $N_\R$, called the \emph{tropical toric variety associated to $\Sigma$}. Its definition is a variant of the above one, where the morphisms may also take negative values. Namely, for a cone $\eta$ in $\Sigma$, introduce
\[\canchart\eta:=\hom_{\RMod}(\eta^\vee,\eR),\] 
and note that $N_\R \subset \canchart \eta$ (since $\hom_{\RMod}(\eta^\vee, \R)\simeq N_\R$).

\smallskip

The tropical toric variety $\TP_\Sigma$ is then obtained by gluing the spaces $\canchart\eta$ via the inclusions $\canchart\delta \subseteq\canchart\eta$ for $\delta \subseteq \eta$ in $\Sigma$. We endow $\TP_\Sigma$ with the induced quotient topology. We note that $\TP_\Sigma$ can be defined alternatively  as the tropicalization of the toric variety $\P_\Sigma$ associated to $\Sigma$.

\smallskip

Since $\eR_+ \subset \eR$, we get an inclusion  $\cancomp\eta \subset \canchart\eta$ and hence an embedding $\cancomp\Sigma \subseteq \TP_\Sigma$. This embedding identifies $\cancomp \Sigma$ as the closure of $\Sigma$ in $\TP_\Sigma$.

\smallskip

The tropical toric variety $\TP_\Sigma$ admits a natural {\em stratification} in terms of the faces of $\Sigma$. Namely, note first that each space $\canchart \eta$ comes with an addition and multiplication by scalars, induced by $\eR$. This allows to define the space $N^\eta_{\infty,\R}:= N_\R + \infty_\eta \subseteq \canchart\eta$. Moreover, we have $N^\delta_{\infty, \R} \simeq N^\delta_\R :=\rquot{N_\R}{N_{\delta,\R}}$ and, in particular, $N^\zerocone_{\infty,\R}=N_\R$. Clearly, $\canchart\eta$ is stratified into a disjoint union of the subspaces $N^\delta_{\infty,\R}$,  $\delta \subseteq \eta$. Altogether, the partial compactification $\TP_\Sigma$ is stratified into a disjoint union of the spaces $N^\eta_{\infty,\R} \simeq N^\eta_\R$, $\eta \in \Sigma$. 

\smallskip

The stratification of $\TP_\Sigma$ induces a natural {\em conical stratification} on the canonical compactification $\cancomp \Sigma$ as follows. Let $\delta \subseteq \eta$ be a pair of faces in $\Sigma$. We denote by $\keg^\delta_{\eta}$ the subset of $\cancomp\eta$ defined as 
\[\keg^\delta_{\eta} :=\cancomp\eta \cap N^{\delta}_{\infty, \R} = \infty_\delta + \eta = \{\infty_\delta + x \mid x\in \eta\}.\]
In the above definition, the sum is taken in $\canchart \delta$ where both $\eta$ and $\infty_\delta$ live via the embeddings $\eta = \hom_{\R_+}(\eta^\vee , \R_+) \subset \hom_{\R_+}(\delta^\vee , \eR) = \canchart\delta$ and $\infty_\delta\in \cancomp\delta \subset \canchart \delta$. By definition, $\keg^\delta_{\eta}$ is a cone in the $\R$-vector space $N^\delta_{\infty, \R}$.

\smallskip

Under the isomorphism $N^\delta_{\infty,\R}\simeq N^\delta_\R$, the cone $\keg^\delta_{\eta}$ is identified precisely with the projection of $\eta$ in the linear space $N^\delta_\R \simeq N^\delta_{\infty,\R}$, which we denote by $\eta^\delta$. We stress that, although $\eta^\delta\simeq \keg^\delta_\eta$, these two cones live in two different spaces, the first one in $N^\delta_\R$ and the second in $N^\delta_{\infty, \R}$.

We denote by $\inn\keg^\delta_{\eta}$ the relative interior of $\keg^\delta_{\eta}$.  For pair of faces $\delta \subseteq \eta$ in $\Sigma$, define $N^\delta_{\infty, \eta,\R}$ as the projection of $N_{\eta, \R}$ in $N^\delta_{\infty, \R}$, that is
\[
N^\delta_{\infty, \eta,\R}:= N_{\eta, \R} + \infty_\delta.
\]
Note that $N^\delta_{\infty, \eta,\R} \simeq \rquot{N_{\eta, \R}}{N_{\delta, \R}}$.

\smallskip

We have the following description of the canonical compactification.

\begin{prop}\label{prop:con-strat} Let $\Sigma$ be a fan in $N_\R$.
\begin{itemize}
\item The canonical compactification $\cancomp\Sigma$ is a disjoint union of \emph{(}open\emph{)} cones $\inn\keg^\delta_{\eta}$ with $\delta \subseteq \eta$ both faces of $\Sigma$. The linear space generated by the cone $\keg^\delta_{\eta}$ is the real vector space $N^\delta_{\infty, \eta,\R}$.
\item For a pair of faces $\delta \subset \eta$ in $\Sigma$, the closure $\cancomp \keg^\delta_{\eta}$ of $\keg^\delta_{\eta}$ in $\cancomp\Sigma$ is the union of all the \emph{(}open\emph{)} cones $\inn\keg^{\delta'}_{\eta'}$ with $\delta \subseteq \delta' \subseteq \eta' \subseteq \eta$.
\end{itemize}
\end{prop}

The open cones $\inn \keg^\delta_{\eta}$ give the (open) strata of the conical stratification of $\cancomp\Sigma$. 

\medskip

The \emph{sedentarity} of a point $x$ in $\cancomp\Sigma$ is defined to be the face $\delta$ of $\Sigma$  such that $x$ lies in the corresponding linear space $N^\delta_{\infty, \R}$.

\subsection{Higher rank canonical compactifications: Motivation} We now extend the above picture to higher rank. The outcome of the constructions in this section is another compactification 
\[
\Sigma \subseteq  \cancomp \Sigma^\trop.
\]
of the fan $\Sigma$ which naturally lives above the canonical compactification $\cancomp \Sigma$. We call this fan the \emph{absolute canonical compactification} or the \emph{tropical compactification} of $\Sigma$. We stress however that this is one of the possible higher rank compactifications, the finest one, appearing in a \emph{bi-indexed} class of higher rank canonical compactifications  $\cancomp\Sigma^{\tropr{a,b}}$of $\Sigma$ (see Section~\ref{sec:botany}). 
\begin{remark}[Refinement relative to behavior both at zero and infinity]
The higher rank canonical compactifications in this paper are all defined relative to the behavior of points in the fan at \emph{infinity}. Remark~\ref{remark:refinement} sketches a refinement  including the behavior relative to \emph{zero}. 
\end{remark}

The idea behind the definition is to refine the canonical compactification $\cancomp \Sigma$ by further separating the special points $\infty_\delta$ living at its boundary. To motivate the definition, consider now a sequence of elements $\phi_t$, $t \in \N$, in the cone $\eta$ of $\Sigma$. We wish to give a meaningful sense to the refined convergence of (a subsequence of) the $\phi_t$'s in the (to-be-defined) refined compactification $\cancomp\Sigma^\trop$. 

\smallskip

We first use the duality to identify $\eta = \hom_\RMod(\eta^\vee, \R_+)$. Extracting a subsequence if necessary, we can suppose that the sequence converges pointwise to an element $\phi_{\infty} \colon \eta^\vee \to \eR_+$. We will repeatedly use the following lemma.
\begin{lem} Let $\tau$ be a polyhedral cone in $N_\R$, and consider a morphism $\theta\colon \tau \to \eR_+$. The locus of points $x$ in $\tau$ where the function $\theta$ takes a finite value is a face of $\tau$. The same statement holds for the locus of points where the function takes value zero.
\end{lem} 
\begin{proof} The set $A$ of points where $\theta(x)<\infty$ is a closed convex set in $\tau$. Moreover, if $A$ intersects the relative interior of a face of $\tau$, it should contain the face entirely. It follows that $A$ is a face of $\tau$. The same argument applies to the second part.
\end{proof}
Applying the first part of the lemma to the limit morphism $\phi_{\infty}$, we obtain a face $\tau$ of $\eta^\vee$. By duality, this face is of the form $\tau = \eta^\vee \cap \delta^{\perp}$ for a face $\delta \subseteq \eta$, and the limit map $\phi_{\infty}$ takes value $\infty$ at all the points $\eta^\vee \setminus \tau $. This implies that $\phi_{\infty}$ is uniquely determined as an element of the dual to the cone $\tau = \eta^\vee \cap \delta^\perp$. Now, we identify $\tau = \eta^\vee \cap \delta^\perp \simeq (\eta^\delta)^\vee$ and realize  $\phi_{\infty}$ as a point living in the cone $\eta^\delta$.  This is precisely the identification $C^\delta_\eta \simeq \eta^\delta$ of the previous section. 

\smallskip 

Let $\delta_{\infty+\fin}$ be the subface of $\eta$ such that the point $\phi_\infty$ lives in the relative interior of $\delta_{\infty+\fin}^\delta \subseteq \eta^\delta$.

At this step, we turn to the points where the value of the limit morphism $\phi_{\infty}$ is $\infty$. Replacing $\phi_t$ with a subsequence, if necessary, we can find a sequence of functions $L_r (t), L_{r-1}(t), \dots, L_1(t)$, $r\in \N$, with the following properties for $j=1, \dots, r$:
\begin{itemize}
\item the ratio $L_{j}(t)/L_{j+1}(t)$ converges to infinity. For $j=r+1$, we set $L_{r+1}(t)\equiv 1$, 

\smallskip

\item the renormalized sequence $\frac 1{L_j(t)}\, \phi_t$ converges to a non-zero, continuous morphism $\phi_{\infty, j} \colon \eta^\vee \to \eR_+$. In particular, $\phi_{\infty, j}$ takes all the finite values in $\R_+$.
\smallskip

\item For each ray $\varrho$ of $\eta^\vee \setminus \delta_{\infty+\fin}^\perp$, there is some $j$ such that $\phi_{\infty, j}\rest\varrho$ is neither zero nor identically $\equiv + \infty$.
\end{itemize}

Applying now the previous lemma to the maps $\phi_{\infty, j}$, we infer the existence of a non-zero face $\delta_j \subseteq \delta_{\infty+\fin} \subseteq \eta$ with the property that $\phi_{\infty, j}$ takes finite values on $\eta^\vee \cap \delta_j^\perp$, and values $\infty$ on $\eta^\vee \setminus \delta_j^\perp$, $j=0, \dots, r$.

Using the properties of $L_j$ and $\phi_{\infty, j}$ listed above we infer the following:

\begin{itemize}
 
 \item $\phi_{\infty,r} =\phi_\infty$ and it lies in the interior of the face $\keg_{\delta_{\infty+\fin}}^{\delta} \simeq \delta_{\infty+\fin}^\delta$ of $\keg_\eta^\delta$.
 
 \item We have $\delta_{j} \subsetneq \delta_{j+1}$, and $\phi_{\infty, j}$ is uniquely determined as a morphism from $(\delta_{j+1}^{\delta_{j}})^\vee \to \R_+$ with only non-zero values taken on non-zero points of this cone. This identifies $\phi_{\infty, j}$ with an interior point of the cone $\delta_{j+1}^{\delta_j}$.
  \item $\phi_{\infty, j}$ for every $0 \le j \le r-1$ is well-defined modulo multiplication by a positive scalar. 

\end{itemize}

The outcome of the above discussion is thus a flag $\pi$ of faces in $\eta$, $r\in \N\cup\{0\}$, 
\[\pi\colon\qquad \delta_0 = \zerocone \subsetneq \delta_1 \subsetneq \delta_2 \subsetneq \dots \subsetneq \delta_r = \delta_\infty \subseteq \delta_{\infty+\fin} \subseteq \eta. \]
In the above inclusions, $\delta_{\infty+\fin}$ is the face of $\eta$ such that $\phi_\infty$ lies in the relative interior of $\delta_{\infty+\fin}^{\delta_r}$. 
 
For a pair of cones $\zeta \subset \delta$, we introduce the \emph{projectivized cone} 
\[\sigma_{\delta}^\zeta :=\rquot{\bigl(\delta^\zeta  - 0\bigr)}{\R_{>0}}  = \phom_{\RMod}\bigl(({\delta}^{\zeta})^\vee, \R_{+} \bigr) = \rquot{\hom_{\RMod}\bigl(({\delta}^{\zeta})^\vee, \R_{+} \bigr) -0}{\R_{>0}}. \]

Using the above discussion, the refined limit of the sequence $\phi_t$ lives in the space 
\begin{equation} \label{eq:DefTropicalCones} \inn\keg^{\trop}_\pi := \inn\sigma_{\delta_{1}}^{\delta_0} \times \inn\sigma_{\delta_{2}}^{\delta_1} \times \dots \times \inn\sigma_{ \delta_{r}}^{\delta_{r-1}} \times \inn \keg^{\delta_r}_{\delta_{\infty+\fin}}
\end{equation}
for $\pi$ a flag of faces of $\eta$ as above, for some $r\in \N$.

\subsection{Tropical compactification $\cancomp\Sigma^\trop$: Definition}The above discussion motivates the following definition. Let $\Sigma$ be a fan in $N_\R$ and denote by $\mathscr{F}(\Sigma)$ the set of all flags of faces in $\Sigma$ of the form
\[\pi\colon\qquad \delta_0 = \zerocone \subsetneq \delta_1 \subsetneq \delta_2 \subsetneq \dots \subsetneq \delta_r \subseteq \delta_{\infty+\fin}, \qquad r\in \N \cup\{0\}. \]
The integer $r \in \N\cup\{0\}$ is called the {\em rank} of the flag $\pi$ and we also set $\delta_\infty := \delta_r$  in this context. For any flag $\pi \in \mathscr{F}(\Sigma)$ of rank $r$, its {\em sedentarity} $\sed(\pi) = \pi_\infty$ is the flag
\begin{equation} 
\pi_\infty\colon\qquad \delta_0 = \zerocone \subsetneq \delta_1 \subsetneq \delta_2 \subsetneq \dots \subsetneq \delta_r = \delta_\infty.
\end{equation}
We introduce the \emph{tropical compactification} of the fan $\Sigma$ as 
\begin{equation} \label{eq:DefSigmaTrop}
\cancomp\Sigma^{\trop}: = \bigsqcup_{\pi \in \mathscr{F}(\Sigma)} \inn\keg^{\trop}_\pi.
\end{equation}
 The fan $\Sigma$ corresponds to flags of trivial sedentarity (equivalently, of rank $r = 0$),
\[
\Sigma = \bigsqcup_{\substack{\pi \in \mathscr{F}(\Sigma) \colon \\ \sed(\pi) = (\zerocone)}} \inn\keg^{\trop}_\pi \:  \subseteq \: \cancomp\Sigma^{\trop}.
\]
There is a natural projection map
 \begin{align*}
\pr \colon \cancomp\Sigma^{\trop} &\longrightarrow   \cancomp \Sigma\\
\thy =(x_1, \dots, x_r, t_\fin) &\longmapsto t_\fin
\end{align*}
which maps each stratum $\inn\keg^{\trop}_\pi$ in $ \cancomp\Sigma^{\trop}$ onto the stratum $\inn \keg^{\delta_\infty}_{\delta_{\infty+\fin}}$ in $ \cancomp \Sigma$ by projecting onto the last coordinate. In what follows, we equip $\cancomp\Sigma^{\trop}$ with a topology such that this projection map is continuous.

\subsection{Tropical compactification: Topology}  \label{ss:TropTopDef}  We begin by specifying the topological relation between the different strata $\inn\keg^{\trop}_\pi$ in \eqref{eq:DefSigmaTrop}. Recall that in the canonical compactification $\cancomp \Sigma$ (see Proposition \ref{prop:con-strat}), a stratum $\inn \keg_{ \eta}^{ \delta}$ is contained in the closure of  $\inn \keg_{ \eta'}^{ \delta'}$ exactly when $\delta' \subseteq \delta \subseteq  \eta \subseteq \eta'$. We introduce the following generalization to flags.

 For two flags $\pi$ and $\tilde \pi$,
 \begin{align*}
 \pi&\colon\qquad \delta_0 = \zerocone \subsetneq \delta_1 \subsetneq \delta_2 \subsetneq \dots \subsetneq \delta_r \subseteq \delta_{\infty+\fin}, \qquad r\in \N \cup\{0\},\\
 \tilde \pi &\colon\qquad \tilde \delta_0 = \zerocone \subsetneq \tilde\delta_1 \subsetneq \tilde\delta_2 \subsetneq \dots \subsetneq \tilde\delta_{\tilde r} \subseteq \tilde \delta_{\infty+\fin}, \qquad \tilde r\in \N \cup\{0\}
 \end{align*}  
 in $\mathscr{F}(\Sigma)$, we say that $\pi$ is {\em finer} than $ \tilde  \pi $ and write $ \tilde  \pi  \preceq \pi$, if
\begin{center} $\sedbis( \tilde  \pi ) \subseteq \sedbis(\pi)$ (as sets) \qquad and \qquad $\delta_{\infty+\fin} \subseteq  \tilde  \delta_{\infty+\fin}$.
\end{center}
The first condition means that $\sedbis(\pi)$ is a refinement of $\sed(\tilde \pi)$, that is, there are integers $1 \le n(1) < n(2) < \dots < n(\tilde r) \le r$ such that
\begin{equation} \label{eq:FormCoarserStrata}
	\tilde \delta_j = \delta_{n(j)}, \qquad j= 1, \dots, \tilde r,
\end{equation}
where $r$ and $\tilde r$ denote the ranks of $\pi$ and $\tilde \pi$, respectively. The relation "$\preceq$" is clearly a partial order on the set of flags $\mathscr{F}(\Sigma)$. In the topology of $\cancomp \Sigma^\trop$, a stratum $\inn \keg^{\trop}_{\pi}$ is contained in the closure of $\inn \keg^{\trop}_{\tilde \pi}$ exactly when $\tilde\pi \preceq \pi$. 

\smallskip

From the perspective of dual cones, the order "$\preceq$" has the following interpretation. For two flags $\pi$, $\tilde \pi$ in $\mathscr{F}(\Sigma)$ with $ \tilde  \pi \preceq   \pi$, we obtain an ascending flag of faces of $\eta:=  {\tilde  \delta_{\infty+\fin}}$,
\[
\zerocone =  \delta_0  \subsetneq   \delta_1 \subsetneq \dots \subsetneq  \delta_r \subseteq { \tilde  \delta}_{\infty+\fin} = \eta.
\]
By duality, this gives a decreasing flag of faces of the corresponding dual cone $\eta^\vee$,
\[
\eta^\vee = { \tau}_0  \supsetneq { \tau}_1 \supsetneq \dots \supsetneq { \tau}_r  \supseteq \tilde { \tau}_{\infty+\fin} = \eta^\vee \cap \eta^\perp,
\]
where  $\tau_i = \eta^\vee  \cap {\delta}_i^\perp$ for $i = 1, \dots, r$. In particular, the dual cone $\eta^\vee$ is decomposed into
\begin{equation} \label{eq:Dual1}
\eta^\vee = (\eta^\vee \setminus \tau_r) \cup { \tau}_r =  \bigcup_{i=1}^r ({\tau}_{i-1} \setminus { \tau}_i ) \cup { \tau}_r.\end{equation}
Since the flag $\tilde \pi$ satisfies $\tilde \pi \preceq \pi$, it yields a coarser decomposition of the dual cone $\eta^\vee$,
\begin{equation}  \label{eq:Dual2}
\eta^\vee = (\eta^\vee \setminus \tau_{n(\tilde{r})}) \cup \tau_{n(\tilde{r})} =  \bigcup_{j=1}^{\tilde{r}} ({\tau}_{n(j-1)} \setminus { \tau}_{n(j)} ) \cup\tau_{n(\tilde{r})}.
\end{equation}

The points of the cones $ \inn\keg^{\trop}_\pi$ and $\inn\keg^{\trop}_{\tilde \pi}$ can be interpreted as morphisms on ${\eta^\vee}$ with particular properties w.r.t.\ the partitions \eqref{eq:Dual1} and \eqref{eq:Dual2}. Namely, let $\thy = (x_1, \dots, x_r, t_\fin) \in \inn\keg^{\trop}_{ \pi}$ and $\shy = (y_1, \dots, y_{\tilde r}, s_\fin)$ in $\inn\keg^{\trop}_{\tilde \pi}$.

\smallskip

(i) By the inclusions ${{\delta}_\fin}^{{\delta}_r} \subseteq {\eta}^{{ \delta}_r} \subseteq \cancomp {\eta}$ and ${\eta}^{{\tilde \delta}_{\tilde r}} \subseteq \cancomp {\eta}$, the last components $t_\fin$ and $s_\fin$ are identified with morphisms $t_\fin, s_\fin\colon \eta^\vee \to \eR_+$ such that
\begin{align*}
t_\fin^{-1}\big (+ \infty\big ) = \eta^\vee \setminus { \tau}_{r} \qquad \text{and}  \qquad s_\fin^{-1}\big (+ \infty\big ) = \eta^\vee \setminus { \tau}_{n(\tilde r)}.
\end{align*}

Notice that the inclusion $\eta^\vee \setminus \tau_{n(\tilde r)} \subseteq \eta^\vee \setminus \tau_{r}$ holds true.

\smallskip

(ii) The representatives $x_i' $ and $y_j'$ of the projective components $x_i \in  \inn\sigma_{\delta_{i}}^{\delta_{i-1}}$ and  $y_j' \in \inn\sigma_{\tilde \delta_{j}}^{\tilde \delta_{j-1}}$ correspond to morphisms $x_i', y_j' \colon \eta^\vee \to \eR_+$ satisfying
\begin{align*}
 x_i' \equiv 0 \text{ on } {\tau}_{i} 
  , \qquad  0 < x_i' < + \infty \text{ on } {\tau}_{i-1} \setminus { \tau}_i , \qquad x_i' \equiv + \infty \text{ on } \eta^\vee \setminus {\tau}_{i-1}
\end{align*}
for every $i= 1,\dots, r,$ and 
\begin{align*}
 y_j' \equiv 0 \text{ on } {\tau}_{n(j)} 
  , \qquad  0 < y_j' < + \infty \text{ on } {\tau}_{n(j-1)} \setminus { \tau}_{n(j)} , \qquad y_j' \equiv + \infty \text{ on } \eta^\vee \setminus {\tau}_{n(j-1)}
\end{align*}
for every $j= 1, \dots, \tilde r$. To avoid the ambiguity in choosing representatives, we will fix a point $a_i$ in each set $\tau_{i-1} \setminus \tau_i$, $i = 1, \dots, r$. Then the elements $x_i' / x_i'(a_i) \in \cancomp \eta$ and $y_j'/y_j'(a_i) \in \cancomp \eta$ only depend on $x_i$ and $y_j$.

\medskip 

After these preparations, we proceed with the definition of the topology on $\cancomp\Sigma^{\trop}$. For each point $\thy \in \cancomp\Sigma^{\trop}$, we will specify a neighborhood base $\mathcal{U}(\thy) = (U)_{U \in \mathcal{U}(\thy)}$. Every set $U$ in the neighborhood base $\mathcal{U}(\thy)$ is of the form
\[
U = \bigsqcup_{\tilde \pi \preceq \pi} U_{\tilde \pi},
\]
where $U_{\tilde \pi}$ is a subset of the (open) stratum $\inn \keg_{\tilde \pi}^\trop$ defined as follows: Fix neighborhoods $V_\fin$ and $V_i$ in $\cancomp \eta$ of the points $t_\fin$ and $x_i' / x_i'(a_i)$,  $i=1, \dots, r$. The set $U_{\tilde \pi}$ consists of all the points $\shy = (y_1, \dots, y_{\tilde r}, s_\fin)$ in $\inn \keg_{\tilde \pi}^\trop$ which satisfy the following set of conditions:

\begin{itemize}
\item [(i)]  $s_\fin \in V_\fin$,
\smallskip

\item [(ii)] for each face $\delta_i$ of the flag $\pi$ with ${\tilde \delta}_{\infty} \subsetneq \delta_i$,
\begin{equation}
	\frac{s_\fin}{s_\fin(a_i)} \in V_i, 
	\end{equation}
\smallskip

\item [(iii)] and for each two consecutive faces $\tilde \delta_{j-1} \subseteq \tilde \delta_j$, $j= 1, \dots, \tilde r$ of $\tilde \pi$ and each face $\delta_i$ of $\pi$ with $\tilde \delta_{j-1} \subsetneq \delta_i \subseteq \tilde \delta_j$,
\begin{equation}
	\frac{y_j'}{y_j'(a_i)} \in V_i. 
\end{equation}
\end{itemize}

The neighborhood base $\mathcal{U}(\thy)$ consists of all sets $U = \bigsqcup_{\tilde \pi \preceq \pi} U_{\tilde \pi}$ of the above form, that is, they are obtained when $V_\fin$ and $V_i$ range over all neighborhoods of $t_\fin$ and $x_i'/x_i'(a_i)$.

\begin{thm} 
There is a unique topology on $\cancomp\Sigma^{\trop}$ such that $\mathcal{U}(\thy)$ is a neighborhood base for any $\thy \in \cancomp\Sigma^{\trop}$. Moreover, this topology is independent of the choice of the points $a_i$ as above.
\end{thm}
\begin{proof} The set systems $\mathcal{N}(\thy):= \bigl\{N \subseteq \cancomp\Sigma^{\trop} | \, N \supseteq U \text{ for some } U \in \mathcal{U}(\thy)\bigr\}$, $\thy \in \cancomp\Sigma^{\trop}$, are easily seen to satisfy the axioms of a neighborhood system and the claim follows.
\end{proof}

We have the following description of the convergence of sequences in $\cancomp \Sigma^\trop$.  When applied to a sequence $(\varphi_n)_{n \in \N}$ in the fan $\Sigma$, we recover the idea of the preceding section. 
\begin{prop} \label{prop:convseq}
Let $\thy \in \cancomp \Sigma^\trop$ and assume $\thy=(x_1,\dots, x_r, t_\fin) \in \inn \keg_\pi^{\trop}$ for the flag $\pi \in \mathscr{F}(\Sigma)$. Suppose that $(\shy_n)_n$ is a sequence in $\cancomp \Sigma^\trop$. Then:

\medskip
(a) If $(\shy_n)_n$ converges to $\thy$ in $ \cancomp \Sigma^\trop$, then almost all $\shy_n$ belong to strata $\inn \keg_{\tilde \pi}^{\trop}$ of flags $\tilde \pi \in \mathscr{F}(\Sigma)$ with $\tilde \pi \preceq  \pi$.

\medskip
(b) Assume that $\shy_n = (y_{1,n}, \dots, y_{\tilde r, n}, s_{\fin, n})$ defines a sequence in the stratum $ \keg_{\tilde \pi}^\trop$ for some fixed flag $\tilde \pi \in \mathscr{F}(\Sigma)$ with $\tilde \pi \preceq \pi$. In particular, $\tilde \pi$ is of the form \eqref{eq:FormCoarserStrata} and, as above, the points $\thy$ and $\shy_n$ can be described via the dual cone of  $\eta := {\tilde \delta_{\infty+\fin}}$.

Then, $\shy_n$ converges to $\thy$ in $\cancomp \Sigma^\trop$ if and only if the following conditions hold:

\begin{itemize}
\item [(i)] $s_{\fin, n}$ converges to $t_\fin$ in $\cancomp{\Sigma}$. In particular, $\lim_{n \to \infty} s_{\fin, n}(a) = + \infty$ for all $a \in \eta^\vee \setminus { \tau}_{r}$.
\item [(ii)] If $j=1, \dots, \tilde r$ and $\delta_i$ is a face in the flag $\pi$ with $\tilde \delta_{j-1} \subsetneq \delta_i \subseteq \tilde \delta_j$, then 
\[
	\lim_{n \to \infty} \frac{y_{j, n}'}{y_{j, n}'(a)} =  \frac{x_i'}{x_i'(a)} \qquad \text{in } \cancomp \eta
\]
for every point $a \in \tau_{i -1} \setminus \tau_i$. In particular, if $\tilde \delta_{j-1} \subsetneq \delta_i \subsetneq \delta_{\hat i} \subseteq \tilde \delta_j$, then
\[
\lim_{n \to \infty} \frac{y_{j, n}'(a)}{y_{j, n}'(\hat a)} = + \infty
\]
for all points $a \in \tau_{i-1} \setminus \tau_i$ and $\hat a\in \tau_{\hat i-1} \setminus \tau_{\hat i}$.

\item [(iii)] If $\delta_i$ is a face in the flag $\pi$ such that $\tilde \delta_\infty \subsetneq \delta_i$, then
\[
	\lim_{n \to \infty} \frac{s_{\fin, n}}{s_{\fin, n}(a)} =  \frac{x_i'}{x_i'(a)} \qquad \text{in } \cancomp \eta
\]
for every point $a \in \tau_{i -1} \setminus \tau_i$. In particular, if $\tilde \delta_\infty \subsetneq \delta_i \subsetneq \delta_{\hat i}$, then
\[
\lim_{n \to \infty} \frac{s_{\fin, n}(a)}{s_{\fin, n}(\hat a)} = + \infty
\]
for all points $a \in \tau_{i-1} \setminus \tau_i$ and $\hat a\in \tau_{\hat i-1} \setminus \tau_{\hat i}$. 
\end{itemize}
 \end{prop}
\begin{proof}
The proposition is clear from the definition of the topology on $ \cancomp \Sigma^\trop$.
\end{proof}

The next proposition clarifies the basic properties of $\cancomp\Sigma^{\trop}$ and the topological relation between the strata (analogous to Proposition~\ref{prop:con-strat} for $\cancomp \Sigma$).

 \begin{prop}\label{prop:con-strat} Let $\Sigma$ be a fan in $N_\R$ and let $\cancomp\Sigma^{\trop}$ be its tropical compactification.
 
 \begin{itemize}
 \item The space $\cancomp\Sigma^{\trop}$ is a second countable, compact Hausdorff space. The convergent sequences are described in Proposition~\ref{prop:convseq}.
 \smallskip
 
 \item  The projection map $\pr \colon \cancomp\Sigma^{\trop} \to \cancomp\Sigma$, $\thy \mapsto t_\fin$, is continuous. It maps each stratum $\inn\keg^{\trop}_\pi$ in $ \cancomp\Sigma^{\trop}$ onto the stratum $\inn \keg^{\delta_\infty}_{\delta_{\infty+\fin}}$ in $ \cancomp \Sigma$ by projecting to the last coordinate. 
 \smallskip
 
 \item The induced topology on the stratum $ \inn\keg^{\trop}_\pi = \inn\sigma_{\delta_{1}}^{\delta_0} \times  \dots \times \inn\sigma_{ \delta_{r}}^{\delta_{r-1}} \times \inn \keg^{\delta_r}_{\delta_{\infty+\fin}}$ of a flag $\pi \in \mathscr{F}(\Sigma)$ coincides with the product topology. Moreover, the closure of the stratum $\inn\keg^{\trop}_\pi $ in $\cancomp \Sigma^\trop$ is the union of all the strata $\inn\keg^{\trop}_{\tilde \pi}$ with $\pi \preceq \tilde \pi$.
 \end{itemize}
\end{prop}
\begin{proof}
The above claims are easily deduced from the corresponding properties of $\cancomp \Sigma$.
\end{proof}

\subsection{Sedentarity stratification of $\Tropcomp\Sigma$} For any point $x$ in $\cancomp \Sigma^\trop$ lying on the stratum $\inn \keg_\pi^\trop$, with the flag 
 \[\pi \colon\qquad \delta_0 = \zerocone \subsetneq \delta_1 \subsetneq \delta_2 \subsetneq \dots \subsetneq \delta_r \subseteq \delta_{\infty+\fin}, \]
 we define the \emph{sedentarity of $x$} denoted by $\sed(x) = \pi_\infty$ as the sedentarity of $\pi$, that is
 \[\pi_\infty\colon\qquad \delta_0 = \zerocone \subsetneq \delta_1 \subsetneq \delta_2 \subsetneq \dots \subsetneq \delta_r = \delta_\infty. \]
A point $x$ is called of \emph{full sedentarity} if $\sed(x) = \pi$, \ie, if $\delta_{\infty+\fin} =\delta_\infty$.

The locus of points in $\cancomp\Sigma^\trop$ of given sedentarity $\pi_\infty = \zerocone\subsetneq \delta_1 \subsetneq \delta_2 \subsetneq \dots \subsetneq \delta_r$ is denoted by $\cancomp\Sigma^\trop_{\pi_\infty}$. It is isomorphic to the product 
\[\cancomp\Sigma^\trop_{\pi_\infty} = \inn\sigma_{\delta_{1}}^{\delta_0} \times \inn\sigma_{\delta_{2}}^{\delta_1} \times \dots \times \inn\sigma_{ \delta_{r}}^{\delta_{r-1}} \times \Sigma^{\delta_r}.\]
We extend the definition to $\pi_\infty =\zerocone$ by 
\[\cancomp\Sigma^\trop_{\zerocone} :=\Sigma\]
and call this the \emph{finite stratum} of the compactification. 

The set of points of full sedentarity in each stratum $\cancomp\Sigma^\trop_{\pi_\infty}$ corresponds to the product 
\[\inn\sigma_{\delta_{1}}^{\delta_0} \times \inn\sigma_{\delta_{2}}^{\delta_1} \times \dots \times \inn\sigma_{ \delta_{r}}^{\delta_{r-1}} \times \zerocone.\]

 \subsection{Botany of higher rank canonical compactifications}\label{sec:botany} As we discuss next, the space $\cancomp \Sigma^\trop$ appears as the "finest one"  in a certain class of compactifications  of a fan $\Sigma$. This class appears to contain several polyhedral analogues of the existing spaces in the literature. Although we will not need this point of view in the sequel to this paper, we include it here because of its conceptual importance and of its potential utility in further developments of tropical and hybrid geometry. \smallskip

Consider a point $\thy$ in $ \cancomp \Sigma^\trop$ which belongs to the stratum $\inn \keg_\pi^\trop$ of a flag $\pi \in \mathscr{F}(\Sigma)$ of rank $r$. Then its $r+1$ coordinates $\thy= (x_1, \dots, x_r, t_\fin)$ describe the asymptotic behavior of an approximating morphism sequence $(\phi_t)_t$ for $t_\fin$ in $\cancomp \Sigma$ on $r+1$ faces of the dual cone $\delta_{\infty+\fin}^\vee$. The idea behind the space $ \cancomp \Sigma^{\tropr{p, q}}$ that we will define below is that it can keep track of the behavior of morphisms on at most $p+q$ faces: namely the $p$ with fastest growth and the $q$ with slowest growth.

\medskip
 
This concept is formalized as follows. Fix $p,q  \in \N\cup\{0\}$ with $1 \le p+q \le d+1$ (recall that $N_\R \simeq \R^d$). Consider first the case $q\neq 0$.

There is a natural forgetful operation $\forget$ on $\cancomp\Sigma^{\trop}$ which only remembers the first $p$ and the last $q$ coordinates of a point $\thy \in\Sigma^{\trop}$. Namely, if $\thy \in \inn \keg_{\pi}^\trop $ for a flag $\pi \in \mathscr{F}(\Sigma)$ of rank $r$, we define
\[ \forget(\thy) = \forget(\underbrace{x_1, x_2, \dots, x_r, t_\fin}_{r+1 \text{ components}}) := \begin{cases} (\underbrace{x_1, \dots, x_{p}}_{p \text{ components}}, \underbrace{x_{r - q + 2},  \dots, x_{r}, t_\fin}_{q \text{ components}}) & \text{if $r + 1 > p + q $} \\ \thy  &\text{if $r + 1 \le p+q $}\end{cases}.
\]
As a set, we introduce the space $ \cancomp \Sigma^{\tropr{p,q}}$ as the image
\[
	\cancomp \Sigma^{\tropr{p, q}} := \forget(\cancomp \Sigma^\trop).
\]
Analogous to \eqref{eq:DefSigmaTrop}, it can be stratified into a disjoint union 
\begin{equation*}
 \cancomp\Sigma^{\tropr{p,q}}  = \bigsqcup_{\pi \in \mathscr{F}^{(p,q)} (\Sigma) } \inn\keg^{\, (p,q)}_\pi,
\end{equation*}
where $ \mathscr{F}^{(p,q)}(\Sigma)$ is the set of all flags $\pi \in \mathscr{F}(\Sigma)$ of rank $r \le p+q$ and 
\begin{equation} \label{eq:pqStrata}
\inn\keg^{\, (p,q)}_\pi := \begin{cases}
\inn \keg_{ \pi}^\trop & \text{if $r \le p+q -1$} \\
\prod_{\substack{j \in [r] \\ j \neq p+1}}
\inn\sigma_{\delta_{j}}^{\delta_{j-1}}   \times \inn \keg^{\delta_r}_{\delta_{\infty+\fin}} &\text{if $r  = p+q$} \end{cases}.
\end{equation}

\smallskip

 The remaining case $q = 0$ corresponds to looking only at the $p$ faces with fastest growth. In this case, the definition is slightly different: the set $\mathscr{F}^{(p, 0)}(\Sigma) $ consists of all flags $\pi \in \mathscr{F}(\Sigma)$ of rank $r \le p -1$ and all flags $\pi \in  \mathscr{F}(\Sigma)$ of length rank $r = p$ with $\delta_p = \delta_{\infty+\fin}$, that is , $\pi \colon \zerocone \subsetneq \delta_1 \subsetneq \dots \subsetneq \delta_p= \delta_{\fin}$. Accordingly, \eqref{eq:pqStrata} is changed to $\inn\keg^{\, (p,q)}_\pi := \prod_{j \in [p]} \inn\sigma_{\delta_{j}}^{\delta_{j-1}}$ in the second case.

\smallskip

Notice that there is a natural forgetful map from $\mathscr{F}(\Sigma)$ to $\mathscr{F}^{(p,q)}(\Sigma)$. By an abuse of the notation, we will denote this map by $\forget \colon \mathscr{F}(\Sigma) \to \mathscr{F}^{(p,q)}(\Sigma)$ as well. With this convention, we arrive at the equality
\[
	\inn\keg^{\,(p,q)}_{\forget(\pi)} = \forget( \inn\keg^{\trop}_{\pi}).
\]

In particular, we have $\cancomp \Sigma^{\tropr{p,q}} = \cancomp{\Sigma}^\trop$ whenever $p+q = d+1$. Moreover, if $p' \le p$ and $q' \le q$, then there is a surjective forgetful map from $\cancomp \Sigma^{\tropr{p, q}}$ to $\cancomp \Sigma^{\tropr{p', q'}}$ . Altogether, we obtain a commutative diagram of maps (see Figure~\ref{fig:CompactificationTower}).

\begin{center}
\begin{figure}[h!] 
\begin{tikzpicture}[scale=0.8, every node/.style={scale=0.9}]

\foreach \x in {1,..., 2}{
	\node at (2*\x, 0) {$\cancomp \Sigma^{\scaleto{\tropr{\x, 0}}{8 pt}}$};
	\draw[->>]  (2*\x, 1.25) -- (2*\x, 0.75) ;
	\draw[->>]  (2*\x + 1.25, 0) -- (2*\x + 0.75, 0) ;
}

\foreach \x in {0,..., 2}{
	\node at (2*\x, 2) {$\cancomp \Sigma^{\tropr{\x, 1}}$};
	\draw[->>]  (2*\x, 3.25) -- (2*\x, 2.75) ;
	\draw[->>]  (2*\x + 1.25, 2) -- (2*\x + 0.75, 2) ;
}

\foreach \x in {0,..., 1}{
	\node at (2*\x, 4) {$\cancomp \Sigma^{\tropr{\x, 2}}$};
	\draw[->>]  (2*\x + 1.25, 4) -- (2*\x + 0.75, 4) ;
		\draw[->>]  (2*\x, 5.25) -- (2*\x, 4.75) ;
}

	\node at (0, 9) {$\cancomp \Sigma^{\tropr{0, d+1}}$};
	\node[right] at (0.75, 9) {$= \cancomp \Sigma^\trop$};

	\node at (0, 7) {$ \cancomp \Sigma^{\tropr{0, d}}$};
	\draw[->>]  (0, 9 - 0.75) -- (0, 9-1.25) ;
	
	\node at (2, 7) {$\cancomp \Sigma^{\tropr{1, d}} $};
	\node[right] at (2.75, 7) {$= \cancomp \Sigma^\trop$};

	\draw[->>]  (1.25, 7) -- (0.75, 7) ;

	\draw[->>, dotted]  (2, 7 - 0.75) -- (2, 7-2.25) ;
	\draw[->>, dotted]  (0, 7 - 0.75) -- (0, 7-2.25) ;
	\node at (9, 0) {$\cancomp \Sigma^{\tropr{d+1, 0}}$};
	\node[right] at (9.75, 0) {$= \cancomp \Sigma^\trop$};

	\node at (7, 0) {$ \cancomp \Sigma^{\tropr{d, 0}}$};
	\draw[->>]  (9 - 0.75, 0) -- (9-1.25, 0) ;
	
	\node at (7, 2) {$\cancomp \Sigma^{\tropr{d, 1}} $};
	\node[right] at (7.75, 2) {$= \cancomp \Sigma^\trop$};

	\draw[->>]  (7, 1.25) -- (7, 0.75 ) ;

	\draw[->>, dotted]  (7 - 0.75, 2) -- (7-2.25, 2) ;
	\draw[->>, dotted]  (7 - 0.75, 0) -- (7-2.25, 0) ;

\end{tikzpicture}
\caption{The spaces $\cancomp \Sigma^{\tropr{p, q}}$ and the corresponding forgetful maps.}\label{fig:CompactificationTower}
\end{figure}
\end{center}

Each space $\cancomp \Sigma^{\tropr{p,q}}$ carries a natural topology, which formalizes the above idea. In terms of sequences, this topology can be described as follows. For two flags $\pi$, $ \tilde  \pi $ in $\mathscr{F}^{(p,q)}(\Sigma)$, we say that $\pi$ is {\em finer} than $ \tilde  \pi $ and write $ \tilde  \pi  \preceq \pi$, if there are flags $ \tilde  \pi^\ast  \preceq \pi^\ast$ in $\mathscr{F}(\Sigma)$ with
\[
\pi = \forget(\pi^\ast) \qquad \text{and} \qquad \tilde \pi = \forget(\tilde \pi^\ast)
\]
where $\forget$ is the corresponding forgetful map. 

One checks that this defines a partial order on $\mathscr{F}^{(p,q)} (\Sigma)$. A sequence $(\shy_n)_n$ of points in a stratum $\inn\keg^{\, (p,q)}_{\tilde \pi}$ can converge only to points $\thy$ in finer strata $\inn\keg^{\,( p,q)}_{ \pi}$, that is, with $\tilde \pi \preceq \pi$. Moreover, the sequence $(\shy_n) = (y_{1,n}, \dots, s_{\fin, n})_n$ converges to such a point $\thy = (x_1, \dots, t_\fin)$ exactly when the following conditions hold. Consider again the case $q\neq 0$ first.

\begin{itemize}
\item [(i)] $s_{\fin, n}$ converges to $t_\fin$ in $\cancomp{\Sigma}$. 

\smallskip

\item [(ii)] Clearly, each of the other (at most $p+q-1$ many) components $x_i$ of $\thy$ corresponds to two consecutive faces $\delta_{i-1} \subsetneq \delta_i$ of $\pi$. Since $\tilde \pi \preceq \pi$, there are two consecutive faces $\tilde \delta_{j-1} \subsetneq \tilde \delta_j$ of $\tilde \pi$ with $\tilde \delta_{j-1} \subseteq \delta_{i-1} \subsetneq \delta_i   \subseteq \tilde\delta_j$ (here we allow $j = \tilde r+1$ and set $\tilde\delta_{\tilde r + 1} := \tilde\delta_{\fin}$). The elements $x_i$ and $y_{j,n}$ can be seen as (equivalence classes of) morphisms on the extended dual cone $\cancomp \eta$ for $\eta := \tilde \delta_{\infty+\fin}$. Denoting by $x_i'$ and $y_{j,n}'$ the representatives, we require that
\[
	\lim_{n \to \infty} \frac{y_{j, n}'}{y_{j, n}'(a)} =  \frac{x_i'}{x_i'(a)} \qquad \text{in } \cancomp \eta
\]
for every point $a \in  \eta^\vee  \cap {\delta}_{i-1}^\perp \setminus  {\delta}_{i}^\perp$ (for $j = \tilde r+1$, we set $y_{j, n} := s_{\fin, n}$).
\end{itemize}

\smallskip

In the case $q=0$, for the space $\cancomp \Sigma^{\tropr{p,0}}$, it might occur that $\shy_n =(y_{1,n}, \dots, s_{\fin, n})$ has a finite component $s_{\fin, n}$, but the point $\thy = (x_1, \dots, x_p)$ not. The condition (i) then is replaced by $s_{\fin, n} \to + \infty$ on $\eta^\vee \setminus \delta_{\infty+\fin}^\perp$ (with $\eta := \tilde{\delta_{\infty+\fin}}$) and says that indeed the morphisms $s_{\fin, n}$ go to infinity on the $p$ faces represented by $\thy$. If both $s_{\fin, n}$ and $\thy$ have no finite components, the condition (i) is simply dropped.

\smallskip

The topology on $\cancomp \Sigma^{\tropr{p, q}}$ can be introduced analogously $\cancomp \Sigma^\trop$ and we omit the details. The following basic properties can be verified.

\begin{prop}\label{prop:con-strat} Let $\Sigma$ be a fan in $N_\R$ and $p, q \in \N \cup\{0\}$ with $1 \le p+q \le d+1$. 
 \begin{itemize}
 \item The space $\cancomp\Sigma^{\tropr{p,q}}$ is a second countable, compact Hausdorff space.
 \smallskip
 
  \item  The forgetful map $\forget \colon \cancomp\Sigma^{\trop} \to \cancomp\Sigma^{\tropr{p,q}}$ is continuous. It maps each stratum $\inn\keg^{\trop}_\pi$ of the tropical compactification $ \cancomp\Sigma^{\trop}$ onto the stratum $\inn \keg^{ \,(p, q)}_{\forget(\pi)}$ in $ \cancomp \Sigma^{\tropr{p,q}}$.
 \smallskip
 
 \item The topology on $\cancomp\Sigma^{\tropr{p,q}} = \forget(\cancomp \Sigma^\trop)$ coincides with the quotient topology induced by the forgetful map $\forget \colon \cancomp \Sigma^\trop \to \cancomp \Sigma^{\tropr{p,q}}$.
 
 \smallskip
 
  \item All maps in the above commutative diagram (see Figure~\ref{fig:CompactificationTower}) are continuous.
 
 \end{itemize}
\end{prop}
 
\begin{remark}
Since all the spaces $\cancomp \Sigma^{\tropr{p, q}}$ are obtained as quotients of $\cancomp \Sigma^\trop$, the latter can be considered as the "finest" of all these compactifications.
\end{remark}

\begin{remark}[A refinement of higher rank compactifications]\label{remark:refinement}
One can actually refine the above constructions by including \emph{vanishing layers} which are defined in the same way as the infinitary layers. That is we will have a second rank taking into account the speed of vanishing for the polyhedral coordinates. 
This viewpoint is particularly adapted to the moduli space setting. Indeed, the moduli space $\mgtrop{\grind{G}}$ of tropical curves of a given combinatorial type $G$ that we will define later is not compact, because of the impossibility to have zero edge lengths in this moduli space. However, adapting the possibility of adding "vanishing layers", this combinatorial moduli space becomes compact. From the geometric point of view, related to the results we will prove in the second part, this space has the same layered expansion for solutions of the Poisson equation, and for the convergence of one-forms. Moreover, it allows to speak about the convergence of the differentials also on vanishing edges. 
\end{remark}

The special interest in the spaces $\cancomp \Sigma^{\tropr{p, q}}$ stems from the fact that they provide a uniform framework for a treatment of canonical compactifications, in the realm of tropical geometry, including some specific constructions that have appeared in the literature before. Namely,

\begin{itemize}
\item The space $\cancomp \Sigma^{\tropr{0,1}} = \cancomp \Sigma$ is the original canonical compactification that we described in the beginning of this section. 
\smallskip

\item The space $\cancomp \Sigma^{\tropr{1,0}}$ can be viewed as a polyhedral analog of the constructions of Boucksom and Jonsson~\cite{BJ17} for complex manifolds using  dual complexes.

\smallskip

\item More generally, the spaces $\cancomp \Sigma^{\tropr{p,0}}$, $1 \le p \le d$, should be regarded as higher rank polyhedral analogs of the construction of~\cite{BJ17}.
\smallskip

\item The spaces $\cancomp \Sigma^{\tropr{p,1}}$, $1 \le p \le d$, are polyhedral analogs of the hybrid spaces  ${B}^{{\hybr{p}}}$ constructed in \cite{AN}.

\smallskip
\item The last line in Figure~\ref{fig:CompactificationTower} is the polyhedral analogue of the tower of compactifications which appeared in~\cite{AN}.
\end{itemize}


\subsection{The tame topology} \label{ss:TropicalTameTopology}
In what follows, it will be convenient to equip the higher rank compactification $\cancomp \Sigma^\trop$ of a fan $\Sigma$ with another, slightly finer, topology. Recall that, broadly speaking, the convergence of points $t \in \Sigma$ to a point $\thy$ in a rank $r$ stratum $\inn\keg_\pi^{\trop}$ means that $t$ goes to infinity with $r$ different speeds on $r$ faces of some dual cone. In applications, it turns out sometime necessary to quantify this layered behavior further. That is, instead of simply saying that on some faces the values approach infinity faster than on others, we would like to capture {\em how much faster} this happens.

\medskip

In order to formalize this, we introduce the following notion. Fix a family $\tameclass$ of functions $\Psi \colon \R_+ \to \R_+$ such that each $\Psi$ satisfies the following properties:
\begin{itemize}
\item [(i)] $\Psi$ is continuous and
\item [(ii)] $\lim_{\lambda \to \infty} \Psi(t) = +\infty$ and 
\item [(iii)] there is a continuous function $\psi \colon \R_+ \to \R_+$ with $\Psi(\lambda) \le \psi (\frac{\lambda}{\lambda'}) \Psi(\lambda')$ for all $\lambda$, $\lambda' \in \R_+$. 
\end{itemize}

\begin{example} \label{polynomial families}
Every function $\Psi_m(\lambda) = \lambda^m$, $m>0$, satisfies the above conditions (one can take $\psi =\max\{1, \Psi_m\}$). In particular, the monomials $\Psi_m$ of integral degree $m$ bounded by a given  number $N$ form a natural family of functions that we denote by $\tameclass_N$. Note that setting $N = \infty$ yields the family $\tameclass_\infty = \{ \lambda^m | \, m \ge 1\}$ of all monomials.
\end{example}

The idea of the following construction is to define a new topology on $\cancomp \Sigma^\trop$, by taking into account the family $\tameclass$. Naively speaking, we would like to refine the original topology on $\cancomp \Sigma^\trop$ by ensuring that the layers at infinity are, in a certain sense, separated from each other by the growth functions $\Psi$ in $\tameclass$. Namely, the convergence of points $t \in \Sigma$ to a limit $\thy$ in a stratum $\inn\keg_\pi^{\trop}$ requires additionally that on each dual face, $t$ approaches infinity faster than every function $\Psi\in \tameclass$ of the values on lower infinitary faces. This property is more carefully explained in Proposition~\ref{prop:convseqtame}.

\smallskip

Given two flags $\tilde \pi =(\tilde \delta_1, \dots,\tilde \delta_{\tilde r}, \tilde\delta_{\infty + \fin})$ and $\pi =(\delta_1, \dots, \delta_r, \delta_\fin)$  in $\mathscr{F}(\Sigma)$, we say that $\pi$ {\em tamely refines} $\tilde \pi$, if $\tilde \pi \preceq \pi$ and $\delta_j =  \tilde  \delta_j$ for all $j$ in $[\tilde r] = [\tilde r(\tilde \pi)]$. Equivalently, $\pi$ is obtained from $\tilde \pi$ by changing the last pieces of $\tilde \pi$ by replacing
\[
	\tilde \delta_{\tilde r}  \subseteq \tilde \delta_{\infty + \fin} \qquad \to \qquad \delta_{\tilde r}  = \tilde \delta_{\tilde r} \subsetneq  \delta_{\tilde{r}+1} \subsetneq \dots \subsetneq \delta_{r-1} \subsetneq \delta_{ r} \subseteq \delta_{\infty + \fin}
\]
with the condition that $\delta_{\infty + \fin} \subseteq  \tilde \delta_{\infty + \fin}$.

Notation as in Section~\ref{ss:TropTopDef}, recall that the original topology on $\cancomp\Sigma^{\trop}$ was given by specifying for each point $\thy \in \cancomp\Sigma^{\trop}$ a neighborhood base $\mathcal{U}(\thy) = (U)_{U \in \mathcal{U}(\thy)}$. For a point $\thy = (x_1, \dots, x_r, t_\fin)$ in the stratum $\inn \keg_{\pi}^\trop$ of some flag $\pi$, every set $U$ in the neighborhood base $\mathcal{U}(\thy)$ was of the form
\[
U = \bigsqcup_{\tilde \pi \preceq \pi} U_{\tilde \pi},
\]
where $U_{\tilde \pi}$ was a subset of the (open) stratum $\inn \keg_{\tilde \pi}^\trop$, constructed from neighborhoods $V_\fin$ and $V_i$ in $\cancomp \eta$ of the points $t_\fin$ and $x_i' / x_i'(a_i)$,  $i=1, \dots, r$ with $\eta = \tilde \delta_{\infty + \fin}$.

\smallskip

In order to define a neighborhood base around $\thy \in \inn \keg_{\pi}^\trop$ in the topology induced by $\tameclass$, we slightly modify the above sets $U_{\tilde \pi}$. Namely, fix a finite subset $S \subseteq \tameclass$ and some number $A \in \N$. For any flag $\tilde \pi$ such that $\pi$ tamely refines $\tilde \pi$, consider the subset $V_{\tilde \pi} \subseteq \inn \keg_{\tilde \pi}^\trop$ given by
\begin{align*}
V_{\tilde \pi} = \Bigl \{ &\shy = (y_1, \dots, y_{\tilde r}, s_\fin) \in \inn \keg_{\tilde \pi}^\trop \st  \,  \frac{s_\fin(a_i)}{\Psi (s_\fin(a_{\hat i}))} > A \\ &\text{for all } \Psi \in S \text{ and faces $\delta_i$, $\delta_{\hat i}$ of $\pi$ with $\tilde \delta_r \subsetneq \delta_i \subsetneq \delta_{\hat i} \subseteq \delta_\infty$}   \Bigr \}.
\end{align*}
In the above definition, the last component $s_\fin$ of $\shy$ is viewed as a morphism on the dual cone $\eta^\vee = (\tilde \delta_{\infty+ \fin})^\vee$ and $a_i$, $i = 1, \dots, r$, are the points fixed in Section~\ref{ss:TropTopDef}.
Combining the previous sets $U_{\tilde \pi}$ with the new sets $V_{\tilde \pi}$, we introduce

\[
V := \bigsqcup_{\tilde \pi\colon \pi \text{ tamely refines } \tilde \pi} U_{\tilde \pi} \cap V_{\tilde \pi}.
\]

The neighborhood base $\mathcal{U}(\thy)$ around a point $\thy$ consists of all sets $V $ of the above form. That is, they are obtained when the neighborhoods $V_\fin$ and $V_i$ defining $U_{\tilde \pi}$ range over all neighborhoods of $t_\fin$ and $x_i'/x_i'(a_i)$, $A$ over all integers, and $S$ over all finite subsets of $\tameclass$.

\smallskip

Analogous to the original topology on the higher rank compactification $\cancomp\Sigma^{\trop}$, we arrive at the following result.

\begin{thm} 
Given a family of functions $\tameclass$ as above, there is a unique topology on $\cancomp\Sigma^{\trop}$ such that $\mathcal{U}(\thy)$ is a neighborhood base for any $\thy \in \cancomp\Sigma^{\trop}$ (which is independent of the choice of the points $a_i$ as above).
\end{thm}

We call the above topology the {\em $\tameclass$-tame topology} on $\cancomp \Sigma^\trop$. The higher rank compactification $\cancomp \Sigma^\trop$ equipped with a $\tameclass$-tame topology is easily shown to be a second countable Hausdorff space. However, we stress that the resulting space is not compact in general.

Moreover, the $\tameclass$-tame topology is finer than the original topology on $\cancomp \Sigma^\trop$. Notice also that the tame topology for $\tameclass_1 = \{ \lambda\}$ is, by definition, strictly finer than the original topology (since the topological relation between the strata at infinity has changed).

\begin{remark} The restriction to {\em tame refinements of flags} is indeed necessary when introducing $\tameclass$-tame topologies. That is, otherwise, the resulting space will not satisfy the axioms of a topological space (even for the basic families $\tameclass_N$, $N \ge 2$, from Example~\ref{polynomial families}).
\end{remark}

Characterizing the convergent sequences, we have the following analog of Proposition~\ref{prop:convseq}. Notice, in particular, that the $\tameclass$-tame topology on $\cancomp \Sigma^\trop$ formalizes the ideas outlined above.

\begin{prop} \label{prop:convseqtame}
Let $\tameclass$ be a family of functions as above and consider the $\tameclass$-tame topology on $\cancomp \Sigma^\trop$. Fix $\thy \in \cancomp \Sigma^\trop$ and assume $\thy=(x_1,\dots, x_r, t_\fin) \in \inn \keg_\pi^{\trop}$ for the flag $\pi \in \mathscr{F}(\Sigma)$. Suppose that $(\shy_n)_n$ is a sequence in $\cancomp \Sigma^\trop$. Then:

\medskip
(a) If $(\shy_n)_n$ converges to $\thy$ in $ \cancomp \Sigma^\trop$, then almost all $\shy_n$ belong to strata $\inn \keg_{\tilde \pi}^{\trop}$ of flags $\tilde \pi \in \mathscr{F}(\Sigma)$ such that $\pi$ tamely refines $\tilde \pi$.

\medskip
(b) Assume that $\shy_n = (y_{1,n}, \dots, y_{\tilde r, n}, s_{\fin, n})$ defines a sequence in the stratum $ \keg_{\tilde \pi}^\trop$ for some fixed flag $\tilde \pi \in \mathscr{F}(\Sigma)$ such that $\pi$ tamely refines $\tilde \pi$. In particular, $\tilde \pi$ is of the form \eqref{eq:FormCoarserStrata} and, as above, the points $\thy$ and $\shy_n$ can be described via the dual cone of  $\eta := {\tilde \delta_{\infty+\fin}}$. Then, $\shy_n$ converges to $\thy$ in the $\tameclass$-tame topology on $\cancomp \Sigma^\trop$ if and only if $\shy_n$ converges to $\thy$ in the original topology on $\cancomp \Sigma^\trop$  (see conditions (i)---(iii) in Proposition~\ref{prop:convseqtame}) and, additionally,

\begin{itemize}
\item [(iv)] If $\delta_i, \delta_{\hat i}$ are faces in the flag $\pi$ such that $\tilde \delta_\infty \subsetneq \delta_i \subsetneq \delta_{\hat i} \subseteq \delta_\infty$, then
\[
\lim_{n \to \infty} \frac{s_{\fin, n}(a)}{\Psi\left( s_{\fin, n}(\hat a)\right)} = + \infty
\]
for all points $a \in \tau_{i-1} \setminus \tau_i$ and $\hat a\in \tau_{\hat i-1} \setminus \tau_{\hat i}$, and all functions $\Psi \in \tameclass$.
\end{itemize}
 \end{prop}
\begin{proof}
The proposition is clear from the definition of the topology on $ \cancomp \Sigma^\trop$.
\end{proof}

\begin{example}[Example~\ref{polynomial families} continued] The monomial families $\tameclass_N$, $N \ge1$, of growth functions from Example~\ref{polynomial families} appear to be particularly natural choices. Convergence in the $\tameclass_N$-tame topology on $\cancomp \Sigma^\trop$ refines the original topology on $\cancomp \Sigma^\trop$ by requiring that the values on dual faces approach infinity faster than the $N$-th power of the values on less important faces. Note that for $\tameclass_\infty$, we require that the values on dominant faces for every polynomial grow faster than the values on lower infinitary faces.
\end{example}

In applications to moduli spaces (of metric graphs and Riemann surfaces), the topology induced by $\tameclass_2 = \{\lambda, \lambda^2\}$ plays a crucial role. For reasons of abbreviation, we introduce the following terminology.
\begin{defi}[Tropical tame topology] \label{def:TropicalTameTopology}
The topology on $\cancomp \Sigma^\trop$ induced by the family $\tameclass_2 = \{\lambda, \lambda^2\}$ is called in the rest of the paper the {\em tame topology} on $\cancomp \Sigma^\trop$.
\end{defi}
If  $t \in \cancomp \Sigma^\trop$ converges to some $\thy \in \cancomp \Sigma^\trop$ in the tame topology on $\cancomp \Sigma^\trop$, we also say that $t$ {\em converges tamely } to $\thy$.

 \subsection{The case of a simplicial cone}
In order to introduce the moduli spaces of tropical curves, we will employ the tropical compactification of simplicial fans. Namely, we consider the standard fans $\eta = \R_+^E$ in $\R^E$  for the edge set $E$ of a stable marked graph $G$, and then glue the respective compactifications over different graphs of a given genus. (Note that here we use $\eta$ both for the cone $\eta$ and the fan defined by the faces of $\eta$.)

 In this case, we can give an alternate description of the tropical compactificiation in terms of ordered partitions. The topology becomes more explicit as well and can be described in terms of the coordinates without using the dual cones. The aim of this section is to discuss this.

 Notations as above, each cone $\sigma \in \eta$ is uniquely determined by a subset of $E$. The flags in $\mathscr F(\eta)$ can be  then identified with a filtrations
 \[\mathscr F_\bullet^\pi : F_0=\emptyset \subsetneq F_1 \subsetneq F_2 \subsetneq \dots \subsetneq F_r \subseteq F_{\infty + \fin} \subseteq E\]
of subsets $F_{\infty + \fin}$ of $E$. 

The above flags are in turn in bijection with ordered partitions $\pi = (\pi_1, \dots, \pi_r, \pi_\fin)$ of subsets $S$ of $E$, where $r\in \N\cup\{0\}$, $\pi_1, \dots, \pi_r$ are all non-empty, and $\pi_\fin$ can be empty.  As in Section~\ref{sec:preliminaries}, we denote by $\pi_\infty =(\pi_1, \dots, \pi_r)$ the sedentarity of the ordered partition. Slightly abusing the notation, we also refer to the set $E_\infty = \pi_1\cup \dots \cup\pi_r$ as the sedentarity set of the ordered partition. The part $\pi_\fin$ is the finitary part of $\pi$. Recall that the set of ordered partitions of a given set $S$ is denoted by $\Piall(S)$, that consisting of ordered partitions of full sedentarity is denoted by $\Pifs(S)$.

Given an ordered partition $\pi = (\pi_1, \dots, \pi_r, \pi_\fin)\in \Piall(S)$, $S \subseteq E$, the corresponding stratum $\inn\keg_\pi^\trop$ of $\cancomp\eta^\trop$ becomes isomorphic to the product 
\[\inn\keg_\pi^\trop = \inn \sigma_{\pi_1} \times \dots \times \inn\sigma_{\pi_r} \times \inn \keg_{\pi_\fin}\]
where for a subset $A \subset E$, we denote by $\keg_{A}$ and $\inn \keg_A$ the cone $\R_{\geq 0}^A$ and its relative interior $\R_{>0}^A$, respectively. Similarly, for a subset $A \subset E$, we denote by $\sigma_A$ and $\inn \sigma_A$ the standard simplex 
\[\sigma_A = \bigl\{x = (x_i)_{i\in A}\in \R^A\,\bigl|\, x_i \geq 0 \,\, \textrm{ and }\,\, \sum_{i\in A} x_i=1\bigr\}\] 
and its relative interior, respectively. Note in particular that $\inn\sigma_A$ is isomorphic to the projectivized open cone $\rquot{\inn \keg_A}{\R_{>0}}$. 

The compactified cone is then given as 
\[\cancomp \eta^\trop  =\bigsqcup_{\substack{F\subseteq E \\ \pi\in \Piall(F)}} \inn \keg_\pi  = \eta \sqcup \bigsqcup_{\substack{F\subseteq E \\ \pi\in \Piall(F), \pi_\infty \neq \varnothing}} \inn \keg_\pi\]

Assume that $\tilde \pi = (\tilde \pi_1, \dots, \tilde \pi_{\tilde r},\tilde \pi_\fin)$ and that $\tilde \pi\preceq \pi$. This implies that 
\[ E_{\tilde \pi_\infty} \subseteq E_{\pi_\infty} \subseteq E_{\pi} \subseteq E_{\tilde \pi}\]
 and that $\pi$ has the following form
\begin{equation} \label{eq:preceqdecomp}
\pi = \Big( \varrho^1 , \varrho^2 , \dots, \varrho^r \,, \varrho \Big ),
\end{equation}
where each $\varrho^i =(\varrho_k^i)_k$, $i \in  [\tilde r]$, is itself an ordered partition of $\tilde \pi_i$ (without finitary part) and the last part $\varrho =(\varrho_1, \dots, \rho_s, \rho_\fin)$ is an ordered partition of a subset $A \subset \tilde \pi_\fin$. We allow that $A=\varnothing$ and in this case, we have $\varrho = \varnothing$. Note that the finitary part of $\pi$ coincides with $\rho_\fin$.

The convergent sequences in the tropical topology are described as follows.

\begin{prop} \label{prop:convseq}
Let $\thy \in \cancomp\eta^{\trop}$ and assume $\thy = (x_1, \dots, x_r, x_\fin)$ belongs to the stratum $\inn\keg_\pi^{\trop}$ for the ordered partition $\pi = (\pi_1, \dots, \pi_r, \pi_\fin) $ in $\Piall(F)$, $F \subseteq E$. Suppose that $(\shy_n)_n $ is a sequence in $\cancomp\eta^{\trop}$. Then:

\medskip
(a) If $(\shy_n)_n$ converges to $\thy$ in $\cancomp\eta^{\trop}$, then almost all $\shy_n$ belong to strata $\inn \keg_{\tilde \pi}^{\trop}$ of ordered partitions $\tilde \pi$ with $\tilde \pi \preceq \pi$.

\medskip
(b) Assume that $(\shy_n)_n \subseteq \inn\keg_{\tilde \pi}^{\trop}$ for a fixed ordered partition $\tilde \pi = (\tilde \pi_1, \dots, \tilde \pi_{\tilde r},  \tilde \pi_\fin) \preceq \pi$ of $ E_{\tilde \pi} \subseteq E$. In particular, $\pi$ is of the form \eqref{eq:preceqdecomp} and  $\shy_n$ can be written as $\shy_n = (y_{1, n}, \dots, y_{\tilde r,n}, y_{\fin, n})$ with components $y_{i,n} \in \inn \sigma_{\tilde \pi_i}$ and $ y_{\fin, n} \in \R^{\tilde \pi_\fin}$. \smallskip

Then $(\shy_n)_n$ converges to $\thy$ in $\cancomp\eta^{\trop}$ if and only if the following conditions hold:

\begin{itemize}
\item [(i)] The finitary part $y_{\fin, n}$ of $\shy_n$ has the limiting behavior (here $A$ is defined in \eqref{eq:preceqdecomp})
\begin{equation} \label{eq:conv1}
\lim_{n\to\infty}y_{\fin, n} (e) = \begin{cases}
x_{\fin }(e) &e \in  \pi_\fin \\
+ \infty & e \in A\setminus \pi_\fin  \\
0 & e \in  \tilde  \pi_\fin \setminus A \end{cases}.
\end{equation}

\item[(ii)]
Let $\pi_j$, $j \neq \fin$, be a non-finitary set of $\pi$ which is contained in $\tilde  \pi_\fin$ (equivalently, $\pi_j = \varrho_k$ for some $k \neq \fin$, see  \eqref{eq:preceqdecomp}). Then
\begin{equation} \label{eq:conv2}
\lim_{n\to\infty} \frac{ y_{\fin,n} \rest{\pi_j}}{ \sum_{e \in \pi_j} y_{\fin,n}(e) } =x_j \qquad \text{ in } \R^{\pi_j}.
\end{equation}

Moreover, if $\pi_j, \pi_{j+1} \subseteq \tilde \pi_\fin$ are consecutive sets of $\pi$ contained in  $\tilde \pi_\fin$, then
\begin{equation} \label{eq:conv3}
\lim_{n\to\infty}\frac{y_{\fin, n}(e)}{y_{\fin, n}(e')}= + \infty 
\end{equation}
for all $e \in \pi_j$ and $e' \in \pi_{j+1}$.

\item [(iii)] For $i  \in [\tilde r]$, let $\pi_j$ be a set of $\pi$ contained in $ \tilde \pi_i$ (i.e.\ $\pi_j = \varrho^i_k$ for some $k$). Then 
\begin{equation} \label{eq:conv3}
\lim_{n\to\infty} \frac{ y_{i,n} \rest{\pi_j}}{\sum_{e \in \pi_j} y_{i,n}(e)} =x_j \qquad \text{ in } \R^{\pi_j}.
\end{equation}
Moreover, if $\pi_j, \pi_{j+1} \subseteq \tilde \pi_i$ are consecutive sets of $\pi$ contained in  $ \tilde \pi_i$, then
\begin{equation} \label{eq:conv4}
\lim_{n\to\infty}\frac{y_{i, n}(e)}{y_{i, n}(e')}= + \infty
\end{equation}
for all $e \in \pi_j$ and $e' \in \pi_{j+1}$.
\end{itemize}

\medskip
(c) Assume, in addition, that $\pi$ is a tame refinement of $\tilde \pi$. Then $(\shy_n)_n$ converges to $\thy$ in the $\tameclass$-tame topology on $\cancomp\eta^{\trop}$ if and only if the above conditions (i)--(iii) hold and

\begin{itemize}
\item [(iv)] in addition, for all consecutive sets $\pi_j, \pi_{j+1} \subseteq \tilde \pi_\fin$ of $\pi$ contained in  $\tilde \pi_\fin$,
\begin{equation} \label{eq:conv3}
\lim_{n\to\infty}\frac{y_{\fin, n}(e)}{\Psi(y_{\fin, n}(e'))}= + \infty 
\end{equation}
for any $\Psi\in \tameclass$,  $e \in \pi_j$, and  $e' \in \pi_{j+1}$.
\end{itemize}

 \end{prop}
\begin{proof}
The claims are simple consequences of the definition of the topology on $\cancomp\eta^{\trop}$.
\end{proof}

\section{Tropical log maps}\label{sec:tropical_log_map}
In this section, we describe natural log maps in the tropical setting. These are of fundamental importance in the following sections as they provide a language to carry on information between different strata in both the tropical and hybrid moduli spaces.  

\smallskip

 For a simplicial cone $\eta \simeq \R_+^n$,  recall that $\cancomp\eta^\trop$ denotes the tropical compactification and $\partial \cancomp\eta^\trop =  \cancomp\eta^\trop \setminus \eta$ is its boundary at infinity. Similarly, we identify the (closed) standard simplex $\sigma \subset \R^{n}_{+}$ with the projectivized cone $\rquot{\R^{n}_+  \setminus \{0\}}{\R_{>0}}$ and denote by $\cancomp \sigma^\trop$ the tropical compactification of $\inn \sigma$. Note that this is the closure of $\inn \sigma$ in $\cancomp\eta^\trop$ and coincides with the locus of points of full sedentarity in $\cancomp \eta^\trop$. The boundary at infinity of $\cancomp \sigma^\trop$ is denoted by $\partial_\infty \cancomp \sigma^\trop := \cancomp \sigma^\trop \setminus \inn \sigma$ .

\smallskip

The {\em tropical log map} of a simplicial cone $\eta\simeq \R_+^n$ is a retraction map 
\[\logtrop =\logtropind{\eta}\colon \cancomp\eta^\trop\setminus \cube_\eta \to \partial \cancomp\eta^\trop
\]
where $\cube_\eta$ is the hypercube in $\eta$ with all coordinates bounded by one (via the isomorphism of $\eta$ with $\R_+^n$).

\smallskip

More generally, for a fan $\Sigma$, the log maps $\logtropind{\eta}$ of the faces $\eta$ can be combined to a global log map
\[\logtrop = \logtrop_\Sigma \colon \cancomp\Sigma^{\trop}\setminus\bigcup_{\eta \in \Sigma} \cube_\eta \to \partial \cancomp\Sigma^{\trop}. \]

The approach also gives a log map on the \emph{tropical compactification} $\cancomp \sigma^\trop$ of the open {\em standard simplex} $\inn \sigma \subseteq \R_+^n$. Viewing the open standard simplex $\inn\sigma$ as the projectivized cone $\R_{> 0}^n /\R_{>0}$, the latter is a retraction map
\[\logtrop\colon \cancomp \sigma^\trop \setminus\{\one_\sigma\} \to \partial_\infty \cancomp \sigma^\trop,\]
where $\one_\sigma$ has projective coordinates $\one_\sigma  = [1:  \dots :1]$ and $\partial_\infty \cancomp \sigma^\trop = \cancomp \sigma^\trop \setminus \inn \sigma$.

\smallskip

The construction of the above maps is described in the following sections.

\smallskip

This later leads in Section~\ref{sec:tropical_moduli} to the definition of a log map on the {\em moduli space of tropical curves} (of genus $g$ with $n$ marked points)
\[\logtrop\colon \mgtrop{\grind{g,n}} \setminus \umggraph{g,n} \to \partial \mgtrop{\grind{g,n}}, \]
$\umggraph{g,n}$ being the \emph{moduli space of metric graphs} with edge lengths all bounded by one, as well as to the definition of log maps on the \emph{moduli space of hybrid curves}, see Sections~\ref{sec:hybrid_log_I} and~\ref{sec:hybrid_log_II}.

\subsection{Permutohedra and combinatorics of tropical compactifications} \label{ss:Preliminariespermutohedra}
Our constructions rely on the combinatorics of the compactifications $ \cancomp \sigma^\trop$ and $\cancomp\eta^\trop$, which we show to have the same combinatorial structure as two particular polytopes, the so-called {\em permutohedron} and the corresponding {\em generalized permutohedron}. In this section, we briefly recall the definition of these polytopes, and for more general definitions and results refer to e.g.~\cite{Pos09, PRW08}.

\smallskip

Let $y = (y_1, \dots, y_n) \in \R^n$. The \emph{permutohedron} $Q(y)$ in $\R^n$ associated to $y$ is the polytope obtained as the convex hull of all the vectors of the form $y_s = (y_{s(1)}, \dots, y_{s(n)})$ for $s \in \mathfrak S_n$ a permutation of the set $[n]$, that is 
\[Q(y) :=\textrm{conv-hull}\Bigl\{(y_{s(1)}, \dots, y_{s(n)})\,\,\bigl|\,\, s\in \mathfrak S_n\Bigr\}.\]
Alternatively, by a theorem of Rado~\cite{Rad52}, $Q(y)$ can be described as follows. Assume without loss of generality that $y_1\geq y_2 \geq \dots \geq y_n$. Then we have 
\begin{align*}Q(y) = \Bigl\{(z_1, \dots, z_n)\,\,\Bigl|\,\, y_1+\dots+y_k &\geq \sum_{i\in S} z_i \,\textrm{ for any set } S \in {[n]\choose k}  \\
& \textrm{ with equality } y_1+ \dots+y_n = z_1+ \dots+z_n \Bigr\}.
\end{align*}
In the case the coordinates of the vector $y$ are all different,  the face structure of $Q(y)$ can be described as follows. Here, we use the terminology and notation we introduced in the previous sections.

\begin{thm} Notations as above, assume that $y_1> \dots >y_n$. The faces of $Q(y)$ are in bijection with the elements of $\Pifs([n])$, that is, with ordered partitions of $[n]$ of full sedentarity. More precisely, there is a bijection between faces of $Q(y)$ of codimension $r$ and ordered partitions $\pi = (\pi_1, \dots, \pi_r)$ of $[n]$ of rank $r$. Under this bijection, the face corresponding to $\pi$ is given by cutting the permutohedron $Q(y)$ by $r$ affine hyperplanes with the equations
\[\sum_{j\in \pi_k} z_j = \sum_{i = |\pi_1|+\dots + |\pi_{k-1}|+1}^{|\pi_1|+\dots + |\pi_k|}y_{i}, \qquad \textrm{for any } k\in [r].\]
 \end{thm}

Let $\cancomp\sigma^\trop$ be the tropical compactification of the open standard simplex $\inn\sigma=\R_{>0}^{n}/\R_{>0}$.  It follows from the above description that 
\begin{cor} The tropical compactification $\cancomp\sigma^\trop$ has the same combinatorics as the permutohedron $Q(y)$ for any vector $y\in \R^n$ with distinct coordinates.  
\end{cor}

Consider now a vector $y=(y_1, \dots, y_n)$ as above with $y_1>\dots>y_n$. Assume all the $y_j$ are strictly positive, and let $Q(y)$ be the corresponding permutohedron. Define the polytope $P(y)$ by taking the convex hull of all the vectors $z=(z_1, \dots, z_n) \in \R_{\geq 0}^n$ with the property that there exists a permutation $s\in \mathfrak S_n$ such that $z_j \leq y_{s(j)}$ for all $j\in[n]$. Then $P(y)$ can be alternatively described as follows\footnote{We thank Raman Sanyal for pointing to us these existing notions in the literature relevant to our constructions.}:

\smallskip

\begin{itemize}
\item $P(y)$ coincides with the anti-blocking polytope (in the sense of Fulkerson~\cite{Ful71}) associated to $Q(y)$ and obtained by taking the intersection of the Minkowski sum $\bigl(-\R_{\geq 0}^n\bigr) + Q(y)$ with the quadrant $\R_{\geq 0}^n$.  

\smallskip

\item $P(y)$ coincides with the polymatroid (introduced by Edmonds~\cite{Edm70}) associated to the non-decreasing submodular function $f\colon 2^{[n]} \to \R_+$ given by 
\[f(S) = y_1 + \dots + y_{|S|}, \textrm{ for all } S\subseteq [n].\]
\smallskip

\item $P(y)$ is a generalized permutohedron in the sense of Postnikov~\cite{Pos09}.
\smallskip
\item $P(y)$ is the convex-hull of $0$ and all the vectors of the form $\nu_\phi = (z_1, \dots, z_n)$ for an injection $\phi\colon [k]\hookrightarrow [n]$, $k\in [n]$, with 
\[z_j =\begin{cases} y_{i} & \textrm{ if } j=\phi(i), \,\, i\in[k]\\
0 & \textrm{ otherwise.}
\end{cases}\]
\end{itemize}
The face structure of $P(y)$ is given by the following theorem.
\begin{thm} Notations as above, assume that $y_1> \dots >y_n$. The faces of $P(y)$ are in bijection with the elements of $\Piall = \bigsqcup_{S\subset [n]}\Piall(S)$, that is, with ordered partitions $\pi=(\pi_1, \dots, \pi_r, \pi_\fin)$ of subsets $S$ of $[n]$. More precisely, there is a bijection between faces of $P(y)$ of codimension $d$ and ordered partitions $\pi = (\pi_1, \dots, \pi_r, \pi_\fin)$ of subsets $S\subset [n]$ such that $r+n-|S|=d$. Under this bijection, the face corresponding to $\pi$ is given by cutting the polytope $P(y)$ by $r$ affine hyperplanes with the equations
\begin{align*}
z_j &=0  \qquad \textrm{for all } j\in [n]\setminus S, \textrm{ and }\\
\sum_{j\in \pi_k} z_j &= \sum_{i = |\pi_1|+\dots + |\pi_{k-1}|+1}^{|\pi_1|+\dots + |\pi_k|}y_{i}, \qquad \textrm{for any } k\in [r].
\end{align*}
\end{thm}
\begin{proof} This can be obtained from the general theory of generalized permutohedra~\cite{Pos09, PRW08}, or by direct verification. We omit the details. 
\end{proof}

Note in particular that the permutohedron $Q(y)$ appears as the \emph{central} closed face of $P(y)$, that is, the closed face of $P(y)$ associated to $\pi = (\pi_1=[n], \pi_\fin =\varnothing)$. Clearly, $Q(y)$ is the disjoint union of all the open faces of $P(y)$ associated to ordered partitions of $[n]$ of full sedentarity.

\medskip

We will refer to $P(y)$ as the \emph{generalized permutohedron} associated to the vector $y$. 

\smallskip
Let $\eta$ be the standard simplical cone $\R_{+}^n$.  We deduce from the above results the following statement.
\begin{cor} The tropical compactification $\cancomp\eta^\trop$ has the same combinatorics as the generalized permutohedron $P(y)$ for any vector $y\in \R^n_{>0}$ with distinct coordinates.  
\end{cor}

\subsection{Exhausting permutohedra $P_t$ and the foliation of $\cancomp\eta^\trop$} \label{ss:PermutohedraExhaustionFoliation}
As the first step, we introduce a specific family $(P_t)_t$ of generalized permutohedra in $\eta$. This leads to an ascending exhaustion and a foliation of $\cancomp\eta^\trop$.

The constructions which follow depend on the choice of $n$ auxiliary functions 
\[\Phi_1, \Phi_2, \dots, \Phi_n \colon [1, +\infty) \to \R_+\]
 that we suppose 
 \begin{itemize}
 \item smooth, increasing each, with $\Phi_j(1)=1$ and $\lim_{t\to \infty}\Phi_j(t)=\infty$, $j\in [n]$, and verifying
 \smallskip
 
 \item (Log-increasing property) $\Phi'_1(t)/\Phi_1(t)\,< \dots <\,\Phi'_n(t)/\Phi_n(t)$, $t\in [1, +\infty)$, with 
 \[\lim_{t \to \infty} \frac{\Phi_j(t)}{\Phi_{j+1}(t)} =0.\]
 \end{itemize}
Note in particular that the above properties imply $\Phi_1(t)< \dots <\Phi_n(t)$, $t\in (1, +\infty)$. 
 
 \smallskip
 
 We introduce as well a parametrized version of the log-increasing property stated above. Let $C>1$ be a constant. 
  \begin{itemize}
 \item (Log-increasing property with parameter $C$) $\Phi'_{i+1}(t)/\Phi_{i+1}(t) \, >\, C\,\Phi'_i(t)/\Phi_i(t)$, $i\in [n-1]$ and $t\in [1, +\infty)$, with 
 \[\lim_{t \to \infty} \frac{\Phi_j(t)^{^{C}}}{\Phi_{j+1}(t)} =0.\]
 \end{itemize}
In particular, we have $\Phi_i(t)^{^{C}}< \Phi_{i+1}(t)$, $i\in [n-1]$ and $t\in (1, +\infty)$. 
\begin{remark} In application to the moduli spaces in the sequel, we always assume that the family satisfies the log-increasing property with parameter $2$.
\end{remark}

 \subsubsection{Vertices of $\cancomp\eta^\trop$ and the curves $\gamma_\phi$} The {\em vertices at infinity of $\cancomp\eta^\trop$} are in bijection with injective maps $\phi\colon [k] \hookrightarrow [n]$, $k\in [n]$. These correspond themselves to ordered partitions of full sedentarity $\pi = (\pi_\infty, \pi_\fin)$ of subsets $S$ of $[n]$ with $\pi_\fin = \varnothing$. The link is given by the relations
\begin{equation} \label{eq:PartitionFromMap}
S = \phi([k]), \qquad  \qquad \pi = (\pi_\infty, \pi_\fin = \varnothing) \in \Pi(S) \quad \text{with} \quad \pi_\infty = \big ( \{\phi(1)\},   \dots, \{\phi(k)\} \big ),
\end{equation}
that is, $\pi_\infty$ consisting of singletons. 

In terms of this notation, denoting by $\hat \pi = (\pi_\infty, \hat \pi_\fin = [n]\setminus S) \in \Piall([n])$, the vertex $v_\phi \in \cancomp\eta^\trop$ associated with $\phi\colon [k] \hookrightarrow [n]$ is the special point in the stratum $\keg_{\hat \pi} = \sigma_{\pi_\infty} \times \keg_{\hat \pi_\fin}$ with finitary part $x_\fin = 0$ in $ \R_{\geq 0}^{\hat \pi_\fin}$. 
\smallskip

For each vertex $v_\phi$, let $\one_\phi$ be the point in the cone $\eta$ with coordinates $\one_\phi \equiv 1$ on $\phi([k])$ and all other coordinates equal to zero. We join the points $\one_\phi$ and $v_\phi$ by the curve 
\begin{equation} \label{eq:DefinitionVerticesCurves}
\gamma_\phi\colon [1, \infty] \to \cancomp\eta^\trop
\end{equation}
which maps $t \in[1,\infty)$ to the point $\gamma_\phi(t) = (x_1(t), \dots, x_n(t))$ in $\eta$ with coordinates 
\begin{equation} \label{eq:DefinitionVerticesCurves2}
x_j(t) :=\begin{cases} \Phi_{n+1 - l}(t) & \textrm{if } j=\phi(l), \, l\in[k],\\
0& \textrm{otherwise},
\end{cases}
\end{equation}
and sends $\infty$ to $\gamma_\phi(\infty) := v_\phi$. The map $\gamma_\phi\colon [1, \infty] \to \cancomp\eta^\trop$ is clearly a continuous curve. 
\subsubsection{The generalized permutohedra $P_t$} For each $t\in [1,\infty)$, we define the polytope $P_t \subset \eta$ as the convex hull of the point $0$ and the curve points $\gamma_\phi(t)$ (see Figure~\ref{fig:Pt}), 
\[
	P_t := \textrm{conv-hull} \Big( \{0 \} \cup \bigl\{ \gamma_\phi(t)\,\bigl | \,  \phi\colon [k] \hookrightarrow [n], \, k \in [n]\bigr \}\Big).
\]
Note that $P_1 = \cube_\eta$ is a hypercube and the polytope $P_t$ for $t>1$ is the generalized permutohedron $P(y_t)$ associated to the vector $y_t = (\Phi_n(t) , \Phi_{n-1}(t), \dots, \Phi_1(t))$ in $\R^n$.

\smallskip 

We denote by $\partial_\infty P_t$ the union of the proper faces of $P_t$ which do not contain the origin. Each such face corresponds to an ordered partition $\pi = (\pi_\infty, \pi_\fin)$ with non-empty sedentarity part $\pi_\infty \neq \varnothing$, which is an ordered partition of a subset $S \subseteq [n]$.

\smallskip

Setting $P_\infty := \cancomp\eta^\trop$ and $\partial_\infty P_\infty = \partial  \cancomp\eta^\trop$, the following properties are easily verified:\begin{itemize} \setlength\itemsep{1 mm}
\item The polytopes $P_t$, $t \in (1,  \infty)$ form an {\em increasing exhaustion} of the cone $\eta$. \newline
(That is, $P_s \subset P_t$ for $s<t$ and $\bigcup_{t \in (1, \infty)} P_t = \eta$).
\item The boundaries $\partial_\infty P_t$, $t\in (1, \infty]$, provide a {\em foliation} of $\cancomp\eta^\trop \setminus \cube_\eta$. \newline
(That is, the sets are pairwise disjoint and $\bigcup_{t \in (1, \infty]} \partial_\infty P_t =  \cancomp\eta^\trop \setminus \cube_\eta$).
\end{itemize}

\begin{center}

\begin{minipage}{0.5 \textwidth}
\centering 
\begin{tikzpicture}[scale = 3]


\filldraw (0,0)  circle (0.4 pt);

\coordinate (A) at ({cos(90)}, {sin(90)});

\coordinate (B) at ({cos(210)}, {sin(210)});

\coordinate (C) at ({cos(330)}, {sin(330)});


\draw[thick] (A)--(B);
\draw[thick] (B)--(C);
\draw[thick] (C)--(A);

\draw[scale=1, domain=0:1, smooth, variable=\s, blue] plot
({(\s*\s*cos(90) + \s*cos(210) + cos(330))/(1+\s+\s*\s) },
{(\s*\s*sin(90) + \s*sin(210) + sin(330))/(1+\s+\s*\s) });

\draw[scale=1, domain=0:1, smooth, variable=\s, blue] plot
({(\s*cos(90) + \s*\s* cos(210) + cos(330))/(1+\s+\s*\s) },
{(\s*sin(90) + \s*\s*sin(210) + sin(330))/(1+\s+\s*\s) });

\draw[scale=1,  domain=0:1, smooth, variable=\s, blue] plot
({(\s*\s*cos(90) + cos(210) + \s*cos(330))/(1+\s+\s*\s) },
{(\s*\s*sin(90) + sin(210) + \s*sin(330))/(1+\s+\s*\s) });

\draw[scale=1, domain=0:1, smooth, variable=\s, blue] plot
({(\s*cos(90) + cos(210) + \s*\s*cos(330))/(1+\s+\s*\s) },
{(\s*sin(90) + sin(210) + \s*\s*sin(330))/(1+\s+\s*\s) });

\draw[scale=1, domain=0:1, smooth, variable=\s, blue] plot
({(cos(90) + \s*\s*cos(210) + \s*cos(330))/(1+\s+\s*\s) },
{(sin(90) + \s*\s*sin(210) + \s*sin(330))/(1+\s+\s*\s) });

\draw[scale=1,  domain=0:1, smooth, variable=\s, blue] plot
({(cos(90) + \s*cos(210) + \s*\s*cos(330))/(1+\s+\s*\s) },
{(sin(90) + \s*sin(210) + \s*\s*sin(330))/(1+\s+\s*\s) });

\def \R {1/4}

\coordinate (Q1) at  ({(cos(90) + \R*\R*cos(210) + \R*cos(330))/(1+\R+\R*\R) },
{(sin(90) + \R*\R*sin(210) + \R*sin(330))/(1+\R+\R*\R) });

\coordinate (Q2) at
({(cos(90) + \R*cos(210) + \R*\R*cos(330))/(1+\R+\R*\R) },
{(sin(90) + \R*sin(210) + \R*\R*sin(330))/(1+\R+\R*\R) });


\coordinate (Q3) at ({(\R*cos(90) + cos(210) + \R*\R*cos(330))/(1+\R+\R*\R) },
{(\R*sin(90) + sin(210) + \R*\R*sin(330))/(1+\R+\R*\R) });
\coordinate (Q4) at ({(\R*\R*cos(90) + cos(210) + \R*cos(330))/(1+\R+\R*\R) },
{(\R*\R*sin(90) + sin(210) + \R*sin(330))/(1+\R+\R*\R) });

\coordinate (Q5) at ({(\R*\R*cos(90) + \R*cos(210) + cos(330))/(1+\R+\R*\R) },
{(\R*\R*sin(90) + \R*sin(210) + sin(330))/(1+\R+\R*\R) });

\coordinate (Q6) at ({(\R*cos(90) + \R*\R*cos(210) + cos(330))/(1+\R+\R*\R) },
{(\R*sin(90) + \R*\R*sin(210) + sin(330))/(1+\R+\R*\R) });

\draw[fill opacity=0.3, fill = teal] (Q1) -- (Q2) -- (Q3) -- (Q4) -- (Q5) -- (Q6) -- (Q1);

\def \R {0.5}

\coordinate (P1) at  ({(cos(90) + \R*\R*cos(210) + \R*cos(330))/(1+\R+\R*\R) },
{(sin(90) + \R*\R*sin(210) + \R*sin(330))/(1+\R+\R*\R) });

\coordinate (P2) at
({(cos(90) + \R*cos(210) + \R*\R*cos(330))/(1+\R+\R*\R) },
{(sin(90) + \R*sin(210) + \R*\R*sin(330))/(1+\R+\R*\R) });


\coordinate (P3) at ({(\R*cos(90) + cos(210) + \R*\R*cos(330))/(1+\R+\R*\R) },
{(\R*sin(90) + sin(210) + \R*\R*sin(330))/(1+\R+\R*\R) });
\coordinate (P4) at ({(\R*\R*cos(90) + cos(210) + \R*cos(330))/(1+\R+\R*\R) },
{(\R*\R*sin(90) + sin(210) + \R*sin(330))/(1+\R+\R*\R) });

\coordinate (P5) at ({(\R*\R*cos(90) + \R*cos(210) + cos(330))/(1+\R+\R*\R) },
{(\R*\R*sin(90) + \R*sin(210) + sin(330))/(1+\R+\R*\R) });

\coordinate (P6) at ({(\R*cos(90) + \R*\R*cos(210) + cos(330))/(1+\R+\R*\R) },
{(\R*sin(90) + \R*\R*sin(210) + sin(330))/(1+\R+\R*\R) });

\draw[dashed] (P1) -- (P2) -- (P3) -- (P4) -- (P5) -- (P6) -- (P1);

\def \R {0.4}

\coordinate (S1) at  ({(cos(90) + \R*\R*cos(210) + \R*cos(330))/(1+\R+\R*\R) },
{(sin(90) + \R*\R*sin(210) + \R*sin(330))/(1+\R+\R*\R) });

\coordinate (S2) at
({(cos(90) + \R*cos(210) + \R*\R*cos(330))/(1+\R+\R*\R) },
{(sin(90) + \R*sin(210) + \R*\R*sin(330))/(1+\R+\R*\R) });


\coordinate (S3) at ({(\R*cos(90) + cos(210) + \R*\R*cos(330))/(1+\R+\R*\R) },
{(\R*sin(90) + sin(210) + \R*\R*sin(330))/(1+\R+\R*\R) });
\coordinate (S4) at ({(\R*\R*cos(90) + cos(210) + \R*cos(330))/(1+\R+\R*\R) },
{(\R*\R*sin(90) + sin(210) + \R*sin(330))/(1+\R+\R*\R) });

\coordinate (S5) at ({(\R*\R*cos(90) + \R*cos(210) + cos(330))/(1+\R+\R*\R) },
{(\R*\R*sin(90) + \R*sin(210) + sin(330))/(1+\R+\R*\R) });

\coordinate (S6) at ({(\R*cos(90) + \R*\R*cos(210) + cos(330))/(1+\R+\R*\R) },
{(\R*sin(90) + \R*\R*sin(210) + sin(330))/(1+\R+\R*\R) });

\draw[dashed] (S1) -- (S2) -- (S3) -- (S4) -- (S5) -- (S6) -- (S1);

\def \R {0.3}

\coordinate (S1) at  ({(cos(90) + \R*\R*cos(210) + \R*cos(330))/(1+\R+\R*\R) },
{(sin(90) + \R*\R*sin(210) + \R*sin(330))/(1+\R+\R*\R) });

\coordinate (S2) at
({(cos(90) + \R*cos(210) + \R*\R*cos(330))/(1+\R+\R*\R) },
{(sin(90) + \R*sin(210) + \R*\R*sin(330))/(1+\R+\R*\R) });


\coordinate (S3) at ({(\R*cos(90) + cos(210) + \R*\R*cos(330))/(1+\R+\R*\R) },
{(\R*sin(90) + sin(210) + \R*\R*sin(330))/(1+\R+\R*\R) });
\coordinate (S4) at ({(\R*\R*cos(90) + cos(210) + \R*cos(330))/(1+\R+\R*\R) },
{(\R*\R*sin(90) + sin(210) + \R*sin(330))/(1+\R+\R*\R) });

\coordinate (S5) at ({(\R*\R*cos(90) + \R*cos(210) + cos(330))/(1+\R+\R*\R) },
{(\R*\R*sin(90) + \R*sin(210) + sin(330))/(1+\R+\R*\R) });

\coordinate (S6) at ({(\R*cos(90) + \R*\R*cos(210) + cos(330))/(1+\R+\R*\R) },
{(\R*sin(90) + \R*\R*sin(210) + sin(330))/(1+\R+\R*\R) });

\draw[dashed] (S1) -- (S2) -- (S3) -- (S4) -- (S5) -- (S6) -- (S1);

\node  at (-0.6,0.5) {\color{black}{$Q_t$}};
\node  at (0.4,0.8) {\color{blue}{$\gamma_\phi$}};

\end{tikzpicture}
  \captionof{figure}{The permutohedra $Q_t$ in dimension $n=2$.}
  \label{fig:Qt}
\end{minipage}\begin{minipage}{0.5 \textwidth}
\centering

\begin{tikzpicture}


\draw [->] (0,0) -- (5,0);
\draw [->] (0,0) -- (0,5);



\draw[scale=1, domain=1:2.2, smooth, variable=\s, blue] plot
({(\s) },{(\s*\s)});

\draw[scale=1, domain=1:2.2, smooth, variable=\s, blue] plot
({(\s*\s)}, {(\s) });


\def \R {2}

\coordinate (Q1) at  ({\R^2}, {(0)}) ;
\coordinate (Q2) at  ({\R^2}, {(\R)}) ;
\coordinate (Q3) at  ({\R}, {(\R^2)}) ;
\coordinate (Q4) at  ({0}, {(\R^2)});

\draw[fill opacity=0.3, fill = teal] (Q1) -- (Q2) -- (Q3) -- (Q4) -- (0,0) -- (Q1);

\def \R {1.8}

\coordinate (Q1) at  ({\R^2}, {(0)}) ;
\coordinate (Q2) at  ({\R^2}, {(\R)}) ;
\coordinate (Q3) at  ({\R}, {(\R^2)}) ;
\coordinate (Q4) at  ({0}, {(\R^2)}); 

\draw[dashed] (Q1) -- (Q2) -- (Q3) -- (Q4) -- (0,0) -- (Q1);

\def \R {1.6}

\coordinate (Q1) at  ({\R^2}, {(0)}) ;
\coordinate (Q2) at  ({\R^2}, {(\R)}) ;
\coordinate (Q3) at  ({\R}, {(\R^2)}) ;
\coordinate (Q4) at  ({0}, {(\R^2)}); 

\draw[dashed] (Q1) -- (Q2) -- (Q3) -- (Q4) -- (0,0) -- (Q1);

\def \R {1.4}

\coordinate (Q1) at  ({\R^2}, {(0)}) ;
\coordinate (Q2) at  ({\R^2}, {(\R)}) ;
\coordinate (Q3) at  ({\R}, {(\R^2)}) ;
\coordinate (Q4) at  ({0}, {(\R^2)}); 

\draw[dashed] (Q1) -- (Q2) -- (Q3) -- (Q4) -- (0,0) -- (Q1);

\draw (0,0) -- (1,0) -- (1,1) -- (0,1) -- (0,0);
\node  at (-1, 0.5) {$P_1 =\cube$ };

\node  at (2.5, 4.5) {\color{blue}{$\gamma_\phi$}};
\node  at (-1, 3.8)  {\color{black}{$P_t$}};

\end{tikzpicture}
  \captionof{figure}{The generalized permutohedra $P_t$ in dimension $n=2$.}
  \label{fig:Pt}
\end{minipage}
\end{center}

\begin{remark} \label{rem:OtherVertexCurves} Imagine we are given $\Phi_1, \dots, \Phi_m$, $m\geq n$, auxiliary functions with the same properties listed in the beginning of this section. The possibility of having more auxiliary functions adds flexibility.
For example, we can choose the vertex curves $\gamma_\phi(t)$ in \eqref{eq:DefinitionVerticesCurves} differently by considering instead the curves $\gamma_\phi^m(t)$ with coordinates $(x^m_1(t) ,\dots, x^m_n(t))$ given by 
\begin{equation} \label{eq:AlternativeVerteCurves}
x^m_j(t) :=\begin{cases} \Phi_{m+1 - l}(t) & \textrm{if } j=\phi(l), \, l\in[k],\\
0& \textrm{otherwise},
\end{cases}
\end{equation}
and with $\gamma^m_\phi(\infty):=v_\phi$. In this notation, we have $\gamma_\phi = \gamma^n_\phi$, where $n = \dim \eta$.

The associated polytopes $P^m_t$, $t>1$, in $\eta$ are again generalized permutohedra and all relevant properties also hold for this sequence. Intersecting the polytopes $P^m_t$ with the coordinate hyperplanes of a subset $S \subseteq [n]$, we get
\[
	 \bigl\{x \in P^m_{t}\, \bigl| \, x(j) = 0 \text{ for } j \notin S \bigr\} = P^m_{S, t} \times 0_{[n]\setminus S}
	 \]
for the corresponding polytope $P^m_{S, t}$ in $\R_+^S$ given by the same auxiliary functions. Hence this extension will appear naturally when we glue the log maps on fans and moduli spaces, even if we start with $n$ auxiliary functions. 
\end{remark}

\subsection{Projectivized version} Analogous to the case of cones, we obtain an increasing exhaustion and a foliation of $\cancomp\sigma^\trop$ using a specific family $(Q_t)_t$ of permutohedra in $\inn \sigma$.

\smallskip
  
Namely, the vertices of the compactified simplex $\cancomp \sigma^\trop$ are in bijection with the permutations $\phi\colon [n] \to [n]$. For any such  permutation $\phi$, consider the continuous curve 
\begin{align*}
\gamma_\phi\colon [1,
\infty] \to \cancomp \sigma^\trop \\
t \mapsto [\Phi_{\phi(n)}(t) : \dots : \Phi_{\phi(1)}(t)] 
\end{align*}
connecting the point $\one_\sigma$ to its vertex $v_\phi$ of $\cancomp \sigma^\trop$.  As for the case of cones treated before, we define the polytopes $Q_t$, $t\in(1, \infty)$, as the convex hull of the curve points $\gamma_\phi(t)$ (see Figure~\ref{fig:Qt}),
\[
	Q_t :=\textrm{conv-hull}\Big( \bigl\{ \gamma_\phi(t) \,\st   \phi\colon [n] \rightarrow [n] \text{ permutation} \bigr\}\Big).
\] 
The convex hull here is taken inside the simplex $\inn\sigma$. Alternatively, we take the cone inside $\eta$ whose rays are given by $\R_{\geq 0} \gamma_\phi$, $\phi$ a permutation of $[n]$, $Q_t$ corresponds then to the projectivization of this cone. 
 
By construction, $Q_t$ is a permutohedron for each $t > 1$. We denote its boundary by $\partial_\infty Q_t$, which is the union of the proper faces of $Q_t$. Each such face corresponds to an ordered partition $\pi = (\pi_1, \dots, \pi_r, \pi_\fin)$ of $[n]$ with $\pi_\fin = \varnothing$ and $r \ge 2$.

\smallskip

Note that the permutohedron $Q_t$ is isomorphic to the central face of the generalized permutohedron $P_t$ under the projection map from $\eta = \R_+^n$ to $\sigma = \rquot{\R^{n}_+ \setminus \{0\}}{\R_{>0}} $.

\smallskip

Denoting again $Q_\infty = \cancomp\sigma^\trop$ and $\partial_\infty Q_\infty = \cancomp\sigma^\trop \setminus \inn \sigma$, we obtain the following result.
\begin{prop} The permutohedra $Q_t$, $t \in (1, \infty)$, form an {\em increasing exhaustion} of $\inn \sigma$ and the boundaries $\partial_\infty Q_t$, $t\in (1, \infty]$, provide a {\em foliation} of $\cancomp \sigma^\trop \setminus\{\one_\sigma\}$. 
\end{prop}
\begin{proof} Consider the standard simplex $\sigma \subseteq \R_{\geq 0}^n$ consisting of all the points whose coordinates sum to one. We view each polytope $Q_t$ inside $\sigma$, as the convex hull of $\frac 1{L(t)} \bigl(\Phi_{\phi(1)}(t), \dots,\Phi_{\phi(n)}(t)\bigr)$, with $L(t) := \sum_{j=1}^n\Phi_{j}(t)$ and $\phi$ running over the permutations of $[n]$.

 We need to prove that $Q_s \subset \inn Q_t$ provided that $s<t$.  By the characterization of the polytopes in terms of inequalities, it will be enough to prove that for each subset $S\subseteq [n]$ of size $|S|=k$, $k\in [n]$, we have 
 \[\frac 1{L(s)}\sum_{j\in S} \Phi_j(s) < \frac 1{L(t)} \bigl(\Phi_n(t) + \dots + \Phi_{n-k+1}(t)\bigl).\]
 This itself reduces to the case $S = \{n-k+1, \dots, n\}$, that is, to 
 \[\frac {L(t)}{L(s)} < \frac{\Phi_n(t) + \dots + \Phi_{n-k+1}(t)}{\Phi_n(s) + \dots + \Phi_{n-k+1}(s)},\]
 which we get from the inequalities 
 \[\Phi_1(t)/\Phi_1(s) < \dots <\Phi_n(t)/\Phi_n(s),\]
 using the log-increasing property of the sequence of auxiliary functions. 
\end{proof}

The constructions extend to the case where we have $m$ auxiliary functions $\Phi_1, \dots, \Phi_m$ and give rise to the polytopes $Q_t^m$, with the same properties listed above.

\subsection{Faces of the polytopes $P_t$ and $Q_t$} Our later constructions rely on the recursive structure of the exhausting polytopes $P_t$ and $Q_t$. As we show next, each boundary face of these polytopes is a product of the same type of polytopes in lower dimension.
\smallskip

For the ease of presentation, we place ourselves in the more general setting where we are given $m$ auxiliary functions $\Phi_1, \dots, \Phi_m$, so $P_t= P^m_t$ and $Q_t = Q^m_t$.

 Let $\face_{\pi, t}$ be a face in $\partial_\infty P_t$ corresponding to an ordered partition $\pi = (\pi_\infty=(\pi_1, \dots, \pi_r), \pi_\fin)$ of $S \subseteq [n]$ with $\pi_\infty \neq \varnothing$. The associated face $\face^\trop_\pi$ in $\cancomp \eta^\trop$ is the product
 \[\face^\trop_\pi =\face_{\pi, \infty} \cong \cancomp\sigma^\trop_{\pi_1}\times  \dots \cancomp\sigma^\trop_{\pi_r}\times \cancomp \keg_{\pi_\fin}^\trop.
 \] 

Similarly, let $\sface_{\pi, t}$ be a face in $\partial_\infty Q^m_t$ associated to an ordered partition $\pi = (\pi_1, \dots, \pi_r, \pi_\fin)$ of $[n]$ with $\pi_\fin = \varnothing$ and $r \ge 2$. Then the corresponding face $\sface^\trop_\pi$ of $\cancomp \sigma^\trop$ is the product
\[\sface^\trop_\pi =\sface_{\pi, \infty} \cong \cancomp\sigma^\trop_{\pi_1}\times \dots \cancomp\sigma^\trop_{\pi_r}.
\]

For a subset $E \subseteq [n]$, we denote by $\cancomp\sigma_E^\trop$ and  $\cancomp \eta_{E}^\trop$ the compactifications of the simplex $\inn \sigma_E \subseteq \R_+^E$ and the cone $\eta_E := \R_+^E$. For each $m\geq k\geq \dim\eta$, we write $Q^k_{E, t}$ and $P^k_{E, t}$, $t\in(1, \infty]$, for the exhausting polytopes in $\cancomp\sigma_E^\trop$ and $\cancomp \eta_{E}^\trop$, respectively, obtained from the curves $\gamma_\phi^k(t)$ in \eqref{eq:AlternativeVerteCurves}, relative to the auxiliary functions $\Phi_1, \dots, \Phi_k$. The polytopes $Q_{E,t}^k$, $t\in(0,+\infty)$, are all viewed inside the standard simplex of the corresponding cones, that is the sum of coordinates of points in $Q^k_{E,t}$ is equal to one.

\begin{lem}\label{lem:RecursiveStructureFaces}
(i) Let $\pi=(\pi_\infty=(\pi_1, \dots, \pi_r), \pi_\fin)$ be an ordered partition of $S \subset [n]$. The face $\face_{\pi,t}$ in the boundary $\partial_\infty P_t$ admits the product decomposition
\[
\face_{\pi, t} = A_{1,t} \times \dots A_{r, t} \times A_{\fin,t} \times 0^{n - |S|}, \qquad t \in (1, \infty),
\]
where $A_{\fin,t} = P_{\pi_f, t}^{m_\fin}$ is the generalized permutohedron in $\R_+^{\pi_\fin}$ for $m_\fin = m - |\pi_\infty|$, and
\begin{align*}
A_{j, t} = \chi_j(t) \cdot Q^{m_j}_{\pi_j,t} & \qquad \text{with}\qquad  m_j =m-(|\pi_1| + \dots + |\pi_{j-1}|), m_1=m,  \\
  \qquad \chi_j(t) := \sum_{i=m_{j+1}+1}^{m_j}\Phi_{i}(t).
\end{align*}
In particular, the projection to $\face_{\pi}^\trop$ gives an isomorphism $\face_{\pi_,t} \cong Q^{m_1}_{\pi_1, t}\times \dots\times Q^{m_r}_{\pi_r,t} \times P^{m_\fin}_{\pi_\fin, t}$.

\smallskip

(ii) Each face $\sface_{\pi,t}$, $\pi=(\pi_\infty=(\pi_1, \dots, \pi_r), \pi_\fin=\varnothing)$ of full sedentarity, in the boundary $\partial_\infty Q_t$, $t \in (1, \infty)$ admits the product decomposition
\[
\sface_{\pi,t} = A_{1,t} \times \dots \times A_{r,t}, 
\]
where 
\begin{align*}
A_{j,t} = \chi_j(t) L(t)^{-1} \cdot Q^{m_j}_{\pi_j, t}  & \qquad \textrm{with} \qquad m_j = m-(|\pi_1| + \dots + |\pi_{j-1}|),\\
\chi_j(t) := \sum_{i=m_{j+1}+1}^{m_j}\Phi_{i}(t) & \qquad \textrm{and} \qquad L(t) = \sum_{i=m-n+1}^{m}\Phi_{i}(t).
\end{align*}
In particular, the natural projection to $\sface_{\pi}^{\trop}$ gives an isomorphism $ \sface_{\pi,t} \cong Q^{m_1}_{\pi_1, t}\times \dots\times Q^{m_r}_{\pi_r,t}$.
\end{lem}
\begin{proof} The first claim follows from the shape of the curves $\gamma^m_\phi$ in \eqref{eq:DefinitionVerticesCurves}. Also, (i) implies (ii), since $\sface_{\pi, t}$ is the image of the corresponding face $\face_{\pi,t}$ of $P_t$ under projectification.
\end{proof}
For a fixed face $\face_{\pi,t}$ or $\sface_{\pi, t}$, we write $\pr_{j,t}$, $j \in [r]$ for the natural isomorphisms $\pr_{j,t}\colon A_{j,t} \to Q^{m_j}_{\pi_j, t}$ (which are given by projectifying).

\subsection{The definition of the log map}
To define the log maps, we first introduce maps from the polytopes $P_t$ and $Q_t$ to $\cancomp\eta^\trop$ and $\cancomp\sigma^\trop$. More precisely, we will construct for each $t\in (1, \infty]$ two homeomorphisms
\[h_t \colon P_t  \to \cancomp\eta^\trop  \qquad \text{and} \qquad \hbar_t \colon Q_t  \to \cancomp\sigma^\trop\]
which preserve the structure of faces. The construction gives a globally continuous family of maps. Namely, define the subsets $\mathcal P \subset \cancomp\eta^\trop \times (1,\infty]$ and $\mathcal Q \subset \cancomp \sigma^\trop \times (1,\infty]$ by
\begin{align} \label{eq:GlobalSets}
\mathcal P = \bigcup_{t \in (1, \infty]} P_t  \times\{t\} \qquad \text{and} \qquad \mathcal Q = \bigcup_{t \in (1, \infty]} Q_t \times \{t\} ,
\end{align}
equipped with their induced topologies as subsets. Then, as we will show, the global functions 
\begin{align} \label{eq:GlobalFunctions}
h\colon \mathcal P \to \cancomp\eta^\trop \qquad \text{and} \qquad \hbar \colon \mathcal Q \to \cancomp\sigma^\trop
\end{align}
which restrict to $h_t$ on $P_t \times\{t\}$ and $\hbar_t$ on $Q_t \times\{t\}$, respectively, are continuous.

\smallskip

Using the foliation property, the log map $\logtrop \colon \cancomp\eta^\trop\setminus \cube_\eta \to \partial \cancomp\eta^\trop$ will be defined by
\[\logtrop\rest{\partial_\infty P_t} := h_t\rest{\partial_\infty P_t}, \qquad t \in (1, \infty]. \]
Similarly, the log map $\logtrop\colon \cancomp \sigma^\trop \setminus\{\one_\sigma\} \to \partial_\infty \cancomp \sigma^\trop$ will be given by
\[\logtrop\rest{\partial_\infty Q_t} := \hbar_t\rest{\partial_\infty Q_t}, \qquad t \in (1, \infty]. \]

\subsection{The homeomorphisms $h_t$ and $\hbar_t$: Outline} Let us briefly sketch the construction of the homeomorphisms $h_t\colon P_t  \to \cancomp \eta^\trop$ and $\hbar_t\colon Q_t \to \cancomp \sigma^\trop$, $t > 1$.

\smallskip

The maps $h_t$ and $\hbar_t $ will be defined {\em recursively}, using induction on the dimension of the respective space. Observe that the polytopes $P_t$ and $Q_t$ are naturally decomposed into regions, induced by the faces of their respective boundary. Namely, we can expand each boundary face $\face_{\pi,t}$ and $\sface_{\pi,t}$ (coming from an ordered partition $\pi$) into the following region
\[
R_{\pi,t} := \bigcup_{s \in (1, t]} \face_{\pi,s}
\]
and obtain a cover of $P_t \setminus \cube_\eta$ and $Q_t \setminus \{\one_\sigma\}$, respectively (see Figure \ref{fig:RegionsPt} and Figure \ref{fig:RegionsQt}). Similarly, the spaces $\cancomp \eta^\trop \setminus \cube_\eta$ and $\cancomp \sigma^\trop  \setminus \{\one_\sigma\}$, respectively, are covered by the regions $R_{\pi, \infty}$ (setting $t = \infty$ in the above definition). We define a map from each region $R_{\pi,t}$ to $R_{\pi,\infty}$ and glue them to a global homeomorphism (on $\cube_\eta$ and $\{\one_\sigma\}$, we define the map as the identity map).

\smallskip

By Lemma~\ref{lem:RecursiveStructureFaces}, the boundary face $\face_{\pi,t}$ is naturally obtained from lower dimensional polytopes of the same type. Hence, using the corresponding maps for lower dimensions, we obtain a natural homeomorphism $g_t$ from $\face_{\pi,t}$ to the face at infinity $\face_{\pi, \infty}$. We then "stretch" $g_t$ to a map $h_t$ between the regions $R_{\pi,t}$ to $R_{\pi,\infty}$. The idea is to think of $R_{\pi,t}$ and $R_{\pi, \infty}$ as products $(1, t]  \times C$ and $(1, \infty] \times C$ for some fixed combinatorial object $C$. Having an auxiliary homeomorphism $\iota_t$ from $ (1, t]$ to $(1, \infty]$, each boundary face $\face_{\pi,s}$ is mapped to $\face_{\pi,\tilde s}$ with $\tilde \iota_t(s)$, in a way dictated by the already existing map $g_t$ on $\face_{\pi, t}$.

\smallskip

We call this procedure the {\em unfolding} of the map $g_t \colon \face_{\pi,t} \to \face_{\pi, \infty}$ to the region $R_{\pi, t}$ (since, similar to an accordion, the region $R_{\pi,t}$ becomes "unfolded to infinity").

\begin{center}
\begin{minipage}{0.45 \textwidth}
\centering

\begin{tikzpicture}


\draw (0,0) -- (5,0);
\draw(0,0) -- (0,5);

\draw[scale=1, domain=1:2.2, smooth, variable=\s, black] plot
({(\s) },{(\s*\s)});


\def \R {2}

\coordinate (Q1) at  ({\R^2}, {(0)}) ;
\coordinate (Q2) at  ({\R^2}, {(\R)}) ;
\coordinate (Q3) at  ({\R}, {(\R^2)}) ;
\coordinate (Q4) at  ({0}, {(\R^2)});

\draw[dashed] (Q1) -- (Q2) -- (Q3) -- (Q4) -- (0,0) -- (Q1);

\def \R {1.8}

\coordinate (Q1) at  ({\R^2}, {(0)}) ;
\coordinate (Q2) at  ({\R^2}, {(\R)}) ;
\coordinate (Q3) at  ({\R}, {(\R^2)}) ;
\coordinate (Q4) at  ({0}, {(\R^2)}); 

\draw[dashed] (Q1) -- (Q2) -- (Q3) -- (Q4) -- (0,0) -- (Q1);

\def \R {1.6}

\coordinate (Q1) at  ({\R^2}, {(0)}) ;
\coordinate (Q2) at  ({\R^2}, {(\R)}) ;
\coordinate (Q3) at  ({\R}, {(\R^2)}) ;
\coordinate (Q4) at  ({0}, {(\R^2)}); 

\draw (Q1) -- (Q2) -- (Q3) -- (Q4) -- (0,0) -- (Q1);

\def \R {1.4}

\coordinate (Q1) at  ({\R^2}, {(0)}) ;
\coordinate (Q2) at  ({\R^2}, {(\R)}) ;
\coordinate (Q3) at  ({\R}, {(\R^2)}) ;
\coordinate (Q4) at  ({0}, {(\R^2)}); 

\draw[dashed] (Q1) -- (Q2) -- (Q3) -- (Q4) -- (0,0) -- (Q1);

\def \R {1.6}
\draw[scale=1, domain=1:\R, smooth, variable=\s, black, name path = A] plot
({(\s) },{(\s*\s)});

\draw[scale=1, domain=1:\R, smooth, variable=\s, black, name path = B] plot
({(0) },{(\s*\s)});
\tikzfillbetween[of=A and B]{red, opacity = 0.3};

\draw[very thick, red] (0, \R * \R) -- (\R, \R*\R);

\node  at (-1, \R*\R) {\color{red}{$\face_{\pi, t}$}};
\node  at (-1, 1.4) {\color{red}{$R_{\pi,t}$}};

\def \R {2.3}

\draw[scale=1, domain=1:\R, smooth, variable=\s, black, name path = A] plot
({(\s) },{(\s*\s)});

\draw[scale=1, domain=1:\R, smooth, variable=\s, black, name path = B] plot
({(0) },{(\s*\s)});
\tikzfillbetween[of=A and B]{ pattern=north east lines, pattern color=gray};

\node  at (-1, 5) {\color{black}{$R_{\pi,\infty}$}};

\def \R {2.5}
\draw[thick, black] (0, \R * \R) -- (2.8, \R*\R);
\draw[dotted, black] (2.8, \R * \R) -- (3.8, \R*\R);
\draw[dashed] (0,2.3) -- (0, \R);

\node  at (-1, \R*\R) {\color{black}{$\face_{\pi, \infty} $}};



\end{tikzpicture}
  
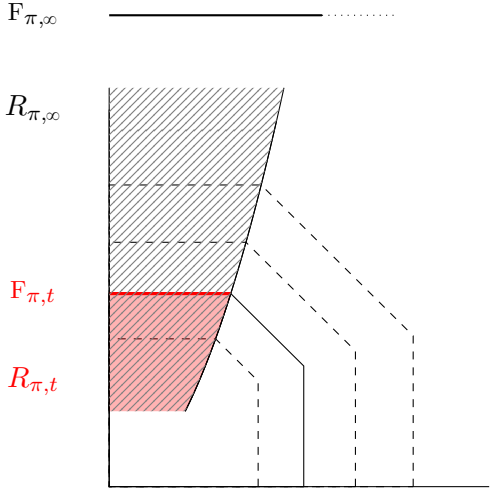
\captionof{figure}{A region $R_{\pi,t}$ of the polytope $P_t$, generated by a boundary face $\face_{\pi,t}$, and the infinite region $R_{\pi, \infty}$.}
  \label{fig:RegionsPt}
\end{minipage} \begin{minipage}{0.45 \textwidth} 

\centering 
\begin{tikzpicture}[scale = 3.3]


\filldraw (0,0)  circle (0.4 pt);

\coordinate (A) at ({cos(90)}, {sin(90)});

\coordinate (B) at ({cos(210)}, {sin(210)});

\coordinate (C) at ({cos(330)}, {sin(330)});


\draw[thick] (A)--(B);
\draw[thick] (B)--(C);
\draw[thick] (C)--(A);

\draw[scale=1, domain=0:1, smooth, variable=\s,] plot
({(\s*\s*cos(90) + \s*cos(210) + cos(330))/(1+\s+\s*\s) },
{(\s*\s*sin(90) + \s*sin(210) + sin(330))/(1+\s+\s*\s) });

\draw[scale=1, domain=0:1, smooth, variable=\s, ] plot
({(\s*cos(90) + \s*\s* cos(210) + cos(330))/(1+\s+\s*\s) },
{(\s*sin(90) + \s*\s*sin(210) + sin(330))/(1+\s+\s*\s) });

\draw[scale=1,  domain=0:1, smooth, variable=\s,] plot
({(\s*\s*cos(90) + cos(210) + \s*cos(330))/(1+\s+\s*\s) },
{(\s*\s*sin(90) + sin(210) + \s*sin(330))/(1+\s+\s*\s) });

\draw[scale=1, domain=0:1, smooth, variable=\s, ] plot
({(\s*cos(90) + cos(210) + \s*\s*cos(330))/(1+\s+\s*\s) },
{(\s*sin(90) + sin(210) + \s*\s*sin(330))/(1+\s+\s*\s) });

\draw[scale=1, domain=0:1, smooth, variable=\s, ] plot
({(cos(90) + \s*\s*cos(210) + \s*cos(330))/(1+\s+\s*\s) },
{(sin(90) + \s*\s*sin(210) + \s*sin(330))/(1+\s+\s*\s) });

\draw[scale=1,  domain=0:1, smooth, variable=\s,] plot
({(cos(90) + \s*cos(210) + \s*\s*cos(330))/(1+\s+\s*\s) },
{(sin(90) + \s*sin(210) + \s*\s*sin(330))/(1+\s+\s*\s) });

\def \R {1/4}

\coordinate (Q1) at  ({(cos(90) + \R*\R*cos(210) + \R*cos(330))/(1+\R+\R*\R) },
{(sin(90) + \R*\R*sin(210) + \R*sin(330))/(1+\R+\R*\R) });

\coordinate (Q2) at
({(cos(90) + \R*cos(210) + \R*\R*cos(330))/(1+\R+\R*\R) },
{(sin(90) + \R*sin(210) + \R*\R*sin(330))/(1+\R+\R*\R) });


\coordinate (Q3) at ({(\R*cos(90) + cos(210) + \R*\R*cos(330))/(1+\R+\R*\R) },
{(\R*sin(90) + sin(210) + \R*\R*sin(330))/(1+\R+\R*\R) });
\coordinate (Q4) at ({(\R*\R*cos(90) + cos(210) + \R*cos(330))/(1+\R+\R*\R) },
{(\R*\R*sin(90) + sin(210) + \R*sin(330))/(1+\R+\R*\R) });

\coordinate (Q5) at ({(\R*\R*cos(90) + \R*cos(210) + cos(330))/(1+\R+\R*\R) },
{(\R*\R*sin(90) + \R*sin(210) + sin(330))/(1+\R+\R*\R) });

\coordinate (Q6) at ({(\R*cos(90) + \R*\R*cos(210) + cos(330))/(1+\R+\R*\R) },
{(\R*sin(90) + \R*\R*sin(210) + sin(330))/(1+\R+\R*\R) });

\draw[dashed] (Q1) -- (Q2) -- (Q3) -- (Q4) -- (Q5) -- (Q6) -- (Q1);

\def \R {0.3}

\coordinate (S1) at  ({(cos(90) + \R*\R*cos(210) + \R*cos(330))/(1+\R+\R*\R) },
{(sin(90) + \R*\R*sin(210) + \R*sin(330))/(1+\R+\R*\R) });

\coordinate (S2) at
({(cos(90) + \R*cos(210) + \R*\R*cos(330))/(1+\R+\R*\R) },
{(sin(90) + \R*sin(210) + \R*\R*sin(330))/(1+\R+\R*\R) });


\coordinate (S3) at ({(\R*cos(90) + cos(210) + \R*\R*cos(330))/(1+\R+\R*\R) },
{(\R*sin(90) + sin(210) + \R*\R*sin(330))/(1+\R+\R*\R) });
\coordinate (S4) at ({(\R*\R*cos(90) + cos(210) + \R*cos(330))/(1+\R+\R*\R) },
{(\R*\R*sin(90) + sin(210) + \R*sin(330))/(1+\R+\R*\R) });

\coordinate (S5) at ({(\R*\R*cos(90) + \R*cos(210) + cos(330))/(1+\R+\R*\R) },
{(\R*\R*sin(90) + \R*sin(210) + sin(330))/(1+\R+\R*\R) });

\coordinate (S6) at ({(\R*cos(90) + \R*\R*cos(210) + cos(330))/(1+\R+\R*\R) },
{(\R*sin(90) + \R*\R*sin(210) + sin(330))/(1+\R+\R*\R) });

\draw[dashed] (S1) -- (S2) -- (S3) -- (S4) -- (S5) -- (S6) -- (S1);

\def \R {0.5}

\coordinate (P1) at  ({(cos(90) + \R*\R*cos(210) + \R*cos(330))/(1+\R+\R*\R) },
{(sin(90) + \R*\R*sin(210) + \R*sin(330))/(1+\R+\R*\R) });

\coordinate (P2) at
({(cos(90) + \R*cos(210) + \R*\R*cos(330))/(1+\R+\R*\R) },
{(sin(90) + \R*sin(210) + \R*\R*sin(330))/(1+\R+\R*\R) });


\coordinate (P3) at ({(\R*cos(90) + cos(210) + \R*\R*cos(330))/(1+\R+\R*\R) },
{(\R*sin(90) + sin(210) + \R*\R*sin(330))/(1+\R+\R*\R) });
\coordinate (P4) at ({(\R*\R*cos(90) + cos(210) + \R*cos(330))/(1+\R+\R*\R) },
{(\R*\R*sin(90) + sin(210) + \R*sin(330))/(1+\R+\R*\R) });

\coordinate (P5) at ({(\R*\R*cos(90) + \R*cos(210) + cos(330))/(1+\R+\R*\R) },
{(\R*\R*sin(90) + \R*sin(210) + sin(330))/(1+\R+\R*\R) });

\coordinate (P6) at ({(\R*cos(90) + \R*\R*cos(210) + cos(330))/(1+\R+\R*\R) },
{(\R*sin(90) + \R*\R*sin(210) + sin(330))/(1+\R+\R*\R) });

\draw[dashed] (P1) -- (P2) -- (P3) -- (P4) -- (P5) -- (P6) -- (P1);

\def \R {0.4}

\coordinate (S1) at  ({(cos(90) + \R*\R*cos(210) + \R*cos(330))/(1+\R+\R*\R) },
{(sin(90) + \R*\R*sin(210) + \R*sin(330))/(1+\R+\R*\R) });

\coordinate (S2) at
({(cos(90) + \R*cos(210) + \R*\R*cos(330))/(1+\R+\R*\R) },
{(sin(90) + \R*sin(210) + \R*\R*sin(330))/(1+\R+\R*\R) });


\coordinate (S3) at ({(\R*cos(90) + cos(210) + \R*\R*cos(330))/(1+\R+\R*\R) },
{(\R*sin(90) + sin(210) + \R*\R*sin(330))/(1+\R+\R*\R) });
\coordinate (S4) at ({(\R*\R*cos(90) + cos(210) + \R*cos(330))/(1+\R+\R*\R) },
{(\R*\R*sin(90) + sin(210) + \R*sin(330))/(1+\R+\R*\R) });

\coordinate (S5) at ({(\R*\R*cos(90) + \R*cos(210) + cos(330))/(1+\R+\R*\R) },
{(\R*\R*sin(90) + \R*sin(210) + sin(330))/(1+\R+\R*\R) });

\coordinate (S6) at ({(\R*cos(90) + \R*\R*cos(210) + cos(330))/(1+\R+\R*\R) },
{(\R*sin(90) + \R*\R*sin(210) + sin(330))/(1+\R+\R*\R) });

\draw (S1) -- (S2) -- (S3) -- (S4) -- (S5) -- (S6) -- (S1);


\draw[scale=1, domain=0:1, smooth, variable=\s, name path = A] plot
({(cos(90) + \s*\s*cos(210) + \s*cos(330))/(1+\s+\s*\s) },
{(sin(90) + \s*\s*sin(210) + \s*sin(330))/(1+\s+\s*\s) });

\draw[scale=1,  domain=0:1, smooth, variable=\s, name path = B] plot
({(cos(90) + \s*cos(210) + \s*\s*cos(330))/(1+\s+\s*\s) },
{(sin(90) + \s*sin(210) + \s*\s*sin(330))/(1+\s+\s*\s) });

\def \R {0.4}
\fill [red, opacity = 0.3,  domain=\R:1, smooth, variable=\s,]

plot ({(cos(90) + \s*\s*cos(210) + \s*cos(330))/(1+\s+\s*\s) },
{(sin(90) + \s*\s*sin(210) + \s*sin(330))/(1+\s+\s*\s) })

-- (0,0) -- 
({(cos(90) + \R*cos(210) + \R*\R*cos(330) + cos(90) + \R*\R*cos(210) + \R*cos(330))/(2+2*\R+2*\R*\R) },
{(sin(90) + \R*sin(210) + \R*\R*sin(330) + sin(90) + \R*\R*sin(210) + \R*sin(330))/(2+2*\R+2*\R*\R) })

-- ({(cos(90) + \R*\R*cos(210) + \R*cos(330))/(1+\R+\R*\R) },
{(sin(90) + \R*\R*sin(210) + \R*sin(330))/(1+\R+\R*\R) });

\fill [red, opacity = 0.3,  domain=\R:1, smooth, variable=\s,]

plot ({(cos(90) +\s* cos(210) +  \s*\s* cos(330))/(1+\s+\s*\s) },
{(sin(90) + \s* sin(210) + \s*\s* sin(330))/(1+\s+\s*\s) })

-- (0,0) -- 
({(cos(90) + \R*cos(210) + \R*\R*cos(330) + cos(90) + \R*\R*cos(210) + \R*cos(330))/(2+2*\R+2*\R*\R) },
{(sin(90) + \R*sin(210) + \R*\R*sin(330) + sin(90) + \R*\R*sin(210) + \R*sin(330))/(2+2*\R+2*\R*\R) })

-- ({(cos(90) +  \R* cos(210) + \R*\R* cos(330))/(1+\R+\R*\R) },
{(sin(90) + \R* sin(210) +\R*\R* sin(330))/(1+\R+\R*\R) });

\coordinate (U) at  ({(cos(90) + \R*\R*cos(210) + \R*cos(330))/(1+\R+\R*\R) },
{(sin(90) + \R*\R*sin(210) + \R*sin(330))/(1+\R+\R*\R) });

\coordinate (V) at ({(cos(90) +  \R* cos(210) + \R*\R* cos(330))/(1+\R+\R*\R) },
{(sin(90) + \R* sin(210) +\R*\R* sin(330))/(1+\R+\R*\R) });

\draw [thick, red] 
(U) -- (V);

\tikzfillbetween[of=A and B]{ pattern=north east lines, pattern color=gray};

\node  at (-0.8, 0.5) {\color{red}{$\sface_{\pi, t}$}};
\node  at (-0.8, 0.2) {\color{red}{$R_{\pi,t}$}};
\node  at (-0.8, 0.8) {\color{black}{$R_{\pi,\infty}$}};

\fill [pattern=north east lines, pattern color=gray, domain=0:1, smooth, variable=\s,]

plot ({(\s*\s*cos(90) + \s*cos(210) + cos(330))/(1+\s+\s*\s) },
{(\s*\s*sin(90) + \s*sin(210) + sin(330))/(1+\s+\s*\s) })

plot ({((1 -\s)*(1-\s)*cos(90) + cos(210) + (1-\s)*cos(330))/(1+(1-\s)+(1-\s)*(1-\s)) },
{((1-\s)*(1-\s)*sin(90) + sin(210) + (1-\s)*sin(330))/(1+(1-\s)+(1-\s)*(1-\s)) })

-- (C);

\fill [blue, opacity = 0.3, domain=\R:1, smooth, variable=\s,]

plot ({(\s*\s*cos(90) + \s*cos(210) + cos(330))/(1+\s+\s*\s) },
{(\s*\s*sin(90) + \s*sin(210) + sin(330))/(1+\s+\s*\s) })

plot ({( (\R + 1 -\s) *(\R + 1 -\s)*cos(90) + cos(210) + (\R + 1 -\s)*cos(330))/(1+(\R + 1 -\s)+(\R + 1 -\s)*(\R + 1 -\s)) },
{((\R + 1 -\s)*(\R + 1 -\s)*sin(90) + sin(210) + (\R + 1 -\s)*sin(330))/(1+(\R + 1 -\s)+(\R + 1 -\s)*(\R + 1 -\s)) })
-- (S5);

\draw [thick, blue] 
(S4) -- (S5);

\node  at (1.2, - 0.2) {\color{blue}{$R_{\bar \pi, t}$}};
\node  at (1.2, - 0.35) {\color{blue}{$\sface_{\bar \pi, t}$}};
\node  at (1.2, - 0.5) {\color{black}{$R_{\bar{\pi},\infty}$}};

\end{tikzpicture}
  
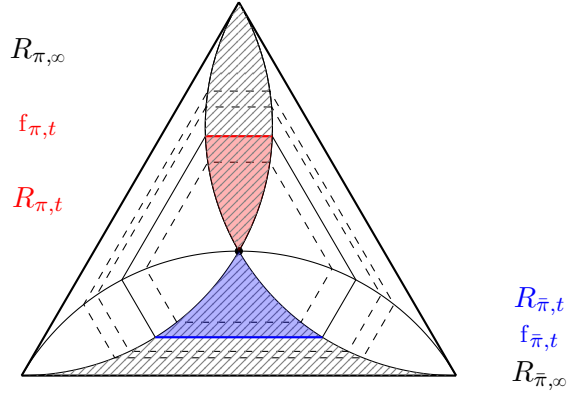
\captionof{figure}{Two regions in the simplex, induced by two faces with partitions $\pi$ and $\bar \pi$.}
  \label{fig:RegionsQt}
\end{minipage}
\end{center}

\subsection{Auxiliary homeomorphism $\iota_t$} All our constructions use the homeomorphisms
\[
\begin{array}{cccc}
\iota_t &\colon [1,t] &\longrightarrow &[1, \infty] \\[2 mm]
&s &\mapsto & (1+\frac{s-1}{t-s}) s
 \end{array}
\]
for $t\in (1, \infty]$. For a fixed $t$, we will sometimes denote by $\tilde s$ the real number $\iota_t(s)$.

\subsection{Construction of the maps $h_t$ and $\hbar_t$} \label{ss:ConstructionHomeomorphisms}
We proceed to introduce the homeomorphisms $h_t\colon P_t \to \cancomp\eta^\trop$ and $\hbar_t \colon Q_t  \to \cancomp \sigma^\trop$ by {\em recursive unfolding}. By construction, the maps will satisfy the following \emph{symmetry property}:

\begin{itemize}
\item for each $s\in (1, t]$,  the image by $h_t$ of $\partial_\infty P_s$ in the compactified cone $\cancomp\eta^\trop$ coincides with $\partial_\infty P_{\tilde s}$ with $\tilde s = \iota_t(s)$.

\item for each $s\in (1, t]$,  the image by $\hbar_t$ of $\partial_\infty Q_s$ in the compactified simplex $\cancomp\sigma^\trop$ coincides with $\partial_\infty Q_{\tilde s}$ with $\tilde s = \iota_t(s)$.
\end{itemize}

\subsubsection{Definition of the map $\hbar_t$ by unfolding}

In order to define the homeomorphism $\hbar_t : Q_t \to \cancomp\sigma^\trop$, we follow the decomposition idea outlined above. Fix a face $\sface_{\pi,t}$ in $\partial_\infty Q_t$ corresponding to an ordered partition $\pi = (\pi_1, \dots, \pi_r, \pi_\fin = \varnothing)$ of $[n]$. Consider the region $R_{\pi, t}$ in $Q_t \setminus \{\one_\sigma\}$ given by the disjoint union
\[
	R_{\pi, t} = \bigcup_{s \in(1, t]} \sface_{\pi, s}.
\]
By Lemma~\ref{lem:RecursiveStructureFaces}, each face $\sface_{\pi,s}$, $s<  t$ is a product $\sface_{\pi,s} = A_{1,s} \times \dots \times A_{r,s}$ and each set $A_{j,s}$ is isomorphic to $Q^{m_j}_{\pi_j,s}$ under the projection map $\pr_{j,s}$. Since the permutohedra $Q^{m_j}_{\pi_j,t}$ are of lower dimension, we have the recursively defined homeomorphisms $\hbar_{j,t}^{m_j} \colon Q^{m_j}_{\pi_j, t} \to \cancomp \sigma_{\pi_j}^\trop$.

\smallskip

We will define $\hbar_t$ on ${R_{\pi,t}}$ by {\em unfolding} the maps $\hbar^{m_j}_{j,t}$ to $R_{\pi, t}$, that is, as the product \smallskip
\[\hbar_t\rest{R_{\pi,t}}  = \hat\hbar^{m_j}_{1, t}\times \dots \times \hat\hbar^{m_r}_{r, t}\]
where each $\hat\hbar^{m_j}_{j,t}$ is obtained by unfolding $\hbar^{m_j}_{j,t}$. 
Namely, for a point $x$ in $\sface_{\pi, s}$, $s \in (1,  t]$,
\[\hat\hbar^{m_j}_{j, t} (x) := \big ( \pr^{-1}_{j,\tilde s}\circ\hbar^{m_j}_{j,t}\circ\pr_{j,s}  \big )  (x_{\pi_j}) \in A_{j, \tilde s}\]
where $\tilde s =\iota_t(s)$ and $x_{\pi_j} := (x_k)_{k \in \pi_j}$. Notice that the latter is well-defined by the above {\em symmetry property} of $\hbar^{m_j}_{j,t}$. Moreover, $\hbar_t^{m_j}$ maps $\sface_{\pi,s}$ bijectively onto $\sface_{\pi, \tilde s}$ and, altogether, we get a bijection $\hbar_t\rest{R_{\pi,t}}  \colon R_{\pi,t} \to R_{\pi, \infty}$.

\smallskip

Setting additionally $\hbar_t(\one_\sigma) = \one_\sigma$, all these maps glue to a global map 
\[
\hbar_t \colon Q_t = \{ \one_\sigma\} \cup \bigcup_{\pi}R_{\pi,t} \to \cancomp\sigma^\trop.
\]
Moreover, $\hbar_t$ satisfies again the above symmetry property, respects the combinatorics of $Q_t$ and $\cancomp\sigma^\trop$ and is a homeomorphism (the latter follows from Theorem~\ref{thm:global_continuity}).

\subsubsection{Definition of the map $h_t$ by unfolding}
In the case of cones, we proceed analogously. Fix a face $\face_{\pi,t}$ in the boundary $\partial_\infty P_t$ with ordered partition $\pi = (\pi_1, \dots, \pi_r, \pi_\fin)$ of $S \subseteq [n]$. Define the region $R_{\pi,t}$ of $\cancomp\eta^\trop\setminus \cube_\eta$ as the union of the faces $\face_{\pi, s}$, $s \in (1, t]$.

\smallskip

Denote by $\hbar^{m_1}_{1,t}, \dots, \hbar^{m_r}_{r,t}$ the maps defined by induction for $Q^{m_1}_{\pi_1, t}, Q^{m_2}_{\pi_2, t},\dots, Q^{m_r}_{\pi_r, t}$. Moreover, let $h_{\fin,t}^{m_\fin} \colon P_{\pi_\fin, t}^{m_\fin} \to \cancomp \eta_{\pi_\fin}^\trop$ be the analogous map for the lower dimensional polytope $P_{\pi_\fin, t}^{m_\fin}$. Here, $m_1, \dots, m_r, m_\fin$ are defined as in Lemma~\ref{lem:RecursiveStructureFaces}.

\smallskip

On the region $R_{\pi,t}$, we define the map $h_t$ as
\begin{equation} \label{eq:htProduct}
h_t\rest{R_{\pi,t}} := \hat{\hbar}^{m_1}_{1,t}\times \dots\times \hat{\hbar}^{m_r}_{r,t}\times \hat{h}^{m_\fin}_{\fin,t} \times 0^{n - |S|}
\end{equation}
where $\hat{\hbar}^{m_j}_{j,t}$ and $\hat{h}^{m_\fin}_{\fin,t}$ are obtained by unfolding the previously defined maps. More precisely, its last part is given by (notice that this map is well-defined by Lemma~\ref{lem:RecursiveStructureFaces})
\[
\hat{h}^{m_\fin}_{\fin,t} (x) := {h}^{m_\fin}_{\fin,t} (x_{\pi_\fin}), \qquad x \in R_{\pi,t},
\]
where $x_{\pi_\fin} := (x_k)_{k \in \pi_\fin}$. Similar to the case of simplices, the unfolded maps $ \hat{\hbar}^{m_j}_{j,t}$, $j \in [r]$, are
\[
\hat\hbar_{j, t}^{m_j} (x) := \big ( \pr^{-1}_{j,\tilde s}\circ\hbar^{m_j}_{j,t}\circ\pr_{j,s}  \big )  (x_{\pi_j}), \qquad x \in \face_{\pi,s},
\]
with $x_{\pi_j} := (x_k)_{k \in \pi_j}$ and $\tilde s =\iota_t(s)$. Here, we have again used the isomorphisms $\pr_{j,s} \colon A_{j,s} \to Q^{m_j}_{\pi_j, s}$ in the decomposition $\face_{\pi,s} = A_{1,s} \times \dots A_{\fin,s} \times 0^{n - |S|}$ of the face $\face_{\pi,s}$ (see Lemma~\ref{lem:RecursiveStructureFaces}). Notice that $h_t$ maps $F_{\pi,s}$ bijectively onto $F_{\pi, \tilde s}$.

\smallskip

Setting additionally $h_t(x) = x$ for $x \in \cube_\eta$, these maps glue to a global map
\[
h_t \colon P_t = \cube_\eta \cup \bigcup_{\pi}R_{\pi,t} \to \cancomp\eta^\trop.
\]
Analogous to the case of simplices, $h_t$ is a homeomorphism with the above symmetry property and respects the combinatorics of $P_t$ and $\cancomp \eta^\trop$ (see Theorem~\ref{thm:global_continuity}).

\subsubsection{Examples} For illustration, we describe the above maps in low dimensions for $\Phi_j(t) =t^j$, $j=1, 2$.
\begin{itemize} \item [(i)] In dimension $n = 1$, we trivially have $\inn {\sigma_1} = \{ 1 \} = \cancomp \sigma_1^\trop$ and $Q_t = \{ 1 \} = \cancomp \sigma_1^\trop$ for all $t \in (1, \infty]$. Hence $\hbar_t(1) := 1 $ is simply the trivial map on $\inn {\sigma_1}$.
\item [(ii)] For the one-dimensional cone $\eta = \R_+$, the generalized permutohedron is the interval $P_t = [0, t]$. Note that $\partial_\infty P_t$ has exactly one face $\face_t = \{t\}$. Using the map $\hbar_t$ of the simplex $\inn \sigma_1$, it follows that $\hbar_t(x) = \iota_t(x)$ for $x$ in $(1, t]  = P_t \setminus\cube$ and $h_t(x) = x$ for $x \in [0, t]$.
\item [(iii)] For the simplex $\inn \sigma_2 \subseteq \R_+^2$, the permutohedron $Q_t$ is the line segment
\[
Q_t =\textrm{conv-hull} \Big \{ \varphi(t) \big (1, t \big), \varphi(t) \big(t, 1 \big) \Big\}, \qquad  \varphi(t) := \frac{1}{1+t}.
\]
The boundary $\partial_\infty Q_t$ has two faces $\sface_{\pi, t} = \{ \varphi(t) \big(1, t \big)\}$ with $\pi = (\{2\}, \{1\}, \varnothing )$ and $\sface_{\bar{\pi}, t} = \{  \varphi(t) \big(t, 1 \big)\}$ with $\bar{\pi} = (\{1\}, \{2\}, \varnothing )$. The map $\hbar_t$ is given by
\[
\hbar_t (x) = \begin{cases} \varphi( \tilde s ) \big(1, \tilde s \big) & \text{if } x = \varphi(s) \big(1, s \big) \text{ for } s \in [1, t] \\[2mm]
\varphi( \tilde s ) \big(\tilde s, 1 \big) & \text{if } x = \varphi(s) \big(s, 1 \big) \text{ for }  s \in [1, t] \end{cases},
\]
where $\tilde s := \iota_t(s)$.  Note that $\hbar_t$ is not linear (as a map between line segments).
\end{itemize}

\subsection{The tropical log map} 
Next, we prove continuity of the global functions defined by $h_t$ and $\hbar_t$ (see \eqref{eq:GlobalFunctions}). Let $\mathcal P$ and $\mathcal Q$ be the global sets in \eqref{eq:GlobalSets}. 

\begin{thm}\label{thm:global_continuity} The globally defined map $h\colon \mathcal P \to \cancomp\eta^\trop$ which restricts to $h_t$ on $P_t \times \{t\}$ is continuous. Similarly, the globally defined map $\hbar\colon \mathcal Q \to \cancomp\sigma^\trop$ is continuous.\smallskip\newline
Here, by convention $h(x, \infty) := x$ for $x \in \cancomp\eta^\trop = P_\infty$ and $\hbar(x, \infty) := x$ for $x \in \cancomp\sigma^\trop = Q_\infty$.
\end{thm}
\begin{proof} 
In the following, we prove the continuity of the global map $\hbar\colon \mathcal Q \to \cancomp\sigma^\trop$. The proof for the global map $h\colon \mathcal P \to \cancomp\eta^\trop$ is analogous. Let $(x_n, t_n)_{n \in \N}$ be a sequence in $\mathcal Q$ converging to a point $(x, t)$ in $\mathcal Q$.  
We prove that $\lim_{n \to \infty} \hbar(x_n, t_n) = \hbar(x,t)$ in $\cancomp \sigma\trop$.

\smallskip

We proceed by case distinction and use induction on the dimension of the simplex $\inn \sigma$. If the limit $x$ is in the interior $Q_{t} \setminus \partial_\infty Q_{t}$ of the polytope $Q_t$, then the convergence reduces to the continuity of the global map for lower dimensional simplices. The same reasoning applies if $x$ and all $x_n$'s are boundary points of the polytopes, that is, $x \in \partial_\infty Q_{t}$ and $x_n \in \partial_\infty Q_{t_n}$ for all $n$.

\smallskip

It remains to treat the case that $x \in \partial_\infty Q_t$ is a boundary point and all $x_n$'s are in the interior $Q_{t_n} \setminus \partial_\infty Q_{t_n}$. Notice that in this case $\hbar(x_n, t_n) \in \inn \sigma$ for all $n$. Let $\sface_{\pi, t}$ be the finest face of $\partial_\infty Q_t$  with $x \in \sface_{\pi, t}$ and $\pi = (\pi_\infty=(\pi_1, \dots, \pi_r), \pi_\fin=\varnothing)$ its ordered partition. By definition,
\[
\hbar(x,t) = \hbar_t(x) = \Big (\hbar^{m_1}_{\pi_1, t} (\bar x_{\pi_1}), \dots, \hbar^{m_r}_{\pi_r, t} (\bar x_{\pi_r}) \Big ) \text{ belongs to } \prod_{j=1}^r \inn \sigma_{\pi_j} \subseteq \cancomp \sigma^\trop,
\]
where $\bar x_{\pi_j}$, $j \in [r]$, is the projection of $x$ to $\inn \sigma_{\pi_j}$. It turns out that, as $n \to \infty$, we will get the asymptotics
\begin{equation} \label{eq:GlobalContinuityQuantitive}
\hbar(x_n,t_n)  = \Big ((1 + o(1)) \,  \frac{\Phi_{m_j}(\iota_{t_n}(s_n))}{\Phi_{m}(\iota_{t_n}(s_n))} \cdot \hbar^{m_j}_{\pi_j, t} (\bar{x}_{\pi_j}) \Big)_{j=1, \dots, r} \, \in \inn \sigma.
\end{equation}

Since $\tilde s_n := \iota_{t_n}(s_n) \to \infty$ as $n \to \infty$, by the relative behavior of the auxiliary functions at infinity, this in particular implies that $\hbar(x_n, t_n) \to \hbar(t,x)$ in $\cancomp \sigma^\trop$, as required.

Now, to prove \eqref{eq:GlobalContinuityQuantitive}, since the expression is quantitative, it can be again verified by induction on the dimension of $\inn \sigma$. For simplicity, assume first that $x_n \in \sface_{\pi, s_n}$ with $s_n < t_n$ (for all $n$) and that $\sface_{\pi, s_n}$ is the finest such face. Using the global continuity property for lower dimensions, we get
\begin{align*}
\hbar(x_n,t_n) &= \Big (\pr_{j,\tilde s_n}^{-1} \circ  \hbar^{m_j}_{\pi_j, t_n} (\bar{x}_{n, \pi_j}) \Big)_{j} = \Big (  \chi_j(\tilde s_n)L(\tilde s_n)^{-1} \cdot \hbar^{m_j}_{\pi_j, t_n} (\bar x_{n,\pi_j}) \Big)_{j} \\ &= \Big ((1 + o(1)) \,  \frac{\Phi_{m_j}(\tilde s_n)}{\Phi_{m}(\tilde s_n)} \cdot \hbar^{m_j}_{\pi_j, t} (\bar x_{\pi_j}) \Big)_j,
\end{align*}
where $\chi_j(\cdot)$ and $L(\cdot)$ are the functions from Lemma~\ref{lem:RecursiveStructureFaces} (ii). Hence, \eqref{eq:GlobalContinuityQuantitive} holds in this situation. The general case can be inferred by using \eqref{eq:GlobalContinuityQuantitive} inductively for lower dimensional simplices. The proof is complete.
\end{proof}

\begin{remark} \label{rem:Retractions}
Note that Theorem~\ref{thm:global_continuity} implies that for $t_0 > 1$, the following maps
\[
\begin{array} {cccccccccc}
G \colon &\cancomp \eta^\trop \times [t_0, \infty] &\to &\cancomp \eta^\trop & & &H \colon &\cancomp \sigma^\trop \times [t_0, \infty] &\to &\cancomp \sigma^\trop  \\[2 mm]
&(x,t)  & \mapsto &h_t^{-1}(x) & & & &(x,t) & \mapsto & \hbar_t^{-1}(x).
\end{array}
\] 
are deformation retractions from $\cancomp \eta^\trop$ and $\cancomp \sigma^\trop$ to the polytopes $P_{t_0}$ and $Q_{t_0}$, respectively.
\end{remark}

Based on the foliation properties, we can introduce the tropical log map as follows.
\begin{defi}[Tropical log map] The map $\logtrop \colon \cancomp \eta^\trop \setminus \cube_\eta \to \partial  \cancomp \eta^\trop$ is defined by 
\[\logtrop\rest{\partial_\infty P_t} := h_t\rest{\partial_\infty P_t}, \qquad t \in (1, \infty]. \]
Analogously, the map $\logtrop \colon \cancomp \sigma^\trop \setminus \{\one_\sigma\} \to \partial_\infty   \cancomp \sigma^\trop$ is defined by setting 
\[\logtrop\rest{\partial_\infty Q_t} := \hbar_t\rest{\partial_\infty Q_t}, \qquad t \in (1, \infty]. \]
\end{defi}

\begin{thm} The tropical log maps $\logtrop$ for $\cancomp \eta^\trop$ and $\cancomp \sigma^\trop$ are retractions. 
\end{thm}
\begin{proof} This follows from the continuity properties  stated in Theorem~\ref{thm:global_continuity}.
\end{proof}

\subsection{Log maps on fans} Finally, we describe how to globalize the above constructions to fans. This will be used later in Section~\ref{sec:tropical_log_map2} to get a log map in the setting of tropical moduli spaces.

\smallskip

Suppose that $\Sigma$ is  a simplicial fan in $\R_+^n$ or an abstract simplicial fan. In the latter case, we choose a generating vector $e_\rho$ for each ray $\rho$ of $\Sigma$ so that we can naturally identify each separate face $\eta $ with $ \R_+^{n_\eta}$ with $n_\eta = \dim\eta$. We choose $m \geq \dim \Sigma$ and auxiliary functions $\Phi_1, \dots, \Phi_m$ as in the previous section. 

 We have constructed a log map using the generalized permutohedra $P^m_t$. More generally, for each cone $\eta$ in $\Sigma$, we consider the sequence of generalized permutohedra $P_t^m$ in $\eta$, and  obtain in this way,  a log map $\logtropind{\eta, m} \colon \cancomp\eta^\trop\setminus \cube_\eta \to \partial \cancomp\eta^\trop$.

As is evident from \eqref{eq:AlternativeVerteCurves} and \eqref{eq:htProduct}, the restriction of the log map $\logtropind{\eta,m}$ on a face $\eta$ to a subface $\tau \subseteq \eta$ coincides with the log map $\logtropind{\tau, m}$ on $\tau$. In particular, the log maps of all the cones $\eta$ in the fan glue together to give a global log map
\[\begin{array}{cccc}\logtrop = \logtropind{m} \colon &\cancomp\Sigma^{\trop}\setminus\bigcup_{\eta \in \Sigma} \cube_\eta &\to &\partial \cancomp\Sigma^{\trop} 
\end{array}. \]

\begin{thm} The tropical log map $\logtrop$ for $\cancomp \Sigma^\trop$ constructed above is a retraction. 
\end{thm}
\begin{proof} This is immediate from the above discussion. 
\end{proof}

\newpage


\part{Tropical curves} \label{part:TropicalCurves}
The aim of this section is to develop a  function theory in higher rank tropical geometry. We define a notion of Laplacian on higher rank tropical curves, formulate the tropical Poisson equation, and define higher rank Green functions associated to layered measures on tropical curves.  With the help of tropical Green functions, and appropriate log maps, we then describe the behavior of Green functions on metric graphs near the boundary of their corresponding moduli spaces, introduced as well in this section.

\section{Moduli space of higher rank tropical curves} \label{sec:tropical_moduli}

This section introduces the \emph{moduli space of stable tropical curves of genus $g$ with $n$ marked points} that we will denote by $\mgtrop{\grind{g,n}}$. 
Let $\mggraph{\grind{g,n}}$ be the moduli space of stable metric graphs of genus $g$ with $n$ markings and denote by $\cancomp{\mggraph{\grind{g,n}}}$ its canonical compactification, both studied in ~\cite{ACP}. The moduli space $\mgtrop{\grind{g,n}}$ defined in this section shall be regarded as \emph{the absolute canonical compactification} of $\mggraph{\grind{g,n}}$, and together  with the space $\cancomp{\mggraph{\grind{g,n}}}$ belongs to the bi-indexed family of higher rank canonical compactifications $\mgbartropr{\grind{g,n}}{a,b}$ of $\mggraph{\grind{g,n}}$ (in which $\cancomp{\mggraph{\grind{g,n}}}$ is identified with $\mgbar_{g,n}^{\tropr{0,1}}$).

\subsection{Stable metric graphs of given combinatorial type} \label{sec:cone_metric_graph_given_genus}
Let $G = (V, E, \genusfunction, \marking)$ be a stable graph of genus $g$ with $n$ marked points. The \emph{open cone metrics on $G$} that we call as well \emph{the open cone of stable metric graphs of combinatorial type $G$} is the open simplicial cone $\inn\eta_{\grind{G}} := \R_{>0}^E$. The closed cone $\eta_{\grind{G}} =\R_{\geq 0}^E$ is the disjoint union of $\inn\eta_{\grind{H}}$ for $H$ a stable graph of genus $g$ with $n$ markings obtained by contracting a subset of edges in $G$.

\subsection{Tropical curves of given combinatorial type}\label{sec:moduli_space_tropical_curves}
For a given stable graph with marking $G = (V, E, \genusfunction, \marking)$, we denote by $\cancomp\eta_{\grind{G}}^\trop$ the tropical compactification of $\eta_{\grind{G}}$. Using the description of the tropical compactifications, and given that the cone $\eta_{\grind{G}}$ is simplicial, we obtain the following description of $\cancomp\eta_{\grind{G}}^\trop$
\[
\cancomp \eta_{\grind{G}}^\trop = \bigsqcup_{F \subseteq E}  \bigsqcup_{\substack{\pi \in \Piall(F)}}  \inn\keg_\pi 
\]
where for $\pi=(\pi_\infty=(\pi_1, \dots, \pi_r), \pi_\fin)$, we have $\inn\keg_\pi = \inn \sigma_{\pi_\infty} \times \inn\keg_{\pi_\fin}$ and $\inn\keg_{\pi_\fin} = \R_{+}^{\pi_\fin}$. As usual, $\pi_\infty$ is the sedentarity part of $\pi$, $r\in \N \cup\{0\}$ is the rank of $\pi$, and the finite part $\pi_\fin$ is allowed to be empty (in this case, $\R_{+}^{\varnothing} =\zerocone =(0)$ by convention).

\smallskip

The \emph{space of tropical curves of combinatorial type $G$} denoted by $\mgtropcombin{\grind{G}}$ is defined by
\[ \mgtropcombin{\grind{G}} := \bigsqcup_{\substack{\pi \in \Piall(E)}}  \inn\keg_\pi 
\]
that we endow with the subspace topology from $\cancomp\eta_{\grind{G}}^\trop$.  We moreover define for an ordered partition $\pi=(\pi_\infty, \pi_\fin)$ of $E$, the subspace
\[\mgtropcombin{\combind{(G, \subface \pi)}} := \bigsqcup_{\substack{\pi' \in \Pihat(E) \\ \pi'\subface \pi}}  \inn\keg_{\pi'}. \]

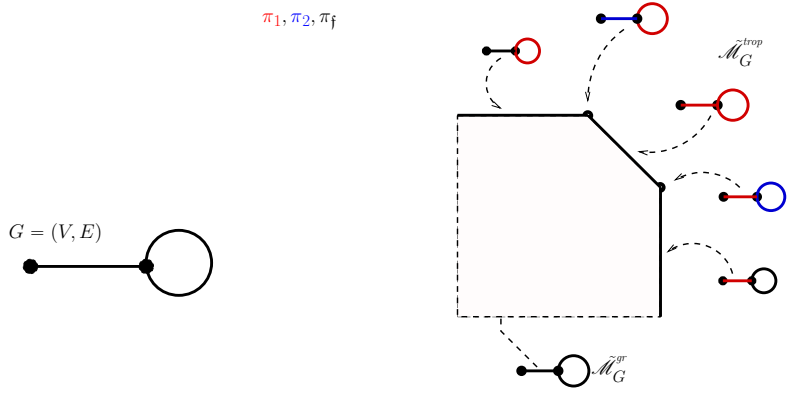
\begin{figure}[!t]
\centering
    \scalebox{.3}{\input{example4.tikz}}
\caption{Six types of tropical curves in the compactified space of tropical curves of the given combinatorial type depicted on the left.}
\label{fig:moduli_comb}
\end{figure}

In the sequel, for  consistency with the above terminology, the stratum $\inn\keg_\pi$ of $\mgtropcombin{\grind{G}}$ representing all tropical curves of combinatorial type $(G, \pi)$ will be denoted by $ \mgtropcombin{\combind{(G,\pi)}}$. In the case $\pi_\infty =\emptyset$, that is for the open cone $\inn\eta_{\grind{G}}$ consisting of metric graphs of combinatorial type $G$ we simply use the notation $\mggraphcombin{\grind{G}}$.

\smallskip

\begin{figure}[!t]
\centering
    \scalebox{.3}{\input{example6.tikz}}
\caption{Tropical curves of full sedentarity of the same combinatorial type in the canonically compactified cone of Figure~\ref{fig:compactified_cone}. The underlying graph $G$ has two vertices and three edges, two loops at vertices, and an extra edge joining them.}
\label{fig:full_sedentarity}
\end{figure}
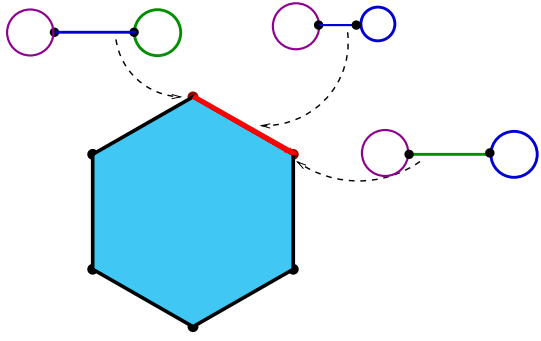

The space $\mgtropcombin{\grind{G}}$ has a \emph{central subspace}, which consists of the points representing tropical curves of full sedentarity. This is a compact space which coincides with the tropical compactification $\cancomp\sigma_{\grind{G}}^\trop$ of the projectivized cone $\inn\sigma_{\grind{G}} = \inn\eta_{\grind{G}} /\R_{+}$. We note likewise that $\cancomp\sigma_{\grind{G}}^\trop$ coincides with the space introduced in~\cite{AN} that we denote $\tilde \mg_{\grind{G},\fsed}^{\trop} = \bigsqcup_{\pi\in \Pi(E)}\inn\sigma_\pi$. The example of a graph with three edges is shown in Figure~\ref{fig:full_sedentarity}.

\smallskip

We therefore get the following stratified description of the compactified cone $\cancomp \eta_{\grind{G}}^\trop$:
\[\cancomp \eta_{\grind{G}}^\trop  =\bigsqcup_{H \subface G} \mgtropcombin{H}\]
where $H\subface G$ is understood in the sense of Section~\ref{sec:preliminaries}. In particular, we have the inclusion $\cancomp \eta_{\grind{H}}^\trop \subseteq \cancomp \eta_{\grind{G}}^\trop$ provided that $H\subface G$. Moreover, we have an inclusion 
\[\bigsqcup_{H \subface G} \cancomp\sigma_{\grind{H}}^\trop \hookrightarrow \cancomp\eta_{\grind{G}}^\trop\]
which identifies $\bigsqcup_{H \subface G} \cancomp\sigma_{\grind{H}}^\trop$ as the full sedentarity part of  $\cancomp\eta_{\grind{G}}^\trop$.

 \subsection{The universal curve $\unicurvetrop{\grind{G}}$ over $\cancomp\eta_{\grind{G}}^\trop$} \label{sec:univresal_curve_G} We denote by $\unicurvetrop{\grind{G}}$ the corresponding \emph{universal tropical curve} which on each space $\inn\eta_{\grind{H}}$ restricts to \emph{the universal tropical curve of combinatorial type $H$}.  With this terminology, we get the compatibility $\unicurvetrop{\grind{G}}\rest{\cancomp\eta_{\grind{H}}^\trop} = \unicurvetrop{\grind{H}}$, for a pair of stable marked graphs $H \subface G$.
 
 We also denote by $\unicurvetrop{\grind{G}}$ and $\unicurvetrop{\combind{(G, \subface \pi)}}$ the corresponding universal curve over $\mgtropcombin{\grind{G}}$ and $\mgtropcombin{\combind{(G, \subface \pi)}}$, respectively. The universal metric graph of combinatorial type over $\mggraphcombin{\grind{G}}$ is denoted by $\unicurvetrop{\grind{G}}$.

\subsection{Moduli space of tropical curves of combinatorial type $G$} Let $\aut(G)$ be the automorphism group of the stable marked graph $G$. The group $\aut(G)$ naturally acts on $\mgtropcombin{\grind{G}}$. We define the \emph{moduli space of tropical curves of combinatorial type $G$} denoted by $\mg^\trop_{\grind{G}}$ by taking the quotient 
\[\mgtrop{\grind{G}} := \rquot{\mgtropcombin{\grind{G}}}{\aut(G)}.\]
We moreover define 
\[\mgtrop{\combind{(G,\subface \pi)}} := \rquot{\mgtropcombin{\combind{(G, \subface \pi)}}}{\aut(G)}.\]

\subsection{Tropical moduli space} \label{ss:TropicalModuliSpace}
We now define the (coarse) moduli space $\mgtrop{\grind{g,n}}$ of tropical curves of genus $g$ with $n$ markings as the direct limit of the diagrams of inclusions
\[\cancomp\eta_{\grind{H}}^\trop \hookrightarrow \cancomp\eta_{\grind{G}}^\trop\]
for pairs $H\subface G$ of stable curves of genus $g$ with $n$ markings, similar to the construction of the tropical moduli space $\mggraph{\grind{g,n}}$ and $\cancomp{\mggraph{\grind{g,n}}}$ in~\cite{ACP}. We endow $\mgtrop{\grind{g,n}}$ with the topology induced by the ones on $\cancomp\eta_{\grind{G}}^\trop$ as the corresponding quotient topology on the limit. 

\smallskip

For each stable graph $G$ of genus $g$ with $n$ marked points, we get a canonical map $\mgtropcombin{\grind{G}} \to \mgtrop{\grind{g,n}}$. The universal tropical curve $\unicurvetrop{\grind{g,n}}$ is defined over these \emph{charts}, so that the maps $\mgtropcombin{\grind{G}}\to \mgtrop{\grind{g,n}}$ shall be regarded as \emph{tropical \'etale charts} underlying the moduli \emph{stack} of stable tropical curves of genus $g$ with $n$ marked points.

\medskip

The following theorem summarizes basic properties of the tropical moduli space $\mgtrop{\grind{g,n}}$.

\begin{thm} Let $g,n$ be a pair of non-negative integers with $3g-3+n>0$, and let $\mgtrop{\grind{g,n}}$ be the moduli space of tropical curves of genus $g$ with $n$ markings defined above. Then we have 
\begin{itemize}
\item The tropical moduli space $\mgtrop{\grind{g,n}}$ is compact.
\smallskip

\item We have a canonical embedding of the moduli space $\mggraph{\grind{g,n}}$ of stable metric graphs of genus $g$ with $n$ markings in $\mgtrop{\grind{g,n}}$. Moreover, under this embedding, $\mggraph{\grind{g,n}}$ becomes an open dense subset of $\mgtrop{\grind{g,n}}$.

\smallskip
\item For each stable marked graph $G$, the closure of the moduli space $\mgtrop{\grind{G}}$ in $\mgtrop{\grind{g,n}}$ is the disjoint union of the strata $\mgtrop{\grind{H}}$ for $H \subface G$.
\smallskip
\item We have a continuous projection map $\mgtrop{\grind{g,n}}$ to the moduli space $\cancomp{\mggraph{\grind{g,n}}}$ from~\cite{ACP}.
\end{itemize}
\end{thm}

\subsection{Family of higher rank canonical compactifications $\mgbartropr{\grind{g,n}}{a,b}$} The moduli space $\mgtrop{\grind{g,n}}$ belongs to the family of canonical compactifications $\mgbartropr{\grind{g,n}}{a,b}$ of $\mggraph{\grind{g,n}}$. The definition of $\mgbartropr{\grind{g,n}}{a,b}$ mimics that of $\mgtrop{\grind{g,n}}$ replacing each $\cancomp\eta_{\grind{G}}^\trop$ with the corresponding compactified cone $\cancomp\eta_{\grind{G}}^{\tropr{a,b}}$, for each stable graph $G$ of genus $g$ with $n$ markings. 
\smallskip

For each pair $(a,b)$ with $a+b \leq N:=3g-3+n$, we thus get the canonical compactification $\mgbartropr{\grind{g,n}}{a,b}$ of $\mggraph{\grind{g,n}}$. The family verifies the following properties.

\begin{itemize}
\item For each pair $a,b$ with $a+b=N+1$, we have $\mgbartropr{\grind{g,n}}{a,b} = \mgtrop{\grind{g,n}}$.
\smallskip

\item For two pairs $(a,b)$ and $(c,d)$ with  $c \le a$ and $d \le b$, we have a surjective forgetful map from $\mgbartropr{\grind{g,n}}{a,b}$ to $\mgbartropr{\grind{g,n}}{c,d}$.
\smallskip

\item We have $\mgbartropr{\grind{g,n}}{0,1} = \cancomp{\mggraph{\grind{g,n}}}$.
\end{itemize}

\smallskip

Altogether  we get a commutative diagram of maps depicted in Figure~\ref{fig:CompactificationTower-moduli}.

\begin{center}
\begin{figure}[h!] 
\begin{tikzpicture}[scale=0.8, every node/.style={scale=0.9}]

\foreach \x in {1,..., 2}{
	\node at (2*\x, 0) {$\mgbartropr{\grind{g,n}}{\x, 0}$};
	\draw[->>]  (2*\x, 1.25) -- (2*\x, 0.75) ;
	\draw[->>]  (2*\x + 1.25, 0) -- (2*\x + 0.75, 0) ;
}

\foreach \x in {0,..., 2}{
	\node at (2*\x, 2) {$\mgbartropr{\grind{g,n}}{\x, 1}$};
	\draw[->>]  (2*\x, 3.25) -- (2*\x, 2.75) ;
	\draw[->>]  (2*\x + 1.25, 2) -- (2*\x + 0.75, 2) ;
}

\foreach \x in {0,..., 1}{
	\node at (2*\x, 4) {$\mgbartropr{\grind{g,n}}{\x, 2}$};
	\draw[->>]  (2*\x + 1.25, 4) -- (2*\x + 0.75, 4) ;
		\draw[->>]  (2*\x, 5.25) -- (2*\x, 4.75) ;
}

	\node at (0, 9) {$\mgbartropr{\grind{g,n}}{0, N+1}$};
	\node[right] at (0.75, 9) {$= \mgtrop{\grind{g,n}}$};

	\node at (0, 7) {$\mgbartropr{\grind{g,n}}{0, N}$};
	\draw[->>]  (0, 9 - 0.75) -- (0, 9-1.25) ;
	
	\node at (2, 7) {$\mgbartropr{\grind{g,n}}{1, N}$};
	\node[right] at (2.75, 7) {$= \mgtrop{\grind{g,n}}$};

	\draw[->>]  (1.25, 7) -- (0.75, 7) ;

	\draw[->>, dotted]  (2, 7 - 0.75) -- (2, 7-2.25) ;
	\draw[->>, dotted]  (0, 7 - 0.75) -- (0, 7-2.25) ;
	\node at (9.2, 0) {$\mgbartropr{\grind{g,n}}{N+1, 0}$};
	\node[right] at (9.95, 0) {$= \mgtrop{\grind{g,n}}$};

	\node at (7, 0) {$\mgbartropr{\grind{g,n}}{N, 0}$};
	\draw[->>]  (9 - 0.75, 0) -- (9-1.25, 0) ;
	
	\node at (7, 2) {$\mgbartropr{\grind{g,n}}{N, 1} $};
	\node[right] at (7.75, 2) {$= \mgtrop{\grind{g,n}}$};

	\draw[->>]  (7, 1.25) -- (7, 0.75 ) ;

	\draw[->>, dotted]  (7 - 0.75, 2) -- (7-2.25, 2) ;
	\draw[->>, dotted]  (7 - 0.75, 0) -- (7-2.25, 0) ;

\end{tikzpicture}
\caption{The spaces $\mgbartropr{\grind{g,n}}{a,b}$ and the corresponding forgetful maps.}\label{fig:CompactificationTower-moduli}
\end{figure}
\end{center}

In this way, the double tower of tropical spaces depicted in Figure~\ref{fig:CompactificationTower-moduli} interpolates between $\mgtrop{\grind{g,n}}$ and $\cancomp{\mggraph{\grind{g,n}}}$. Moreover, the family should be compared with the tower of hybrid spaces
\begin{align}\label{eq:tower-moduli}
\mgbar_{\grind{g,n}} \longleftarrow \mgbarhybr{\grind{g,n}}{1} \longleftarrow \mgbarhybr{\grind{g,n}}{2} \longleftarrow \dots \longleftarrow \mgbarhybr{\grind{g,n}}{N}= \mg_{\grind{g,n}}^\hyb,
\end{align}
for $N = 3g-3+n$ that we introduced in~\cite{AN}. The last line in Figure~\ref{fig:CompactificationTower-moduli} is precisely the tropical analogue of the tower \eqref{eq:tower-moduli}.
We thus refer to $\mgbartropr{\grind{g,n}}{k}:=\mgbartropr{\grind{g,n}}{k,0}$ as \emph{the rank $k$ tropical compactification of $\mggraph{\grind{g,n}}$} and call it the \emph{compactified moduli space of tropical curves of rank bounded by $k$} (because of the possibility of presence of length zero edges in the intermediate rank). The \emph{moduli space of tropical curves of rank bounded by $k$} will form an open subset of $\mgbartropr{\grind{g,n}}{k}$ and is denoted by $\mgtropr{g,n}{k}$.

\begin{remark} Our notations are consistent as in the case $k=3g-3+n$, we have $\mgtropr{\grind{g,n}}{3g-3+n} = \mgbartropr{\grind{g,n}}{3g-3+n}$.
\end{remark}

\subsection{The tropical log map revisited} \label{sec:tropical_log_map2} Define \emph{the boundary at infinity} of $\mgtrop{\grind{g,n}}$ by
\[
\partial_\infty \mgtrop{\grind{g,n}}:= \mgtrop{\grind{g,n}} \setminus \mggraph{\grind{g,n}}.
\] 
Let $\umggraph{\grind{g,n}}$ be the subspace of $\mggraph{\grind{g,n}}$ parametrizing the stable marked graphs whose edge lengths are all bounded by one. That is for each stable marked graph $G$, we consider the image of the hypercube $\cube_{\eta_{\grind{G}}} \subset \eta_{\grind{G}}$ under the map $\eta_{\grind{G}} \to \mggraph{\grind{g,n}}$, and define $\umggraph{\grind{g,n}}$ as the union of these images.

\smallskip

Setting $N=3g-3+n$, picking $N$ auxiliary functions $\Phi_1, \dots, \Phi_N$ as in Section~\ref{sec:tropical_log_map} and applying the constructions of that section to each cone $\eta_{\grind{G}}$, we get a tropical log map
 \[
 \begin{array}{cccc}
 \logtrop = \logtropind{\grind{G,N}}  &\colon \eta_{\grind{G}} \setminus \cube_{\eta_{\grind{G}}} &\to & \partial_\infty \cancomp\eta_{\grind{G}}^\trop
 \end{array}.
\]
We now observe that by construction, and by compatibility of the polytopes  $P_t^N$ and $Q_t^N$ from Section~\ref{sec:tropical_log_map} under restriction to faces, the above log maps are compatible with the inclusion of cones and boundaries at infinities $\eta_{\grind{H}} \hookrightarrow \eta_{\grind{G}}$ and $\partial_\infty \cancomp\eta_{\grind{H}}^\trop \hookrightarrow \partial_{\infty}\cancomp\eta_{\grind{G}}^\trop$, for $H \subface G$. That is, we get a commutative diagram
\[
\begin{tikzcd}
\eta_{\grind{H}} \setminus \cube_{\eta_{\grind{H}}} \arrow[r,"\logtropdiag"] \arrow[hookrightarrow]{d} & \partial_\infty \cancomp\eta_{\grind{H}}^\trop\arrow[hookrightarrow]{d}\\
\eta_{\grind{G}} \setminus \cube_{\eta_{\grind{G}}} \arrow[r,"\logtropdiag"]  & \partial_\infty \cancomp\eta_{\grind{G}}^\trop.
\end{tikzcd}
\]

Passing to the limit, we therefore get a well-defined tropical log map
\[\logtrop =  \logtropind{\grind{N}}\colon \mgtrop{\grind{g,n}} \setminus \umggraph{\grind{g,n}} \to \partial_\infty \mgtrop{\grind{g,n}}. \]

\smallskip

Analogously, for every fixed stable marked graph $G$ with edge set $E$, we can define a tropical log map
\[\logtrop =\logtropind{\grind{N}} \colon \mgtrop{\grind{G}} \setminus \umggraph{\grind{G}} \to \partial_\infty \mgtrop{\grind{G}}, \]
with $\umggraph{\grind{G}}$ the moduli space of metric graphs of combinatorial type $G$ and all edge lengths bounded by one.

\smallskip

In summary, we obtain the following result.

\begin{thm} The above  defined  tropical log maps $\logtrop$ on $\mgtrop{\grind{g,n}}$ and $\mgtrop{\grind{G}}$  are retractions to the respective boundaries at infinity. 
\end{thm}

\section{Potential theory on metric graphs} \label{sec:potential_theory_metric_graphs}
In this section, we provide basic information on the canonical measure, Laplacian and canonical Green function on metric graphs. For the remainder of this section, let $\mgr$ be the metric graph obtained as the metric realization of a pair $(G, l)$ consisting of a finite graph $G=(V, E)$ and an edge length function $l  \colon E \to (0, + \infty)$. We denote by $h$ the genus of the graph $G$. 

Unless otherwise stated, we assume that $G$ (and hence $\mgr$) is connected. 

\subsection{The canonical measure on an (augmented) metric graph}

We first recall what we mean by the \emph{canonical measure} on metric graphs. This measure was introduced by Zhang in \cite{Zhang} as a refinement of the constructions carried out by Chingburg and Rumley in\cite{CR93}, and has several different representations (see, e.g., \cite{BF11, SW19, AN}). The following definition uses the \emph{cycle space} of $G$ and highlights the analogy to the Arakelov--Bergman measure on Riemann surfaces (see \eqref{eq:cannonical_measure_RS}). In the following we do not assume the connectivity of $G$.

\medskip

Let $H := H_1(G, \Z)$ be the first homology of $G$. By definition, we have an exact sequence 
 \[ 0\to H_1(G, \Z) \longrightarrow \Z^E \longrightarrow \Z^{V} \to 0.
 \]
 For each edge $e$ of the graph, denote by $\innone{e}{\cdot\,,\cdot}$ the bilinear form on $\R^E$ defined by
 \begin{equation} \label{eq:EdgePairing}
 \innone{e}{x,y} := x_e y_e
 \end{equation}
 for any pair of elements $x = (x_f)_{f\in E}, y=(y_f)_{f\in E} \in \R^E$. We denote by $q_e$ the corresponding quadratic form on $\R^E$. Notice that the edge length function $l \colon E \to (0, + \infty)$ defines an inner product on $\R^E$ given by 
 \begin{equation} \label{eq:InnerProduct}
 	\innone{l}{x,y} := \sum_{e \in E} l(e) ( x\,, y)_e = \sum_{e \in E} l(e) x_e y_e, \qquad x,y \in \R^E.
 \end{equation}
The corresponding quadratic form on $\R^E$ is denoted by $q_l$. After fixing a basis $\gamma_1, \dots, \gamma_h$ of $H$, the quadratic form $q_e$ restricted to $H_1(G, \R)$ is identified with an $h \times h$ symmetric matrix $\rmM_e$ of rank at most one. Thinking of elements of $H_1(G, \R)$ as column vectors, we have
\begin{equation}\label{eq:1}
 q_e(x) = \prescript{\transpose}{}{x}\rmM_e x.
\end{equation}
The restriction of $q_l$ to $H_1(G, \R)$ can be identified with the $h\times h$ symmetric matrix 
 \begin{equation} \label{eq:M_l}
 	\rmM(\mgr) := \rmM_l  = \rmM(G, l) = \sum_{e\in E}l(e)  \rmM_e.
 \end{equation}
that we call the \emph{period matrix} of the metric graph $\mgr$ (more precisely, with respect to the model $(G, l)$).

\smallskip

The {\em Zhang measure} on the metric graph $\mgr$ is the measure
\begin{equation} \label{eq:FosterCoefficient}
\mu_\Zh :=   \sum_{e \in E} \frac{\mu (e)}{\ell(e)} d\lambda_e,
\end{equation}
where $\lambda_e$ denotes the uniform Lebesgue measure on the edge $e$ and the {\em Foster coefficients} $\mu(e)$, $e \in E$, are given by  (see \cite[Theorem 5.1]{AN})
 \[\mu(e) =  \sum_{i,j=1}^h \rmM^{-1}_l (i,j) \,  l(e) \gamma_i(e) \gamma_j(e), \]
with respect to the fixed basis $\gamma_1,\dots,\gamma_h$ of $H_1(G, \R)$. It is easy to see that $\mu_{\Zh}$ has total mass $\graphgenus$. We define the {\em canonical  measure} $\mu^\can$ of $\mgr$ by normalizing $\mu_\Zh$, that is, 
\[\mu^\can := \frac 1{\graphgenus} \mu_\Zh = \frac 1{\graphgenus}\sum_{e \in E} \frac{\mu (e)}{\ell(e)} d\lambda_e.\] 
In particular, our definition of the canonical measure differs by the normalization factor $h$ from some other works (see, e.g., \cite{AN, BF11}). Moreover, we highlight that {\em we allow $G$ to be  disconnected in the above definition}.
 
\medskip

The above notion of canonical measure extends naturally to augmented metric graphs. Namely, assume that we have further given a genus function $\genusfunction \colon V \to \N \cup\{0\}$ for $G$.  Denote by $g$ the genus of the augmented graph $(G, \genusfunction)$, that is
\[
	g = h + \sum_{v \in V} \genusfunction(v).
\]

 The canonical measure associated to the augmented metric graph $\mgr$, metric realization of the triple $(G, \genusfunction, l)$, is the measure $\mu^{\can}$ on $\mgr$ given by
\[
	\mu^\can = \frac{1}{g} \left( \mu_\Zh + \sum_{v \in V}  \genusfunction(v) \delta_v \right)
\]
where $\mu_\Zh$ refers to the Zhang measure of the (unaugmented) metric graph $\mgr$, as defined in~\eqref{eq:FosterCoefficient}.  By definition, $\mu^\can$ has again total mass equal to one.

\subsection{The Laplacian on a metric graph} \label{ss:LaplacianMetricGraph}

For the purpose of the present paper, we define the Laplacian on functions $f \colon \mgr \to \C$ satisfying the following properties:
\begin{itemize}
\item [(i)] $f$ is continuous on the metric graph $\mgr$,
\item [(ii)] $f$ is piecewise $\mathcal{C}^2$ (i.e., piecewise $\mathcal{C}^2$ on each edge $e \in E$), 
\item [(iii)] and $f''$ belongs to $L^1(\mgr)$. 
\end{itemize}
The space of all such functions is sometimes called the {\em Zhang space} \cite{BR10, Cinkir, Zhang} and we denote it by $D(\Deltaind{\mgr})$ (or simply $D(\Delta)$).

\smallskip

For any function $f \in D(\Deltaind{\mgr})$, we define its {\em Laplacian} as the following measure on $\mgr$
\[
	\Deltaind{\mgr} f = - f'' d \lambda - \sum_{x \in \mgr} \Big( \sum_{\nu \in T_x} \slp_\nu f(x) \Big) \delta_x.
\]
In the above definition, $d \lambda$ denotes the natural Lebesgue measure on the metric graph $\mgr$. Moreover, for every point $x \in \mgr$, $T_x$ is the set of tangential directions at $x$ and $\slp_\nu f(x)$ is the slope of $f$ at $x$ in the unital direction of $\nu$. 
 Note that the tangential directions at a vertex $v \in V$ correspond to the incident (half)-edges. Therefore, for an incident edge $e$ to the vertex $v$, we denote by $\slp_e f(v)$ the corresponding tangential direction. 

\smallskip

At a point $x \in \mgr$ in the interior of an edge $e \in E$, there are precisely two tangential directions. By the regularity properties of a function $f \in D(\Deltaind{\mgr})$, the above derivatives exist for all $x \in \mgr$. Moreover, the sum $\sum_{\nu \in T_x} \partial_\nu f(x)$ vanishes identically on all open sets $U$ where  $f$ is $\mathcal{C}^2$. In particular, $\Deltaind{\mgr}f$ is of the simple form 
\begin{equation} \label{eq:SpecialMus}
	\mu = g d \lambda + \sum_{x \in A} a_x \delta_x,
\end{equation}
where $g \in L^1(\mgr)$ is a piecewise continuous function, the set $A \subseteq \mgr$ is finite, and $a_x \in \C$, $x \in A$, are coefficients. Altogether, we obtain a map
\begin{equation} \label{eq:LaplacianMap}
\begin{array}{cccc}
\Deltaind{\mgr} \colon & D(\Deltaind{\mgr})  &\longrightarrow & \mathcal{M}(\mgr)
\end{array}
\end{equation}
from the Zhang space $D(\Deltaind{\mgr})$ to the space $\mathcal{M}(\mgr)$ of complex-valued Borel measures $\mu$ on $\mgr$ with finite total variation, $|\mu|(\mgr) < \infty$. Alternatively, $\Deltaind{\mgr} f$ can be described as the unique measure in $\mathcal{M}(\mgr)$ such that
\begin{equation} \label{eq:GraphLaplacianDistribution}
\int_\mgr \varphi(y) \, \Delta_\mgr f = - \int_\mgr f(y) \varphi''(y) \, d\lambda(y)
\end{equation}
for every test function $\varphi\colon \mgr \to \C$ in the space 
\[
\mathcal{C}^2(\mgr) = \Big \{\varphi \in \mathcal{C}(\mgr) \st \, \text{$\varphi$ is edgewise in $\mathcal{C}^2([0, l(e)])$ and $\sum_{e \sim v}  \slp_e \varphi (v) = 0$ for all $v \in V$} \Big \}.
\]
In particular, \eqref{eq:GraphLaplacianDistribution} implies that the measure $\Deltaind{\mgr} f$ has zero total mass for every $f \in D(\Delta)$. Hence, the image of the map $\Deltaind{\mgr}$ in \eqref{eq:LaplacianMap} is contained in
\begin{equation} \label{eq:NiceMeasures}
	\widetilde{\mathcal{M}}^0(\mgr) := \bigl\{\mu \in \mathcal{M}(\mgr) \st \, \text{$\mu$ is of the form \eqref{eq:SpecialMus} and $\mu(\mgr) = 0$} \bigr\}.
\end{equation}

\medskip 

It is well-known (see, e.g., \cite{BR10}) that for any measure $\mu \in \widetilde{\mathcal{M}}^0(\mgr)$, the {\em Poisson equation}
\begin{equation} \label{eq:LaplacianEquationGraph}
\Delta_\mgr f = \mu
\end{equation}
has a  solution $f \in D(\Deltaind{\mgr})$ which is unique up to additive constants.  In particular, imposing the value of $f$ to be zero at a given point $x \in \mgr$ uniquely determines $f$.

\subsection{The metric graph $\jvide$-function} \label{ss:LaplacianMetricGraph}

It follows from the above discussion that for every three points $p, q, x \in \mgr$, the following system of equations
\begin{equation} \label{eq:JFunctionGraph}
\begin{cases}
\Deltaind{\mgr} \, \jfunc{p \tiret q, x}(\cdot) = \delta_p - \delta_q \\
 \jfunc{p \tiret q, x} (x) = 0
 \end{cases}
\end{equation}
has a unique solution $y \mapsto \jfunc{p \tiret q, x} (y)$ which belongs to $D(\Deltaind{\mgr})$. By definition, $\jfunc{p \tiret q, x}(\cdot)$ is a piecewise affine function on $\mgr$ for all $p, q, x \in \mgr$. We call the mapping
\[
(x,y, p, q) \mapsto  \jfunc{p \tiret q, x} (y), \qquad x,y, p,q \in \mgr,
\]
the {\em $\jvide$-function} of the metric graph $\mgr$.

\begin{remark}
We stress that the above notion differs slightly from the standard $j$-function (or {\em potential kernel}) used in the literature (see, e.g., \cite{BR10, CR93, Cinkir}). Indeed, usually the $j$-function is defined as the mapping
\[
(s, t, \zeta) \mapsto j_{\zeta}(s,t), \qquad s,t,\zeta \in \mgr,
\] defined by the equations
\begin{equation*}
\Deltaind{\mgr} \, j_{\zeta}(\cdot,t) = \delta_t - \delta_\zeta, \qquad j_{\zeta}(\zeta,t) = 0.
\end{equation*}
Clearly, the connection between these two notions is given by the equation
\begin{equation} \label{}
\jfunc{p \tiret q, x} (y) = j_q(y, p) - j_q(x, p), \qquad x,y,p,q \in \mgr.
\end{equation}
However, our notion is more closely related to the $j$-function on Riemann surfaces and this is advantageous for the considerations in the final part of this paper.
\end{remark}
For later considerations, we also mention the following useful equality. If $\mu \in \widetilde{\mathcal{M}}^0(\mgr)$ is a measure of zero mass and $x \in \mgr$, then the solution to the Poisson equation
\begin{equation*} \label{eq:Poisson_metric_graphs}
\begin{cases}
\Deltaind{\mgr} f = \mu \\
f (x) = 0
 \end{cases}
 \end{equation*}
 can be written in terms of the $\jvide$-function as
\begin{equation} \label{eq:SolutionFormulaJFunction}
f(y) = \int_\mgr \jfunc{s \tiret q, x} (y) \, d\mu(s) = - \int_\mgr \jfunc{p \tiret s, x} (y) \, d\mu(s), \qquad y\in \mgr,
\end{equation}
where $p,q \in \mgr$ are arbitrary (see \cite[Chapter 3]{BR10}).
 
\begin{remark} \label{rem:DistributionalLaplacian}
In the definition of the Laplacian $\Deltaind{\mgr}$, we have restricted to the function space $D(\Deltaind{\mgr})$ only for the sake of a clear exposition. For instance, in \cite{BR07, BR10} the Laplacian is an $\mathcal{M}(\mgr)$-valued operator on the larger space $\operatorname{BDV}(\mgr)$ of functions with bounded differential variation.
 
 On the other hand, one might also pursue a more distributional approach, where the Laplacian takes values in the dual $\mathcal{C}^2(\mgr)^\dual$ (possibly after equipping $\mathcal{C}^2(\mgr)$ with a suitable topology). Following \eqref{eq:GraphLaplacianDistribution}, the Laplacian of every function $f \in L^1(\mgr)$ could be introduced as the functional
 \begin{equation} \label{eq:LaplacianDistr}
 (\Deltaind{\mgr} f ) (\varphi) =  - \int_\mgr f(y) \varphi''(y) \, d\lambda(y), \qquad \varphi \in \mathcal{C}^2(\mgr).
 \end{equation}
This actually allows to treat the Poisson equation \eqref{eq:LaplacianEquationGraph} for general measures $\mu \in \mathcal{M}(\mgr)$ of mass zero (see \cite{BR07, BR10} for details). In particular, the function $f$ defined by \eqref{eq:SolutionFormulaJFunction} belongs to $\operatorname{BDV}(\mgr)$ and is the unique solution of \eqref{eq:LaplacianEquationGraph} in the sense of \cite{BR07, BR10}.  Moreover, it is the unique $L^1$-function such that $\Deltaind{\mgr}f = \mu$ in the dual space $\mathcal{C}^2(\mgr)^\dual$, that is verifying Equation~\eqref{eq:LaplacianDistr}. Of course, uniqueness is understood as up to additive constants here.
 
 \end{remark}

\subsection{The canonical Green function on (augmented) metric graphs}
The {\em canonical Green function} on an (augmented) metric graph $\mgr$ is the unique function $\gri{\mgr} \colon \mgr \times \mgr \to \R$ such that for each point $p \in \mgr$, the function $\gri{\mgr}(p, \cdot)$ belongs to $D(\Delta_\mgr)$ and satisfies
\begin{equation} \label{eq:DefGFGraph}
\Delta \gri{\mgr}(p, \cdot)=  \delta_p - \mu^{\can}, \qquad \int_\mgr \gri{\mgr}(p, y) \, d\mu^{\can}(y) = 0.
\end{equation}
If there is no risk of confusion, we will sometimes abbreviate and write $\grg(p,y)$ instead of $\gri{\mgr}(p,q)$. Notice that $\gri{\mgr}(p, \cdot)$ is a piecewise quadratic function on $\mgr$ for each $p \in \mgr$, see e.g., \cite{Cin14}. 

From the above discussion, we also deduce the following representation
\begin{equation} \label{eq:GreenByJFunction}
\gri{\mgr}(p,y) = \int_\mgr \jfunc{p\tiret q,x}(y) \, d\mu^{\can}(q)  - \int_\mgr \int_\mgr \jfunc{p\tiret q,x}(y)  \, d\mu^{\can}(q) \, d\mu^{\can}(p), \qquad p,y \in \mgr,
\end{equation}
which holds for all points $x \in \mgr$.

\smallskip

\section{Function theory on tropical curves} \label{sec:tropical_function_theory}
The aim of this and the forthcoming Section~\ref{sec:hybrid_laplacian} is to introduce a framework for doing function theory in the layered setting and in the presence of sedentarities on tropical and hybrid curves. This will be used to define Laplace operators on tropical and hybrid curves (the analog of the operator $\Delta = \frac1{\pi i}\partial \bar{\partial}$ on Riemann surfaces) and to explain how they arise as the (weak measure-theoretic) limit of the Laplacians on degenerating metric graphs and Riemann surfaces, respectively.

\subsection{$\R$-divisors and one-forms on graphs} We start by recalling the definition of {\em $\R$-divisors} and {\em one-forms} on graphs and collect some elementary properties.

\smallskip

Let $G = (V, E)$ be a finite graph. An \emph{$\R$-divisor} $D$ on $G$ is an element of the real vector space generated by the vertices of the graph. The generator associated to the vertex $v$ of the graph is denoted by $(v)$. An $\R$-divisor $D$ can be thus written as a formal sum
\[
	D = \sum_{v \in V} D(v) (v) 
\]
where $D(v) \in \R$ is called the \emph{coefficient} of $v$ in $D$. In this paper we only consider divisors with real coefficients, so we  drop the mention of $\R$ and simply talk about divisors. Divisors with integral coefficients arise in connection with algebraic geometry of tropical and hybrid curves, and will be considered in our forthcoming work~\cite{AN-AG-hybrid}.

The {\em degree} of a divisor $D$ is given by
\[
	\deg(D) := \sum_{v \in V} D(v).
\]
We denote the space of divisors on $G$ and its subspace of degree zero divisors by $\Div(G)$ and $\Div^0(G)$, respectively.

\medskip

Let $\EE$ be the set of all the possible orientations of edges in $E$ so that for each edge $\{u,v\}$, there are two possible oriented edges $uv$ and $vu$. By an abuse of the notation, we use $e$ when referring both  to elements of $E$ and $\EE$. In the latter case, when $e =uv$, we denote by $\bar e$ the same edge with the reverse orientation, so $\bar e = vu$. For an oriented edge  $e=uv$, we call $u$ the \emph{tail} and $v$ the \emph{head} of $e$ and denote them by $\tail_e$ and $\head_e$, respectively.   By an \emph{orientation} of the edges of $G$, we mean a map $\oo\colon E \to \EE$ which is a section of the forgetful map $\EE \to E$, that is, for an edge $e=\{u,v\}$, $\oo(e)$ is either $e$ with the orientation $uv$ or $e$ with the orientation $vu$.

\smallskip

The {\em space of one-forms} on $G$ denoted by $C^1(G, \R)$ consists of all the maps $\omega\colon \E \to \R$ which verify the equation 
\[\omega(e) = - \omega(\bar e), \qquad e\in \EE.\]
It is a vector space which can be identified with $\R^E$ upon the choice of a  fixed orientation $\oo\colon E \to \EE$.  
An element $\omega \in C^1(G, \R)$ is called a {\em one-form} on $G$. We consider the usual {\em boundary operator} $\partial \colon C^1(G, \R) \to \Div(G)$ given by
\[
\partial \omega(v) := \sum_{v \in V} \Big( \sum_{\substack{e \in \EE \\
\head_e =v}} \omega(e)\Big ) (v).
\]
A one-form $\omega$ on $G$ is called {\em harmonic} if it has vanishing boundary
\[
	\partial \omega = 0,
\]
and the {\em space of harmonic one-forms} on $G$ is denoted by
\[
\Omega^1(G) := \Bigl\{ \omega \in C^1(G, \R) \,\st \, \omega \text{ is a harmonic one-form} \Bigr\}.
\]

\subsubsection{Hodge decomposition} Suppose further that the graph $G=(V,E)$ is equipped with an edge length function $l \colon E \to (0, + \infty)$.  Having this additional structure,  the {\em differential} of a function $f \colon V \to \R$ is the one-form $df \in C^1(G, \R)$ given by
\[
	(df) (e) = \frac{f(v) - f(u)}{l(e)}, \qquad e= uv \in \EE.
\]

\smallskip

A one-form $\omega \in C^1(G, \R)$ arising in this way is called {\em exact on} $(G, l)$ and the {\em space of exact forms} is denoted by 
\[\exact(G, l) := \Bigl\{ df \, \st \,f \colon V \to \R \Bigr\}.\]

\smallskip

Recall from \eqref{eq:InnerProduct} that the edge lengths on $G$ define a natural inner product on $C^1(G, \R)$ given by
\[
\innone{l}{\alpha, \beta}= \frac 12 \sum_{e \in \EE} l(e) \alpha(e) \beta(e).
\]

In terms of this pairing, a one-form $\alpha$ on $G$ is exact if and only if
\[
	\innone{l}{\alpha, \omega}= 0
\]
for all harmonic one-forms $\omega \in \Omega^1(G)$. In particular, the space of one-forms $C^1(G, \R)$ has the classical orthogonal Hodge decomposition 
\begin{equation} \label{eq:OrthDecomposition}
C^1(G, \R) = \Omega^1(G) \oplus \exact(G, l).
\end{equation}
with respect to the inner product defined by \eqref{eq:InnerProduct}.

\smallskip

Notice that there is a natural identification between the space of harmonic one-forms and the elements of $H_1(G, \R)$. This identification is consistent with the above decomposition via \emph{integration along edges}, which for an edge $e \in \EE$, sends $\alpha \in C^1(G, \R)$  to $\int_e \alpha = \alpha(e)l(e)$.
In particular, each cycle $\gamma \in H_1(G, \R)$ defines via integration an element of the dual $C^1(G, \R)^\dual$
\[\alpha \to \int_\gamma \alpha = \frac 12 \sum_{e\in \EE} \int_{e} \alpha, \qquad \alpha \in C^1(G, \R),\]
The the dual space can be identified with $C^1(G, \R)$ via the pairing $\innone{l}{\,,}$. Since the integration of exact one-forms along cycles vanishes, this element is necessarily in $\Omega^1(G, \R)$ and is thus a harmonic one-form on $G$. By a slight abuse of the notation, we will denote it with the same letter $\gamma$. Equivalently, in what follows, an expression of the form $\innone{l}{\gamma, \alpha}$ for $\gamma \in H_1(G, \R)$ and $\alpha \in C^1(G, \R)$ is viewed as referring to the integration pairing $\int_{\gamma}\alpha$ in the metric realization $\mgr$ of $(G, l)$.

\subsubsection{Orthogonal projection to the space of harmonic one-forms} \label{ss:ProjectionHarmonicOneForms}

In the upcoming sections, the orthogonal projection $\projhar$ from $C^1(G, \R)$ onto its subspace $\Omega^1(G)$ of harmonic forms will play a crucial role. For later use, we recall here an explicit expression for this orthogonal projection. 

\smallskip

Fix a basis $\gamma_1, \dots, \gamma_h$ of $H_1(G, \Z)$ (as always, $h$ is the genus of $G$) and consider, as~\eqref{eq:M_l}, the associated $h\times h$ period matrix $\rmM_l$ defined as 
\[\rmM_l = \sum_{e\in E} l(e) \rmM_e.\]
 Fix an orientation $\oo\colon E \to \EE$ and let $E_\oo \subset \EE$ be the image of $\oo$. We introduce another matrix $P_l$ of dimension $h \times |E|$ by setting
\begin{equation} \label{eq:MatrixP}
	P_l(\gamma_i, e) := \innone{l}{\gamma_i, e} = l(e) \gamma_i(e), \qquad i=1, \dots, h, \, e\in E_\oo,
\end{equation}
for any cycle  in our fixed basis of $H_1(G, \R)$, and any oriented edge in $E_\oo$. By what preceded, the matrix $P_l$ corresponds simply to the integration pairing between cycles and one-forms on $(G,l)$.

\smallskip

 We have the following explicit expression of the orthogonal projection.

\begin{prop} \label{prop:OrthoProjection}
Notations as above, the orthogonal projection $\projhar\colon C^1(G, \R) \to \Omega^1(G)$ is given by
\begin{equation} \label{eq:OrthoProjection}
\projhar (\alpha) = \sum_{i=1}^h (\rmM_l^{-1} P_l \alpha)_i \, \gamma_i, \qquad \alpha \in C^1(G, \R).
\end{equation}
\end{prop}
\begin{proof}  
 It suffices to show that the one-form
\[
\widehat{\alpha} = \alpha - \sum_{i=1}^h (\rmM_\ell^{-1} P_l \alpha)_i  \, \gamma_i
\]
satisfies $\innone{l}{\gamma,  \widehat{\alpha}} = 0$ for all cycles $\gamma \in H_1(G, \R)$. However,
\[
\innone{l}{\gamma_i,  \widehat{\alpha}} = \innone{l}{\gamma_i, \alpha} - \sum_{j=1}^h (\rmM_l^{-1} P_l \alpha)_i M_l (i,j) = \innone{l}{\gamma_i, \alpha} - (P_l \alpha)_i = 0
\]
for any cycle $\gamma_i$ in the fixed basis of $H_1(G, \R)$, and hence, the claim follows.
\end{proof}

\subsection{Spaces of functions on tropical curves}  \label{ss:FunctionTheoryCurve} Let $\curve$ be a tropical curve of rank $r$ with underlying model $(G, \pi, l)$, with the ordered partition $\pi =(\pi_\infty, \pi_\fin)$, $\pi_\infty = (\pi_1, \dots, \pi_r)$, and  the length function $l\colon E \to (0, + \infty)$ satisfying the normalization property $\sum_{e\in \pi_j} l(e)=1$ for all $j \in [r]$, in the infinitary part. Let $E_\infty =E\setminus \pi_\fin$, and denote by $\Gamma^j$ the $j$-th graded minor of $\curve$ obtained as the metric realization of the pair $(\grm{\pi}{j}(G), l_j)$.

\smallskip

A \emph{real} (\emph{resp. complex}) \emph{valued function} $\lf$ on $\curve$ is an $r+1$-tuple $\lf = (f_1, \dots, f_r, f_\fin)$ consisting of functions $f_j \colon \Gamma^j \to \R$ (resp. $f_j \colon \Gamma^j \to \C$) for each $j=1, \dots, r, \fin$. The function $\lf =(f_j)_j$ is called \emph{continuous, piecewise smooth, $L^2$, etc.} if the corresponding functions $f_j$ on $\Gamma^j$ verify the same property, that is, if they are continuous, piecewise smooth, $L^2$, etc., respectively.

\subsection{Pullback of functions on tropical curves to metric graphs} \label{sss:PullbackExplanation} Notations as in the previous section, the curve $\curve$ lies in the stratum $\mgtropcombin{\combind{(G,\pi)}}$ of  $\mgtropcombin{\grind G}$.
 Let $\unicurvetrop{\combind{(G,\subface \pi)}}$ be the universal tropical curve over $\mgtropcombin{\combind{(G, \subface \pi)}}$. (Recall that
 $\mgtropcombin{\combind{(G, \subface \pi)}}$ is the union of all the strata $\mgtropcombin{\combind{G,\pi'}}$ for $\pi' \subface \pi$.)  For a function $\lf = (f_1, \dots, f_r, f_\fin)$ on $\curve$, we introduce a {\em pullback} $\lf^*$ of $\lf$ to the universal tropical curve $\unicurvetrop{\combind{(G,\subface \pi)}}$ as follows.

\smallskip

By definition, functions on a tropical curve are automatically layered into $r+1$ components, and different graded minors of the tropical curve live in different infinities, and components are separated via these different infinities. Nevertheless, in order to define the pullback, we adapt the following  perspective which will allow to simplify the notations.

Considering the edge length function $l_j\colon \pi_j \to (0, + \infty)$, $j\in [r] \cup\{\fin\}$, we denote by $\Gamma$ the element of $\mggraphcombin{\grind{G}}$ with edge length function given by the $l_j$. (More generally, for any $\pi'$ with $\pi'\subface \pi$, we get an element in the stratum $\mgtropcombin{\combind{(G,\pi')}}$ with the edge length functions given by $l_j$s.) We first \emph{pullback} the functions $f_j$ to $\Gamma$ and then use these functions to propagate the definition of the pullback over the full stratum $\mgtropcombin{\grind{G}}$ (more generally, over the stratum  $\mgtropcombin{\combind{(G,\pi')}}$). 

To do this, observe that functions on a graded minor $\Gamma^j$ can be identified with a special class of functions on $\Gamma$, namely, those functions $f\colon \Gamma \to \C$ satisfying the following properties:
\begin{itemize}
\item $f$ is constant on all edges of larger index layers $\pi_i$, $i>j$ (with our convention that $\fin >j$ for $j\in [r]$) and
\item $f$ is affine on all edges of lower index layers $\pi_i$, $i<j$.
\end{itemize}
Indeed, such a function $f$ clearly defines by restriction a function on $\Gamma^j$. On the other hand, every function on $\Gamma^j$ can be first extended to (the metric realization of) the subgraph $G_\pi^j$,  by composing with the projection map $\kappa \colon G_\pi^j \to \grm{\pi}{j}(G)$, and then to $\Gamma$ by linear interpolation on remaining edges. 

In conclusion, each piece $f_j$ of a function $\lf = (f_j)_j$ on a given tropical curve $\curve$ can be identified with a function on $\Gamma$ constant on higher index layers and affine on lower index layers. By a slight abuse of notation, we will denote this function by $f_j$ as well and stress the difference explicitly whenever there is a risk of confusion.

\smallskip

Let now $t$ be a point in $\mggraphcombin{\grind{G}}$ with $\umgr_t$ the corresponding metric graph with the same underlying graph as the tropical curve $\curve$ and with the edge length function  $\ell_t:E \to (0, + \infty)$. In this situation, we introduce the {\em pullback} of  the tropical function $\lf = (f_1, \dots, f_r, f_\fin)$ to $\umgr_t$ as follows.

\smallskip

For each $j=1,  \dots, r$, we first define a natural pullback $f_j^\ast \colon \mgr_t \to \C$ of $f_j$ to $\mgr_t$. By the discussion which preceded,  $f_j$ can be viewed as a function on the metric graph $\Gamma \in \mggraphcombin{\grind{G}}$. Moreover, there is a natural map from $\umgr_t$ to $\Gamma$  which is given on each edge by rescaling the edge lengths. We define $f_j^\ast$ as the corresponding pullback. More formally, on each edge $e \in E$ in $\umgr_t$, the pullback $f_j^\ast$ is given by
\[
f_j^\ast \rest e (x) =  f_j \Big ( \frac{l(e)}{\ell_t(e)} x \Big),
\]
where $x$ is the parametrization of the (oriented) edge $e=uv$ in $\mgr_t$ (after identification with the interval $[0,\ell_t(e)]$). 

\noindent Altogether, we introduce the pull-back of $\lf = (f_j)_j$ as the sum
\[
 \lf^\ast := f_\fin^\ast+ \sum_{j=1}^r L_j f_j^\ast  ,
\]
where the scaling factors $L_j=L_j(t)$, $j=1,\dots, r$, are given by $L_j(t) = \sum_{e \in \pi_j} \ell_t(e)$ (see \eqref{eq:DefGrLength}).

Note in particular that, along a sequence of metric graphs $\umgr_t$ degenerating to the tropical curve $\curve$, the pullbacks $\lf_t^\ast$ to $\umgr_t$ of a fixed function $\lf=(f_j)_j$ on $\curve$ become unbounded (since the $L_j$'s go to infinity).

\smallskip

The above definition naturally extends to all the strata $\mgtropcombin{\combind{(G,\pi')}}$ for any $\pi' \subface \pi$, and leads to the definition of a pullback $\lf^*$ over $\mgtropcombin{\combind{(G,\subface \pi)}}$. For any point $\thy \in \mgtropcombin{\combind{(G, \subface \pi)}}$, we denote by $\lf_\thy^*$ the pullback of $\lf$ to $\hcurve^\trop_\thy$, the tropical curve given by the point $\thy$.

\subsection{Harmonically arranged functions on tropical curves} We now define the important class of \emph{harmonically arranged functions} on tropical curves. As we will see later, this type of functions is precisely the one which captures the asymptotic behavior of solutions to the Poisson equation $\Delta f = \mu$ on degenerating metric graphs.

\smallskip

Notations as above, let $\curve$ be a tropical curve with underlying graph $G = (V, E)$, edge length function $l: E \to (0, + \infty)$, the ordered partition $\pi=(\pi_1, \dots, \pi_r, \pi_\fin)$,  and the corresponding layered metric graph $\Gamma \in \mggraph{\grind{G}}$ represented by the metric realization of $(G, l)$, as in the previous section. We adapt the convention from the previous section and view a function defined on the graded minor $\Gamma^j$ as being defined on the whole metric graph $\Gamma$. 

\medskip

A function $f_j$ on the $j$-th graded minor $\Gamma^j$ is called \emph{harmonic in lower index layers} or simply \emph{lower harmonic} if it is harmonic up to layer $j-1$. Hereby we mean that, when interpreted as a function on $\Gamma$, the differential $d f_j $ on $\Gamma$  is harmonic on lower index layers, that is, more precisely,
\[
 d f_j \rest{\pi_i} \in \Omega^1(\grm{\pi}{i}(G))
\]
for all $i <j$.
Equivalently, this amounts to requiring that for all $i<j$ and all vertices $w$ of the $i$-th graded minor $\Gamma^i$,
\[
\sum_{u \in \proj_{i}^{-1}(w)} \sum_{\substack{{e = uv} \\ e \in \pi_i}} \frac{f_j(v) - f_j(u)}{l_i(e)} = 0,
\]
where $\proj_i \colon V \to V(\grm{\pi}{i}(G))$ is the contraction map. Note that this makes sense since the restriction of the function $f_j$ on the edges of $\pi_i$, $i<j$, is affine linear, so $df\rest{\pi_i}$ restricts to a one-form on $\gr^i_\pi(G)$.

\smallskip

Finally, a function $\lf =(f_j)_j$ on the tropical curve $\curve$ is called \emph{harmonically arranged} if all the components $f_j$'s are lower harmonic.

\begin{figure}[!t]
\centering
    \scalebox{.4}{\input{example9.tikz}}
\caption{A tropical curve of full sedentarity with $\pi_\infty=(\pi_1, \pi_2, \pi_3)$, and its corresponding graded minors. The edges in each layer have the same length.  The function $\lf=(f_1,f_2, f_3)$ on $\curve$ has components $f_1$, $f_2$, and $f_3$ defined on the first, second, and third layers, respectively.} 
\label{fig:harmonically_arranged}
\end{figure}
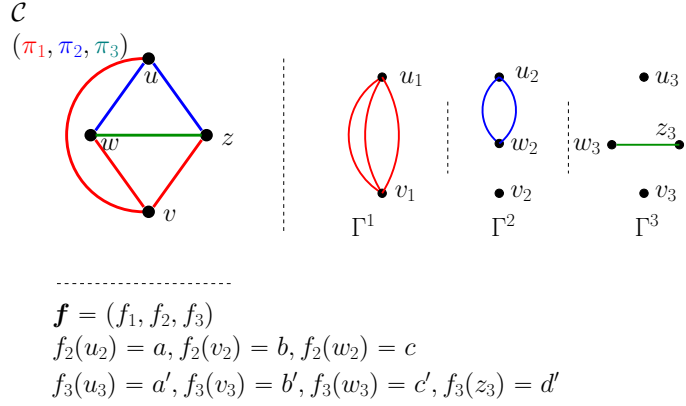

\begin{example}[Harmonically arranged functions on the tropical curve~\ref{fig:harmonically_arranged}] \label{example:harmonicall_arranged_function} An example of a tropical curve of rank three with full sedentarity and a function $\lf=(f_1, f_2, f_3)$ is depicted in Figure~\ref{fig:harmonically_arranged}. The function $\lf$ is harmonically arranged if $f_2$ and $f_3$ are lower harmonic. Assume $f_2(u_2)=a$, $f_2(v_2)=b$ and $f_2(w_2)=c$ where $u_2, v_2$, and $w_2$ are the vertices of $\Gamma^2$. The function $f_2$ viewed on $\Gamma$ takes values $f_2(u)=a, f_2(v)=b$, and $f_2(w)=f_2(z)=c$. The differential $df_2$ takes values $df_2(uv) = 3(b-a)$ and $df_2(wv) = df_2(zv)=3(b-c)$.  It follows that $f_2$ is lower harmonic if $b-a+2(b-c)= 3b - a -2c =0$, that is, if $b = \frac{a+2c}3$. Consider now $f_3$ and assume $f_3(u_3) = a', f_3(v_3)= b', f_3(w_3) = c'$, and $f_3(z_3) =d'$, where $u_3, v_3, w_3, z_3$ are the vertices of $\Gamma^3$ (they coincide with those of $\Gamma$). Then, $df_3$ on the edges of $\pi_1$ takes values $df_3(uv) = 3(b'-a')$ and $df_3(wv) = 3(b'- c')$, $df_3(zv)=3(b'-d')$, and takes values $df_3(uw)=2(c'-a')$, $df_3(uz)=2(d'-a')$, on the edges of $\pi_2$. It follows that $f_3$ is lower harmonic if $3 b' = a' + d'+ c'$ and $2a' = c' + d'$. 
\end{example}

\subsection{Harmonic rearrangements}\label{sec:harmonic_rearrangement} We  describe in this section how an arbitrary function $\lf$ on a tropical curve $\curve$ can be transformed into a harmonically arranged function $\tilde{\lf}$ by what we call a \emph{rearrangement} of $\lf$, as described below.

\smallskip

Notations as in the previous section, let $\pi = (\pi_\infty, \pi_\fin)$, $\pi_\infty=(\pi_1, \dots, \pi_r)$, be the ordered partition underlying the tropical curve $\curve$. Let $\lf =(f_1, \dots, f_r, f_\fin)$ be a function on $\curve$. \smallskip

First we choose a constant $\const{H}{j}$ for each connected component $H$ of each graded minor $\Gamma^j$, and then we set
\[
\tilde{f}_j\rest{H} := f_j\rest H + \const{H}{j}
\]
for each such connected component $H$. A function $\tilde{\lf} = (\tilde{f}_j)_j$ obtained in this way will be called a \emph{rearrangement} of $\lf$; it is called a \emph{harmonic rearrangement} if the resulting function $\tilde{\lf}$ is harmonically arranged.

\begin{remark} The term \emph{rearrangement} shall be understood as the \emph{relative rearrangement} of the restrictions to connected components of graded minors. Simplifying, we drop the term relative. Note that, modifying a harmonic rearrangement $\lf$ by the same constant on each component of a layer results in another harmonic rearrangement. That is, modifying $\lf$ with a \emph{global constant} $\constwi \in \R^{r+1}$ (or $\C^{r+1}$, if the functions are complex valued) preserves the harmonically arranged property.
\end{remark}

The following  important result shows the relevance of this concept to our setting. 

\begin{thm}[Harmonic rearrangement of tropical functions] \label{thm:HarmonicExtension}
Every function $\lf = (f_j)_j$ on $\curve$ has a harmonic rearrangement. Moreover, this is unique up to a global constant $\constwi \in \R^{r+1}$, $\constwi \in \C^{r+1}$ if functions are complex valued. That is, for two such harmonic rearrangements $\tilde{\lf} = (\tilde f_j)_j$ and $\tilde{\lf}'=(\tilde f'_j)_j$, the difference $\tilde{f}_j - \tilde{f}_j'$ is constant on $\Gamma^j$ for all $j$.
\end{thm}
The proof of this theorem is given in the next section.

\begin{example}[Example~\ref{example:harmonicall_arranged_function} continued] Consider the example given in Figure~\ref{fig:harmonically_arranged}. The graded minor $\Gamma^2$ has two components. The metric graph $\Gamma^2$ has two connected components and a rearrangement of $f_2$ consists in changing the values of $f_2$ by adding  constants $\constwi^2_1$ and $\constwi^2_2$ to each of the two components. This results in $\tilde f_2$ taking values $\tilde f_2(u_2) = a + \constwi^2_{1}, \tilde f_2(w_2) = c + \constwi^2_1$, and $\tilde f_2(v_2) = b+\constwi^2_{2}$. By Example~\ref{example:harmonicall_arranged_function}, $\tilde f_2$ becomes lower harmonic if the equation 
$3 b+3 \constwi^2_2  = a+ 2c + 3\constwi^2_1$ is verified. That is, if   $\constwi^2_2 - \constwi^2_1 = \frac{a+2c}3 - b$. 

\smallskip

Similarly, $\Gamma^3$ has three components, and a rearrangement of $f_3$ consists in adding a constant $\constwi^3_1$, $\constwi^3_2$, and $\constwi^3_3$ to each of these three components. This results in  $\tilde f_3$ taking values 
\[\tilde f_3(u_3) = a' + \constwi^3_{1},  \,\,  \tilde f_3(w_3) = c' + \constwi^3_2, \,\, \tilde f_3(z_3) = d' + \constwi^3_2, \quad \textrm{and}\qquad \tilde f_3(v_3) = b'+\constwi^3_{3}.\]
In particular, combined with Example~\ref{example:harmonicall_arranged_function}, $\tilde f_3$ becomes lower harmonic if 
\[\constwi^3_1 - \constwi^3_2 = \frac{c'+d' -2a'}2 \qquad \textrm{and} \qquad 3 \constwi^3_3 - 2 \constwi^3_2 - \constwi^3_1 = a'+c'+d' - 3b'. \]
\end{example}

\subsection{Extension lemma} We will deduce Theorem~\ref{thm:HarmonicExtension} from the following lemma. 

\begin{lem}[Extension lemma] \label{lem:extension_lemma}
Let $\curve$ be a tropical curve with underlying graph $G=(V, E)$,  the ordered partition $\pi = (\pi_\infty, \pi_\fin)$ of $E$, $\pi_\infty=(\pi_1, \dots, \pi_r)$, and the length function $l$. Then, every exact one-form $\alpha_{\pi_j}\in \exact(\grm{\pi}j(G), l_j)$ on the graded minor $(\grm{\pi}j(G), l_j)$ has a unique extension $\alpha_G \in C^1(G, \R)$  such that the following properties hold:

\begin{itemize}
\item [(i)] $\alpha_G$ is exact on $(G, l)$,
\item [(ii)] the restriction $\omega_i := \alpha_G\rest{\pi_i}$ is harmonic on $\grm{\pi}{i}(G)$ for all $i< j$, and
\item [(iii)] $\alpha_G\rest{\pi_i} = 0$ for all $i > j$. 
\end{itemize}
\end{lem}

In the rest of this section, we give a proof of this lemma. For further applications later, we will need to provide an explicit way of constructing $\alpha_G$. This will be done in the next section.

\smallskip

We begin by proving the uniqueness statement.

\begin{proof}[Proof of the uniqueness in Lemma~\ref{lem:extension_lemma}]
 First of all, observe that every cycle $\gamma$ in a graded minor $\grm{\pi}{i}(G)$ is a contraction of some cycle $\tilde \gamma \in H_1(G, \R)$ on the full graph $G$ using only edges of layers $\pi_k$ with $k \ge i$, see Section~\ref{sec:admissible_basis_layered_graphs} and \cite{AN}. More formally, for every $\gamma \in H_1( \grm{\pi}{i}(G) , \R)$, we can fix some cycle $\tilde \gamma \in H_1(G, \R)$ such that 
 
 \begin{itemize}
 \item $\tilde \gamma(e) = \gamma(e)$ for all oriented edges $e$ whose underlying edge belongs to $\pi_i$, and
 \item $\tilde \gamma(e) = 0$ for all oriented edges in $\pi_{1} \cup \dots \cup \pi_{i-1}$. 
 \end{itemize}

Let now $\alpha_G \in C^1(G, \R)$ be an extension of $\alpha_{\pi_j}$ with the properties stated in the lemma, assume that $i<j$, and let $\omega_i$ be the restriction of $\alpha_G$ to $\pi_i$. Since $\alpha_G$ is exact on $\Gamma$, we get 
\begin{equation} \label{eq:Exactness}
\innone{l}{\tilde\gamma , \alpha_G} = 0.
\end{equation}
Separating $\alpha_G$ into its restrictions to different layers, and using that $\gamma$ is defined on $\Gamma^i$, this implies the equation
\begin{equation} \label{eq:CycleEquation}
\innoneind{l}{i}{\gamma, \omega_i} = - \innoneind{l}{j}{\tilde \gamma\rest{\pi_j}, \alpha_{\pi_j}} - \sum_{i<k<j} \innoneind{l}{k}{\tilde \gamma\rest{\pi_k}, \omega_k}.
\end{equation}
In particular, for $i = j-1$, we obtain the equations
\[
\innoneind{l}{j-1}{\gamma, \omega_{j-1}} = - \innoneind{l}{j}{\tilde\gamma\rest{\pi_j}, \alpha_{\pi_j}}.
\]
By assumption, the one-form $\omega_{j-1}$ is harmonic on $\grm{\pi}{j-1}(G)$. We conclude from the above equations for all cycles $\gamma$ in $H_1(\grm{\pi}{i}(G,\R))$ that $\omega_{j-1}$ is uniquely determined by $ \alpha_{\pi_j}$. Using now \eqref{eq:CycleEquation} inductively, it follows that $\omega_i$'s, $i < j$, are all determined uniquely. Since the restriction of $\alpha_G$ to higher indexed layers is null, we conclude that $\alpha_G$ is unique.
\end{proof}

We now prove the existence of a one-form $\alpha_G$ on $G$ with the above properties.

\begin{proof}[Proof of the existence in Lemma~\ref{lem:extension_lemma}] We use the fact for any two lifts $\tilde\gamma$ and $\tilde \gamma'$ of a cycle $\gamma$ in $\grm{\pi}{i}(G)$, the restriction of the difference $\tilde \gamma - \tilde \gamma'$ to $\pi_j$ is a cycle in $\grm{\pi}{j}(G)$. It follows the pairing $\innoneind{l}{j}{\tilde\gamma\rest{\pi_j}, \alpha_{\pi_j}}$ is a linear form on $H_1(\grm{\pi}{i}(G),\R)$.  Notice now that \eqref{eq:CycleEquation} inductively defines elements $\omega_i \in \Omega^1(\grm{\pi}{i}(G))$ for $i<j$.  Introducing the one-form $\alpha_G$ on $G$ as
\[
 \alpha_G\rest{\pi_i} := \begin{cases} \omega_i, &\text{if } i < j, \\
 \alpha_{\pi_j}, &\text{if } i = j \\
 0,  &\text{if } i > j \end{cases}, 
\] it easily follows that \eqref{eq:Exactness} holds true for all cycles $\gamma \in H_1( \grm{\pi}{i}(G) , \R)$ of all minors $i=1,\dots, r$. However, ranging over all cycle spaces $ H_1( \grm{\pi}{i}(G) , \R)$ for $i=1,\dots, r, \fin$, the extended cycles $\gamma$ generate the cycle space $H_1(G, \R)$, see Section~\ref{sec:admissible_basis_layered_graphs}. Hence, $\alpha_G$ is exact on $G$ and the proof is complete.
\end{proof}

\subsubsection{Proof of the harmonic rearrangement theorem}
We are now ready to prove Theorem~\ref{thm:HarmonicExtension}.

We first prove the existence. Without loss of generality, we can assume all components $f_j$ of $\lf=(f_1, \dots, f_r, f_\fin)$ are real-valued.  We define an exact one-form $\alpha_{\pi_j}$ on $(\grm{\pi}{j}(G), l_j)$ as the differential
\[
\alpha_{\pi_j}(uv) := d f_j (u) = \frac{f_j(v) - f_j(u)}{l_j(e)}, \qquad e \in \pi_j,
\]
where $uv$ is an orientation of the edge $e \in \pi_j$ with extremities $u$ and $v$ in the minor $\grm{\pi}{j}(G)$. Let $\alpha_G$ be the exact one-form on $G$ given by Lemma~\ref{lem:extension_lemma}. Since $\alpha_G$ is exact, there is a continuous function $g_j$ on $\Gamma$, affine linear on edges of $G$, such that $d g_j = \alpha_G$. Moreover, by the properties of $\alpha_G$, the function $g_j$ is constant on the edges of higher index layers $\pi_i$, $i>j$ (with the convention that $\fin> j$, $j\in [r]$). It follows that $g_j$ gives rise to a well-defined function on the graded minor $\Gamma^j$ that we denote by $\tilde{f}_j$.

 By construction, $\tilde{f}_j$ is harmonic on lower index layers. Moreover, since $\alpha_G\rest{\pi_j} =\alpha_{\pi_j}$, the difference  $f_j -\tilde{f}_j $ is constant on each connected component of the $j$-th minor $\Gamma^j$. This implies that $\tilde{\lf} := (\tilde{f}_1, \dots, \tilde f_r, \tilde f_\fin)$ is a harmonic rearrangement of $\lf$.

\smallskip

It remains to prove the uniqueness. Let $\tilde{f}'_j$ be the $j$-th piece of a second harmonic rearrangement $\tilde \lf'$ of $\lf$. We view $\tilde f'_j$ as a function defined on the vertices of $G$, and consider the differential $\alpha' := d \tilde{f}'_j$ as an exact one-form in $\exact(G, l)$. Since $\tilde f'_j$ is harmonic in lower layers and since $\alpha'$ coincides on $\pi_j$ with $\alpha_{\pi_j}$, as defined above, it follows that $\alpha'$ has all the properties of Lemma~\ref{lem:extension_lemma}. This implies that $\alpha'$ coincides with $\alpha_G$  from Lemma~\ref{lem:extension_lemma}. This implies $\alpha_G - \alpha' =0$, and so $\tilde f_j'-\tilde f_j$ is a constant function, as desired.
 The proof is complete.
\qed

\subsection{Explicit sum-product  description of $\alpha_G$ in the extension lemma} \label{sec:explicit_extension_form} For future use in Sections~\ref{sec:tropical_laplacian} and~\ref{sec:GraphPeriodMatrices}, we now explain how to express the extension $\alpha_G$ of an exact one-form $\alpha_{\pi_j}$ in $\exact(\grm{\pi}{j}(G), l_j)$, obtained in Extension Lemma~\ref{lem:extension_lemma}, in closed form. For this we need to introduce a few additional notation. 

\smallskip

 We adapt the notations of Section~\ref{sec:admissible_basis_layered_graphs} and fix an admissible basis $\gamma_1, \dots, \gamma_h$ of the cycle space $H_1(G, \Z)$. By definition, the contracted cycles $\proj_j(\gamma_k)$, $k \in J_\pi^j$, form a basis of $H_1(\grm{\pi}{j}(G), \Z)$ for every $j = 1, \dots, r, \fin$.
 
 We denote by
\begin{equation}
\rmM_{j} =\rmM_j^\pi= \rmM_{\Gamma^j}:=\rmM \big (\grm{\pi}{j}(G), l_j \big) \, \in \R^{h^{j}_\pi \times h^{j}_\pi} 
\end{equation}
the corresponding graph period matrix \eqref{eq:M_l} for the $j$-th graded minor $\Gamma^j$ of $\curve$.

\smallskip

To a pair of indices $(j, k)$ in $ \{1, \dots, r, \fin\}$, we also associate the $h_\pi^j \times h_\pi^k$ matrix $T_{jk}$ given by
\begin{equation} \label{eq:TMatrix}
T_{jk} (m,n) := \sum_{e \in \pi_{\max\{j,k\}}} l(e) \gamma_m (e) \gamma_n (e),
\end{equation}
where $\gamma_m$ and $\gamma_n$, $m \in J_\pi^j$ and $n \in J_\pi^k$, are cycles in the fixed admissible basis of $H^1(G, \Z)$. 
(Note  that with this notation, $T_{jj} = M_{j}$ for every $j \in [r]\cup \{\fin\}$.)

\subsubsection{Strictly increasing sequences and their corresponding matrices}
Let $j$ and $k$ be two elements of $[r]\cup \{\fin\}$ with $j \le k$. A \emph{strictly increasing sequence from $j$ to $k$}  is a sequence 
\[p:\qquad i_0 = j< i_1 < \dots< i_{s-1} < i_s= k.\]
The integer $s$ is called the \emph{length} of the increasing sequence  and is denoted by $\abs p$.

For every pair $(j,k)$ in $[r]\cup\{\fin\}$ with $j<k$, we denote by $\mathcal{P}_{jk}$ the set of strictly increasing sequences from $j$ to $k$, that is,
\[
\mathcal{P}_{jk}:= \Bigl \{ p = (i_0, \dots, i_{|p|}) \st  \, \, i_0 = j, \, i_{|p|} = k \text{ and }  i_0 < i_1 < \dots < i_{|p|} \Bigr\}.
\]
Note that in the case $j = k$, the set $\mathcal{P}_{jj}$ consists of a single sequence $(j)$, which has length zero. 

\smallskip

Each increasing sequence $p \in \mathcal{P}_{jk}$ of length $s =\abs p$ gives rise to an $h_\pi^j \times h_\pi^k$ matrix denoted by $A_p=A_p(\curve) \in \R^{h_\pi^j \times h_\pi^k}$ and defined as the product 
\begin{align} \label{eq:AsmyptoticsApTropical}
A_p=A_p(\curve)  := (-1)^{s}  \,   \rmM_{j}^{-1} \, \bigl(T_{i_0 \, i_1} \, \rmM_{i_1}^{-1}\bigr) \, \bigl(T_{i_1 \, i_2} \, \rmM_{i_2}^{-1}\bigr) \, \dots \,\bigl(T_{i_{s-1} \, i_{s}} \, \rmM_{i_s}^{-1}\bigr).
\end{align}
For the trivial sequence $p = (j)$, we simply obtain $A_{(j)} =  \rmM_{j}^{-1}$.

\subsubsection{Sum-product formula} The following important result provides the explicit form of $\alpha_G$. 

\begin{prop}[Sum-product formula] \label{prop:harmonic_extension_matrices}
Let $\curve$ be a tropical curve with underlying graph $G$, ordered partition $\pi$, and length function $l$. Assume $\alpha_{\pi_j}$ an exact one-form on some graded minor $(\grm{\pi}{j}(G), l_j)$. Denote by $\alpha_G \in C^1(G, \R)$ the extension given in Lemma~\ref{lem:extension_lemma} and denote by $\omega_i$ the restriction of $\alpha_G$ to $\pi$, $i=1, \dots, j-1$, which is harmonic on $\grm{\pi}{i}(G)$.

Then, for any $i \in \{1, \dots, j-1\}$, $\omega_i$ is  the linear combination of the cycles $\kappa_i (\gamma_m)$ in $H_1(\grm{\pi}{i}(G), \R)$, $m\in J_\pi^i$, given by
\[
	\omega_i = - \sum_{m \in J_\pi^i} \Big (  \sum_{k \ge i} \, \sum_{p \in \mathcal{P}_{ik}} A_p(\curve)   \rmW|_{J_\pi^k}  \Big )(n) \,  \kappa_i (\gamma_m)
\]
where the vector $\rmW \in \R^h$ is defined as
\[
\rmW(a) := \innoneind{\ell}{j}{\gamma_a|_{\pi_j}, \alpha_{\pi_j}}, \qquad a \in \{1, \dots, h\}.
\]
\end{prop}

\begin{proof}
Recall that in the proof of Lemma~\ref{lem:extension_lemma}, the extension $\alpha_G$ is constructed inductively, using the equation \eqref{eq:CycleEquation}. The claimed explicit form follows by rewriting the equation \eqref{eq:CycleEquation} in terms of the above matrices and the vector $\rmW$.
\end{proof}

\subsection{One-forms on tropical curves}
In this subsection, we introduce the space of exact one-forms on tropical curves of arbitrary  sedentarity. Moreover, we show how it appears as the limit of the space of exact forms on degenerating metric graphs in the compactified moduli space $\mgtrop{\grind{g}}$ of tropical curves of a given genus.

\smallskip

Suppose $\curve$ is a tropical curve  with underlying combinatorial graph $G$, the ordered partition $\pi = (\pi_\infty=(\pi_1, \dots, \pi_r), \pi_\fin)$ and edge length function $l \colon E \to (0, + \infty)$. We define the \emph{space of one-forms} on $\curve$ denoted by $C^1(\curve, \R)$ as the product
\[
	C^1(\curve, \R) := C^1(\grm{\pi}{1}(G), \R)  \times \dots \times C^1(\grm{\pi}{r} (G), \R) \times C^1(\grm{\pi}{\fin}(G), \R).
\] 
Notice that, formally, we have an isomorphism
\[
C^1(\curve, \R) \cong C^1(G, \R). 
\]
In particular, every one-form on $G$ can be identified with a one-form on $\curve$ and vice versa. However, the graded structure of $\curve$ is reflected in the definition of the {\em space of exact one-forms}, which we  introduce as
\[
	\exact(\curve) := \exact(\grm{\pi}{1}(G), l_1) \times \dots \times  \exact(\grm{\pi}{r}(G), l_r) \times \exact(\grm{\pi}{\fin}(G), l_\fin).
\]
As for graphs, $\exact(\curve)$ coincides with the space of differentials of edge-wise affine linear functions on the tropical curve $\curve$. More precisely, for an edge-wise linear function $\lf = (f_1, \dots, f_r, f_\fin)$ on $\curve$ given by the restrictions $f_j \colon V(\grm{\pi}{j}(G)) \to \R$,
we define the differential $\ld(\lf)$ as the following one-form on $\curve$
\[
	\ld(\lf) := \bigl( d (f_1), \dots, d ( f_r), d( f_\fin)  \bigr ) 
\]
where each $d (f_j) \in \exact(\grm{\pi}{j}(G), l_j)$ is the exact one form associated to $f_j$. 
Clearly, a one-form $\lalpha \in C^1(\curve, \R)$ is exact if and only if it arises as the differential of a function $\lf = (f_1, \dots, f_r, f_\fin)$ on $\curve$.

\smallskip

The following proposition shows that the exact one-forms on $\curve$ are precisely the limits of exact one-forms on degenerating metric graphs, and that limits of exact one-forms on tropical curves remain exact.  This observation is the main motivation behind our terminology.

\begin{prop} \label{prop:exact_forms_limits}
Suppose $\curve =\hcurve^\trop_{\thy_0}$ for a point $\thy_0 \in \mgtrop{\combind{(G, \pi)}}$ is a tropical curve with underlying combinatorial graph $G=(V, E)$, ordered partition  $\pi=(\pi_\infty=(\pi_1, \dots, \pi_r), \pi_\fin)$ of rank $r$, and edge length function $l=\ell_{\thy_0}\colon E\to (0, + \infty)$.

\begin{itemize}
\item [(i)] Let $\umgr_t$ be a sequence of metric graphs in $\mggraph{\grind{G}}$ degenerating to $\curve$, that is $t$ tends to $\thy$ in the tropical moduli space $\mgtrop{\grind{G}}$.  Assume further that $\alpha_t$ are exact one-forms on $\umgr_t$ and the sequence $\alpha_t$ converges to some $\lalpha \in C^1(\curve, \R) \cong C^1(G, \R)$ in the sense that
\[
	\lalpha (e) := \lim_{t \to \thy} \alpha_t (e)
\]
exists for all edges $e \in E$. Then, the limit one-form $\lalpha$ is an exact one-form on the tropical curve $\curve$.
\smallskip

More generally, let $\hcurve_\thy$ be a sequence of tropical curves in $\mgtrop{\combind{(G, \subface \pi)}}$ degenerating to $\curve$, that is, $\thy$ tends to $\thy_0$ in $\mgtrop{\combind{(G, \subface \pi)}}$. If a sequence $\lalpha_t$ of exact forms on $\hcurve_\thy$ converges to $\lalpha$, then $\lalpha$ is exact. 

\medskip

\item [(ii)] Conversely, suppose $\lalpha$ is an exact one-form on the tropical curve $\curve$. For a metric graph $\mgr=\umgr_t$, $t\in \mggraph{\grind{G}}$, of combinatorial type $G$ with edge length function $\ell=\ell_t$, define the exact one-form 
\[
\alpha_t := \projexact \lalpha,
\]
where $\projexact$ denotes the orthogonal projection of $C^1(G, \R)$ onto its subspace $\exact(G, \ell)$ of exact one-forms, and $\lalpha_t$ is seen in $C^1(G, \R)$. Then, as $\mgr$ degenerates to $\curve$, that is, as $t$ converges to $\thy$, the exact one-forms $\alpha_t$ converge to $\lalpha$.

\smallskip
 More generally, for a tropical curve $\hcurve_\thy$, $\thy \in \mgtrop{\combind{(G, \subface \pi)}}$, with ordered partition $\pi'\subface \pi$ and edge length function $\ell_\thy$, let $\projexact$ be the orthogonal projection of $C^1(G, \R)$ onto its subspace $\exact(\hcurve_\thy)$ of exact one-forms  on $\hcurve_\thy$ and set 
 \[
\lalpha_t := \projexact \lalpha.
\]
 Then, as $\thy$ converges to $\thy_0$ in $\mgtrop{\combind{(G, \subface \pi)}}$, $\lalpha_\thy$ converges to $\lalpha$. 
\end{itemize}
\end{prop}
\begin{proof} To prove the first claim, notice as before that for any cycle $\gamma \in H_1(\grm{\pi}{j}(G), \R)$ in the $j$-th graded minor, we can find an extension of $\gamma$ to a cycle $\tilde \gamma \in H_1(G;\R)$  in $G$ using only edges of higher index layers. Taking into account the definition of the topology on $\mgtrop{\grind{G}}$, and the exactness of $\alpha_t$ which implies $\innoneind{\ell}{t}{\tilde \gamma, \alpha_t}$, we arrive at the limit equation
\[
0 =  \lim_{t \to \thy} \frac{ \innoneind{\ell}{t}{\tilde \gamma, \alpha_t}}{L_j(t)} = \lim_{t \to \thy} \sum_{i \ge j} \sum_{e \in \pi_i} \frac{\ell_t(e)}{L_j(t)} \tilde \gamma(e) \alpha_t(e) = \innoneind{l}{j}{\gamma, \lalpha\rest{\pi_j}}.
\]
Here, $L_j(t) = \sum_{e\in \pi_j} \ell_t(e)$. Since $\gamma \in H_1(\grm{\pi}{j}(G), \R)$ was arbitrary, we infer that $\lalpha\rest{\pi_j}$ is an exact one-form on $(\grm{\pi}{j}(G), l_j)$, and thus, altogether, $\lalpha$ is exact on the tropical curve $\curve$. The same proof works for the more general case of degenerations of tropical curves in $\mgtrop{\combind{(G, \subface \pi)}}$.
\smallskip

It remains to prove the second claim. We fix an admissible basis for $H = H_1(G, \R)$ in the terminology of Section~\ref{sec:admissible_basis_layered_graphs}. Taking into account the orthogonal decomposition \eqref{eq:OrthDecomposition} and Proposition~\ref{prop:OrthoProjection}, it suffices to show that 
\begin{equation*} \label{eq:ProofExactFormsLimit}
\rmM_\ell^{-1} P_\ell  \,\lalpha = o(1)
\end{equation*}
as $\mgr$, equivalently, $(G, \ell)$, degenerates to $\curve$  in the tropical moduli space. (Here, $\rmM_\ell$, $P_\ell$ are the matrices in \eqref{eq:OrthoProjection}). From Theorem~\ref{thm:PeriodAsymptotics} that we will prove in Section~\ref{sec:GraphPeriodMatrices}, we get
\[
	(\rmM^{-1}_\ell)_{ij} = O\left( L_{\min\{i, j\}}^{-1}(\ell) \right ) 
\]
as $\mgr$ degenerates to $\curve$. Moreover, we claim that for any $j=1, \dots r, \fin$, the restricted vector
\[
	\rmW_j := \Big(P_\ell \, \lalpha \Big)\rest{J^j_\pi} \in \R^{h^j_\pi}
\]
has the following order
\[\rmW_j = o(L_{j}(\ell)).
\] 
Indeed, since the fixed basis of $H = H_1(G, \R)$ is admissible, for any $k\in J^j_\pi$, the restriction of the cycle $\gamma_k$ has zero coefficients on the edges of $\pi_i$, $i <j$ and the restriction to $\pi_j$ belongs to the homology group $H_1(\grm{\pi}{j}(G), \R)$. Now, since $\lalpha$ is an exact form on $\curve$, it follows that
\[
\rmW_j =  \innone{\ell}{\gamma_k, \lalpha} = \sum_{i \ge j} L_i(\ell) \Big ( \innoneind{l}{i}{\gamma_k |_{\pi_i}, \lalpha|_{\pi_i}}  + o(1)  \Big ) = L_j(\ell) o(1) + O(L_{j+1}(\ell)) = o(L_j(\ell))
\]
as $\mgr$ degenerates to $\curve$ (where for $j=r$, we use our convention $r+1 =\fin$). Altogether, the $i$-th part of the column vector $\rmM_\ell^{-1} P_\ell \lalpha\in \R^h$ (i.e., when restricting to coordinates $k$, with $k \in J^i_\pi$) has order
\begin{align*}
(\rmM_\ell^{-1} P_\ell \lalpha)_i &= \sum_{j\in[r]\cup\{\fin\}}  (\rmM^{-1}_\ell)_{ij}\rmW_j = \sum_{j\in [r]\cup\{\fin\}}  L^{-1}_{\min\{i, j\}}(\ell) o(L_j(\ell)) = o(1)
\end{align*}
and this completes the proof.

The proof in the general case is similar.\end{proof}

\subsection{Tame families of tropical functions} \label{ss:TameFunctions} 
In the following, we introduce two notions of tameness for functions on families of tropical curves: the {\em stratumwise tameness} and {\em weak tameness}. 

\medskip

Notations as in previous sections, let $G = (V,E)$ be a fixed (augmented) graph (with marking) and consider the space $\mgtropcombin{\grind{G}}$ of tropical curves of combinatorial type $G$. Let $\unicurvetrop{\grind{G}}$ be the corresponding universal family of tropical curves, and denote by $\unicurvetrop{\thy}$ the fiber at the point $\thy \in \mgtropcombin{\grind{G}}$, that is, the tropical curve represented by $\thy$.

\subsubsection{Weak tameness} \label{ss:TropicalWeakTame}

Let $\lf_\thy$, $\thy \in \mgtropcombin{\grind{G}}$ be a {\em family of tropical functions} on $\unicurvetrop{\grind{G}}$. That is, for each $\thy \in \mgtropcombin{\grind{G}}$, we have a function $\lf_\thy$ on the tropical curve $\unicurvetrop{\thy}$ represented by $\thy$.

\smallskip

By the constructions in Section~\ref{sec:tropical_log_map2}, we have the tropical log map
\[
 \begin{array}{cccc}
 \Lognoind = \logtropind{\combind{G,m}}  &\colon \mgtropcombin{\grind{G}} \setminus {\umggraphcombin{\grind{G}}} &\to &  \partial_\infty \mgtropcombin{\grind{G}}
 \end{array}
\]
to the boundary at infinity $\partial_\infty \mgtropcombin{\grind{G}}  =  \mgtropcombin{\grind{G}} \setminus  \mggraphcombin{\grind{G}}$, where we can either take $m = |E(G)|$ or $m = 3g - 3$ ($g$ is the genus of $G$).

\smallskip

 The boundary $ \partial_\infty \mgtropcombin{\grind{G}} = \mggraphcombin{\grind{G}}  \setminus \mggraphcombin{\grind{G}}$ is stratified as $\partial_\infty \mgtropcombin{\grind{G}}  = \bigsqcup_\pi  \mgtropcombin{\combind{(G, \pi)}}$, which induces the following decomposition of the domain of definition $\mgtropcombin{\grind{G}} \setminus {\umggraphcombin{\grind{G}}}$ of the log map $\Lognoind$,
\begin{equation} \label{eq:DecompositionFromTropicalLogMap}
\mgtropcombin{\grind{G}} \setminus {\umggraphcombin{\grind{G}}} = \bigsqcup_\pi \inn R_\pi, \qquad \qquad \inn R_\pi := \Lognoind^{-1}\Big ( \mgtropcombin{\combind{(G, \pi)}} \Big ).
\end{equation}
Here the union is taken over all ordered partitions $\pi$ with non-empty sedentarity on $E$.

\smallskip

Assume further that $\layh_\thy$, $\thy \in  \mgtropcombin{\grind{G}}$, are tropical functions satisfying the following properties:
\begin{itemize}
\item[(i)] for each $\thy \in \mgtropcombin{\grind{G}}$, the function $\layh_\thy$ is a tropical function \emph{on the tropical curve} $\unicurvetrop{\shy}$ for $\shy = \Lognoind(\thy)$.
\item[(ii)] The equality  $\layh_{\thy} = \lf_\thy$ holds for all $\thy \in \partial_\infty  \mgtropcombin{\grind{G}}  $.
\item[(iii)] Consider the boundary stratum $\mgtropcombin{\combind{(G, \pi)}}$ for an ordered partition $\pi =(\pi_1, \dots, \pi_r, \pi_\fin)$ on $E$. By property (i), the tropical function $\layh_\thy$ is of the form $\layh_\thy = (h_{\thy, 1}, \dots, h_{\thy, r}, h_{\thy, \fin})$  for each base point $\thy$ in the region $\inn R_\pi$.  Then all components $h_{\thy, j}$ , $j=1, \dots, r, \fin$,  depend continuously on the parameter $ \thy \in \inn R_\pi$.
\end{itemize}

For any point $t \in \mggraphcombin{\grind{G}}$, denote by $\layh^\ast_t$ the pull-back of $\layh_t$ from the tropical curve $\unicurvetrop{\shy}$, $\shy = \Lognoind(t)$, to the metric graph $\unicurvetrop{t}$.

\begin{defi}[Weak tameness] \label{def:WeakGlobalTameTropical}
A family of tropical functions $\lf_\thy$, $\thy \in \mgtropcombin{\grind{G}}$, is {\em weakly tame} if there are tropical functions $\layh_\thy$, $\thy \in \mgtropcombin{\grind{G}}$, as above such that, in addition, the difference 
\begin{equation} \label{eq:Differencev}
 \lf_t  -  \layh^\ast_t , \qquad t \in \mggraphcombin{\grind{G}},
 \end{equation}
goes to zero upon approaching the boundary $\partial_\infty \mgtropcombin{\grind{G}}$ in the tame topology.
\end{defi}
More precisely, hereby we mean that for any point $\shy \in \partial_\infty   \mgtropcombin{\grind{G}}$, we have
\[
	\lim_{t \to \shy} \sup_{x \in \unicurvetrop{t}} \Big |  \lf_t  (x)  - \layh^\ast_t (x) \Big | = 0
\]
if $t \in \mggraphcombin{\grind{G}}$ converges to $\shy$ in the tame topology on $\mgtropcombin{\grind{G}}$. By the definition of the pullbacks, it implies that the function $\lf_t$, $t \in \mggraphcombin{\grind{G}}$, has the following \emph{asymptotic expansion}
\[
\lf_t  = \sum_{j=1}^{r_\thy} L_j(t)  h_{ \thy, j}^*  + h_{\thy, \fin}^*  + o(1), \qquad \thy = \Lognoind (t),
\]
close to a boundary point $\shy \in \partial_\infty   \mgtropcombin{\grind{G}}$ in the tame topology on $ \mgtropcombin{\grind{G}}$. The number of terms $r_\thy$ in the additive expansion is the rank of the tropical curve $\unicurvetrop{\thy}$, $\thy = \Lognoind  (t)$, and hence \emph{depends on the point} $t \in \mggraphcombin{\grind{G}}$.


\subsubsection{Stratumwise tame} 
Let $\lf = (\lf_\thy)_\thy$ be a family of tropical functions on $\unicurvetrop{\grind{G}}$, in the sense of the previous section. 

Consider a boundary stratum $\mgtropcombin{\combind{(G, \pi)}}$ associated to an ordered partition $\pi =(\pi_\infty, \pi_\fin)$ of $E$ with non-empty $\pi_\infty =(\pi_1, \dots,  \pi_r)$. Denote by $\pr_\pi \colon  \mggraphcombin{\grind{G}} \to \mgtropcombin{\combind{(G, \pi)}}$ the natural projection map, which, recall, sends every point $t \in  \tilde \mg_{\grind{G}}$ to
\begin{equation} \label{eq:StratumwiseProjection}
\pr_\pi(t) = \Big( [t\rest{\pi_1}], \dots, [t\rest{\pi_r}], t\rest{\pi_\fin}  \Big )
\end{equation}
Here, the brackets denotes projectivization. Alternatively, $\pr_\pi(t)$ is the point representing the conformal equivalence class of the layered metric graph $(\mgr_t, \pi)$.

\medskip

We call the family $\lf = (\lf_\thy)_\thy$ {\em tame at the stratum} $\mgtropcombin{\combind{(G, \pi)}}$ if for any point $\thy$ in the \emph{closure in $\mgtropcombin{\grind{G}}$} of $\mgtropcombin{\combind{(G, \pi)}}$, we have
\begin{equation} \label{eq:StratumwiseTameness}
	\lim_{t \to \thy} \sup_{x \in \unicurvetrop{t}} \big |  \lf_t  (x)  - \lf_{\pr_\pi(t)}^\ast  (x) \big | = 0
\end{equation}
as $t \in  \mggraphcombin{\grind{G}}$ converges to $\thy$ in the tame topology on $\mgtropcombin{\grind{G}}$.

\begin{defi}
A family of functions $\lf_\thy$, $\thy \in \mgtropcombin{\grind{G}}$, is called {\em stratumwise tame} if it is tame at all boundary strata $\mgtropcombin{\combind{(G, \pi)}}$.
\end{defi}

\begin{remark} \label{rem:TamenessVS} We make the following important comment. Using the properties of the tropical logarithm map $\logtrop = \Lognoind$, it can be shown that every continuous, stratumwise tame family of functions $\lf_\thy$, ${\thy \in \mgtropcombin{\grind{G}}}$, is also weakly tame. The validity of this connection is the main reason for allowing the point $\thy$ in \eqref{eq:StratumwiseTameness} to belong to the {\em closure of the stratum} $\mgtropcombin{\combind{(G, \pi)}}$. Indeed, it can be shown that the implication might become wrong if just $\mgtropcombin{\combind{(G, \pi)}}$ is considered.
\end{remark}

\subsection{Layered measures on layered measurable spaces} \label{ss:LayeredMeasures}
The aim of the present section is to introduce the notion of measures on tropical curves. This will be needed in order to define a measure valued Laplacian on tropical curves. Since we will need the same set-up when defining hybrid Laplacian on hybrid curves, we present the results in  a more general framework of \emph{topological spaces $X$ endowed with a layring}; the application we have in mind is in the case where $X$ is either a tropical or a hybrid curve.

\smallskip

In the following we suppose given 
\begin{itemize}
\item A measurable space $X$, and
\item A descending filtration 
\[\mathfrak{X} : \qquad X^1 = X \supsetneq X^2 \supsetneq \dots \supsetneq X^r \supseteq X^\fin\] of $X$, and
\item for each $j = 1, \dots, r, \fin$, a partition $X^j = \bigsqcup_{i=1}^{n_j} X^{j}_i$ for measurable subsets $X_i^j \subseteq X^j$ such that these partitions are globally compatible: each $X^{j+1}_i$ lies entirely in a part $X^j_k$ in the partition of $X^j$. The collection of all these partitions is denoted by $\pi$.
\end{itemize}
In the following, we will refer to such an object as a \emph{layered measurable space}. 

\smallskip

We then  define the graded minors $\grm{\pi}{j}X$ of $X$ as the quotient topological space $\grm{\pi}{j}X := X^j  / \sim$ obtained by contracting the subspace $X^{j+1}_i$, $i=1, \dots, n_{j+1}$, to a point.  (Here, as before, we set $r+1 =\fin$. So for $r=0$, we only have one term $X^\fin =X$.)

\smallskip

In the case $X$ is a topological space, we assume that the filtration consists of closed subsets of $X$ and take the connected components for the decomposition of the $X^j$'s as the partition.

\medskip

 In the following, a measure on a topological space $X$ means a complex-valued Borel measure $\mu$ on $X$ with finite total variation $|\mu| (X) < \infty$. We denote the vector space of all such measures by $\mathcal{M}(X)$ and the subspace of measures with total mass $\mu(X) = 0$ by $\mathcal{M}^0(X)$.

\medskip

 A \emph{layered measure} $\lmu = \lmu_\pi$ on a layered measurable space as above is an $r+1$-tuple $\lmu=(\mu_1, \dots, \mu_r, \mu_\fin)$ consisting of measures $\mu_j$ on the graded minors $\grm{\pi}{j}X$, for $j=1, \dots, r, \fin$.  The set of all layered measures on $X$ is denoted by $\mathcal{M}_\pi(X)$.
 
\medskip

For a given measurable space $X$ and a measure $\mu$ on $X$, we consider its \emph{mass function} $\mass$ which associates to any connected component $H$ of $X$ the mass $\mu(H)$ of $H$. 
 
We generalize the mass function to the layered setting as follows. The mass function $\mass$ of a layered measure $\lmu$ on a layered topological space is defined on connected components of the graded minors $\grm{\pi}{j}X$, $j=1, \dots, r, \fin$. It assigns to each connected component $H$ of a graded minor $\grm{\pi}{j}X$ the value
\[
\mass(H) := \mu_j(H) - \mu_{j-1}(\{x_H\}),
\]
where $x_H$ is the point of the graded minor $\grm{\pi}{j-1}X$ associated to $H$. (For $j=\fin$, we set $\fin-1=r$.)

\smallskip

We say that a layered measure $\lmu$ on $X$ has \emph{mass zero} if the function $\mass$ is null. The set of layered measures of mass zero on $X$ is denoted by $\mathcal{M}^0_\pi(X)$.

\medskip

We now clarify the link between layered measures and usual (non-layered) measures on a topological space $X$.

To every measure $\mu$ on $X$, we can naturally associate a layered measure $\lmu = (\mu_1, \dots, \mu_r, \mu_\fin)$ defined as follows. Consider the filtration $\mathfrak X$ on $X$ given above. Then, we set
\begin{equation} \label{eq:PushoutMeasure}
\mu_j:= \proj_{j*}\Bigl(\mu\rest{X^j}\Bigr)
\end{equation}
where $\proj_j: X^j \to\grm{\pi}j X$ is the natural projection map. In other words, $\mu_j$ is the push-out of the restriction of $\mu$ to the subspace $X^j$ on the graded minor $\grm{\pi}jX$.

\begin{prop} \label{prop:MassZeroMeasures}
Notations as above, the layered measure $\lmu$ has mass zero if and only if the measure $\mu$ on $X$ has mass zero. Moreover, \eqref{eq:PushoutMeasure} defines a bijection between $\mathcal{M}^0(X)$ and $\mathcal{M}_\pi^0(X)$. 
\end{prop}
\begin{proof}
All claims are direct consequences of the definition. For instance, the measure $\mu$ on $X$ corresponding to a layered measure $\lmu=(\mu_1, \dots, \mu_r, \mu_\fin)$ is explicitly given  by
\begin{equation} \label{eq:Gr/NGrMea}
\mu = \iota_{\fin\,*}(\mu_\fin)  + \sum_{j=1}^r \iota_{j*}\bigl(\mu_{j}\rest{X^j \setminus X^{j+1}}\bigr),
\end{equation}
where $\iota_j$ is the inclusion $\iota_j\colon X^j \hookrightarrow X$, and $\iota_{j*}$ is the pushout of the measure to $X$. 
\end{proof}

\begin{remark} 
By the previous proposition, when assuming mass zero, we can actually identify layered measures and usual (non-layered) measures on a layered topological space. However, the layered perspective turns out to be more convenient in context with the Poisson equation on degenerating metric graphs and tropical curves.
\end{remark}

We will apply the above discussion to the setting of a tropical curve $\curve$ (and later to the setting of hybrid curves) with underlying ordered partition $\pi$. In this case, a layered measure is the data of a tuple $\lmu=(\mu_1, \dots, \mu_r, \mu_\fin)$ consisting of measures $\mu_j$ on the graded minors $\Gamma_\pi^j$, for $j=1, \dots, r, \fin$.  The set of all layered measures on $\curve$ is denoted by $\mathcal{M}_\pi(\curve)$, and those of mass zero with $\mathcal{M}_\pi^0(\curve)$. Moreover, $\widetilde{\mathcal{M}}_\pi^0(\curve)$ is the space of layered measures $\lmu = (\mu_j)_j$ in $\mathcal{M}_\pi^0(\curve)$ such that $\mu_j$ is of the form \eqref{eq:SpecialMus} for all $j$ (see \eqref{eq:NiceMeasures}).

\section{Tropical Laplacian and Poisson equation} \label{sec:tropical_laplacian}

In this section, we define the tropical Laplace operator and show how it arises as the limit of the Laplace operator on degenerating metric graphs. We further study the behavior of solutions to the tropical Poisson equation on degenerating families of tropical curves.

\subsection{Definition} Consider a tropical curve $\curve$ lying on the stratum $\mgtrop{\combind{(G,\pi)}}$ of the tropical moduli space $\mgtrop{\grind{g,n}}$, with underlying graph $G=(V, E)$, edge length function $l$, and ordered partition $\pi=(\pi_\infty, \pi_\fin)$, $\pi_\infty=(\pi_1, \dots, \pi_r)$, of rank $r$. Denote by $\Gamma^j$, $j\in[r] \cup\{\fin\}$, the graded minors of $\curve$. We suppose that the edge length function $l$ is normalized, that is $l_j = l\rest{\pi_j}$ has total length one for $j\in[r]$. 

\smallskip

The \emph{tropical Laplacian} $\Deltatrop=\Deltatrop_{\grind{\curve}}$ on $\curve$ is a measure-valued operator which maps piecewise smooth functions $\lf =(f_1, \dots, f_r, f_\fin)$ on $\curve$ to layered measures.  It is defined as the sum
\begin{equation} \label{eq:TropicalLaplacian}
\Deltatrop (\lf) = \Deltatrop_1(f_1) + \Deltatrop_2(f_2) + \dots \Deltatrop_r(f_r) + \Deltatrop_\fin(f_\fin)
\end{equation}
where for each $j=1, \dots, r, \fin$, and for a function $f_j$ on the graded minor $\Gamma^j$, the layered measure $\Deltatrop_j(f_j)$ on the tropical curve $\curve$ is given by
\[
\Deltatrop_j(f_j)  := \Bigl(0, \dots, 0, \Delta_j (f_j),  \divind{j}{j+1}(f_j), \dots,  \divind{j}{r}(f_j),\divind{j}{\fin} (f_j) \Bigr).
\]
Here, 
\begin{itemize}
\item $\Delta_j (f_j) \in \widetilde{\mathcal{M}}^0(\Gamma^j)$ denotes the Laplacian on the $j$-th graded minor $\Gamma^j$, and
\smallskip

\item for the higher graded minors $\Gamma^i$ with $i>j$, the point measure $\divind{j}{i} (f_j)$ is defined as
\[
\divind{j}{i}(f_j)  := -  \sum_{\substack{e \in \pi_j} } \sum_{\substack{v \in e}} \slp_e f_j (v) \, \delta_{\proj_i(v)},
\]
where, $\proj_i$ denotes the projection map $\proj_i:V\to V(\Gamma^i)$, and $\slp_e f_j(v)$ is the slope of $f_j$ at $v$ along the unit tangent direction at the incident edge $e$ to $v$.
\end{itemize}

\medskip

The following proposition is immediate from the definition of $\Deltatrop$.
\begin{prop} \label{eq:LaplacianMassZero}
For any function $\lf$ on a tropical curve $\curve$, the measure $\Deltatrop(\lf)$ has mass zero. 
\end{prop}

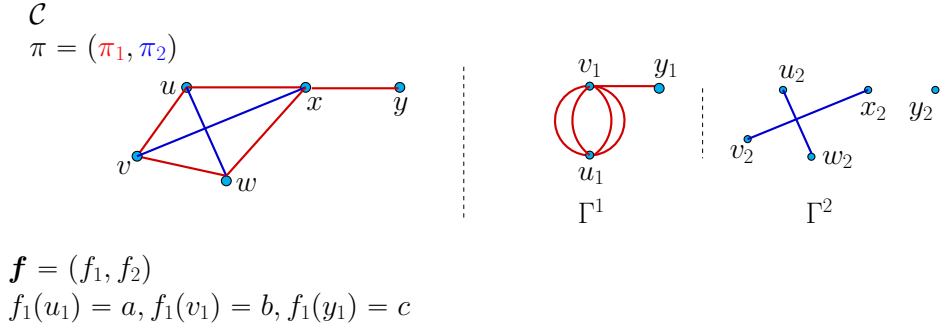
\begin{figure}[!t]
\centering
    \scalebox{.45}{\input{example10.tikz}}
\caption{A tropical curve of full sedentarity with $\pi_\infty=(\pi_1, \pi_2)$, and its corresponding graded minors. The edges in each layer have the same length.} 
\label{fig:kite_Laplacian}
\end{figure}

\begin{example} Consider the tropical curve depicted in Figure~\ref{fig:kite_Laplacian}. Consider the function $\lf=(f_1,f_2)$ on $\curve$ with components $f_1$ and $f_2$ defined on the first and second minors, respectively. Assume $f_1$ is affine linear taking values $a,b,c$ on the vertices $u_1, v_1$, and $y_1$ of $\Gamma^1$, respectively. Then, we have $\Deltatrop(\lf) = (\mu_1, \mu_2)$ with 
\[\mu_1 = 5\cdot(4b-4a) \delta_{u_1} + 5\cdot(4a+c -5b) \delta_{v_1} + 5\cdot(a -c) \delta_{y_1}, \qquad \textrm{and}\] 
\[\mu_2 = \Delta_2(f_2) + 5\cdot(2b-2a) \delta_{u_2} + 5\cdot (2a -2b) \delta_{v_2} + 5\cdot (2b-2a) \delta_{w_2} + 5\cdot (2a- 3b + c) \delta_{x_2} + 5\cdot (b-c)\delta_{y_2}. \]
\end{example}
\subsection{Tropical Laplacian as a weak limit of metric graph Laplacians} \label{ss:TropicalLaplacianAsLimit}

It turns out that the tropical Laplacian $\Deltatrop$ appears naturally as the limit of the Laplacian on degenerating metric graphs. In order to make this statement precise, we first need to formalize how tropical curves can be seen as limits of metric graphs. We can and will work with the space $\mgtropcombin{\grind{G}}$ introduced in Section~\ref{sec:tropical_moduli}, which provides an \'etale chart for the stratum $\mgtrop{\grind{G}}$ of $\mgtrop{\grind{g,n}}$. Recall that, as a set, $\mgtropcombin{\grind{G}}$ is defined as
\begin{equation} \label{eq:TildeMgTrop}
\mgtropcombin{\grind{G}} =  \bigsqcup_{\pi\in\widehat \Pi(E)} \inn\keg_\pi =\bigsqcup_{\pi\in \widehat\Pi(E)}   \inn\sigma_{\pi_\infty} \times \inn\keg_{\pi_\fin},
\end{equation}
with $\inn\sigma_{\pi_\infty} = \inn\sigma_{\pi_1}\times \dots\times \inn\sigma_{\pi_r}$ provided that $\pi_\infty=(\pi_1, \dots, \pi_r)$.

\smallskip

In terms of this topology, there is a rather simple description of the degeneration of metric graphs to tropical curves. More precisely, let $(t_n)_n$ be a sequence of points in $\mggraphcombin{\grind{G}}$ representing metric graphs $\mgr_n$ with edge function $\ell_n$ on $E$, and let $\thy$ be a point of $\mgtropcombin{\grind{G}}$ representing a tropical curve $\curve$ with the edge length function $l$. Then, $t_n$ converges to $\thy$ in $\mgtropcombin{\grind{G}}$ (and so in $\mgtrop{\grind{g,n}}$) exactly when
\begin{itemize}
\item [(i)] for all edges $e \in E_\infty =\pi_1 \sqcup \dots \sqcup \pi_r$,
\[
\lim_{n \to \infty} \ell_n(e) = + \infty,
\] 
\item [(ii)] and for all $i=1, \dots, r$ (see \eqref{eq:DefGrLength} for the definition of $L_i$),
\[
	\lim_{n \to \infty} \frac{L_{i+1}(n)}{L_i(n)}  = 0,
\]
	with $L_\fin \equiv 1$,

\item [(ii)] and for each edge $e$ in each layer $\pi_i$, $i=1, \dots, r, \fin$,
\[
\lim_{n \to \infty} \frac{\ell_n(e)}{L_i(n)} = l (e).
\]
\end{itemize}
Moreover, note that if $t_n$ converges to $\thy$ in the tame topology on $\mgtropcombin{\grind{G}}$, then in addition to (i)--(iii), we get
\begin{itemize}
\item [(iv)] for all $i=1, \dots, r$,
\[
	\lim_{n \to \infty} \frac{L_{i+1}(n)}{L_i(n)^2}  = 0.
\]
\end{itemize}

\smallskip

\begin{thm} \label{thm:GraphLaplacianConvergence}
Suppose $\curve = \hcurve_\thy^\trop$ is a tropical curve corresponding to the point $\thy$ lying on the stratum $\mgtrop{\combind{(G, \pi)}}$, $\pi\in \Piall(E)$, of $\mgtrop{\grind{G}}$. Let $\lf = (f_j)_j$ be a  function on $\curve$. Let $\umgr_t$, $t\in \mggraph{\grind{G}}$, be the metric graph of combinatorial type $G$ associated to the point $t$, and let $f_t = \lf^\ast$ be the pull-back of $\lf$ to $\umgr_t$. Then, as $\umgr_t$ degenerates to $\curve$ as $t$ tends to $\thy$ in $\mgtrop{\grind{G}}$, we get
\[
\Delta(f_t)  \to \Deltatrop(\lf)
\]
in the weak sense.
\end{thm}

In this setting, convergence in the weak sense means, taking the universal curve $\hcurveg{G}^\trop$ over the chart $\mgtropcombin{\grind{G}}$ of $\mgtrop{\grind G}$, that for any continuous function $h \colon \hcurveg{G}^\trop \to \R$, we have
\[
 \int_{\umgr_t} h(s) \,  \Delta f_t \quad  \longrightarrow  \quad \int_{\hcurve_{\thy}^\trop} h(s)  \, \Deltatrop \lf
\]
as $t \in \mggraphcombin{\grind{G}}$ tends to $ \thy$ in $\mgtropcombin{\grind{G}}$. On the right hand side, $\Deltatrop \lf$ is understood as a measure of mass zero on the tropical curve $\curve=\hcurve_\thy^\trop$ (see Proposition~\ref{prop:MassZeroMeasures} and Proposition~\ref{eq:LaplacianMassZero}).

\begin{proof} We will work in the chart $\mgtropcombin{\grind{G}}$ and denote by $\ell_t$ the edge length function of $\umgr_t$. Taking into account the (additive) definition of $\Deltatrop = \Deltatrop_1 + \dots \Deltatrop_r+\Deltatrop_\fin$, it suffices to prove the claim in the case that $\lf = (0, \dots, 0, f_j, 0, \dots)$. 

\smallskip

First of all, notice that the Laplacian of $f_j^\ast$ is the measure on $\umgr_t$ given by
\[
- \Delta f_j^\ast = \sum_{e \in \pi_j}  \frac{l(e)^2}{\ell_t (e)^2} \, \frac{d^2 f_j\rest{e}}{dx^2} \Big (\frac{l(e)}{\ell_t(e)} x \Big ) +  \sum_{i \le j} \sum_{e \in \pi_i} \sum_{v \in e} \frac{l (e)}{\ell_t(e)} \slp_e f_j (v) \, \delta_{v}.
\]
Integrating against a function $h \colon \hcurveg{G}^\trop \to \R $, we arrive at the expression
\begin{align*}
 - \int_{\umgr_t} h \,  \Delta f_j^\ast &= \sum_{e \in \pi_j}  \frac{l (e)^2}{\ell_t(e)} \int_0^1 h\rest{e} (l (e) x) \frac{d^2 f_j\rest{e}}{dx^2} (l (e) x) \, d x \\
 &+  \sum_{i \le j} \sum_{e \in \pi_i} \sum_{v \in e} \frac{l(e)}{\ell_t(e)} \slp_e f_j (v) \, h(v).
\end{align*}
Taking into account the topology on $\mgtropcombin{\grind{G}}$, we conclude that $\ell_t (e) / L_i(t)$ goes to $l(e)$ for all $e  \in \pi_i$ and $L_j(t)  / L_i(t)$ goes to zero for all $i < j$. Hence after multiplying with $L_j$ and taking the limit, we end up with 
\begin{align*}
-\lim_{t\to \thy} \int_{\umgr_t} h \,  \Delta f_t &= -\lim_{t\to \thy} L_j(t)\int_{\umgr_t} h \,  \Delta f^\ast_j =
\sum_{e \in \pi_j} \int_0^{l(e)} h\rest{e} (x) \frac{d^2 f_j \rest{e}}{dx^2} ( x) \, d x + \sum_{e \in \pi_j} \sum_{v \in e} \slp_e f_j (v) \, h(v) \\
&= -\int_{\hcurve_{\thy}^\trop} h \, \Deltatrop \lf, 
\end{align*}
where the last equality follows from \eqref{eq:Gr/NGrMea}.
\end{proof}

\begin{remark}\label{rem:justification_tropical_Laplacian} Theorem~\ref{thm:GraphLaplacianConvergence} is in fact the main motivation for our definition of the Laplacian $\Deltatrop$ on a tropical curve $\curve$. A naive approach to introducing such a Laplacian would be, informally speaking, to set
\[
	\Deltaind{\curve} (\lf) := \lim \Deltaind{\mgr} f_\mgr,
\] 
where $f_\mgr$ is a suitable extension of a function $f \colon \Gamma \to \R$ to metric graphs $\mgr$ degenerating to $\curve$. However, as the proof of Theorem~\ref{thm:GraphLaplacianConvergence} shows, the corresponding limit will be the zero measure on the tropical curve $\curve$. On the other hand, multiplying the function $f$ by the lengths $L_r$, we will get a well-defined and non-zero limit contribution precisely on the edges on the last layer $\pi_r$. Multiplying by $L_{r-1}$, there is a limit contribution on the edges on the first layer $\pi_{r-1}$, but an infinite contribution on the edges of $\pi_r$ - except if the function $f$ is constant on the edges of $\pi_r$. Proceeding like this, we conclude that, in order to get well-defined limits on all edges of $\curve$, we should look at functions $f_j$ which are constant on edges of layers $\pi_i$ with $i >j$ and rescale them by $L_j$.

This observation motivates our definition of functions $\lf = (f_j)_j$, their rescaled pullbacks $\lf^\ast$, and the definition of the tropical Laplacian $\Deltatrop$.
\end{remark}

\subsection{Tropical Poisson equation}
In the rest of this section, we investigate the Poisson equation on degenerating families of tropical curves. 

\smallskip

Suppose $\mgr$ is a metric graph with underlying combinatorial graph $G = (V,E)$. Let $\mu$ and $ \nu$ be Borel measures on $\mgr$ such that
\begin{equation} \label{eq:AssGraphMeasures}
\mu \in \widetilde{\mathcal{M}}_0(\mgr) \text{, i.e.,} \,\, \mu(\mgr) = 0,  \qquad \qquad \nu \in {\mathcal{M}}(\mgr) \text{ and } \nu(\mgr) = 1.
\end{equation}
In the above situation, we sometimes call the triple $(\mgr, \mu, \nu)$ a {\em measured metric graph}.

By the discussion in Section~\ref{sec:potential_theory_metric_graphs}, the associated {\em Poisson equation}
 \begin{equation}  \label{eq:PoissonMetricGraphs}
 \begin{cases}
 \Deltaind{\mgr} f = \mu \\
\int_{\mgr} f \, d \nu = 0
 \end{cases}
\end{equation}
has a unique solution $f \colon \mgr \to \C$, and $f$ belongs to the Zhang space $D(\Deltaind{\mgr})$.

\begin{remark}
In fact, the smoothness assumption $\mu \in \widetilde{\mathcal{M}}(\mgr)$ can be omitted. However, then solutions may no longer exist in the Zhang space and the Laplacian $\Delta$ has to be considered on the larger space $\operatorname{BDV}(\mgr)$ (see Remark~\ref{rem:DistributionalLaplacian}). With this modification, the subsequent discussion extends verbatim to the case of a general Borel measure $\mu \in \mathcal{M}(\mgr)$ on $\mgr$ of mass $\mu(\mgr) = 0$. However, the definition of $\Delta$ on $\operatorname{BDV}(\mgr)$ is slightly involved, and, in order to streamline the exposition, we do not develop this here.
\end{remark}

We consider the following question:
\begin{question}\label{question:limit_Poisson-metric_graph} What is the behavior of the solution $f$ when the measured metric graph $(\mgr, \mu_\mgr, \nu_\mgr)$ degenerates to a measured tropical curve $(\curve, \mu, \nu)$?
\end{question}

The aim of this section and the next one is to provide an answer to the above question. 

In order to do so, we first need to introduce a {\em Poisson equation on tropical curves}, and we have all the needed ingredients at our disposal. Suppose that $\curve  = (G, \pi, l)$ is a tropical curve of rank $r$. Let further $\lmu$ and $\lnu$ be two layered Borel measures on $\curve$ such that
\begin{equation} \label{eq:AssTropMeasures}
\lmu \in \widetilde{\mathcal{M}}^0(\curve),  \text{ and } \qquad \lnu \in \mathcal{M}(\curve) \text{ and } \mass_\lnu(\Gamma^1) = 1, \mass_\lmu(H) = 0 \text{ for all } H \neq \Gamma^1.
\end{equation}
That is, $\lnu$ \emph{has total mass one}. 
\medskip

We call the triplet $(\curve, \lmu, \lnu)$ a {\em measured tropical curve}. We consider the following {\em tropical Poisson equation}
 \begin{equation}  \label{eq:PoissonTropicalCurve}
 \begin{cases}
 \Deltatrop \lf = \lmu \\
 \lf \text{ is harmonically arranged} \\
\int_{\curve} \lf \, d \lnu = 0.
 \end{cases}
\end{equation}
 The third condition means that
\[
\int_\Gamma f_j \, d\nu = 0 \qquad \text{for all }j=1, \dots, r, \fin,
\]
where $\lnu$ is identified with a (non-layered) measure $\nu$ on $\Gamma$ of total mass $\nu(\Gamma) = 1$ (which is a consequence of the assumption \eqref{eq:AssTropMeasures} on $\lnu$), and every function $f_j$ on $\Gamma^j$ is identified with a function on $\Gamma$ as we did in formulating the harmonic rearrangement theorem. We have the following.

\begin{thm} The tropical Poisson equation \refeq{eq:PoissonTropicalCurve} has a unique solution $\lf = (f_1, \dots, f_r, f_\fin)$ on $\curve$.
\end{thm}
\begin{proof} Proceeding recursively, from lower index layers to upper index ones, the equation $\Deltatrop \lf = \lmu$ is shown to have a family of solutions by using the same result for Poisson equations on metric graphs (in our situation applied to minors of $\curve$). Two different solutions are related by a rearrangement.
The second property, that of being harmonically arranged, then shows that the solution is unique up to a global constant, that is, up to an additive constant on each graded minor $\Gamma^j$. The last condition, the normalization property, fixes all these constants. 
\end{proof}

\subsection{Poisson equation on degenerating metric graphs: tameness results} We will now formalize the degeneration of a measured metric graph $(\mgr, \mu, \nu)$ by working in the chart $\mgtropcombin{\grind{G}}$. Let $(\lmu_\thy)_\thy$ and $(\lnu_\thy)_\thy$ be continuous families of measures on $\unicurvetrop{\grind{G}} / \mgtropcombin{\grind{G}}$ (in the sense of Section~\ref{ss:FamiliesMeasuredSpaces}) satisfying \eqref{eq:AssTropMeasures}. For each point $\thy \in\mgtropcombin{\grind{G}}$, denote by $\lf_\thy$ the solution to the Poisson equation \eqref{eq:PoissonTropicalCurve} on the associated measured tropical curve $(\curve_\thy,  \mu_\thy, \nu_\thy)$.

\smallskip

We start by stating the following basic result whose proof will be given in Section~\ref{sec:proof_proposition_leading_term}. It shows that the leading term under a degeneration $(\mgr, \mu, \nu)\to (\curve, \lmu, \lnu)$ is given by the solution on the first minor.

\begin{prop} \label{prop:leadingterm}
Let $(\lmu_\thy)_\thy$ and $(\lnu_\thy)_\thy$ be continuous families of measures on $\unicurvetrop{\grind{G}} / \mgtropcombin{\grind{G}}$ satisfying \eqref{eq:AssTropMeasures}. For $\thy \in  \mgtropcombin{\grind{G}}$, let $h_\thy \colon \unicurvetrop{\thy} \to \R$ be the (non-layered) function 
\[
h_\thy = \begin{cases}L(t)^{-1} f_t, &t \in \inn\eta^\trop_{\grind{G}}  \\ f_{\thy, 1}, &\thy \in \partial_\infty \mgtropcombin{\grind{G}} \end{cases}
\]
on its fiber $\unicurvetrop{\thy}$ in the universal curve $\unicurvetrop{\grind{G}}$. (Recall that $\partial_\infty \mgtropcombin{\grind{G}}$ is the set of all the points with non-empty sedentarity.) Then the collection of functions $h_\thy$, ${\thy \in \mgtropcombin{\grind{G}}}$, forms a continuous family of functions on $\unicurvetrop{\grind{G}} / \mgtropcombin{\grind{G}}$. 
\smallskip

In fact, the following strong phenomenon happens. If $\shy$ converges to $\thy$ in $\mgtropcombin{\grind{G}}$, then 
\[
\lim_{\shy \to \thy} \| h_{\shy} - h_{\thy} \|_{\infty, \unicurvetrop{\shy}} = 0.
\]
\end{prop}
In the above equation, in order to take the $\infty$-norm, we take the pull-back of functions $h_\shy$ from the fiber $\unicurvetrop{\thy}$ to $\unicurvetrop{\shy}$. 

\medskip

Our aim next is to obtain a refined description using solutions to the tropical Poisson equation. We will provide three types of results. In a sense, all these results are sharp, and together, they lead to an essentially complete solution to Question~\ref{question:limit_Poisson-metric_graph}. This being said, we will complement them later when discussing Green functions.

\medskip

In our first stated theorem, we will assume that the degeneration $(\mgr, \mu, \nu)\to (\curve, \lmu, \lnu)$ happens sufficiently \emph{fast} (formulated in the statement). We then get the following result.

\begin{thm}[Stratumwise tameness of solutions] \label{thm:main_graph_strong}
Let $(\lmu_\thy)_\thy$ and $(\lnu_\thy)_\thy$ be continuous families of measures on $\unicurvetrop{\grind{G}} / \mgtropcombin{\grind{G}}$ satisfying \eqref{eq:AssTropMeasures}. Suppose that, for every boundary stratum $\mgtrop{\combind{(G, \pi)}}$, we have
\begin{equation} \label{eq:condition_strong_tame}
L(t) \big (\lmu_t - \lmu_{\pr_\pi(t)}^\ast\big) \qquad \text{and} \qquad L(t) \big (\lnu_t - \lnu_{\pr_\pi(t)}^\ast\big)
\end{equation}
converge weakly to zero provided that the point $t \in \mggraphcombin{\grind{G}}$ tends to a point $\thy$ in the closure of $\mgtropcombin{\combind{(G, \pi)}}$ in the tame topology of $\mgtropcombin{\grind{G}}$ .

\smallskip

Then, the solutions $\lf_\thy$, ${\thy \in \mgtropcombin{\grind{G}}}$, to \eqref{eq:PoissonTropicalCurve} form a stratumwise tame family of functions.
\end{thm}

\medskip

In the above, $\lmu_{\pr_\pi(t)}^\ast$ and $\lnu_{\pr_\pi(t)}^\ast$ denote the natural pull-backs of the measures $\lmu_{\pr_\pi(t)}$ and $\lnu_{\pr_\pi(t)}$ to the fiber at $t$ in the universal curve $\unicurvetrop{\grind{G}}$, that is to the associated metric graph. Moreover,  the weak convergence to zero in \eqref{eq:condition_strong_tame} is understood in the distributional sense, that is, we require
\[
\lim_{t \to \thy} L(t) \int_{\mgr_t} f(q) \, d(\lmu_t - \lmu_{\pr_\pi (t)}^\ast)(q)  = 0
\]
for all continuous functions $f \colon \unicurvetrop{\grind{G}} \to \R$, and similar for the measures $\nu_\thy$.

\smallskip

We emphasize that, by definition, the stratumwise tameness of the solutions implies the following (see Section \ref{ss:TameFunctions}): if a sequence of metric graphs $(\mgr_n)_n$ converges to a tropical curve $\curve= (G, \pi,\ell)$ in the sense that $\lim_{n \to \infty} t_n = \thy$ tamely for associated points in $\mgtropcombin{\grind{G}}$, then
\[
	\lim_{n \to \infty} \sup_{p \in \mgr_n} \big |  \lf_{t_n}  (p)  - \lf_{\pr_\pi(t_n)}^\ast  (p) \big | = 0,
\]
for the pullback $ \lf_{\pr_\pi (t_n)}^\ast$ obtained by using the projection map $\pr_\pi \colon \mggraphcombin{\grind{G}} \to \mgtropcombin{\combind{(G, \pi)}}$.

\smallskip

Hence, the above theorem provides a full description of the asymptotic of the solutions in this case.

\medskip

For the sake of illustration, we restate the result in the case of discrete measures.

\begin{thm}[Stratumwise tameness for discrete measures] \label{thm:MainGraphDiscrete}
Let $D \in \Div^0(G)$ be a degree zero divisor and let $v \in V$ be a fixed vertex. We view $D$ and $v$ as measures of total mass zero and one, respectively. For any point $\thy \in \mgtropcombin{\grind{G}}$, let $\lf_\thy$ be the unique solution to the equations
\begin{equation}\label{eq:JDFctTrop}
\begin{cases}
 \Deltatrop (\lf) = D \\
 \lf \text{ is harmonically arranged} \\
 f_j (\proj_j(v)) = 0, \qquad j=1, \dots, r, \fin,
 \end{cases}
\end{equation}
on its associated tropical curve  $\curve_\thy$. 

Then the solutions $\lf_\thy$, ${\thy \in \mgtropcombin{\grind{G}}}$, form a stratumwise tame family of functions.
\end{thm}

We will complement this statement in Theorem~\ref{thm:layered_expansion_tropical_laplacian} below, where we give an expansion of the solutions around each point $\curve$ of the tropical moduli space $\mgtropcombin{\grind{G}}$. 
\medskip

In the above results, we are assuming some conditions on the speed of convergence of the measures. 
Dropping the condition \eqref{eq:condition_strong_tame}, we get the following most general result.

\begin{thm}[Global weak tameness of solutions] \label{thm:MainGraphWeak}
Let $(\mu_\thy)_\thy$ and $(\nu_\thy)_\thy$ be continuous families of measures on $\unicurvetrop{\grind{G}}/\mgtropcombin{\grind{G}}$ satisfying \eqref{eq:AssTropMeasures}. Then, the solutions $\lf_\thy$, $\thy \in  \mgtropcombin{\grind{G}}$, to the tropical Poisson equation form a weakly tame family of functions on $\mgtropcombin{\grind{G}}$.
\end{thm}

We make a couple of remarks. 

\begin{remark}[Sharpness of the results]
\begin{itemize}
\item [(i)] We stress that, informally speaking, stratumwise tameness fails if the convergence of the metric graph measures $\mu_t$ and $\nu_t$ to tropical curve measures $\lmu_\thy$ and $\nu_\thy$ is not fast enough, cf. Condition \eqref{eq:condition_strong_tame} in Theorem~\ref{thm:main_graph_strong}. In fact, since this convergence can be arbitrarily slow, one might expect the weak tameness result to be the only general statement.  Moreover, by choosing $\mu$ and $\nu$ as discrete measures with masses within the edges of the first layer $\pi_1$, one can show that the condition for stratumwise tameness in \eqref{thm:main_graph_strong} is in a certain sense sharp. 
\smallskip

\item [(ii)] On the other hand, even in the situation where the convergence assumption \eqref{eq:condition_strong_tame} fails, it is sometimes possible to find the full asymptotic of solutions by using a higher order approximation on the measure. We demonstrate this approach later by the examples of Green functions associated to the Lebesgue and the canonical measures (see Proposition~\ref{thm:Lebesgue_Green_functions} and Proposition~\ref{prop:simple_layering_Green_functions}).
\end{itemize}
\end{remark}

The rest of this section contains the proofs of the above results. It is divided into two parts: treating first the case of discrete measures, we prove Theorem~\ref{thm:MainGraphDiscrete} and analyze the $j$-function under degeneration (see Lemma~\ref{lem:JfunctionTame}). The remaining results then follow from the integration formula \eqref{eq:SolutionFormulaJFunction} and its tropical curve analog (see Lemma~\ref{lem:SolutionFormulaJFunctionTropical}).

\smallskip

In order to prove Theorem \ref{thm:MainGraphDiscrete}, we establish a suitable local result which provides a local expansion of the solution of Poisson equation around a fixed tropical curve.

To state the theorem, let $\curve$ be a fixed tropical curve in the boundary stratum $\mgtropcombin{\combind{(G, \pi)}}$ for an ordered partition $\pi = (\pi_\infty =(\pi_1, \dots, \pi_r), \pi_\fin)$. For a metric graph $\mgr$ with underlying model $G$ with edge lengths $\ell: E \to (0, + \infty)$, let $\curve_\ell$ be its conformally equivalent tropical curve of combinatorial type $(G, \pi)$, that is, $\curve_\ell$ is the conformal class of the layered metric graph $(G, \ell, \pi)$. Let $D \in \Div^0(G)$ be a degree zero divisor and let $v \in V$ be a fixed vertex. Denote by $\lf_\ell = (f_{\ell,1}, \dots, f_{\ell,r}, f_{\ell,\fin})$ the solution of the Poisson equation \eqref{eq:JDFctTrop} on the tropical curve $\curve_\ell$.

\begin{thm}[Local expansion of solutions]\label{thm:layered_expansion_tropical_laplacian} Let $\curve$ be a fixed tropical curve.  As the metric graph $\mgr$ degenerates to $\curve$, that is, as $t$ converges to $\thy$ in $\mgtropcombin{\grind{G}}$ for associated points in $\mgtropcombin{\grind{G}}$, the function $f_\mgr$ has the following asymptotic expansion
\begin{equation} \label{eq:ExpansionConsideringLayered}
	f_\mgr =  L_1(\ell) f_{\ell,1}^\ast + \dots + L_r(\ell) f_{\ell,r}^\ast + f_{\ell,\fin}^\ast + O \Big ( \max_{j=1, \dots, r, \fin} \frac{L_{j+1}(\ell)^2}{L_j(\ell)}\Big),
\end{equation}
where the error term is of the claimed order in the sup-norm on $\mgr$, $\|\psi \|_{\mgr, \infty} := \sup_{p \in \mgr} |\psi(p)|$.

Moreover, $\lf_\ell = (f_{\ell,1}, \dots, f_{\ell,r}, f_{\ell,\fin}) $ converges uniformly to $\lf_\curve = (f_{1}, \dots, f_{r}, f_{\fin})$, that is, $f_{\ell,j}$ converges uniformly to $f_j$ for all $j=1, \dots, r, \fin$.
\end{thm}

By our convention, $L_{\fin}(\ell)=1$ and $L_{\fin + 1}(\ell)=0$. The symbol $h^\ast$ for a function on the $j$-th graded minor of $\curve$ refers to the pullback function on any metric $\mgr$ of the same combinatorial type (see Section~\ref{sss:PullbackExplanation}).

\subsection{The case of discrete measures: proof of Theorems~\ref{thm:layered_expansion_tropical_laplacian} and~\ref{thm:MainGraphDiscrete}} In this section we treat the case of discrete measures. Let $G = (V, E)$ be a finite graph. Fix further a degree zero divisor $D \in \Div^0(G)$ and a vertex $v \in V$. Suppose $\curve$ is a tropical curve with underlying combinatorial graph $G$, ordered partition $\pi = (\pi_1, \dots, \pi_r, \pi_\fin)$ and edge length function $l = (l_j)_{j\in[r]\cup\{\fin\}}$.

\subsubsection{Asymptotic analysis of differentials}
In order to prove Theorem~\ref{thm:layered_expansion_tropical_laplacian} we undertake a careful asymptotic analysis of the differentials $d(f_\mgr)$ of the solutions $f_\mgr$ the Poisson equation \eqref{eq:JDFctTrop}, when the metric graph $\mgr$ varies. We start with the following basic observation.

\begin{prop} \label{thm:ContForms} Let $\curve$ be a fixed tropical curve.
As the metric graph $\mgr$ degenerates to $\curve$ in $\mgtropcombin{\grind{G}}$, the differentials $d(f_\mgr)$ converge to $\ld(\lf_\curve)$. More precisely, 
\[
	\lim_{\mgr \to \curve} d(f_\mgr) (e) = \ld(\lf_\curve) (e) \qquad \text{for all edges } e \in E
\]
where the limit is taken in $\mgtropcombin{\grind{G}}$ for points representing $\mgr$ and $\curve$.

More generally, denoting by $\lf_\thy$ the solution to the Poisson equation \eqref{eq:JDFctTrop} for $\thy \in \mgtropcombin{\grind{G}}$, the evaluation maps $\thy \mapsto \ld(\lf_\thy)(e)$, for $e\in E$, are  all continuous on $\mgtropcombin{\grind{G}}$. 
\end{prop}

\begin{proof} We drop the index $\curve$ from $\lf_\curve$. 
The theorem is a consequence of Proposition~\ref{prop:exact_forms_limits}. Namely, for each metric graph $\mgr$ of combinatorial type $G$ with edge length function $\ell$ close to the tropical curve $\curve$ in the tropical moduli space, we can consider $\ld(\lf)$ as a one-form on $\mgr$. By Proposition~\ref{prop:exact_forms_limits}, the projection $\projexact(\ld(\lf))$ onto the subspace of exact one-forms on $(G,\ell)$ converges to $\ld(\lf)$. On the other hand, since $\lf$ is a solution of \eqref{eq:JDFctTrop}, we have $\partial \ld(\lf) = D$ (on $G$) and this implies that  $\projexact(\ld (\lf)) = d f_\mgr$. The second statement follows similarly by applying  Proposition~\ref{prop:exact_forms_limits} to a sequence of tropical curves converging to $\curve$.
\end{proof}

The next lemma refines the asymptotics of the differentials $d(f_\mgr)$. Theorem \ref{thm:layered_expansion_tropical_laplacian} then follows immediately by integration.

Recall that for a metric graph $\mgr$ of combinatorial type $G$ with edge length function $\ell$, we denote by $\curve_\ell$ the tropical curve given by the conformal equivalence class of the triple $(G, \pi, \ell)$, and we denote by  $\lf_\ell = (f_{\ell,1}, \dots, f_{\ell,r}, f_{\ell,\fin})$ the corresponding solution of the Poisson equation \eqref{eq:JDFctTrop}. 
 
\begin{lem}[Local Expansion Lemma for Differentials] \label{lem:DerivativeAsymptotics}
Notations as above, if the edge $e$ belongs to the $j$-th set $\pi_j$, $j=1,\dots,r, \fin$, then
\[
	d f_\mgr (e) =  \sum_{k=1, \dots, r,\,\fin} L_k(\ell) d \Big( f_{\ell,k}^\ast \Big )  (e) + R_e(\ell)
\]
as $\mgr$ degenerates to the tropical curve $\curve$, where the error term $R_e(\ell)$ is of order
\begin{equation} \label{eq:DerivativeErrors}
R_e(\ell) = O \Bigl ( \max_{k < j} \frac{L_{k+1}}{L_{k}} (\ell) + \max_{k \ge j} \max_{s \le k} \frac{L_{k+1}L_{s+1}}{L_j L_s} (\ell) \Bigr ).
\end{equation}
\end{lem}
\begin{proof}
By definition, $d f_\mgr$ is the unique exact one-form on $\mgr$ which solves the equation 
\[
	\partial \omega = D \qquad \text{on }\mgr. 
\]
On the other hand, the differential equation $\Deltatrop \bigl(\lf_l\bigr) = D$ implies that, considering $\ld\bigl(\lf_l\bigr) $ as a one-form on $G$, we have
\[
	\partial \ld(\lf_l) = D \qquad \text{on }\mgr.
\]
as observed by considering the last component of the tropical Laplacian $\Deltatrop$ given in Equation~\eqref{eq:TropicalLaplacian}). In particular, it follows from the above two equations that
\[
d f_\mgr = \ld( \lf_\ell ) - \projhar (\ld( \lf_\ell )) = \ld( \lf_\ell )  - \sum_{n=1}^h (\rmM_l^{-1} P_\ell \, \ld( \lf_\ell ) )_n \, \gamma_n
\]
where $\projhar$ denotes the orthogonal projection from $C^1(G, \R)$ onto its subspace $\Omega^1(G, \R)$, for the inner product $\innone{\ell}{.\,,.}$, $\rmM_\ell$ is the period matrix, and $P_\ell$ is the matrix from Proposition~\ref{prop:OrthoProjection}. We assume further that we are using a fixed admissible basis relative to the layering $\pi$ on $G$ and use freely the notations of Section~\ref{sec:admissible_basis_layered_graphs}, e.g., $h_\pi^j$ will denote the genus of the $j$-th minor $\Gamma^j$.

\smallskip

The rest of the proof consists in a careful analysis of all above terms, relying on the results we will prove in Section~\ref{sec:GraphPeriodMatrices} for period matrices of degenerating graphs. We will thus refer to relevant parts of that section for the used results. 

\smallskip

We start by analyzing the column vector $\rmW_\ell :=  P_\ell  \,  \ld( \lf_\ell)$, see also the proof of Proposition~\ref{prop:exact_forms_limits}. Since $\ld( \lf_\ell)$ is exact on each layer $\pi_j$, the restricted vectors
\[
	\rmW_{\ell,j} := \rmW_\ell\rest{J_\pi^j}=\Big( \rmW_{\ell} (n) \Bigr)_{n \in J_\pi^j} \in \R^{h_\pi^j}
\]
have the following growth orders (here, by convention $L_{r+1}(\ell) = L_\fin(\ell)=1$ and $L_{\fin +1} = 0$)
\begin{equation} \label{eq:OrdersW}
\rmW_{\ell,j} = O\bigl(L_{j+1}(\ell) \bigr), \qquad j = 1,\dots, r, \fin.
\end{equation}
as $\mgr$ degenerates to $\curve$. 

\smallskip

Combining this with Theorem~\ref{thm:PeriodAsymptotics}$(i)$, which provides the precise asymptotics of the inverse of the metric graph period matrices, we can analyze the entries of the vector $\rmM_\ell^{-1} \rmW_\ell$ and conclude that, if the index $n$ belongs to the $j$-th set $J_\pi^j$ for some $j=1, \dots,r, \fin$, then
\[
\bigl(\rmM_\ell^{-1} \rmW_\ell\bigr)(n) = O\bigl(\max_{k \le j} \frac{L_{k+1}(\ell)}{L_k(\ell)}\bigr), \qquad {n \in J_\pi^j}.
\]

\smallskip

Assume now that the edge $e \in E$ belongs to the $j$-th layer $\pi_j$, $j \in [r]\cup\{\fin\}$. Since our fixed basis of $H_1(G, \Z)$ is admissible, we have 
\begin{align*}
d f_\mgr(e) &= \ld( \lf_\ell )(e) - \sum_{k \le j} \sum_{n \in J_\pi^k} \bigl(\rmM_\ell^{-1} \rmW_\ell\bigr)(n) \, \gamma_n(e) \\
		&= \ld( \lf )(e) -  \sum_{n \in J_\pi^j} \bigl(\rmM_\ell^{-1} \rmW_\ell\bigr)(n) \, \gamma_n (e) + O\bigl(\max_{k < j} \frac{L_{k+1}(\ell)}{L_k(\ell)}\bigr).
\end{align*}
On the other hand, it follows from Theorem~\ref{thm:PeriodAsymptotics}$(i)$ and \eqref{eq:OrdersW} that the $j$-th part of the vector $\rmM_\ell^{-1} \rmW_\ell$ is given by
\begin{align*}
\bigl( \rmM_\ell^{-1} \rmW_\ell \bigr)\rest{J_\pi^j} &:= \bigl((\rmM_\ell^{-1} \rmW_\ell) (n) \bigr)_{n \in J_\pi^j} \\ 
&= \sum_{k=1, \dots,r, \fin}  (\rmM_\ell^{-1})_{jk} \rmW_{\ell,k}  =  \sum_{k=j, \dots,r, \fin} (\rmM_\ell^{-1})_{jk} \rmW_{\ell,k} + O\bigl(\max_{k < j} \frac{L_{k+1}(\ell)}{L_k(\ell)}\bigr),
\end{align*}
where $(\rmM_\ell^{-1})_{jk}$ denotes the $(j,k)$-th block in the block decomposition of $\rmM_\ell^{-1}$ using admissible basis of $H_1(G)$ (see Section~\ref{sec:AsymptoticsGraphPeriods}). Moreover, Theorem~\ref{thm:PeriodAsymptotics}$(i)$ together with \eqref{eq:OrdersW} implies that
\[
(\rmM_\ell^{-1})_{jk} \rmW_{\ell, k}  = \sum_{p \in \mathcal{P}_{jk}} A_p(\mgr) \rmW_{\ell,k}  + \frac{L_{k+1}(\ell)}{L_j(\ell)} O\bigl( \max_{s \le k} \frac{L_{s+1}(\ell)}{L_s(\ell)} \bigr )
\]
for every $k \in \{j, \dots, r, \fin\}$. Here, we are using the notations of Section~\ref{sec:explicit_extension_form}. Comparing the matrices $A_p(\mgr)$ and $A_p(\curve_l)$ (see \eqref{eq:AsmyptoticsAp} and \eqref{eq:AsmyptoticsApTropical}), we see that for every strictly increasing sequence $p \in \mathcal{P}_{jk}$, we have
\[
A_p(\mgr) \rmW_{\ell, k}  = \frac{1}{L_j(\ell)} A_p(\curve_\ell) \rmW_{\ell, k}  +  \frac{L_{k+1}(\ell)}{L_j(\ell)} O\bigl( \max_{s \le k} \frac{L_{s+1}(\ell)}{L_s(\ell)} \bigr )
\]
(up to the desired orders, the $(\rmM(\mg))_{mn}$'s in the product can be replaced by the $T_{mn}$'s). The claimed expansion now follows from the sum-product formula established in Proposition~\ref{prop:harmonic_extension_matrices}.
\end{proof}

\subsubsection{Proof of Theorem~\ref{thm:layered_expansion_tropical_laplacian}} The asymptotic expansion \eqref{eq:ExpansionConsideringLayered} follows by integrating the local expansion given in  Lemma~\ref{lem:DerivativeAsymptotics} along paths of the metric graph $\mgr$, starting from the fixed vertex $v$. Note that all functions $f_{\ell,j}^\ast$ take value zero at $v$ by \eqref{eq:JDFctTrop}. To conclude, it only remains to notice that the product $\ell(e) R_e(\ell)$ has order $\max_{k} {L_{k+1}^2(\ell)}/{L_k(\ell)}$ for every edge $e \in E$.

\smallskip

The uniform convergence of the functions $f_{\ell, j}$ to $f_j$ is obtained by combining Proposition~\ref{thm:ContForms}, Local Expansion Lemma~\ref{lem:DerivativeAsymptotics} and the expression in Proposition~\ref{prop:harmonic_extension_matrices}. 
\qed

\subsubsection{Stratumwise tameness in the discrete case}

Combining the local asymptotics from Theorem~\ref{thm:layered_expansion_tropical_laplacian} with the tame topology on $\mgtropcombin{\grind{G}}$, we arrive at the result on the stratumwise tameness.

\begin{proof}[Proof of Theorem~\ref{thm:MainGraphDiscrete}]
By definition, we must prove the following. Suppose that $\mgtropcombin{\combind{(G, \pi)}}$ is a boundary stratum in $\mgtropcombin{\grind{G}}$. If $(t_n)_n$ is a sequence in $\tilde \mg_{\grind{G}}$ (representing metric graphs $\mgr_n$) such that $\lim_{n \to \infty} t_n = \thy$ in the tame topology on $\mgtropcombin{\grind{G}}$ for some point $\thy$ in the closure of  $\mgtropcombin{\combind{(G, \pi)}}$ (representing a tropical curve $\curve$ of type $(G,  \overline{\pi})$, for an ordered partition $\overline \pi \in \Piall(E)$), then
\begin{equation} \label{eq:StatementProofMainGraphDiscrete}
\lim_{n \to \infty} \big \| f_{t_n} - \lf_{\pr_\pi (t_n)}^\ast  \big \|_{\infty, \mgr_n} = 0,
\end{equation}
where $\pr_\pi \colon \tilde \mg_{\grind{G}} \to \mgtropcombin{\combind{(G, \pi)}}$ is the projection map defined in Section~\ref{ss:TameFunctions}. 

Suppose first that the limit point $\thy$ belongs to $\mgtropcombin{\combind{(G, \pi)}}$, that is, its tropical curve $\curve$ has the same partition $\bar \pi = \pi$. Then, \eqref{eq:StatementProofMainGraphDiscrete} is precisely the statement of Theorem~\ref{thm:layered_expansion_tropical_laplacian} (since, in terms of the notation there, 
$\lf_{\ell_n} = \lf_{\pr_\pi (t_n)}$). Notice also that the error term in \eqref{eq:ExpansionConsideringLayered} goes to zero, since we require that $t_n \to \thy$ in the tame topology of $\mgtropcombin{\grind{G}}$.

In the general case, the ordered partition $\bar \pi$ of $\curve$ has to be a refinement of $\pi$, that is, $\pi \preceq \bar \pi$. In this case, \eqref{eq:StatementProofMainGraphDiscrete} is obtained inductively, using the result for the lower-dimensional spaces $\tilde \mg_{\grm{\pi}{j}(G)}^\trop$ of the graded minors $\grm{\pi}{j}(G)$ of $(G, \pi)$. We omit the details.
\end{proof}

\subsection{Asymptotic description of the $\jvide$-function} \label{ss:AsymptoticsJFunc}
We use the results of the previous section to provide an asymptotic description of the $\jvide$-function on degenerating metric graphs (see \eqref{eq:JFunctionGraph} for definitions).

 \smallskip

Let $\mgr$ be a metric graph over $G$ with edge length function $\ell \colon E \to (0, +\infty)$.  If $p, q $ are points on $\mgr$, then the $j$-function of $\mgr$ gives the solution $y \in \mgr \mapsto \jfunc{p\tiret q, v}(y)$ of the equations
 \begin{equation}
\begin{cases}
\Deltatrop_\mgr \, f  = {\delta}_p -{\delta}_q \\
f (x) = 0
 \end{cases}
\end{equation}
(see Section~\ref{ss:LaplacianMetricGraph} for details).

\smallskip

We introduce the following generalization to the tropical setting. Suppose that $\curve$ is a tropical curve with underlying graph $G=(V,E)$, ordered partition $\pi = (\pi_1, \dots, \pi_r, \pi_\fin)$ and edge length function $l$, normalized on each infinitary layer $\pi_j$, $j=1, \dots, r$. Let $\Gamma$ be the canonical representative in the conformal equivalence class of $\curve$. Then we can identify every point $p \in \Gamma$ naturally with a (layered) Dirac measure $\bm{\delta}_p$ on $\curve$. In particular, for three points $p,q,x \in \Gamma$, the following system of equations 
\begin{equation} \label{eq:LayeredJFunction}
\begin{cases}
\Deltatrop_\curve \, \lf = \bm{\delta}_p - \bm{\delta}_q \\
 \lf \text{ is harmonically arranged} \\
 \int_{\curve} \lf \, d \bm{\delta}_x = 0
 \end{cases}
\end{equation}
has a unique solution $\lf = (f_1, \dots, f_r, f_\fin)$ on $\curve$. We denote the latter by $\ljfuncbis{ p \tiret q, x}$. After pulling back the components $\ljfuncbis{ p \tiret q, x, j} = f_j$ from the graded minors $\Gamma_j$, $j=1, \dots, r, \fin$ to $\Gamma$, we can also view the latter as a function $y \in \Gamma \mapsto (\ljfuncbis{ p \tiret q, x, j}(y))_j \in \R^{r+1}$. Altogether, this gives a map 
\[
(p,q,x,y) \in \Gamma^4 \quad \mapsto \quad (\ljfuncbis{ p \tiret q, x, j}(y))_{j}  \in \R^{r+1},
\]
which we call the \emph{$j$-function of the tropical curve $\curve$}. In Theorem~\ref{thm:TropicalHPvsJFunction} we will prove that, analogous to the case of metric graphs \eqref{eq:MGHPvsJ},  the $\jvide$-function on tropical curves can be expressed in terms of the height pairing on tropical curves.

\smallskip

Consider now a combinatorial graph $G = (V,E)$ and a fixed boundary stratum $ \mgtropcombin{\combind{(G, \pi)}}$ (associated to an ordered partition $\pi$) in $\mgtropcombin{\grind{G}}$. Recall that the projection map $\pr_\pi \colon \tilde \mg_{\grind{G}} \to \mgtropcombin{\combind{(G, \pi)}}$ assigns to every metric graph $\mgr$ over $G$ with length function $\ell \colon E \to (0, + \infty)$ the tropical curve $\curve_\ell =(G, \pi, \ell)$ of combinatorial type $(G, \pi)$, whose edge length function $\ell_j$ given by the (normalization of the) restriction $\ell\rest{\pi_j}$ on each $\pi_j$.
By using homothecies between edges, any three points $p,q,x \in \mgr$ can be naturally viewed as points on the canonical representative $\Gamma = \Gamma(\curve_\ell)$ of $\curve_\ell$, with normalized edge length function on each layer. We denote by $\ljfunc{\ell, p \tiret q, x}^\ast$ the pull-back of the solution to \eqref{eq:LayeredJFunction} on $\curve_\ell$ to $\mgr$.

\begin{lem} \label{lem:JfunctionTame}
Let $\mgtropcombin{\combind{(G, \pi)}}$ be a boundary stratum associated to an ordered partition $\pi$ of $E$. Suppose $\curve$ is a tropical curve of type $(G, \overline{\pi})$ with $\pi \preceq \overline{\pi}$.  As the metric graph $\mgr$ degenerates to $\curve$  in the tame topology on $\mgtropcombin{\grind{G}}$, \begin{equation} \label{eq:JFunctionTameUniform}
	\sup_{p,q,x,y  \in \mgr} \big | \jfunci{p \tiret q, x} (y)-  \ljfunc{\ell, p \tiret q, x}^\ast(y) \big | = o(1).
\end{equation}
\end{lem}
\begin{proof}
It suffices to treat the case that $x = v$ is a vertex of $G$, since $\jfunci{p \tiret q, x} (y) = \jfunci{p \tiret q, v} (y)  - \jfunci{p \tiret q, v} (x)$ and
\[
\ljfunc{\ell, p \tiret q, x}^\ast(y) = \ljfunc{\ell, p \tiret q, v}^\ast(y) - \ljfunc{\ell, p \tiret q, v}^\ast(y)
\]
by additivity. Moreover, since trivially $\jfunci{p \tiret q, v} (y) = \jfunci{p \tiret r, v} (y)  -  \jfunci{ q \tiret r, v} (y)$, we can also assume in \eqref{eq:JFunctionTameUniform} that $p = u$ is a vertex of $G$. Moreover, if $q \in \mgr$ belongs to an edge $e = ab$ (and hence can be seen as a number in the interval $[0, \ell(e)]$), then one easily verifies that
\begin{equation} \label{eq:JFunctionMiddleEdges}
\jfunci{u \tiret q, v} (y) = \big ( 1 - \frac{q}{\ell(e)} \big ) \jfunci{u \tiret a, v} (y) + \frac{q}{\ell(e)} \jfunci{u \tiret b, v} (y) + h_e(y),\qquad y \in \mgr,
\end{equation} where the function $h_e \colon \mgr \to \R$ is given by (see also \cite[Lemma 4]{Faltings21})
\[
h_e(y) = \begin{cases} 0, & y \notin e \\
- \min \Big \{ \big ( 1 - \frac{q}{\ell(e)} \big ) y,  \frac{q}{\ell(e)} (\ell(e) - y) \Big \}, & y \in e \cong [0, \ell(e)] \end{cases}
\]
It follows that in proving \eqref{eq:JFunctionTameUniform}, the supremum over all points $p, q \in \mgr$ can be replaced by the maximum over all pairs of vertices $v_1, v_2 \in V$. However, for each divisor $v_1-v_2 \in \Div^0(G)$, we can apply Theorem~\ref{thm:MainGraphDiscrete} and hence the proof is complete.
\end{proof}

Let $G = (V, E)$ be a fixed finite graph. Suppose that $(\mgr, \mu,  \nu)$ is a measured metric graph over $G$.  By \eqref{eq:SolutionFormulaJFunction}, the solution $f \colon \mgr \to \C$ to the Poisson equation \eqref{eq:PoissonMetricGraphs} can be written as 
\begin{equation} \label{eq:SolutionFormulaJFunctionGeneral}
f(y) = \int_\mgr  \jfunci{p \tiret q, x}(y) \, d \mu(p)  - \Big ( \int_\mgr \int_\mgr  \jfunci{p \tiret q, x} (r) \, d \mu(p) \, d \nu(r) \Big ) , \qquad y \in \mgr,
\end{equation}
where $q$, $x \in \mgr$ are arbitrary fixed points.
 Thus, in order to prove Theorem~\ref{thm:main_graph_strong} and Theorem~\ref{thm:MainGraphWeak}, it only remains to integrate the asymptotics for the $j$-function from above. We will formalize this in the following sections.


\subsection{Proof of Proposition~\ref{prop:leadingterm}}  \label{sec:proof_proposition_leading_term}
In this section, we present the proof of Proposition~\ref{prop:leadingterm}, which actually follows from the basic estimates in Proposition~\ref{thm:ContForms}.

\begin{proof}[Proof of Proposition~\ref{prop:leadingterm}]
Let $(t_n)_n$ be a sequence in $\mgtropcombin{\grind{G}}$ which converges in $\mgtropcombin{\grind{G}}$ to some point $\thy \in \mgtropcombin{\grind{G}}$.  For simplicity, assume first that all $t_n$'s belong to $\mggraphcombin{\grind{G}}$ (and hence represent metric graphs $\mgr_n$) and that $\thy$ belongs to some boundary stratum $\mgtropcombin{\combind{(G, \pi)}}$ (and hence represents a tropical curve $\curve =(G, \pi, l)$).

Fix a vertex $v$ of $G=(V,E)$. For a point $p \in \mgr_n$, we will denote by $\overline{p} \in \Gamma^1$ the image of $p$ in the first graded minor $\Gamma^1$ of $\curve$. Let $\jfunci{p \tiret v, v}(\cdot)$ be the $j$-function on $\mgr_n$ and $\jfunci{\Gamma^1, \overline{p} \tiret \overline{v}, \overline{v}}$ the $j$-function on $\Gamma^1$.  From Proposition~\ref{prop:leadingterm}, it follows that 
\begin{equation} \label{eq:BasicJApproximation}
\lim_{n \to \infty} \sup_{p, y \in \mgr_n} |L(t)^{-1} j_{p \tiret v, v}(y) - j_{\Gamma^1, \overline{p} \tiret \overline{v}, \overline{v}}(\overline{y})| = 0.
\end{equation}
The uniformity in $p$ follows as in the proof of Lemma~\ref{lem:JfunctionTame} from \eqref{eq:JFunctionMiddleEdges}. 

Since the measure sequences  $(\mu_n)_n$ and $(\nu_n)_n$ are weakly convergent, their total variations are uniformly bounded,  $\sup_n |\mu_n| + |\nu_n| < \infty$. Combining \eqref{eq:BasicJApproximation} with   the solution formula \eqref{eq:SolutionFormulaJFunctionGeneral}, it follows that
\[
\lim_{n \to \infty} \sup_{y \in \mgr_n} |h_{t_n}(y) - g_{t_n}^\ast(y)| = 0,
\]
where $g_{t_n}^\ast$ is the pull-back of the solution $g_{t_n} \colon \Gamma^1 \to \C$ on the first minor $\Gamma^1$ to
\[
	\Delta_{\Gamma^1}( g_{t_n}) = \tilde \mu_{t_n}  \qquad \text{and} \qquad \int_{\Gamma^1} g_{t_n} \, d \tilde \nu_{t_n}=0,
\]
and $\tilde \mu_{t_n}$, $\tilde \nu_{t_n}$ are the natural push-outs of $\mu_{t_n}$, $\nu_{t_n}$ from $\mgr_n$ to $\Gamma^1$. However, since $\tilde \mu_{t_n}$ and $\tilde \nu_{t_n}$ converge weakly to $\mu_{\thy}$ and $\nu_{\thy}$, the solutions $g_{t_n}$ converge uniformly to $f_{1, \thy}$ on $\Gamma^1$ (this follows again from \eqref{eq:SolutionFormulaJFunctionGeneral}, see also \cite[Proposition 3.11]{BR10}).

Hence the claim holds true if all $t_n$'s belong to the open part $\mggraphcombin{\grind{G}}$ of $\mgtropcombin{\grind{G}}$. In general, the convergent sequence $(t_n)_n$ can belong to some boundary stratum $\mgtropcombin{\combind{(G, \pi)}}$. However,  using the second claim in Proposition~\ref{thm:ContForms}, this case can be treated analogously.
\end{proof}

\subsection{Integral representation of tropical solutions using the tropical $\jvide$-function}
As an auxiliary result, we next extend the representation \eqref{eq:SolutionFormulaJFunctionGeneral} to tropical curves (for simplicity, we only treat the case that $q =x = v$ is a vertex in \eqref{eq:SolutionFormulaJFunctionGeneral}).

\begin{lem}\label{lem:SolutionFormulaJFunctionTropical}
Let $(\curve, \lmu, \lnu)$ be a measured tropical curve of rank $r$. Fix  a vertex $v \in V(G)$ and for each point $p \in \Gamma$, denote by $\lz_{p} = (\zeta_{p,1}, \dots, \zeta_{p,r}, \zeta_{p,\fin} )$ the solution of the Poisson equation \eqref{eq:LayeredJFunction} for $q=x= v$, that is, $\lz_p = \ljfuncbis{p\tiret v, v}$.

Then the $j$-th piece $f_j \colon \Gamma^j \to \C $ of the solution $\lf$ to \eqref{eq:PoissonTropicalCurve} is given by
\begin{equation} \label{eq:SolutionFormulaJFunctionTropical}
f_j (y) = \int_{\Gamma} \zeta_{p,j}(y) \, d\lmu (p) - \int_{\Gamma} \int_{\Gamma} \zeta_{p,j}(r) \, d\lmu (p) d\lnu (r), \qquad y \in \Gamma_{j}, \, j=1, \dots, r, \fin,
\end{equation}
where $\lmu$ and $\lnu$ are seen as measures on the space $\Gamma$ and $\zeta_{p,j}$ as a function on $\Gamma$.
\end{lem}
\begin{proof}
We only have to prove the function $\lf = (f_1,  \dots, f_r, f_\fin)$, with $f_j \colon \Gamma^j \to \C$ defined by Equation~\eqref{eq:SolutionFormulaJFunctionTropical} above, satisfies the conditions in the Poisson equation~\eqref{eq:PoissonTropicalCurve}. The last two properties in \eqref{eq:PoissonTropicalCurve} follow immediately from \eqref{eq:SolutionFormulaJFunctionTropical} and the fact that $\lz_{p}$ is harmonically arranged on $\curve$ for all $p$. In order to prove that $\Deltatrop(\lf) = \lmu$, we must show that
\[
	\int_{\Gamma^{j}} f_j \, \Delta_{\Gamma^{j}}(\varphi) = \int_{\Gamma^{j}} \varphi \, d \mu_j - \sum_{i < j} \int_{\Gamma^{j}} \varphi \, \divind{j}{i}(f_i)
\]
for every $j=1, \dots, r, \fin$ and every continuous, edgewise smooth function $\varphi\colon \Gamma^{j} \to \R$. For $j=1$, the Poisson equation~\eqref{eq:LayeredJFunction} and a change of order in the integration yields 
\begin{align*}
\int_{\Gamma^{1}} f_1\, \Delta_{\Gamma^{1}}(\varphi) &= \int_{\Gamma} \Big ( \int_{\Gamma^{1}} \zeta_{p,1}(s)  \, \Delta_{\Gamma^{1}}(\varphi)(s) \Big ) d \lmu (p) = \int_{\Gamma}  \Big ( \int_{\Gamma^{1}} \varphi(s) \, (\delta_{\bar p} - \delta_{\bar v}) (s) \Big ) d \lmu (p) \\ &=  \int_{\Gamma}   \varphi(\bar p) \, d \lmu (p) =  \int_{\Gamma^1}   \varphi \, d \mu_1,
\end{align*}
where, we recall, for any $p \in \Gamma$, $\bar p$ denotes its natural image in the first minor $\Gamma^1$ under contraction.
Using the exact same steps (take into account the shape of the layered measure $\bm{\delta}_p - \bm{\delta}_q$ in \eqref{eq:LayeredJFunction}), one verifies that for $j >1$ and a test function $\varphi \colon \Gamma^{j} \to \R$,
\begin{align*}
\int_{\Gamma^{j}} f_j \, \Delta_{\Gamma_{j}}(\varphi) &= \int_{\Gamma^{j}} \varphi \, d \mu_j   -  \sum_{i < j} \int_{\Gamma^{j}} \varphi \, d \mu_{i<j},
\end{align*}
where $\mu_{i<j}$ is the discrete measure on $\Gamma^j$ associated with the degree-zero divisor 
\[
 \int_{\Gamma} \divind{i}{j}(\zeta_{p,i}) \, d \mu(p) = \sum_{v \in V(\Gamma^j)} a_v \cdot \delta_v
\]
supported on vertices of $\Gamma^j$, and obtained by integrating the coefficients of $\divind{i}{j}(\kappa_{p,i})$ with respect to the measure $d \mu(p)$. However, the coefficient at a vertex $v$ of the $j$-th minor $\Gamma^j$ is precisely
\[
a_v = \int_{\Gamma} \sum_{\substack{u \in V(G) \\ \zeta_j(u)  = v}} \sum_{\substack{e \in \pi_i \\ e \sim v}} \slp_e(\zeta_{p,i})(u) \, d\mu(p) = \sum_{\substack{u \in V(G) \\ \proj_j(u)  = v}} \sum_{\substack{e \in \pi_i \\ e \sim v}} \slp_e(f_i)(u)  =  \divind{i}{j}(f_i)(v).
\]
This implies that $\mu_{i<j} =  \divind{i}{j}(f_i)$ for all $i<j$, and completes the proof.
\end{proof}

\subsection{One more auxiliary lemma} We will deduce both theorems  Theorem~\ref{thm:main_graph_strong} and Theorem~\ref{thm:MainGraphWeak} from a slightly technical lemma proved in this section. 

\smallskip

Given a measured metric graph $(\mgr, \mu, \nu)$, we denote by $f_\mgr$ the solution to the Poisson equation \eqref{eq:PoissonMetricGraphs}.

As before, consider a fixed boundary stratum $\mgtropcombin{\combind{(G, \pi)}}$ associated to an ordered partition $\pi$. The projection map $\pr_\pi \colon \tilde \mg_{\grind{G}} \to \mgtropcombin{\combind{(G, \pi)}}$ assigns to every metric graph $\mgr$ over $G$ with length function $\ell \colon E \to (0, + \infty)$ the tropical curve $\curve_\ell$ of combinatorial type $(G, \pi)$ whose edge length function is given by the point $\pr_\pi(\ell) \in \mgtropcombin{\grind{G}}$. Using the natural maps between tropical curves over the same graph $G$, we can push-out $\mu$ and $\nu$ on $\mgr$ to two layered measures $\lmu$ and $\lnu$ on $\curve_\ell$. Clearly, $(\curve_\ell, \lmu, \lnu)$ is a measured tropical curve. For a Borel measure $\omega$ on a metric graph $\mgr$, we denote by $|\omega| (\mgr)$ the total variation of $\omega$.

We denote as before by $\lf_\ell = (f_{\ell, 1}, \dots, f_{\ell, r}, f_{\ell, \fin})$ the unique  function on $\curve_\ell$ to equation \eqref{eq:PoissonTropicalCurve}. Moreover, $\lf_\ell^\ast = L_1(\ell) f_{\ell,1}^\ast + L_2(\ell) f_{\ell,2}^\ast +\dots + f_{\ell,\fin}^\ast$ is the pull-back of $\lf_\ell$ to $\mgr$. Note that $\curve_\ell$ is of type $(G, \pi)$.

\begin{lem} \label{lem:TechnicalLemmaGeneralMeasures}
Notations as above, suppose that for each metric graph $\mgr$ with edge length function $\ell \colon E \to (0, + \infty)$, we have two Borel measures $\mu_\mgr$ and $\nu_\mgr$ satisfying \eqref{eq:AssGraphMeasures}, and denote by $\curve_\ell$ the tropical curve of type $(G, \pi)$ associated to $\mgr$. Let $\curve$ be a tropical curve of type $(G, \overline{\pi})$ for some $\pi \preceq \overline{\pi}$. As the metric graph $\mgr$ degenerates to $\curve$ in the tame topology on $\mgtrop{\grind{G}}$, we get
\[
	\sup_{y  \in \mgr} | f_\mgr (y)-  \lf_\ell^\ast(y)| = \Big (|\mu_\mgr| (\mgr)  + |\mu_\mgr| (\mgr) \cdot |\nu_\mgr| (\mgr)  \Big )  o(1).
\]
\end{lem}
\begin{proof}
Fix an arbitrary vertex $v$ of $G$. Combining Lemma~\ref{lem:JfunctionTame} with the representation \eqref{eq:SolutionFormulaJFunctionGeneral}, we conclude that
\[ \| f_{\mgr} - (L_1(\ell)  h_{\ell,1} + \dots + h_{l, \fin}) \|_\infty = \Big (|\mu_\mgr| (\mgr)  + |\mu_\mgr| (\mgr) \cdot |\nu_\mgr| (\mgr)  \Big )  o(1)
\]
as the metric graph $\mgr$ degenerates tamely to $\curve$, where the functions $h_{\ell,j} \colon \mgr \to \C$, $j=1, \dots, r,\fin$ are given by
\[
h_{\ell,j} (y) = \int_{\mgr} \zeta_{p,j}^\ast(y) \, d\mu_\mgr (p) - \int_{\mgr} \int_{\mgr} \zeta_{p,j}^\ast(r) \, d\mu_\mgr (p) d\nu_\mgr(r) \qquad y \in \mgr,
\]
and we abbreviate notation by writing $\lz_{p} = (\zeta_{p,1}, \dots, \zeta_{p,\fin} )$ for the function $\lz_{p} = \ljfunc{\ell, p \tiret v, v} $ on $\curve_\ell$ from Lemma~\ref{lem:JfunctionTame}. In particular, the error term in Lemma~\ref{lem:JfunctionTame} results in the $o(1)$ term on the RHS, since $\mgr$ degenerates to $\curve$ in the tame topology on $\mgtropcombin{\grind{G}}$. It remains to notice that, by Lemma~\ref{lem:SolutionFormulaJFunctionTropical}, the $j$-th function $h^\ast_{\ell,j}$ is precisely the pull-back of $f_{\ell,j}$ to the metric graph $\mgr$, that is, $h_{\ell,j} = f_{l,j}^\ast$.
\end{proof}

\subsection{Proof of Theorem~\ref{thm:main_graph_strong}} We have now all the ingredients to prove the stratumwise tameness of the solutions. 

\smallskip

Fix a boundary stratum ${\mgtropcombin{\combind{(G, \pi)}}}$ in $\mgtropcombin{\grind{G}}$ and let $\thy$ be a point in the closure of the stratum $\tilde{\mgtrop{\combind{(G, \overline{\pi})}}}$ in $\mgtropcombin{\grind{G}}$ (representing a tropical curve $\curve$). Suppose that $(t_n)_n$ is a sequence in $\mggraph{\grind{G}}$ (representing metric graphs $\mgr_n$) with $\lim_{n \to \infty} t_n = \thy$ in the tame topology on $\mgtropcombin{\grind{G}}$. We have to prove that
\[
\lim_{n \to \infty} \| f_{t_n} - \lf_{\pr_\pi(t_n)}^\ast\|_{\infty, \mgr_n} = 0.
\]
Denote by $\tilde f_{t_n}$ the solution on $\mgr$ to the Poisson equation for the measures $\tilde \lmu_{t_n} :=  \lmu_{\pr_\pi(t_n)}^\ast$ and $\tilde \lnu_{t_n} := \lnu_{\pr_\pi(t_n)}^\ast$, the pull-backs of $\lmu_{\pr_\pi(t_n)}$ and $\lnu_{\pr_\pi(t_n)}$ from the fiber $\unicurvetrop{\pr_\pi(t_n)}$ to $\mgr_n$. Since $\lim_n \|  \tilde f_{t_n} - \lf_{\pr_\pi(t_n)}^\ast\|_\infty = 0$ by Lemma~\ref{lem:TechnicalLemmaGeneralMeasures}, it remains to show that $\lim_n \| f_{t_n} - \tilde f_{t_n}\|_\infty = 0$.

Fix a vertex $v \in V(G)$ and let $\zeta_p := \jfunci{p\tiret v, v}(\cdot)$ be the corresponding $j$-function on $\mgr_n$. Using the solution formula \eqref{eq:SolutionFormulaJFunctionGeneral}, we conclude that, for any $y \in \mgr_n$,
\begin{equation} \label{eq:SolutionsDifferenceStratumwiseTame}
f_{t_n}(y) - \tilde f_{t_n}(y) = g_{t_n}(y) + \int_\mgr \int_\mgr \zeta_p(r) \, d \tilde \lmu_{t_n}(p) d \tilde \lnu_{t_n}(r) - \int_\mgr \int_\mgr \zeta_p(r) \, d  \lmu_{t_n}(p) d \lnu_{t_n}(r),
\end{equation}
where $g_{t_n} \colon \mgr_n \to \C$ is the solution to the Poisson equation $\Delta g = \lmu_{t_n} - \tilde \lmu_{t_n}$ on $\mgr_n$ with $g(v) = 0$. Due to our assumption \eqref{eq:condition_strong_tame}, we can apply Proposition~\ref{prop:leadingterm} to $g_{t_n}$ and infer that $ \lim_n \| g_{t_n} \|_{\infty, \mgr_n} = 0$. Applying the solution formula \eqref{eq:SolutionFormulaJFunctionGeneral} once again, we see that the remaining term in \eqref{eq:SolutionsDifferenceStratumwiseTame} can be written as
\[
R_n(y) = - \int_\mgr g_{t_n}(r) \, d \tilde \lnu_{t_n}(r) + \int_\mgr h_{t_n}(r) \, d ( \tilde \lnu_{t_n} - \lnu_{t_n}) (r), \qquad y \in \mgr_n,
\]
where $h_{t_n} \colon \mgr_n \to \C$ is the solution to $\Delta h = \mu_{t_n}$ with $h_{t_n} (v) = 0$. Due to our assumption \eqref{eq:condition_strong_tame}, we can apply Proposition~\ref{prop:leadingterm} to $h_{t_n}$ and conclude that the second term goes to zero. Moreover, since the measures $\lnu_{t_n}$ converge weakly to $\lnu_\thy$, their total variations are uniformly bounded. In particular, since $\lim_n \| g_n\|_\infty = 0$, we get that $\sup_{y \in \mgr_n} | R_n(y) |  \to 0$ for $n \to \infty$. Hence the claim is proven.
\qed

\subsection{Proof of Theorem~\ref{thm:MainGraphWeak}}  In what follows, we use the same notation as in Section~\ref{ss:TropicalWeakTame} and consider the tropical log map
\[
 \begin{array}{cccc}
 \Lognoind \colon  &\mgtropcombin{\grind{G}} \setminus {\umggraphcombin{\grind{G}}} &\to &  \partial_\infty \mgtropcombin{\grind{G}}
  \end{array}
\]
to the boundary at infinity $\partial_\infty \mgtropcombin{\grind{G}}  =  \mgtropcombin{\grind{G}} \setminus  \mggraphcombin{\grind{G}}$. Recall that its domain decomposes into
\[
\mgtropcombin{\grind{G}} \setminus {\umggraphcombin{\grind{G}}} = \bigsqcup_{\substack{\pi = (\pi_\infty, \pi_\fin) \in \Piall(E) \\ \pi_\infty \neq \varnothing}  } \inn R_\pi, \qquad \qquad \inn R_\pi := \Lognoind^{-1}\Big ( \mgtropcombin{\combind{(G, \pi)}} \Big ).
\]
Let $\lf_\thy$, $\thy \in  \mgtropcombin{\grind{G}}$, be the family of tropical solutions from Theorem~\ref{thm:MainGraphWeak}

\begin{proof}[Proof of Theorem~\ref{thm:MainGraphWeak}]
In the following, we define tropical functions $\layh_{\thy}$, $\thy \in  \mgtropcombin{\grind{G}} \setminus {\umggraphcombin{\grind{G}}}$, with the properties in Definition~\ref{def:WeakGlobalTameTropical}. We specify the functions $\layh_\thy$ separately on each region $\inn R_\pi$ associated to an ordered partition $\pi = (\pi_\infty, \pi_\fin)$ of $E$ with $\pi_\infty = (\pi_1, \dots, \pi_r) \neq \varnothing$. For $t \in \inn R_\pi \cap \mggraphcombin{\grind{G}}$, consider the two points $\thy := \Lognoind(t)$ and $\thy' := \pr_\pi(t)$ in the boundary stratum $\mgtropcombin{\combind{(G, \pi)}}$. We pushout the measures  $\lmu_t$ and $\lnu_t$ from the metric graph $\unicurvetrop{t}$ to (layered) measures $\lmu_t'$ and $\lnu_t'$ on the tropical curve $\unicurvetrop{\thy'}$ (by using homotecies on intervals), and consider the tropical Poisson equation for $(\unicurvetrop{\thy'},\lmu_t', \lnu_t')$. The tropical function $\layh_t = (h_{t,1}, \dots, h_{t,r}, h_{t, \fin})$ on $\unicurvetrop{\thy}$ is defined as the pullback of the solution from $\unicurvetrop{\thy'}$ to $\unicurvetrop{\thy}$ (that is, all components $h_{t, j}'$ are pulled back by using homotecies). For $\thy$ belonging to $\inn R_\pi \setminus \mggraphcombin{\grind{G}} = \mgtropcombin{\combind{(G, \pi)}}$, we set $\layh_\thy := \lf_\thy$. 

We proceed to verify the properties in Definition~\ref{def:WeakGlobalTameTropical}. The dependence $\thy \mapsto h_{\thy, j}$ for $j=1, \dots, r, \fin$ is continuous on $\inn R_\pi$. This follows from the continuity of the maps $\Lognoind$ and $\pr_\pi$, and the measures $\lmu_\thy$ and $\lnu_\thy$, combined with the following property: if $(\curve_n, \lmu_n, \lnu_n)$ is a sequence of measured tropical curves of type $(G, \pi)$ which converges to a measured tropical curve $(\curve, \lmu, \lnu)$ of the same type $(G, \pi)$, meaning that $\lim_n \thy_n = \thy$ for corresponding points in $\mgtropcombin{\combind{(G, \pi)}}$ and $\lim_n \lmu_n = \lmu$, $\lim_n \lnu_n = \lnu$ weakly, then the corresponding solutions to the Poisson equation satisfy $\lim_n \| \layg_n - \layg\|_\infty = 0$. The latter claim is analogous to Lemma~\ref{lem:TechnicalLemmaGeneralMeasures}, and can be deduced by using the convergence of the tropical $j$-functions, namely the components $\zeta_{p,j}$ in Lemma~\ref{lem:SolutionFormulaJFunctionTropical}, obtained from Lemma~\ref{thm:ContForms}, and then, using the integration formula \eqref{eq:SolutionFormulaJFunctionTropical}.

Finally we prove that, if a sequence $(t_n)_n$ in $\mggraphcombin{\grind{G}}$ converges tamely to a point $\shy \in \partial_\infty \mgtropcombin{\grind{G}}$, then
\[
\lim_{n \to \infty} \sup_{x \in \unicurvetrop{t_n}}  \big | f_{t_n}(x)  - \layh_{t_n}^\ast(x)  \big |  = 0.
\]
Without loss of generality, we may assume that all $t_n$'s belong to the same region $\inn R_\pi$ for some ordered partition $\pi$. However, then the claim follows directly from Lemma~\ref{lem:TechnicalLemmaGeneralMeasures} (since the measures are weakly convergent, their total variations remain uniformly bounded).
\end{proof}

\section{Tropical Green functions} \label{sec:TropicalGreenFunctions}
In this section, we introduce {\em tropical Green functions} and employ them in describing the asymptotic of Green functions on degenerating metric graphs.

\subsection{Metric graph Green functions} We start by recalling basic definitions in the setting of metric graphs. Let $\mgr$ be a metric graph and $\mu \in \widetilde{\mathcal{M}}(\mgr)$ a measure of mass $\mu(\mgr) = 1$. The {\em Green function} associated to $\mu$ is the unique function $\grg_{\mu} \colon \mgr \times \mgr \to \R$ such that for each point $x \in \mgr$, the function $\grg_{\mu}(x, \cdot)$ belongs to $D(\Delta_\mgr)$ and satisfies
\begin{equation} \label{eq:DefGFGraph}
\Delta \grg_\mu (x, \cdot)=  \delta_x - \mu, \qquad \int_\mgr  \grg_\mu (x, y) \, d\mu(y) = 0.
\end{equation}
 If there is no risk of confusion, we sometimes abbreviate and write $\grg(x,y)$ instead of $ \grg_\mu(x,y)$.

\smallskip

As follows from \eqref{eq:SolutionFormulaJFunction}, $ \grg_\mu(x,y)$ has the following representation using the $j$-function of $\mgr$,
\begin{equation} \label{eq:GreenByJFunction}
 \grg_\mu(x,y) = \int_\mgr \jfunci{x\tiret q,z}(y) \, d\mu (q)  - \int_\mgr \int_\mgr \jfunci{x\tiret q,z}(r)  \, d\mu(q) \, d\mu (r), \qquad x,y \in \mgr
\end{equation}
which holds for all points $z \in \mgr$.

\begin{remark} 
The regularity assumption $\mu \in \widetilde{\mathcal{M}}(\mgr)$ can be omitted. In fact, our results hold true for general, non-negative Borel measures $\mu$ with $\mu(\mgr) = 1$ (one only has to allow functions in the class $\operatorname{BDV}(\mgr)$ in order to define the Green function, see Remark~\ref{rem:DistributionalLaplacian}). 
\end{remark}
\smallskip

\begin{remark}[Spectral theory on weighted metric graphs] \label{rem:SpectralTheory}
 The measure-valued Laplacian $\Delta$ gives rise to a self-adjoint operator in the $L^2$-space $L^2(\mgr, \mu)$, whose domain of definition and action are given by (for simplicity, we assume that $\mu$ is a nonnegative measure of the form \eqref{eq:SpecialMus} with total mass $\mu(\mgr) = 1$)
\begin{align*} \label{eq:DefSaOperator}
\dom (\bH) &= \bigl \{f \in H^1(\mgr) \st \, \Delta f = g \cdot d \mu \text{ for some } g \in L^2(\mgr, \mu) \bigr \}   \\
 \bH f &:=  g.
\end{align*}
In the above, $H^1(\mgr)$ denotes the first order Sobolev space on $\mgr$,
\[
H^1(\mgr) = \Bigl\{ \, f \in \mathcal{C}(\mgr) \st \, f \rest{e} \text{ belongs to } H^1(e) \text{ for all edges $e$}\,\Bigr\},
\]
where $H^1(e)$ is the standard first order Sobolev space on the edge $e$, identified with the interval $\mathcal{I}_e$ of given length, and the equality $\Delta f = g \cdot d \mu$ means that
\[
 \sum_{e \in E} \int_e f' \varphi' \, dx_e = \int_{\mgr} \varphi g \, d\mu
\]
for all continuous, edgewise smooth functions $\varphi \colon \mgr \to \C$. Notice in particular that the domain $\operatorname{dom} (\bH)$ embeds naturally into $L^2(\mgr, \mu)$ (since $\Delta f = g \cdot d \mu$, the function $f$  is uniquely determined by harmonic extension on the complement of the support of $\mu$).  

The self-adjoint operator $\bH$ is non-negative and has purely discrete spectrum. The smallest eigenvalue is given by $\lambda = 0$ and its eigenspace consists of all constant functions. The image of $\bH$ is contained in the subspace
\[
L_0^2(\mgr, \mu) = \Big \{ f \in L^2(\mgr, \mu) \,\st \, \int_\mgr f d \mu = 0 \Big \}.
\]
Consider the restriction of $\bH$ to $L^2_0(\mgr, \mu)$, that is, the self-adjoint operator $\bH_0 f := \bH f$ with domain $\dom (\bH_0 ) = \dom(\bH ) \cap L_0^2(\mgr, \mu)$ in the Hilbert space $L_0^2(\mgr, \mu)$. Then the inverse of $\bH_0$ is precisely the integral operator associated with the Green function $ \grg_\mu \colon \mgr \times \mgr \to \R$. That is, the operator $J \colon L_0^2(\mgr, \mu) \to L_0^2(\mgr, \mu)$ given by
\[
 (J f) (x) = \int_\mgr   \grg_\mu(x, y) f(y) \, d\mu(y),  \qquad x \in \mgr,
\]
is bounded and $J \bH_0 f = f$ for all $f \in \dom (\bH_0)$. Moreover, for any $f \in L^2_0(\mgr, \mu)$, the image $Jf$ belongs to $\dom(\bH_0)$ and the equality $\bH_0 J f = f$ holds.

In particular, the results in the present section can be interpreted as degeneration results on the corresponding resolvents. We will address the spectral theoretic implications of our results in a separate publication. For further results on the resolvents of (more general) Schrödinger-type operators on metric graphs with shrinking edges (formulated in more classical spectral theoretic terms), we refer to \cite{bls, Borisov, cacciapuoti} and the references therein.
\end{remark}

\subsection{Definition of tropical Green functions} Using the tropical Poisson equation, we now extend the definition of Green functions to tropical curves. Let $\curve=(G, \pi, l)$ be a tropical curve with underlying graph $G$, edge length function $l \colon E(G) \to (0,+\infty)$, and ordered partition $\pi=(\pi_\infty=(\pi_1, \dots, \pi_r), \pi_\fin)$. Recall also that  a layered measure on $\curve$ is a tuple $\lmu=(\mu_1, \dots, \mu_r, \mu_\fin)$ of measures $\mu_j$ on the graded minors $\Gamma_j$, $j=1, \dots, r, \fin$ (see Section~\ref{ss:LayeredMeasures}). 
Let $\lmu = (\mu_j)_j$ be a non-negative measure of total mass one on $\curve$. That is, each $\mu_j$ is a non-negative Borel measure on the respective minor $\Gamma_j$ and the mass function $\mass=\mass_{\lmu}$ takes value one on the first minor $\Gamma^1$, and value zero on any other component of any graded minor. 

\smallskip

As discussed in Section~\ref{sss:PullbackExplanation}, the normalized edge length function $l \colon E(G) \to (0, + \infty)$ defines a canonical layered metric graph $\Gamma$ in the equivalence class of the tropical $\curve$.  Any point $x \in \Gamma$ can be identified with a layered Dirac measure $\bm{\delta}_x$ on $\curve$ (notice that $\bm{\delta}_x$ has total mass one). In particular, the following Poisson equation (see \eqref{eq:PoissonTropicalCurve})
\begin{equation} \label{eq:TropicalGreenFunction}
\begin{cases}
 \Deltatrop \lf = \bm{\delta}_x - \lmu, \\[1 mm]
 \lf \text{ is harmonically arranged}, \\[1mm]
 \int_\Gamma \lf \, d \lmu = 0
 \end{cases}
\end{equation}
has a unique solution on the tropical curve $\curve$ which we denote by 
\[
\lgri{\curve} (x, \cdot) = \big ( \gri{\curve, 1}(x, \cdot), \dots  ,\gri{\curve, r}(x, \cdot), \gri{\curve, \fin}(x, \cdot) \big).
\]
Each component $\gri{\curve, j}(c, \cdot)$ can be viewed as a function on the space $\Gamma$ by linear interpolation  (see Section~\ref{sss:PullbackExplanation}). Therefore, we get a function 
\[
\begin{array}{cccc}
\lgri{\curve, \lmu} \colon &\Gamma \times \Gamma  &\longrightarrow & \R^{r+1} = \R^{[r]}\times \R^\fin.
\end{array}
\]

\begin{defn}[Green function on a tropical curve]\rm
Notations as above, the {\em Green function} for the measure $\lmu$ on the tropical curve $\curve$ of rank $r$ is the $\R^{r+1}$-valued function
\[
\begin{array}{cccc}
\lgri{\curve, \lmu}  \colon &\Gamma \times \Gamma  &\longrightarrow & \R^{r+1} = \R^{[r]}\times \R^\fin
\end{array}
\]
whose $r+1$ components are defined by
\[
	\big (\lgri{\curve, \lmu} (x,y)  \big )_j = \gri{\curve, j}(x,y)
\]
for every $j=1, \dots, r,\fin$ and any pair of points $x,y$ on $\Gamma$. 
\end{defn}

\begin{remark}[Weight of tropical Green functions] The above procedure generalizes to any layered metric graph $\mgr = (G, \pi, \ell)$ in the conformal equivalence class of a tropical curve $\curve$. However, if the edge length function $\ell$ of $\mgr$ is given by $\ell\rest{\pi_j} = \alpha_j l_j$, for positive real numbers $\alpha_j$, for all layers $\pi_j$, then the components of the resulting Green function on $\mgr$ are simply scaled pull-backs of the Green function on $\Gamma$, that is,
\[
\gri{\mgr, j}(p, q) = \alpha_j \, \gri{\curve, j}^{\! \! \! \! \! \ast}(p, q), \qquad p, q \in \mgr_\pi
\]
for $j=1, \dots, r$.  In other words, $\lgri{\curve,\lmu}$ lives in the space of continuous functions of \emph{weight one} on the tropical curve $\curve$.

More generally, we can define \emph{the space of continuous functions of weight $k$} on a tropical curve $\curve$ as those functions $\lf = (f_1, \dots, f_r, f_\fin)$ whose components rescale by $\alpha^k$ under the rescaling by $\alpha$ of the edge length function $l_j$ in the conformal equivalence class.

The layered metric graph $\Gamma$  provides a natural represent of the conformal equivalence class of the curve $\curve$, and we formulate our definitions in terms of this space. \end{remark}

\subsection{Tame behavior of tropical Green functions} In order to formalize the degeneration of a metric graph $\mgr$, we fix a combinatorial graph $G = (V,E)$ and consider the corresponding space $\mgtropcombin{\grind{G}}$ (see \eqref{eq:TildeMgTrop}). Recall that each point $\thy \in \mgtropcombin{\grind{G}}$ represents a tropical curve, which coincides with the fiber $\unicurvetrop{\thy}$ over $\thy$ in the universal family $\unicurvetrop{\grind G}$ over $\mgtropcombin{\grind{G}}$. Moreover, as a topological space, the fibers $\unicurvetrop{\thy}$ are all naturally homeomorphic to the metric graph $\Gamma$.

\smallskip In what follows, we consider a continuous family of measures $\lmu_\thy$, $\thy \in \mgtropcombin{\grind{G}}$, on the universal curve $\unicurvetrop{\grind G} /\mgtropcombin{\grind{G}}$ with the properties mentioned earlier. More precisely, for each $\thy \in \mgtropcombin{\grind{G}}$, we have a layered measure $\lmu_\thy = (\mu_{\thy,j})_j$ on the tropical curve $\unicurvetrop{\thy}$ such that $\lmu_\thy$ has total mass one, each component $\mu_{\thy,j}$ is a non-negative measure of the form \eqref{eq:SpecialMus} on the corresponding minor, and the family of measured spaces $(\unicurvetrop{\thy}, \lmu_\thy)$, $\thy \in \mgtropcombin{\grind{G}}$, is continuous in the sense of Section~\ref{ss:FamiliesMeasuredSpaces}.

\smallskip

Before stating our main results, we briefly discuss notations. Notice that all fibers of the universal curve $\unicurvetrop{\grind G} /\mgtropcombin{\grind{G}}$ are naturally homeomorphic to each other (by using the obvious homothecies between edges). In particular, given a function $f_\thy$ on some fiber $\unicurvetrop{\thy}$, $\thy \in \mgtropcombin{\grind{G}}$, more precisely, on the canonical layered metric graph associated to $\unicurvetrop{\thy}$, and another base point $\shy \in \mgtropcombin{\grind{G}}$, we can naturally pull back $f_\thy$ to a function $f_\thy^\ast$ on $\unicurvetrop{\shy}$, that is, on the canonical layered metric graph associated to $\unicurvetrop{\thy}$. This notion obviously extends to functions defined on products of fibers $\unicurvetrop{\thy}\times\unicurvetrop{\thy}$.

\smallskip

For any point $\thy \in \mgtropcombin{\grind{G}}$, we denote by $\lgri{\thy} = (\gri{\thy, 1}, \dots, \gri{\thy, r}, \gri{\thy, \fin})$ the Green function on the tropical curve $\unicurvetrop{\thy}$ for the measure $\lmu_\thy$.  If the point $\thy$ represents a metric graph $\mgr$, or equivalently, $\thy = t$ belongs to  $\mggraphcombin{\grind{G}}$, then we also write $\grg_t\colon \mgr\times \mgr \to \R$ or $\grg_\mgr \colon \mgr\times \mgr \to \R$ instead of $\grg_\thy$.

\smallskip

The following analog of Proposition~\ref{prop:leadingterm} provides the leading term of the Green function under degeneration.

\begin{prop}[Leading terms] \label{prop:leadingtermGreen}
Let $(\lmu_\thy)_\thy$ be a continuous family of measures on $\unicurvetrop{\grind{G}} / \mgtropcombin{\grind{G}}$ as above. For $\thy \in  \mgtropcombin{\grind{G}}$, let $h_\thy \colon \unicurvetrop{\thy} \times \unicurvetrop{\thy}  \to \R$ be the function sending $x,y \in \unicurvetrop{\thy}$ to
\[
h_\thy(x,y) = \begin{cases}\frac{1}{L(t)} \, \grg_{t}(x,y)  &t \in \tilde  \mg_{\grind{G}}  \\ \grg_{\thy, 1}(x,y) &\thy \in \partial \mgtropcombin{\grind{G}}  \end{cases},
\]
where $L(t) := \sum_{e \in E(G)} t(e)$. Then the functions $h_\thy$, ${\thy \in \mgtropcombin{\grind{G}}}$, form a continuous family of functions over $\mgtropcombin{\grind{G}}$. More precisely, if $\shy$ converges to $\thy$ in $\mgtropcombin{\grind{G}}$, then
\[
\lim_{\shy \to \thy} \, \sup_{x,y} | h_{\shy}(x,y) - h_{\thy}^\ast(x,y) | = 0.
\]
\end{prop}
\begin{proof}
First of all, notice that by Proposition~\ref{prop:leadingterm}, we already know that $\lim_{\shy \to \thy} \, \sup_{y} | h_{\shy}(x,y) - h_{\thy}(x,y) | = 0$ for every fixed $x$. Hence it only remains to prove that this convergence is uniform in $x$.  However, by \cite[Proposition 3.3]{BR10}, we know that for fixed points $p,q,x$ on a metric graph $\mgr$, the $j$-function $y \mapsto \jfunci{p\tiret q,x}(y)$ has bounded derivatives $|\partial_y \jfunc{p\tiret q,x}(\cdot)| \le 1$. In particular,  \eqref{eq:GreenByJFunction} implies that for any measure $\mu$ on $\mgr$, the Green function satisfies
\[
\sup_{y \in  \mgr} |\grg_\mu(x, y) - \grg_\mu(x', y)| = \sup_{y \in  \mgr} |\grg_\mu(y, x) - \grg_\mu(y, y')| \le |\mu|(\mgr) d(x,x'),
\]
where $d(x,x')$ denotes the distance of the points $x, x' \in \mgr$ in the path metric of $\mgr$. The claim now follows easily by standard arguments, choosing a large finite set of points $x$ on the respective fibers, and then using the above estimate for all other points $x'$. 
\end{proof}

Refining the asymptotics in Proposition~\ref{prop:leadingtermGreen}, we arrive at the following layered expansion of the Green function, when approaching a tropical curve.

\smallskip

Suppose that $\curve = (G, \pi, l)$ is a fixed tropical curve in the boundary stratum $\mgtropcombin{\combind{(G, \pi)}}$ of $ \mgtropcombin{\grind{G}}$, corresponding to an ordered partition $\pi = (\pi_\infty, \pi_\fin)$, $\pi_\infty=(\pi_1, \dots, \pi_r, \pi_\fin)$. The projection map $\pr_\pi \colon \tilde \mg_{\grind{G}} \to \mgtropcombin{\combind{(G, \pi)}}$ assigns to every metric graph $\mgr$ of combinatorial type $G$ with length function $\ell \colon E \to (0, + \infty)$ a tropical curve denoted $\curve_\ell$ of type $(G, \pi)$ whose edge lengths are given by the point $\pr_\pi(\ell) \in \mgtropcombin{\grind{G}}$. 

\begin{thm}[Local expansion of tropical Green functions] \label{thm:GreenFunctionConvergence} Let $(\lmu_\thy)_\thy$ be a continuous family of measures on $\unicurvetrop{\grind{G}} / \mgtropcombin{\grind{G}}$ as above. Assume that the metric graph $\mgr$ degenerates tamely to the tropical curve $\curve$ of rank $r$.
 Then the Green function $\grg_\mgr  \colon \mgr \times \mgr \to \R$ has the expansion
\begin{equation} \label{eq:GreenFunctionConvergence}
	\gri{\mgr} (x,y) =  \sum_{k=1}^r L_k(\ell) \grg_{\ell,k}^\ast (x,y) + \grg_{\ell,\fin}^\ast(x, y) + o(1) 
\end{equation}
where the $o(1)$ term goes to zero uniformly in $x,y$ and $\lgri{\ell} = (\gri{\ell,1},  \dots, \gri{\ell, r}, \gri{\ell, \fin})$ is the Green function on the tropical curve $\curve_\ell$ for the measure obtained by pushing out $\mu_\ell = \mu_t$ to the tropical curve $\curve_\ell$. Moreover, $\gri{\ell,k}(x, y)$ converges to the $k$-th component $\gri{\curve, k}$ of the Green function $\lgri{\curve}$ on $\curve$, that is,
\[
\sup_{x,y \in \mgr} | \grg_{\ell,k}^\ast(x, y)  - \gri{\curve,k}(x, y) | \to 0, \qquad k=1, \dots, r, \fin,
\]
as $\mgr$ degenerates tamely to the tropical curve $\curve$ in $\mgtropcombin{\grind{G}}$. 
\end{thm}
\begin{proof}
Recall that, by \eqref{eq:GreenByJFunction}, the Green function can be obtained by integrating the $j$-function. Moreover, the $j$-function has a layered expansion by Lemma~\ref{lem:JfunctionTame}. However, by Lemma~\ref{lem:SolutionFormulaJFunctionTropical}, integrating each summand in the expansion separately, we recover precisely the pull-back of the Green function $\gri{\ell}$ on the tropical curve $\curve_\ell$ to $\mgr$. The claim is proven.
\end{proof}

\begin{thm}[Weak tameness of tropical Green functions]
The family of Green functions $\lgri{\thy}$, $\thy \in \mgtropcombin{\grind{G}}$, is weakly tame under the assumptions of Proposition~\ref{prop:leadingtermGreen}. 
\end{thm}
That is, the analog of Theorem~\ref{thm:MainGraphWeak} holds true for Green functions.  More precisely, the tropical log map $\Lognoind \colon \mgtropcombin{\grind{G}} \to  \partial_\infty \mgtropcombin{\grind{G}}$induces the decomposition $\mgtropcombin{\grind{G}}  = \bigsqcup_\pi \inn R_\pi$ for $\inn R_\pi = \Lognoind^{-1}(\mgtropcombin{\combind{(G, \pi)}})$ (see~\eqref{eq:DecompositionFromTropicalLogMap}).
We push out the measure $\lmu_\thy$ on fibers of the universal curve $\unicurvetrop{\grind{G}}$ over $\inn R_\pi$ to the fibers over $ \mgtropcombin{\combind{(G, \pi)}} $ by the projection map $\pr_\pi \colon\inn R_\pi \to\mgtropcombin{\combind{(G, \pi)}}$. For $t \in \inn R_\pi$, let $\layh_{t}$ be the Green function for the pushed-out measure on the tropical curve $\unicurvetrop{\thy'}$, $\thy' = \pr_\pi(t)$. Using homotecies of edges, we can view alternatively view $\layh_{t}$ as a function $\layh_{t} \colon \unicurvetrop{\thy} \times \unicurvetrop{\thy} \to \R^{r_\thy+1}$  on the tropical curve $\unicurvetrop{\thy}$, $\thy = \Lognoind(t)$. Clearly, for $\thy \in \mgtropcombin{\combind{(G, \pi)}}$, we recover the original Green function $\lgri{\thy} = \layh_{\thy}$ of the measure $\lmu_{\thy}$. Moreover, the family $\layh_{\thy}$,  $\thy \in  \mgtropcombin{\grind{G}}$, satisfies the conditions in Definition~\ref{def:WeakGlobalTameTropical}.

\begin{proof} The proof is analogous to the proof of Theorem~\ref{thm:GreenFunctionConvergence}, using weak tameness of the solutions of Poisson equations, Theorem~\ref{thm:MainGraphWeak}.
\end{proof}

If we assume that the measures $\mu_\thy$ converge sufficiently fast when approaching a boundary stratum $\mgtropcombin{\combind{(G, \pi)}}$ in $\mgtropcombin{\grind{G}}$, we recover the following analog of Theorem~\ref{thm:main_graph_strong}.

\begin{thm}[Stratumwise tameness of tropical Green functions] \label{thm:GraphStrongGreenFunction}
Let $(\lmu_\thy)_\thy$ be a continuous family of measures on $\unicurvetrop{\grind{G}} / \mgtropcombin{\grind{G}}$ as above. Suppose that, in addition, for every boundary stratum $\mgtropcombin{\combind{(G, \pi)}}$, we have
\begin{align} \label{eq:ConditionStrongGreenFunction}
L(t) \big (\lmu_t - \lmu_{\pr_\pi(t)}^\ast\big), \qquad  \text{where } L(t) := \sum_{e \in E} t(e),
\end{align}
goes to zero weakly if $t \in \tilde \mg_{\grind{G}}$ converges tamely in $\mgtropcombin{\grind{G}}$ to a point $\thy$ in the closure of $\mgtropcombin{\combind{(G, \pi)}}$. Here, $\lmu_{\pr_\pi(t)}^\ast$ denotes the natural pull-backs of the measures $\lmu_{\pr_\pi(t)}$ to the fiber of $t$ in $\unicurvetrop{\grind{G}}$.

Then the associated Green functions $\lgri\thy$, ${\thy \in \mgtropcombin{\grind{G}}}$, form a stratumwise tame family of functions. That is, if $t \in \tilde  \mg_{\grind{G}}$ converges in the tame topology to a point $\thy$ in the closure of the stratum $\mgtropcombin{\combind{(G, \pi)}}$, then we get
\[
\grg_t(x,y) =  \sum_{k=1}^r L_k(t) \grg_{\pr_\pi(t),k}^\ast (x,y) + \grg_{\pr_\pi(t),\fin}^\ast(x, y) + o(1)  
\]
where the $o(1)$ term goes to zero uniformly in $x$ and $y$.
\end{thm}
\begin{proof}
The theorem is easily deduced from Theorem~\ref{thm:main_graph_strong} together with Theorem~\ref{thm:GreenFunctionConvergence}.
\end{proof}
We conclude the section with a few examples. We study Green functions for three particular measures:
\emph{standard (uniform) discrete measures, Lebesgue measures, canonical measures}.

\subsection{Green functions of discrete measures}
For any metric graph $\mgr$ with model $G = (V,E)$, represented by a point $t \in  \mggraphcombin{\grind{G}}$, consider the measure on $\mgr$
\[
\mu_t  := \frac{1}{\abs V}\sum_{v \in V} \delta_v
\]
that we call the standard discrete measure. By definition, $\mu_t$ is non-negative and has mass one on $\mgr$. 

\begin{remark}[Remark~\ref{rem:SpectralTheory} continued]
By identifying $\C^V \cong L^2(\mgr, \mu_t)$, the self-adjoint operator associated with the measure $\mu_t$ (see Remark~\ref{rem:SpectralTheory}) coincides with the {\em discrete graph Laplacian} $\bH _t  \colon \C^V  \to \C^V $, which maps $\varphi \in \C^V$ to 
\[
(\bH _t  \varphi) (v) =  \sum_{u \sim v} \Big ( \sum_{e = uv} \frac{ 1}{t(e)} \Big )  (\varphi(u) - \varphi(v)), \qquad v \in V.
\]
More precisely, the self-adjoint operator associated to $\mu_t$ is unitarily equivalent to $|V| \cdot \bH_t$.
\end{remark}
Altogether, we end up with a family of measures $\mu_t$, $t \in \mggraphcombin{\grind{G}}$, on the metric graph part of the space $\mgtropcombin{\grind{G}}$. In order to apply our results, we first have to extend this family to a continuous family of measures on $ \mgtropcombin{\grind{G}}$. Let $\curve$ be a tropical curve of rank $r$, represented by a point $\thy \in \mgtropcombin{\combind{(G,\pi)}}$, for an ordered partition $\pi=(\pi_\infty=(\pi_1, \dots, \pi_r), \pi_\fin)$ of $E$. Define $\lmu_\curve = \lmu_\thy = (\mu_{\thy, j})_j$ as the layered measure on $\curve$ given by
\begin{align*}
&\lmu_\thy = (\mu_1,  \dots, \mu_r, \mu_\fin), &\mu_{\thy, j} = \sum_{u \in  V(\Gamma^j)} \frac{\abs{\kappa_j^{-1}(u)}}{\abs V} \cdot  \delta_u, \quad j = 1, \dots, r, \fin,
\end{align*}
where $\kappa_j \colon V \to V(\Gamma^i)$ is the contraction map. Clearly, we have the following property.

\begin{prop}
The measures $\mu_\thy$, $\thy \in \mgtropcombin{\grind{G}}$, form a continuous family of measures on $\unicurvetrop{\grind{G}} / \mgtropcombin{\grind{G}}$. In fact, $\lmu_t = \lmu_{\pr_\pi(t)}^\ast$ holds true for any ordered partition $\pi$ and $t \in  \tilde  \mg_{\grind{G}}$.
\end{prop}

Hence all results of the previous section can be applied. In particular, Theorem~\ref{thm:GraphStrongGreenFunction} yields the full asymptotic of the metric graph Green functions under a tame degeneration.

\subsection{Green functions of Lebesgue measures} \label{ss:ExampleLesbesgue}
Recall that any metric graph $\mgr$ carries a natural Lebesgue measure $\lambda_\mgr$, which is simply the sum of the Lebesgue measures $d x_e$ of its edges (the edges are identified with intervals of respective lengths).

\smallskip

For a metric graph $\mgr$ of combinatorial type $G = (V,E)$, represented by $t \in  \mggraphcombin{\grind{G}}$, the renormalized measure
\begin{align} \label{eq:LebesgueGraph}
&\mu_t = \frac{1}{L(\mgr)} \cdot \lambda_{\mgr}
\end{align}
has mass total mass one on $\mgr$. Let $ \grg_t \colon \mgr \times \mgr \to \R$ be the corresponding Green function. 

\begin{remark}[Remark~\ref{rem:SpectralTheory} continued]
From the spectral theoretic point of view (see Remark~\ref{rem:SpectralTheory}), the self-adjoint operator $\bH$ associated in $L^2(\mgr, \mu_t)$ with the measure $\mu_t$ is the so-called {\em Kirchhoff Laplacian} on $\mgr$, multiplied by the factor $L(\mgr)$ (see, e.g., \cite{BerkolaikoKuchment} and the references therein).
\end{remark}
As before, we can easily extend the measures $\mu_t$, $t \in \mgtropcombin{\grind{G}}$, to a continuous family of measures on $ \mgtropcombin{\grind{G}}$. Let $\curve$ be a tropical curve of positive rank $r$, represented by a point $\thy \in \mgtropcombin{\combind{(G,\pi)}}$, $\pi = (\pi_\infty=(\pi_1, \dots, \pi_r), \pi_\fin))$.  Define $\lmu_\thy = (\mu_{\thy, j})_j$ as the (layered) measure on $\curve$ given by
\begin{align} \label{eq:LebesgueTropical}
&\mu_{\thy,1} = \lambda_{\Gamma^1} \text{ on } \Gamma^1,  &\mu_{\thy,j} \equiv 0 \text{ on } \Gamma_j \text{ for } j = 2, \dots, r, \fin.
\end{align}
By definition, $\mu_\thy$ has total mass one on $\curve$. The following statement is easy to check.
\begin{prop}
The measures $(\mu_\thy)_{\thy \in \mgtropcombin{\grind{G}}}$ form a continuous family of measures on $\unicurvetrop{\grind{G}} / \mgtropcombin{\grind{G}}$.
\end{prop}
In particular, both Proposition~\ref{prop:leadingtermGreen} and Theorem~\ref{thm:GreenFunctionConvergence} can  be applied. However, it turns out that condition \eqref{eq:ConditionStrongGreenFunction} in Theorem~\ref{thm:GraphStrongGreenFunction} fails.

\smallskip

Fix again a tropical curve $\curve=(G, \pi, l)$, represented by a point $\thy$ in the stratum $\mgtropcombin{\combind{(G, \pi)}}$ for a partition $\pi= (\pi_\infty=(\pi_1, \dots, \pi_r), \pi_\fin)$. For $n \in \{2, \dots,r, \fin \}$, introduce the following measure
\[
\nu_\thy^n = \sum_{e \in \pi_n} dx_e - \lambda_{\Gamma^n} (\Gamma^n) \cdot  \sum_{e \in \pi_1} dx_e.
\]
on the space $\Gamma$. Since $\nu_\thy^n$ has mass zero on $\Gamma$, it can be identified with a layered measure $\lnu_\thy^n$ on $\curve$ (see Proposition~\ref{prop:MassZeroMeasures}). Clearly, the first component of $\lnu_\thy^n$ is given by
\[
\omega_\thy^n = \omega_\curve^n = \sum_{v \in V(\Gamma^1)} \Big ( \sum_{\substack{e \in \pi_n\colon \\ \kappa_1(e) = v}} \ell(e) \Big )  \cdot \delta_v- \lambda_{\Gamma^n} (\Gamma^n) \cdot \lambda_{\Gamma^1},
\]
which is a measure of mass zero on $\Gamma^1$ ($\kappa_1$ is the contraction map from $\Gamma$ to $\Gamma^1$).

\smallskip

For any point $t \in \mggraphcombin{\grind{G}}$ representing a metric graph $\mgr$, one easily verifies the following equality
\begin{equation} \label{eq:HigherAsymptoticsLebesgue}
L(t) \big (\mu_t - \lmu_{\pr_\pi(t)}^\ast \big ) =   \sum_{j \in \{2, \dots, r, \fin \}} {L_j(t)} \big ( \lnu_{\pr_\pi(t)}^n \big )^\ast \text{ on } \mgr,
\end{equation} 
where again $\ast$ denotes the pullback of measures. In particular, condition \eqref{eq:ConditionStrongGreenFunction} fails. However, the higher order terms in \eqref{eq:HigherAsymptoticsLebesgue} lead to the following complete asymptotics.

\begin{thm}[Full asymptotics of tropical Lebesgue Green functions] \label{thm:Lebesgue_Green_functions} Notations as in Theorem~\ref{thm:GreenFunctionConvergence}, assume that the metric graph $\mgr$ with edge length function $\ell$ degenerates tamely to the tropical curve $\curve$ of rank $r$ lying on the boundary stratum $\mgtropcombin{\combind{(G, \pi)}}$. Consider the corresponding Green functions $\gri{\mgr}$ and $\lgri{\curve} = (\gri{\curve,j})_j$ for the measures $\mu_\mgr$ and $\lmu_\curve = (\mu_{\curve,j})_j$ defined by \eqref{eq:LebesgueGraph} and \eqref{eq:LebesgueTropical}, respectively. Then, we have 
\[
	\grg_\mgr (x,y) =  \sum_{k=1}^r L_k(\mgr) \grg_{\ell,k}^\ast (x,y) + \grg_{\ell,\fin}^\ast(x, y) + o(1),
\]
where the $o(1)$ goes to zero uniformly in $x,y$ and $\lgri{\ell} = (\gri{\ell,j})_j$ is the Green function on the tropical curve $\curve_\ell$ for the measure obtained by pushing out $\mu_\mgr$ to $\curve_\ell$.

Moreover, comparing the Green functions $\lgri{\ell}$ with the Green function $\lgri{\curve_\ell}$ on $\curve_\ell$, we have
\[
\gri{\ell, j} = \gri{\curve_\ell, j} +\frac{o(1)}{ L_j(\ell)},
\]
for all $j \in \{2, \dots, r, \fin\}$ and
\begin{align*}
\gri{\ell,1}(x,y) = \gri{\curve_\ell, 1} (x,y) + \sum_{n \in \{2, \dots,r, \fin\}} \frac{L_j(\ell)}{L_1(\ell)} \, \Big(  \psi_{\curve_\ell, n} (y) - \int_{\Gamma_{\ell}^1}  \gri{\curve_\ell, 1}(x,y) \, d \omega_{\curve_\ell}^n(y) \Big)  + \frac{o(1)}{ L_1(\ell)},
\end{align*}
where $\psi_{\curve_\ell, j} \colon\Gamma_{\ell}^1 \to \R$ is the solution of the Poisson equation~\eqref{eq:PoissonMetricGraphs} for the measures $\mu = \omega_{\curve_\ell}^n$ and $\nu = \mu_{\curve_\ell, 1}$ on the first minor $\Gamma_{\ell}^1$ of $\curve_\ell$.
\end{thm}
\begin{proof}
The claim is a consequence of Theorem~\ref{thm:GreenFunctionConvergence} together with \eqref{eq:HigherAsymptoticsLebesgue}. Indeed, all other contributions obtained from integration in Lemma~\ref{lem:SolutionFormulaJFunctionTropical} go to zero uniformly, if we assume that  $\mgr$ degenerates to $\curve$ in the tame topology.
\end{proof}

\subsection{Canonical Green functions} \label{sec:canonical_green_functions_tropical}
We finally investigate the behavior of the canonical Green function on degenerating metric graphs, which will be related to the Arakelov Green functions in the third part of our paper. 

We begin by introducing a notion of {\em canonical measure for tropical curves}. For tropical curves of full sedentarity (that is, $\pi_\fin = \varnothing$), a notion of canonical measure was already defined in \cite{AN}. In what follows, we adapt the definition from \cite[Section 5]{AN} to the general case. The following example shows this is necessary for describing the limits of canonical measure on metric graphs, when the edge lengths tend to infinity.

\begin{example} \label{ex:CounterExampleCanonicalMeasure}
Consider a metric graph with model $\mgr=(V,E, \ell)$ consisting of two vertices connected by three parallel edges $e_i$, $i \in [3]$ (see Figure~\ref{fig:BB} in the introduction). Denoting by $\ell_i = \ell(e_i)$ the length of $e_i$, $i \in [3]$, the Foster coefficients (see~\cite{AN} or discussion in Section~\ref{sec:canonical_green_functions_tropical}) are given by
\begin{align*}
& \text{ $ \mu(e_1) = \frac{\ell_1 \ell_2 + \ell_1 \ell_3}{\ell_1 \ell_2 + \ell_2 \ell_3 + \ell_1 \ell_3}$}  &\text{ $ \mu(e_2) = \frac{\ell_2 \ell_3 + \ell_1 \ell_2}{\ell_1 \ell_2 + \ell_2 \ell_3 + \ell_1 \ell_3}$}  &&\text{ $ \mu(e_3) = \frac{\ell_1 \ell_3 + \ell_2 \ell_3}{\ell_1 \ell_2 + \ell_2 \ell_3 + \ell_1 \ell_3}$}
\end{align*}
Suppose that all edge lengths go to infinity, that is, $\ell_i \to \infty$ for all $i$. Then, in the sense of the compactified moduli space $\cancomp{\mggraph{\grind g}}$, the metric graph  $\mgr$ degenerates to a metric graph $\mgr_\infty$ over $G = (V,E)$ with all edges $e_i$ having length $\ell_{\mgr_\infty} (e_i) = + \infty$. However, the limit of the masses of edges, as $\mgr$ degenerates to $\mgr_\infty$ in $\cancomp{\mggraph{\grind g}}$ does not exist.

To see this, note that different ways of approaching $\mgr_\infty$ in $\cancomp{\mggraph{\grind g}}$ yield different limit mass distributions. For instance, if $\ell_1 = \ell_2 = \ell_3 \to \infty$, then the limit of the Foster coefficients are
\[
\mu_\infty(e_1) = \mu_\infty(e_2) = \mu_\infty(e_3) = 2/3.
\]
On the other hand, if all $\ell_i$'s go to infinity, but $\ell_1 \gg \ell_2 = \ell_3$, the limits are given by
\[
\mu_\infty(e_1) = 1,  \qquad \mu_\infty(e_2) = \mu_\infty(e_3) = 1/2.
\]
In particular, the canonical measure \emph{cannot be extended continuously to the compactification} $\cancomp{\mggraph{\grind g}}$. That is, there is no notion of a canonical measure $ \mu^\can_{\mgr_\infty}$ for the boundary points $\mgr_\infty$ of  $\cancomp{\mggraph{\grind g}}$ such that $\mu^\can_\mgr \to \mu^\can_{\mgr_\infty}$ whenever a metric graph  $\mgr$ converges to $\mgr_\infty$ in $\cancomp{\mggraph{\grind g}}$. 
\end{example}

Let $\curve$ be a tropical curve with underlying graph $G = (V,E)$, layering $\pi = (\pi_\infty, \pi_\fin)$, $\pi_\infty=(\pi_1, \dots, \pi_r)$ and edge length function $l \colon E \to (0, + \infty)$.  The canonical measure on $\curve$ is the layered measure $\lmu^{\can}=(\mu^{\can}_{1}, \dots, \mu^{\can}_{r}, \mu^{\can}_{\fin})$ whose components $\mu^{\can}_{j}$, $j=1, \dots r, \fin$ are the measures
\[
{\mu}^{\can}_{j} = \frac{1}{h} \Big(\mu^j_{\Zh} + \sum_{v \in V^j} h_v \cdot \delta_v \Big ) \qquad \text{on } \Gamma^j.
\]
In the above definition...
\begin{itemize}
\item $\graphgenus$ is the genus of the underlying graph $G$,
\smallskip

\item $\mu^j_{\Zh}$ is the Zhang measure of the $j$-th graded minor $\Gamma^j$,  and
\smallskip

\item for every vertex $v \in V^j$, $\graphgenus_v$ is the genus of the connected component of $\Gamma^{j+1}$ associated with $v$ (by convention, $\graphgenus_v := 0$ for all vertices $v$ of the last minor $\Gamma^\fin$).
\end{itemize}
By definition, $\lmu^{\can}$ has total mass one on $\curve$. In the case $\curve = \mgr$ is a metric graph, $\lmu^{\can}$ obviously coincides with the canonical measure of $\mgr$, that is, $\lmu^{\can}=\mu^{\can}_{\fin} = \frac 1{\graphgenus}\mu_{\Zh}$.

\begin{remark} \label{rem:GlobalCanonicalMeasure} Notice that $\lmu^{\can}$ can be identified with the measure ${\mu}^{\can}$ on $\Gamma$ given by
\begin{equation} \label{eq:CanonicalMeasureWholeSpace}
{\mu}^{\can}= \frac{1}{h}\sum_{j \in \{1, \dots, r, \fin\}} \frac{ \mu_{\Gamma_j}(e)  }{l(e)} \, dx_e,
\end{equation}
where $\mu_{\Gamma^j}(e)$, $e \in \pi_j$ denote the Foster coefficient in the minor $\Gamma^j$ (see \eqref{eq:FosterCoefficient}). More precisely, the layered measure $\lmu^{\can}$ is obtained from ${\mu}^{\can}$ by the restriction-push-out procedure in \eqref{eq:PushoutMeasure}.

If $\pi_\fin = \varnothing$, then $\curve$ is a tropical curve in the sense of \cite{AN} and ${\mu}^{\can}$ coincides with the canonical measure introduced in \cite[Section 5.4]{AN}. 
\end{remark}

We consider the variation of canonical measures on families of tropical curves. Fix a finite graph $G=(V,E)$ and consider the space  $\mgtropcombin{\grind{G}}$. For any base point $\thy \in \mgtropcombin{\grind{G}}$, representing a tropical curve of combinatorial type $G$, we equip its fiber $\unicurvetrop{\thy}$ in the universal family $\unicurvetrop{\grind{G}}$ with the canonical measure $\lmu_\thy = \lmu_{\thy}^{\can} = \lmu_{ \curve}^{\can}$. Notice that the fibers above the subspace $\mggraphcombin{\grind{G}}$ correspond to metric graphs equipped with the Zhang measure.

\smallskip

We have the following basic continuity result (see Section~\ref{ss:FamiliesMeasuredSpaces} for definitions).

\begin{thm} \label{thm:CanMeasuresCont}
Let $G= (V,E)$ be a finite graph. Then the canonical measures $\mu_{\thy}^{\can}$, $\thy \in \mgtropcombin{\grind{G}}$, form a continuous family of measures on $\unicurvetrop{\grind{G}} / \mgtropcombin{\grind{G}}$. 
\end{thm}

\begin{proof}
The proof of \cite[Theorem 6.1]{AN} carries over line to line to the present setting and we omit the details. 
\end{proof}

In what follows, we extend Theorem~\ref{thm:CanMeasuresCont} to the moduli space of tropical curves $\mgtrop{\grind{g}}$. The strata of $\mgtrop{\grind{g}}$ correspond to {\em augmented layered graphs} $(G, \pi, \genusfunction)$, where $\genusfunction\colon V(G) \to \N \cup\{0 \}$ is a genus function on the vertex set. Hence we first define a canonical measure in this setting.

\smallskip 

Let be $\curve=(G, \pi, \ell)$ a tropical curve of rank $r$ and $\genusfunction\colon V(G) \to \N \cup\{0 \}$. The associated canonical measure is the layered measure $\lmu^\can=(\mu^\can_1, \dots, \mu^\can_r, \mu^\can_\fin)$ on $\curve$, whose $j$-th component $\mu^\can_j$ is given by
\[
\mu^\can_j =\frac 1g \Big( \mu^j_{\Zh} + \sum_{v \in V^j}  \big ( \graphgenus_v +\sum_{\substack{u \in V(G) \colon \\ \cont{j} (u) = v}}  \genusfunction(u) \big ) \cdot \delta_v  \Big) \qquad \text{on } \Gamma_j,
\]
where $\graphgenus$ is the genus of $G$, $g = \graphgenus + \sum_{u \in V} \genusfunction(u)$ is the genus of the tropical curve, and $\mu^j_{\Zh}$ is the Zhang measure on $\Gamma^j$ and $\graphgenus_v$ is the genus of the component of $\Gamma^{j+1}$ which has been contracted to the vertex $v$. Note that $\lmu^\can$ is a layered measure of total mass one on $\curve$.

\smallskip

Consider the moduli space $\mgtrop{\grind{g}}$ of genus $g$ tropical curves. For any base point $\thy \in \mgtrop{\grind{g}}$, representing (an equivalence class of) a marked tropical curve $\curve$, with combinatorial type an augmented graph $G$, the fiber $\unicurvetrop{\grind g}\rest{\thy}$ in the universal family $\unicurvetrop{\grind{g}}$ over $\mgtropcombin{\grind{G}}$ can be equipped with the canonical measure $\mu_\thy = \mu_{\thy}^{\can} = \mu_{ \curve}^{\can}$. Using Theorem~\ref{thm:CanMeasuresCont}, one can easily deduce the following result (see Section~\ref{ss:FamiliesMeasuredSpaces}). 

\begin{thm} \label{thm:CanMeasuresContModuli}
The canonical measures $\mu_{\thy}^{\can}$, $\thy \in \mgtropcombin{\grind{G}}$, form a continuous family of measures on the universal curve $\unicurvetrop{\grind g}$ over $\mgtropcombin{\grind{G}}$. 
\end{thm}

Finally, we investigate the corresponding Green functions on degenerating metric graphs.
\begin{defn}[Canonical Green function on a tropical curve]\rm
The {\em canonical Green function} on a tropical curve $\curve$ is the Green function  $\lgri{\curve} = \lgri{\curve}^{\can} = (\gri{\curve,j})_j$ of its canonical measure $\lmu^{\can}$.
\end{defn}
By Theorem~\ref{thm:CanMeasuresCont}, both Proposition~\ref{prop:leadingtermGreen} and Theorem~\ref{thm:GreenFunctionConvergence} apply. Hence we obtain the leading term and a layered expansion for the canoncial Green function on degenerating metric graphs in terms of Green functions on tropical curves.

However, the convergence of canonical measures under a degeneration $\mgr \to \curve$ of a metric graph $\mgr$ to a tropical curve $\curve$ is not very fast. Similar as for the Lebesgue measure, the condition \eqref{eq:ConditionStrongGreenFunction} fails and hence the strong convergence result, Theorem~\ref{thm:GraphStrongGreenFunction}, cannot be applied. On the other hand, recall that in Proposition~\ref{thm:Lebesgue_Green_functions}, we could recover the full asymptotic of the Green function using a higher order approximation of the measure (see \eqref{eq:HigherAsymptoticsLebesgue}). Due to their more complicated structure, this approach is more difficult to adapt for canonical measures. However, if the layering of the limit tropical curve $\curve$ has the simple form $\pi = (\pi_1, \pi_\fin)$, that is if $\curve$ is of rank one, there is a short combinatorial solution. We finish the section with a detailed discussion of this case.

\medskip
We recall first a representation of the canonical measure in terms of spanning trees. For a fixed graph $G = (V,E)$, let 
\[
	\mathcal{T} = \mathcal{T}(G) =\big \{ T =(V(T), E(T) | \, \text{T is a spanning tree of $G$} \big \}.
\]
be the set of spanning trees of $G$. Given some edge length function $\ell \colon E \to (0, + \infty)$, the weight of a spanning tree $T \in \cT$ is defined by
\[
	\omega(T) = \prod_{e\in E \setminus E(T)} \ell(e).
\]
In particular, we can associate to any metric graph $\mgr$ or tropical curve $\curve$ over $G$, the weight $\omega(T)$ of $T$ with respect to  $\mgr$ or $\curve$. If we want to stress the geometrical object behind, we will write $\omega_\mgr(T)$ and $\omega_\curve(T)$.

The canonical measure on a metric graph $\mgr$ of combinatorial type $G$ with edge length function $\ell \colon E(G) \to (0, +\infty)$ can be written as
\[
\mu^{\can} = \sum_{e \in E(G)} \frac{\mu(e)}{\ell(e)} \cdot dx_e
\]
where $dx_e$ denotes the uniform Lebesgue measure on the edge $e$ and the edge masses $\mu(e)$, $e \in E$, also called {\em Foster coefficients}, can be expressed as
\[
\mu(e) = {\sum_{\substack{T \in \cT(G)\colon \\ e \notin E(T)}} \omega(T)} / {\sum_{T \in \cT(G)} \omega(T)}.
\]

Consider a layering of $G$ having the form $\pi = (\pi_1, \pi_\fin)$. For a spanning tree $T \in \cT(G)$, let
 \[
 x_j(T) := \abs{E(T) \cap \pi_j}, \qquad j \in \{1, \fin\}.
 \]
 be the number of edges of $T$ in the respective layers.
 It is easily shown that for every $T \in \cT(G)$,
 \[
 x_1(T) + x_\fin(T) = \abs{E} -  h, \qquad x_1(T) \ge \abs{\pi_1} -  h^1, \qquad x_\fin(T) \le \abs{\pi_\fin}-  h^\fin,
 \]
where $h^1, h^\fin$ are the genera of the minors $\grm{\pi}{1}(G)$ and $\grm{\pi}{\fin}(G)$, and $h = h^1 + h^\fin$ is the genus of $G$. Using the layering $\pi$, we define two special classes of spanning trees on $G$. Consider first
\[
\cT_\pi(G) =\bigl  \{ T \in \cT(G) \st \, x_i(T) = \abs{\pi_i} -  h^i \text{ for } i \in \{1, \fin \} \bigr  \}.
\]
We show in~\cite{AN} that $\cT_\pi(G)$ can be identified with the product
\[
\cT_\pi(G) = \cT(\grm{\pi}{1}(G)) \times \cT(\grm{\pi}{\fin}(G)).
\]
The spanning trees in $\cT_\pi(G)$ are called \emph{layered}.

The second class of trees is defined by
\[
\cT_{1}(G) = \big  \{ T \in \cT(G) \st \, x_1(T) = \abs{\pi_1} -  h^1 + 1 \big  \},
\]
so every tree $T \in \cT_{1}(G)$ has exactly $\abs{\pi_\fin} - h^\fin -1$ edges in $\pi_\fin$.

\medskip

Notice that for the metric graph $\mgr$ of combinatorial type $G$ with edge length function $\ell$, the following scaling property holds
\begin{equation} \label{eq:ScalingTrees}
\omega_\mgr (T) = L_1(\ell)^{\abs{\pi_1}  - x_1(T)} \, \, \omega_{\curve_\ell}(T), \qquad T \in \cT(G),
\end{equation}
where $\curve_\ell$ is the tropical curve of type $(G, \pi)$ with normalized edge lengths in the first layer representing the layered metric graph $(\mgr, \pi)$.

In particular, one obtains the following lemma.
\begin{lem} \label{lem:FosterApprox}
Notations as above, assume that the metric graph $\mgr$ over $G$ degenerates to a tropical curve $\curve$ of type $(G, \pi)$ with $\pi = (\pi_1, \pi_\fin)$. Then, for all edges $e \in E$,
\[
\mu(e) = \mu_{\curve_\ell}(e) + \Big (\sum_{T \in \cT_\pi} \omega(T) \Big)^{-1} \Big (\sum_{\substack{T \in \cT_{1}\colon \\ e \notin T}} \omega_{\curve_\ell} (T) - \mu_{\curve_\ell}(e) \sum_{T \in \cT_{1}} \omega_{\curve_\ell}(T)   \Big ) + O(L_1(\ell)^{-2}),\]
where $\mu_{\curve_\ell}(e) = \mu_{\curve_\ell}^{\can}(e)$ denotes the mass of the edge $e$ in the canonical measure on the tropical curve $\curve_\ell$.
\end{lem}
\begin{proof}
The proof is a straightforward consequence of the scaling property \eqref{eq:ScalingTrees}. First of all, as $\mgr$ degenerates to $\curve$ in $\mgtropcombin{\grind{G}}$, the weights of spanning trees have orders
\[ L_1(\ell) ^{h_1 - \abs{\pi_1}} \, \, \omega(T) = \begin{cases} O(1) &\text{if } T \in \cT_\pi  \\   O(L_1(\ell)^{-1}) &\text{if } T \in \cT_{1} \\ O(L_1(\ell)^{-2})
 &\text{if } T \notin \cT_{1} \cup \cT_\pi \end{cases}
 \]
This implies that for every edge $e \in E$, as $\mgr$ degenerates to $\curve$ in $\mgtropcombin{\grind{G}}$,
\[
\mu(e) = {\sum_{\substack{T \in \cT_{\pi} \cup \cT_{1} \colon \\ e \notin E(T)}} \omega(T)} / {\sum_{T \in \cT_{\pi} \cup \cT_{1}} \omega(T)} + O(L_1(\ell)^{-2}).
\]
The above expression follows easily after taking into account that
\begin{align*}
\mu_{\curve_\ell} (e) &=  {\sum_{\substack{T \in \cT(\grm{\pi}{1}(G)) \colon \\ e \notin E(T)}} \omega_{\Gamma^1}(T)  \cdot \sum_{\substack{T \in \cT(\grm{\pi}{\fin}(G)) \colon \\ e \notin E(T)}} \omega_{\Gamma_\fin}(T) } / {\sum_{T \in \cT(\grm{\pi}{1}(G)) } \omega_{\Gamma^1}(T)  \cdot \sum_{T \in \cT(\grm{\pi}{\fin}(G)) } \omega_{\Gamma_\fin}(T) }   \\ &= {\sum_{\substack{T \in \cT_{\pi} \colon \\ e \notin E(T)}} \omega(T)} / {\sum_{T \in \cT_{\pi}} \omega(T)}.
\end{align*}
\end{proof}

We reformulate the above lemma in terms of measures. Fix a tropical curve $\curve$ with the combinatorial triplet $(G, \pi, l)$ and layering $\pi = (\pi_1, \pi_\fin)$. For a spanning tree $T \in \cT_{1}$, consider the following measure on $\Gamma$,
\[
\mu_T = \frac{\omega(T)}{\sum_{T \in \cT_\pi} \omega(T)} \, \sum_{e \in E} \big ( \chi_{T^c} (e) - \mu_{\curve}(e)  \big ) \frac{dx_e}{l(e)} ,\]
where $\chi_{T^c}(\cdot)$ is the indicator function of $E(T)^c  = E \setminus E(T)$. Since $\mu_T$ has total mass zero on $\Gamma$, it can be identified with a (layered) measure $\lmu_T$ of mass zero on $\curve$ (see Poposition~\ref{prop:MassZeroMeasures}). Taking the sum over all spanning trees in the set $\cT_1(G)$,
\[
\lnu = \lnu_{\curve} = \sum_{T \in \cT_{1}} \lmu_T
\]
is again a measure of mass zero on $\curve$. The first component $\nu_1$ of $\lnu = (\nu_1, \nu_\fin)$  is the measure
\[
\nu_1 = \nu_{\curve, 1} = \sum_{e \in \pi_1} \nu_1(e) \frac{dx_e}{l(e)} + \sum_{v \in V^1} \nu_1(v) \delta_v,
\]
on the first minor $\Gamma^1$ of $\curve$, where the edge and vertex masses are
\begin{align*}
\nu_1(e) &= \sum_{T \in \cT_{1}} \frac{\omega_\curve(T)}{\sum_{T \in \cT_\pi} \omega_\curve(T) } \, \Big (\chi_{T^c}(e) - \mu_\curve(e) \Big ),  && e \in \pi_1 \\
\nu_1(v) &= \sum_{T \in \cT_{1}} \frac{\omega_\curve(T)}{\sum_{T \in \cT_\pi} \omega_\curve(T) } \Big  (\abs{E(T)\cap H_v} - h_v \Big ),  &&v \in V(\grm{\pi}{1}(G))
\end{align*}
where $H_v$ is the connected component of $G_\pi^2$ associated with a vertex $v \in V^1$ and $h_v$ is its genus.

\medskip

We denote by $\psi_\curve \colon  \Gamma^1 \to \R$ the solution on $\Gamma^1$ to the Poisson equation
\[ 
 \begin{cases}
 \Deltaind{\Gamma^1} f = \nu_1 \\
\int_{\Gamma^1} f \, d\mu^{\can}_{1} = 0
 \end{cases}
\]
($\lmu^{\can} =\lmu_{\curve}^{\can}   = (\mu^{\can}_{1}, \mu^{\can}_{\fin})$ is the canonical measure of $\curve$).

Consider now a metric graph $\mgr$ with edge length function $\ell$ and the corresponding tropical curve $\curve_\ell$. By Lemma~\ref{lem:FosterApprox}, the difference measure
\[
L_1(\ell) \big ( \mu_{\mgr}^{\can} - (\lmu_{\curve_\ell}^{\can})^\ast - \frac{1}{L_1(\ell)} \lnu_{\curve_\ell}^\ast \big ) 
\]
goes to zero weakly if the metric graph $\mgr$ degenerates to $\curve$ in $\mgtropcombin{\grind{G}}$ (again, $\ast$ denotes the pull-back of measures from $\curve_\ell$ to $\mgr$).

\begin{prop} \label{prop:simple_layering_Green_functions}
Notations as in Theorem~\ref{thm:GreenFunctionConvergence}, assume that the metric graph $\mgr$ degenerates to the tropical curve $\curve$ in the boundary stratum $\mgtropcombin{\combind{(G, \pi)}}$ for an ordered partition of the form $\pi = (\pi_1, \pi_\fin)$. Let $\lgri{\ell} = (\gri{\ell, 1},  \gri{\ell, \fin})$ be the Green function on the tropical curve $\curve_\ell$ in the expansion \eqref{eq:GreenFunctionConvergence} associated to the push-out of the canonical measure on $\mgr$.  

Comparing $\gri{\ell}$ to the canonical Green function $\lgri{\curve_\ell}=\lgri{\curve_\ell}^\can = (\gri{\curve_\ell,1}, \gri{\curve_\ell, \fin})$ on $\curve_\ell$, we have 
\[
\gri{\ell, \fin}(x,y) = \gri{\curve_\ell, \fin}(x,y) +o(1),
\]
and
\begin{align*}
g_{\ell,1}(x,y) = \gri{\curve_\ell, 1} (x,y) + \frac{1}{L_1(\ell)} \Big(  \psi_{\curve_\ell} (y) - \int_{\Gamma_{\ell}^1}  \gri{\curve_\ell,1} (x,y) \, d \nu_{\curve_\ell, 1}(y) \Big)  + \frac{o(1)}{ L_1(\ell)},
\end{align*}
where $\Gamma_{\ell}^{1}$ is the first minor of $\curve_\ell$ and the $o(1)$ terms are uniform in $x$ and $y$.
\end{prop}
\begin{proof}
The claim follows from Theorem~\ref{thm:GreenFunctionConvergence} together with \eqref{eq:HigherAsymptoticsLebesgue}. Indeed,  if we integrate the sum of terms obtained from Lemma~\ref{lem:SolutionFormulaJFunctionTropical}, then all other contributions go to zero uniformly.
\end{proof}


\section{Asymptotics of height pairing on degenerating metric graphs} 
The aim of this section is to answer the following question:

\begin{question} What is the limit behavior of the height pairing $\hp{\mgr}{\cdot\,,\cdot}$ on a metric graph $\mgr$ degenerating to a tropical curve $\curve$ in the tropical  moduli space $\mgtrop{\grind{g}}$?
\end{question}
After recalling basic facts about the height pairing on metric graphs in Section~\ref{ss:HeightPairingGraphs},  we introduce in Section~\ref{ss:HeightPairingTropicalCurves} a notion of height pairing on a tropical curve $\curve$ of arbitrary sedentarity. This tropical height pairing $\lhp{\curve}{\cdot\,,\cdot}$ takes values in $\R^{r+1}$, where $r$ is the rank of the tropical curve $\curve$, and each component corresponds to a height pairing on the respective graded minor.

In Section~\ref{ss:HeightPairingDegeneration}, we then use this new notion to describe the asymptotics of the height pairing $\hp{\mgr}{\cdot \,, \cdot}$ on degenerating metric graphs.

\subsection{Height pairing on metric graphs} \label{ss:HeightPairingGraphs}
We begin by recalling basic definitions and properties in connection with the {\em height pairing} on compact metric graphs.

\smallskip

Let $G = (V, E)$ be a finite graph together with an edge length function $l \colon E \to (0, + \infty)$. Let $D \in \Div^0(\mgr)$ be a degree zero $\R$-divisor on the associated metric graph $\mgr$, that is,
\[
D = \sum_{x \in \mgr} D(x) (x)
\]
is a formal combination of finitely many points $x$ on $\mgr$, with real coefficients $D(x)$, $x \in \mgr$ such that $\sum_x D(x) = 0$. 
Viewing $D$ as a discrete measure of zero total mass on $\mgr$, we consider the Poisson equation
 \begin{equation} \label{eq:LapDiv}
	\Delta f = D
\end{equation}
on the metric graph $\mgr$, which has a (piecewise affine) solution $f \colon \mgr \to \R$, unique up to an additive constant. Identifying an edge $e\in E$ with an interval $\Ical_e = [0, l(e)]$ parametrized by $x_e$, and setting $f_e = f\rest{\Ical_e}$, it follows that the derivative 
\begin{equation} \label{eq:ContinuousDifferential}
f_e'(x_e) = \frac{df_e}{dx_e}(x_e), \qquad x_e \in \Ical_e,
\end{equation}
is a (possibly discontinuous) piecewise constant function on $e$ and only depends on the divisor $D \in \Div^0(\mgr)$. The {\em height pairing} of two degree zero divisors $D_1, D_2 \in \Div^0(\mgr)$ is defined as
\begin{equation} \label{eq:DefPairing}
	\hp{\mgr}{D_1, D_2} := \int_\mgr f_1' \cdot f_2' := \sum_{e \in E} \int_{\Ical_e}  f_{1, e}'(x_e)  f_{2, e}' (x_e) \, dx_e
\end{equation}
with $f_i$ a solution to the Poisson equation for $D_i$, $i=1,2$.

If the underlying combinatorial graph $G$ is understood from the context, we sometimes denote the height pairing by $\hp{l}{D_1, D_2}$.

\medskip

In the special case that the divisor $D$ is supported on the vertices of $G$, that is, $D$ belongs to $\Div^0(G) \subseteq \Div^0(\mgr)$, we can express the above in terms of one-forms on $G$ (see Section~\ref{sec:tropical_function_theory}). By the above discussion, the {\em differential} $\alpha_D$ of the solution $f \colon \mgr \to \R$ to \eqref{eq:LapDiv},
\begin{equation} \label{eq:SlopesDivisor}
	\alpha_D (e) : = df\rest V (e) = \frac{f(v) - f(u)}{l(e)}, \qquad e=uv \in \EE,
\end{equation}
only depends on the original divisor $D\in \Div^0(G)$. By definition, $\alpha_D$ is an exact one-form on the pair $(G, l)$ and so belongs to $\exact(G, l)$. The orthogonal decomposition~\eqref{eq:OrthDecomposition} implies that $\alpha_D$ is the unique exact one-form on $(G, l)$ satisfying the equation
\begin{equation}  \label{eq:BdOmega}
\partial \alpha = D, \qquad \alpha \in C^1(G, \R).
\end{equation}
It follows that the height pairing of two degree zero divisors $D_1, D_2 \in \Div^0(G)$ is equal to
\begin{equation} \label{eq:DefPairing2}
	\hp{\mgr}{D_1, D_2} = \innone{l}{\alpha_{D_1}, \alpha_{D_2}}= \sum_{e \in E} l(e)  \alpha_{D_1} (e)  \alpha_{D_2} (e).
\end{equation}

\begin{remark} \label{rem:ContinuousOneForms} Introducing the notion of \emph{one-forms on metric graphs}, one can recover an expression similar to \eqref{eq:DefPairing2} for general divisors supported on $\mgr$. In this language, \eqref{eq:ContinuousDifferential} corresponds to the differential $df$ of the function $f \colon \mgr \to \R$, which is the one-form on $\mgr$ given by $df\rest{e} = f'_e(x_e) dx_e$ on every edge $e$. The definition \eqref{eq:BdOmega} then can be written as a pairing $\hp{\mgr}{D_1, D_2} = \innone{\mgr}{d f_1 , d f_2}$. However, since we do not need this perspective in the sequel, we omit a discussion of one-forms on metric graphs in order to shorten the exposition.
\end{remark}

We also recall another expression for the height pairing using the graph period matrix $\rmM_l$.
\begin{prop}[Hodge theoretic description of the height pairing] For any degree zero divisor $D \in \Div^0(G)$, consider the equation \eqref{eq:BdOmega}.  If the one-form $\alpha \in C^1(G, \R)$ solves \eqref{eq:BdOmega}, then
\begin{equation} \label{eq:HeightPairingAlternative}
 \hp{l}{D, D}= \| \alpha \|^2_l - \alpha^T P_l^T \rmM_l^{-1} P_l \alpha,
\end{equation}
where $\rmM_l$ is the period matrix, and $P_l$ is the matrix defined right before Proposition~\ref{prop:OrthoProjection}.
\end{prop}

\begin{proof} We give the proof for the sake of completeness. By definition, the exact one-form $\alpha_D$ is a solution of \eqref{eq:BdOmega}. Moreover, for any other solution $\alpha$ of \eqref{eq:BdOmega}, the difference $\alpha - \alpha_D$ is a harmonic one-form on $G$. By the orthogonality of the decomposition \eqref{eq:OrthDecomposition}, we infer that
\[
	\alpha_D = \alpha - \projhar(\alpha),
\]
where $\projhar$ denotes the orthogonal projection from $C^1(G, \R)$ to $\Omega^1(G)$. 
It remains to apply Proposition~\ref{prop:OrthoProjection} and take into account once again the orthogonal decomposition \eqref{eq:OrthDecomposition}.
\end{proof}

Finally, recall also that the $\jvide$-function of the metric graph $\mgr$ (see Section~\ref{ss:LaplacianMetricGraph}) can be expressed in terms of the height pairing. More precisely, the following equality holds
\begin{equation} \label{eq:MGHPvsJ}
\jfunci{p \tiret q, x}(y) =  \hp{\mgr}{p-q, \, y -x}
\end{equation}
for all points $p,y,x,y$ on $\mgr$. Indeed, we have
\[
  \hp{\mgr}{p-q, \, y -x} = \sum_{e \in E} \int_{\Ical_e} f_{p\tiret q, e}' f_{y\tiret x, e}' = \int_\mgr  f_{p\tiret q} \, \Delta f_{y\tiret x} =  f_{p\tiret q}(y) - f_{p\tiret q}(x) = \jfunci{p \tiret q, x}(y), 
\]
where the second equality is obtained from integration by parts. Here, $f_{p\tiret q}$ and $f_{y\tiret x}$ are the solutions to the Poisson equation for divisors $p-q$ and $y-x$, respectively.

\subsection{Height pairing on tropical curves} \label{ss:HeightPairingTropicalCurves}

We now come to the definition of height pairing on a more general tropical curve. Suppose $\curve$ is a tropical curve of rank $r$ with underlying combinatorial graph $G$, ordered partition $\pi = (\pi_1, \dots, \pi_r, \pi_\fin)$ and edge length function $l$ with restriction $l_j$, $j\in [r] \cup\{\fin\}$, that we suppose normalized on the infinitary layers $\pi_j$, $j\in [r]$. Let $\Gamma$ be the associated metric graph.

\smallskip

Consider a degree zero divisor $D \in \Div^0(\curve)$ on $\curve$. By this, we mean a formal combination 
\[ 
D = \sum_{x \in \Gamma} D(x) (x)
\]
of finitely many points $x$ on $\Gamma$, with real coefficients $D(x)$, $x \in \Gamma$ such that $\sum_x D(x) = 0$. As before, we identify each divisor with a point measure, and we sometimes use the same letter for these two objects.

As discussed in Section~\ref{ss:LayeredMeasures}, the degree zero divisor $D \in \Div^0(\curve)$ induces a layered measure $\lmu_D$ of mass zero on $\curve$ (which we sometimes denote by $D$ as well). The $j$-th piece $\mu_j$ of $\lmu_D = (\mu_1, \dots, \mu_r, \mu_\fin)$ is given by
\begin{align*}
&\mu_j=  \sum_{v \in V^j} \Big ( \sum_{u \in \kappa_j^{-1} (v)} D(u) \Big ) \cdot \delta_v  + \sum_{e \in \pi_j} \sum_{x \in e \setminus V}  D(x) \cdot \delta_x, &&j =1,\dots,r, \fin,
\end{align*}
where $V^j=V(\grm{\pi}{j}(G))$ and $\kappa_j \colon V \to V^j$ is the projection map on vertices. In particular, $\mu_j$ is a point measure (divisor) on the $j$-th graded minor $\Gamma^j$.

\medskip

Analogous to the definition of the height pairing on metric graphs (see Section~\ref{ss:HeightPairingGraphs}), we consider the Poisson equation
\begin{equation} \label{eq:TropicalPairingEquation}
	\Deltatrop (\lf) = D
\end{equation}
on the tropical curve $\curve$. Since $D$ is a measure of total mass zero on $\curve$, \eqref{eq:TropicalPairingEquation} has a (piecewise linear) solution $\lf = (f_1, \dots, f_r, f_\fin)$. Any solution $\lf$ defines a collection of degree zero divisors $D^j$ on the graded minors $\Gamma^j$ of $\curve$, given by
\begin{align} \label{eq:GraphHPInducedDivisor}
D^j = \Delta_{\Gamma^j} (f_j) = \mu_j - \sum_{k < j} \divind{k}{j}(f_k), &&j =1,\dots,r, \fin,
\end{align}
where again we identify divisors with point measures. The divisors $D^j$, $j =1,\dots,r, \fin$, depends only on the divisor $D \in \Div^0(\curve)$, and not the choice of the solution $\lf$ to \eqref{eq:TropicalPairingEquation}.

\smallskip

With these preparations, we can now state the definition of the height pairing on tropical curves. In the following, we view $\R^{r+1}$ as the product of $\R^{[r]}$ and $\R^{\fin}$. By convention, in case that $\pi_\fin$ is empty, the value taken by the functions in $\R^\fin$ is zero.
\begin{defn}[Height pairing on a tropical curve]\rm
Notations as above, the {\em height pairing} on a tropical curve $\curve$ of rank $r$ is the $\R^{r+1}$-valued bilinear map 
\[
\begin{array}{cccc}
 \lhp{\curve}{\cdot\,, \cdot} \colon &\Div^0(\curve) \times \Div^0(\curve)  &\longrightarrow & \R^{r+1} = \R^{[r]}\times \R^\fin
\end{array}
\]
whose components are defined by
\[
	 \lhp{\curve, j}{D, D}  =  \hp{\Gamma^{\grind{j}}}{D^j, D^j}    = \hp{l_{\grind{j}}}{D^j, D^j}  \]
for  $j\in[r]\cup\{\fin\}$, and for any degree zero divisor $D \in \Div^0(\curve)$.
\end{defn}

\begin{remark} If the divisor is supported on the vertex set $V$, that is, $D$ belongs to $\Div^0(G) \subseteq \Div^0(\curve)$, then the height pairing can be expressed in terms of one-forms on $\curve$. The differential $\ld(\lf)$ of a solution $\lf$ to \eqref{eq:TropicalPairingEquation} is an exact one-form on $\curve$, which we denote by $\lalpha_D = (\alpha_{D}^1, \dots, \alpha_{D}^r, \alpha_{D}^{\fin})$. Here, the associated differential $\alpha_D^j$, $j\in [r]\cup\{\fin\}$, is taken in the $j$-th graded minor $\Gamma^j$, and is given by 
\begin{equation} \label{eq:SlopesDivisor}
	\alpha_D^j := d f_j, \qquad j=1, \dots, r,\fin.
\end{equation}
In terms of the one-form $\lalpha_D$, the components of the height pairing can be written as
\[
	 \lhp{\curve, j}{D, D} =  \innone{l_{\grind{j}}}{\alpha_D^j , \alpha_D^j}= \sum_{e \in \pi_j} l_j(e)  |\alpha_D^j(e)|^2, \qquad j\in[r] \cup\{\fin\}.
\]
 The latter representation is the analog of \eqref{eq:DefPairing2}. Using the notion of one-forms on metric graphs, one can derive a similar expression for general divisors on the tropical curve $\curve$, however, we omit a discussion (see Remark~\ref{rem:ContinuousOneForms}).
\end{remark}

Recall from Section~\ref{sec:tropical_moduli} that for any graph $G=(V,E)$   and for any ordered partition $\pi = (\pi_\infty, \pi_\fin)$, $\pi_\infty=(\pi_1, \dots, \pi_r)$ of rank $r$ in $\Piall(E)$, we have a projection map 
\[\mgtrop{\combind{(G,\subface \pi)}} \to \mgtrop{\combind{(G, \pi)}}.\] 
This projection map on the open part $\mggraph{\grind{G}}$, consisting of metric graphs of combinatorial type $G$, sends a metric graph $\mgr \in \mggraph{\grind{G}}$ corresponding to the pair $(G,\ell)$ with an edge length function $\ell\colon E \to (0, + \infty)$, to the  tropical curve $\curve_\ell$, the conformal equivalence class of $\mgr$. The canonical element representing $\curve_\ell$ is the layered metric graph associated to the triple $(G, \pi, \tilde \ell)$ obtained by taking the graph $G$, the ordered partition $\pi$ and the edge lengths $\tilde \ell$ given by the normalization of the restrictions $\ell_j=\ell\rest{\pi_j}$, that is,
\[
	\tilde \ell_{j} =\tilde\ell\rest{\pi_j} =  \frac{\ell\rest{\pi_j}}{L_j(\ell)}, \qquad j=1, \dots, r,
\]
and 
$\tilde \ell_{\fin} = \ell\rest{\pi_\fin}$. (Here, we recall, $L_j(\ell) = \sum_{e\in \pi_j} \ell(e)$.)

\subsection{Degeneration of the height pairing} \label{ss:HeightPairingDegeneration}

We are now in position to state the main result of this section.

Let $\curve$  be a tropical curve of rank $r$ with underlying graph $G=(V, E)$, ordered partition $\pi = (\pi_\infty=(\pi_1, \dots, \pi_r), \pi_\fin)$ and edge length function $l$. Suppose that $\mgr$ is a metric graph with underlying graph $G$ and edge length function $\ell$, and consider the corresponding tropical curve  $\curve_\ell$  in $\mgtrop{\combind{(G, \pi)}}$.

Let $D \in \Div^0(G)$ be a degree zero divisor on the combinatorial graph $G$. 

\begin{thm}\label{thm:HeightPairingGraphs} Notations as above, as the metric graph $\mgr$ degenerates to $\curve$ (for the respective points in $\mgtrop{\grind{G}}$), the height pairing has the asymptotic expansion 
\begin{equation}  \label{eq:HeightPairingGraphAsymptotics}
\langle D, D \rangle_{\mgr} = \sum_{j=1, \dots, r, \fin} L_j(\ell) \lhp{\curve_\ell,j}{D, D} + O \Bigl ( \max_{j \in \{1, \dots, r, \fin\}} \frac{ L_{j+1}^2}{ L_j} (\ell) \Bigr ).
\end{equation}
Moreover, the graded pieces of the height pairing converge, that is,
\begin{equation} \label{eq:HeightPairingGraphPieces}
 \lhp{\curve_\ell,j}{D, D} = \lhp{\curve,j}{D, D} + o(1)
\end{equation}
for all $j = 1,\dots, r, \fin$.
\end{thm}

We will conclude Theorem~\ref{thm:HeightPairingGraphs} from our results on the degeneration of the solutions of the Poisson equation. In order to apply this procedure, we first need to establish a connection between the height pairing and the Poisson equation on tropical curves. The following theorem relates the height pairing to the tropical $j$-function (see Section~\ref{ss:AsymptoticsJFunc}) and is the precise analog of \eqref{eq:MGHPvsJ}

\begin{thm}[Tropical height pairing and $j$-function] \label{thm:TropicalHPvsJFunction} Let $\Gamma$ be the layered metric graph representing the tropical curve $\curve$ with normalized edge lengths in each infinity layer. Let $p,x,y$ be points on $\Gamma$. Then the tropical height pairing
\[
y \in \Gamma \mapsto \lhp{\curve}{p-q \,, y-x} \in \R^{r+1}
\]
coincides with the tropical $j$-function $\ljfuncbis{p\tiret q, x}$, that is, it is the solution of the tropical Poisson equation
\begin{equation} \label{eq:TropicalHPvsJFunction}
\begin{cases}
 \Deltatrop (\lf) = \bm{\delta}_p - \bm{\delta}_q, \\[1 mm]
 \lf \text{ is harmonically arranged}, \\[1mm]
\int_\curve \lf \, d\bm{\delta}_x= 0.
 \end{cases}
\end{equation}
More precisely, the functions $f_i \colon \Gamma \to \R$, $i = 1,\dots, r, \fin$, given by
\begin{align*}
&f_i(y) :=  \lhp{\curve, i}{p-q, y-x},\qquad y \in \Gamma,
\end{align*}
are precisely the components of the solution $\lf$ (viewed as functions on the space $\Gamma$).
\end{thm}

Denote by $D^k_{p \tiret q}$ and $D^k_{y \tiret x}$ the degree zero divisors on $\Gamma^k$, $k \in [r] \cup\{\fin\}$, induced from $D_{p\tiret q} = p-q$ and $D_{y\tiret x} = y-x$ by \eqref{eq:GraphHPInducedDivisor}. We will reduce the theorem to the situation where we consider vertices $u$ and $v$ which form extremities of an edge. 

Let thus $e=uv$ be an (oriented) edge of $G$ in $\E$ with $\head_e =v$ and $\tail_e =u$. 
We  provide an explicit sum-product expression for the divisor $D^k_{v \tiret u}$, in terms of an admissible basis $\gamma_1, \dots, \gamma_h$ of the cycle space $H^1(G, \Z)$ and the matrices $A_p(\curve)$ introduced in Section~\ref{sec:explicit_extension_form}. 

For $j\in[r]\cup\{\fin\}$, we denote $\rmM_{j}$ and $P_{j}$ the (period) matrices from Proposition~\ref{prop:OrthoProjection} for the graded minor $(\grm{\pi}j(G), l_j)$, and let $\proj_j \colon C^1(G, \R) \to C^1(\gr_\pi^j(G), \R) $ be the natural contraction map.

In the following we denote by $\chi_e$ the one-form in $C^1(G, \R)$ which takes values $\chi_e(e') = \delta_{e,e'}$ with $\delta_{e,e'} =\pm1$ depending on whether $e'=e$ or $\bar e$, respectively, and otherwise, $\delta_{e,e'}=0$.
\begin{prop}[Sum-product description of the divisor $D^k_{u\tiret v}$]\label{prop:sum-product-Dk} Assume the edge $e = uv$ belongs to $\pi_j$, $j\in[r]\cup\{\fin\}$. Then, we have $D^k_{v\tiret u} =0$ for $k<j$, and $D^j_{v\tiret u} = \proj_j(v) -\proj_j(u)$. For all $k>j$, we have 
\[
D^k_{v\tiret u} = - \partial \Big (\sum_{j \le i < k} \, \sum_{p \in \mathcal{P}_{j i}}\,  \sum_{n \in J_\pi^{k}}  (A_p(\curve)^T P_{j} \cdot \chi_e) (n) \,\proj_k(\gamma_n)  \Big ). 
\]
\end{prop}
\begin{proof} The first and second assertion follow from the fact that $\proj_k(u) = \proj_k(v)$ for $k<j$. We need the last statement for $k>j$. Let $\omega_i$, $i\ge j$ be the exact one-form on $(\grm{\pi}i(G),l_i)$ with $\partial \omega_i = D^i_{v \tiret u}$. 
 
By Proposition~\ref{prop:OrthoProjection}, for $i=j$, we have 
\[
\omega_j = \chi_{e} - \sum_{n \in J_\pi^j} (\rmM_j^{-1} P_{j} \cdot \chi_e)(n) \, \kappa_j(\gamma_n).
\]
 We obtain the expression
\[
D^{j+1}_{v\tiret u} =  - \partial  \Big ( \sum_{n \in J_\pi^{j}}  (\rmM_{j}^{-1} P_{j} \cdot \chi_e)(n) \, \proj_{j+1}(\gamma_n)  \Big ).
\]
Applying now Proposition~\ref{prop:OrthoProjection} this time for the graded minor $(\grm{\pi}{j+1}(G), l_{j+1})$, we infer that
\[
\omega_{j+1} = - \sum_{n \in J_\pi^j}  (\rmM_{j}^{-1} P_{j} \cdot \chi_e)(n) \, \proj_{j+1}(\gamma_n) + \sum_{n \in J_\pi^{j+1}}  (\rmM_{j+1}^{-1} T_{j+1, j}\rmM_{j}^{-1} P_{j} \cdot \chi_e)(n) \, \kappa_{j+1}(\gamma_n).
\]
Here $T_{j+1, j}$ is the $h_\pi^{j+1} \times h_\pi^j$-matrix \eqref{eq:TMatrix} appearing in the definition of the matrices $A_p(\curve)$ (see \eqref{eq:AsmyptoticsApTropical}). Proceeding by induction and using the notation introduced in \eqref{eq:AsmyptoticsApTropical}, we conclude that for any $k > j$, we have
\[
D^k_{v\tiret u} = - \partial \Big (\sum_{j \le i < k} \, \sum_{p \in \mathcal{P}_{j i}}\,  \sum_{n \in J_\pi^{i}}  (A_p(\curve)^\transpose P_{j} \cdot \chi_e) (n) \,\proj_k(\gamma_n)  \Big ) 
\]
and the corresponding exact one-form $\omega_{k}$ is given by
\[
\omega_{k} = - \sum_{j \le i \le k} \, \sum_{p \in \mathcal{P}_{j i}} \, \sum_{n \in J_\pi^{i}}  (A_p(\curve)^\transpose P_{j} \cdot \chi_e) (n) \,\proj_k(\gamma_n).
\]
\end{proof}

\begin{proof}[Proof of Theorem~\ref{thm:HeightPairingGraphs}]
We use the above notations and denote by $D^k_{p \tiret q}$ and $D^k_{y \tiret x}$ the degree zero divisors on $\Gamma^k$, $k=1,\dots, r,\fin$, induced from $D_{p\tiret q} = p-q$ and $D_{y\tiret x} = y-x$ by \eqref{eq:GraphHPInducedDivisor}. By the definition of the tropical height pairing,
\begin{equation} \label{eq:HPSpelledOut}
f_i (y) = \hp{\Gamma^{\grind{i}}}{D^i_{p \tiret q}, \, D^i_{y \tiret x}}
\end{equation}
for all $i=1,\dots, r, \fin$ and $y \in \Gamma$.  The last condition in \eqref{eq:TropicalHPvsJFunction} is trivially satisfied, since
\[
f_i(x) = \hp{\Gamma^{\grind{i}}}{D^i_{p \tiret q}, \, D^i_{x \tiret x}}  = 0.
\]
Suppose that $e = uv$ is an edge in the $i$-th layer $\pi_i$. Taking into account \eqref{eq:JFunctionMiddleEdges}, one easily shows that for a point $y$ on $e$,
\begin{equation} \label{eq:LinearityDivisor}
D^j_{u\tiret x} = D^j_{y\tiret x}\quad \text{for $j <i$}, \qquad \text{and} \qquad D^j_{u-x} = D^j_{y-x} + |y-u|D_{j,e}  \quad \text{for $j >i$}
\end{equation}
where $D_{j,e}$ is a certain degree zero divisor on $\Gamma^j$ and $|y-u|$ is the distance of $y$ to $u$ on the edge $e = [0, l(e)]$. The latter proves that each function $f_i \colon \Gamma \to \R$, $i=1,\dots, \fin$, is linear on edges of lower indexed layers $\pi_j$, $j <i$, and constant on edges of higher indexed layers $\pi_j$, $j >i$. In particular, $f_i \colon \Gamma \to \R$ is the pull-back $f_i = g_i^\ast$ of a function $g_i \colon \Gamma^i \to \R$.

Moreover, for two points $y,y'$ on some connected component of the $i$-th graded minor $\Gamma^i$,
\[
g_i(y) - g_i(y') = \hp{\Gamma^{\grind{i}}}{D^i_{p \tiret q}, \, y - y'}.
\]
Using \eqref{eq:MGHPvsJ} for the graded minor $\Gamma^i$, it follows that $\Delta g_i = D^i_{p \tiret q}$.

It remains to verify the second condition in \eqref{eq:TropicalHPvsJFunction}, that is, to show that each $f_i$ is lower harmonic. Fix some  $k = 1, \dots, r, \fin$ and consider the function $f_k \colon \Gamma \to \R$.  Let $\pi_j$ be another layer with $j < k$. For an edge $e = uv$ in $\pi_j$, we have to analyze
\begin{equation} \label{eq:TropicalHPvsJFunction2}
df_k(e) = \frac{f_k(v) - f_k(u)}{l_j(e)} = \frac{1}{l_j(e)} \hp{\Gamma^{\grind{k}}}{D^k_{p \tiret q}, \, D^k_{v \tiret u} }.
\end{equation}
We will use the explicit sum-product description of the divisor $D^k_{v\tiret u}$ given in Proposition~\ref{prop:sum-product-Dk}. Using the notations in the proof of that result,
  we are interested in computing
\[
df_k(e) = \frac{1}{l(e)} \hp{\Gamma^{\grind{k}}}{D^k_{p \tiret q}, \, D^k_{v \tiret u} } = (\alpha, \omega_k)_{{\Gamma^{\grind{k}}}},
\]
where $\alpha \in C^1(\gr_\pi^k(G), \R)$ is the differential of $f_k$ seen as a function on $\Gamma^k$, that is, $\alpha = d_{\Gamma^k} f_k$. For fixed $i \le k$ and a strictly increasing sequence $p \in \mathcal{P}_{ji}$, we can write
\[
\sum_{n \in J_\pi^{i}} ( A_p(\curve)^\transpose P_{j} \cdot \chi_e)(n)  \, (\alpha, \proj_k(\gamma_n))_{{\Gamma^{\grind{k}}}} =  \chi_e^\transpose P_{j}^\transpose A_p(\curve) \cdot \rmW\rest{J_\pi^{i}} =  - l(e) \sum_{n \in J_\pi^{i}} (A_p(\curve) \cdot \rmW\rest{J_\pi^{i}}) (n) \, \gamma_n(e),
\] 
in terms of the column vector $\rmW \in \R^h$ given by $\rmW(n) := (\alpha, \proj_k(\gamma_n))_{\Gamma^{\grind{k}}}$ for $n \in [h]$.  Applying the sum-product Proposition~\ref{prop:harmonic_extension_matrices} to the one-form $\alpha$ on $\Gamma^k$, we see that $f_k$ is harmonically arranged.
\end{proof}


\begin{remark}[Harmonic rearrangement using height pairing] \label{rem:TechnicalRemark}
The proof of Theorem~\ref{thm:TropicalHPvsJFunction} reveals that harmonic rearrangement can be described in terms of the divisors \eqref{eq:GraphHPInducedDivisor} appearing in the height pairing. Suppose that $f_k \colon \Gamma^k \to \R$ is a function defined on the $k$-th minor $\Gamma^k$, $k \in [r]\cup \{\fin\}$,  of a tropical curve $\curve$. Let $g_k$ be a harmonic rearrangement of $f_k$. Then the slope of $g_k$ on an edge $e = uv$ in a layer $\pi_i$, $i < k$, can be written as
\begin{equation} \label{eq:SlopeHarmonicExtensions}
\frac{g(v) - g(u)}{l_i(e)}= \frac{f_k(D^k_{v \tiret u})}{l_i(e)} := \frac{1}{l_i(e)} \sum_{x \in \Gamma^k} D^k_{v\tiret u}(x) f_k(x)
\end{equation}
where $D^k_{v\tiret u}$ is the divisor on $\Gamma^k$  induced from $D_{v\tiret u} = v-u$ by \eqref{eq:GraphHPInducedDivisor}.
\end{remark}

\begin{proof}[Proof of Theorem~\ref{thm:HeightPairingGraphs}]
By the bilinearity of the height pairing, it suffices to treat the case when $D$ is of the form $D = p-q$ for two vertices $p,q \in V$. However, in this case the claim follows from Theorem~\ref{thm:TropicalHPvsJFunction} and Theorem~\ref{thm:layered_expansion_tropical_laplacian}.
\end{proof}

\begin{remark}
Note that the degeneration of the height pairing is a consequence of our results on the Poisson equation (see Section~\ref{sec:tropical_laplacian}). In fact, one can also prove Theorem~\ref{thm:HeightPairingGraphs} by different methods, and then deduce the results in Section~\ref{sec:tropical_laplacian}. However, since our main interest lies in the Poisson equation, the first approach seemed more streamlined to us. On the other hand, for Riemann surfaces, it is more difficult to analyze the solutions to the Poisson equation directly. In this case, we will follow precisely the second strategy, and derive the asymptotics of general solutions from the degeneration of the height pairing. The results of this section will be crucial in providing the precise description of the height pairing close to the hybrid limit.  
\end{remark}

\section{Period matrices of degenerating metric graphs} \label{sec:GraphPeriodMatrices}

In the present section, we collect results on the asymptotics of period matrices of degenerating metric graphs.  Assume that $G = (V,E)$ is a fixed combinatorial graph of genus $h$ together with a fixed basis $\gamma_1, \dots, \gamma_h$ of the cycle space $H_1(G, \Z$). Recall that for any metric graph $\mgr$, arising as a metric realization of a pair $(G, \ell)$ for some edge length function $\ell \colon E \to (0, + \infty)$, the corresponding graph period matrix $\rmM_\ell \in \R^{h \times h}$ (see \eqref{eq:M_l}) is given by 
\[
\rmM_\ell (i,j) = \innone{\mgr}{\gamma_i, \gamma_j} = \sum_{e \in E} \ell(e) \gamma_i(e) \gamma_j(e), \qquad i,j=1,\dots, h.
\]
In the present section, we discuss the asymptotic behavior of $\rmM_\ell$ (and its inverse) as $\mgr$ degenerates to a tropical curve $\curve$.

\subsection{Asymptotics of the period matrix} \label{sec:AsymptoticsGraphPeriods}

More precisely, let $\pi = (\pi_\infty=(\pi_1, \dots, \pi_r), \pi_\fin)$ be a layering on the edges of $G$ and suppose $\curve$ is a tropical curve of combinatorial type $(G, \pi)$ whose class is given by the edge length function $l$, that we can assume normalized on each layer. Assume further that our fixed basis of the cycle space $H_1(G, \Z)$ is admissible with respect to $\pi$. Since the basis of  $H_1(G, \Z)$ is decomposed into $r+1$ parts, this induces a $(r+1, r+1)$ block decomposition of the graph period matrix $\rmM_\ell$ associated to any edge length function $\ell$ on $E$. In other words, the $(m,n)$-th block $(\rmM_\ell)_{mn}$ is simply the $\R^{h_\pi^m \times h_\pi^n}$ matrix given by
\begin{equation} \label{eq:BlockFormula}
(\rmM_\ell)_{mn} (i,j) = \sum_{e \in \pi_{\max\{m,n\}} \cup \dots \cup \pi_\fin} \ell(e) \gamma_i(e) \gamma_j(e) 
\end{equation}
for all indices $m, n \in \{1, \dots, r, \fin\}$ and cycles $\gamma_i$, $\gamma_j$ with $i$ and $j$ in $J_\pi^m$ and $J_\pi^n$, respectively. 

\smallskip

The key advantage of choosing an admissible basis is a convenient description of the block-by-block asymptotics using graded minors. As is clear from \eqref{eq:BlockFormula}, if $\mgr$ degenerates to $\curve$, then the corresponding blocks have asymptotic growth rates
\begin{equation} \label{eq:AsymGraphPeriod1}
(\rmM_\ell)_{mn} = O \Big ( \max\{L_m(\ell), L_{n}(\ell) \} \Big ) 
\end{equation}
for all indices $m, n \in \{1, \dots, r, \fin\}$. Here, $L_m(\ell) = L(\ell_m)$ denotes the total length of the $m$-th graded minor $\mgr^m$, $m = 1, \dots, r, \fin$ of $\mgr$ (see \eqref{eq:DefGrLength} for the definition).

\smallskip

Moreover, since our fixed basis of $H_1(G, \Z)$ is admissible, the contracted cycles $\proj_m(\gamma_i)$ with $i \in J_\pi^m$ form a basis of $H_1(\gr_\pi^m(G), \Z)$ for every $m = 1, \dots, r, \fin$. We denote by
\[
	\rmM_{\ell_m} := \rmM \big (\gr_\pi^m(G), \ell\rest{\pi_m} \big) \in \R^{h_\pi^m \times h_\pi^m}
\]
the corresponding graph period matrix \eqref{eq:M_l} for the $m$-th graded minor $\mgr^m$ of $\mgr$. In terms of this notation, the on-diagonal entry for ever $m \in \{1, \dots, r, \fin\}$ is given by
\begin{equation} \label{eq:AsymGraphPeriod2}
(\rmM_\ell)_{mm} = \rmM_{\ell_m} + O(L_{m+1}(\ell))  
\end{equation}
as $\mgr$ degenerates to $\curve$. Moreover, taking into account the topology on $\mgtrop{\grind{g}}$,
\begin{equation} \label{eq:AsymGraphPeriod4}
 \rmM_{\ell_m}  = L_m( \ell ) \Big( \rmM_{l_m} + o(1) \Big )
\end{equation}
as $\mgr$ degenerates to $\curve$, where
\begin{equation} \label{eq:graded_matrix}
\rmM_{l_m} :=\rmM \big (\gr_\pi^m(G), l \rest{\pi_m} \big) \, \in \R^{h^{m}_\pi \times h^{m}_\pi} 
\end{equation}
is the graph period matrix \eqref{eq:M_l} of the $m$-th graded minor $\Gamma^m$ of $\curve$, again for the basis of $H_1(\grm{\pi}{m}(G), \Z)$ given by the projections $\proj_m(\gamma_i)$, $i \in J^m_\pi$.

\subsection{Asymptotics of the inverse period matrix}
Our next step is to obtain the asymptotics of the inverse matrix $\rmM_\ell^{-1}$. In order to formulate our results, we first need to introduce to a few notions. As above, let $\curve = (G, \pi, l)$ be a tropical curve with ordered partition $\pi = (\pi_1, \dots, \pi_r, \pi_\fin)$ and $\mgr = (G, \ell)$ a metric graph.

\smallskip

For every pair of indices $m, n \in \{1, \dots, r, \fin \}$ with $m \le n$, denote by $\mathcal{P}_{mn}$ the set of (strictly) increasing sequences from $m$ to $n$,
\[
\mathcal{P}_{mn} = \Bigl \{\, p = (i_0, \dots, i_{|p|}) \in \{1, \dots, r, \fin \}^{|p|} \, \, \st \, \, i_0 = m, \, i_{|p|} = n \text{ and }  i_0 < i_1 < \dots < i_{|p|} \,\Bigr\}.
\]
If $m = n$, then we simply have $\mathcal{P}_{mm} = \bigl\{ (m) \bigr\}$. For each $p \in \mathcal{P}_{mn} $, define the matrix $A_p(\mgr) \in \R^{h_\pi^m \times h_\pi^n}$ by
\begin{align} \label{eq:MatrixApGraph}
A_p(\mgr)  := (-1)^{|p|}  \,   \rmM({\ell_m})^{-1} \,  \prod_{k=1}^{|p|} (\rmM(\ell))_{i_{k-1}, \, i_k} \, \,  \rmM({\ell_{i_k}})^{-1}.
\end{align}
Again, if $m = n$, we simply obtain $A_{(m)}(\mgr) = \rmM_{\ell_m}^{-1}$.

Notice that the matrix $A_p(\mgr)$ looks similar to the matrix $A_p(\curve)$ of a tropical curve $\curve$ (see \eqref{eq:AsmyptoticsApTropical} and Proposition~\ref{prop:harmonic_extension_matrices}). In fact, the following asymptotic formula holds
\begin{equation} \label{eq:AsmyptoticsAp}
A_p(\mgr)  = \frac{1}{L_{m} (\ell) } \Big ( A_p(\curve) + o(1) \Big ), \qquad p \in \mathcal{P}_{mn},
\end{equation}
as $\mgr$ degenerates to $\curve$.

In terms of the above notions, we arrive at the following asymptotic description of the inverse graph period matrix.

\begin{thm} \label{thm:PeriodAsymptotics}
Let $\curve = (G, \pi, l)$ be a tropical curve and $\mgr = (G,\ell)$ a metric graph. Then:
\begin{itemize}
\item [(i)]  As $\mgr$ degenerates to $\curve$, the $(m,n)$-th block of $\rmM_\ell^{-1}$ has asymptotic behavior
\[
(\rmM_\ell^{-1})_{mn} = \sum_{p \in \mathcal{P}_{mn}} A_p(\mgr) + \frac{1}{L_m(\ell)} O\Big (\max_{k\le n} \frac{L_{k+1}}{L_k} (\ell) \Big ) 
\]
for every $m, n \in \{1, \dots, r, \fin\}$ with $m \le n$. In particular, 
\[
(\rmM_\ell^{-1})_{mn} =  \frac{1}{L_{\min\{m,n\}}(\ell)} (\mB_{mn} + o(1))  
\]
for some matrices $\mB_{mn} \in \R^{h^{m}_\pi \times h^{n}_\pi}$ and every $m, n \in \{1, \dots, r, \fin\}$.
\item[(ii)] The on-diagonal limit matrices $\mB_{m m }$ are given by
\[
	\mB_{m m } =  \rmM_{l_m}^{-1}, \qquad m \in \{1, \dots, r, \fin\},
\]
with $\rmM_{l_m} \in \R^{h^{m}_\pi \times h^{m}_\pi}$ the matrices of the graded minors of $\curve$ (see \eqref{eq:graded_matrix}). 
\item [(ii)] The off-diagonal matrix $\mB_{m n }$ are given
\[
\mB_{m n } = \sum_{p \in \mathcal{P}_{mn}} A_p(\curve)
\]
for every $m <n$. Moreover, $\mB_{mn} = \mB_{nm}^\transpose$ if $n < m$.

\end{itemize}
\end{thm}
\begin{proof}
This follows from \eqref{eq:AsymGraphPeriod1}, \eqref{eq:AsymGraphPeriod2} and \eqref{eq:AsmyptoticsAp}  together with Lemma~\ref{lem:inverse_lemma} and Lemma~\ref{lem:InvLem2}.
\end{proof}

\subsection{The inverse lemma(s)} \label{ss:InverseLemmas}
In the present section, we collect several auxiliary statements on asymptotics of inverse block matrices. We have the following basic result on the on-diagonal blocks (see \cite[Lemma 9.3]{AN}).

\begin{lem} \label{lem:inverse_lemma} Let $X$ be a topological space, $ \thy \in X$ a fixed point and $y_1, \dots, y_r\colon X \setminus \{\thy\} \to \C$  functions such that
\begin{equation} \label{eq:inverse_lemma_hyp}
	\lim_{t \to \thy} \frac{y_{j+1}(t)}{y_j(t)} = 0, \qquad j = 1, \dots, r-1.
\end{equation}
Suppose $\rmM \colon  X \setminus \{\thy\} \to \C^{N \times N}$ is a matrix-valued function. Assume that $\rmM(t)$ has an $(r,r)$ block decomposition (here, $N = \sum_{k=1}^r N_k$ for some $N_k \in \N$)
\[
	\rmM(t) = \Big(\rmM_{mn} (t) \Big )_{1 \le m,n \le r}
\]
where, as $t$ goes to $\thy$ in $X$, the blocks are asymptotically given by 
\[
	\rmM_{mn}(t) = \mA_{kl}(t)  + O\Big ( \frac{y_{\max\{m,n \} +1}(t)}{y_{\max\{m,n\}}(t)}\Big ) 
\]
for matrix-valued functions ${\mA}_{mn}  \colon  X \setminus \{\thy\} \to \R^{N_m \times N_n}$ such that
\[
	\lim_{t \to \thy} \, \frac{1}{y_{\max\{m,n \}}(t)} {\mA}_{mn}(t)  = \widehat{\mA}_{mn}
\]
exists for all $m,n$ and all diagonal limit matrices $\widehat{\mA}_{mm}$, $m=1,\dots, r$ are invertible.

Then, as $t$ goes to $\thy$ in $X$, the inverse $\rmM(t)^{-1}$ has the $(r,r)$ block decomposition
\[
	\rmM(t)^{-1} = \Big( y_{\min\{m,n\}}(t)^{-1} \big (\mB_{mn} + o(1) \big) \Big)_{1 \le m,n \le r}
\]
for some matrices $\mB_{mn} \in \C^{N_m \times N_n}$, and the diagonal terms are given by
\[
	\mB_{mm} = \widehat{\mA}_{mm}^{\, -1}, \qquad m=1, \dots r.
\]
\end{lem}

In our applications, we need information on the off-diagonal entries and hence a refined version of Lemma~\ref{lem:inverse_lemma}. The following lemma suffices for our purpose.

\begin{lem} \label{lem:InvLem2}
Let $\rmM \colon  X \setminus \{\thy\} \to \C^{N \times N}$ be a matrix-valued function satisfying the assumptions in Lemma~\ref{lem:inverse_lemma}. Notations as above, for every $m \le n$, the corresponding block of the inverse matrix is of the form
\[
(\rmM(t)^{-1})_{mn} = \sum_{p \in \mathcal{P}_{mn}} A_p(t) + \frac{1}{y_m(t) } O\Big (\max_{k\le n} \frac{y_{k+1}}{y_k} (t) \Big ) ,
\]
where, as before, 
\[
\mathcal{P}_{mn} = \Bigl \{\, p = (i_0, \dots, i_{|p|}) \in \N^{|p|} \,\, \st  \, \, i_0 = m, \, i_{|p|} = n \text{ and }  i_0 < i_1 < \dots < i_{|p|} \,\Bigr\}
\]
and for each $p \in \mathcal{P}_{mn} $, the matrix $A_p\colon X \setminus \{\thy\} \to \C^{N_m \times N_n} $ is given by
\begin{align*}
A_p(t)  = (-1)^{|p|}  \,   \mA_{mm}^{-1}(t) \,   \left ( \prod_{k=1}^{|p|} \rmM_{i_{k-1} \, i_k}(t) \,  \mA_{i_k i_k}^{-1}(t)  \right ). 
\end{align*}
\end{lem}

Before giving the proof of Lemma~\ref{lem:InvLem2}, we first obtain another auxiliary statement. Notations as above, assume that $\rmM \colon  X \setminus \{\thy\} \to \C^{N \times N}$ is a matrix-valued function with an $(r,r)$-block decomposition satisfying the assumptions of Lemma~\ref{lem:inverse_lemma}. Then for any $j=1, \dots, r$, we denote by 
\begin{equation} \label{eq:InductiveMatrices}
\rmM_j(t) := \big ( \rmM_{m, n}(t) \big)_{j \le m,n  \le r} 
\end{equation}
the matrix-valued function $\rmM_j \colon  X \setminus \{\thy\} \to \C^{N \times N}$ obtained by taking the lower right blocks of $M(t)$ starting from index $j$. Notice that each $\rmM_j(t)$ is a square matrix of dimension $N_\pi^j = \sum_{j \le k } N_k$ and $\rmM_1(t) = \rmM(t)$. Moreover, we denote by 
\[
Q_j(t) := \big ( \rmM_{j, k} \big)_{j< k} = \Big ( \begin{matrix} \rmM_{j, j+1}(t) & \rmM_{j, j+2}(t) & \dots & \rmM_{n, r}(t) \end{matrix} \Big ) \in \R^{N_j \times N_\pi^{j+1}}.
\]

We have the following results on the inverses of the matrices $\rmM_j$, $j= 1, \dots, r$.
\begin{lem} \label{lem:InvLem2Preparation}
Let $\rmM \colon  X \setminus \{\thy\} \to \C^{N \times N}$ be a matrix-valued function satisfying the assumptions in Lemma~\ref{lem:inverse_lemma}. Then, for every $j  =1, \dots, r$ the inverse matrix is given by
\begin{equation}\label{eq:InvLem2Preparation}
\rmM_{j}(t)^{-1} = 
 \left( \begin{matrix}\quad  \quad  ( \id + O  \big ( \frac{y_{j+1}}{y_j}   \big)  ) \mA_{jj}^{-1} \quad \quad   & - (\id +  O \big ( \frac{1}{y_j}   \big ) )  \,  \mA_{jj}^{-1} \, Q_j   \, \rmM_{{j+1}}^{-1} \,  (\id +  O \big ( \frac{y_{j+1}}{y_j} \big )) \\
 \quad\quad   O\Big(  \frac{y_{j+1}}{y_j}  \Big)  \quad \quad  &  \rmM_{{j+1}}^{-1}  (\id +  O \big ( \frac{y_{j+1}}{y_j} \big )) \end{matrix} \right ) 
\end{equation}
as $t$ goes to $\thy$ in $X$. 
\end{lem}
\begin{proof}
Recall first that the inverse of a symmetric $2 \times 2$- block matrix is expressed in terms of the Schur complement, that is,
\[
\left( \begin{matrix} A & B \\ B^T & D \end{matrix} \right )^{-1} = \left( \begin{matrix}  A^{-1} + A^{-1} B S^{-1} B^T A^{-1} & - A^{-1} B S^{-1}   \\ - S^{-1} B^T A^{-1}  & S^{-1} \end{matrix} \right )
\]
where $S =   D - B^T A^{-1} B$. We prove the claim by applying the Schur complement formula to the matrix
\[
\rmM_{j}(t) = 
 \left( \begin{matrix} \mA_{jj}(t) + O(y_{j+1}(t)) & Q_j(t) \\ Q_j^T (t) &  \rmM_{j+1}(t) \end{matrix} \right ) = \left( \begin{matrix} A & B \\ B^T & D \end{matrix} \right ).
\]
To proceed further, we first need the asymptotics of the Schur complement matrix $S$. Since $\mA_{jj} ^{-1} = O(1/y_j)$m it follows that
\begin{align*}
S  &= \rmM_{{j+1}} - Q_j^T  \, \big ( \mA_{jj}  + O(y_{j+1})  \big )^{-1}  \, Q_j  \\ &=  \Big [ \id -    Q_j^T  \, O(1)  \, \mA_{jj}^{-1}  \, Q_j \, \rmM_{{j+1}}^{-1}\Big ] \,  \rmM_{{j+1}}.
\end{align*}
Moreover, applying Lemma \ref{lem:inverse_lemma} to the matrix $\rmM_{{j+1}}$, we get
\begin{align*}
&\rmM_{j+1}^{-1} =  \left ( L_{\min\{m, n\}}^{-1}  \right ) _{m, n \ge j+1},
&&Q_j  = \Big (   \begin{matrix} O(y_{j+1}) & O(y_{j+2}) & \dots&& O(y_r)  \end{matrix} \Big ), 
\end{align*}
and hence
\begin{equation} \label{eq:ProductAsympt}
\rmM_{{j+1}}^{-1}(t) \, Q_j^T(t) = O(1)
\end{equation}
as $t$ goes to $\thy$ in $X$. Altogether, we conclude that
\[
S^{-1} =\rmM_{{j+1}}^{-1}  \Big [ \id + O(y_{j+1} / y_j) \Big ]^{-1}  =   \rmM_{{j+1} }^{-1} \Big [ \id + O(y_{j+1} / y_j) \Big ].
\]
and this proves the asymptotic of the lower right entries. Moreover, we clearly have that
\[
A^{-1} = \Big (\id+ O(y_{j+1}/y_j) \Big)  \mA_{jj}^{-1}
\]
and the claims concerning the lower left and upper right terms follow as well. Finally,
\[
A^{-1} B S^{-1} B^T  = O(y_{j+1}/ y_j)
\]
and this gives the missing statement on the upper left entry.
\end{proof}

\begin{proof}[Proof of Lemma~\ref{lem:InvLem2}]
We will prove the statement by induction on the number $r$ of blocks in each row/column of the matrix $\rmM(t)$. If $r=1$, then the claim is trivial since the matrix has no off-diagonal blocks. Hence we can assume that $\rmM(t)$ has $r>1$ and we have already proven the statement for matrices with $r-1$ blocks.

Consider now a fixed pair $(m, n)$ of indices with $m \le n$. In case that $m \ge 2$ for the row index $m$, 
it follows from Lemma \ref{lem:InvLem2Preparation} that, as $t$ goes to $\thy$ in $X$,
\[
(\rmM^{-1})_{mn} = (\rmM_2^{-1})_{mn} \, \big   (\id +  O ( y_{2} / y_1 )  \big),
\]
where $\rmM_2$ is the corresponding $(r-1, r-1)$ block matrix introduced in \eqref{eq:InductiveMatrices}. On the other hand,  the induction hypothesis can be applied to $\rmM_2$ and since
\[
\rmM_{ij}(t) \,  \mA_{jj}^{-1} (t)= O(y_j) O (1/y_j) = O(1) 
\]
for all $j\ge i$, we conclude that altogether $A_p = O(1/y_m)$ for all $p \in \mathcal{P}_{mn}$ and
\begin{align*}
(\rmM^{-1})_{mn} &=  \Big ( \sum_{p \in \mathcal{P}_{mn}} A_p  + \frac{1}{y_m} O (\max_{k\le n} {y_{k+1}}/{y_k} )\Big )    \big   (\id +  O ( y_{2} / y_1 )  \big) \\
&=   \sum_{p \in \mathcal{P}_{mn}} A_p  + \frac{1}{y_m} O (\max_{k\le n} {y_{k+1}}/{y_k} )   +  O(1/y_m) O( y_{2} / y_1 ) \\
&= \sum_{p \in \mathcal{P}_{mn}} A_p  + \frac{1}{y_m} O (\max_{k\le n} {y_{k+1}}/{y_k} ). 
\end{align*}
It remains to treat the case that $m=1$. Again, we will use Lemma \ref{lem:InvLem2Preparation}. For the upper left block, that is, the case $m= n=1$, we immediately see that $\mA_{1}^{-1} = O(1/y_1)$ and
\[
(\rmM^{-1})_{11} = \big (\id + O(y_2 / y_1) \big ) \mA_{11}^{-1} =  \mA_{11}^{-1} + \frac{1}{y_1} O (\max_{k\le n} {y_{k+1}}/{y_k} ).
\]
Turning to the case that $n > 1$, notice that $Q_1  \rmM_{{2}}^{-1} = O(1)$ (see \eqref{eq:ProductAsympt}) and hence the upper right entry in \eqref{eq:InvLem2Preparation} has asymptotic behavior
\begin{align*}
(\id +  O \big ( y_{2} / y_1    \big ) )  \,  \mA_{1}^{-1} \, Q_1   \, \rmM_{{2}}^{-1} \,  (\id +   O \big ( y_{2} / y_1   \big )) &=  \mA_{1}^{-1} \, Q_1   \, \rmM_{{2}}^{-1} + O(y_2 / y_1^2) \\ &=  \mA_{1}^{-1} \, Q_1   \, \rmM_{{2}}^{-1} + \frac{1}{y_1} O (\max_{k\le n} {y_{k+1}}/{y_k} )
\end{align*}
as $t$ goes to $\thy$ in $X$. On the other hand, we can apply Lemma~\ref{lem:inverse_lemma} to the matrix $\rmM_2(t)$ and conclude that the $(m,n)$-th entry of the product $ T_1   \, \rmM_{{2}}^{-1}$ is given by
\begin{align*}
\Big ( T_1   \, \rmM_{{2}}^{-1}\Big )_{mn}  = \sum_{a=2}^r \rmM_{1 a}  \, (\rmM_{{2}}^{-1})_{a n} = \sum_{a=2}^n \rmM_{1 a}  \, (\rmM_{{2}}^{-1})_{an} + O(y_{n+1} / y_{n}).
\end{align*}
Furthermore, the induction hypothesis applied to $\rmM_2(t)$ implies that
\begin{align*}
 \sum_{a=2}^n \rmM_{1 a}  \, (\rmM_{{2}}^{-1})_{an} &= \rmM_{an}  \, (\rmM_{{2}}^{-1})_{nn} +  \sum_{a=2}^{n-1} \sum_{p \in \mathcal{P}_{an}} \rmM_{1a}  \, (A_p +  \frac{1}{y_a} O( \max_{k \le r} y_{k+1}/y_k)) \\
 &= \rmM_{1n}  \, (\rmM_{{2}}^{-1})_{nn} +  \sum_{a=2}^{n-1} \sum_{p \in \mathcal{P}_{an}} \rmM_{1a}  \, A_p +   O( \max_{k \le r} y_{k+1}/y_k)).
\end{align*} 
On the other hand, we also have that
\begin{align*}
\rmM_{1n} (\rmM_{{2}}^{-1})_{nn} &= \rmM_{1n} (\id + O(y_2 / y_1)) \dots  (\id + O(y_{n+1} / y_n))   \mA_{nn}^{-1} \\
&= \rmM_{1n}  \mA_{nn}^{-1} + O( \max_{k \le r} y_{k+1}/y_k)).
\end{align*}
Hence, altogether, we get
\begin{align*}
 (\mA_{11}^{-1} \, Q_j   \, \rmM_{{2}}^{-1})_{mn} & =  \mA_{11}^{-1}  \rmM_{1n}  \mA_{nn}^{-1} +  \sum_{a=2}^{n-1} \sum_{p \in \mathcal{P}_{an}}  \mA_{11}^{-1} \rmM_{1a} A_p + \frac{1}{y_1} O( \max_{k \le r} y_{k+1}/y_k)) \\
 &= - \sum_{p \in \mathcal{P}_{1n}} A_p + \frac{1}{y_1} O( \max_{k \le r} y_{k+1}/y_k).
\end{align*}
Combining the above statements, the claim is proven.
\end{proof}


\newpage


\part{Hybrid curves} \label{part:HybridCurves}
The last part of our paper is centered around the {\em Poisson equation on hybrid curves}. We develop a function theory on hybrid curves, formulate a Poisson equation and study the behavior of its solutions, introduce hybrid Green functions, and prove our main theorem on the tameness of Arakelov Green functions. 

\section{Hybrid log maps I: towards the hybrid boundary} \label{sec:hybrid_log_I}

The  problem of {\em constructing log maps on families} plays a central role in the present paper. Loosely speaking, the problem can be formulated as follows.

Consider a stratified space $X = X^\ast \sqcup \partial_\infty X$ and another space $\rsf \to X$ lying above $X$.  In situations of interest, the points $t$ of $X^\ast$ represent classical geometric objects, the points $\thy$ of $\partial_\infty X = X \setminus X^\ast$ represent limiting geometric objects and for any $t \in X$, the fiber $\rsf_t$ in the family $\rsf$ is isomorphic to the geometric object it represents. In order to compare classical with limiting geometry, we want to select for each point $t \in X^\ast$ an appropriate point $\thy \in\partial_\infty X$, and also construct a map between the corresponding geometric objects. 

More formally, we would like to construct a suitable commutative diagram of maps:
 \begin{equation} \label{eq:LogMapPrelude}
\begin{tikzcd}
 \rsf \rest{X} \arrow[d]\arrow[r, "\Lognoind_\rsf"] & \arrow[d] \rsf \rest{\partial_\infty X}\\
X^\ast \arrow[r, "\Lognoind_X"] &\partial_\infty X
\end{tikzcd}
\end{equation}
In Part~\ref{part:TropicalCurves} of our paper, $X = \mgtropcombin{\grind{G}}$ and $X^\ast =\mggraphcombin{\grind{G}}$ is the subset representing metric graphs, and $\rsf$ is the family of tropical curves over $\mgtropcombin{\grind{G}}$. As maps on the base, we used a tropical log map $\logtrop{}$ and the stratumwise projection maps $\pr_\pi \colon \mggraphcombin{\grind{G}} \to \mgtropcombin{\combind{(G, \pi)}}
 $. The map on the family $\rsf /\mgtropcombin{\grind{G}} $ is simply given by homotecies of edges.

\smallskip
In applications to Riemann surfaces, $X = B^\hyb$ is the base of the hybrid versal deformation space, associated to a versal deformation of a stable marked Riemann surface, $X^\ast = B^\ast$ is the subset representing smooth Riemann surfaces, and $\rsf^{\hyb}$ is the versal family of hybrid curves over $B^\hyb$ (see Section~\ref{sec:hybrid_log_II} for definitions and further details). In this framework, we would like to map a point $t \in B^\ast$ (representing a smooth Riemann surface) to a point $\thy \in B^\hyb \setminus B^\ast$ (representing a hybrid curve), and find a map from the fiber $\rsf^{\hyb}_t$ over $t$ to the fiber $\rsf^{\hyb}_\thy$.

\medskip

In what follows, we will construct log maps for this setting. In the present section, we define the log map on the base $B^\hyb$, which applies locally to the more general setting, namely the hybrid spaces associated with a complex manifold and a simple normal crossing divisor in~\cite{AN}. In Section~\ref{sec:hybrid_log_II}, we then specialize to hybrid versal deformation spaces and construct the second map on the versal family $\rsf^\hyb$.

\subsection{Hybrid spaces of higher rank} We start by recalling the definition of the hybrid spaces of higher rank introduced in our paper~\cite{AN}. To a given complex manifold $B$ and a simple normal crossing divisor  $D$, we associate a hybrid space $B^{\hyb}$ by enriching the points $x$ of the divisor $D$ with additional simplicial coordinates.  Specializing to the case of a polydisc $B = \Delta^N$ with the SNC divisor given by the coordinate axes, we obtain a suitable base space for the hybrid versal deformation space $\rsf^\hyb/B^\hyb$, associated to the versal deformation of a stable marked Riemann surface (see Section~\ref{sec:hybrid_log_II} for details). By further gluing these spaces, one obtains the moduli spaces $\mghyb{\grind{g,n}}$ of hybrid curves with $n$ marked points \cite{AN}.

\smallskip

Let $N$ be a positive integer and let $B$ be an $N$-dimensional complex manifold. Let $E$ be a subset of $[N]$ and consider a divisor  $D = \bigcup_{e \in E} D_e$ which we assume simple normal crossing, that is, we require that $(D_e)_{e \in E}$ is a finite family of smooth, connected and closed submanifolds of codimension one in $B$ such that for any subset $F \subseteq E$, the intersection 
\begin{equation} 
 D_F := \bigcap_{e \in F} D_e
\end{equation}
is either empty or a smooth submanifold of codimension $|F|$ with only finitely many connected components (all of the same dimension). We call the intersection $D_F$ the \emph{stratum} associated to the subset $F$ of $E$. Furthermore, we define the \emph{inner stratum} of $F$ to be the subset
\begin{equation} \label{eq:inncomp}
{\inn D_F} := D_F \setminus \bigcup_{e \notin F} D_e = \bigl\{t \in B \,\st \ E_t = F \bigr\}.
\end{equation}
Here, for $t \in B$, we set
\begin{equation}
E_t := \bigl\{e \in E\,\st \, t \in D_e \bigr\}.
\end{equation}
Note that in particular, we recover for $F = \varnothing$, 
\[
 {\inn D_\varnothing} = B \setminus D =: B^\ast.
\]
Moreover, the inner strata \eqref{eq:inncomp} provide a partition of $B$,
\begin{equation} \label{eq:partitionB}
B = \bigsqcup_{\substack{ F \subseteq E }} {\inn D_F}.
\end{equation}

Recall from Section~\ref{subsec:OrdPart} that $\Pifs(F)$ is the set of ordered partitions with full sedentarity of a subset $F \subseteq E$. Each element $\pi \in \Pifs(F)$ is an ordered sequence $\pi = (\pi_i)_{i=1}^r$, for some $r\in \N$, consisting of non-empty, pairwise disjoint subsets $\pi_i \subseteq F$  with $\bigsqcup_{i=1}^r \pi_i = F.$ The integer $r$ is called the rank of $\pi$. In other words, $\Pifs(F)$ is the subset of $\Piall(F)$ consisting of those $\pi = (\pi_\infty, \pi_\fin)$ with $\pi_\fin =\emptyset$. Recall also that, by definition, $F = \varnothing$ has $\pi_\varnothing = \varnothing$ as its only ordered partition.

For an ordered partition $\pi = (\pi_i)_{i=1}^r$ in $\Pifs(F)$, $F \subseteq E$, the \emph{hybrid stratum} $D_\pi^\hyb$ is given by
\begin{equation} \label{eq:hybridstratum}
{D}^{\hyb}_\pi :=  \inn \sigma_{\pi} \times \inn D_\pi  =  \inn \sigma_{\pi} \times \inn D_F 
\end{equation}
where for a finite set $A$, $\sigma_A$ denotes the standard simplex in $\R^A$, $\inn \sigma_A$ its relative interior, and
\begin{align*}
\inn \sigma_\pi := \inn\sigma_{\pi_1} \times \dots \times \inn\sigma_{\pi_r},  \qquad \inn D_\pi := \inn D_F = \inn D_{\pi_1 \cup \dots \cup \pi_r}.
\end{align*}
For $\pi = \pi_\varnothing$ (i.e., the empty ordered partition of $F= \varnothing$), we set
\[
	{D}^{\hyb}_{\pi_\varnothing} := \inn D_{\varnothing} = B^\ast.
\]
The \emph{hybrid space} $B^{\hyb}$ is the disjoint union
\begin{equation} \label{eq:Bhyb}
B^{\hyb} := \bigsqcup_{\substack{  F \subseteq E }} \bigsqcup_{\substack{ \pi \in \Pifs(F)}} {D}^{\hyb}_\pi = B^\ast \sqcup \bigsqcup_{\substack{ \varnothing \subsetneq F \subseteq E }} \bigsqcup_{\substack{ \pi \in \Pifs(F)}} {D}^{\hyb}_\pi
\end{equation}
equipped with the \emph{hybrid topology}. We call  the stratum $B^\ast$ the \emph{finite part} of of $B^\hyb$; its complement $\partial_\infty B^\hyb := B^\hyb \setminus B^\ast$ is called the \emph{hybrid} \emph{boundary (at infinity)}. 

\smallskip

We refer to~\cite{AN} for the definition and discussion of the hybrid topology. In Proposition~\ref{prop:HybridTopologyViaLogMap} below, we will give an alternative characterization in terms of locally defined maps and higher rank compactifications of simplicial fans. It turns out that the hybrid topology can be viewed as a mixture of the complex topology on $B$ and the tropical topology from Section~\ref{sec:higher_rank_compactifications}.   In later applications to Riemann surfaces, we will use a refined topology, the {\em hybrid tame topology} (see Section~\ref{ss:HybridTameTopology}). The latter mixes the complex topology on $B$ and the tropical tame topology (see Section~\ref{sec:higher_rank_compactifications}).

\medskip

In the sequel, elements of $B^{\hyb}$ will be usually written as pairs $\thy = (x,t)$, meaning that $t \in {\inn D}_\pi$ and $x \in {\inn \sigma}_\pi$ for some ordered partition $\pi$ in $\Pifs := \bigsqcup_{F \subseteq E} \Pi(F)$.

\subsection{Adapted coordinates} \label{ss:AdaptedCoordinates} An \emph{adapted coordinate neighborhood} or a \emph{system of local parameters} for $D$ around a point $t \in B$ is a pair $(U,z)$ where $U$ is an open neighborhood of $t$ in $B$ and $z=(z_i)_{i=1}^N$ are local coordinates on $U$ with
\smallskip

\begin{itemize}
\item [(i)] $|z_i|< 1$ on $U$ for all $i = 1,\dots,N$, and
\smallskip

\item [(ii)] $D_e \cap U = \varnothing$ for all $e \notin E_t$, and
\smallskip

\item [(iii)] $D_e \cap U =\bigl\{s \in U\,| \; z_e(s) = 0 \bigr\}$ for all $e\in E_t$. 
\end{itemize}

In the following, we will usually fix the point $t$ and a system of local parameters $(z_i)_{i\in [N]}$. Replacing $B$ with $(U,z)$ and $E =E_t$, we can assume that $B = \Delta^N$ is a polydisc, $t= 0  \in \Delta^N$, and $D =\bigcup_{e\in E} D_e$ with $D_e =\{s\in B \,|\, z_e(s) =0\}$, for some subset $E \subseteq [N]$.

Here and below, $\Delta$ denotes the open unit disk in $\C$ and $\inn\Delta =\Delta \setminus \{0\}$ the punctured disk.

 \subsection{The map $\LOG$: hybrid spaces and tropical compactifications} \label{ss:LogMapHybridTropicalCompactification}
We adapt the set-up presented in Section~\ref{ss:PermutohedraExhaustionFoliation} and fix a choice of $n$ auxiliary functions 
\[\Phi_1, \Phi_2, \dots, \Phi_n \colon [1, +\infty) \to \R_+\]
 that we suppose 
 \begin{itemize}
 \item smooth, increasing each, with $\Phi_j(1)=1$ and $\lim_{t\to \infty}\Phi_j(t)=\infty$, $j\in [n]$, and verifying
 \smallskip
 
 \item (Log-increasing property with parameter 2) $\frac{\Phi'_{j+1}(t)}{\Phi_{j+1}(t)}\,> \, 2\,\frac{\Phi'_j(t)}{\Phi_j(t)}$, $j\in[n-1]$ , $t\in [1, +\infty)$, with 
 \[\lim_{t \to \infty} \frac{\Phi_j(t)^2}{\Phi_{j+1}(t)} =0.\]
 \end{itemize}

 \label{sec:basic_log_map} In the situation of the previous section, (i.e., $B =  \Delta^N$ is a polydisc), consider the cone $\eta := \R_{+}^E$ and its tropical compactification $\cancomp\eta^\trop$ (see Section~\ref{sec:higher_rank_compactifications}).  
 
 \smallskip

There is a natural logarithmic map from the finite part $B^\ast = \prod_{e\in E} \inn\Delta_e \times \prod_{e\in E^c}\Delta_e$ to $\eta$,
\begin{align} \label{eq:LogComplexToCone}
\begin{array}{cccc}
\LOG=\LOG_E \colon &B^\ast &\to &\eta\\
&(z_e)_{e\in [N]} &\mapsto  &(-\log\abs{z_e})_{e\in E},
\end{array}
\end{align}
sending the coordinates $z_e$, $e\in E$, to $-\log\abs{z_e}$, and forgetting the other coordinates $z_e$, for $e$ in $E^c := [N] \setminus E$. Note that the image of $\Log$ coincides with the open stratum $\R_{+}^E$ of $\cancomp\eta^\trop$ . We will extend this map to a stratumwise defined map
\begin{equation} \label{eq:LogHybridCones}
\LOG \colon  B^\hyb = \bigsqcup_{\pi} {D}^{\hyb}_\pi \to \cancomp\eta^\trop,
\end{equation}
which sends every hybrid stratum ${D}^{\hyb}_\pi$ of $B^\hyb$ onto a tropical stratum $\inn\keg_{\hat \pi}$ in $\cancomp\eta^\trop$. Recall that each ordered partition $\hat \pi = (\pi_\infty, \pi_\fin)$ of $E$ yields a stratum $\inn\keg_{\hat \pi}$ of $\cancomp\eta^\trop$ of the form
\[
\inn\keg_{\hat \pi} = \inn \sigma_{\pi_\infty} \times \inn\keg_{\pi_\fin}.
\]
Moreover, any ordered partition of full sedentarity $\pi$ of some $F\subseteq E$ naturally gives an ordered partition $\hat \pi = (\pi_\infty, \pi_\fin)$ of $E$, by setting
\begin{equation} \label{eq:HybridTropicalPartitions}
\hat \pi = (\pi_\infty, \pi_\fin), \qquad \pi_\infty = \pi, \qquad \pi_\fin = E \setminus F.
\end{equation}
On the stratum ${D}^{\hyb}_\pi$, we define the map $\LOG$ in \eqref{eq:LogHybridCones} as
\begin{align*}\begin{array}{cccc}
\LOG_\pi \colon & D_\pi^\hyb = \inn \sigma_\pi \times \inn D_F & \to & \inn\keg_{\hat \pi} = \inn \sigma_{\pi_\infty} \times \inn\keg_{\pi_\fin} \\[2mm]
	&\thy = \Big ( x, (z_e)_{e \in [N]} \Big) & \mapsto &\Big (x, (-\log\abs{z_e})_{e\in \pi_\fin} \Big )
\end{array}.
\end{align*}
Using logarithmic polar coordinates, the open unit disc $\inn\Delta$ is decomposed as $\inn\Delta = \R_{+} \times \realtor$, where $\realtor = \R / 2 \pi$, and we can write
\begin{equation} \label{eq:HybridStratumbyTropicalCone}
 D_\pi^\hyb  =  \inn \sigma_\pi \times \inn D_F  = \inn \sigma_\pi \times \inn \keg_{\pi_\fin} \times  \realtor_{\pi_\fin} \times \prod_{e \in E^c} \Delta_e = \inn\keg_{\hat \pi} \times \realtor_{\pi_\fin} \times \prod_{e \in E^c} \Delta_e
 \end{equation}
where for any subset $A \subset E$, the \emph{real torus} $\realtor_A$ is given by 
\[\realtor_A :=   (\rquot \R{2\pi\R})^A. 
\]
 In these coordinates, $\LOG_\pi$ is a simple projection map.
\begin{prop} \label{prop:ContLogMap} The map $\LOG \colon B^\hyb \to \cancomp\eta^\trop$ is continuous. It maps each stratum ${D}^{\hyb}_\pi$ of $B^\hyb$ surjectively onto the stratum $\inn\keg_{\hat \pi}$ of $\cancomp\eta^\trop$ for $\hat \pi$ given by \eqref{eq:HybridTropicalPartitions}.
\end{prop}

\begin{proof} This follows directly  from the definition of the tropical and hybrid topologies on $\cancomp\eta^\trop$ and $B^\hyb$, respectively.
\end{proof}

\subsection{The stratumwise hybrid log map $\loghyb_\pi$}
\label{sec:log-hybrid-stratumwise}
Notations as above, we proceed with the definition of suitable log maps for a polydisc $B = \Delta^N$. Recall that, when discussing tropical curves, we used two different log maps: the stratumwise projection map $\pr_{\pi}$ and the (global) tropical log map $\logtrop$  (see, e.g., Section~\ref{ss:TameFunctions}).

In what follows, we combine these tropical maps with the above map $\LOG \colon  B^\hyb  \to \cancomp\eta^\trop$, and obtain two different hybrid log maps:
\begin{itemize}
\item for each stratum $D_\pi^\hyb$, a {\em stratumwise hybrid log map} $\loghyb_\pi \colon B^\ast \to D_\pi^\hyb$ mapping the finite part $B^\ast$ to $D_\pi^\hyb$, and
\item the {\em (global) hybrid log map} $\loghyb \colon \mathscr U \to \partial_\infty B^\hyb$, which is a retraction from an open neighborhood $ \mathscr U$ of $\partial_\infty B^\hyb$ to the hybrid boundary $\partial_\infty B^\hyb$
\end{itemize}
We begin with the construction of the first map $\loghyb_\pi$. Fix a hybrid stratum $D_\pi^\hyb$, associated to an ordered partition $\pi = (\pi_j)_{j=1}^r$ of full sedentarity on a subset $F \subset E$. By separation of coordinates, the polydisc $B$ is a product of two polydiscs $B = B_F \times B_{F^c}$, where $F^c = [N] \setminus F$. Consider the logarithm map
\begin{align*}
\begin{array}{cccc}
\Log_F \colon &B^\ast_{F} &\to &\inn \eta_{F}\\[2mm]
&(z_e)_{e\in F} &\mapsto &(-\log\abs{z_e})_{e\in F}
\end{array},
\end{align*}
where $\inn \eta_{F} =\R_{>0}^F$.  The {\em stratumwise log map} $\loghyb_\pi$  is given  by
\begin{equation*}
\begin{array}{cccc}
\loghyb_\pi \colon &B^\ast &\to &D^\hyb_\pi = \inn \sigma_\pi \times \inn D_F \\[2mm]
&\big( z_F , z_{F^c} \big ) &\mapsto & \big (\pr_{\pi}\circ \Log_F(z_F), z_{F^c} \big)
\end{array}
\end{equation*}
where $\pr_{\pi} \colon \inn \eta_F \to \inn \sigma_\pi$ is the projection map to $\inn\sigma_\pi$ from \eqref{eq:StratumwiseProjection}. That is,
\[
\pr_\pi\big( (x_e)_{e \in F} \big) = \Big( [ (x_e)_{e \in \pi_1}  ], \dots, [ (x_e)_{e \in \pi_r} ]  \Big ), 
\]
where the brackets denote projectivization. 

\subsection{The map $\loghyb$ to $\partial_\infty B^\hyb$} \label{ss:GlobalHybridLogMap} Notations as in the previous section, in this section we define the \emph{hybrid log map} $\loghyb$. The latter is a retraction $\loghyb \colon \mathscr U \to \partial_\infty B^\hyb$ for an open neighborhood $\mathscr U$ of the hybrid boundary  $\partial_\infty B^\hyb$ to  $\partial_\infty B^\hyb$.

\medskip
 
As before, let $\eta =\R_{\geq 0}^E$, and consider the above map $\LOG \colon  B^\hyb  \to \cancomp\eta^\trop$ (see \eqref{eq:LogHybridCones}). Since $\Log$ is continuous (see Proposition~\ref{prop:ContLogMap}), the following set
\[
{\mathscr U} = \LOG^{-1}(\cancomp\eta^\trop \setminus \cube_\eta),
\]
is an open neighborhood of $\partial_\infty B^\hyb$ in $B^\hyb$. Here $\cube_\eta$ is the closed unit hypercube in $\eta$ (see Section~\ref{sec:tropical_log_map}). The set $\mathscr{U}$ provides the domain of definition for the hybrid log map $\loghyb$.%

\smallskip

As in the previous section, for an ordered partition $\pi=(\pi_1, \dots, \pi_r) \in\Pifs(F)$ of full sedentarity, $F\subset E$, let $\hat\pi = (\pi_\infty, \pi_\fin) \in \Piall(E)$ be the ordered partition of $E$ given by \eqref{eq:HybridTropicalPartitions}.

Using logarithmic polar coordinates, we decompose as before the hybrid stratum $D_\pi^\hyb$  as
\begin{equation} \label{eq:WriteHybridStratum} D_\pi^\hyb = \inn \sigma_\pi \times \inn D_F  = \inn \sigma_\pi \times \inn\keg_{\pi_{\fin}}\times \realtor_{\pi_\fin} \times \prod_{j\in E^c }\Delta_j.
\end{equation}

\subsubsection{Projections to the sedentary and finitary parts}

Consider the tropical compactification $\cancomp \eta^\trop$ of $\eta$. 
We define the \emph{projection-to-the-sedentarity-part} map  
\[\forget_{\sed}\colon \partial_\infty \cancomp\eta^\trop\to \sqcup_{\pi} \inn\sigma_\pi,\]
the union being taken over all ordered partitions $\pi \in \Pifs(F)$ of full sedentarity for non-empty subsets $F \subseteq E$, as the map which forgets the finitary part, that is: 

On each stratum $\inn \keg_{\hat \pi} = \inn \sigma_{\pi_\infty} \times \inn \keg_{\pi_\fin}$, for ordered partition $\hat\pi =(\pi_\infty, \pi_\fin)\in \Piall(F)$, for non-empty $F\subseteq E$, the map is given by the forgetful projection map  $\inn \sigma_{\pi_\infty} \times \inn \keg_{\pi_\fin} \to \inn \sigma_{\pi_\infty}$.  

\smallskip

Similarly, we define the \emph{projection-to-the-finitary-part} map  
\[\forget_{\fin}\colon \partial_\infty \cancomp\eta^\trop\to \sqcup_{\pi_\fin} \inn\keg_{\pi_\fin},\]
the union being taken over all subsets $\pi_\fin \subseteq E$.

\subsubsection{Definition of the hybrid log map}

With these preparations, we are now ready to define the hybrid log map $\loghyb$.  

\smallskip

First of all, on the boundary $\partial_\infty B^\hyb$, we set $\loghyb$ to be the identity.

\smallskip

Let now $\underline z=(z_e)_{e\in [N]}$ be a point in $\mathscr{U} \setminus \partial_\infty B^\hyb$, and let $\underline x := \LOG(\underline z) \in \eta$. The image $\logtrop(\underline x)$ lies on some stratum $\inn \keg_{\hat \pi}$ of $\cancomp \eta^\trop$, associated to an ordered partition  $\hat \pi = (\pi_\infty, \pi_\fin) = (\pi_1, \dots, \pi_r, \pi_\fin)$ in $\Piall(E)$. Here $\logtrop \colon \eta \setminus\cube_\eta \to \partial_\infty\cancomp \eta^\trop$ is the tropical log map defined in Section~\ref{sec:tropical_log_map}.

\smallskip

We define the image $\loghyb(\underline z)$ as the point of $\partial_\infty B^\hyb$ lying on the stratum 
\[D^\hyb_{\pi_\infty} = \inn \sigma_{\pi_\infty} \times \inn D_F  = \inn \sigma_{\pi_\infty} \times \prod_{e\in \pi_\fin}\inn \Delta_e \times \prod_{e\in E^c} \Delta_e, \qquad F = \pi_1 \cup \dots \cup \pi_r, \] 
defined by
\begin{equation}\label{eq:loghyb}
\loghyb(\underline z) := \forget_{\sed}\circ\logtrop(\underline x) \times  \Bigl(\forget_{\fin}\circ \logtrop(\underline x) , (\arg(z_e))_{e\in \pi_f}\Bigr) \times (z_e)_{e\in E^c}. \end{equation}
That is, the  coordinates of $\loghyb(\underline z)$ w.r.t. the decomposition \eqref{eq:WriteHybridStratum} are given by
\begin{align*} 
\forget_{\sed}\circ\logtrop(\underline x) \in \inn\sigma_{\pi_\infty}, \qquad \Bigl(\forget_{\fin}\circ \logtrop(\underline x), (\arg(z_e))_{e\in \pi_f}\Bigr) \in  \inn \keg_{\pi_\fin}\times \realtor_{\pi_\fin}, \qquad (z_e)_{e\in E^c} \in \prod_{e\in E^c}\Delta_e.
\end{align*}

\begin{thm}\label{thm:cont-hybrid} Notations as above, the map $\loghyb$ defines a retraction from $\mathscr U$ to the boundary at infinity $\partial_\infty B^\hyb$.
\end{thm}

 \begin{proof} We need to show that the map $\loghyb$ is continuous. 
The coordinates in the first term in \eqref{eq:loghyb} are continuous as the tropical logarithmic map is. The third term in \eqref{eq:loghyb} is obviously continuous. It remains to check that the middle term behaves continuously. This follows from the definition of the hybrid topology: the middle term is obviously continuous on each stratum of $B^\hyb$. Moreover, moving from a stratum associated to an ordered partition to the one associated to a refined partition leads to forgetting some of the coordinates of the middle term. The theorem follows. 
\end{proof}

\subsection{The tame topology} \label{ss:HybridTameTopology}
In this section, we return to the general setting of a complex manifold $B$, together with a fixed SNC divisor $D$.

We first characterize the hybrid topology on $B^\hyb$ in terms of log maps to tropical spaces and the projection map to $B$. This confirms the intuitive idea that the hybrid topology mixes complex and tropical topologies.

\begin{prop}\label{prop:HybridTopologyViaLogMap} The hybrid topology on $B^\hyb$ is the weakest topology such that
\begin{itemize}
\item the natural projection map $B^\hyb \to B$ is continuous and
\smallskip

\item for all adapted coordinate neighborhoods $(U,z)$ around any $t \in B$, the logarithmic map $\LOG \colon U^\hyb \to \cancomp\eta_t^\trop$ is continuous, for $\eta_t = \R_{\geq 0}^{E_t}$ and $\cancomp\eta_t^\trop$ the absolute canonical compactification of $\eta_t$ endowed with the tropical topology. 
\end{itemize}
\end{prop}
Here, the map $\LOG \colon U^\hyb \to \cancomp\eta_t^\trop$ is defined in \eqref{eq:LogHybridCones} (recall that $U$ is a subset of a polydisc).
\begin{proof}
The claim is a reformulation of the definition of the topology on $B^\hyb$, and can be verified in a straightforward manner. We thus omit the details.
\end{proof}

Recall that in Definition~\ref{def:TropicalTameTopology}, we introduced another topology on the higher rank compactifications of fans, which is called the (tropical) tame topology. 
Mimicking the characterization of the hybrid topology on $B^\hyb$ in Proposition~\ref{prop:HybridTopologyViaLogMap}, we introduce the hybrid analogue of the tame topology by pulling back the tame topology from the spaces $\cancomp\eta_t^\trop$. The following definition obviously extends to the $\tameclass$-tame topology for any class of functions $\tameclass$ considered in Section~\ref{sec:higher_rank_compactifications}.

\begin{defi}[Tame topology] \label{def:TameHybridTopology}
The (hybrid) {\em tame topology} on $B^\hyb$ is the weakest topology such that
\begin{itemize}
\item the natural projection map $B^\hyb \to B$ is continuous and
\smallskip

\item for all adapted coordinate neighborhoods $(U,z)$ around every $t \in B$, the logarithmic map $\LOG \colon U^\hyb \to \cancomp\eta_t^\trop$ is continuous for the tame topology on $\cancomp\eta_t^\trop$ (here $\eta_t = \R_+^{E_t}$).
\end{itemize}
\end{defi}
If a sequence $(t_n)_n$ in $ B^\hyb$ converges to some $\thy \in B^\hyb$ in the tame topology on $B^\hyb$, we also say that $t_n$ {\em converges tamely } to $\thy$ or that $(t_n)_n$ is {\em tamely convergent} to $\thy$.

\smallskip

The space $B^\hyb$ equipped with the tame topology is a second countable Hausdorff space. However, it is not compact in general. Moreover, the tame topology is finer than the original topology on $B^\hyb$.

\smallskip

In later applications, we will often consider a tamely convergent sequence of points $(t_n)_n$ in the finite part $B^\ast$. For better illustration, we describe this convergence in the next proposition.

\begin{prop}
Let $\thy \in B^{\hyb}$ and assume that $\thy=(t,x) \in D_\pi^{\hyb}$ for the ordered partition $\pi = (\pi_1, \dots, \pi_r)$ in $\Pifs(E_t)$. Suppose that $(t_n)_n$ is a sequence in $B^\ast$.

\smallskip

Then, $(t_n)_n$ converges tamely to $\thy$ if and only if the following conditions hold:
\begin{itemize}
\item [(i)] $t_n$ converges to $t$ in $B$.
\smallskip

\item [(ii)] For every $k= 1, \dots, r$  and every coordinate $e \in \pi_k$,
\begin{equation}  \label{eq:conv1} \lim_{n\to\infty}\frac{\log\abs{z_e(t_n)}}{  \sum_{f \in \pi_k}  \log\abs{z_{f}(t_n)}}=x(e).
\end{equation}
\smallskip

Moreover, for coordinates $e \in \pi_k$ and $e' \in \pi_{k'}$ with $k < k'$,
\begin{equation} \label{eq:conv2}
\lim_{n\to\infty}\frac{\bigl(\log\abs{z_{e'}(t_n)}\bigr)^2}{\log\abs{z_e(t_n)}}=0.
\end{equation}
\smallskip

Here $(U,z)$ is a fixed (or equivalently, any) adapted coordinate neighborhood of $t$.
\end{itemize}
\end{prop}

\section{Hybrid log maps II: from Riemann surfaces to hybrid curves} \label{sec:hybrid_log_II}

In this section, we complete the construction of log maps on families of hybrid curves. As outlined in the beginning of Section~\ref{sec:hybrid_log_I},  we extend the hybrid log map of the previous section to versal families of hybrid curves (see \eqref{eq:LogMapPrelude} for details).

We start by recalling basic facts on versal families and versal deformation spaces associated with stable Riemann surfaces (Section~\ref{sec:deformations}), as well as the hybrid deformation space $B^\hyb$ and the family of hybrid curves $\rsf^\hyb \to B^\hyb$  introduced in~\cite{AN}  (Section \ref{sec:hybrid_family}). We then introduce what we call \emph{adapted sets of coordinates} for versal families of Riemann surfaces (Section~\ref{sec:AdaptedCoordinatesFamily} and Section~\ref{sec:plumbing}). Using these coordinates, we construct the hybrid log maps on the hybrid family $\rsf^\hyb$ (Section~\ref{sec:LogSurfacesHybridCurves} and Section~\ref{sec:hybrid_log_map_rs}).

 \subsection{Deformations of stable curves with markings} \label{sec:deformations}

 Let  $(S_0, q_1, \dots, q_n)$ be a fixed stable curve of arithmetic genus $g$ with $n$ marked points. We allow $n=0$,  in which case, we just get a stable curve without any markings.
 
  Let $G=(V,E,\genusfunction,\marking)$ be the dual graph of $S_0$ and for each vertex $v\in V$, denote by $C_v$ the normalization of the corresponding irreducible component of $S_0$. Recall that the genus function $\genusfunction: V \to \mathbb Z_{\geq0}$ assigns to each vertex $v$ of $G$ the genus $\genusfunction(v)$ of $C_v$. By the genus formula,
\[g = \graphgenus + \sum_{v\in V} \genusfunction(v),\]
where $\graphgenus = |E| - |V|+1$ denotes the genus of the graph $G$.

Recall that the marking function $\marking: [n] \to V$ associates to each element $i \in [n]$ a vertex of the graph. More precisely, the vertex $\marking(i) \in V$ represents the component on which the $i$th marked point $q_i$ lies. The marking function $\marking$ gives rise to a counting function $\countmarking: V \to \N \cup\{0\}$, which assigns to every vertex $v\in V$ the number $\countmarking(v)$ of elements $i\in[n]$ with $\marking(i) = v$. If $n=0$, the marking function is void. The stability condition means that the inequality
\[2\genusfunction(v) - 2 +\deg(v)+\countmarking(v) >0\]
holds at every vertex $v$ of the dual graph.

\medskip

For each edge $e =uv$ of the dual graph, $S_0$ has a node $p_e$. We denote by $p^e_u$ and $p^e_v$ the corresponding points in $C_u$ and $C_v$, respectively. For each vertex $v$, let $\mathcal A_v$ be the collection of the marked points $q_j$ lying on $C_v$ and all the points $p^e_v$, for edges $e$ incident to $v$.   We get for each $v$ in $G$, the marked Riemann surface $(C_v, \mathcal A_v)$ of genus $\genusfunction(v)$ with $\deg(v)+\countmarking(v)$ marked points.

 \smallskip
 
For the marked curve  $(S_0, q_1, \dots, q_n)$,  there exists a  versal formal deformation space $\pr \colon \widehat{\rsf} \to \widehat{B}$ over a formal disc $\widehat{B}$ of dimension $N:=3g-3+n$ such that the fiber at $0$ of the family is isomorphic to $S_0$. Moreover, we have sections $\sigma_i \colon \widehat{B} \to \widehat{\rsf}$ such that $\sigma_i(0)=q_i$. 

For each edge $e\in E$, we get a local parameter $z_e$  and the corresponding formal divisor $\widehat{D}_e\subset \widehat{B}$. The divisors $\widehat D_e$, $e\in E$, meet transversally. In addition to these local parameters indexed by the edges of the graph, there are $N- |E| = 3g-3+n-|E|$ extra local parameters. These parameters correspond to those deformations of the stable marked curve which preserve the dual marked graph, that is, to deformations of the marked curves $(C_v, \mathcal A_v)$, $v\in V$. (Note that $3g-3+n-|E| = \sum_{v\in V} \bigl(3\genusfunction(v)-3+\deg(v) +\countmarking(v)\bigr)$.)

\smallskip

Using formal schemes,  the above local picture leads to analytic deformations $\rsf \to B$ of stable Riemann surfaces with $n$ markings  over a polydisc $B = \Delta^N$ of dimension $N$, where $\Delta$ is the open unit disk in $\C$. In addition, we get analytic divisors  $D_e \subset B$ given by the equations
\[	D_e = \{ z \in  \Delta^N | z_e = 0 \}\]
for $e\in E$; $D_e$ is the locus of all points $t$ in $B$ such that  the corresponding fiber $\pr^{-1}(t)$ in the family $\pr: \rsf \to B$ has a node $p_e(t)$ corresponding to $e$.  For each edge $e\in E$, we thus get a section $p_e \colon D_e \to \rsf$ of the family, $t\mapsto p_e(t)$, defined over $D_e = \{z_e=0\}$. 

\smallskip

In what follows, we denote by $B$ the base of the versal deformation of  the stable marked curve $(S_0, q_1, \dots, q_n)$, and define $B^\ast : = B \setminus \bigcup_{e \in E} D_e$ and $\rsf^\ast := \rsf\rest{B^\ast}=\pr^{-1}(B^\ast)$ to be the locus of points with smooth fibers and their corresponding family. The fiber over  $t \in B$ is denoted $\rsf_t$.

\subsection{Hybrid versal deformation space and hybrid versal curve} \label{sec:hybrid_family} Notations as in the previous section, consider the analytic versal family $\rsf \to B$ associated to a marked stable Riemann surface $(S_0, q_1,  \dots, q_n)$. Let $B^{\hyb}$ be the hybrid space obtained from $B$ and the divisor $D$ (see Section~\ref{sec:hybrid_log_I}). $B^\hyb$ will then serve as the base space for the hybrid versal deformation $\pr^\hyb\colon \rsf^\hyb \to B^\hyb$ associated to the stable Riemann surface  $(S_0, q_1, \dots, q_n)$. 

The {\em hybrid versal curve} $\rsf^\hyb$ is obtained from $\rsf$ by replacing non-smooth fibers $\rsf_t$ in the versal curve $\pr\colon \rsf \to B$ by hybrid curves. Taking into account the hybrid stratification of $B^\hyb$,
\[
B^{\hyb} = B^\ast \sqcup \bigsqcup_{\substack{ \varnothing \subsetneq F \subseteq E }} \bigsqcup_{\substack{ \pi \in \Pi(F)}} {D}^{\hyb}_\pi,
\]
the versal curve $\rsf^\hyb$ inherits a similar decomposition
\begin{equation} \label{eq:S_hyb_decomposition}
\rsf^\hyb := \rsf^\ast \sqcup \bigsqcup_{\substack{ \varnothing \subsetneq F \subseteq E }} \bigsqcup_{\pi \in \Pi(F)} \rsf_\pi^\hyb.
\end{equation}
In this decomposition, each $\rsf_\pi^\hyb \to D_\pi^\hyb$ is a family of hybrid curves defined over $D_\pi^\hyb$.
That is, for a hybrid point ${\thy} = (t, x) \in D_\pi^{\hyb}$, the fiber $\rsf_{\thy}^\hyb$ is the hybrid curve obtained by replacing each singular point in $\rsf_t$ with an interval of length given by the respective coordinate in $x \in \inn \sigma_{\pi}$. Notice that, topologically, any fiber $\rsf_{\thy}^\hyb$ is homeomorphic to a metrized complex. However, certain analytic quantities associated to the hybrid curve $\rsf_{\thy}^\hyb$ depend on the layered structure (see e.g. the definition of the canonical measure and Green function in Section~\ref{sec:ArakelovGreenfunctionsSurfaces}).

The hybrid versal family $\rsf^\hyb$ is equipped with a hybrid topology, which is again closely related to the decomposition \eqref{eq:S_hyb_decomposition}. We refer to~\cite{AN} for details underlying the construction of $\rsf^\hyb$ and its topology.

\subsection{Adapted coordinates on versal deformation families} \label{sec:AdaptedCoordinatesFamily}
Consider a stable marked  Riemann surface $(S_0, q_1, \dots, q_n)$ and the framework described in the previous sections.

We denote the coordinates in $B = \Delta^N$ by $\underline z = (z_e)_{e\in [N]}$. They are naturally decomposed into $\underline z = \underline z_E\times \underline z_{E^c}$, where $E^c = [N] \setminus E$.

\smallskip

An \emph{adapted set of coordinates} or {\em adapted coordinate system} for the versal family $\pr \colon \rsf \to B$ is the choice of...

\begin{itemize}
\item[(1)] an adapted set of coordinates $(z_e)_{e\in [N]}$ in the sense of Section~\ref{ss:AdaptedCoordinates}, giving a biholomorphism $B \simeq \Delta^N$ with $z_e$, $e\in E$, corresponding to the smoothing of the node $p_e$.

\smallskip

\item[(2)] for each edge $e\in E$ with extremities $u, v \in V$, an open subset $W_e$ around the image $p_e(D_e)$ of the section $p_e \colon D_e \to \rsf$,  and holomorphic coordinates $z_u^e, z_v^e$ defined on (a neighborhood of) the closure $\overline W_e$, such that $\overline W_e = U_e$ for the hypersurface  
\[U_e : =\Bigl\{\bigl(\underline z = (z_e)_{e\in[N]}, z^e_u, z^e_v\bigr) \, \st\, z_u^e z_v^e =z_e\Bigr\}, \qquad U_e \subset B \times \bar \Delta^2\]
and $\pr(\underline z,  z^e_u, z^e_v) = \underline z$ in these coordinates. Here $\bar \Delta$ is the closed unit disk in $\C$.  

\smallskip

Note that $W_e$ is naturally fibered over $B$; the fiber over $t \in B$ is denoted by $W_{e,t}$.

We require in addition that the open sets $W_e$, $e\in E$, be disjoint and to contain no marked points, that is, $q_i(t) \notin W_{e,t}$ for all $i \in [n]$ and $t \in B$.

\smallskip

\item[(3)] for each point $t\in B$, consider the stable Riemann surface $\rsf_{t}$. Let $\rsf_{t} \setminus \bigsqcup_{e\in E} W_{e,t}$ be the closed subset obtained by removing the fibers $W_{e, t}$, $e \in E$ of $W_e$ at point $t$.

Consider the point $t_0 \in B$ with coordinates $\underline z_E(t_0) =0$ and  $\underline z_{E^c}(t_0) = \underline z_{E^c}(t)$. Its fiber $\rsf_{t_0}$ in the versal family is a stable Riemann surface  with nodes $p_e(t_0)$, $e \in E$. The normalization $\widetilde{\rsf}_{t_0}$ is a disjoint union of marked curves $C_{v,t_0}$, $v\in V$. For each vertex $v$ and each incident edge $e=vu$, remove an open disk of radius one around the point $p^e_v$ in $C_{v,t_0}$, in the coordinates given by $z^e_v$ in $W_{e,t_0}$.  Let $Y_{v,t_0}$, $v\in V$, be the corresponding marked Riemann surfaces with boundaries.

Then we require that
\[\rsf_{t} \setminus \bigsqcup W_{e, t} = \bigsqcup_{v\in V} Y_{v,t_0}, \qquad t\in B.\]
\end{itemize}

\begin{thm}\label{thm:adapted-coordinates} For any stable marked curve $S_0$, an adapted set of coordinates exists.
\end{thm} 
The statement of the theorem can be seen as a rectification theorem, arising in the study of holomorphic vector fields, applied to the versal deformation space and the versal family.

\smallskip

Before proving Theorem~\ref{thm:adapted-coordinates}, we discuss its consequences. The key advantage of adapted coordinates is to provide a {\em decomposition of the fibers} $\rsf_t$ which reflects the structure of the dual graph $G=(V,E)$ of the central fiber $S_0$.

\smallskip

More precisely, we can assign to each base point $t\in B$ a collection of marked Riemann surfaces $C_{v, t}$, $v \in V$. Namely, consider the base point $t_0 \in B$ with coordinates $\underline z_E(t_0) =0$, and $\underline z_{E^c} (t_0) = \underline z_{E^c}(t)$. Its fibre $\rsf_{t_0}$ has nodes $p_e(t_0)$ corresponding to the edges $e \in E$. The marked Riemann surface $C_{v, t}$ is the component in the normalization $\widetilde\rsf_{t_0}$ representing $v$, with marking $\mathcal A_v$ given by all points $p^e_v$ for incident half-edges $e$ to $v$ and all markings $q_i(t_0)$ lying on $C_{v, t}$.

Note that  for the edge $e=uv$, $W_{e, t_0}$ consists of two open unit discs $ \{|z_v^e(t)| < 1 \}$, $ \{|z_u^e(t)| < 1 \}$ glued together at their origins (see property (2)). Clearly, in the normalization procedure $W_{e,t_0}$ is split into two open unit disks $ \{|z_v^e(t)| < 1 \}$, $ \{|z_u^e(t)| < 1 \}$, contained in $C_{v,t}$ and $C_{u,t}$, respectively. Defining the marked Riemann surfaces with boundaries
\[Y_{v,t} = C_{v,t}\setminus \bigsqcup_{e\sim v} \big  \{|z_v^e(t)| < 1 \big \},\]
we arrive at the decomposition (see property (3) above)
\begin{equation} \label{eq:AdaptedCoordinates} \rsf_{t}  = \bigsqcup_{v\in V} Y_{v,t} \sqcup \bigsqcup_{e\in E} W_{e, t}.
\end{equation}
Note also that the fiber $W_{e,t}$ over $t \in B$ can be explicitly described as
\begin{equation} \label{eq:Cylinder}
W_{e,t} = \Bigl\{ (z^e_u, z^e_v) \in \Delta \times \Delta \, \st\, z_u^e z_v^e =z_e(t) \Bigr \}.
\end{equation}
The representations \eqref{eq:AdaptedCoordinates} and \eqref{eq:Cylinder} are crucial in the construction of log maps below.

\begin{remark} \label{rem:SetsYvt} We will use the notations $Y_{v, t}$ and $C_{v,t}$ in order to emphasize the dependency on the point $t$. Note however that the marked curve $C_{v, t}$ and its closed subset $Y_{v, t}$ only depend on the parameters $\underline z_{E^c}(t)$.
\end{remark} 
 
\subsection{Proof of Theorem~\ref{thm:adapted-coordinates}}\label{sec:plumbing}
The proof is based on the analytic construction of Deligne--Mumford compactifications in~\cite{HK14}, using the plumbing coordinates~\cite{Mas76,Wol90}. 

\subsubsection{The parameterized analytic cylinder}
Consider the (open and closed) cylinders
\begin{align*} C &= \Big\{(\zeta_1,\zeta_2, z) \in \C^3 \,\st \,\zeta_1\zeta_2=z, \,\, |\zeta_1|, |\zeta_2|, |z|<1 \Bigr\}, \\
\overline C &= \Bigl\{(\zeta_1,\zeta_2, z) \in \C^3 \,\st\,\zeta_1\zeta_2=z, \,\, |\zeta_1|, |\zeta_2| \leq 1, |z|< 1 \Bigr\}.
\end{align*}
$C$ and $\bar C$ are naturally fibered over the open unit disk $\Delta$ by sending
\begin{align*}
\big (\zeta_1,\zeta_2, z \big ) \mapsto z.
\end{align*}
The fibers of $C$ and $\overline C$ at $z \in \Delta$ are denoted by $C_z$ and $\overline C_z$, respectively. 
\subsubsection{The plumbing family}

Let $B = \Delta^N =B_E \times B_{E^c}$ be the base space for the versal deformation space of the stable marked Riemann surface $(S_0, q_1,\dots, q_n)$.
Let $\rsf \to B$ be the analytic versal family of stable marked Riemann surfaces. Recall that $(C_v, \mathcal A_v)$, $v \in V$, denote the irreducible components of $S_0$ with their markings,

\smallskip

The \emph{plumbing family} $\plumb \to B$ is constructed as follows. Consider first the family of stable marked Riemann surfaces $\pr_2 \colon \rsf_{E^c}:=\rsf\rest{B_{E^c}} \to B_{E^c}$ consisting of those Riemann surfaces with the same dual graph as $S_0$, and which are obtained by deforming the marked Riemann surfaces $(C_v, \mathcal A_v)$ in their respective moduli spaces.  

For each edge $e\in E$, we have a section $p_e$ of $\pr_2$ given by the node corresponding to the edge $e$ of the stable Riemann surfaces in the family.

Consider the parametrized closed cylinder 
\[\overline C_e = \Bigl\{(z_u^e,z_v^e, z_e) \in \C^3 \,\st\,z_u^e z_v^e=z_e, \,\, |z_u^e|, |z_v^e| \leq 1, |z_e| < 1 \Bigr\}\]
where $u$ and $v$ are the two extremities of the edge $e$. We choose an open neighborhood $W_e$ of the node $p_e$ in the family $\rsf_{E^c}$, $e\in E$, and an isomorphism $\overline W_e \to B_{E^c} \times \overline C_{e,0}$. Here $\overline C_{e,0}$ is the fiber of $\rho\colon \overline C_e \to \Delta$ at zero. Moreover, we impose that $W_e$ are disjoint from the markings. In the following, for an edge $e$ with the extremities $u$ and $v$, we thus identify each fiber of $\overline W_e$ with $\overline C_{e,0}$ with corresponding coordinates $z^e_u, z_v^e$ with $|z^e_u|, |z^e_v| \leq 1$.

By an abuse of the notation, we denote by $\rsf_{E^c}/B$ the family of marked Riemann surfaces over $B$ obtained by pulling back the family $\rsf_{E^c}/B_{E^c}$ under the projection map $B \to B_{E^c}$. Similarly, we denote by $W_e$ the open neighborhood of $p_e$ in $\rsf_{E^c}/B$ which is isomorphic to  $B \times C_{e,0}$.

For each edge $e\in E$ with extremities $u$ and $v$, let $\overline W'_e$ be the subset of $\overline W_e$ in $\rsf_{E^c}/B$ defined by 
\[W'_e := W_e \setminus \Bigl\{((\underline z_E, \underline z_{E^c}),z_u^e,z_v^e, 0) \in B\times \C^3 \,\st\,z_u^ez_v^e=0, \,\, |z_u^e|, |z_v^e| < |z_e| \Bigr\}\]
where $\underline z = \underline z_E \times \underline z_{E^c} \in B_E \times B_{E^c}$.  That is, we consider a point $\underline z_E \times \underline z_{E^c} \in B_E \times B_{E^c}$ with $|z_e|<1$ for all $e \in E^c$, and from the fiber $\rsf_{E^c, \underline z}$ at $\underline z =\underline z_E\times \underline z_{E^c}$, we remove the two open disks of radius $|z_e|$ around $p_e$ in each of the two branches of $\rsf_{E^c, \underline z }$ around the node $p_e$, $e\in E$.

Moreover, we define $\rsf_{E^c}'/B$ in a similar way as the subset of $\rsf_{E^c}/B$ obtained by removing the two open disks of radius $|z_e|$ around $p_e$ in each of the two branches of $\rsf_{E^c, \underline z}$ around the node $p_e$, for each edge $e$ and each point $\underline z =\underline z_E\times \underline z_{E^c} \in B_E \times B_{E^c}$. That is,
\[\rsf'_{E^c} := \rsf_{E^c} \setminus \bigcup_{e\in E}\Bigl\{(\underline z,z_u^e,z_v^e, 0) \in B\times \C^3 \,\st\,z_u^ez_v^e=0, \,\, |z_u^e|, |z_v^e| < |z_e| \Bigr\}\]

The plumbing family of Riemann surfaces $\plumb/B$ is obtained by plugging in the cylinders $\bar C_{z_e}$, $e \in E$, into $\rsf_{E^c}'/B$. More precisely, for a point $\underline z=\underline z_E \times \underline z_{E^c} \in B_E \times B_{E^c}$, we define 
\[\plumb_{\underline z} := \rquot{\Bigl( \rsf'_{E^c,\underline z} \sqcup \bigsqcup_{e \in E} \bar C_{z_e}\Bigr)}{\sim}\] 
with  $\sim$ corresponding to the identification for $e\in E$ of
\[(z^e_u, z_v^e, 0)\in W'_{e, (\underline z, t)} \qquad \textrm{with}\qquad \begin{cases}\bigl(z^e_u, z_e/z^e_u\bigr) \in \bar C_{z_e} & \textrm{if } z^e_u \neq 0\\
\bigl(z_e/z^e_v , z_v^e \bigr) \in \bar C_{z_e} & \textrm{if } z^e_v \neq 0
\end{cases}, \]
with $u,v$ the extremities of the edge $e$.

The result of~\cite{HK14} we need can be stated as follows. 

\begin{thm}[Hubbard-Koch~\cite{HK14}]\label{thm:hk} There exists a biholomorphic map $\phi \colon B \to B$ which induces an isomorphism between $\plumb/B$ and $\phi^*\rsf/B$, that is, such that we get a commutative diagram
\[
\begin{tikzcd}
\plumb \arrow[r, "\sim"] \arrow[rd] & \arrow[d,  "\phi^*\pi"]\phi^*\rsf \arrow[r] & \arrow[d, "\pi"]\rsf\\
& B \arrow[r,"\phi"]& B
\end{tikzcd}
\]
\end{thm}

\begin{proof}[Proof of Theorem~\ref{thm:adapted-coordinates}]
The adapted coordinates are precisely $z_e$, $e\in [N]$, and the neighborhood $W_e$ of $p_e(D_e)$ is the union of $C_{z_e}$, $z_e\in\Delta$, which is isomorphic to the parametrized cylinder $C$ with coordinates $z^e_u$ and $z^e_v$ verifying $z^e_uz_v^e =z_e$. By the construction, the third property in the definition of the adapted coordinates is also verified. The theorem follows.  
\end{proof}

 \subsection{Mapping Riemann surfaces to hybrid curves} \label{sec:LogSurfacesHybridCurves}
 Notations as above, consider the hybrid versal family $\rsf^\hyb \to B^\hyb$ associated to a marked stable Riemann surface $(S_0, q_1, \dots, q_n)$. In what follows, we construct a map from a smooth Riemann surface to a hybrid curve. Given a smooth Riemann surface and a hybrid curve, which appear as fibers $\rsf_t$ and $\rsf^\hyb_\shy$ in the hybrid versal family $\rsf^\hyb$, we construct a map from $\rsf_t$ to $\rsf^\hyb_\shy$ by using an adapted system of coordinates.  
  
\smallskip 
 
 More precisely, suppose that $t \in B^\ast$ and $\shy \in \partial_\infty B^\hyb$ are two points in the finite part and the hybrid boundary of $B^\hyb$, respectively. By definition, $\shy$ belongs to the hybrid stratum $D_\pi^\hyb$ of an ordered partition $\pi$ of full sedentarity on some $F \subseteq E$, and hence can be written as $\shy = (s, l(e)_{e \in F})$ with $s \in \inn D_F$ and $l\colon E \to (0, + \infty)$ the edge length function normalized in each layer.  Suppose further that $t$ and $\shy$ are compatible in the sense that
  \begin{equation} \label{eq:CompatibleFiberMap}
 	z_e(t) = z_e(s), \,e \in E^c \qquad \text{and} \qquad \arg(z_e(t)) =\arg(z_e(s)), \, e \in E \setminus F.
 \end{equation}
 Under these conditions, we construct a map
 \begin{equation} \label{eq:Maptshy}
 	\psi^{t}_\shy \colon \rsf_t \longrightarrow \rsf^\hyb_\shy
\end{equation}
between the corresponding fibers $\rsf_t = \rsf^\hyb_t$ and $\rsf^\hyb_\shy$ in the family $\rsf^\hyb$. The construction is based on the fiber representations \eqref{eq:AdaptedCoordinates} and \eqref{eq:Cylinder}. Recall that
 \[
 	\rsf_t^\hyb = \rsf_t = \bigsqcup_{v \in V} Y_{v,t} \sqcup \bigsqcup_{e \in E} W_{e,t}.
 \]
 Moreover, the fiber $\rsf_\shy^\hyb$ can be decomposed in a similar way as
 \[
 \rsf_\shy^\hyb = \bigsqcup_{v \in V} Y_{v,s} \sqcup \bigsqcup_{e \in E \setminus F} W_{e,s} \sqcup \bigsqcup_{e=uv \in F} \Delta^e_u \cup I_e \cup \Delta^e_v.
 \]
 where $\Delta^e_u = \{0 \le |z^e_u| < 1\}$, $\Delta^e_v = \{0 \le |z^e_v| < 1\}$ denote open unit disks around the points $p^e_u$ and $p^e_v$  and $I_e$ is an interval of length $l(e)$ connecting $p_u^e$ and $p_v^e$.
 The map $\psi^{t}_\shy \colon \rsf_t \to \rsf^\hyb_\shy$ will send the different types of regions in the above decompositions to each other.

 \subsubsection{The map $\psi^t_\shy$ on ${Y_{v,t}}$} \label{ss:LogMapSmoothCylinders} By the assumption \eqref{eq:CompatibleFiberMap} and  Remark~\ref{rem:SetsYvt}, we have $Y_{v, t} = Y_{v,s} = Y_{v, t_0}$ for the point $t_0 \in B$ with coordinates $z_e(t_0) = 0$, $e \in E$ and $z_e(t_0) = z_e(t) =   z_e(s)$,  $e \in E^c$. Hence on those parts, we can simply define $\psi^t_\shy$ as the identity map,
 \[
\psi^t_\shy \rest{Y_{v, t}} := \id \colon Y_{v, t} \rightarrow Y_{v, s}.
\] 
\subsubsection{The map $\psi^t_\shy$ on ${W_{e,t}}$ for edges $e$ in $E \setminus F$} Note that for all edges $e \in E \setminus F$, we can define a natural map 
\[\psi^t_\shy \colon W_{e, t}\to W_{e,s}\]
 between the cylinders $W_{e,t}$ and $W_{e,s}$. By definition, for any fixed $e \in E \setminus F$,
 \[W_{e, t} =\Bigl\{ (z_u^e, z^e_v) \in \Delta^2 \,|\, z_u^e z_v^e =z_e(t)\Bigr\}, \qquad W_{e, s} =\Bigl\{ (\xi_u^e, \xi^e_v) \in \Delta^2 \,|\, \xi_u^e \xi_v^e =z_e(s)\Bigr\}.\]
 Moreover, $\arg(z_e(s)) =\arg(z_e(t))$ by assumption \eqref{eq:CompatibleFiberMap}. Let  $\alpha_e  = \log|z_e(s)|/\log(|z_e(t)|)$. Then we define $\psi^t_\shy\rest{W_{e, t}}$ as the \emph{real homothecy of factor $\alpha_e$ in logarithm polar coordinates, preserving the angles}, that is a $(z_u^e, z^e_v)$ of $W_{e, t}$ is mapped to the point $(\xi_u^e, \xi^e_v)$ of $W_{e, s}$ with 
 \begin{align*}
 \arg(\xi^e_u) = \arg(z^e_u), &\qquad \log|\xi_u^e| = \alpha_e \log|z^e_u|, \\
 \arg(\xi^e_v) = \arg(z^e_v),&\qquad  \log|\xi_v^e| = \alpha_e \log|z^e_v|.
 \end{align*}

 \subsubsection{Restriction of $\psi^t_\shy$ to ${W_{e,t}}$ for edges in $F$} \label{sec:hybrid_log_map_einf} It remains to define the map $\psi^t \colon \rsf^\hyb_t \to \rsf^\hyb_\shy$ on the cylinders $W_{e,t}$ for $e \in F$.

 In order to define a map $\psi^t_\shy \colon W_{e,t} \to \Delta^e_u \cup I_e \cup \Delta^e_v$, we decompose the cylinder $W_{e,t}$ into three parts $W_{e,t} = A^e_{u,t} \cup B_{e,t} \cup A^e_{v,t}$. Let $\varrho_e, \tau_e \colon B^\ast \to \R_{+}$ be the auxiliary functions
  \begin{equation} \label{eq:Defrhoe}
  \varrho_e(t) := - \frac{1}{\log|z_e(t)|} = \frac{1}{\ell_e(t)}, \qquad \tau_e(t) := - \log\varrho_e(t) = - \log| \log|z_e(t)|| = \log\ell_e(t).
 \end{equation}
 (Here, $\ell_e(t) = -\log\abs{z_e(t)}$, for $e\in F$.) Note that
 \begin{align*}
 \varrho_e(t)\to 0 \qquad \text{and} \qquad \tau_e(t) \rho_e(t) \to 0,
 \end{align*}
 if $z_e(t) \to 0$. 

\smallskip

Consider the two regions $A^e_{u,t}$ and $A_{v,t}^e$ in $\overline W_{e,t}$ defined by
  \begin{align} \label{eq:DefAeu} \begin{split} A^e_{u,t} :&=\Bigl\{ (z_u^e, z^e_v) \,\st \, z_u^e z^e_v =z_e(t),\, |z_v^e|, |z_u^e| \leq 1 \, \textrm{and} \,  |z_u^e| > \varrho_e(t) \Bigr\} \\
  &= \Bigl\{ (z_u^e, z^e_v) \,\st \, z_u^e z^e_v =z_e(t),\, |z_v^e|, |z_u^e| \leq 1 \, \textrm{and} \,  -\log\abs{z_u^e} < \tau_e(t) \Bigr\}.\end{split} \\
A^e_{v,t} :&=\Bigl\{ (z_u^e, z^e_v) \,\st\, z_u^e z^e_v =z_e(t),\, |z_v^e|, |z_u^e| \leq 1 \, \textrm{and} \,  |z_v^e| > \varrho_e(t) \Bigr\}\\\begin{split}
&=\Bigl\{ (z_u^e, z^e_v) \,\st\, z_u^e z^e_v =z_e(t),\, |z_v^e|, |z_u^e| \leq 1 \, \textrm{and} \,  -\log\abs{z_v^e} < \tau_e(t) \Bigr\}\end{split}
\end{align}
 Let moreover
 \begin{align*} B_{e,t} :&= W_{e, t} \setminus \Bigl(A^e_{u,t}\sqcup A^e_{v,t} \Bigr) = \Bigl\{ (z_u^e, z^e_v) \,\st\, z_u^e z^e_v =z_e(t),\, |z_v^e|, |z_u^e| \leq 1 \, \textrm{and} \,  |z_v^e|, |z_u^e| \leq \varrho_e(t) \Bigr\}\\
 &= \Bigl\{ (z_u^e, z^e_v) \,\st\, z_u^e z^e_v =z_e(t),\, |z_v^e|, |z_u^e| \leq 1 \, \textrm{and} \,  -\log|z_v^e|, -\log|z_u^e| \geq \tau_e(t) \Bigr\}
 \end{align*}

 \begin{figure}[!t]
\centering
   \scalebox{.4}{\input{logmap.tikz}}
\caption{The map $\psi^t_\shy$. The edge set $E$ has size two, $F \subset E$ is of size one.}
\label{fig:log-map}
\end{figure}
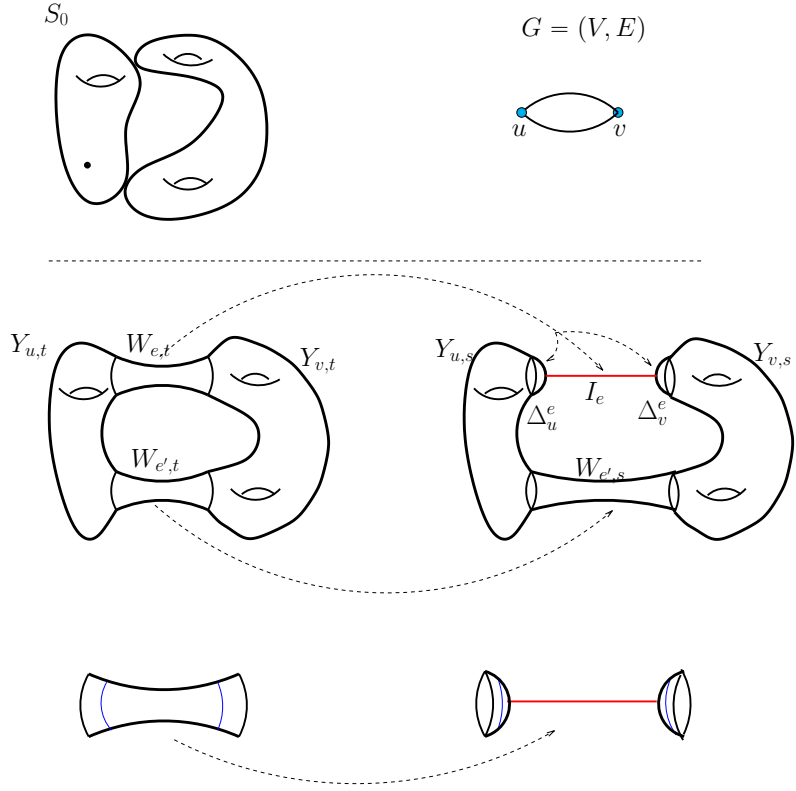

 We will introduce $\psi^t_{\shy}\rest{W_{e, t}}$ by defining its restriction to $A^e_{u,t}, A^e_{v,t}$ and $B_{e,t}$. More precisely, we set
 \begin{align*}
 	\psi^t_\thy(z^e_u, z^e_v) &= z^ e_u \in \Delta^e_u \qquad \text{ for } (z^e_u, z^e_v) \in A^e_{u,t}, \\
	\psi^t_\thy(z^e_u, z^e_v) &= z^ e_v \in \Delta^e_v \qquad \text{ for } (z^e_u, z^e_v) \in A^e_{v,t},
 \end{align*}
 and for $(z^e_u, z^e_v) \in B_{e,t}$, we define
 \[
 \psi^t_\thy(z^e_u, z^e_v) = \frac{l(e)}{\log|z_e(t)| } \bigl(\log|z^e_u(t)|, \log|z^e_v(t)| \bigr ) \in I_e,
 \]
 where we are parametrizing the interval $I_e$ in the form
 \[I_{e} = \Bigl\{(x^e_u, x^e_v) \,\st\, x^e_u+ x^e_v =l(e), \, x_u^e, x_v^e \geq 0\Bigr\}.\]
 
 Note that the map $\psi^t_\shy \colon W_{e,t} \to \Delta^e_u \cup I_e \cup \Delta^e_v$ is \emph{neither continuous nor surjective}. In fact, the images of the respective parts are given by
 \[
 \psi^t_\shy(A^e_{u,t}) = \Bigl \{ z^e_u \in \Delta^e_u \st \, |z^e_u| > \varrho_e(t) \Bigr \}, \qquad  \psi^t_\shy(A^e_{v,t}) = \Bigl \{ z^e_v \in \Delta^e_v \,\st \, |z^e_v| > \varrho_e(t) \Bigr \},
 \]
and
 \[
\psi^t_\shy(B_{e,t}) = \Bigl\{(x^e_u, x^e_v) \in I_e \,\st\, x^e_u, x^e_v \ge l(e)\tau_e(t)\rho_e(t) = l(e)\cdot\frac{\log\ell_e(t)}{\ell_e(t)} \Bigr\}.
 \]
 
  \begin{remark}\label{rem:discontinuity}
 The above map $\psi^{t}_\shy \colon \rsf_t \to \rsf^\hyb_\shy$ clearly depends on the choice of adapted coordinates on the versal family $\rsf \to B$. Note as well that $\psi^{t}_\shy$ is continuous away from the inner cycles of the regions $A^{e}_{u,t}$.
 \end{remark}
 
 \subsubsection{An alternative construction} \label{sss:ContinuousLogMap}
 In the following, we outline how to define another map
 \begin{equation}  \label{eq:ContLog}
 \tilde \psi^t_\shy \colon \rsf_t \to \rsf^\hyb_\shy
  \end{equation} which is both continuous and surjective. The alternative definition will be only used in Section~\ref{ss:HybridLaplacianAsLimit} in order to prove the hybrid Laplacian is a weak limit of the Laplacian on close-by Riemann surfaces (see however Remark~\ref{rem:validity-weak-convergence}).  So, unless otherwise explicitly stated, we will use the construction given in the previous section which is the adapted one in the situation of interest in this part, that is, in providing approximate solutions to the Poisson equation on Riemann surfaces close to the boundary of the moduli spaces, \cf\, Remark~\ref{rem:two-log-maps}.

  \smallskip
  
  We will only modify the map on the sets $W_{e,t}$ for $e \in F$. First of all, there is a natural symmetric homeomorphism $\varphi \colon \psi^t_\shy(B_{e,t}) \to I_e$. This map is defined by 
 \begin{align*}
 \varphi \colon \psi^t_\shy(B_{e,t}) &\longrightarrow I_e \\
 (x^e_u, x^e_v) &\longmapsto \frac{l(e)}2\cdot \Bigl(1+\frac{x^e_u-x^e_v}{l(e) - 2 l(e)\tau_e(t)\rho_t(t)}\,, 1+\frac{x^e_v-x^e_u}{l(e) - 2 l(e)\tau_e(t)\rho_t(t)}\Bigr)
 \end{align*} 
  
\smallskip
  
On the subregion $B_{e,t} \subseteq W_{e,t}$, the map $\tilde \psi^t_\shy$ is the composition $\tilde \psi^t_\shy := \varphi \circ \psi^t_\shy$. In terms of the original coordinates, we send $(z^e_u, z^e_v) \in B_{e,t}$ to the point $\tilde \psi^t_\shy (z^e_u, z^e_v) = (\tilde x^e_u, \tilde x^e_v)$ in $I_e$  given by
  \[\tilde x^e_u =  \frac{l(e)}{\ell_e(t)-2 \tau_e(t)} (- \log|z^e_u| -\tau_e) \qquad \textrm{and} \qquad \tilde x^e_v =  \frac{l(e)}{\ell_e(t)-2\tau_e(t)} (- \log|z^e_v| -\tau_e).\]

Moreover, using polar coordinates, the image $\psi^t_\shy(A^e_{u,t})$ and the closed disc $\bar \Delta^e_u \subseteq \curve$ can be written as products
\[
\psi^t_\shy(A^e_{u,t}) \cong [\varrho_e, 1 ] \times \realtor  \qquad \bar \Delta^e_u = \{0\} \cup (0,1] \times \realtor.
\]
 Using these coordinates, we define the map
 \[\tilde \psi^t_\shy \rest{A^e_u} \colon A^e_{u,t} \to \bar \Delta^e_u\]
 by sending
  \begin{align*}
 (z^e_u, z^e_v)&\mapsto \Big ( \frac{|z^e_u|- \varrho_e(t)}{1 - \varrho_e(t)}, \arg(z^e_u) \Big ).
 \end{align*}

 We use the analog definition for the restriction 
  \[ \tilde \psi^t_\shy \rest{A^e_{v,t}} \colon A^e_{v,t} \to \bar \Delta^e_v.\]

\subsection{The hybrid log maps on $\loghyb$ and $\loghyb_\pi$ on versal families}\label{sec:hybrid_log_map_rs}

Notations as in the previous sections, consider the hybrid versal family $\rsf^\hyb \to B^\hyb$ associated with a marked stable Riemann surface $(S_0, q_1, \dots, q_n)$.

Recall that in Section~\eqref{sec:hybrid_log_I}, we defined two types of log maps on the hybrid base $B^\hyb$: for each boundary stratum $D_\pi^\hyb$, the stratumwise log map $\loghyb_\pi \colon B^\ast \to D_\pi^\hyb$; and the (global) hybrid log map $\loghyb \colon \mathscr U \to \partial_\infty B^\hyb$, which is a retraction from an open neighborhood $ \mathscr U$ of $\partial_\infty B^\hyb$ to the hybrid boundary $\partial_\infty B^\hyb$.

\smallskip

Using the construction from the previous section, we extend the log maps on the base to the hybrid versal family. More precisely, we define two maps
\begin{align*}
\loghyb_{\pi} \colon \rsf^\ast \to \rsf_\pi^\hyb &&\loghyb \colon \rsf^\hyb\rest{\mathscr U}  \to  \rsf^\hyb\rest{\partial_\infty B^\hyb}
\end{align*}
where $\rsf^\ast=\rsf^\hyb\rest{B^\ast}$ and $\rsf^\hyb_\pi =\rsf^\hyb\rest{D^\hyb_\pi}$, such that the following two diagrams commute:
 \begin{equation} \label{eq:CDHybridLogs}
\begin{tikzcd}
 \rsf^\ast \arrow[d]\arrow[r, "\loghybdiagpi{\pi}"] & \arrow[d] \rsf_\pi^\hyb\\
B^\ast\arrow[r, "\loghybdiagpi{\pi}"] &D_\pi^\hyb
\end{tikzcd}
\qquad  \qquad \qquad
\begin{tikzcd}
 \rsf^\hyb\rest{\mathscr U} \arrow[d]\arrow[r, "\loghybdiag"] & \arrow[d] \rsf^\hyb\rest{\partial_\infty B^\hyb}\\
\mathscr U\arrow[r, "\loghybdiag"] &\partial_\infty B^\hyb 
\end{tikzcd}
\end{equation}
To ensure commutativity of the diagrams, we should define for each point  $t \in B^\ast$ a map
 \[\loghyb_{\pi,t} \colon \rsf_{t} \to \rsf^\hyb_{\shy},\qquad \shy = \loghyb_\pi(t) \in D^\hyb_\pi\]
 and for each point $t \in \mathscr U$ a map
 \[\loghyb_{t} \colon \rsf^\hyb_{t} \to \rsf^\hyb_{\shy'}, \qquad \shy' = \loghyb(t) \in \partial_{\infty} B^\hyb.\]
 The log maps are then given by setting $\loghyb_{\pi}(q) = \loghyb_{\pi,t}(q)$ and $\loghyb (q) = \loghyb_{t}(q)$ for $q \in \rsf^\hyb_t$.
  
  \smallskip
 
Note that for given $t \in B^\ast$, its images $\shy = \loghyb_\pi(t)$ and $\shy' =  \loghyb(t)$ satisfy the conditions~\eqref{eq:Maptshy} (see Section~\ref{sec:log-hybrid-stratumwise} and Section~\ref{ss:GlobalHybridLogMap}).  Hence we can set
 \[
 \loghyb_{\pi,t} =\Psi^t_\shy \colon \rsf_{t} \to \rsf^\hyb_{\shy},
 \]
 where $\Psi^t_\shy \colon \rsf_{t} \to \rsf^\hyb_{\shy}$ is the map defined in the previous section (see \eqref{eq:Maptshy}).
 
 Moreover, recall that $\loghyb(\thy) = \thy$ for $\thy \in \partial_\infty B^\hyb$. Hence we can define
 \begin{align*}
  \loghyb_{\pi,t} &=\Psi^t_{\shy'} \colon \rsf_{t} \to \rsf^\hyb_{\shy'}, &&t \in B^\ast\cap \mathscr{U}, \\[2mm]
  \loghyb_{\pi,\thy} &=\operatorname{Id} \colon \rsf_{\thy}^\hyb \to \rsf^\hyb_{\thy}, &&\thy \in \partial_\infty B^\hyb.
 \end{align*}
Note that $\loghyb \colon \rsf^\hyb\rest{\mathscr U}  \to  \rsf^\hyb\rest{\partial_\infty B^\hyb}$ retracts  $\rsf^\hyb\rest{\mathscr U}$ to $ \rsf^\hyb\rest{\partial_\infty B^\hyb}$.

 \begin{thm} The hybrid log maps $\loghyb_{\pi} \colon \rsf^\ast \to \rsf_\pi^\hyb$ and $\loghyb \colon \rsf^\hyb\rest{\mathscr U}  \to  \rsf^\hyb\rest{\partial_\infty B^\hyb}$ are well-defined and continuous away from a subspace of codimesion one. Moreover, the diagrams in \eqref{eq:CDHybridLogs} commute.
 \end{thm}
\begin{proof}
This follows from the explicit form of the log maps described above. Note that the maps are continuous away from fibrations by circles indexed by edges and their extremities, namely, along the inner cycles of the regions $A^e_{u,t}$, which are the locus of discontinuities of auxiliary functions $\psi_{e,t}$, see Remark~\ref{rem:discontinuity}.
\end{proof}


\section{Hybrid functions and hybrid Laplace operator}
\label{sec:hybrid_laplacian}

In this section, we introduce a framework for doing \emph{function theory on hybrid curves} (Sections~\ref{ss:FunctionTheoryHybridCurve} and~\ref{ss:TameHybridFunctions}), define a notion of \emph{Laplacian on hybrid curves} (Section~\ref{ss:HybridLaplacian}), and explain how the latter can be viewed as a limit of the Laplacian on degenerating Riemann surfaces (Section~\ref{ss:HybridLaplacianAsLimit}).

\subsection{Function theory on hybrid curves}   \label{ss:FunctionTheoryHybridCurve}
Let $\curve=\curve^\hyb$ be a hybrid curve of rank $r$ with underlying graph $G=(V, E)$, the ordered partition $\pi = (\pi_1, \dots, \pi_r)$ of $E$, and the edge length function $l\colon E \to (0, + \infty)$. Let $\curve^\trop$ be the underlying tropical curve, which is of full sedentarity, and $\Gamma$ the normalized layered metric graph in the conformal equivalence class of $\curve^\trop$. Recall that this means the sum of edge lengths in each layer is equal to one.  For each vertex $v \in V$, let $C_v$ be the corresponding Riemann surface. 
We call $\pi_\smallc := \bigsqcup_{v\in V} C_v$ the \emph{finitary} or \emph{complex part} of $\curve^\hyb$. 

\smallskip

\subsubsection{Basic definitions}
A  (\emph{complex-})\emph{valued hybrid function} on $\curve$ is an $(r+1)$-tuple
\[
\lf = (f_1, \dots, f_r, f_\smallc)
\]
consisting  of functions $f_j\colon \Gamma^j \to \C$ on the graded minors $\Gamma^j$, $j\in[r]$, and a function $f_\smallc \colon\pi_\smallc \to \C$ on the complex part $\pi_\smallc$.  The \emph{tropical part} of $\lf$ is the function $\lf^\trop := (f_1, \dots, f_r)$ on the underlying tropical curve $\curve^\trop$. We refer to  $f_\smallc$ as the \emph{complex part} of $\lf$. 

 A hybrid function $\lf$ on $\curve^\hyb$ is called \emph{continuous, smooth}, etc., if the corresponding functions on $\Gamma^j$, $j\in[r]$, and on $\pi_\smallcc$ are continuous, smooth, etc., respectively. 

\medskip

In Section~\ref{ss:HybridLaplacian}, we will define a hybrid Laplacian on sufficiently regular hybrid functions. In this context, we consider hybrid functions $\lf = (f_1, \dots, f_r, f_\smallcc)$ such that each component $f_j$, $j \in [r]$, belongs to the Zhang space $D(\Delta_{\Gamma^{\grind{j}}})$ on $\Gamma^j$ (see Section~\ref{sec:potential_theory_metric_graphs}). That is, we assume that
\begin{itemize}
\item[(i)] $f_j$ is continuous on $\Gamma^j$,
\item[(ii)] $f_j$ is piecewise $\mathcal C^2$ on $\Gamma^j$, and
\item[(iii)] $f_j''$ belongs to $L^1(\Gamma^j)$.
\end{itemize}
Furthermore, we suppose that the complex part $f_\smallc$ is either $\mathcal C^2$, or, more generally, has finitely many logarithmic poles on $\pi_\smallc$. Hereby we mean that there are at most finitely many points $p_1, \dots, p_n$ on $\pi_\smallcc = \bigsqcup_{v\in V} C_v$ such that in a local holomorphic coordinate $z_j$ around each $p_j$,
\begin{equation} \label{eq:RegularizedFunction}
f_\smallcc = r_j \log|z_j| + f^\reg_{\smallc}
\end{equation}
for a real number $r_j$ and a (locally defined) $\mathcal{C}^2$ function $f^\reg_{\smallc}$, and that $f_\smallcc$ is $\mathcal{C}^2$ everywhere else on $\pi_\smallcc$. The space of all functions $f_\smallcc$ on $\pi_\smallcc$ with these properties is denoted by $D_{\log}(\Delta_\smallcc)$. Allowing logarithmic poles in the complex part $f_\smallcc$ turns out to be crucial in context with the hybrid Poisson equation (see Section~\ref{sec:hybrid_green}).

\smallskip 

The \emph{hybrid Arakelov--Zhang space} $D(\Deltahybind{\curve})$ is defined as the set  all hybrid functions $\lf = (f_1, \dots, f_r, f_\smallc)$ such that $f_j \in D(\Gamma^j)$ for $j\in [r]$, and the complex part $f_\smallc$ is a $\mathcal C^2$ function. We also define \emph{the extended hybrid Arakelov--Zhang space} $D_{\log}(\Deltahybind{\curve})$ as the set of all hybrid functions $\lf $ with $f_j \in D(\Gamma^j)$ for $j\in [r]$, and $f_\smallc \in D_{\log}(\Delta_\smallcc)$.

\subsubsection{Regularization of complex functions with logarithmic poles} \label{ss:RegularizationComplexPart}
As outlined above, we often consider hybrid functions $\lf$ whose complex parts $f_\smallcc$ have logarithmic poles. In particular, $f_\smallcc$ might have poles in the attachement points $p^e_v$, and hence the values there are not necessarily defined. However, for us, it will be crucial to assign values also in these points.

\smallskip 

We resolve this issue by the following \emph{regularization} procedure. Consider a function $f_\smallcc$ in $D_{\log}(\Delta_\smallcc)$. Let $p^e_v$, $e \sim v$ be one of the marked points on a component $C_v$ of the hybrid curve $\curve$. Choosing a local holomorphic coordinate $z^e_v$ around $p^e_v$ on $C_v$, the function $f_\smallcc$ is locally of the form \eqref{eq:RegularizedFunction}. In particular, \emph{the regularized value} $f^\reg_{\smallcc}(p^e_v)$ is well-defined. We extend $f_\smallcc$ to the point $p^e_v$ by setting
\begin{equation} \label{eq:RegularizedValue}
f_\smallc(p^e_v) :=  f^\reg_{\smallc}(p^e_v) := \lim_{p \to p^e_v} f_\smallcc(p) - r^e_v\, \log|z^e_v(p)|.
\end{equation}
In case that $p^e_v$ is not a pole, we simply recover the original value of the function.

\smallskip

We stress that \emph{the regularized values depend on the choice of the local coordinates}. In the case where the hybrid curve $\curve$ arises as the limit of a family of Riemann surfaces, we will fix them globally by choosing an adapted system of coordinates. We address this issue later.

\subsubsection{Pullback of hybrid functions to metrized complexes} \label{ss:PullbackHCtoMC}

Similar to defining pullbacks of functions on tropical curves to metric graphs, we now introduce \emph{pullbacks of hybrid functions on hybrid curves to metrized compelxes}.

\smallskip

Let $\mc$ be a metrized complex with the same underlying stable Riemann surface as the hybrid curve $\curve$. Denote by $\ell \colon E \to (0, + \infty)$ the edge length function in $\mc$ and let $\mgr$ be the corresponding metric graph. Let $\lf = (\lf^\trop, f_\smallc)$ be a hybrid function on $\curve$ with the tropical part $\lf^\trop = (f_1, \dots, f_r)$ a function on $\curve^\trop$. 

\smallskip

For each $j\in[r]$, we define the pullback $f_j^\ast \colon \mc \to \C$ by taking first the pullback of $f_j$ to $\mgr$, as in Section~\ref{sec:tropical_function_theory}, and then by composing it with the projection map $\mc \to \mgr$ (which contracts each Riemann surface $C_v$ to the vertex $v$). 

\smallskip

The pullback of $f_\smallc$ to $\mc$ is the function 
\[f_{\smallc}^{\ast}:\mc \to \C\]
obtained by extending $f_\smallcc$ to the edges through \emph{linear interpolation}. More precisely, on each edge $e=uv\in E$, the pullback $f_{\smallc }^{\ast}$ is given by
\begin{equation}
f_{\smallcc}^\ast \rest e (x) = f_\smallc (p^e_u) + \frac{f_\smallcc(p^e_v) - f_\smallcc(p^e_u)}{\ell(e)}x.
\end{equation} 
Here $x$ is the parametrization of the (oriented) edge $e=uv$ in  $\mgr$ (after identification with the interval $[0,\ell(e)]$) and $p^e_u$ and $p^e_v$ are the corresponding marked points on $C_u$ and $C_v$, respectively.

\smallskip
Note that the above can only be applied if the complex part $f_\smallcc$ has values defined in all of the marked points $p^e_v$ . However, we can use the regularization from Section~\ref{ss:RegularizationComplexPart} to extend the definition to the case where $f_\smallcc$ has logarithmic poles.
\smallskip

More precisely, suppose that the complex part $f_\smallcc$ belongs to $D_{\log}(\Delta_\smallcc)$. Fix further local coordinates $z^e_v$ around the marked points $p^e_v$, $e \sim v$ on the components $C_v$, $v \in V$ of $\curve$. We define the pullback of $f_\smallc$ to $\mc$ as  the function $f_{\smallc}^{\ast}:\mc \to \C$ obtained by \emph{linear interpolation of the regularization} of $f_\smallc$ on the edges. That is, on each edge $e=uv\in E$,
\begin{equation}
f_{\smallcc}^\ast \rest e (x) = f_\smallcc^\reg (p^e_u) + \frac{f_\smallcc^\reg(p^e_v) - f_\smallcc^\reg(p^e_u)}{\ell(e)}x,
\end{equation}
where $f_\smallcc^\reg (p^e_u)$ is the regularized value \eqref{eq:RegularizedValue}. Clearly, the pullback $f_\smallcc^\ast$ depends on the choice of the local coordinates $z^e_v$.

\smallskip

Recall also that to a hybrid curve $\curve$, we can associate naturally a metrized complex $\mccan$ with normalized edge length function on each layer. The latter provides a canonical representative in the respective conformal equivalence class. 
The above in particular leads to the definition of a pullback of $\lf$ to $\mccan$. Note that the metric graph underlying $\mccan$ is $\Gamma$, the canonical representative of $\curve^\trop$. In the sequel, by an abuse of the notation we sometime drop the $\ast$ and simply denote by $f_j\colon \mccan \to \C$ the corresponding pullback on the metrized complex $\mccan$, for $j\in [r] \cup\{\smallcc\}$.

\subsubsection{Harmonically arranged hybrid functions} \label{ss:HarmonicallyArrangedHybridFunctions} Notations as above, a hybrid function $\lf = (\lf^\trop, f_\smallc)$ on $\curve$ is called \emph{harmonically arranged} if the following conditions hold:

\begin{itemize}
\item the tropical part $\lf^\trop$ is harmonically arranged on $\curve^\trop$, \ie, all functions $f_j$, $j\in [r]$ are lower harmonic in $\curve^\trop$, and
\item the complex part $f_{\smallc}$ is harmonic on all tropical layers. Hereby we mean that for all $i\in[r]$, and for any vertex $w$ of the graded minor $\Gamma^i$, we have 
\[
 \sum_{u \in \kappa_i^{-1}(w)} \sum_{\substack{{e = uv} \\ e \in \pi_i}} \frac{f_\smallc(p^e_v) - f_\smallc(p^e_u)}{l(e)} = 0,
\]
where $\kappa_i \colon V \to V(\Gamma^i)$ is the contraction map.
\end{itemize}
Equivalently, the last condition says that the differential of the pullback $f_{\smallc}^{\ast}$ of $f_\smallc$ is harmonic on all tropical layers, meaning that
\[
 d f_{\smallc}^{\ast}\rest{\pi_i} \in \Omega^1(\grm{\pi}{i}(G))
\]
for all $i=1, \dots, r$, where
\[
d f_{\smallc}^{\ast} (e)  = \frac{f_\smallc(p^e_v) - f_\smallc(p^e_u)}{l(e)}, \qquad e=uv \in \pi_i.
\]
Note that this definition can only be applied if the complex part $f_\smallcc$ has values defined in all of the marked points $p^e_v$ . However, in case that $f_\smallcc$ has logarithmic poles, we can use the regularization from Section~\ref{ss:RegularizationComplexPart} to define the property.

\smallskip

More formally, let $\lf = (\lf^\trop, f_\smallc)$ be a hybrid function whose complex part $f_\smallcc$ belongs to $D_{\log}(\Delta_\smallcc)$. Fix a local coordinate $z^e_v$ around each marked point $p^e_v$, $e \sim v$ on each component $C_v$, $v \in V$ of the hybrid curve $\curve$.

We say that $\lf$ is \emph{harmonically arranged}, with respect to the coordinate system $(z^e_v)$, if the tropical part $\lf^\trop$ is harmonically arranged on $\curve^\trop$ and the regularization of the complex part $f_\smallcc$ is harmonic on all tropical layers (in the above sense).

More precisely, we require that the one-form $\alpha_\smallc \in C^1(G, \R)$ defined by
\[
\alpha_\smallc (e) = \frac{f^\reg_\smallcc (p^e_v) -  f^\reg_\smallcc (p^e_w)}{l(e)}, \,\, e=vw \in \E,
\]
belongs to the space of harmonic one-forms $\Omega^1(\grm{\pi}{j}G)$ for any $j
\in [r]$. Equivalently,
\begin{equation} \label{eq:RegEq}
\sum_{e= vw \in \pi_j} \alpha_{\smallc}(e) =0
\end{equation}
for all vertices $v\in V(\grm{\pi}{j}G)$. 

\smallskip
Note that for functions with singularities, the above property depends on the choice of the local coordinates $z^e_v$. However, it is preserved under the following basic changes.

\begin{remark} Consider a map $a \colon \E \to \C^*$, $e\mapsto a_e$, which verifies $a_e = a_{\bar e}^{-1}$ for all oriented edge $e\in \E$. Changing local parameters $z^e_u$  to $a_e z^e_u$ for all $e\sim v$  does not change the one form $\alpha_{\smallc}$. Indeed, the functions $f_j$ are all affine linear on each edge of $G$, and so we have $c_{u}^e =-c_v^e$ for all edges $e=uv$ in the graph. It follows that the one-form $\alpha_\smallc$ remains unchanged under this choice of local coordinates. 
\end{remark}

\subsubsection{Harmonic rearrangement of hybrid functions} Similar to the case of tropical curves, hybrid functions on hybrid curves can be transformed in an essentially unique manner into harmonically arranged ones.
 \smallskip

 Let $\lf = (\lf^\trop, f_\smallc)$ be a hybrid function on $\curve$ with tropical part $\lf^\trop=(f_1,\dots, f_r)$ a function on the tropical curve $\curve^\trop$. In the following discussion, we allow that the complex part $f_\smallcc$ has finitely many logarithmic poles.

 We then construct a new hybrid function $\tilde{\lf} = (\tilde \lf^\trop, \tilde f_\smallc)$ by first taking a rearrangement $\tilde \lf^\trop$ of $\lf^\trop$ on $\curve^\trop$, and then choosing a constant $\consts{\smallc}{v}$ for each $v\in V$, and finally setting
\[\tilde f_{\smallc{_{,v}}} = f_{\smallc{_{,v}}} + \consts{\smallc}{v}
\]
on the corresponding component $C_v$.

We call a function $\tilde{\lf} = (\tilde{f}_1, \dots, \tilde f_r, \tilde f_\smallc)$ obtained in this way a \emph{rearrangement} of $\lf$, and a \emph{harmonic rearrangement} if $\tilde{\lf}$ is harmonically arranged.

\smallskip

We have the following analogue of the rearrangement property in the hybrid setting.

\begin{thm}[Harmonic rearrangement of hybrid functions] \label{thm:HarmonicExtensionHybrid}
Every hybrid function $\lf = \bigl(f_1, \dots, f_r, f_\smallc\bigr)$ on a hybrid curve $\curve$ has a harmonic rearrangement. The latter is unique up to a global constant on $\curve$, that is, an additive constant on each layer of the tropical curve, and  an additive constant on the disjoint union $\pi_\smallcc = \bigsqcup_{v \in V} C_v$.

Namely, for two such rearrangements $\tilde{\lf}$ and $\tilde{\lf}'$, the difference of the $j$th components $\tilde{f}_j - \tilde{f}_j'$ is constant on $\Gamma^j$ for all $j\in [r]$, and the difference of the complex parts is constant on $\pi_\smallc$.
\end{thm}
\begin{proof} The proof is similar to that of Theorem~\ref{thm:HarmonicExtension}. 
\end{proof}

\subsubsection{Pullback of hybrid functions to Riemann surfaces} \label{sec:pullback_rs} Consider the setting from Section~\ref{sec:hybrid_log_II}. Namely, let $\rsf \to B$ be the analytic versal deformation space of a stable marked Riemann surface $(S_0, q_1, \dots, q_n)$, and let  $\rsf^\hyb \to B^\hyb$ be the corresponding family of hybrid curves. We assume that an adapted coordinate system is fixed on $\rsf \to B$.  In the following, we explain how to pullback hybrid functions living on the hybrid curves in the boundary $\partial_\infty B^\hyb$ to the smooth Riemann surfaces in  the family $\rsf$. These constructions are based on the log maps introduced in Section~\ref{sec:hybrid_log_II}. 

\medskip

Let $\shy = (s, x)$ be a point in the hybrid boundary $\partial_\infty B^\hyb$, belonging to the hybrid stratum $D_\pi^\hyb$ of an ordered partition $\pi=(\pi_1, \dots, \pi_r) \in \Pi(F)$ of some $F \subseteq E$. Consider a hybrid function $\lf = (f_1, \dots, f_r, f_\smallc)$ on the hybrid curve $\rsf^\hyb_\thy$. We stress that in what follows, the complex part $f_\smallcc$ is allowed to have logarithmic poles, that is, we allow that  $f_\smallcc \in D_{\log}(\Delta_\smallcc)$.

Viewing the fiber $\rsf^\hyb_\shy$ as the canonical metrized complex $\mccan_\shy$ with normalized edge lengths on each layer $\pi_j$, $j\in[r]$, the construction in Section~\ref{ss:PullbackHCtoMC} gives a pulled back function $f_k^\ast \colon \rsf^\hyb_\shy \to \C$ for each piece $f_k$, $k\in [r] \sqcup\{\smallc\}$ of $\lf$.  Recall that in Section~\ref{sec:LogSurfacesHybridCurves}, for every $t \in B^\ast$ such that 
\begin{equation} \label{eq:CompatibleFiberMap2}
z_e(t) = z_e(s), \,e \in E^c, \qquad \text{and} \qquad \arg(z_e(t)) =\arg(z_e(s)), \, e \in E \setminus F,
\end{equation}
we constructed a log map $\psi^t_\shy$ from the smooth Riemann surface $\rsf_t$ to the fiber $\rsf^\hyb_\shy$. 

\smallskip

For any $t \in B^\ast$ satisfying \eqref{eq:CompatibleFiberMap2}, we define the pullback
\[
\lf^*_{t}\colon \rsf_{t} \to \C
\]
of the hybrid function $\lf$ to the smooth Riemann surface $\rsf_t$ by
\begin{equation} \label{eq:SimplestPullback}
\lf^*_{t} := \sum_{j=1}^r L_j(t) \, \big( f_j^*\circ \psi^t_\shy  \big) + f_\smallc^*\circ \psi^t_\shy,
\end{equation}
where the scaling factors $L_j(t)$, $j=1, \dots, r$ are defined by $L_j(t) = - \sum_{e\in \pi_j} \log|z_e(t)|$.

Altogether, the hybrid function $\lf$ becomes propagated from the hybrid curve $\rsf^\hyb_\shy$ to all smooth Riemann surfaces $\rsf_t$ in the family $\rsf^\ast  \to B^\ast_\shy$, where 
\[B^\ast_\shy= \Bigl\{ s \in B^\ast \,\st \, s \text{ satisfies \eqref{eq:CompatibleFiberMap2}} \Bigr\}.\]
 This procedure can be viewed as the analog of the considerations in Section~\ref{sss:PullbackExplanation} for tropical curves and metric graphs.

\medskip

Using the hybrid log maps from Section~\ref{sec:hybrid_log_map_rs}, the above construction can be \emph{globalized} as follows. Given a family $\mathscr{F}$ of hybrid functions, living on the hybrid curves over the hybrid boundary $\partial_\infty B^\hyb$, we obtain a pulled-back family $\mathscr{F}^\ast$ of functions living on the smooth Riemann surfaces over $B^\ast$. Recall that in Section~\ref{sec:hybrid_log_map_rs}, we constructed two log maps, the stratumwise log map $\loghyb_\pi$ (for each boundary stratum $D_\pi^\hyb$) and the (global) hybrid log map $\loghyb$, resulting in the commutative diagrams (see \eqref{eq:CDHybridLogs} for details)
 \begin{equation} \label{eq:CDHybridLogs}
\begin{tikzcd}
 \rsf^\ast \arrow[d]\arrow[r, "\loghybdiagpi{\pi}"] & \arrow[d] \rsf_\pi^\hyb\\
B^\ast\arrow[r, "\loghybdiagpi{\pi}"] &D_\pi^\hyb
\end{tikzcd}
\qquad  \qquad \qquad
\begin{tikzcd}
 \rsf^\hyb\rest{\mathscr U} \arrow[d]\arrow[r, "\loghybdiag"] & \arrow[d] \rsf^\hyb\rest{\partial_\infty B^\hyb}\\
\mathscr U\arrow[r, "\loghybdiag"] &\partial_\infty B^\hyb.
\end{tikzcd}
\end{equation}

Fix a boundary stratum $D_\pi^\hyb$ in $B^\hyb$, associated to an ordered partition $\pi = (\pi, \dots, \pi_r)$ on some subset $F \subseteq E$. Consider a family $\mathscr{F}_\pi = (\lf_\thy)_{\thy \in D_\pi^\hyb}$ of hybrid functions over $D_\pi^\hyb$. Hereby we mean that for each $\thy \in D_\pi^\hyb$, we have given a hybrid function
\[
	\lf_\thy = (f_{\thy, 1}, \dots, f_{\thy, r}, f_{\thy , \smallcc})
\]
on the hybrid curve $\rsf_\thy^\hyb$ over $\thy$. From the family of hybrid functions $\mathscr{F}_\pi$, one obtains a family $\mathscr{F}_\pi^\ast = (\lf_t^\ast)_{t \in B^\ast}$ of functions on the smooth Riemann surfaces $\rsf_t$, $t \in B^\ast$, as follows. For each $t \in B^\ast$, we define the function $\lf_t^\ast \colon \rsf_t \to \C$ by setting
\begin{equation} \label{eq:StratumwisePullback}
\lf^*_{t} := \sum_{j=1}^r L_j(t) \Big( f_{\thy, j}^*\circ \loghyb_{\pi}  \Big) + f_{\thy,  \smallcc}^*\circ \loghyb_{\pi}, \qquad \thy = \loghyb_{\pi}(t), 
\end{equation}
and $L_j(t) = - \sum_{e \in \pi_j} \log|z_e(t)|$. Since $\loghyb_{\pi} = \psi^t_\thy$ on the Riemann surface $\rsf_t$  for $ \thy = \loghyb_{\pi}(t)$, this is actually just the pullback of the hybrid function $f_\thy$  to $\rsf_t$ in the sense of \eqref{eq:SimplestPullback}.

\smallskip

Similarly as above, we can define a pullback by considering the global hybrid log map $\loghyb$. Given a family $\mathscr{F}_\infty = (f_\thy)_{\thy \in \partial_\infty B^\hyb}$ of hybrid functions over the hybrid boundary $\partial_\infty B^\hyb$, we obtain a family $\mathscr{F}_\pi^\ast = (f_t^\ast)_{t \in B^\ast}$ of functions on the smooth Riemann surfaces $\rsf_t$, $t \in \mathscr U$. For each $t \in \mathscr U$, the function $\lf_t^\ast \colon \rsf_t \to \C$ is given by
\begin{equation} \label{eq:GlobalPullback}
\lf^*_{t} := \sum_{j=1}^r L_j(t) \Big( f_{\thy, j}^*\circ \loghyb  \Big) + f_{\thy, \smallcc}^*\circ \loghyb, \qquad \thy = \loghyb (t),
\end{equation}
where $r$ is the rank of the hybrid curve $\rsf^\hyb_\thy$ and $L_j(t) = - \sum_{e \in \pi_j} \log|z_e(t)|$, $j=1,\dots,r$ for each $\pi_j$ in the ordered partition $\pi=(\pi, \dots, \pi_r)$ associated with $\thy$. Note that the number of terms in the additive expansion depends on the stratum $D_\pi^\hyb$ of the point $\thy$. 

\smallskip

We refer to the families $\mathscr{F}_\pi^\ast$ and $\mathscr{F}_\infty^\ast$ as the \emph{pullbacks} of the families $\mathscr{F}_\pi$ and $\mathscr{F}_\infty$, by the maps $\loghyb_{\pi}$ and $\loghyb$, respectively. The two notions of pullbacks are obviously distinct.

\smallskip

Note that, in general, the pullbacks $f_t^\ast$ of a continuous hybrid function $\lf$ can be discontinuous along the inner cycles of the regions $A^e_v$, since the map $\psi^t_\shy \to  \rsf^\hyb_\shy$ is discontinuous there. Note also that, if the complex part $f_\smallcc$ has a logarithmic pole at a marked point $p^e_v$, then the pullback is bounded on the region $W_{e,t}$ (since $\psi^t_\shy \to  \rsf^\hyb_\shy$ is not surjective).

 \begin{remark} \label{rem:two-log-maps}
In all the above constructions, we can replace the discontinuous maps $\psi^t_\shy \to  \rsf^\hyb_\shy$ by the continuous map $\tilde \psi^t_\shy \to  \rsf^\hyb_\shy$ (and also make the respective change in the definition of the log maps $\loghyb_\pi$ and $\loghyb$). As a result, the pullbacks $f_t^\ast$ of continuous hybrid functions $\lf$ are continuous on the smooth Riemann surfaces $\rsf_t$.

However, it turns out that this is not suitable to obtain approximations of the solutions to the Poisson equation $\Delta f_t = \mu_t$  on degenerating Riemann surfaces $\rsf_t$ . Naively, one can justify this by observing that the pull-back of a function $\lf$ in $D_{\log}(\Deltahybind{\curve})$ (these appear in context with the hybrid Poisson equation)  to $\rsf_t$ can have poles along the inner cycles of the regions $A^e_v$. However, the solution of the Poisson equation $\Delta f_t = \mu_t$ on the Riemann surface $\rsf_t$ is smooth (if $\mu_t$ is smooth). 
 \end{remark}


\subsection{Tameness of families of hybrid functions} \label{ss:TameHybridFunctions} In the following, we introduce two notions of tameness for functions on families of hybrid curves: the {\em stratumwise tameness} and {\em weak tameness}. These definitions provide the hybrid analog of the corresponding notions for tropical curves in Section~\ref{ss:TameFunctions}.

\smallskip

The idea behind these notions is to give a mathematical meaning to the convergence of functions $f$ on smooth Riemann surfaces $S$ to hybrid functions $\lf$ on hybrid curves $\curve$. Since in contrast to a usual function $f \colon S \to \C$, a hybrid function $\lf =(f_1, \dots, f_r, f_\smallcc)$ has several components, classical continuity notions obviously do not apply. Vaguely speaking, the tameness of a family of functions $(f_\thy)_{\thy \in B^\hyb}$ on a hybrid deformation family $\rsf^\hyb \to B^\hyb$ means that, when a Riemann surface $\rsf_t$ degenerates to a hybrid curve, the function $f_t$ admits a \emph{logarithmic asymptotic expansion in terms of the components of hybrid curve functions}. This idea is formalized by using the notion of pullbacks from the previous section.

\smallskip

Notations as above, let $\rsf^\hyb \to B^\hyb$ be a hybrid versal family associated to a stable marked curve $S_0$. Assume further that an adapted system of coordinates on $\rsf \to B$ is given.

\subsubsection{Stratumwise tameness}

Let $D_\pi^\hyb$ be a hybrid stratum in the hybrid base space $B^\hyb$, associated to an order partition $\pi= (\pi_1, \dots, \pi_r)$ on some subset $F \subseteq E$.

\smallskip
Consider a  \emph{family of hybrid functions}  $\mathscr{F} = (\lf_\thy)_\thy$,  $\thy \in B^\hyb \cup D_\pi^\hyb$.  That is, for each point $\thy \in B^\hyb \cup D_\pi^\hyb$, we have given a function $\lf_\thy$ on the hybrid curve $\rsf^\hyb_\thy$.

Suppose that the functions $\lf_\thy$ live fiberwise in $D_{\log}(\Deltahyb)$ of each hybrid curve.

\smallskip

Using the hybrid log map $\loghyb_\pi$, we pullback the family $\mathscr{F}_\pi := \mathscr{F}\rest{D_\pi^\hyb}$ from the boundary stratum $D_\pi^\hyb$ to a family of functions $\mathscr{F}_\pi^\ast =(\lf^\ast_t)_{t \in B^\ast} $ on $\rsf^\ast$ (see Section~\ref{sec:pullback_rs}). 

\smallskip

A family of hybrid functions $\mathscr{F} = (\lf_\thy)_{\thy \in B^\hyb \cup D_\pi^\hyb}$ is \emph{tame at the stratum $D_\pi^\hyb$} if the difference
\begin{equation} \label{eq:DifferenceStratum}
\lf_t - \lf^\ast_{t}, \qquad  t\in B^\ast,
\end{equation}
goes to zero uniformly upon approaching the closed boundary stratum $\overline{D}_\pi^\hyb$ in the tame topology. That is, for each point $\thy$ in the closure $\overline{D}_\pi^\hyb$, we have that
\[
\sup_{x \in \rsf_t} \big |\lf_t(x) - \lf_t^\ast(x) \big| \to 0
\]
 if $t \in B^\ast$ tamely converges to $\thy$ in $B^\hyb$.

\smallskip

The above property means that the pullback family $\mathscr{F}_\pi$ approximates the original family $\mathscr{F}$ close to the boundary stratum $D_\pi^\hyb$. By \eqref{eq:StratumwisePullback}, it implies the following asymptotic expansion
\[
\lf_t = \sum_{j=1}^r L_j(t) \Big( f_{\thy, j}^*\circ \loghyb_{\pi}  \Big) + f_{\thy,  \smallcc}^*\circ \loghyb_{\pi} + o(1), \qquad \thy = \loghyb_{\pi}(t),
\]
of the function $\lf_t$, $t \in B^\ast$, close to a hybrid point $\shy \in D_\pi^\hyb$ in the tame topology on $B^\hyb$.

\begin{defi}[Stratumewise tameness]
Let $\mathscr{F} = (\lf_\thy)_{\thy \in B^\hyb}$ be a family of hybrid functions on $\rsf^\hyb$ living fiberwise in $D_{\log}(\Deltahyb)$.
We say that $\lf$ is \emph{stratumwise tame} (or simply tame), if the restriction of $\mathscr{F}$  to $B^\ast \cup D_\pi^\hyb$ is tame at all boundary strata $D_\pi^\hyb$ of $B^\hyb$.   
\end{defi}

\begin{remark} Analogous to the tropical setting (see Remark~\ref{rem:TamenessVS}), it can be shown that every continuous, stratumwise tame family of functions $\lf_\thy$, ${\thy \in B^\hyb}$, is also weakly tame in the sense of the next section. Again, the validity of this connection is the main reason for allowing the point $\thy$ in \eqref{eq:DifferenceStratum} to belong to the {\em closure of the stratum} $D_\pi^\hyb$. Indeed, it can be shown that the implication might become wrong if just the open stratum $D_\pi^\hyb$ is considered.
\end{remark}


\subsubsection{Weak tameness} \label{ss:HybridWeakTame}

Let $\mathscr{F} = (\lf_\thy)_\thy$ be a {\em family of hybrid functions} on $\rsf^\hyb$. That is, for each $\thy \in B^\hyb$, we have a function $\lf_\thy$ on the hybrid curve $\rsf^\hyb_\thy$. Suppose that the functions $\lf_\thy$ live fiberwise in the extended Arakelov--Zhang space $D_{\log}(\Deltahyb)$ of each hybrid curve.

\smallskip
Consider the global log maps $\Lognoind := \loghyb$ defined on the family $\rsf^\hyb$ and the base space $B^\hyb$. The stratification $\partial_\infty B^\hyb = \bigsqcup_\pi D_\pi^\hyb$ of the hybrid boundary $\partial_\infty B^\hyb$ induces the following decomposition of the domain of definition $\mathscr{U}$ of the hybrid log map $\Lognoind$,
\[
\mathscr{U} = \bigsqcup_\pi \inn R_\pi, \qquad \qquad \inn R_\pi := \Lognoind^{-1}(D_\pi^\hyb),
\]
where the union is over all ordered partitions $\pi \in \Pifs(F)$ on subsets $F \subseteq E$.

\smallskip

Assume further that $\layh_\thy$, $\thy \in B^\hyb$, are hybrid functions with the following properties:
\begin{itemize}
\item[(i)] for each $\thy \in B^\hyb$, the function $\layh_\thy$ is a hybrid function \emph{on the hybrid curve} $\rsf^\hyb_\shy$ for $\shy = \Lognoind(\thy)$.
\item[(ii)] The equality  $\layh_{\thy} = \lf_\thy$ holds for all $\thy \in \partial_\infty B^\hyb$.
\item[(iii)] Consider the boundary stratum $D_\pi^\hyb$ for an ordered partition $\pi =(\pi_1, \dots, \pi_r) \in \Pifs(F)$ on some subset $F \subseteq E$. By property (i), the hybrid function $\layh_\thy$ is of the form $\layh_\thy = (h_{\thy, 1}, \dots, h_{\thy, r}, h_{\smallcc, 1})$   for each base point $\thy$ in the region $\inn R_\pi$.  Then for every $j=1, \dots, r, \smallcc$,  the $j$-th component $h_{\thy, j}$ of $\layh_{\thy}$ depends continuously on $\inn R_\pi $.
\end{itemize}

By the continuity in property (iii), we mean the following. Using the log map $\Lognoind \colon \inn R_\pi \to D_\pi^\hyb$,  we get $r$ families $\mathscr{G}_1, \dots, \mathscr{G}_r$ of metric graphs over $\inn R_\pi$ (the fiber $\mathscr{G}_{j, \thy}$ is the $j$-th graded minor of the hybrid curve $\rsf^\hyb_\shy$ for $\shy = \Lognoind(\thy)$); and a family $\~{\rsf}_\pi$ of disjoint Riemann surfaces over $\inn R_\pi$  (the fiber $\~{\rsf}_{\pi, \thy}$ is the normalization of the stable Riemann surface underlying $\rsf^\hyb_\shy$ for $\shy = \Lognoind(\thy)$).

The continuity means that the functions $h_{\thy, j}$, $\thy \in \inn R_\pi$ form a continuous function on $\mathscr{G}_j$, $j = 1, \dots, r$; moreover, the complex components $h_{\thy, \smallcc}$, $\thy \in \inn R_\pi$, form a continuous function on $\~{\rsf}_{\pi}$. In case that the functions $h_{\thy, \smallcc}$ have logarithmic poles, these conditions are slightly relaxed: we require that there are continuous sections $p^1_\thy, \dots, p^m_\thy$ of $\~{\rsf}_{\pi} / \inn R_\pi$ such that the poles of $h_{\thy, \smallcc}$ are contained in $\{ p^1_\thy, \dots, p^m_\thy\}$, the functions $h_{\thy, \smallcc}$, $\thy \in \inn R_\pi$, form a continuous function on $\~{\rsf}_{\pi}$ away from the image of these sections. 

\smallskip

For a base point $t \in B^\ast$, denote by $\layh^\ast_t$ the pull-back of $\layh_t$ from the hybrid curve $\rsf^\hyb_\shy$, $\shy = \Lognoind(t)$, to the smooth Riemann surface $\rsf_t$ (see \eqref{eq:SimplestPullback}).

\begin{defi}[Weak tameness] \label{def:WeakGlobalTameHybrid}
A family of hybrid functions $\lf_\thy$, $\thy \in B^\hyb$, is {\em weakly tame} if there are hybrid functions $\layh_\thy$, $ \thy \in B^\hyb$, as above and such that, in addition, the difference 
\begin{equation} \label{eq:Differencev}
 \lf_t  -  \layh^\ast_t , \qquad t \in B^\ast,
 \end{equation}
goes to zero upon approaching the hybrid boundary $\partial_\infty B^\hyb$ in the tame topology.
\end{defi}
More precisely, we require that for any point $\shy \in \partial_\infty  B^\hyb$,
\[
	\lim_{t \to \shy} \sup_{x \in \rsf^\hyb_t} \big |  \lf_t  (x)  - \layh^\ast_t (x) \big | = 0
\]
if the point $t \in B^\ast$ converges to $\shy$ in the tame topology on $B^\hyb$. By \eqref{eq:GlobalPullback}, weak tameness implies that the function $\lf_t$, $t \in B^\ast$, has the following asymptotic expansion
\[
\lf_t  = \sum_{j=1}^{r_\thy} L_j(t) \Big( h_{ \thy, j}^*\circ\Lognoind   \Big) + h_{\thy, \smallcc}^*\circ \Lognoind  + o(1), \qquad \thy = \Lognoind (t),
\]
close to a hybrid point $\shy \in \partial_\infty  B^\hyb$  in the tame topology on $B^\hyb$. Note that the number of terms $r_\thy$ in the additive expansion is the rank of the hybrid curve $\rsf^\hyb_\thy$, $\thy = \Lognoind  (t)$, and hence \emph{depends on the point} $t \in B^\ast$.


\subsection{Hybrid Laplacian} \label{ss:HybridLaplacian} In this section, we introduce the hybrid Laplacian.

\smallskip

Recall first that, on a smooth Riemann surface $S$, the Laplacian is the following operator $\Delta = \Delta_\smallcc =\frac 1{\pi i}\partial_z \bar\partial_z$ mapping smooth functions $f \colon S \to \C$ to smooth $(1,1)$-forms on $S$. We usually consider a slightly weaker definition and view $\Delta$ as an operator with values in the space of measures on $S$,  defined on $D_{\log}(\Delta_\smallcc)$, the space of $\mathcal{C}^2$-functions with finitely many logarithmic singularities. For any function $f \in D_{\log}(\Delta_\smallcc)$ with logarithmic poles $p_1, \dots, p_n \in S$,
\[
\Delta f  = \Delta_\smallcc f \rest{S \setminus\{p_1, \dots, p_n\}} - \sum_{j=1}^n r_j \delta_{p_j},
\]
where $r_j$ is the coefficient from \eqref{eq:RegularizedFunction}. In particular, the equality
\[
\int_S  h \, \Delta f = \int_S f \, \Delta h
\]
holds true for all smooth functions $h \colon S \to \R$.
\smallskip 

Consider a hybrid curve $\curve$ with underlying graph $G=(V,E)$, ordered partition $\pi = (\pi_1, \dots, \pi_r)$ of full sedentarity on the edge set $E$, and edge length function $l\colon E \to (0, + \infty)$.  Denote by $\Gamma^j$, $j\in [r]$,  the graded minors of $\curve^\trop$ which are metric graphs with normalized edge lengths.
Recall also that we introduced a notion of layered measures on $\curve$ in Section~\ref{ss:LayeredMeasures}.

\medskip

\medskip

The \emph{hybrid Laplacian} $\Deltahyb=\Deltahybind{\curve}$ on $\curve$ is a measure-valued operator which maps hybrid functions $\lf =(f_1, \dots, f_r, f_\smallc)$ in the extended Arakelov--Zhang space $D_{\log}(\Deltahybind{\curve})$ to layered measures on $\curve$.  It is defined as the sum
\begin{equation} \label{eq:HybridLaplacian}
\Deltahyb(\lf) = \Deltahyb_1(f_1) + \Deltahyb_2(f_2) + \dots \Deltahyb_r(f_r) + \Deltahyb_\smallcc(f_\smallc)
\end{equation}
where the components in the above sum are given as follows. For $j\in [r]$, and for a function $f_j$ on the graded minor $\Gamma^j$, the layered measure $\Deltatrop_j(f_j)$ on the hybrid curve $\curve$ is given by
\[
\Deltatrop_j(f_j)  := \Bigl(0, \dots, 0, \Delta_j (f_j),  \divind{j}{j+1}(f_j), \dots,  \divind{j}{r}(f_j),\divind{j}{\smallcc} (f_j) \Bigr).
\]
Here, as in Section~\ref{sec:tropical_laplacian}, $\Delta_j (f_j) \in \mathcal{M}^0(\Gamma^j)$ denotes the Laplacian on the $j$-th graded minor $\Gamma^j$. For the higher minors $\Gamma^i$ with $i>j$, $i\in[r]$, the point measure $\divind{j}{i} (f_j)$ on $\Gamma^i$ is defined by
\[
\divind{j}{i}(f_j)  := -  \sum_{\substack{e \in \pi_j} } \sum_{\substack{v \in e}} \slp_e f_j (v) \, \delta_{\proj_i(u)},
\]
where, as before, $\proj_i$ is the projection map $\proj_i:V\to V(\Gamma^j)$, and $\slp_e f_k(u)$ is the slope of $f_k$ at $u$ along the unit tangent direction at the incident edge $e$.

Finally, for $j\in [r]$,
\begin{equation} \label{eq:DefineDivisorsHybridLaplacian}
\divind{j}{\smallcc}(f_j)  := -  \sum_{\substack{e \in \pi_j} } \sum_{\substack{v \in e}} \slp_e f_j (v) \, \delta_{p^e_v},
\end{equation}
and 
\[\Deltahyb_{\smallcc}(f_\smallc) := \Bigl(0, \dots, 0, \Delta_\smallcc (f_\smallc)\Bigr)\] 
where $\Delta_\smallcc =\frac 1{\pi i}\partial_z \bar\partial_z$ is the Laplacian on the complex part $\pi_\smallcc = \bigsqcup_v C_v$.

\smallskip

\noindent The following proposition is immediate from the definition of the hybrid Laplacian $\Deltahyb$.
\begin{prop} \label{prop:hybridLaplacianMassZero}
For any function $\lf=(\lf^\trop, f_\smallc)$ on a hybrid curve $\curve$, the measure $\Deltahyb(\lf)$ is a layered measure of mass zero on $\curve$. 
\end{prop}

In Section~\ref{sec:hybrid_green} we will discuss the Poisson equation on hybrid curves.

\subsection{Hybrid Laplacian as a weak limit of Laplacians on Riemann surfaces} \label{ss:HybridLaplacianAsLimit}
In this section, we show that the hybrid Laplace operator $\Deltahyb$ arises as the limit of the Laplace operator on degenerating Riemann surfaces (see Section~\ref{ss:TropicalLaplacianAsLimit} for the parallel result on the Laplacian on tropical curves). Notations as in the previous sections, consider the versal deformation family $\rsf \to B$ of a stable marked Riemann surface $S_0$, together with a fixed adapted system of coordinates. Let $\rsf^\hyb \to B^\hyb$ be the corresponding family of hybrid curves.

\smallskip

Consider a fixed point $\thy$ in the hybrid boundary $\partial_\infty B^\hyb$ and the corresponding hybrid curve $\curve := \rsf_\thy^\hyb$, with underlying ordered partition $\pi = (\pi_1, \dots, \pi_r) \in \Pi(F)$ on some $F \subseteq E$. Let  $\lf = (f_1, \dots, f_r, f_\smallc)$ be a hybrid function on $\curve$ in the Arakelov--Zhang space $D(\Deltahybind{\curve})$.

\smallskip
By the construction in Section~\ref{sec:pullback_rs} (see \eqref{eq:SimplestPullback}), we can propagate $\lf$ to the smooth Riemann surfaces $\rsf_t$, with $t$ belonging to 
\[
B_\thy^\ast = \{t \in B^\ast | \, \text{$t$ satisfies \eqref{eq:CompatibleFiberMap2}} \},
\]
\emph{however, using the continuous log maps $\tilde \psi^t_\thy \colon \rsf_t \to \rsf^\hyb_\thy$ from Section~\ref{sss:ContinuousLogMap} instead of $\psi^t_\thy$}. That is, for any $t \in B^\ast_\thy$ the pullback $f^*_{t}\colon \rsf_{t} \to \C$ is given by
\[
f^*_{t} := \sum_{j=1}^r L_j(t) \, \big( f_j^*\circ \tilde \psi^t_\thy  \big) + f_\smallc^*\circ \tilde \psi^t_\thy,
\]
where the scaling factors $L_j(t)$, $j=1, \dots, r$, are defined by $L_j(t) = - \sum_{e\in \pi_j} \log|z_e(t)|$.
 
 \smallskip
 
 We have the following result on the convergence of the Laplacian $\Delta$ on Riemann surfaces to the hybrid Laplacian $\Deltahyb$. Analogous to the tropical Laplacian (see Theorem~\ref{thm:GraphLaplacianConvergence}), it justifies our definition of hybrid functions $\lf = (f_1, \dots, f_r, f_\smallc)$, their rescaled pullbacks $\lf^\ast$ and the definition of the hybrid Laplacian $\Deltahyb$.

\begin{thm} \label{thm:RiemannSurfaceLaplacianConvergence}
As $\rsf_t$ degenerates to $\curve$, that is, as $t \in B^\ast_\thy$ tends to $\thy$ in $B^\hyb$,
\[
\Delta(f^\ast_t)  \to \Deltatrop(\lf)
\]
in the weak sense.
\end{thm}

Convergence in the weak sense is defined as before. Namely, we require for any continuous function $h\colon \rsf^\hyb\rest{B^\ast_\thy \cup \{\thy\}} \to \R$ the convergence 
\begin{equation} \label{eq:ExplainLaplacianConvergence}
 \int_{\rsf_t} h(s) \,  \Delta f_t \quad  \longrightarrow  \quad \int_{\curve} h(s)  \, \Deltahyb \lf
\end{equation}
as $t \in B^\ast_\thy$ tends to $ \thy$. On the right hand side,  $\Deltahyb \lf$ is understood as a measure of mass zero on the hybrid curve $\curve$. (see Proposition~\ref{prop:MassZeroMeasures} and Proposition~\ref{prop:hybridLaplacianMassZero}).

\begin{remark}[Validity of the theorem for the former pullbacks]\label{rem:validity-weak-convergence} For ease of formulation, in Theorem~\ref{thm:RiemannSurfaceLaplacianConvergence}, we have restricted to the continuous map $\tilde \psi^t_\thy$ instead of $\psi^t_\thy$  and hybrid functions $\lf$ in $D(\Deltahyb)$ instead of $D_{\log}(\Deltahyb)$. More generally, we can consider the pullbacks from \eqref{eq:SimplestPullback},
\[
f^*_{t} := \sum_{j=1}^r L_j(t) \, \big( f_j^*\circ \psi^t_\thy  \big) + f_\smallc^*\circ  \psi^t_\thy,
\]
for any $t \in B^\ast_\thy$ and allow that $\lf$  belongs to $D_{\log}(\Deltahyb)$, that is, the complex part $f_\smallcc$ may have logarithmic poles.

However, then the pullbacks $f_t^\ast$ are no longer continuous on their respective fibers $\rsf_t$. As a result, the Laplacian $\Delta f_t^\ast$ does not exist as a measure on $\rsf_t$. Namely, integration by parts in $\int_{\rsf_t} f_t^\ast \Delta h$ gives also contributions of the derivatives of $h$, integrated on the cycles $\realtor^e_v$ between the regions $B_{e,t}$ and $A^e_{v,t}$.

Note however that $\Delta f_t^\ast$ still exists as a distribution defined on smooth functions $h \colon \rsf_t \to \C$. Moreover, it turns out that \eqref{eq:ExplainLaplacianConvergence} remains true if the test function $h \colon \rsf^\hyb \to \C$ is sufficiently regular and the left hand side is interpreted as a pairing between smooth functions and distributions on $\rsf_t$. More precisely, we additionally require that the restrictions $h_t := h \rest{\rsf_t}$ are smooth on each smooth fiber $\rsf_t$, $t \in B^\ast_\thy$, and their derivatives satisfy $\frac{\partial h_t}{ \partial |z^e_u|} \le C$ on all boundary cycles $\realtor^e_u$, for some constant $C >0$ independent of $t \in  B^\ast_\thy$.
\end{remark}
In the following proof, we use the notation from Section~\ref{sss:ContinuousLogMap}.

\begin{proof}[Proof of Theorem~\ref{thm:RiemannSurfaceLaplacianConvergence}] 
Let $l\colon E \to (0, + \infty)$ be the normalized edge length function of the hybrid curve $\curve$. Recall that, by definition, the Laplacian of a sufficiently smooth function $f$ on a Riemann surface $S$ is the measure $\Delta f$  such that
\begin{equation} 
	\int_S h \, \Delta f = \frac{1}{\pi i}\, \int_S f  \, \partial \bar \partial h  
\end{equation}
for all smooth functions $h \colon S \to \R$. Recall also that in local coordinates, the $2$-form $\frac{1}{\pi i} \partial \bar \partial h$ can be written as
\begin{align*}
\frac{1}{\pi i} \partial \bar \partial h = - \frac{1}{2 \pi} \left ( \frac{d^2 h}{dx^2} + \frac{d^2 h}{dy^2}  \right) dx \wedge dy =  \frac{1}{2 \pi}  \left ( \frac{d^2 h }{d \rho^2} + \frac{d^2 h }{d\theta^2}  \right) d\rho \wedge d\theta
\end{align*}
in Cartesian coordinates $(x,y)$ and logarithmic polar coordinates $(s, \theta)$, that is for $s = - \log|z|$, $z=x+i y$.

\smallskip

Taking into account the (additive) definition of $\Deltahyb = \Deltahyb_1 + \dots \Deltahyb_r+\Deltahyb_\smallc$, it suffices to prove the claim in the case that $\lf = (0, \dots, 0, f_j, 0, \dots, 0)$ for $j\in[r]\cup\{\smallcc\}$. We abbreviate $ \ell_e(t) = - \log|z_e(t)|$ for $e \in E$ and $L_j(t) =\sum_{e\in \pi_j}\ell_e(t)$ for $j \in [r]$.

\smallskip

(i) Assume first that $j \in [r]$. Since the function $f_j$ is constant on all connected components of $\curve \setminus \bigsqcup_{e \in \pi_1 \cup \dots \cup \pi_j} I_e$, the pull-back $f_j^\ast$ to the Riemann surface $\rsf_t$ is constant on all components of $\rsf_t  \setminus \bigsqcup_{e \in \pi_1 \cup \dots \cup \pi_j} B_e$. We use log polar coordinates $\bigl((\rho^e_u,\rho^e_v), \theta\bigr)$ with $\rho^e_u = - \log|z^e_u|, \rho^e_v = - \log|z^e_v|$ on the sets $B_e$. The pull-back $f_j^\ast$ is equal to $f_j(\tilde x^e_u, \tilde x^e_v)$ for $(\tilde x^e_u, \tilde x^e_v)$ given in~\ref{sss:ContinuousLogMap} by
  \[\tilde x^e_u =  \frac{l(e)}{\ell_e(t)-2 \tau_e(t)} (\rho^e_u -\tau_e) \qquad \textrm{and} \qquad \tilde x^e_v =  \frac{l(e)}{\ell_e(t)-2\tau_e(t)} (\rho^e_v -\tau_e).\]
on each set $B_e$. Hence, $f_j^\ast$ is radially symmetric.  Moreover, if $e\in \pi_i$ with $i<j$, then $f_j^\ast$  is affine linear in $\rho^e_u$ .

\smallskip

By the definition of the hybrid Laplacian, and using Stokes' theorem, we infer that $\Delta(f_j^\ast)$ is the following measure on $\rsf_t$:
\[
 \Delta f_j^\ast = \sum_{e \in \pi_j}  \frac{1}{2 \pi}  \frac{l(e)^2}{\bigl(\ell_e -2 \tau_e \bigr)^2} \, \frac{d^2 f_j\rest{e}}{{{d\tilde x^{e}_v}^2}} \Big ( \tilde x^{e}_v \Big ) \cdot \lambda_{B_e} +  \sum_{e \in \pi_1 \cup \dots \cup \pi_j} \sum_{v \in e} a^e_v \cdot \delta_{\combind{\realtor^e_v}}
\]
where $\lambda_{B_e}$ is the measure associated with the $2$-form $d \rho^e_v \wedge d \phi$ on $B_e \cong [\tau_e, \ell_e - \tau_e] \times \realtor$, $\realtor^e_{v}$ is the circle between $A^e_v$ and $B_e$ (corresponding to $|z^e_v| = \varrho_e$) and $\delta_{\combind{\realtor^e_{v}}}$ is the uniform measure on $\realtor^e_{v}$ of total mass $2 \pi \rho_e$ with coefficient 
\[
a^e_v = - \frac{1}{2 \pi \varrho_e} \frac{l(e)}{\ell_e - 2 \tau_e} \, \, \slp_e f_j (v).
\]  
Integrating against a continuous function $h \colon \rsf^\hyb\to \R $, one computes that
\begin{align*}
\int_{\rsf_t} h \,  \Delta f_j^\ast &= \sum_{e \in \pi_j}   \frac{ - l(e)}{\ell_e -2\tau_e} \int_{I_e}  \left ( \frac{1}{2\pi} \int_{\combind{\realtor}} h \rest{B_e} \big (\varphi^{-1}(\tilde x^e_v), \theta \big ) \, d\theta \right ) \frac{d^2 f_j\rest{e}}{{d{\tilde x}^e_v}^2} (\tilde x^e_v) \,  d \tilde x^e_v \\
&+ \sum_{e \in \pi_1 \cup \dots \cup \pi_j} \sum_{v \in e} a^e_v \int_{\combind{\realtor^e_{v}}} h (\theta)  \, \delta_{\combind{\realtor^e_{v}}}  (\theta).
\end{align*}
(Here we, for simplicity, assume that the function $f_j \colon \Gamma^j \to \C$ is $\mathcal{C}^2$ on all edges $e$ of $\pi_j$; if $f_j$ is only piecewise $\mathcal{C}^2$ on $e$, we would have to add additional integrals of $h$ over spheres $\realtor$, corresponding to the positions where the derivative $f_j'$ is discontinuous.)
\smallskip

Recall that, as $t \to \thy$, we have that $\ell_e(t) \to \infty$ for $e \in \pi_1 \cup \dots \cup \pi_r$ and
\begin{equation} \label{eq:RelHybrid}
\ell_e(t) = L_i(t) \big (l_e + o(1) \big ) \quad \text{$e \in \pi_i$, $i \in [r]$}, \qquad L_i (t) / L_j(t) = o(1) \quad \text{for $i > j$}.
\end{equation}
Integrating against a continuous function $h \colon \rsf^\hyb \to \R$, we arrive at the limit
\begin{align*}
\lim_{t \to \thy} \int_{\rsf_t} h  \, \Delta (L_j f_j^\ast) &= \sum_{e \in \pi_j} \int_{I_e} - \frac{d^2 f_j}{dy^2} h(y) dy + \sum_{e \in \pi_j} \sum_{v \in e} - \slp_e f_j (v) h(v) = \int_{\curve} h \, \Deltahyb \lf,
\end{align*}
where the last equality follows from \eqref{eq:Gr/NGrMea} and the definition of the hybrid Laplacian.

\smallskip

(ii) It remains to treat the case where $\lf = (0, \dots, 0,  f_\smallc)$. Using Stokes' theorem once again, we infer that the Laplacian $\Delta(f_\smallc^\ast)$ is the following measure on $\rsf_t$,
\begin{equation} \label{eq:LaplacianComplexPullback}
\Delta(f_\smallc^\ast) = \sum_{v \in V} \frac{\partial \bar \partial}{\pi i}  f_\smallc \rest{Y_{v,t}} + \sum_{e \in \pi_\fin} \frac{ \partial \bar \partial}{\pi i} f_\smallc^\ast\rest{W_{e,t}} + \sum_{v \in e} \widetilde{a}^e_v \cdot \delta_{\combind{{\widetilde \realtor}^e_{v}}}  + \sum_{j=1}^r \sum_{e \in \pi_j} \sum_{v \in e} a^e_v \cdot  \delta_{\combind{\realtor^e_{v}}}  + \widetilde{a}^e_v \cdot \delta_{\combind{{\widetilde \realtor}^e_{v}}},  \end{equation}
where for any vertex $v$ on an edge $e \in E$, ${\widetilde \realtor}^e_{v}$ denotes the cycle between the two regions $Y_{v,t}$ and $W_{e,t}$ (corresponding to $|z^e_v|=1$), $\delta_{\combind{{\widetilde \realtor}^e_{v}}}$ is the uniform measure on ${\widetilde \realtor}^e_{v}$ of total mass $2 \pi$ and the density functions $a^e_v \colon \realtor^e_v \to \C$ and $\widetilde a^e_v \colon \widetilde \realtor^e_v \to \C$ are given by
\begin{align*}
	{a}^e_v (\theta) &= -\frac{1}{2 \pi \varrho_e} \Big (\frac{l(e)}{\ell_e - 2 \tau_e} \slp_e f_\smallcc(v) + \frac{\varrho_e}{1-\varrho_e} \frac{\partial f_\smallc}{\partial r} (0, \arg(\theta)  ) \Big ),  &&\theta \in  {\realtor}^e_{v}, \\[1mm]
	 \widetilde{a}^e_v (\theta) &= \frac{1}{2\pi} \frac{\varrho_e}{1 - \varrho_e} \, \frac{\partial f_\smallc}{\partial r} (1, \arg(\theta)), &&\theta \in {\widetilde \realtor}^e_{v} ,
\end{align*}
for $e \in E \setminus \pi_\fin$, and
\[
{a}^e_v (\theta)  = \frac{1}{2 \pi} \Big (\frac{\log|z_e(s)|}{\log|z_e(t)|}  - 1\Big ) \frac{\partial f_\smallc}{\partial r} (1, \arg(\theta)) , \qquad \theta \in {\widetilde \realtor}^e_v \]
for $e \in \pi_f$. Here we have used polar coordinates $\rho^e_v= -\log|z^e_v|$ and $\theta = \arg(z^e_v)$ on the region $W_{e,t}$. Taking a continuous function $h \colon \rsf^\hyb \to \R$, one verifies that
\begin{align*}
\lim_{t \to \thy} \sum_{v \in V}  \int_{Y_{v,t}} h \, \partial \bar \partial f_\smallc +    \sum_{e \in \pi_\fin} \int_{W_{e,t}} h  \, \partial \bar \partial f_\smallc = \sum_{v \in V} \int_{C_v} h \, \partial \bar \partial f_\smallc = \pi i \int_{\curve} h \, \Deltahyb \lf,
\end{align*}
where $C_{v}$, $v \in V$ denote the Riemann surface components of the hybrid curve $\curve$, and  the last equality follows from \eqref{eq:Gr/NGrMea}. However, by the asymptotics obtained in \eqref{eq:RelHybrid},
\[
	\sup_{\theta \in {\widetilde S}^e_{v}} \big |\widetilde{a}^e_v (\theta) \big | = o(1) ,\quad  e \in E \qquad \text{and} \qquad \sup_{\theta \in {S}^e_{v}} \big |{a}^e_v (\theta) \big | = \frac{o(1)}{\varrho_e(t)}  ,\quad  e \in E \setminus \pi_\fin
\]
as $t \to \thy$ in $B^{\hyb}$. It follows that the integrals of $h$ over the remaining measures in \eqref{eq:LaplacianComplexPullback} vanish in the limit $t \to \thy$. The proof is complete.

\end{proof}


\section{Hybrid Poisson equation}\label{sec:hybrid_green} In this section, we introduce the \emph{Poisson equation on hybrid curves} and study the associated solutions. We furthermore state our main theorems on the relationship of these solutions to limits of solutions of the Poisson equation on degenerating families of Riemann surfaces.  The proof of these theorems will be given in the upcoming sections.


\subsection{Poisson equation on a hybrid curve} \label{eq:GeneralHybridCurve}
We introduce the following general \emph{Poisson equation on hybrid curves}.

\smallskip

Let $\curve$ be a hybrid curve of rank $r$ with underlying partition $\pi$.  Consider two layered Borel measures $\lmu$ and $\lnu$ on $\curve$ such that $\lmu$ is of mass zero, and $\lnu$ has mass one (meaning that it has mass one on the first graded minor $\Gamma^1$ and mass zero on all the components of any higher indexed graded minor, including on each Riemann surface $C_v$, $v\in V(G)$). Altogether, we refer to the triplet $(\curve, \lmu, \lnu)$ as a \emph{bimeasured hybrid curve}.

\smallskip

In order to have a well-posed equation in the extended Arakelov--Zhang space $D_{\log}(\Deltahyb)$, we require that the pieces $\mu_j$, $j \in [r]$, of the measure $\lmu$ are of the form \eqref{eq:SpecialMus} on the metric graphs $\Gamma^j$. Moreover,  we suppose that the complex parts of the measures are of the form
\begin{equation} \label{eq:NiceHybridMeasures}
\mu_\smallcc = \alpha_\mu + D_\mu, \qquad \nu_\smallcc = \alpha_\nu + D_\nu
\end{equation}
for some continuous $(1,1)$-form $\alpha_\mu$ and $\alpha_\nu$ and some (real) divisors $D_\mu$ and $D_\nu$ on the components $C_v$, $v \in V$. In addition, we assume that these divisors have disjoint supports.

\smallskip

If $\curve = S$ is a smooth Riemann surface, we recover the familiar setting of a smooth Riemann surface $S$ equipped with two measures $\mu$ and $\nu$. The associated Poisson equation
\[
 \begin{cases}
\Delta f = \mu \\
\int_S f \, d\nu = 0
\end{cases}
\]
is \emph{well-posed}, that is, it has a unique solution $f$ on $S$. Indeed, since $S$ is connected, solutions to $\Delta f = \mu$ are unique up to a constant, which is fixed by the second equation.

\medskip

We extend this to general bimeasured hybrid curves $(\curve, \lmu, \lnu)$. Consider first the equation
\begin{equation} \label{eq:PoissonHybridCurveNotUnique}
\Deltahyb \lf = \lmu
\end{equation}
on the hybrid curve $\curve$. By definition of the hybrid Laplacian, \eqref{eq:PoissonHybridCurveNotUnique} is a coupled system of Poisson equations
\begin{equation} \label{eq:IterativeLaplacian}
\begin{cases}
\Delta_{\Gamma^1}(f_1)= \mu_1 & \text{on $\Gamma^1$} \\[1mm]
\Delta_{\Gamma^2}(f_2)= \mu_2 - \divind{1}{2}(f_1)  & \text{on $\Gamma^2$}  \\[1mm]
\dots \\[1mm]
\Delta_\smallcc(f_\smallcc) = \mu_\smallcc - \sum_{j=1}^r  \divind{j}{\smallcc}(f_j)  & \text{on $\pi_\smallcc = \bigsqcup_v C_v$}
\end{cases}
\end{equation}
on the metric graphs $\Gamma^j$, $j \in [r]$, and on $\pi_\smallcc$. The next lemma ensures existence of solutions.

\begin{lem} \label{lem:ExistenceHybridSolutions}
The Poisson equation \eqref{eq:PoissonHybridCurveNotUnique} has a solution $\lf$ in the extended Arakelov--Zhang space $D_{\log}(\Deltahyb)$, for every layered measure $\lmu$ on $\curve$ of mass zero satisfying \eqref{eq:NiceHybridMeasures}.
\end{lem}
\begin{proof}
Iteratively solving the above equations for $f_1, \dots, f_\smallcc$, one obtains a solution. \end{proof}

However, since the spaces $\Gamma^j$, $j \neq 1$ and $\pi_\smallcc$ are in general not connected, the solutions to \eqref{eq:PoissonHybridCurveNotUnique} are far from being unique. The $j$-th part $f_j$, $j \in \{1, \dots, r, \smallcc\}$ of a solution $\lf $ is unique up to the choice of 
\[
n_j = 1 + \sum_{k < j} \abs{E\big (\gr_\pi^k(G)\big)} - \operatorname{genus}\big (\gr_\pi^k(G) \big ) = 1 + \sum_{k < j} \abs{\pi_k} - h_\pi^k
\]
constants, one for each connected component of $\Gamma^j$, $j \in [r]$, and of $\pi_\smallcc$ ($\Gamma^j$ has precisely $n_j$ connected components). In order to fix this, we impose additional conditions on $\lf$ based on harmonically arranged functions (see Section~\ref{ss:HarmonicallyArrangedHybridFunctions}). Note also that $n_j -1$ is precisely the number of linearly independent exact one-forms on the metric graph $\Gamma^j$. It turns out that this way of normalizing the hybrid solutions describes the behavior of  solutions to the Poisson equation on degenerating families of Riemann surfaces (see Theorem~\ref{thm:TameContinuityHybridCurves} and Theorem~\ref{thm:MainGeneralMeasures} below).

\smallskip

Let $(\curve, \lmu, \lnu)$ be a bimeasured hybrid curve. Fix moreover a local parameter $z^e_v$ in $C_v$ around the point $p^e_v$ for each $v \in V$ and $e \sim v$. We consider the following {\em hybrid Poisson equation} 
 \begin{equation}  \label{eq:PoissonHybridCurve2}
 \begin{cases}
 \Deltahyb \lf = \lmu  \\[1 mm]
 \lf \text{ is harmonically arranged} \\[1mm]
\int_{\curve} \lf \, d \lnu = 0
 \end{cases}
\end{equation}
in the extended Arakelov--Zhang space $D_{\log}(\Deltahyb)$.

\smallskip

The latter two conditions mean the following. By \eqref{eq:NiceHybridMeasures}, the complex part $f_\smallcc$ of any solution to $\Deltahyb \lf = \lmu$ is a $\mathcal{C}^2$-function on $\pi_\smallcc$, except for finitely many logarithmic singularities. Using regularization in the local parameters $z^e_u$, we can define the property of being harmonically arranged for $\lf$ (see Section~\ref{sec:hybrid_laplacian}). Moreover, we can identify each piece $f_j$ of the solution $\lf$ with a function $f_j \colon \mccan \to \C$ on the distinguished metrized complex $\mccan$ associated to $\curve$, that is, the one with normalized edge length on each layer  (see Section~\ref{sec:hybrid_laplacian}).   The third condition in \eqref{eq:PoissonHybridCurve2} means that
\[
\int_{\mccan} f_j \, d\nu = 0 \qquad \text{for all }j=1, \dots, r, \smallcc.
\]
Here $\nu$ is the measure of total mass one on the metrized complex $\mccan$ obtained from the layered measure $\lnu$. Note that, in general, the following integrals are different.
\[
\int_{\mccan} f_j \, d\nu \neq \int_{\Gamma^j} f_j \, d\nu_j.
\]

\begin{thm}[Existence and uniqueness of solutions of the hybrid Poisson equations] \label{thm:ExUniqHybridPoisson} The hybrid Poisson equation \eqref{eq:PoissonHybridCurve2} has a unique solution $\lf$ in the extended Arakelov--Zhang space $D_{\log}(\Deltahyb)$ for every bimeasured curve $(\curve, \lmu, \lnu)$. 
\end{thm}
\begin{proof}
The theorem is a direct consequence of Lemma~\ref{lem:ExistenceHybridSolutions} and Theorem~\ref{thm:HarmonicExtensionHybrid}.
\end{proof}

\begin{remark}[Dependency on the regularization] 
We stress that the solution to the hybrid Poisson \eqref{eq:PoissonHybridCurve2} \emph{depends on the choice of local parameters} $z^e_u$, $e \sim v$, since the regularized function values $f_\smallcc^\reg(p^e_v)$ do. However, the solutions $\lf$ and $\tilde \lf$ for two different choices $z^e_u$, and $\tilde z^e_u$ have the same tropical part, $\lf^\trop = \tilde \lf^\trop$. Their complex parts differ by a constant on each Riemann surface component of $\curve$, that is, $f_\smallcc - \tilde f_\smallcc$ is constant on punctured surfaces $C_v \setminus \bigl\{p^e_v\,\st \, e \sim v\bigr\}$ for all $C_v$, $v \in V$.  \end{remark}


\subsection{Poisson equation on degenerating families of Riemann surfaces: discrete measures}

Notations as above, let $\rsf^\hyb/B^\hyb$ be a hybrid versal family parameterized by the hybrid basis $B^\hyb$ associated to a stable marked curve $(S_0, q_1, \dots, q_n, p)$, and assume that an adapted system of coordinates $(z_e)_{e\in [N]}$ and $z^e_v$, for $e\sim v$, is given. The markings naturally produce continuous sections $q_i(\thy)$ and $p(\thy)$ of $\rsf^\hyb \to B^\hyb$. Namely, they map a hybrid base point $\thy = (l,t) \in B^\hyb$ to the respective marked points on the stable marked Riemann surface $\rsf_t$, which lie in the smooth part of $\rsf_t$ and can hence be viewed naturally as points on $\rsf^\hyb_\thy$.

Let $D$ be a divisor of degree zero with constant coefficients supported on marked points $q_1, \dots, q_n$. That is, consider a family of divisors $D_\thy = \sum_{j=1}^n d_j q_j(\thy)$, $\thy\in B^\hyb$, for constants $d_1, \dots, d_n$ which sum up to zero.

\smallskip

For each point $\thy\in B^\hyb$, consider the Poisson equation
\begin{equation}  \label{eq:DiscrPoissonRiemannSurfaces}
\begin{cases}
 \Deltatrop (\lf) = D_\thy \\[1 mm]
 \lf \text{ is harmonically arranged} \\[1mm]
 f_j (\proj_j(p_\thy)) = 0, \qquad j=1, \dots, r, \smallcc,
 \end{cases}
\end{equation}
on the hybrid curve $\rsf^\hyb_\thy$. The second conditions means that $\lf$ is harmonically arranged with respect to the
\emph{local charts on $\rsf^\hyb_\thy$ induced by the adapted system of coordinates on $\rsf^\hyb / B^\hyb$.} Moreover, $\proj_j$ is the projection from $V$ to $V^j = V(\gr^j_\pi(G))$, for $j\in [r]$, and $\proj_\smallcc$ is identity map on the disjoint union of the curves $C_v$, $v\in V$. 

\smallskip

Let $\lf_\thy$, $\thy \in B^\hyb$, be the unique solution to the above Poisson equation.

\begin{thm}[Tameness of solutions of the Poisson equation: discrete measures] \label{thm:TameContinuityHybridCurves} The solutions $\lf_\thy$, $\thy\in B^\hyb$, to \eqref{eq:DiscrPoissonRiemannSurfaces} form a stratumwise tame family of hybrid functions on $\rsf^\hyb$.
\end{thm} 

\begin{proof} This result follows from the asymptotic description of the $\jvide$-function in Theorem~\ref{thm:TamenessJFunc}.
\end{proof}


\subsection{Poisson equation on degenerating families of Riemann surfaces: general measures} \label{ss:SurfacesPoissonGeneralLimit}
We generalize the above theorem to continuous families of measures. In this case, we can prove the weak tameness. Notations as above, let $\rsf^\hyb/B^\hyb$ be a hybrid versal family associated to a stable marked curve $(S_0, q_1, \dots, q_n)$ together with a fixed system of adapted coordinates.

\smallskip 

Consider two families of (layered) measures $\lmu_\thy$ and $\lnu_\thy$ on $\rsf^\hyb_\thy$, $\thy \in B^\hyb$. Assume that $\lmu_\thy$ and $\lnu_\thy$ form continuous families of measures on $\rsf^\hyb / B^\hyb$ (satisfying for each $\thy$ the conditions in Section~\ref{eq:GeneralHybridCurve}). On each hybrid curve $\rsf_\thy^\hyb$, $\thy \in B^\hyb$, we consider the Poisson equation
 \begin{equation}\label{eq:poisson_general_measure}
 \begin{cases}
 \Deltatrop \lf_\thy = \lmu_\thy  \\[1 mm]
 \lf_\thy \text{ is harmonically arranged} \\[1mm]
\int_{\rsf_\thy^\hyb} \lf_\thy \, d \lnu_\thy = 0.
 \end{cases}
\end{equation} 
In order to interpret the second and third equation, we use the  \emph{local charts on $\rsf^\hyb_\thy$ induced by the adapted system of coordinates on $\rsf^\hyb / B^\hyb$.}

\smallskip

Moreover, we assume that the two measures satisfy the following \emph{uniform growth bounds}. For an edge $e \in E$, consider the cylinder $W_{e,t} \subseteq \rsf_t$ (see \eqref{eq:AdaptedCoordinates}) and its subregion
\[
\~{W}_{e,t} = \Big \{ (\underline{z}, z^e_u, z^e_v) \in W_{e,t} \, \st \, |z^e_u| \le |\log |z_e(t)||^{-1/2} \text{ and }  |z^e_v| \le |\log|z_e(t)||^{-1/2} \Big \}.
\]
Recall moreover that, in Section~\ref{sec:hybrid_log_map_einf}, the cylinder $W_{e,t}$ was decomposed into
\[
W_{e,t} = A^e_{u,t} \sqcup B_{e,t} \sqcup A^e_{v,t}.
\]
Note that $\~{W}_{e,t}$ intersects all three sets in this decomposition. We impose the following additional assumptions on the measures $\lmu_\thy$ and $ \lnu_\thy$.

\smallskip

\begin{itemize}
\item[$(1)$] (Uniformly bounded variation) The measures $\lmu_\thy$ and $\lnu_\thy$ are uniformly bounded in total variation.
\smallskip
\item[$(2)$] (Bounded densities) The densities of the measures $\mu_t= \lmu_t$ and $\nu_t = \lnu_t$ on the smooth fibers $\rsf_t$, $t \in B^\ast$, satisfy the following uniform bounds.

 \smallskip
 
(2a) In the coordinate $z^e_u$ on the cylinder $W_{e,t}$, the measures $\mu_t$ and $\nu_t$ are given by 
\[
\mu_t \rest{W_{e,t}}= f_t \, \, \sqrt{-1} dz^e_u \wedge d \overline{z^e_u}, \qquad \qquad \nu_t \rest{W_{e,t}} = h_t  \, \, \sqrt{-1} dz^e_u \wedge d \overline{z^e_u}
\] for continuous functions $f_t, h_t \colon W_{e,t} \to \R$. Moreover, for $t \in B^\ast$ close to a point $s \in D = B\setminus B^*$,
\[
 \max\big \{ |f_t(q)|, |h_t(q)| \big \} \le C \frac{1}{|\log|z_e(t)||} \frac{1}{|z^e_u(q)|^2}, \qquad q \in \widetilde{W}_{e,t},
 \]
 and
\[
 	 \max\big \{ |f_t(q)|, |h_t(q)| \big \}  \le C \Big (1 + \frac{1}{|\log|z_e(t)||} \frac{1}{|z^e_u(q)|^2} \Big )  , \qquad q \in A^e_{u,t},
\]
with a constant $C > 0$ independent of $t$.
 
 \smallskip
 
 (2b) Furthermore, consider a fixed point $p$ in the smooth part of a singular fiber $\rsf_s$ over $s \in D= B\setminus B^*$. Let $U$ be a small open neighborhood of $p$ in the family $\rsf \to B$, parametrized by $z_1, \dots, z_N, z$, where  $\underline z = (z_1, \dots, z_N)$ are the adapted coordinates from the base $B= \Delta^N$ and $z$ represents an additional parameter on the fibers of $\rsf \to B$.
 
 Then on the subset $U_t = U \cap \rsf_t$ of $\rsf_t$, $t \in B^\ast$, the measures $\mu_t$ and $\nu_t$ are given by 
\[
\mu_t \rest{U_t} = f_t \, \, dz \wedge d \overline z, \qquad \qquad \nu_t \rest{U_t} = h_t \, \, dz \wedge d \overline z,
\]
for continuous functions $f_t, h_t \colon U_t \to \R$.  Moreover, for $t \in B^\ast$ close to $s$ in $B$,
 \[
 	 \max\big \{ |f_t(q)|, |h_t(q)| \big \}  \le C, \qquad q \in U_t,
\]
with a constant $C >0$ independent of $t$.

\smallskip

\item[$(3)$] (Support of the point masses) The set of marked points of Riemann surfaces is partitioned into two parts $Q_\mu=\{q_1, \dots, q_a\}$ and $Q_\nu=\{q_{a+1}, \dots, q_{n}\}$ such that 
for each $\thy \in B^\hyb$, the discrete parts of the measures $\mu_\thy$ and $\nu_\thy$ on the corresponding metrized complex $\mccan_\thy$ lie in the disjoint parts $Q_\mu$ and $Q_\nu$ of the partition, respectively. 
\end{itemize}

\begin{remark}
As follows from Lemma~\ref{lem:CanonicalMeasureBound}, these properties hold in context with canonical measures.
\end{remark}

By Theorem~\ref{thm:ExUniqHybridPoisson}, the hybrid Poisson equation \eqref{eq:poisson_general_measure} has a unique solution $\lf_\thy$ for every base point $\thy \in B^\hyb$. The main theorem on continuity of the solutions reads as follows.
\begin{thm}[Weak tameness of solutions of the Poisson equation] \label{thm:MainGeneralMeasures} Notations as above, assume that the measures $\lmu_\thy$ and $\lnu_\thy$, $\thy \in B^\hyb$, additionally verify the properties (1)-(2)-(3).

Then, the family of solutions $\lf_\thy$, $\thy \in B^\hyb$, is weakly tame.
\end{thm}

More precisely, we will establish the following result. Fix a boundary stratum $D_\pi^\hyb$ in $B^\hyb$. For a point $t\in B^\ast$, let $\thy = \loghyb_\pi(t)$, and consider the hybrid log map
\[\loghyb_{\pi,t}\colon \rsf_t \to \rsf_\thy^\hyb.\]

Denote by $\loghyb_{t_*}(\mu_t)$ and $\loghyb_{t_*}(\nu_t)$ the push-out measures on $\rsf^\hyb_\thy$, and let $\hat\lf_t$ be the unique harmonically arranged solution to the Poisson equation
 \begin{equation} \label{eq:JDFcthyb}
\begin{cases}
 \Deltatrop (\hat \lf_t) = \loghyb_{t_*}\mu_t \\[1 mm]
 \hat \lf_t \text{ is harmonically arranged} \\[1mm]
\int_{\rsf^\hyb_\thy} \hat \lf_t \, \,  d  \loghyb_{t_*}\nu_t  = 0
 \end{cases}
\end{equation}
on the hybrid curve $\rsf^\hyb_\thy$, $\thy = \loghyb_\pi(t)$.

\begin{thm} \label{thm:WeakConvergence} Assume the conditions in Theorem~\ref{thm:MainGeneralMeasures} and suppose that $t \in B^\ast$ approaches $D_\pi^\hyb$, that is, $t$ converges tamely to a point $\shy$ in $D_\pi^\hyb$.

Then, the difference $f_t -\hat\lf^\ast_t$ goes to zero uniformly on $\rsf_t$.
\end{thm}

The proofs of Theorem~\ref{thm:MainGeneralMeasures} and Theorem~\ref{thm:WeakConvergence} are similar to the proof for the asymptotics of the Green function. They are given in the final Section~\ref{sec:Green_functions_asymptotics}.

\begin{remark}Under additional conditions on the speed of convergence of $\mu_t$ and $\nu_t$ to their hybrid limit measures (similar to \eqref{eq:condition_strong_tame}), the family of solutions $\lf_\thy$, $\thy \in B^\hyb$,  to \eqref{eq:poisson_general_measure} is stratumwise tame. However, we do not develop these results here.
\end{remark}


\section{Hybrid canonical measures and Arakelov Green functions} \label{sec:ArakelovGreenfunctionsSurfaces}
After recalling the definition of the canonical measure and Arakelov Green function on Riemann surfaces, we introduce the hybrid canonical measure and its corresponding Arakelov Green function on hybrid curves. The latter controls the asymptotics of the  Arakelov Green function on Riemann surfaces.

\subsection{Canonical measure on hybrid curves}
In this section, we recall the definition of the canonical measure on hybrid curves.

\subsubsection{Riemann surfaces} Let $S$ be a compact Riemann surface of positive genus $g$. Denote by $\Omega^1(S)$ the vector space of holomorphic one-forms $\omega$ on $S$. It has complex dimension $g$ and admits a natural hermitian inner product  defined by 
\[ \innone{S}{\omega_1, \omega_2} :=\frac i{2} \int_S \omega_1 \wedge \bar{\omega}_2.\]
Choosing an orthonormal basis of $\eta_1, \dots,\eta_g$ of  $\Omega^1(S)$, the canonical measure of $S$ is defined by
\[\mu^\can = \frac{i}{2g}\sum_{j=1}^g \eta_j \wedge \bar{\eta}_j,\]
which is a positive density measure of total mass one on $S$. 

An alternate description of $\mu^\can$ can be given as follows. Let $a_1, \dots, a_g, b_1, \dots, b_g$ be a symplectic basis of $H_1(S, \Z)$. This means for the intersection pairing $\intprod{}{\cdot\,,\cdot}$ between 1-cycles in $S$, we have 
\[\intprod{}{a_i, a_j} = \intprod{}{b_i, b_j} =0, \qquad \textrm{and} \qquad \intprod{}{a_i, b_j} =\delta_{i,j}. \]
Let $\omega_1, \dots, \omega_g$ be an adapted basis of $\Omega^1(S)$ in the sense that 
\[\int_{a_j} \omega_i =\delta_{i,j}\]
and define the period matrix of $S$ by
\[\Omega := \Bigl(\Omega_{i,j}\Bigr) = \Bigl(\int_{b_j}\omega_i\Bigr)_{i,j=1}^g.\]
The imaginary part  $\Im(\Omega)$ of $\Omega$ is a symmetric positive definite real matrix, and we denote by $\Im(\Omega)^{-1}$ its inverse. Then we have 
\begin{equation} \label{eq:cannonical_measure_RS}
\mu^\can = \frac {\sqrt{-1}}{2g} \sum_{i,j=1}^g \Im(\Omega)^{-1}_{i,j} \, \omega_i \wedge \bar{\omega}_j
\end{equation}
with $\Im(\Omega)^{-1}_{i,j} $ denoting the $(i,j)$-coordinate of $\Im(\Omega)^{-1}$.

\subsubsection{Hybrid curves} Consider now a hybrid curve $\curve$ with underlying graph $G=(V, E)$, ordered partition $\pi=(\pi_1, \dots, \pi_r)$ of the edge set $E$, and the edge length function $l$ with the normalization property that the sum of edge lengths in each layer is equal to one. Denote by $\Gamma^j$, $j\in [r]$, the graded minors of $\curve^\trop$, the corresponding tropical curve to $\curve$. For each vertex $v$, let $C_v$ be the corresponding smooth compact Riemann surface of genus $\genusfunction(v)$.

The canonical measure $\lmu^\can$ on $\curve^\hyb$ is by definition the measure which restricts to the canonical measure of the tropical curve $\curve^\trop$ on the intervals $\Ical_e$, for each edge $e\in E$, and which coincides with the rescaled Arakelov--Bergman measure $\frac1g \mu_{\Ar,v}$ on each Riemann surface $C_v$ of positive genus $\genusfunction (v) >0$. The restriction to Riemann surfaces $C_v$ of genus zero is the null measure.

We recall the definition of the canonical measure on a tropical curve.  Denote by $\mu_\Zh^j$ the Zhang measure of $\Gamma^j$.  By an abuse of the notation, we denote by $\mu_\Zh^j$ the corresponding measure on $\curve^\trop$ with support in the intervals $\Ical_e$ for $e\in \pi_j$. We can write 
\[\mu_\Zh^j = \sum_{e\in \pi_j} \frac{\mu^j(e)}{l(e)} d\theta_e\]
for the uniform Lebesgue measure $d\theta_e$ on the interval $\Ical_e$ and the Foster coefficients $\mu^j(e)$, $e \in \pi_j$ of the metric graph $\Gamma^j$. The canonical measure $\lmu^{\can}$ of the tropical curve $\curve^\trop$ is the measure given by
\[\lmu^\can = \frac 1g \Bigl(\,\sum_{v\in V} \genusfunction(v) \delta_v + \sum_{j=1}^r\mu^j_\Zh\,\Bigr).
\]
The canonical measure of the hybrid curve $\curve$ is given by 
\[\lmu^\can := \frac 1g \Bigl(\,\sum_{v\in V} \mu_{\Ar, v} + \sum_{j=1}^r\mu^j_\Zh\,\Bigr).
\]
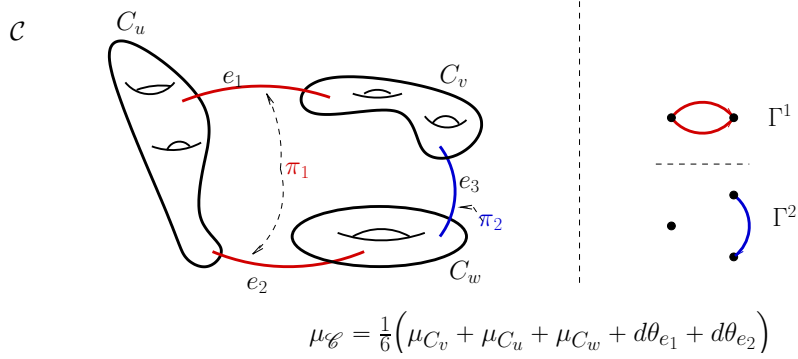
\begin{figure}[!t]
\centering
   \scalebox{.4}{\input{canonical_measures_2.tikz}}
\caption{Example of a hybrid curve and its corresponding canonical measure}
\label{fig:hybrid_canonical_measure}
\end{figure}

The following is straightforward.

\begin{prop} Let $\forget: \curve \to \curve^\trop$ be the natural forgetful projection map which contracts each Riemann surface $C_v$ to the vertex $v$ of $\curve^\trop$. The push-out of the canonical measure on the hybrid curve $\curve$ by $\forget$ coincides with the canonical measure of the tropical curve $\curve^\trop$.
\end{prop}

\subsection{Hybrid Green functions} 
In this section we define the canonical Green function on hybrid curves. 

\subsubsection{Green functions on Riemann surfaces} Consider first the case of a Riemann surface $S$ of genus $g$. Let $\mu$ be a Borel measure on $S$ of mass $\mu(S) = 1$, which for simplicity we assume to be given by a continuous $(1,1)$-form.

The Green function $\grs \colon (S \times S) \setminus \mathrm{diag}_S \to \R$  of $S$ is the function defined by the conditions
\begin{align*}
\frac{1}{\pi i}\partial_z \partial_{\bar{z}}\, \grs(p, \cdot) &= \delta_p - \mu   , \qquad \int_S \grs(p, q) \, d\mu(q) = 0\end{align*}
which hold for all $p \in S$. The first equation has to be interpreted in a distributional sense and means that
\[
	\frac{1}{\pi i} \int_S \grs(p, \cdot) \partial_z \partial_{\bar{z}} f =  - f(p) + \int_S f(q)  \, d \mu(q), \qquad f \in \mathcal{E}^0 (S),
\]
where $\mathcal{E}^0 (S)$ denotes the space of $C^\infty$-functions on $S$.

In the case $\mu = \mu_\Ar$ is the canonical measure on $S$, we call the corresponding Green function the Arakelov Green function of $S$, and denote it by $\grs_\Ar$ or simply $\grs$ if the measure is understood from the context.

\subsubsection{Green functions on hybrid curves}

We follow the notation of the previous section and consider a hybrid curve $\curve$. Assume that around each point $p^e_v$, $v\in V, e\in E, e\sim v$, we have  a holomorphic coordinate $z^e_v$. We consider the canonical layered metrized complex representative $\mccan$ of $\curve$ with normalized edge lengths on each layer. Let $\lmu$ be a layered measure on $\curve$ which is a continuous (1,1)-form on the components $C_v$ of the hybrid curve, $v\in V$.

 For each point $p\in \curve$, consider the Poisson equation for a hybrid function $\lf=(f_1, \dots, f_r, f_\smallc)$
 \begin{equation} \label{eq:green_general_measure_hybrid}
\begin{cases}
 \Deltatrop (\lf) = \delta_{p} - \lmu \\[1mm]
\text{$\lf$ is harmonically arranged} \\[1mm]
 \int_{\curve} \lf d\lmu= 0.
  \end{cases}
\end{equation}
Recall that the last equation means 
\[\int_{\mccan} f_j d \mu=0, \qquad j\in [r]\cup\{\smallcc\}
 \]
 where $f_j$ here refers to the natural extension of $f_j$ to $\mccan$ (by linear interpolation on edges of lower indexed layers). Moreover, for $j=\smallcc$, the extension of $f_\smallc$ refers to that of the regularization $f^\reg_\smallc$, defined after regularization of $f_\smallc$, given by local coordinates $z^e_v$.  The expression $\delta_p -\lmu$ is viewed as a layered measure on $\curve$. By Theorem~\ref{thm:ExUniqHybridPoisson}, the hybrid Poisson equation~\eqref{eq:green_general_measure_hybrid} is well-posed and has a unique harmonically arranged solution.

\begin{defn}[Hybrid canonical Green function] Consider the canonical measure $\lmu=\lmu^\can$ on $\curve$.  The \emph{hybrid canonical Green function} of $\curve$ is the function $\lgri{\curve}\colon \curve \times \curve \to \R^r\times \eRm$ whose restriction $\lgri{\curve}(p, \cdot)$, for $p\in \curve$, is the solution of the above equation~\eqref{eq:green_general_measure_hybrid}.
\end{defn}
Here and below, $\eRm = \R \cup \{- \infty\}$.

The following proposition gathers the main properties of the hybrid Green function.
\begin{prop} \label{prop:BasicPropertiesHybridGreenFunctions} Let $\curve$ be a hybrid curve.
\begin{itemize}
\item For any point $p$, the restriction $\lgri{\curve}(p, \cdot)$ belongs to the extended hybrid Arakelov--Zhang space $D_{\log}(\Deltahybind{\curve})$. Moreover, the only possible singularities of $\lgri{\curve}(p, \cdot)$ are at the points $p^e_v$, $v\in V, e\in E, e\sim v$, or at $p$ when $p$ is a complex point of $\curve$.

\item The tropical part $\lgri{\curve}^\trop$ coincides with the tropical Green function $\lgri{\curve^\trop}$.
\item The  hybrid Green function $\lgri{\curve}$ is symmetric, that is, we have $\lgri{\curve}(p,q) = \lgri{\curve}(q,p)$ for any pair of points $p,q$ in $\curve$.
\end{itemize} 
\end{prop}
Note that, by our convention, the expression $\lgri{\curve}\colon \curve \times \curve \to \R^r\times \eRm$, and the symmetry property stated above, are understood after naturally extending  the components of $\lgri{\curve}(p, \cdot)$ on the canonical metrized complex representative $\mccan$ in the conformal equivalence class of $\curve$. 

\subsection{Weak tameness of hybrid Green functions}

We have the following tameness result on the Arakelov Green function.

\begin{thm}[Weak tameness of hybrid canonical Green functions] \label{thm:WeakTamenessArakelovGreenFunction}
Notations as above, the family of functions $\lgri{\rsf_\thy}(p_\thy, \cdot)$ is a weakly tame family of functions on $\rsf^\hyb/B^\hyb$.
\end{thm}  

In the next section, we complement Theorem~\ref{thm:WeakTamenessArakelovGreenFunction} by a precise description of the asymptotics of the Arakelov Green function (see Theorem~\ref{thm:MainArakelovIntro} and Theorem~\ref{thm:MainArakelovIntro2}).

\smallskip

The proof of Theorem~\ref{thm:WeakTamenessArakelovGreenFunction} is given in Section~\ref{ss:ProofMainGreenFunction}.


\section{Asymptotics of the Arakelov Green function} \label{sec:MainArakelovDetails}
This section contains our main results on the \emph{asymptotics of the Arakelov Green function}.

\smallskip

In Section~\ref{ss:MainArakelovDetails}, we study the behavior of the Green function $\gri{\rsf_t}(p_t, \cdot)$, when the one-marked Riemann surface $(\rsf_t, p_t)$ degenerates to a one-marked hybrid curve $(\curve, p)$. We obtain an asymptotic description in terms of hybrid Green functions (Theorem~\ref{thm:MainArakelovDetails}).

In Section~\ref{ss:MainArakelovDetails}, we proceed to study the Green function $\gri{\rsf_t}(\cdot, \cdot)$, now viewed as a function of two variables, when the Riemann surface $\rsf_t$ degenerates to a hybrid curve $\curve$. In contrast to the first case, we allow the first variable $p$ to be any point on $\rsf_t$, instead of a smooth section $p_t$ staying away from the appearing nodes. We describe the asymptotics in terms of hybrid Green functions and an additional explicit correction term (Theorem~\ref{thm:MainArakelovDetails2}). This result hence answers our original Question~\ref{question:main-surfaces2} in the introduction.

Theorem~\ref{thm:MainArakelovDetails} and Theorem~\ref{thm:MainArakelovDetails2} refine Theorem~\ref{thm:MainArakelovIntro} and Theorem~\ref{thm:MainArakelovIntro2}, which were stated in the introduction.

\subsection{Asymptotics of the Arakelov Green function in $\mghyb{{\combind{g,1}}}$} \label{ss:MainArakelovDetails}

\smallskip

Fix a point $\shy$ in $\mghyb{{\combind{g,1}}}$ representing a hybrid $\curve$ having one marked point. Denote by $(S_0, p_0)$ the underlying stable Riemann surface with one marked point $p_0$ and let $G = (V,E)$ be the underlying stable graph.   Let $\rsf \to B$ be a versal deformation family for $(S_0, p_0)$. 

\smallskip

We \emph{fix an adapted system of coordinates} on the versal family.
 Recall in particular that the coordinates of $B \cong \Delta^N$, $N = 3g-2$ are naturally decomposed into $\underline z = \underline z_E\times \underline z_{E^c}$, where $E^c = [N] \setminus E$. Moreover, the fiber $\rsf_t$ of a point $t \in B^\ast$ decomposes into (see \eqref{eq:AdaptedCoordinates} for details)
\begin{equation} \label{eq:AdaptedDecompositionFinals}
 \rsf_{t}  = \bigsqcup_{v\in V} Y_{v,t} \sqcup \bigsqcup_{e\in E} W_{e, t}.
\end{equation}

The fiber $\rsf_t$ over a point $t \in B^\ast$ is isomorphic to the Riemann surface represented by $t$; moreover, there is a section $p_t \in \rsf_t$, $t \in B$, representing the marking. Due to the choice of the adapted system of coordinates, the point $p_t$, $t \in B^\ast$ always lies in the set $Y_{v,t}$ for the vertex $v$ whose associated component $C_v$ in $S_0$ contains $p_0$.

\smallskip

Denote by $\rsf^\hyb \to B^\hyb$ the corresponding family of hybrid curves. Since the section $p$ lies always in the smooth part of the respective stable Riemann surfaces, it extends naturally to a continuous section $p_\thy$, $\thy \in B^\hyb$, of the hybrid family $\rsf^\hyb \to B^\hyb$. The fixed hybrid point $\shy$ lies in the hybrid stratum $D_\pi^\hyb$, associated to the ordered partition $\pi = (\pi_1, \dots, \pi_r)$ on the edge set $E$ which underlies the hybrid curve $\curve$. It has the form $\shy = (l,0)$ for the complex coordinate $s = 0$, and the simplicial coordinates $l = (l_e)_{e \in E}$ representing the edge lengths of the hybrid curve $\curve$ (satisfying the layerwise normalization conditions). \smallskip

We are interested in \emph{describing the Arakelov Green function} $\gri{t}(p_t, \cdot) \colon \rsf_t \setminus \{p_t\} \to \R$ as the marked Riemann surface $(\rsf_t, p_t)$ degenerates to the marked hybrid curve $(\rsf^\hyb_\shy, p_\shy)$, that is, as $t \in B^\ast$ converges tamely to $\shy$ in $B^\hyb$. For each base point $\thy \in B^\hyb$, we consider the canonical Green function $\lgri{\thy}(p_\thy, \cdot) := \lgri{\rsf_\thy}(p_\thy, \cdot)$. Here and below, $\lgri{\thy}(p_\thy, \cdot) \colon \rsf^\hyb_\thy \setminus \{p_\thy\} \to \R^{r_\thy} \times  \eRm$ is viewed as a function defined on the fiber $\rsf^\hyb_\thy$ (which coincides with the canonical metrized complex $\Sigma_\thy$ of the hybrid curve represented by $\thy$). Recall also that $\eRm := \R \cup \{- \infty\}$.

\medskip

Consider the log map for the stratum $D_\pi^\hyb$, given by the commutative diagram 
\[
\begin{tikzcd}
 \rsf^\ast \arrow[d]\arrow[r, "\loghybdiagpi{\pi}"] & \arrow[d] \rsf_\pi^\hyb\\
B^\ast\arrow[r, "\loghybdiagpi{\pi}"] &D_\pi^\hyb.
\end{tikzcd}
\]
Recall that, on the base $B^\hyb$, the log map sends every point $t \in B^\ast$ representing the marked surface $(\rsf_t, p_t)$ to a point $\thy = \loghyb_\pi(t)$ in $D_\pi^\hyb$, representing the marked hybrid curve $(\rsf^\hyb_\thy, p_\thy)$ of type $(G, \pi)$. The log map on the family is given by an additional map $ \loghyb_\pi \colon \rsf_t \to \rsf^\hyb_\thy$ between the fibers. By the properties of the hybrid log map, the marked point $p_t \in \rsf_t$ is mapped to the marked point $p_\thy \in \rsf^\hyb_\thy$, that is, $ \loghyb_\pi  (p_t) = p_\thy$.

\smallskip

For every $t \in B^\ast$, we push out the canonical measure $\mu_t$ from the Riemann surface $\rsf_t$ to a measure $\widetilde{\mu}_t$ on the hybrid curve $\rsf^\hyb_\thy$, $\thy = \loghyb_\pi(t)$, via the fiberwise map $ \loghyb_\pi \colon \rsf_t \to \rsf^\hyb_\thy$. Let $\widetilde{\lgr}_t(p_t, \cdot)$ be the Green function of the measure $\widetilde{\mu}_t$ on $\rsf^\hyb_\thy$. Denote by $\grihat{t,1}(p_t,\cdot), \dots, \grihat{t,\smallcc}(p_t,\cdot)$ the pullbacks of its components $\gritilde{t,1}(p_\thy,\cdot), \dots, \gritilde{t,\smallcc}(p_\thy,\cdot)$ to $\rsf_t$ (see Section~\ref{sec:pullback_rs}).

\medskip

In the formulation of our result, we use the following \emph{regularization procedure}. Consider a hybrid curve fiber $\rsf^\hyb_\thy$, $\thy \in B^\hyb$, in the family $\rsf^\hyb \to B^\hyb$. Let $f$ be a function in the domain $D_{\log}(\Delta_\smallcc)$ for the disjoint union of components $\pi_\smallcc  = \bigsqcup_u C_u$ of $\rsf^\hyb_\thy$. That is, $f$ is a complex valued $\mathcal{C}^2$-function on $\pi_\smallcc $, except  for finitely many logarithmic poles. As discussed in Section~\ref{ss:PullbackHCtoMC}, we can identify $f$ with a function from the fiber $\rsf^\hyb_\thy$ to $\C$. Around every attachment point $p^e_u = p^e_u(\thy)$ on a component $C_u$ of $\rsf^\hyb_\thy$, we consider the open disc
\[
\Delta_1(p^e_u) := \bigl\{q \in C_u \, \st \,  \text{$q$ is close to $p^e_u$ and } |z^e_u(q)| < 1 \bigr\}
\]
induced by the adapted coordinates, naturally viewed as a subset of $\rsf^\hyb_\thy$.  We denote by
\begin{equation} \label{eq:OpenPunctures}
\inn \Delta_1 (p^e_u) := \bigl  \{\, q \in \Delta_1(p^e_u)\, \st \, 0 < |z^e_u(q) | < 1 \, \bigr  \}
\end{equation}
the corresponding open punctured disc. (More generally, for $1\geq \delta>0$, we denote by $\Delta_\delta (p^e_u)$ the open disc of radius $\delta$ around $p^e_u$, and by $\inn \Delta_\delta (p^e_u)$ the corresponding punctured disc.)

 On  $\inn \Delta_1 (p^e_u)$, the function $f$ is of the form
\[
f(q) = c \log| z^e_u(q)| + h(q) 
\]
where $h$ is a continuous function on $\Delta_1(p^e_u)$. The (\emph{continuous}) \emph{regularization of $f$} is the function  $f^\reg$ on $\rsf^\hyb_\thy$ obtained by modifying $f$ to $f - c \log|z^e_u(q)|$ on each open disc $\Delta_1(p^e_u)$ around each attachment point $p^e_u$. Note that, with this definition, $f^\reg$ is indeed a continuous function on the fiber $\rsf^\hyb_\thy$ (except possibly for finitely many remaining logarithmic poles in the interior of the components $C_u$). Note also that the values of $f$ and $f^\reg$ coincides on the boundary of the unite disc $\Delta_1(p^e_u)$.  

\smallskip 

From this construction, we obtain the regularizations $\gritildereg{t, \smallcc}$, $t \in B^\ast$, and $\grireg{\thy, \smallcc}$, $\thy \in B^\hyb$, of the complex parts in the Green functions. Note that $\gritildereg{t, \smallcc}$ and $\grireg{\thy, \smallcc}$ are continuous functions on $\rsf^\hyb_{\thy'} \setminus\{ p_t\}$, $\thy' = \loghyb_\pi(t)$, and $\rsf^\hyb_\thy \setminus \{ \thy\}$, respectively. Moreover, the pole coefficients in the regularization procedure come from the graph components in the Green function. That is, for an attachment point $p^e_u$ and a point $p$ on a fiber $\rsf^\hyb_\thy$, $\thy \in B^\hyb$, define the coefficients
\begin{align} \label{eq:PoleCoeff1}
c_\thy (p, p^e_u) &:= - \sum_{j=1}^r \divind{j}{\smallcc} \big  (\gri{\thy, j}(p, \cdot) \big) \, (p^e_u).
\end{align}
Moreover, for $t \in B^\ast$, an attachment point $p^e_u$ and a point $p$ on the fiber $\rsf^\hyb_\thy$, $\loghyb_\pi(t) = \thy$, set
\begin{align}  \label{eq:PoleCoeff2}
\~ c_t (p, p^e_u) &:= - \sum_{j=1}^r \divind{j}{\smallcc} \big  (\gritilde{t, j}(p, \cdot) \big) \, (p^e_u),
\end{align}
Here $\divind{j}{\smallcc}(\cdot)$ denotes the divisors from~\eqref{eq:DefineDivisorsHybridLaplacian}. Then the coefficients used in the regularization are precisely $c_\thy(p_\thy, p^e_u)$ for $\grireg{\thy, \smallcc}$ and $\~ c_t (p_\thy, p^e_u)$, $\thy = \loghyb_\pi(t)$, for $\gritildereg{t, \smallcc}$.

\medskip

We obtain the following asymptotic expansion of the Arakelov Green function in terms of hybrid Green functions.
Note that this result gives a precise version of Theorem~\ref{thm:MainArakelovIntro} in the introduction.

\begin{thm} \label{thm:MainArakelovDetails} Suppose that $t \in B^\ast$ converges tamely to $\shy$ in $B^\hyb$. Then:
\begin{enumerate}\item [(i)]The following logarithmic expansion holds for $y \in \rsf_t$, 
\begin{equation} \label{eq:MainArakelovDetails(i)}
\gri{t}(p_t, y) = \sum_{j =1}^r  L_j(t) \grihat{t,j}(p_t,y)  + \grihat{t, \smallcc}(p_t, y) + o(1),
\end{equation}
where $L_j(t) = \sum_{e \in \pi_j} - \log|z_e(t)|$ and the $o(1)$ term goes to zero uniformly for $y \in \rsf_t$. 

\smallskip

\item  [(ii)]  The pullbacks $\grihat{t,j} (p_t, \cdot) \colon \rsf_t \to \eRm$ of $\gritilde{t, j} (p_\thy, \cdot) \colon \rsf^\hyb_\thy \to \eRm$ converge in the following way.

If the point $t \in \loghyb_\pi^{-1}(\shy)$ converges tamely to $\shy$ in $B^\hyb$, then for every $j \in [r]$,
\begin{equation}  \label{eq:MainArakelovDetails(iia)}
\gritilde{t, j}(p_\shy, y) = \gri{\shy, j}(p_\shy, y) + o(1), \qquad y \in \rsf^\hyb_\shy,
\end{equation}
and the regularizations of the complex parts satisfy
\begin{equation} \label{eq:MainArakelovDetails(iib)}
\gritildereg{t, \smallcc} (p_\shy, y) = \grireg{\shy, \smallcc} (p_\shy, y) + o(1),  \qquad y \in \rsf^\hyb_\shy,
\end{equation}
where the $o(1)$ terms go to zero uniformly for $y \in \rsf^\hyb_\shy$.

Moreover, the coefficients of the poles converge, that is, $\~c_t(p_\shy, p^e_u) \to c_\shy (p_\shy, p^e_u)$ for all attachment points $p^e_u$.
\end{enumerate}
\end{thm}

In order to define the $o(1)$ notation for functions with values in $ \eRm = \R \cup \{-\infty\}$, we employ the convention that $a + \infty :=  \infty$ for all $a \in \R$, and $\infty - \infty := 0$. In other words, both sides in \eqref{eq:MainArakelovDetails(i)} take the value $- \infty$ simultanuously and the difference goes to zero uniformly everywhere else; and similar for the functions in \eqref{eq:MainArakelovDetails(iib)}.

\smallskip

The proof of Theorem~\ref{thm:MainArakelovDetails} is given in Section~\ref{ss:ProofMainArakelovDetails}.

\begin{remark}  We briefly explain how the above translates into the relationship between the functions $\grihat{t, j}(p_t, y)$ and the pullbacks $\mathrm{g}^*_{_{\thy, \smallcc}}(p_t, y)$, as stated in Theorem~\ref{thm:MainArakelovIntro}  (2). If $t \in \loghyb_\pi^{-1}(\shy)$ converges tamely to $\shy$ in $B^\hyb$, then by Theorem~\ref{thm:MainArakelovDetails}(iii),
\[
\grihat{t, \smallcc}(p_t, y) - \mathrm{g}^*_{_{\shy, \smallcc}}(p_\shy, y) = o(1) , \qquad y\in  \rsf_t,
\]
where $o(1)$ term goes to zero uniformly for $y$ in the set $\rsf_t \setminus \bigsqcup_{u, e} A^e_{u,t}$ (here, the subset $A^e_{u,t} \subseteq \rsf_t$ for $u \in V$, $e \sim u$ is defined in \eqref{eq:DefAeu}). This formalizes the assumption made in Theorem~\ref{thm:MainArakelovIntro} that \emph{$y \in \rsf_t$ remains uniformly separated from the appearing nodes}. Moreover, on each set $A^e_{u,t}$, we have
\[
\grihat{t, \smallcc}(p_t, y) - \mathrm{g}^*_{_{\shy, \smallcc}}(p_\shy, y) = \Big( \~ c_t (p_\thy, p^e_u) - c(p_\shy, p^e_u) \Big) \log|z^e_u(y)|  + o(1), \qquad y \in A^e_{u,t},
\]
where the $o(1)$ term goes to zero uniformly in $y$.
\end{remark}

\subsection{Asymptotics of the Arakelov Green function in $\mghyb{g}$} \label{ss:MainArakelovDetails2}
We proceed to give an asymptotic description of the Arakelov Green function in the hybrid moduli space $\mghyb{g}$. We follow the discussion in the previous section.

\smallskip

Let $\curve$ be a hybrid curve with underlying stable Riemann surfaces $S_0$, normalized edge length function $l \colon E \to (0, + \infty)$ on the edge set $E$ of the dual graph $G = (V,E)$, and ordered partition $\pi = (\pi_1, \dots, \pi_r) \in \Pifs(E)$. Consider the hybrid family $\rsf^\hyb \to B^\hyb$ associated to a versal family $\rsf \to B$ for $S_0$, equipped with adapted coordinates. The hybrid curve $\curve$ is represented by a point $\shy \in B^\hyb$, which lies in the hybrid stratum $D_\pi^\hyb$ and has the form $\shy = (l, 0)$.

\smallskip

For a base point $\thy \in B^\hyb$,  let $\lgr_\thy$ be the hybrid canonical Green function on $\rsf^\hyb_\thy$. Recall that we view the latter as a function $\lgr_\thy \colon \rsf^\hyb_\thy \times  \rsf^\hyb_\thy \to \R^r \times  \eRm$, where $\rsf^\hyb_\thy$ is the fiber over $\thy$ in the family (which coincides with the canonical metrized complex of the hybrid curve represented by $\thy$). 

\smallskip

Consider again the hybrid log map $\loghyb_\pi$ associated to the stratum $D_\pi^\hyb$ (see the previous section).
For $t \in B^\ast$, we pushout the canonical measure $\mu_t$ from the smooth Riemann surface $\rsf_t$ to a measure $\~{\mu}_t$ on  $\rsf^\hyb_\thy$, $\thy = \loghyb_\pi(t)$. Denote by $\~{\lgr}_t \colon \rsf^\hyb_\thy\times \rsf^\hyb_\thy \to \R^r \times \eRm $ the associated Green function on  $\rsf^\hyb_\thy$. We pullback the components of $\~{\lgr}_t$ to the smooth Riemann surface $\rsf_t$ via the log maps on the family $\rsf^\hyb \to B^\hyb$. That is, we define 
\[
\grihat{t,j}(p,y) := \gritilde{t,j} \Big ( \loghyb_\pi(p), \loghyb_\pi(y) \Big), \qquad p,y \in \rsf_t,
\]
for $j \in \{1, \dots, r, \smallcc\}$ and obtain the pullbacks $\grihat{t,j} \colon \rsf_t \times \rsf_t \to \R$, $j \in [r]$ and $\grihat{t,\smallcc} \colon \rsf_t \times \rsf_t \to  \eRm$.

\smallskip

In order to state our result, we define the following \emph{regularizations of the complex parts of the Green functions}. Let $c_\thy (p, p^e_u)$ and $\~ c_t (p, p^e_u)$ be the coefficients introduced in \eqref{eq:PoleCoeff1} and \eqref{eq:PoleCoeff2}, respectively. Recall that $c_\thy (p, p^e_u)$ is precisely the coefficient of the pole of $\gri{\thy, \smallcc} (p, \cdot)$ at the node $p^e_u$, and similar for $\~ c_t (p, p^e_u)$.

\smallskip

For $\thy \in \partial_\infty B^\hyb$, we define the regularization $\grireg{\thy, \smallcc} \colon \rsf^\hyb_\thy \times \rsf^\hyb_\thy \to \eRm$ as
\[
\grireg{\thy, \smallcc}(p,y) := \gri{\thy, \smallcc}(p,y) - \sum_{p^e_u} \bm{1}_{\inn\Delta_1{(p^e_u)}}(p) \, c_\thy (y, p^e_u) \log|z^e_u(p)| + \bm{1}_{\inn\Delta_1{(p^e_u)}} (y) \, c_\thy (p, p^e_u) \log|z^e_u(y)|.
\]
Here the sum is taken over all attachment points $p^e_u$, $e = uv$ on $\rsf^\hyb_\thy$, and $z^e_u$ is the adapted coordinate on the component $C_u$ of $\rsf^\hyb_\thy$, and $\bm{1}_{\inn \Delta_1{(p^e_u)}}$ is the indicator function of the punctured disc $\inn\Delta_1{(p^e_u)} = \{ 0 < |z^e_u| < 1\}$ around $p^e_u$ in $C_u$ (see \eqref{eq:OpenPunctures}).

\smallskip

Analogously, for $t \in B^\ast$, the regularization $\gritildereg{t, \smallcc} \colon \rsf^\hyb_\thy \times \rsf^\hyb_\thy \to \eRm$, $\thy = \loghyb_\pi(t)$, is given by 
\[
\gritildereg{t, \smallcc}(p,y) := \gritilde{t, \smallcc}(p,y) - \sum_{p^e_u} \bm{1}_{\inn \Delta_1{(p^e_u)}}(p) \, \~{c}_t (y, p^e_u) \log|z^e_u(p)| + \bm{1}_{\inn \Delta_1{(p^e_u)}} (y) \, \~{c}_t (p, p^e_u) \log|z^e_u(y)|.
\]
The above corresponds to removing the poles of the complex parts of Green functions at the points $p^e_u$. However, note that we still have $\gritildereg{t, \smallcc}(p,p) = - \infty$ and $\grireg{\shy, \smallcc}(p,p) = - \infty$ for points $p$ in the interior of the Riemann surface components of hybrid fibers $\rsf^\hyb_\thy$.

\smallskip

We obtain the following asymptotic description of the Arakelov Green $\gri{t}(\cdot, \cdot)$, $t \in B^\ast$, as the smooth Riemann surface $\rsf_t$ degenerates to the hybrid curve $\curve = \rsf^\hyb_\shy$. Note that, in contrast to Theorem~\ref{thm:MainArakelovDetails}, the asymptotics contains an additional correction term $\varepsilon_\pi(p,y)$. Moreover, the result gives a precise version of Theorem~\ref{thm:MainArakelovIntro2} in the introduction.

\begin{thm} \label{thm:MainArakelovDetails2} Suppose that $t \in B^\ast$ converges tamely to $\shy$ in $B^\hyb$. Then,
\begin{itemize} 
\item [(i)] The following expansion holds for $p, y \in \rsf_t$, 
\begin{equation}  \label{eq:MainArakelovDetails2(i)}
\gri{t}(p, y) =  \sum_{j =1}^r  L_j(t) \grihat{t,j}(p,y)  + \grihat{t, \smallcc}(p, y)  - \varepsilon_\pi(p,y) + o(1),
\end{equation}
where the $o(1)$ term goes to zero uniformly for $p, y \in \rsf_t$. Here $L_j(t) = \sum_{e \in \pi_j} - \log|z_e(t)|$, $j \in [r]$, and $\varepsilon_\pi(p,y)$ is a correction term defined below.

\smallskip

\item  [(ii)] The pullbacks $\grihat{t,j} \colon \rsf_t \times \rsf_t \to \eRm$ of $\gritilde{t, j} \colon \rsf^\hyb_\thy \times  \rsf^\hyb_\thy\to  \eRm$ converge in the following way.

If $t \in \loghyb_\pi^{-1}(\shy)$ converges tamely to $\shy$ in $B^\hyb$, then for every $j \in [r]$,
\begin{equation}  \label{eq:MainArakelovDetails2(iia)}
\gritilde{t, j}(p, y) = \gri{\shy, j}(p, y) + o(1), \qquad p,y \in \rsf^\hyb_\shy,
\end{equation}
and the regularizations of the complex parts satisfy
\begin{equation}  \label{eq:MainArakelovDetails2(iib)}
\gritildereg{t, \smallcc} (p, y) = \grireg{\shy, \smallcc} (p, y) + o(1), \qquad p,y \in \rsf^\hyb_\shy,
\end{equation}
where the $o(1)$ terms go to zero uniformly for $p,y \in \rsf^\hyb_\shy$.

Moreover, the coefficients of the poles converge, that is, $\~c_t(p, p^e_u) \to c_\shy (p, p^e_u)$ for all $p \in  \rsf^\hyb_\shy$ and all attachment points $p^e_u$.
\end{itemize}
\end{thm}
As before, in order to define the $o(1)$ notation for functions with values in $ \eRm = \R \cup \{-\infty\}$, we employ the convention that $a + \infty :=  \infty$ for all $a \in \R$, and $ \infty - \infty := 0$. In other words, both sides in \eqref{eq:MainArakelovDetails2(i)} take the value $- \infty$ simultanuously and the difference goes to zero uniformly everywhere else; and similar for the functions in \eqref{eq:MainArakelovDetails2(iib)}.

\smallskip
The proof of Theorem~\ref{thm:MainArakelovDetails2} is given in Section~\ref{ss:ProofMainArakelovDetails2}.

\begin{remark}[The correction term $\varepsilon_\pi$] Before defining the correction term $\varepsilon_\pi$, we briefly explain its meaning. Recall that the Green function on $\rsf_t$ has a logarithmic pole $\gri{t}(p, y) \approx - \log|p-y|$ if, $p$ and $y$ become close to each other on the fiber $\rsf_t$. This holds in particular, if $p$ and $y$ are close and lie in the middle of some cylinder $W_{e,t}$ of $\rsf_t$ (for the decomposition of $\rsf_t$ with respect to the adapted coordinates \eqref{eq:AdaptedDecompositionFinals}). However, the hybrid log maps contract the complex cylinder $W_{e,t}$ to an interval in the metric graph, losing the information on the appearance of this pole. The correction term $\varepsilon_\pi(p,y)$ computes the logarithmic pole (in suitable coordinates, provided by the adapted system) and adds it back. \end{remark}

We proceed to introduce the above \emph{correction term}. For two points $p \neq x$ on a smooth fiber $\rsf_t$, $t \in B^\ast$, and an ordered partition $\pi$ of full sedentarity on a subset $E_\pi \subseteq E$, define
\begin{equation} \label{eq:CorrectionTwoPoints}
\varepsilon_\pi (p,x) := \sum_{e \in E_\pi} \varepsilon_e(p,x),
\end{equation}
where for each edge $e \in E$, the contribution $\varepsilon_e(p,x)$ is given in adapted coordinates as follows. We use the notations from Section~\ref{sec:hybrid_log_map_einf} and consider the auxiliary function
  \[
  \varrho_e(t) := - \frac{1}{\log|z_e(t)|}, \qquad t \in B^\ast.
 \]
The cylinder $W_{e,t} \subseteq \rsf_t$ associated to $e \in E$  decomposes into $W_{e,t} = B_{e, t} \cup A^e_{u,t} \cup A^e_{v,t}$, where 
\begin{align*}
A^e_{u,t}  &= \Big \{ (\underline z, z^e_u, z^e_v) \in W_{e,t} \, \st \, |z^e_u| > \varrho_e(t)  \Big \}, &A^e_{v,t}  = \Big \{ (\underline z, z^e_u, z^e_v) \in W_{e,t} \big | \, |z^e_v| > \varrho_e(t)  \Big \}, \\
B_{e,t} &=\Big \{ (\underline z, z^e_u, z^e_v) \in W_{e,t} \, \st | \, |z^e_u|, |z^e_v| < \varrho_e(t)  \Big  \}.
\end{align*}
The following six cases, vaguely speaking, formalize the idea that "$p$ and $x$ are two points in the middle of the cylinder $W_{e,t}$ and close to each other":
\begin{itemize}
\setlength\itemsep{0.2 mm}
\item[(i)] $p,x \in B_{e,t}$, and $|z^e_u(p)|  > |z^e_u(x)|$ and  $|z^e_u(x)| / |z^e_u(p)| > \varrho_e(t)$,
\item[(ii)] $p,x \in B_{e,t}$ and $|z^e_v(p)|  > |z^e_v(x)|$ and  $|z^e_v(x)| / |z^e_v(p)| > \varrho_e(t)$,
\item[(iii)] $p \in A^e_{u,t}$, $x \in B_{e,t}$, and $\varrho_e(t)^{1/2}  > |z^e_u(p)|$ and $ |z^e_u(x)| > \varrho_e(t)^{3/2}$
\item[(iv)] $ x\in A^e_{u,t}$, $p \in B_{e,t}$, and $\varrho_e(t)^{1/2}  > |z^e_u(x)|$ and  $ |z^e_u(p)| > \varrho_e(t)^{3/2}$
\item[(v)]  $x \in A^e_{v,t}$, $p \in B_{e,t}$, and $\varrho_e(t)^{1/2}  > |z^e_v(x)|$ and $|z^e_v(p)| > \varrho_e(t)^{3/2}$
\item[(vi)] $p \in A^e_{v,t}$, $x \in B_{e,t}$, and $\varrho_e(t)^{1/2}  > |z^e_v(p)|$ and $|z^e_v(x)| > \varrho_e(t)^{3/2}$.
\end{itemize}
The contribution $\varepsilon_e(p,x)$ of the edge $e \in E$ is then defined as
\begin{align*}
\varepsilon_e(p,x) := \begin{cases}
\log  \Big | \frac{z^e_u(x)}{z^e_u(p)} - 1 \Big  |, &\text{if $p,x \in W_{e,t}$ and (i), (iii) or (iv) holds,} \\[3mm]
\log  \Big | \frac{z^e_v(p)}{z^e_v(x)} - 1 \Big  |, &\text{if $p,x \in W_{e,t}$ and (ii), (iv) or (vi) holds,} \\[3mm]
0, & \text{else}.
\end{cases}
\end{align*}


\section{Hybrid height pairing and hybrid $\jvide$-function} \label{sec:JfunctionHeightPairing}
In this section, we study the limit of the height pairing on degenerating Riemann surfaces. We introduce a \emph{hybrid height pairing} and a notion of \emph{hybrid $\jvide$-function} on hybrid curves. We furthermore formulate our results on the degeneration of the height pairing and $\jvide$-function to their hybrid counterparts. The proof of these results will be given in the forthcoming sections.
 
 \smallskip
 
 We begin by recalling the definition of the height pairing on smooth Riemann surfaces, together with its relationship to the Poisson equation and the $\jvide$-function (see Section~\ref{ss:HPBasics}). Section~\ref{ss:HeightPairingHybridCurves} and Section~\ref{ss:HybridJFunction} contain the crucial definition of the height pairing and $j$-function on hybrid curves. The latter is based on a regularization procedure described in Section~\ref{ss:RegularizedHeightPairing}. In Section~\ref{ss:ResultsHybridHeightPairing}, we state our results on the degeneration of the height pairing and $\jvide$-function on Riemann surfaces. Furthermore, in preparation for the proofs, we recall an expression of the height pairing in terms of mixed Hodge structures called biextensions~\cite{Hain}, see Section~\ref{ss:LogFormOmega}--~\ref{ss:RegularizedIntegralsPaths}.   
 
 \subsection{Height pairing and $\jvide$-function on Riemann surfaces} \label{ss:HPBasics}
Let $S$ be a smooth compact Riemann surface of genus $g$. For a divisor $D$ on $S$, we denote by $\supp{D}$ the support of $D$. 

\smallskip
 
The residue $\Res_p(\alpha)$ of a one-form $\alpha$ around a point $p$ on $S$ is defined as the integral $\int_{\gamma}\alpha$ for a small positively oriented cycle $\gamma$ around $p$. The residue divisor $\Res(\alpha)$ is defined by $ \Res(\alpha):= \sum_{p\in S} \Res_p(\alpha) \,p$.

\smallskip

 For any degree zero divisor $D$ (with real coefficients), there exists a unique logarithmic one-form $\alpha\ind{D}$, a one-form of the third kind, with residue divisor $\Res(\alpha\ind{D}) = D$ and with periods all in $\R$, see e.g.~\cite{FK92, GH14}. By this, we mean that for any cycle $\gamma$ on $S$ which does not contain any point in the support of $D$, the integral $\int_\gamma \alpha\ind{D}$ is real. 
 
 \smallskip

For a degree zero divisor $D$, there exists as well a real one-chain $\gamma\ind{D}$ with $\partial(\gamma\ind{D}) = D$. Note that $\gamma\ind{D}$ is not unique: two such choices $\gamma\ind{D}$ and $\gamma'\ind{D}$ differ by an element of $H_1(S, \R)$.

\smallskip

Consider two degree zero divisors $D_1$ and $D_2$ on $S$ with disjoint support, $\supp{D_1}\cap \supp{D_2} =\varnothing$. 
The \emph{height pairing} $\hp{S}{D_1, D_2}$ between $D_1$ and $D_2$  is defined by
\[\hp{S}{D_1, D_2} : = 2  \pi \Im\bigl(\int_{\gamma\indbis{D_2}}\alpha\indbis{D_1}\bigr).\] 
In the above expression, we assume $\gamma\indbis{D_2}$ does not pass through any point of $\supp{D_1}$ so that the integral takes a finite value. By the choice of the form $\alpha\indbis{D_1}$, the height pairing is well-defined and does not depend on the choice of $\gamma\indbis{D_2}$: for two choices $\gamma\indbis{D_2}$ and $\gamma'\indbis{D_2}$, the difference between  $\int_{\gamma\indbis{D_2}}\alpha\indbis{D_1}$ and $\int_{\gamma'\indbis{D_2}}\alpha\indbis{D_1}$ is a period of $\alpha\indbis{D_1}$, which is real.

\medskip

Alternatively, the height pairing can be described in the following way. For a degree zero divisor $D$, there exists a function $\jfunc{D}: S \setminus \supp D \to \R$ such that 
\begin{align}\label{eq:hd}
\frac{1}{\pi i}\partial_z \partial_{\bar{z}}\, \jfunc{D} &= \delta\ind{D},\end{align}
where, as before, $\delta\ind{D}$ is the measure $\delta_D = \sum_{p\in S} D(p)\delta_p$. Again, the equation is interpreted in a distributional sense
\[
	\frac{1}{\pi i} \int_S \jfunc{D}\partial_z \partial_{\bar{z}} f =  - f(D) = -\sum_{p\in S} D(p)f(p), \qquad f \in \mathcal{E}^0 (S).
\]
The function $\jfunc D$ is unique up to addition by a constant. The height pairing between two degree zero divisors with disjoint support $D_1$ and $D_2$ can be written as 
\begin{equation} \label{eq:CHPvsJFunc}
\hp{S}{D_1, D_2} = \jfunci{D_1}(D_2) = - \frac{1}{\pi i} \int_S \jfunci{D_1}\partial_z \partial_{\bar{z}} \jfunci{D_2}.
\end{equation}
In particular, we get from this description the symmetry $\hp{S}{D_1, D_2} = \hp{S}{D_2, D_1}$.  

\smallskip

To make the link to the first definition, note that the function $\jfunc{D}$ can be obtained as follows, see e.g.~\cite{Lang}. Fix a point $p \in S\setminus \supp D$. For any point $q \in S\setminus \supp D$, choose a path $\gamma$ in $S$ from $p$ to $q$ disjoint from the support of $D$, and set 
\[ \jfunc{D}(q):= 2 \pi \Im\bigl(\int_{\gamma}\alpha\ind{D}\bigr).\]
This function verifies Equation~\eqref{eq:hd}.

\smallskip

The following special case will be of particular interest to us. Let $x, p, q$ be three distinct points on $S$. Consider the degree zero divisor $D= p - q$ on $S$.  By the above discussion, there exists a unique smooth function $\jfunc{p\tiret q}: S \setminus \supp D \to \R$ such that 
\begin{align*}\label{eq:hd-pq} 
\Delta \, \jfunc{p \tiret q,x} &= \delta_p -\delta_q, &
\jfunc{p\tiret q, x}(x)=0.\end{align*}

We call the function $\jfunc{p\tiret q,x}\colon S \setminus \supp D \to \R$ the $\jvide$-function of $S$. As discussed above, its values coincide with the height pairing
\begin{equation} \label{eq:HPJFuncSmooth}
 \jfunc{p\tiret q, x}(y) = \hp{S}{p-q, y-x}.
 \end{equation}
Viewed as a function on $S$,  the $\jvide$-function $y \mapsto \jfunc{p\tiret q, x}(y)$ belongs to $D_{\log}(\Delta)$.

\subsection{Regularized height pairing on a Riemann surface} \label{ss:RegularizedHeightPairing}
Let $S$ be a smooth, compact Riemann surface. Consider two divisors $D_1$ and $D_2$ with supports $\suppdiv_1, \suppdiv_2 \subset S$, respectively. We do not assume that $D_1$ and $D_2$ have disjoint support. Fixing a system of local parameters $z_x$, $x\in \suppdiv_1\cup \suppdiv_2$, we define the following regularized version of height pairing between $D_1$ and $D_2$.
\smallskip

We perturb slightly the support $\suppdiv_{1}$ of $D_1$ by replacing each point $x$ of $\suppdiv_1 \cap \suppdiv_2$ with the point $x+\zeta_x$ for small $\zeta_x$ in the local coordinate $z_x$. Denote by $D'_1$ the divisor obtained by this modification,  
\[D'_{1} = \sum_{x \in \suppdiv_1 \setminus \suppdiv_2} D_1(x) x + \sum_{x \in \suppdiv_1 \cap \suppdiv_2} D_1(x) (x + \zeta_x)\]
and define the \emph{regularized height pairing}
\[\hp{S}{D_1, D_2}' : = \lim_{\substack{\zeta_x \to 0 \\ \forall x\in \suppdiv_1\cap \suppdiv_2}} \hp{S}{D'_{1}, D_2} + \sum_{x\in \suppdiv_1 \cap \suppdiv_2 } D_1(x) D_2(x) \log(|\zeta_x|).  \]

\begin{prop} The limit in the definition of the regularized height pairing exists.
\end{prop} 
\begin{proof} This follows, for instance, from the description of the height pairing using the Green function $\gri{S}$. Since $\partial_z \partial_{\bar{z}} \sum_x D_1(x) \gri{S}(x,z) = \pi i \delta\indbis{D_1}$, we have
\[\hp{S}{D_1, D_2} = \sum_{x,y}D_1(x)D_2(y)\gri{S}(x, y).\]
The regularized height pairing can be described as 
\[ \hp{S}{D_1, D_2}' = \sum_{x\neq y}D_1(x)D_2(y)\gri{S}(x, y) + \sum_{x\in \suppdiv_1 \cap \suppdiv_2} \grprimei{S}(x,x)\]
where $\grprimei{S}(x,x) = \lim_{\zeta_x \to 0} \gri{S}(x, x+\zeta_x) + \log(|\zeta_x|)$. Note that the latter limit exists by the description of the Green function on a smooth Riemann surface, see e.g.~\cite{Lang}. 
\end{proof}
\begin{remark}
We stress once again that the regularized height pairing depends on the choice of local coordinates $z_x$ around the points $x \in \suppdiv_1 \cap \suppdiv_2$.
\end{remark}

For later use, we also derive the following explicit expression for the regularized height pairing on the Riemann sphere. Consider the Riemann sphere $\C\P^1$ with two distinguished points $0=[0:1]$ and $\infty =[1:0]$ in projective coordinates. We choose a global coordinate $z$ on $\C\P^1$ inducing coordinates around any point of $\C\P^1$. 
  
  \begin{prop}\label{prop:explicit_regularization_sphere} Let $D_1$, $D_2$ be divisors of degree zero on $\C\P^1$. Then the regularized height pairing (using the global coordinate $z$) is explicitly given by
  \[\hp{\C\P^1}{D_1, D_2}'  = -\sum_{\substack{x,y\in \C \\ x\neq y}} D_1(x) D_2(y) \log|x-y|.\]
  \end{prop}
  \begin{proof} 
 Recall that, by the  Poincar\'e--Lelong formula, $\frac{i}{\pi} \partial_z \partial_{\bar{z}} \log |f| = \div(f) $ for every meromorphic function $f$ on a Riemann surface $S$. Hence given a degree zero divisor $D$ on $S$,
 \[\hp{}{D, \div(f)} = -\log|f(D)|\]
 if the supports $|D|$ and $|\div(f)|$ are disjoint. Since on the Riemann sphere $\C\P^1$ all degree zero divisors  are principal, the result follows by regularizing in the uniform coordinate $z$.
  \end{proof}

\subsection{Height pairing on hybrid curves} \label{ss:HeightPairingHybridCurves} We introduce a hybrid notion of \emph{height pairing on hybrid curves}. Suppose $\curve$ is a hybrid curve of rank $r$ with underlying stable Riemann surface $S_0$, combinatorial graph $G=(V, E)$, ordered partition $\pi=(\pi_1, \dots, \pi_r)$ of  $E$, and (normalized) edge length functions $l = (l_1, \dots, l_r)$. Let $C_v$, $v\in V$, be the Riemann surface components and $\curve^\trop$ the underlying tropical curve of $\curve$. The associated metrized complex is denoted by $\mccan$.

\smallskip

The \emph{hybrid height pairing} on $\curve$ mixes the tropical height pairings on its graded minors $\Gamma^j$ and the regularized height pairing on its complex part $\pi_\smallcc$. The latter requires to \emph{fix a system of local parameters} $z^e_v$ around the points $p^e_v$, $e \sim v$ on each component $C_v$, $v\in V$.

\medskip

Let $ D \in \Div^0(\curve)$ be a degree zero divisor on $\curve$. Hereby we mean that
\[
D = \sum_{x \in \mccan} D(x) (x)
\]
is a formal combination of finitely many points $x$ on the metrized complex $\mccan$ with real coefficients $D(x)$, $x \in \mccan$, such that $\sum_x D(x) = 0$.

\smallskip

Clearly, $D$ induces a layered measure $\lmu\ind{D}$ of total mass zero on $\curve$ that, by an abuse of the notation, we sometimes denote by $D$ as well. The $j$-th piece $\mu_j$ of the layered measure $\lmu\ind{D}$ equals
\begin{align*}
&\mu_j=  \sum_{v \in V^j} \Big ( \sum_{u \in \kappa_j^{-1} (v)} D(u) \Big ) \,  \delta_v  + \sum_{x \in \Gamma^j\setminus V^j} D(x), &&j \in [r],
\end{align*}
where $V^j=V(\grm{\pi}{j}(G))$ and $\kappa_j \colon V \to V(\grm{\pi}{j}(G))$ is the projection map on vertices. In particular, $\mu_j$ is a point measure (divisor) on the $j$-th graded minor $\Gamma^j$.

\smallskip

The last piece $\mu_\smallcc$ of $\lmu\ind{D}$ is a point measure supported on $\pi_\smallcc = \bigsqcup_v C_v$.  It is of the form
\begin{align*}
&\mu_\smallcc= \sum_{v \in V} D_v = \sum_{v \in V} \Big ( \sum_{x \in C_v} D_v(x) \delta_x \Big ),
\end{align*}
where for each vertex $v \in V$, the divisor $D_v$ on $C_v$ is given by $D_v = \sum_{x\in C_v} D(x) x$.

\smallskip

Since $D$ is a measure of total mass zero on $\curve$, the hybrid Poisson equation
\begin{equation} \label{eq:HybridPairingEquation}
	\Deltahyb (\lf) = D
\end{equation}
has a solution $\lf = (f_1, \dots, f_r, f_\smallcc)$  on the hybrid curve $\curve$ (see Lemma~\ref{lem:ExistenceHybridSolutions}). 

\smallskip

The solution $\lf$ defines a collection of degree zero divisors $D^j$ on the graded minors $\Gamma^j$, $j \in [r]$, of $\curve$, given by
\begin{equation}
D^j = \Delta_{\Gamma^j} (f_j) = \mu_j - \sum_{k < j} \divind{k}{j}(f_k)
\end{equation}
where again we identify divisors with point measures.  
Furthermore, we obtain a degree zero divisor $D^\smallcc_v$ on each component $C_v$, $v\in V$, given by
\begin{equation} \label{eq:ComplexDivisorsHybridCurve}
D^\smallcc_v := \Delta_{\smallcc} (f_\smallcc)\rest{C_v} = D_v - \sum_{j=1}^r \divind{j}{\smallcc}(f_j) \rest{C_v}, \qquad \qquad D^\smallcc :=  \sum_{v \in V} D^\smallcc_v  = \Delta_{\smallcc} (f_\smallcc).
\end{equation}
The divisors $D^j$ and $D^\smallcc_v$ depend only on $D$ and not on the choice of the solution $\lf$ to \eqref{eq:HybridPairingEquation}.

\begin{defn}[Height pairing on a hybrid curve] \label{def:HybridHeightPairing} \rm 
Notations as above, the {\em height pairing} on the hybrid curve $\curve$ of rank $r$ is the $\R^{r+1}$-valued bilinear map 
\[
\begin{array}{cccc}
 \lhp{\curve}{\cdot\,, \cdot} \colon &\Div^0(\curve) \times \Div^0(\curve)  &\longrightarrow & \R^{r+1} = \R^{[r]}\times \R^\smallcc
\end{array}
\]
whose $r+1$ components are defined by
\[
	 \lhp{\curve, j}{D, D} =  \hp{\Gamma^j}{D^j, D^j}, \qquad j\in [r],
\]
and 
\[
 \lhp{\curve, \smallcc}{D, D} =  \sum_{v\in V}\hp{C_v}{D^\smallcc_v, D^\smallcc_v}'.
\]
\end{defn}
Note that the height pairing extends naturally to a bilinear pairing between a pair of degree zero divisors.

\begin{remark} \label{rem:TropicalVSHybridHeightPairing}
The height pairing on a hybrid curve $\curve$ of rank $r$  refines the height pairing on its underlying tropical curve $\curve^\trop$ (see Section~\ref{ss:HeightPairingTropicalCurves}). Any degree zero $D \in \Div^0(\curve)$ induces a degree zero divisor $D^\trop$ on $\curve^\trop$ by means of the contraction map $\kappa \colon \mccan \to \Gamma$. The first $r$ components of the hybrid height pairing $\lhp{\curve}{D \,, D}$ form the tropical height pairing
\[
\lhp{\curve^\trop}{D^\trop \,, D^\trop} \in \R^r,
\]
since the tropical part $\lf^\trop$ of  $\lf$ in \eqref{eq:HybridPairingEquation} solves the equation $\Deltatrop(\lf^\trop) = D^\trop$ on $\curve^\trop$.
\end{remark}

\subsection{The hybrid $\jvide$-function} \label{ss:HybridJFunction} 

\medskip

We extend the definition of the $\jvide$-function to hybrid curves as follows. Let $\curve$ be a hybrid curve of rank $r$ and take three points $x, p,q$ living on the associated metrized complex $\mccan$. Then, assuming given local parameters around all the points $p^e_v$, $v\in V, e\in E$, with $e\sim v$, as well as local parameters around $p$ and  $q$, we define the function $\ljfuncbis{p\tiret q, x} \colon \mccan \to \R^{r+1}$ by 
\[
\ljfuncbis{p\tiret q, x}(y) := \lhp{\curve}{p-q, y-x}'.
\]

As in the classical setting of Riemann surfaces, we obtain a solution to a Poisson equation.

\begin{thm}[Hybrid height pairing and $j$-function] \label{thm:HybridHPvsJFunction}
Notations as above, let $p,q,x$ be pairwise distinct points on the metrized complex $\mccan$ associated to $\curve$. Then, the $j$-function is the solution of the hybrid Poisson equation
\begin{equation} \label{eq:SystemHybridJFunction}
\begin{cases}
 \Deltatrop (\lf) = \bm{\delta}_p - \bm{\delta}_q \\[1 mm]
 \lf \text{ is harmonically arranged} \\[1mm]
 \int_\curve \lf \, d\bm{\delta}_x =0.
 \end{cases}
\end{equation}
That is, the functions $f_i \colon \mccan \to \R$ given by $f_i(y) := \ljfuncbis{p\tiret q, x}(y)_i$, $y \in \mccan$, are precisely the components of the solution $\lf$ (viewed as functions on the metrized complex $\mccan$).
\end{thm}
\begin{proof} The proof is analogous to the one of Theorem~\ref{thm:TropicalHPvsJFunction}. For the first $r$ functions $f_1, \dots, f_r$, the claim follows directly from Theorem~\ref{thm:TropicalHPvsJFunction}, since the tropical part $\layg^\trop$ of the solution $\layg$ to \eqref{eq:SystemHybridJFunction} satisfies the Poisson equation \eqref{eq:TropicalHPvsJFunction} on the underlying tropical curve $\curve^\trop$ of $\curve$, and the first $r$ component of $\lhp{\curve}{p-q, y-x}'$ form the height pairing on $\curve^\trop$ (see Remark~\ref{rem:TropicalVSHybridHeightPairing}). 
It remains to prove the respective properties of the function $f_\smallcc \colon \mccan \to \R$. For two points $y, y'$ on the same Riemann surface component $C_u$ of $\curve$,
\[
f_\smallcc(y) - f_\smallcc(y') = \hp{C_u}{D^\smallcc_{p\tiret q}, \, y- y'}',
\]
where $D^\smallcc_{p\tiret q}$ is the divisor on $C_u$ from \eqref{eq:ComplexDivisorsHybridCurve} for $D = p-q$ by \eqref{eq:ComplexDivisorsHybridCurve}. It follows that $f_\smallcc\rest{C_u}$ is a solution to $\Delta g = D^\smallcc_{p\tiret q, u}$ on $C_u$ (see \eqref{eq:CHPvsJFunc}) and the values at the nodes $p^e_u$, $e \sim u$, are precisely the regularized values. Using \eqref{eq:LinearityDivisor}, one easily checks that $f_\smallcc$ is linear on the edges of $\mccan$.

We are left to prove that $f_\smallcc$ is harmonically arranged. In analogy to Remark~\ref{rem:TechnicalRemark}, a function $h \colon \pi_\smallcc \to \R$ is harmonically arranged exactly when
\[
{h(p^e_u) - h(p^e_v)} = h(D^\smallcc_{v\tiret u}) \quad \text{for all edges $e = uv$},
\]
where $D^\smallcc_{v\tiret u} = \sum_{w \in V} D^\smallcc_{v \tiret u, w}$ is the divisor on $\pi_\smallcc$ obtained from $y-x$ by \eqref{eq:ComplexDivisorsHybridCurve}. However, 
\[
{f_\smallcc(p^e_u) - f_\smallcc(p^e_v)} =  \sum_{w \in V} \hp{C_w}{D^\smallcc_{p\tiret q, w}, \, D^\smallcc_{v\tiret u, w}}' =f_\smallcc(D^\smallcc_{v\tiret u})
\]
by the very definition of $f_\smallcc$ and hence the proof is complete.
\end{proof}

\subsection{Tameness of the height pairing and $\jvide$-function} \label{ss:ResultsHybridHeightPairing}

Let $(S_0, p_0, q_0, x_0, y_0)$ be a stable Riemann surface with four marked points. Consider a versal deformation space $\rsf \to B$, $B=\Delta^{N}$ with $N =3g+1$, equipped with a fixed adapted system of coordinates. Let $\rsf^\hyb$ be the corresponding family of marked hybrid curves over $B^\hyb$.

By marking, we mean that for each point $\thy$ in  $B^\hyb$, we consider the corresponding marked points $p_\thy, q_\thy, x_\thy, y_\thy$ on the fiber $\rsf^\hyb_\thy$, which are all distinct and lie on the smooth part of $\rsf_t$ (with $t$ the point in $B$ underlying $\thy = (l,t)$). Altogether, we obtain sections $x, y, p,q$ of $\rsf^\hyb \to B^\hyb$ given by the marked points. 

\smallskip

As the next theorem shows, the hybrid height pairing describes the limit of the height pairing on degenerating families of Riemann surfaces.

\begin{thm}[Stratumwise tameness of the height pairing] \label{thm:tameness-height-pairing1} The function
\[
\thy \in B^\hyb \, \mapsto \, \lhp{\rsf^\hyb_\thy}{x_\thy - y_\thy, p_\thy - q_\thy} \in \R^{r_\thy +1}
\]
on $B^\hyb$ is stratumwise tame ($r_\thy$ is the rank of the hybrid curve $\rsf^\hyb_\thy$).
\end{thm}
Stratumwise tameness here means that the associated family $(\lf_\thy)_{\thy \in B^\hyb}$ of constant functions on $\rsf^\hyb \to B^\hyb$ is stratumwise tame (each piece $f_{j, \thy}$ of $\lf_\thy$ is a constant function on the respective minor). Equivalently, for every hybrid point $\shy$ in the closure $\bar D_\pi^\hyb$ of a boundary stratum $D_\pi^\hyb$,
\begin{align} \label{eq:StratumwiseTamenessHeightPairingSpelldOut}
&\Big | \hp{\rsf_t}{x_t -  y_t, p_t - q_t } - \sum_{j \in \{1, \dots, r_\thy, \smallcc\}} L_j(t) \lhp{\rsf^\hyb_\thy, j}{x_\thy - y_\thy, p_\thy - q_\thy} \Big | \to 0, &\thy = \loghyb_\pi(t),
\end{align}
as $t \in B^\ast$ converges tamely to $\shy$ in $B^\hyb$. 

\medskip

Similarly, let  $(S_0, p_0, q_0, x_0)$ be a stable Riemann surface with three marked points. Consider a versal deformation space $\rsf \to B$, $B=\Delta^{N}$, with $N =3g$, equipped with a fixed adapted system of coordinates, and let $\rsf^\hyb \to B^\hyb$ be the associated family of marked hybrid curves. For each $\thy \in B^\hyb$, let $\ljfuncbis{p_\thy - q_\thy, x_\thy}$ be the solution of the hybrid Poisson equation \eqref{eq:SystemHybridJFunction} for the marked points $p_\thy, q_\thy, x_\thy$ on the hybrid curve $\rsf^\hyb_\thy$.   

\smallskip

The following result strengthens Theorem \ref{thm:tameness-height-pairing1}.

\begin{thm}[Stratumwise tameness of the $\jvide$-function] \label{thm:TamenessJFunc}
The functions $\ljfuncbis{p_\thy - q_\thy, x_\thy}$, $\thy \in B^\hyb$, form a stratumwise tame family of functions.
\end{thm}
Theorem~\ref{thm:TamenessJFunc} implies that $\jfunc{p_t\tiret q_t,x_t}$ is asymptotically given by the pullback $\ljfuncbis{p_\thy - q_\thy, x_\thy}^\ast$, that is,
\begin{align*}
&\sup_{y \in \rsf_t}  \Big | \jfunc{p_t\tiret q_t,x_t}(y) - \ljfuncbis{p_\thy - q_\thy, x_\thy}^\ast(y)  \Big | \to 0,  &\thy = \loghyb_\pi(t), 
\end{align*}
 if $t \in B^\ast$ converges tamely to a point $\shy$ in a hybrid stratum $D_\pi^\hyb$ of $B^\hyb$.

\medskip
 
In fact we will establish the following stronger result. Let $\rsf^\hyb / B^\hyb$ be a hybrid versal family associated to a stable marked Riemann surface $S_0$, equipped with an adapted system of coordinates. Fix a boundary stratum $D_\pi^\hyb$ associated to an ordered partition $\pi = (\pi_1, \dots, \pi_r)$ of a subset $E_\pi \subseteq E$. Let $\Lognoind = \loghyb_\pi$ be the  stratumwise log map for $D_\pi^\hyb$.

\begin{thm} [Uniform stratumwise tameness of the height pairing]  \label{thm:FinalTamenessHeightPairing}
Suppose that $t \in B^\ast$ converges tamely to a point $\shy$ lying in the closure $\bar D_\pi^\hyb$ of $ D_\pi^\hyb$. Then for $p,q,x,y \in \rsf_t$,
\begin{align*}
 \hp{\rsf_t}{p - q , x - y } = & \sum_{j \in \{1, \dots, r, \smallcc\}}^r L_j(t) \,  \lhp{\rsf^\hyb_\thy, j}{\Lognoind(p) - \Lognoind(q) , \Lognoind(
 x) -  \Lognoind(y)} + \varepsilon_{\pi}(p,q,x,y) + o(1) 
\end{align*}
where  $\thy = \Lognoind (t)$, the correction term $\varepsilon_\pi(p,q,x,y)$ is defined below, and the $o(1)$ term goes to zero uniformly with respect to $p,q,x,y \in \rsf_t$
\end{thm}
The \emph{correction term} $\varepsilon_\pi (p,q,x,y)$ appearing in Theorem~\ref{thm:FinalTamenessHeightPairing} is defined as follows.  For four points $p,q,x,y$ on a smooth fiber $ \rsf_t$, $t \in B^\ast$ (with $\{p, q\} \cap \{x,y\} = \varnothing$),
\begin{equation} \label{eq:CorrectionFourPoints}
\varepsilon_\pi (p,q,x,y) := \varepsilon_\pi (p,x) + \varepsilon_\pi (q,y) - \varepsilon_\pi (p,y) - \varepsilon_\pi (q,x),
\end{equation}
where for an ordered partition $\pi$ of full sedentarity on $E_\pi \subseteq E$ and $p',x' \in \rsf_t$, we denote by $\varepsilon_\pi (p', x') = \sum_{e \in E_\pi} \varepsilon_e(p',x')$ the correction term introduced in \eqref{eq:CorrectionTwoPoints}.

\medskip

One can view the statement of Theorem~\ref{thm:FinalTamenessHeightPairing} as the tameness of the corrected pairing $\hp{\rsf_t}{p-q,x-y} - \varepsilon_\pi(p,q,x,y)$ at the stratum $D_\pi^\hyb$. For this reason, we refer to it shorthand as the \emph{uniform (stratumwise) tameness of the height pairing}.

\medskip

In the rest of this section we recall the Hodge-theoretic interpretation of the height pairing~\cite{Hain, ABBF}. Based on this, we prove Theorem~\ref{thm:tameness-height-pairing1} in Section~\ref{sec:ProofHeightPairing}, and Theorem~\ref{thm:FinalTamenessHeightPairing} in Section~\ref{sec:HybridJFctAsymptotics}. Theorem~\ref{thm:TamenessJFunc} follows by applying  Theorem~\ref{thm:FinalTamenessHeightPairing} to $p= p_\thy$, $q = q_\thy$, $x = x_\thy$ and taking into account the connection between the hybrid height  pairing and $j$-function (see Section~\ref{ss:HybridJFunction}).

\subsection{The logarithmic 1-form $\omega_{p\tiret q}$} \label{ss:LogFormOmega} Let $S$ be a Riemann surface of genus $g$ and let $p,q$ be two points on $S$. Fix further a symplectic basis $a_1, \dots, a_g$, $b_1, \dots, b_g$ of $H_1(S, \Z)$ (see Section \ref{sec:preliminaries}) that we lift to elements of the same name in $H_1(S \setminus\{p,q\}, \Z)$. In particular,  for $j= 1,\dots, g$, $a_j$ does not intersect $p$ and $q$. Let $\Omega$ be the period matrix of $S$ with respect to this basis, so that we have 
\[\Omega  = \Bigl(\int_{b_i}\omega_j\Bigr)\]
for holomorphic one forms $\omega_1, \dots, \omega_g$ verifying $\int_{a_i}\omega_j =\delta_{i,j}$.

\smallskip

  We define $\omega_{p\tiret q}$ to be the unique logarithmic one-form on $S$ with $\Res(\omega_{p\tiret q}) = p-q$ and
\begin{equation} \label{eq:CondAlphaInt}
	\int_{a_i} \omega_{p\tiret q} = 0, \qquad i=1,\dots,g.
\end{equation}

\begin{prop}\label{prop:two-forms} 
Notations as above, let $p, q$ be two points on $S$. Then, we have 
\begin{equation} \label{eq:RelationOmegaAlpha}
	\alpha_{p \tiret q} = \omega_{p \tiret q} - \sum_{j=1}^g  ( \Im(\Omega)^{-1} \Im(\rmW))_j \, \omega_j,
\end{equation}
where the column vector $\rmW = (\rmw_j)_{j=1}^g \in \C^g$ is given by
 \[
	\rmw_j =  \int_{b_j} \omega_{p \tiret q}, \qquad j = 1,\dots,g.
\]
\end{prop}

\begin{proof}
Since $\omega_{p\tiret q}$ and $\alpha_{p\tiret q}$ have the same residue, their difference is a holomorphic one-form on $S$. In particular, we can write
\[
	\alpha_{p\tiret q} = \omega_{p\tiret q} + \sum_{j=1}^g \lambda_j \omega_j
\]
for some coefficients $\lambda_j \in \C$, $j= 1,\dots,g$. It remains to show that $\lambda_j$ is given by \eqref{eq:RelationOmegaAlpha}. By integration over the cycles $a_j$, $j=1,\dots,g$, and using that $\alpha_{p\tiret q}$ has only real periods, the condition \eqref{eq:CondAlphaInt} implies that
\[
0 = \Im \big (  \int_{a_j} \alpha_{p\tiret q}  \big ) =   \Im \big (  \sum_{k=1}^g \lambda_k  \int_{a_j} \omega_k ) = \Im(\lambda_j).
\]
This shows that $\lambda_j$ are all real. Integrating now over the cycles $b_j$, $j=1, \dots, g$, leads to
\begin{align*}
0 &= \Im \Big (  \int_{b_j} \alpha_{p\tiret q}  \Big ) = \Im(\rmw_j) + \sum_{k=1}^g \Im \Big ( \lambda_k \Omega_{j,k} \Big ) = \Im(\rmw_j) + \sum_{k=1}^g \lambda_k \Im( \Omega_{j,k}).
\end{align*}
Since $\Im(\Omega)$ is symmetric, $\Im( \Omega_{j,k}) = \Im(\Omega)_{j,k}$, and \eqref{eq:RelationOmegaAlpha} follows immediately.
\end{proof}

\subsection{Biextensions}\label{sec:biextensions} 

Let $S$ be a Riemann surface of genus $g$, and let $x,y, p,q$ be four distinct points on $S$. The degree zero divisors $x-y$ and $p-q$ give rise to a mixed Hodge structure $M={}_{p\tiret q}M_{x\tiret y}$, called a biextension~\cite{Hain}. The mixed Hodge structure $M$ has graded pieces of  weights $-2, -1, 0$ given by
\begin{equation}\label{eq:biextension-graded-pieces}
\grm{-2}\weight\, M=\Z(1), \quad \grm{-1}\weight\,M=H^1(S, \Z(1)), \quad \grm{0}\weight M=\Z(0).
\end{equation}
The first and last isomorphisms are given by the divisors $p-q$ and $x-y$, respectively. 
\smallskip

\begin{defn}\rm The biextension $M={}_{p\tiret q}M_{x\tiret y}$ associated to the divisors $p-q$ and $x-y$ is the relative cohomology group $H^1(S \setminus \{x,y\}, \{p,q\}; \Z(1))$. 
\end{defn}

To prove that it has weight filtration verifying the conditions of biextensions, we first note that there exists a short exact sequence 
\[ 0 \rightarrow  H^1(S, \Z(1))\rightarrow H^1(S\setminus \{x,y\}, \Z(1)) \rightarrow   \Z(0) = \ker\Bigl(\Z x\oplus \Z y \to \Z\Bigr)  \rightarrow 0,
 \] 
 which is obtained from the long exact sequence of cohomology groups for the inclusion $S\setminus \{x,y\} \hookrightarrow S$.

In the same way, we get an exact sequence of mixed Hodge structures 
\[
0\rightarrow \Z(1) = \mathrm{coker}\Bigl( \Z(1) \to \Z(1) p \oplus \Z(1) q \Bigr)\rightarrow H^1(S, \{p,q\};\Z(1)) \rightarrow H^1(S, \Z(1)) \to 0.
\]

\smallskip

The biextension $M={}_{p\tiret q}M_{x\tiret y} = H^1(S \setminus{x,y}, \{p,q\}, \Z(1))$ fits into the diagram  
\begin{equation}\label{eq:biextension}
\xymatrix{ & & 0  & 0  & \\ 
& & \Z(0) \ar@{=}[r]\ar[u] & \Z(0) \ar[u] & \\
0 \ar[r] & \Z(1) \ar[r] &  M={}_{p\tiret q}M_{x\tiret y} \ar[u] \ar[r] & H^1(S, \{x,y\}) \ar[u] \ar[r]& 0 \\
0 \ar[r] & \Z(1) \ar@{=}[u] \ar[r] & H^1(S,\{p,q\};\Z(1)) \ar[u] \ar[r] & H^1(S, \Z(1)) \ar[u] \ar[r]& 0\\
& & 0 \ar[u] & 0\ar[u] &
}
\end{equation}
The weight filtration is given by 
\begin{equation*}
W_{-3}=(0) \subset W_{-2}=\Z(1) \subset W_{-1}=H^1(S, \{p,q\}; \Z(1)) \subset W_0=M,
\end{equation*} and satisfies Equation~\ref{eq:biextension-graded-pieces}, that is,
$$
\grm{-2}\weight M =\Z(1), \quad \grm{-1}\weight M =H^1(S, \Z(1)), \quad \grm{0}\weight M=\Z(0). 
$$

\subsection{Period matrix of the biextension  ${}_{p\tiret q}M_{x\tiret y}$} \label{PeriodMatrixBiextension}
In this section, we describe the period matrix of the biextension  mixed Hodge structure ${}_{p\tiret q}M_{x\tiret y}$ given by a configuration of four points $p,q,x,y$ on a Riemann surface $S$ of genus $g$. This will be used in the next section to describe the period map for the variation of biextensions over $B^\ast$.

\medskip

 Let $ a_1,\dots,a_g,b_1,\dots,b_g $ be a symplectic basis for $H_1(S, \Z)$ that  we lift to $H_1\bigl(S\setminus \{p,q\}, \{x,y\}; \Z\bigr)$. Let $\gamma_p$ and $\gamma_q$ be two small cycles around $p$ and $q$, respectively, with orientation induced by that of the surface. Note that $\gamma_p+\gamma_q =0$ in $H_1\bigl(S\setminus \{p,q\}, \{x,y\}; \Z\bigr)$. Let $\gamma\indbis{x\tiret y}$ be a path on $S$ with boundary $x-y$ such that 
\[\intprod{}{a_j , \gamma\indbis{x\tiret y}} = 0, \quad j=1, \dots, g.\]
We view $\gamma\indbis{x\tiret y}$ as an element of $H_1\bigl(S\setminus \{p,q\}, \{x,y\}; \Z\bigr)$. Note that the homology class of $\gamma\indbis{x\tiret y}$ is unique up to a linear combination of elements $a_1, \dots, a_g$. 

\smallskip

 There exists a unique collection of holomorphic one-forms $\omega_1, \dots, \omega_g$ on $S$, forming a basis of $H^0(S, \omega_S)$, such that 
 \[\int_{a_i} \omega_j = \delta_{i,j}, \quad i,j\in\{1,\dots, g\}.\]
  The period matrix of $S$ is then 
  \[\Omega = \Bigl( \int_{b_i} \omega_j\Bigr)_{i,j\in [g]}.\]

Let as in the previous section, $\omega_{p\tiret q}$ be the unique section of $\omega_S(p+q) = \omega_S(\log(p+q))$, seen as a logarithmic one-form, with residue divisor $p-q$ and with vanishing periods 
\[\int_{a_j} \omega_{p\tiret q} = 0, \quad \forall j=1, \dots, g.\]

\begin{defn} \rm The \emph{period matrix of the biextension mixed Hodge structure ${}_{p\tiret q}M_{x\tiret y}$}, with respect to the symplectic basis $a_1, \dots, a_g, b_1, \dots, b_g$ and $\gamma_{x\tiret y}$, is the $(g+1)\times (g+1)$ matrix $\bipm$ given by 
\[\bipm : = \begin{pmatrix} \Omega & \rmW \\
\rmZ & \rho
\end{pmatrix}\]
where \begin{itemize}
\item the $g\times g$ matrix $\Omega = \Bigl(\int_{b_j} \omega_i\Bigr)$ is the period matrix of $S$;
\item $\rmW = \Bigl(\int_{b_j}\omega_{p\tiret q}\Bigr)_{j=1}^g$ is a column matrix in $\C^g$,
\item $\rmZ = \Big( \int_{\gamma\indbis{x\tiret y}} \omega_j\Bigr)_{j=1}^g$ is a row matrix in $\C^g$, and
\item $\rho = \int_{\gamma\indbis{x\tiret y}}\omega_{p\tiret q}$ is a scalar in $\C$.
\end{itemize}  
\end{defn}

We have the following Hodge-theoretic interpretation fo the height pairing~\cite{Hain, ABBF}.
\begin{prop} \label{prop:HeightPairingMatrix} The height pairing $\hp{S}{x-y, p-q}$ is given by  
\[ \hp{S}{x-y, p-q} = 2\pi \big (  \Im(\rho) -  \Im(\rmZ) \Im(\Omega)^{-1} \Im(\rmW) \big ).\]

\end{prop}
\begin{proof} This follows from Proposition~\ref{prop:two-forms}. 
Indeed, by that result, we have 
\[\alpha_{p\tiret q} = \omega_{p\tiret q} - \sum_{j=1}^g \Bigl(\Im(\Omega)^{-1} \Im(\rmW)\Bigr)_j\omega_j.\]
Now, by definition of the height pairing, we have 
\begin{align*}
\frac{1}{2\pi}\hp{S}{x-y, p-q} &= \Im\bigl(\int_{\gamma\indbis{x\tiret y}} \alpha_{p\tiret q}\Bigr) = \Im\Bigl(\int_{\gamma\indbis{x\tiret y}}\omega_{p\tiret q}\Bigr) - \sum_{j=1}^g \Bigl(\Im(\Omega)^{-1} \Im(\rmW)\Bigr)_j\Im\Bigl(\int_{\gamma\indbis{x\tiret y}}\omega_j \Bigr)\\
&= \Im(\rho) - \Im(\rmZ) \Im(\Omega)^{-1} \Im(\rmW).
\end{align*}
\end{proof}

This brings the proof of Theorems~\ref{thm:tameness-height-pairing1} and~\ref{thm:FinalTamenessHeightPairing} to the problem of understanding the asymptotics of the period matrix.


\subsection{Regularized integrals of forms over paths and the regularized height pairing} \label{ss:RegularizedIntegralsPaths}
For our later steps, in dealing with the asymptotic description of the height pairing, it will be convenient to provide an analogue of Proposition~\ref{prop:HeightPairingMatrix} for the regularized height pairing. Let $S$ be a smooth Riemann surface of genus $g$. 

\smallskip

We start by introducing a \emph{regularized pairing between one-chains and logarithmic one-forms}.

\smallskip

Suppose $\gamma$ is a one-chain on $S$  with boundary divisor $D_\gamma$, and $\omega$ is a differential of the third kind on $S$  with residue divisor $D_\omega$. Fix a local parameter $\zeta_x$ around each point $x$ in the intersection of the supports $\suppdiv_\gamma = \supp{D_\gamma}$ and $\suppdiv_\omega = \supp{D_\omega}$. We now define a new "shorter" one-chain $\gamma'$ by taking a small disc of radius $|\zeta_x|$ around each $x \in \suppdiv_\gamma \cap\suppdiv_\omega$ (in the local coordinate $\zeta_x$), and removing all parts of the one-chain $\gamma$ inside this open disc.

\smallskip

We define the \emph{regularized integral} of $\omega$ over $\gamma$ as the limit
\[
\hp{S}{ \gamma, \omega}' :=  \lim_{\substack{\zeta_x \to 0 \\ \forall x\in  \suppdiv_\gamma \cap\suppdiv_\omega }} \Im \Big( \int_{\gamma'}  \omega  \Big ) + \frac{1}{2\pi} \sum_{x\in  \Sigma_\gamma \cap\Sigma_\omega } D_\gamma(x) D_\omega(x)  \log|\zeta_x| . 
\]
Note that taking the imaginary part of the integral $\int_\gamma \omega$ ensures that the limit exists.

The analogue of Proposition~\ref{prop:HeightPairingMatrix} now reads as follows.

\smallskip

Fix four points $x,y,p,q$ on $S$ with $x \neq y$ and $p \neq q$.  Let $ a_1,\dots,a_g,b_1,\dots,b_g $ be a symplectic basis for $H_1(S, \Z)$ and denote by $\Omega$ its period matrix. As in the previous section, let  $\omega_{p\tiret q}$ be the unique logarithmic one-form on $S$ having residue divisor $p-q$ and vanishing periods on the $a_j$'s. Moreover, let $\gamma_{x\tiret y}$ be a path in $S$ from $x$ to $y$ such that $\langle \gamma_{x\tiret y}, a_j \rangle = 0$ for all $j$.

\begin{lem} \label{lem:RegularizedHeightFormula}
The regularized height pairing $\hp{S}{x-y, p-q}'$ is given by  
\[\hp{S}{x-y, p-q}' = 2\pi \Big ( \hp{S}{ \gamma_{x \tiret y}, \omega_{p \tiret q} }' -  \Im(\rmZ) \Im(\Omega)^{-1} \Im(\rmW) \Big ),\]
\end{lem}
where
\begin{itemize}
\item $\rmW = \Bigl(\int_{b_j}\omega_{p\tiret q}\Bigr)_{j=1}^g$ is a column matrix in $\C^g$,
\item $\rmZ = \Big( \int_{\gamma\indbis{x\tiret y}} \omega_j\Bigr)_{j=1}^g$ is a row matrix in $\C^g$, and
\end{itemize}  
\begin{proof}
This is clear from Proposition~\ref{prop:HeightPairingMatrix}, the definition of the regularized height pairing and the definition of the regularized integrals.
\end{proof}


\section{Degenerations of Riemann surfaces} \label{sec:monodromy}

For later use,  we gather here some well-known results about degenerations of Riemann surfaces, slightly generalized in few cases to our setting, and refer to~\cite{AN, ABBF, ACG, Hof84} for more details.

We resume the discussion, and use the terminology, from Section~\ref{sec:deformations}. Let  $(S_0, q_1, \dots, q_\nmark)$ be a fixed stable curve of arithmetic genus $g$ with $\nmark$ marked points. 
  Let $G=(V,E,\genusfunction,\marking)$ be the dual graph of $S_0$ with the genus function $\genusfunction: V \to \mathbb Z_{\geq0}$ and the marking function $\marking: [\nmark] \to V$. For each vertex $v\in V$, we denote by $C_v$ the normalization of the corresponding irreducible component of $S_0$. 
  For each edge $e =uv$, we denote by $p^e_u$ and $p^e_v$ the corresponding points in $C_u$ and $C_v$, respectively. For each vertex $v$, let $\mathcal A_v$ be the collection of all the points $p^e_v$ for incident edges $e$ to $v$. 
  \smallskip

In what follows, we will denote by $B$ the base of the versal deformation $\rsf$ of  the stable marked curve $(S_0, q_1, \dots, q_\nmark)$, and define $B^* : = B \setminus \bigcup_{e \in E} D_e$ and $\rsf^\ast := \pr^{-1}(B^\ast)$ to be the locus of points with smooth fibers in the family. Moreover, we denote by $\LS$ the relative homology $H_1(\rsf_t, \{\sigma_1(t), \dots, \sigma_\nmark(t)\}, \Z)$, with $\sigma_j(t)$ the corresponding sections of $\rsf \to B$ given by the markings. This provides a local system on $B^*$, and we denote by $\LS_t$ this homology group, which is the fiber of $\LS$ at $t$.

\medskip

We fix a base-point $\bp \in \base^*$. The fundamental group $\pi_1\bigl( \base^* , \bp\bigr)$ is naturally isomorphic to $\Z^E$: for each edge $e\in E$, we have a generator $\lambda_e$ given by a simple loop $\lambda_e\subset \base^*$ based at $\bp$ which turns once around the divisor $D_e$, counter-clockwise with respect to the parameter $z_e$ in the corresponding disk, and is moreover contractible in the space $\base \setminus \bigcup_{e'\neq e}D_{e'}$.

\subsection{Monodromy action} The action of $\lambda_e$ on $\LS_\bp$ can be described by Picard-Lefschetz theory as follows. For any oriented edge $e \in \E$ of the graph, there is a canonical vanishing cycle element $a_e$ in $\LS_\bp$ associated to the singular point $p_e$ of $S_0$. For the same edge with opposite orientation $\bar e$, we have $a_{\bar e} =-a_e$. The monodromy action of $\lambda_e$ is then given by
\begin{align}\label{eq:PL}
&\lambda_e: \LS_\bp \longrightarrow \LS_\bp\\
\beta &\mapsto \beta-\langle \beta,a_e\rangle a_e,
\end{align}
for any $\beta$ in $\LS_\bp$, where $\langle\cdot\,, \cdot \rangle$ is the intersection pairing between one-cycles in $\LS_\bp$. Note that this description implies that $\lambda_e =\lambda_{\bar e}$, for two oriented edges with the same underlying edge, that is, the definition of $\lambda_e$ is independent of the orientation of the edge $e$.

\subsection{Specialization} The inclusion $\Bigl(S_0, \{q_1, \dots, q_\nmark\}\Bigr) \hookrightarrow \Bigl(\rsf, \{\sigma_1, \dots, \sigma_\nmark\}\Bigr)$ admits a deformation retraction to $(S_0,\{q_1, \dots, q_\nmark\})$. Using this, we get the specialization map
\begin{equation*}
\mathrm{sp} \colon \LS_\bp \longrightarrow \LS_0
\end{equation*}
where $\LS_0$ denotes $H_1(S_0, \{q_1, \dots, q_\nmark\}, \Z)$. 

\smallskip

As an application of the Clemens-Schmid exact sequence in asymptotic Hodge theory, we infer that the specialization map above is surjective and its kernel $A\subset \LS_\bp$ corresponds to the subspace spanned by the vanishing cycles $a_e$. This gives an exact sequence
\begin{equation*}
0\to A \to \LS_\bp \xrightarrow{\mathrm{sp}} \LS_0 \to 0.
\end{equation*}

The inclusion of curves $C_v \hookrightarrow S_0$ induces a surjection map $\LS_0 \to H_1(G, \suppdiv, \Z)$ whose kernel we denote by $\LS^\nr_0$, the homology group relative to the markings of the normalized curve obtained from $S_0$. Here, $\suppdiv$ is the support of the counting function, that is, the set of all vertices $v$ such that $C_v$ contains at least one marked point. (In particular, if $\nmark=0$, $\suppdiv =\varnothing$.)  

\noindent More precisely, $\LS^\nr_0$ has the following description. We have a short exact sequence of the form 
\[0 \to \bigoplus_{v\in V} H_1\bigl(C_v, \suppdiv_v, \Z\bigr) \hookrightarrow H_1(S_0, \{q_1, \dots, q_\nmark\}, \Z) \to H_1(G, \suppdiv, \Z) \to 0,\]
where for $v\in V$ a vertex of the graph, $\suppdiv_v \subseteq\{q_1, \dots, q_\nmark\}$ denotes the set of all marked points of $S_0$ which lie on the component $C_v$. This identifies $\LS^\nr_0$ with the direct sum $\bigoplus_{v\in V} H_1\bigl(C_v, \suppdiv_v, \Z\bigr)$.

\medskip

We define now 
\[A':=A + \mathrm{sp}^{-1}\Bigl(\LS^\nr_0\Bigr)\subseteq \LS_\bp.\]

It follows from the two proceeding exact sequences that
\begin{equation}
 \label{eq:aux}
\LS_\bp/A' \simeq H_1(G,\suppdiv, \Z) .
\end{equation}

\smallskip

The subspace of vanishing cycles $A\subset \LS_\bp$ has rank $\graphgenus$, the genus of the graph $G$, and is isotropic, meaning that for any pair of elements $a, b \in H$, $\langle a, b\rangle = 0$. Moreover, we have $\langle A,A'\rangle =0$.
We thus get a pairing $A \times \bigl(\LS_\bp/A'\bigr) \rightarrow \Z$. 

\smallskip

We will provide later in this section all the relations between the vanishing cycles in terms of the \emph{oriented cuts} in the dual graph, as well as an interesting family of bases for the space of vanishing cycles defined in terms of the spanning trees of the graph $G$.

\subsection{Logarithmic  monodromy action} Let $N_e := \id -\lambda_e$. By \eqref{eq:PL}, for any
 $\beta \in \LS_\bp$, we have
$N_e(\beta)=\langle \beta,a_e\rangle a_e$. The image of $N_{e}$ is contained
in $A$ and we get $N_e\circ N_e =0$.
We infer that $N_e=-\log(\lambda_e)$ and $N_{e}$ is vanishing on $A'$.

\medskip

It follows that $N_e : \LS_\bp  \to A$ passes to the quotient by $A'$ and induces 
\begin{equation} \label{eq:nm}
H_1(
G,\suppdiv, \Z) \simeq \LS_\bp/A' \xrightarrow{N_e} A \to 
(\LS_\bp/A')^\vee \simeq H_1(G, \suppdiv, \Z)^\vee.
\end{equation}

\begin{prop} The induced bilinear form on $H_1(G, \suppdiv, \Z)$ given by the composition of the maps in \eqref{eq:nm}
  coincides with the bilinear form $\langle\cdot\,,\cdot\rangle_e$ restricted to $H_1(G, \suppdiv,\Z)$.
\end{prop}

Recall that $H_1(G, \suppdiv,\Z) \subseteq \Z^E$ consists of all the elements $\gamma\in \Z^E$ with $\supp{\partial(\gamma)} \subseteq \suppdiv$.

\begin{proof} This is a slight generalization of the same result in the case $\suppdiv =\varnothing$ proved in~\cite{ABBF}. For any $\beta \in \LS_\bp = H_1\bigl(\rsf_{\bp}, \{\sigma_1(\bp), \dots, \sigma_\nmark(\bp)\}, \Z\bigr)$, 
the image $\mathrm{sp}(\beta)$   in the quotient
  $\LS_\bp/A' \simeq H_1(G, \suppdiv, \Z)$
 can be identified with the 1-chain $\gamma$ in the graph with boundary in $\suppdiv$. Choose an orientation of the edges and write $\gamma=\sum_e n_ee$. To find the coefficients $n_e$, we consider 
  $\langle \gamma, a_e\rangle =\langle \beta, a_e\rangle$. For an edge $e' \neq e$, there is no intersection between $e' $ and $a_e$. Moreover, the edge $e$ is transversal to $a_e$. This leads to the equation $n_e= \langle \beta,a_e\rangle$.
The quadratic form on $H_1\bigl(\rsf_{\bp},\{\sigma_1(\bp), \dots, \sigma_\nmark(\bp)\}, \Z\bigr)$ associated to $N_e$
  sends $\beta$ to $\langle \beta, \langle \beta,a_e\rangle a_e\rangle =
  n_e^2 = q_e(\gamma)$, where $q_e$ is the quadratic form on  $H_1(G, \suppdiv, \Z)$  associated to $\langle\cdot\,,\cdot\rangle_e$, and the
  proposition follows.
\end{proof}

\subsection{Equations for vanishing cycles} \label{ss:EquationsVanishingCycles} Consider the graph $G=(V, E)$. For any subset $W \subset V$, let $\E(W, V \setminus W)$ denote the \emph{oriented cut} defined by $W$, that is the set of all the oriented edges $e \in \E$ with tail $\tail_e$ in $W$ and head $\head_e$ in $V\setminus W$. 

\smallskip

Let $\Z^\E$ be the free abelian group generated by the oriented edges $e\in \E$. Consider the map 
\[\Z^\E \to A, \qquad e \mapsto a_e.\]

\begin{prop}[Equations for vanishing cycles]\label{eq:equations-vanishing-cycles} The kernel of the above map is generated by the elements 
\[e + \bar e, \qquad \forall \, e\in \E\]
and
\[\sum_{e\in \E(W, V\setminus W)} e, \qquad \forall\,W \subset V.\]
\end{prop}
In other words, the equations 
\[\sum_{e\in \E(W, V\setminus W)} a_e =0, \qquad W \subset V,\]
combined with the equations $a_e = -a_{\bar e}$ for $e\in \E$ generate all the relations between vanishing cycles. 

We need the following proposition. 
\begin{lem} For a subset $W \subset V$ of vertices, consider the union of vanishing cycles $a_e$ for $e\in E(W, V \setminus W)$. The surface $\rsf_\bp\setminus \Bigl(\bigcup_{e\in E(S, S^c)}a_e\Bigr)$ obtained from $\rsf_\bp$ by cutting along the cycles $a_e$ has two connected components.  
\end{lem} 
\begin{proof}
The connected components of the surface $\rsf_\bp\setminus \Bigl(\bigcup_{e\in E(W, V\setminus W)}a_e\Bigr)$  are in bijection with that of the surface $S_0 \setminus\{p_e\}_{e\in E(W, V\setminus W)}$ from which the proposition follows. 
\end{proof}

\begin{proof}[Proof of Proposition~\ref{eq:equations-vanishing-cycles}]
Consider the vanishing cycles on $\rsf_\bp$ associated to the edges of an oriented cut $\E(W, V\setminus W)$ in the graph. By the previous lemma, cutting $S_\bp$ along these cycles results in a surface with two components. One of these components have boundary with induced orientation being precisely that of the vanishing cycles coming from the oriented cut. This gives the equation $\sum_{e\in \E(W, V\setminus W)} a_e =0$ in homology. The other equation $a_e = -a_{\bar e}$ was already observed before. Since $A \simeq \Z^h$ and the quotient of $\Z^\E$ with the above relations has rank $h$, the surjectivity of the map $\Z^\E \to A, e \mapsto a_e$, implies the proposition.
\end{proof}

\subsection{Spanning trees and bases of the space of vanishing cycles} \label{sec:basis_spanning_tree}
Proposition~\ref{eq:equations-vanishing-cycles} provides interesting bases for the space of vanishing cycles.  Consider a spanning tree $T$ in the graph $G$. For any edge $e$ in $E \setminus E(T)$, we give an orientation to $e$, and let $a_e$ be the corresponding vanishing cycle. Then,
\begin{prop}[Basis of the space of vanishing cycles associated to a spanning tree] 
Notations as above, the set $a_e$, $e\in E \setminus E(T)$, provides a basis of $A$.
\end{prop} 
\begin{proof} Since we have the right number of elements, that is $\graphgenus$, it will be enough to show that any element $a_e$ for $e\in \E(T)$ can be described in this basis. 

Removing $e$ from the tree $T$ gives two connected components. Let $W$ be the set of vertices of the component which gives $e$ in the oriented cut $\E(W, V\setminus W)$. Then, using the equations for vanishing cycles, we get 
\[a_e = -\sum_{\substack{f\in \E(W, V\setminus W)\\ f\neq e}} a_f.\] 
To conclude, note that for each each such $f$, either $a_f$ or $-a_f=a_{\bar f}$ is in 
is the basis. 
\end{proof}

For what we do in the sequel, it will be sometimes more convenient and transparent to fix a spanning tree and consider the corresponding basis of the space of vanishing cycles. 

\subsection{Admissible symplectic basis} \label{ss:AdmissibleSymplecticBasis} For the details of the materials which follow, we refer~\cite{AN}.

We numerate the vertices of the graph $G$ as $v_1, \dots, v_n$, and define the sequence $g_0, \dots, g_n$ by 
\[g_i = \graphgenus+  \genusfunction(v_1) + \dots +\genusfunction(v_i).\]
In particular, $g_0 = \graphgenus$ and $g_n =g$

For each $r=0, 1, \dots, n$, let $I_r = \{g_{r-1}+1, \dots, g_r\}$.

An \emph{admissible symplectic basis} $a_1,\dotsc,a_{g},b_1,\dotsc,b_{g}$ for
 $H_1(\rsf_{\bp}, \Z)$ is a symplectic basis with the following properties
 
 \begin{enumerate}
 \item[(1)] $a_1, \dots, a_\graphgenus$ form a basis of the space of vanishing cycles
 $A$, and 
 \smallskip
 
 \item[(2)] for any $r=1, \dots, n$, the collection of elements $a_{j}$ and $b_j$ for $j\in I_r$ provides a (symplectic) basis of $H_1(C_{v_r}, \Z)$ in 
 \[H_1(C_{v_r}, \Z) \hookrightarrow H_1(S_0, \Z) \simeq H_1(\rsf_{\bp}, \Z) /A.\]
 \end{enumerate} 
 
When the underlying dual graph has a layering $\pi = (\pi_1, \dots, \pi_r)$, that is, when working with hybrid curves, we impose the following additional admissibility condition.

\begin{itemize}
\item[(3)] The elements $b_1, \dots, b_{\graphgenus}$ form an admissible basis in the sense of Section~\ref{sec:admissible_basis_layered_graphs} for the homology group $H_1(G, \Z)$.
\end{itemize}

 \smallskip

 Let now $\suppdiv_\bp \subset \rsf_\bp$ be a finite set of points. We get an exact sequence
 \[
 0 \to \Z \to \Z^{\suppdiv_\bp} \to H_1(\rsf_\bp \setminus \suppdiv_\bp, \Z) \to H_1(\rsf_{\bp}, \Z) \to 0.
 \]
 For any point $p\in \suppdiv_\bp$, we choose a small cycle $\gamma_p$ around $p$ with orientation compatible with that of the surface. We lift the elements $a_1, \dots, a_g, b_1, \dots, b_g$ to $H_1(\rsf_\bp \setminus \suppdiv_\bp, \Z)$. 
 
 For any oriented edge $e\in \E$, the vanishing cycle $a_e$ canonically lifts to $H_1(\rsf_\bp \setminus \suppdiv_\bp, \Z)$. It can be further decomposed as a linear combination of the basis elements $a_i$, $i\in[\graphgenus]$, and $\gamma_p$, for $p\in \suppdiv_\bp$, 
\begin{equation}\label{eq:vanishing_cycle}
a_e = \sum_{i=1}^\graphgenus c_{e,i}a_i + \sum_{p\in \suppdiv_{\bp}} d_p \gamma_p.
\end{equation}
Moreover, the collection of numbers $d_p$ is unique up to addition by some constant. 

\medskip

One practical choice of the symplectic basis in (1) is the one described in Section~\ref{sec:basis_spanning_tree}. Fix a spanning tree $T$ of the graph $G$, and enumerate the edges of $G$ which are not in $T$ by $e_1, \dots, e_{\graphgenus}$ endowed each with a choice of an orientation. Consider the symplectic basis obtained by taking $a_j := a_{e_j}$, the vanishing cycle corresponding to the edge $e_j$. The advantage of this symplectic basis is that it allows to give a simple description of the coefficients in Equation~\refeq{eq:vanishing_cycle}. Indeed, for any oriented edge of the graph $G$, if $e$ is not an edge of the spanning tree, then 
\begin{equation}\label{eq:vanishing_cycle2}
a_e = \pm a_i,
\end{equation}
 for $i$ such that $e = e_i$ or $e = \bar e_i$. If $e$ is an orientation of an edge of the spanning tree, then removing $e$ from the tree gives a partition of $V$ into sets $W$ and $V\setminus W$ such that $e$ is the only edge of the tree between $W$ and $V \setminus W$, with $e$ having tail in $W$ and head in $W \setminus V$. We get the equation 
 \[\sum_{f\in \E(W, V\setminus W)} a_f + \sum_{p\in \suppdiv_{\bp,\mathrm{w}}} \gamma_p =0. \]
 Here $\suppdiv_{\mathrm{w}, \bp}$ is the set of all points $p\in \suppdiv_\bp$ which lie in the component of $\rsf_\bp \setminus \bigcup_{f\in \E(W, V\setminus W)} a_e$ corresponding to the vertex set $W$. 
 
 \medskip
 
\begin{figure}[!t]
\centering
   \scalebox{.3}{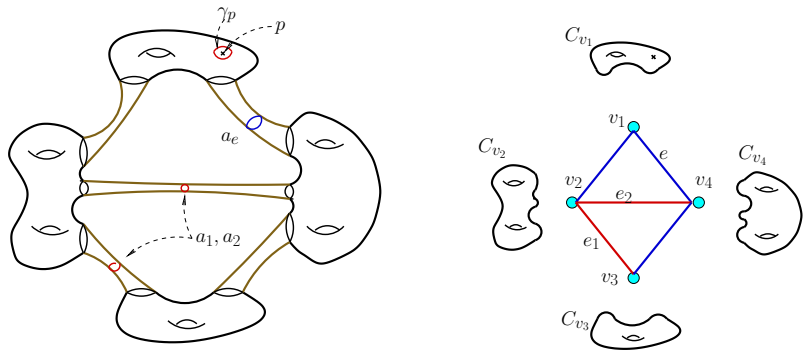}
\caption{Example of an equation in the space of vanishing cycles. The blue edges form a spanning tree. The red edges give rise to the basis $\{a_1, a_2\}$ for the space of vanishing cycles. The vanishing cycle $a_e$ associated to the edge $e$ satisfies the equation $a_e = - a_1 - a_2 -\gamma_p$ (for appropriate choice of orientations). The subset $W$ is composed of the two vertices $v_3$ and $v_4$.}
\label{fig:tree-vc}
\end{figure}

Therefore, we get the equation
 \begin{equation}\label{eq:vanishing_cycle3}
a_e = - \sum_{\substack{f \in \E(W, V\setminus W)\\ f\neq e}} a_f - \sum_{p\in \suppdiv_{\bp, 1}} \gamma_p.
\end{equation}
Note that the first sum above involves only oriented edges $f$ such that either $f= e_i$ or $f=\bar e_i$ for an $i\in[h]$ so that $a_f =\pm a_i$ for an element of the basis, see Figure~\ref{fig:tree-vc}.  

\smallskip

 The elements $b_j$, $j\in[h]$, has as well an explicit description. For any edge $e_j$ not in the spanning tree, with its chosen orientation, adding $e_j$ to the spanning tree creates a unique oriented cycle with orientation consistent with that of $e_j$. This is the element $b_j$.
In order for the elements $b_1, \dots, b_\graphgenus$ to verify (3), we  have to assume moreover that $T$ is a layered spanning tree, and using the notations of Section~\ref{sec:admissible_basis_layered_graphs}, we require that the edges $e_1, \dots, e_\graphgenus$ are enumerated so that the edges $e_i$ for $i\in J^j_\pi$ are in $\pi_j$. 

\section{The period map of a family of marked Riemann surfaces} \label{sec:period}
In this section, we introduce the \emph{period map for the variation of biextension mixed Hodge structures} associated to degree zero divisors on families of Riemann surfaces  $\rsf^\ast\to B^\ast$, and give a description of the action of the monodromy. A part of this discussion is already contained in~\cite{ABBF}. In the next Section~\ref{sec:period-map-asymptotics}, we will describe the asymptotics of the period map by using the nilpotent orbit theorem. Relying on a connection between the period map and the height pairing (see Section~\ref{ss:ImaginaryPartHeightPairing}), this approach eventually leads to the proofs of Theorem~\ref{thm:tameness-height-pairing1} and Theorem~\ref{thm:FinalTamenessHeightPairing} (see Section~\ref{sec:ProofHeightPairing} and Section~\ref{sec:HybridJFctAsymptotics}, respectively).

 \medskip

Let $(S_0, q_1, \dots, q_\nmark)$ be a stable Riemann surface with $\nmark$ marked points and with dual graph $G = (V, E)$. Consider the versal deformation space $\rsf/B$, $B=\Delta^{N}$ with $N =3g - 3 +\nmark$, and let $B^\ast = B \setminus \cup_{e\in E} D_e$. Shrinking the polydisk if necessary and making a choice of local parameters around $0$ for $D_e$, we can write  $B^\ast \simeq (\Delta^*)^E \times \Delta^{N -|E|}$.

\smallskip
Let further $p,q,x,y$ be four disjoint holomorphic sections of $\rsf \to B$, whose images $p_t$, $q_t$, $x_t$, $y_t$ lie in the smooth part of the surface $\rsf_t$ for all $t \in B$. For instance, if $S_0$ has at least  $\nmark \ge 4$ markings, then choosing four distinct marked points $q_i$, the associated sections of $\rsf \to B$ have this property.

 Let $\widetilde {B^\ast}$ be the universal cover of $B^\ast$. We get an isomorphism $\widetilde{B^\ast} \simeq \H^E \times \Delta^{N -|E|}$ where $\H$ is the Poincar\'e half-plane, and the map 
\[
\widetilde{B^\ast} \to B^\ast
\]
is given by sending $\zeta_e\in \H$ to $z_e=\exp(2\pi i \zeta_e) \in \Delta^*$.  Denote by $\widetilde{\rsf}^*$ the family of marked Riemann surfaces over $\widetilde B^\ast$ induced from the family $\rsf^\ast/B^\ast$

\medskip

In addition to the above, we also fix a spanning tree $T = (V(T), E(T))$ of the dual graph $G = (V,E)$ of $S_0$. Moreover, if an ordered partition is given, for example, if we are interested in a particular stratum in $B^\hyb$, we assume furthermore that this spanning tree is layered.  All choices of vanishing cycles, admissible basis, etc., are taken with respect to this spanning tree (see Section~\ref{sec:monodromy}). As we will see later, this approach leads to a particularly transparent expression for the mondromy action.

\subsection{The period map} \label{ss:PeriodMap} Notations as above, consider the family of stable marked Riemann surfaces $\rsf \to B$ with markings $p,q, x,y$. We get a variation of mixed Hodge structures ${}_{p_t\tiret q_t}M_{x_t\tiret y_t}$ over $B^\ast$ (see Section~\ref{sec:biextensions}). 
Being of geometric origin, this is an admissible variation of mixed Hodge structures. If there is no risk of confusion, we denote the corresponding local system over $B^\ast$ by
\[
\localsystem = H_1\Bigl(\rsf \setminus\{p,q\}, \{x,y\}; \Z\Bigr),
\]
the fiber over $t$ being $\localsystem_t = {}_{p_t\tiret q_t}M_{x_t\tiret y_t}$. In this section, we describe the period map for $\localsystem$. 
 
 \medskip
 
Let $\bp$ be a base point in $B^\ast$ and denote by $\tilde \bp \in \widetilde B^\ast$ a lift of  $\bp$.  We choose an admissible symplectic basis  (see Section~\ref{ss:AdmissibleSymplecticBasis})
\begin{equation*}
a_1,\dots,a_g,b_1,\dotsc,b_g \in H_1(\rsf_\bp,\Z)
\end{equation*} 
and lift it to $\localsystem_\bp$. To facilitate the computation of the monodromy, we assume that the admissible basis is induced from the fixed spanning tree $T$ of $G$ (see Section~\ref{ss:AdmissibleSymplecticBasis}). As we recall, we enumerate the edges in $E \setminus E(T)$ by $e_1,\dots, e_\graphgenus$ and endow each with a orientation. The space of vanishing cycles $A$ is then generated by vanishing cycles  $a_1,\dots, a_h$, with $a_j = a_{e_j}$. Moreover, the vanishing cycles $a_j$ canonically lift to $\localsystem_\bp$.

\smallskip

We lift the other classes $a_{j}$ and $b_{j}$, $j=h+1,\dots ,g$ to elements of
$\localsystem_\bp$ by choosing loops that
do not meet the points $p_\bp, q_\bp, x_\bp, y_\bp$. We denote these new classes also by $a_{j}$ and $b_{j}$. Moreover, we choose an element $\gamma_{x_\bp\tiret y_\bp}$ in $\localsystem_\bp$ with boundary $x_\bp -y_\bp$ which does not contain $p_\bp, q_\bp$ and verifies $(\gamma_{x_\bp\tiret y_\bp}, a_j)= (\gamma_{x_\bp\tiret y_\bp}, b_j) =0$, $j\in[g]$.  In addition, we take small cycles $\gamma_{p_\bp}$ and $\gamma_{q_\bp}$ around $p_\bp$ and $q_{\bp}$ with the induced orientation from the surface. We have $\gamma_{p_\bp} + \gamma_{q_\bp}=0$, and $a_1, \dots, a_g, b_1, \dots, b_g, \gamma_{x_\bp\tiret y_\bp}, \gamma_{p_\bp}$ form a basis of $\localsystem_{\bp}$.

\smallskip

The local system $\widetilde \localsystem$, pullback of $\localsystem$ to the universal cover, is trivial over $\widetilde B^\ast$. We can
thus spread out the above basis of $\localsystem_\bp$ to a basis 
\begin{equation*}
a_{1,\tilde t}, \dots, a_{g,\tilde t}, b_{1,\tilde t}, \dots,
b_{g,\tilde t}, \gamma_{x_{\tilde t}\tiret y_{\tilde t}}, \gamma_{p_{\tilde t}}
\end{equation*}
of $\widetilde \localsystem$, for any $\tilde t\in \widetilde B^\ast$.
Since the space of vanishing cycles $A$ is isotropic, the cycles $a_{i,\tilde t}$ are invariant under the monodromy and only depend on the point $t \in B^\ast$. We hence denote them by $a_{i,t}$.  The loops $\gamma_{p_{\tilde t}}$ depend trivially on the point $t$. Summarizing, we get the basis 
\begin{equation*}
a_{1,t}, \dots, a_{g,t}, b_{1,\tilde t}, \dots,
b_{g,\tilde t}, \gamma_{x_{\tilde t}\tiret y_{\tilde t}}, \gamma_{p_{t}}
\end{equation*}
of $\localsystem_{\tilde t}$ for $ \tilde t \in \widetilde B^\ast$ lying above $t\in B^\ast$.

\smallskip

By the admissibility of the variation of Hodge structures, the corresponding holomorphc bundle
$\localsystem \otimes_{\C}{\mathcal O}_{B^\ast}$ can be extended to a
holomorphic vector bundle over $B$. In addition, the Hodge filtration 
$F^{0}W_{-1}\localsystem_{\C}$ can be extended to a coherent subsheaf of it.
We infer the existence of a collection of holomorphic $1$-forms
$\{\omega_{i}\}_{i=1, \ldots, 
  g}$ on $\rsf$ such that, for each $t\in
B^\ast$, the restrictions $\bigl\{\omega_{i, t}\coloneqq
\omega_{i}\rest{\rsf_{t}}\bigr\}_{i=1}^g$ 
form a  basis  of the
holomorphic differentials on $\rsf_t$ and, moreover, we have
\begin{equation} 
 \int_{a_{i, t}}
\omega_{j, t}=\delta_{i, j}. 
\end{equation}
For each point $\tilde t$, the period matrix of the Riemann surface $\rsf_{\tilde t}$ is given by
\[
\Omega_{\tilde t} = \Bigl(\int_{b_{j,\tilde t}} \omega_{i,t}\Bigr)_{i,j \in [g]}.
\]
By the same reasoning as above, we find a relative logarithmic one-form  $\omega_{p\tiret q}$ on $\rsf$ such that for each $t\in B^\ast$, the restriction $\omega_{p_t\tiret q_t}$ has residue divisor $p_t - q_t$ and verifies 

\begin{equation}
  \int_{a_{i,t}}\omega_{p_t\tiret q_t}=0,\quad i=1,\dots,g.
\end{equation}

Let $\bipm_{\tilde t}$ the period matrix associated to the biextension ${}_{p_{\tilde t} \tiret q_{\tilde t}}M_{x_{\tilde t}\tiret y_{\tilde t}}$ with respect to the basis $a_{1,t}, \dots, a_{t,g}, b_{\tilde t, 1}, \dots, b_{\tilde t,g}, \gamma_{x_{\tilde t}\tiret y_{\tilde t}}, \gamma_{p_t}$ (see Section~\ref{PeriodMatrixBiextension}). By definition,
\[\bipm_{\tilde t}  = \begin{pmatrix} \Omega_{\tilde t} & \rmW_{\tilde t} \\
\rmZ_{\tilde t} & \rho_{\tilde t}
\end{pmatrix}\]
where \begin{itemize}
\item the $g\times g$ matrix $\Omega_{\tilde t}$ is the period matrix of $\rsf_{\tilde t}$ and belongs to the Siegel domain    
\begin{equation*}
\H_g := \Bigl\{\text{ $g\times g$ complex symmetric matrix } \Omega\,\,|\
\rm{Im}(\Omega)>0\,\Bigr\},
\end{equation*}
\item $\rmW_{\tilde t} = \Bigl(\int_{b_{\tilde t,j}}\omega_{p_t \tiret q_t}\Bigr)_{j=1, \dots, g}$ is a column matrix in $\C^g$,
\item $\rmZ_{\tilde t} = \Big( \int_{\gamma_{x_{\tilde t} \tiret y_{\tilde t}}} \omega_j\Bigr)_{j=1,\dots, g}$ is a row matrix in $\C^g$, and
\item $\rho = \int_{\gamma_{x_{\tilde t}\tiret y_{\tilde t}}}\omega_{p_t\tiret q_t}$ is a scalar in $\C$.
\end{itemize}

The preceding discussion directly implies the following proposition.  

\begin{prop} The period map of the variation of mixed Hodge
  structures ${}_{x\tiret y}M_{p\tiret q}$ is given by
\begin{align}
\widetilde\Phi \colon \widetilde B^\ast &\longrightarrow  \mathbb
H_{g}\times \row_{g}(\C) \times \col_{g}(\C) \times \C \nonumber \\
\tilde{t} &\longmapsto \begin{pmatrix} \Omega_{\tilde t} & \rmW_{\tilde t}  \\
\rmZ_{\tilde t} & \rho_{\tilde t}
\end{pmatrix} \cong \Bigl( \Omega_{\tilde t}, \rmZ_{\tilde t}, \rmW_{\tilde t}, \rho_{\tilde t}\Bigr).
\end{align}
\end{prop}

The group $\mathrm{Sp}_{2g}(\R)$ acts on $\H_g$ by 
$$\begin{pmatrix}A & B \\ C & D\end{pmatrix}\Omega = (A\Omega+B)(C\Omega+D)^{-1}. 
$$
The Siegel moduli space is defined as the quotient  $\SiegAbel_g= \mathrm{Sp}_{2g}(\Z)\,\backslash\, \H_g$ and 
 parametrizes principally polarized abelian
varieties of dimension $g$.   The period map for the local system $H_1(\rsf_{t}, \Z)$ is the map 
\[\widetilde B^\ast \to \mathbb H_g, \quad \tilde t \to \Omega_{\tilde t}\] 
and the image of the monodromy map lives in $\mathrm{Sp}_{2g}(\Z)$~\cite{ABBF, AN}. 

Let $\group$ be the subgroup of $\gl_{2g+2}(\R)$ whose elements are matrices of the form 
 \begin{align*}
 \begin{pmatrix} 1 & \lambda_1 & \lambda_2 & \alpha \\0 & A & B & \mu_1 \\ 0 & C & D & \mu_2 \\ 0 & 0 & 0 & 1\end{pmatrix},\,\,\quad \lambda_1, \lambda_2 \in \text{Row}_g(\R),\ \mu_1, \mu_2 \in 
\text{Col}_g(\R), \alpha \in \C,\ \begin{pmatrix}A & B \\ C & D\end{pmatrix}\in \mathrm{Sp}_{2g}(\R),
 \end{align*}
and denote by $\group(\Z)$ the subgroup of $\group$ consisting of those matrices with integer coordinates. 

The group $\group$ naturally acts on the space $\mathbb
H_{g}\times \mathrm{Col}_{g}(\C) \times \mathrm{Row}_{g}(\C) \times \C$ and turns it into a homogeneous space.
The image of the monodromy map for the variation of biextension mixed Hodge structures $M$ is in $\group(
\Z)$. This follows from the explicit calculation of the monodromy given in the next section (see also\cite{ABBF}).

\subsection{Logarithmic monodromy} \label{ss:MarkedMonodromy} In this section we describe the action of the logarithm of the monodromy $N_e$ on the coordinates of the period matrix $\bipm_{\bp}$. 

\smallskip

Recall that our choice $a_1, \dots, a_{\graphgenus}$ of the basis of vanishing cycles was induced by the edges of the complement of the spanning tree $T$. For each oriented edge $e\in T$,  the graph $T-e$ obtained by deletion of $e$ has two connected components with vertex sets $W$ and $W^c=V\setminus W$ such that $e$ has tail in $W$ and head in $W^c$. Hence we can expand $a_e$ in the basis $a_1, \dots, a_g, b_1, \dots, b_g, \gamma\indbis{x\tiret y}, \gamma\indbis p$ of $\localsystem_\bp$ as
\begin{align}\label{eq:expansion-vanishing-cycle}
a_e = -\sum_{\substack{f \in \E(W, V\setminus W)\\ f\neq e}}a_f + \epsilon\indbis{e, p} \gamma\indbis{p} + \epsilon\indbis{x\tiret y} \gamma\indbis{x\tiret y}.
\end{align}
Since $a_e$ has zero intersection with all the elements in $A$, it follows that the coefficient $\epsilon\indbis{x\tiret y}$ of $\gamma_{x\tiret y}$ is zero. Note that we have $f = e_i$ or $\bar e_i$, for $i\in [h]$, so $a_f = \pm a_i$.

\smallskip

By the Picard-Lefschetz theorem, we get the following description of the monodromy.

\smallskip

For any point $z$ on the smooth part of $S_0$, denote by $\marking(z)$ the vertex of $G$ whose component $C_{\marking(z)}$ contains $r$. Let $\varrho_{x\tiret y}$ be the unique oriented path in the tree $T$ which goes from $\marking(y_0)$ to $\marking(x_0)$. Similarly, define $\varrho_{p\tiret q}$. Consider the cycles $\gamma_1, \dots, \gamma_h$ in $H_1(G, \Z)$ associated to the elements $b_{1,\bp}, \dots, b_{h,\bp}$. Recall that $\gamma_i$ corresponds to the unique oriented cycle in the graph obtained by adding the edge $e_i$ to the tree $T$, and which has orientation compatible with that of $e_i$, $i\in[h]$ (see Section~\ref{ss:AdmissibleSymplecticBasis}).

\begin{prop}
For all edges in the tree $T$, and all $j \in [\graphgenus]$,
\begin{align}
  N_{e}(b_{j, \tilde \bp}) &= \gamma_j(e) a_e \\
  N_{e}(\gamma_{x_{\tilde \bp} -y_{\tilde \bp} })
  &= \varrho\indbis{x\tiret y}(e)a_e.
\end{align}
Moreover, for every edge $e=e_i$ not in the tree $T$ and all $j \in [\graphgenus]$,
\begin{align}
  N_{e_i}(b_{j, \tilde \bp}) &= \delta_{i,j} a_i \\
  N_{e_i}(\gamma_{x_{\tilde \bp} -y_{\tilde \bp} })
  &= 0.
\end{align}
\end{prop}

\begin{proof} This follows from the Picard-Lefschetz equations~\eqref{eq:PL}.
\end{proof}

The forms $\omega _{j}$, $j\in[g]$, and $\omega_{p\tiret q}$ are defined
globally, and they are invariant under monodromy. By definition, they satisfy the equations 
\begin{align}
  &\int_{a_{e}}\omega _{j}=\gamma_j(e), &&\int_{a_{e}}\omega_{p_{\bp}\tiret q_\bp} =\varrho\indbis{p\tiret q}(e)
\end{align}
for all vanishing cycles $a_e$, $e \in E$.

We deduce the following result. For an edge $e \in E$, let $\innone{e}{\cdot\,, \cdot}$ be the bilinear pairing on $C_1(G, \R)$ mapping two one-chains $\varrho_1, \varrho_2$ to  $\innone{e}{\varrho_1, \varrho_2} := \varrho_1(e) \varrho_2(e)$.

\begin{prop} For all edges $e \in E$ and all $i,j \in [\graphgenus]$,
\begin{align*}
  N_e\bigl(\int_{b_{i,\tilde \bp}}\omega_{j,\bp}\bigr)
  &= \innone{e}{\gamma_i, \gamma_j}\\
  N_e(\int_{\gamma_{x_{\tilde \bp}\tiret y_{\tilde \bp}}}\omega_{j,\bp})  
  &= \innone{e}{\gamma_j ,\varrho\indbis{x\tiret y}}\\ 
  N_e(\int_{b_{i, \tilde \bp}}\omega_{p_\bp\tiret q_\bp})
  &= \innone{e}{\gamma_i, \varrho\indbis{p\tiret q}} \\ 
  N_e(\int_{\gamma_{x_{\tilde \bp}\tiret y_{\tilde \bp}}} \omega_{p_\bp\tiret q_\bp})
  &= \innone{e}{\varrho\indbis{x\tiret y}, \varrho\indbis{p\tiret q}}.
  \end{align*}
\end{prop}
\smallskip

We deduce from the above calculations that the logarithm of the monodromy of the element $\lambda_e$ is the element $N_e$ of the Lie
algebra of $\group$ given by 
\begin{equation}
N_{e}=
  \begin{pmatrix} 0 & 0 & \widetilde \rmW_{e} &
    \rho_e  \\0 &
    0 & \widetilde \rmM_{e} & \widetilde \rmZ_{e} \\
    0 & 0 & 0 & 0 \\ 0 & 0 & 0 & 0\end{pmatrix}, 
\end{equation}
where
\begin{itemize}
\item $\widetilde \rmM_{e}$ is the
$\graphgenus \times \graphgenus$ matrix $\rmM_{e} =\Bigl( \innone{e}{\gamma_i,\gamma_j}\Bigr)_{i,j
\in[\graphgenus]}$, filled to a $g \times g$ matrix by adding zeros;
\item $\widetilde \rmW_e$ is the row matrix  $\rmW_e = \bigl( \innone{e}{\gamma_j, \varrho\indbis{x\tiret y}}\bigr)_{j\in [\graphgenus]} $ of size $\graphgenus$, filled to a row matrix of size $g$ by adding zeros;
\item $\widetilde \rmZ_e$ is the column matrix  $\rmZ_e = \bigl( \innone{e}{\gamma_i, \varrho\indbis{p\tiret q}}\bigr)_{i \in[\graphgenus]}$ of size $\graphgenus$, filled to a column matrix of size $g$ by adding zeros;
\item and $\rho_e$ is the scalar $\rho_e= \innone{e}{\varrho\indbis{x\tiret y}, \varrho\indbis{p\tiret q}}$. 
\end{itemize}
Note in particular that all entries are integers.

\smallskip

It follows that the image of the monodromy lives in $\group(\Z)$ and the
map $\widetilde{\Phi}$ descends to $B^\ast$ by taking the quotient, resulting in  the following commutative diagram.
\begin{equation}
\label{diag-phi}
\begin{CD} 
\widetilde{B^\ast} @>\widetilde\Phi>> 
\H_{g}\times \row_{g}(\C) \times \col_{g}(\C) \times \C  \\
@VVV @VVV \\
B^\ast @>\Phi>> 
\group(\Z)\backslash \Bigl(\H_{g} \times \row_{g}(\C) \times \col_{g}(\C) \times \C  \Bigr).
\end{CD}
\end{equation}

\subsection{The imaginary part of the period matrix and the height pairing} \label{ss:ImaginaryPartHeightPairing} As we discuss next, the imaginary part of the  period matrix (taken componentwise) is closely related to the height pairing $\hp{\rsf_t}{p_t - q_t, x_t - y_t}$.

 \smallskip
 
By the above description of the monodromy action, the imaginary part $\Im(\bipm_{\tilde t})$ is invariant under the monodromy. This means that $\Im(\bipm_{\tilde t})$ descends to a map from $B^\ast$ to the space of $(g+1)\times (g+1)$ real matrices $Q$ which satisfy that the submatrix $(Q(i,j))_{i,j \in [g]}$ is symmetric and positive definite. For a point $t\in B^\ast$, we denote the corresponding matrix by
 \[
 \Im(\bipm_t) = \begin{pmatrix} \Im( \Omega_t) & \Im(  \rmW_t ) \\
\Im(\rmZ_t)  & \Im(\rho_t) \end{pmatrix}.
 \] 

By Proposition~\ref{prop:HeightPairingMatrix}, the height pairing can be expressed as follows.
\begin{lem} \label{lem:HeightPairingByPeriodMap}
In terms of the imaginary part of the period matrix, the height pairing equals
\[
\hp{\rsf_t}{p_t - q_t, x_t - y_t} = 2\pi \big (  \Im(\rho_t) -  \Im(\rmZ_t) \Im(\Omega_t)^{-1} \Im(\rmW_t) \big ) , \qquad t \in B^\ast.
\]
\end{lem}

\subsection{Asymptotics of the holomorphic forms $\omega_{j,t}$} \label{ss:AsymptoticsHolomorphicForms} We briefly discuss the structure of the holomorphic one-forms $\omega_j$, $j \in [g]$. The following description was already used in~\cite{AN}. 

\smallskip

Notations as before, let $\pr \colon \rsf \to B$ be the family of stable Riemann surfaces. Denote by $\pr^*D$ the simple normal crossing divisor in $\rsf$ defined by the pullback of divisors $\pr^{-1}(D_e)$ for $e\in E$. We consider the sheaf $\omega_{\rsf/B}\bigl(\log(\pr^*D)\bigr)$ of holomorphic forms with logarithmic singularities along the divisor $\pr^*D$. The holomorphic forms $\omega_1, \dots, \omega_g$, defined in terms of the symplectic basis, form a basis of the space of global sections $H^0\Bigl(B, \pr_*\omega_{\rsf/B}\bigl(\log(\pr^*D)\bigr)\Bigr)$. The restriction of $\omega_j$ to the fiber $\rsf_t$ is denoted by $\omega_{j,t}$.

\smallskip
The limit forms $\omega_{j, 0}$, $j \in [g]$ on $S_0$ are global sections of the twisted sheaf $\omega_{S_0}(p_0+q_0)$. They can be described as follows. Recall that our admissible basis was defined in terms of the oriented edges $e_1, \dots, e_h$ of $G$, which are not in the fixed spanning tree $T$. Recall as well that $\gamma_1, \dots, \gamma_\graphgenus$ is the corresponding basis of $H_1(G, \Z)$. That is, $\gamma_j$ is the unique oriented cycle in the graph obtained by adding the oriented edge $e_j$ to the tree.

\smallskip

For each graph cycle $\gamma_j$, $j\in [h]$ and Riemann surface component $C_v$, $v \in V$ of $S$, consider the degree zero divisor $D_{j,v} \in \Div^0(C_v)$ given by
\begin{equation} \label{eq:DivisorFromCycle}
D_{j,v} =  \sum_{e = vv'\in \E} \gamma_j(e) p^e_v.
\end{equation}
Let $\omega_{j,v}$ be the meromorphic one-form on $C_v$ associated to $D_{j,v}$. That is, $\omega_{j,v}$ has residue divisor $D_{j,v}$ and $\int_{a_i} \omega_{j,v} = 0$ all $a$-cycles $a_i$ in the fixed symplectic basis on $C_v$.

\begin{prop} \label{prop:AsymptoticsHolomorphicForms}
The limit forms $\omega_{j, 0} := \omega_j \rest{S_0}$ are given as follows.
\begin{itemize}
\item [(i)] For $j \in [\graphgenus]$, the form $\omega_{j,0}$ is the unique section of $\omega_{S_0}$ having the above restrictions $\omega_{j, v}$ to the components $C_v$. 

\item [(ii)] For $r \in [n]$, the restrictions of the forms $\omega_{j, 0}$, $j \in I_r$, to $C_{v_r}$ form the basis of $\Omega^1(C_{v_r})$ corresponding to the symplectic basis of $H_1(C_{v_r}, \Z)$ given by the cycles $a_i$, $b_i $ for $ i \in I_{r}$. Their restrictions to any other component $C_v \neq C_{v_r}$ are identically zero.
\end{itemize}
\end{prop}

\begin{proof} To prove the first claim on $\omega_j$, $j \in [h]$, note first that on every component $C_v$,
\[\int_{a_i} \omega_{j,0} = 0 \]
for all $a$-cycles $a_i$ in the fixed symplectic basis of $H_1(C_v, \Z)$. Suppose that the vertex $v$ lies on the cycle $\gamma_j$. Then there exist two oriented edges $e'$ and $e''$ in $\gamma_j$ with tail $\tail(e') =v$ and head $\head(e'') =v$, and  $D_{j,v} = p^{e'}_v - p_{v}^{e''}$.  We claim that for all $t \in B^\ast$,
\[
\int_{a_{e'}, t} \omega_{j,t} = 1.
\]
If $e'$ belongs to the spanning tree $T$, then $e' = e_j$ and the claim is trivial.  Otherwise,
\[
a_{e'} = -\sum_{\substack{e \in \E(W, V\setminus W)\\ e\neq e'}}a_e + \epsilon_{e', p} \gamma_p + \epsilon\indbis{x\tiret y} \gamma\indbis{x\tiret y}
\]
by Equation~\eqref{eq:expansion-vanishing-cycle}.  However, then the $j$-th edge $e_j$ of $E \setminus E(T)$ belongs to $\E(W, V\setminus W)$, which implies the desired property. Analogously, we obtain that $\int_{a_{e''}, t} \omega_{j,t} = -1$ for all $t \in B^\ast$. By continuity, this proves that $\omega_{j,v} \rest{C_v} = \omega_{j,v}$  for all $v \in V$ lying on the cycle $\gamma_j$. Similarly, one shows that $\omega_{j,0} \equiv 0$ on all other components $C_v$.

The second claim can be deduced easily from the equality
\[
\int_{a_j, t} \omega_i =\delta_{i,j},
\]
which holds for all $i, j \in [g]$ and all $t \in B^\ast$.

\end{proof}

\subsection{Asymptotics of the meromorphic forms $\omega_{p_t \tiret q_t}$} In a similar way, one can describe the limits of the meromorphic forms $\omega_{p_t \tiret q_t}$.

 \smallskip
 
 Recall that, by spreading the first part of symplectic basis $a_1, \dots, a_g$ over $B^\ast$, we had obtained a holomorphic section $\omega_{p \tiret q}$ of the sheaf $\omega_{\rsf^\ast/B^\ast}(p+q)$. Here, by an abuse of the notation, $p$ and $q$ denote the restrictions $p_{|B^\ast}$ and $q_{|B^\ast}$, respectively.  On each fiber $\rsf_t$, for $t\in B^\ast$, $\omega_{p \tiret q}$ has restriction equal to the logarithmic one-form $\omega_{p_t\tiret q_t}$ with residue divisor $p_t-q_t$ and vanishing periods with respect to $a_{1,t}, \dots, a_{g,t}$ (see Section~\ref{ss:PeriodMap}).

\begin{prop} Notations as above, $\omega_{p\tiret q}$ extends holomorphically to a global section of the relative sheaf $\omega_{\rsf/B}(p+q)$.  
\end{prop} 

 We denote this section again by $\omega_{p \tiret q}$ and its restriction to $\rsf_t$, $t \in B$, by $\omega_{p_t\tiret q_t}$. The limit logarithmic form $\omega_{p_0\tiret q_0}$ on $S_0$ is a global section of $\omega_{S_0}(p_0+q_0)$. Since our admissible basis was induced by the spanning tree $T$ of $G = (V,E)$, we get the following nice description.

\smallskip

Recall that $\varrho\indbis{p \tiret q}$ denotes the unique oriented path from $\bar q_0$ to $\bar p_0$ in the spanning tree $T$. Here, for a point $x$ lying on the smooth part of $S_0$, $\bar x$ is the vertex $v$ of $G$ with $x\in C_v$. 

\smallskip

On each component $C_v,$, $v \in V$, consider the degree zero divisor $D_{p_0 \tiret q_0,v}$ given by
\[
D_{p_0 \tiret q_0 ,v} = \delta_{v, \bar p_0}  p_0 - \delta_{v,\bar q_0} q_0 + \sum_{e = vv'\in \E} \varrho_{x\tiret y}(e)p^e_v.
\]

Let $\omega_{p_0 \tiret q_0, v}$ be the unique section of the twisted canonical sheaf $\omega_{C_v}(D_{p\tiret q,v})$ with residue divisor $D_{p_0 \tiret q_0,v}$ and vanishing periods on all the elements $a_j$ in the fixed symplectic basis of $H_{1}(C_v, \Z)$. For any vertex $v$ not appearing on the path $\varrho_{p\tiret q}$, we have $\omega_{p_0 \tiret q_0, v} \equiv 0$. We have the following proposition. 

\begin{prop}
The restriction  $\omega_{p_0\tiret q_0} := \omega_{p \tiret q }\rest{S_0}$ is the unique section of $\omega_{S_0}(p_0+q_0)$ having the above restrictions $\omega_{p_0 \tiret q_0, v}$ to the components $C_v$.  
\end{prop}

\begin{proof}
By definition of $\omega_{p\tiret q}$ on $B^\ast$, we get for any edge $e$ not belonging to the spanning tree
\[\int_{a_{t,e}} \omega_{p_t\tiret q_t}  =0, \qquad t \in B^\ast.\]
By continuity, we infer that the restriction of $\omega_{p_0\tiret q_0}$ to $C_v$ and $C_u$ is holomorphic at the points $p^e_v$ and $p^e_v$, respectively, where $u$ and $v$ are the two extremities of the edge $e$.

On each surface $\rsf_t$,  we choose a small cycle $\gamma_{p_t}$ and $\gamma_{q_t}$ around $p_t$ and $q_t$, respectively, with orientation compatible with that of the surface.

 Now we use Equation~\eqref{eq:vanishing_cycle3} for an oriented edge $e$ of the spanning tree. Let $W$ and $W^c=V\setminus W$ be the vertex set of the two connected components of $T \setminus e$ with the tail of $e$ in $W$. Then we get
 \begin{align*}\int_{a_{e,t}} \omega_{p_t \tiret q_t} &= - \sum_{\substack{e_i \in \E(W, W^c) \\ i\in[\graphgenus]}} \int_{a_{t, e_i}} \omega_{p_t\tiret q_t}- \chi\ind{W}\bigl(\bar p_0\bigr) \int_{\gamma_{p_t}} \omega_{p_t\tiret q_t} - \chi\ind{W}\bigl(\bar q_0\bigr) \int_{\gamma_{q_t}}\omega_{p_t\tiret q_t}\\
 &= - \chi\ind{W}\bigl(\bar p_0\bigr) +  \chi\ind{W}\bigl(\bar q_0\bigr) =- \varrho\indbis{p\tiret q}(e).
 \end{align*}
 Here, $\chi\ind{W}\colon V \to \{0,1\}$ is the characteristic function of $W$. The claim now easily follows.
 \end{proof}


\section{Asymptotics of the period map} \label{sec:period-map-asymptotics}
 
In this section, we continue our study of the period map associated with a family of Riemann surfaces $\rsf \to B$ with four markings $p, q,x,y$. Using the nilpotent orbit theorem ~\cite{KNU, Pearl, Pearl2}, we describe the asymptotics of the period map close to zero. In particular, considering the imaginary parts, we obtain an approximation of the form
\[
\Im(\bipm_t) = \widehat \Ap(t) + o(1) = \widehat \Apc(t) + \~\rmM_\ell(t)  + o(1), \qquad \text{as $t \in B^\ast$ converges to $0$},
\]
where $ \Ap(t)$ is a certain approximative matrix with \emph{complex part} $\widehat \Apc(t)$ and \emph{graph part} $\rmM_\ell(t)$. This approximation result is the main ingredient in the proofs of Theorem~\ref{thm:tameness-height-pairing1} and Theorem~\ref{thm:FinalTamenessHeightPairing} (see Section~\ref{sec:ProofHeightPairing} and Section~\ref{sec:HybridJFctAsymptotics}).

\smallskip

Since this section is a direct continuation of the previous one, we use the same notations and notions, unless otherwise explicitly stated.

\subsection{Asymptotics of the period map}
\label{sec:asympt-peri-map}

In this section, we apply the nilpotent orbit theorem for admissible variations of mixed Hodge structures~\cite{KNU, Pearl, Pearl2} to describe the asymptotics of the period map.

\smallskip

 We will adapt the notations introduced previously and write the coordinates corresponding to the edge set as $
 z_E$ and the other ones by $z_{E^c}$. The coordinates of a point $t$ in $B^\ast$ become $
z_E \times z_{E^c}$. Abusing the notation, we use $z_e$ both for edges $e\in E$ and for the other non-edge coordinates, that is for $e\in E^c$. Whether we are considering edges or non-edges will be clear from the context. We denote the coordinates in the
 universal cover $\widetilde{B^\ast}$ by $\zeta_{e}$. The cover map  $\widetilde{B^\ast}\to B^\ast$ in these coordinates is given by 
 \begin{equation}
   \label{eq:16}
   z_{e}=
   \begin{cases}
     \exp(2\pi i \zeta_{e}),&\text{ for }e\in E,\\
     \zeta_{e},&\text{ for }e\not \in E.
   \end{cases}
 \end{equation}

\medskip

Consider the \emph{twisted period map}
\begin{equation}
  \label{eq:18}
\widetilde\Psi(\zeta)=\exp(-\sum_E \zeta_eN_e)\widetilde\Phi(\zeta)  
\end{equation}
which takes values in the compact dual $\widecheck{\pdomain}$ of the period domain $\pdomain = \sH_g \times \row(\C^g) \times \col(\C^g)\times \C$. The space $\widecheck{\pdomain}$ contains $\pdomain$ as
an open subset and the action of the group $\group(\Z)$ can be extended to $\widecheck \pdomain$. 

Since $\widetilde \Psi$ is invariant  under the transformation
$\zeta_e \mapsto \zeta_e+1$, it descends to a map $\Psi: B^\ast \to
\widecheck{\pdomain}$.

\smallskip

The following result is the Nilpotent Orbit Theorem that
we need. Recall that
$\Delta $ is a disk of small radius and we denote $B=\Delta^{\dimms}$, with $N = 3g-3 + n$.

\begin{thm}[Nilpotent orbit theorem] \label{thm:nilpotentorbit}
Shrinking the radius of $\Delta $ if necessary, the map $\Psi$ extends
to a holomorphic map
\begin{displaymath}
  \Psi: B \longrightarrow \widecheck \pdomain. 
\end{displaymath}
Moreover, there exists a constant $\constant_{0}$ such that for all  $t\in B$ , we have 
\begin{displaymath}
  \exp(\sum_{e\in E} \zeta_eN_e)\Psi(t)\in \pdomain
\end{displaymath}
provided that we have $\Im(\zeta_{e})\ge \constant_{0}$ for any $e\in
E$.

\smallskip

In addition, there are constants $C, \beta >0$ for which the following estimation holds. For any $t \in B$, there exists a small open disc $U_t \subseteq B$ containing $t$ so that 
\[\mathrm{dist}\Bigl(\widetilde \Phi(\tilde s) ,  \exp(\sum_{e\in E} \zeta_e(\tilde s)N_e)\Psi(t)\Bigr) \leq C \sum_{e\in E} \Im(\zeta_e(\tilde s))^\beta \exp(-2\pi \Im(\zeta_e(\tilde s)))\]

for any $s\in U_t\setminus D$ and $\tilde s$ above $s$ in $\widetilde B^\ast$ provided that $\Im(\zeta_e(\tilde s)) \geq \constant_0$ for all $e\in E$.
\end{thm}

\begin{proof} A proof of this theorem for admissible variations of mixed Hodge structures can be found in~\cite{KNU, Pearl, Pearl2}. \end{proof}

We consider now the period map \begin{align*}
\widetilde\Phi \colon \widetilde{B^\ast} &\longrightarrow \pdomain = \sH_g \times \row(\C^g) \times \col(\C^g) \times \C \\
\tilde{t} &\longmapsto \bipm_{\tilde t} = \begin{pmatrix} \Omega_{\tilde t} & \rmW_{\tilde t} \\
\rmZ_{\tilde t} & \rho_{\tilde t}.
\end{pmatrix}
\end{align*}

\begin{equation}
\ExtM_e:=
  \begin{pmatrix} \widetilde \rmM_e & \widetilde \rmW_{e}\\
   \widetilde \rmZ_e & \rho_e
   \end{pmatrix}, 
\end{equation}
where we use the notations from Section~\ref{ss:MarkedMonodromy}.

\begin{prop} Notations as above, we get
\[\widetilde \Psi(\tilde t) = \bipm_{\tilde t} - \sum_{e\in E} \zeta_e(\tilde t) \ExtM_e.\] 
\end{prop}
\begin{proof} This can be obtained from the form of the monodromy matrix $N_e$ and the action of the group $\group$ on $\pdomain$.
\end{proof}
It follows  that we can write $ \widetilde \Psi(\tilde t) =: \ELambda_t \in \pdomain$ for any point $t\in B^\ast$ and $\tilde t$ a point of $\widetilde B^\ast$ above $t$. Moreover, the family of matrices $\ELambda_t$ can be extended holomorphically to a family  over $B$ by $\ELambda_t \coloneqq \Psi(t)$ for the extension $\Psi\colon B \to \widecheck \pdomain$.  In addition, for any point $t\in B$, we get 
\[ \exp(\sum_{e\in E} \zeta_eN_e)\Psi(t) = \ELambda_t + \sum_{e\in E} \zeta_t\ExtM_e \in \pdomain
\]
 provided that $\Im(\zeta_e) \geq \constant_0$ for all $e\in E$. Moreover, for constants $C, \beta$ as in the statement of the theorem, for any point $t\in B$,  we obtain the estimate 
\[\mathrm{dist} \Bigl(\Omega_{\tilde s}, \ELambda_s + \sum_{e\in E} \zeta_e(\tilde s) \ExtM_e\Bigr) \leq C \sum_{e\in E} \Im(\zeta_e(\tilde s))^\beta \exp(-2\pi \Im(\zeta_e(\tilde s)))\]
 for any point $\tilde s \in \widetilde B^\ast$ living above a point $s$ of $U_t \setminus D$, for a small neighborhood $U_t$ of $t$ in $B$. 
  
 \medskip

In applications to the height pairing, we are interested in the imaginary part of the period matrix. As discussed in Section~\ref{ss:ImaginaryPartHeightPairing}, the latter is monodromy invariant and defines a function $t \mapsto \Im(\bipm_t)$ on $B^\ast$. The above estimate has the following immediate consequence.

\begin{thm} \label{thm:nilpotentorbit2} For any $t\in B$, there exists a neighborhood $U_t$ of $t$ such that 
\[\mathrm{dist} \Bigl(\Im(\bipm_s), \Im(\ELambda_t) + \frac{1}{2\pi} \sum_{e\in E} \ell_e(s) \ExtM_e\Bigr) \leq C \sum_{e\in E} \ell_e(s)^\beta \exp(-2\pi \ell_e(s))\]
for all points $s\in U_t \setminus D$ with  $\ell_e(s) := -\log|z_e(s)|$ for all $e\in E$. 
\end{thm}

\subsection{The limit period matrix $\ELambda_0$}

The limit period matrix $\ELambda_0$ depends on the choice of parameters $z_e$, and is well-defined only up to a sum of the form $\sum_{e \in E} \lambda_e\ExtM_e$. In the following discussion, we fix an adapted system of parameters on the base $B$ and on the family $\rsf \to B$. 

\smallskip

 First, we write 
\[\ELambda_0 = \begin{pmatrix} \Lambda_0 & \rmW_0 \\ 
\rmZ_0 & \rho_0\end{pmatrix}\]
and note that $\Lambda_0$ is the limit of the twisted period matrix $\Omega_{\tilde t} - \sum_{e \in E} \zeta_e \widetilde \rmM_e$.

\smallskip

In what follows, we describe the limit matrix $\ELambda_0$. We begin with the submatrix $\Lambda_0$ and refine the results from our previous work~\cite{AN}.
\begin{prop} \label{prop:LPM1} The limit matrix $\Lambda_0$ verifies the following properties. 
\begin{itemize}
\item [(i)] Let $1 \le r \le n$. Then the $\genusfunction(v_r)\times \genusfunction(v_r)$ matrix $\Lambda_0[I_r, I_r]$ coincides with the period matrix $\Omega_{v_r}$ of the component $C_{v_r}$ with respect to the symplectic basis $a_i$, $b_i$ for  $i \in I_r$ and holomorphic forms $\omega_{j,0}$, $j \in I_r$, restricted to the component $C_{v_r}$.
\item  [(ii)]  All the matrices  $\Lambda_0[I_r, I_{r'}]$ for non-zero values $1 \le r, r' \le n$ with $r'\neq r$  are vanishing. 
\end{itemize}
 \end{prop}
  (For a $g\times g$ matrix $P$, recall $P[I_r, I_{r'}]$ is the matrix with rows in $I_r$ and columns in $I_{r'}$.)

 \begin{proof} This follows immediately from the discussion in Section~\ref{ss:AsymptoticsHolomorphicForms}, Proposition~\ref{prop:AsymptoticsHolomorphicForms} and the fact that $t \mapsto \int_{b_i, t} \omega_{j, t}$ defines a holomorphic function on $B$ for all $i,j \in [g] \setminus [\graphgenus]$.
 \end{proof}
It follows that the matrix $\Lambda_0$ has the following form
 \begin{equation} \label{eq:Lambda_0}
 \Lambda_0 = \left(
 \begin{matrix}
\Omega_G & B_{v_1} & B_{v_2} & B_{v_3} & \cdots & B_{v_n}\\
B_{v_1}^\transpose & \Omega_{v_1} & 0 & 0 & \cdots & 0\\
B_{v_2}^\transpose & 0 & \Omega_{v_2} & 0 & \ddots & 0\\ 
\vdots & \vdots & \vdots & \ddots & \ddots &0 \\
B^\transpose_{v_{n-1}} & 0 & \cdots & \cdots & \Omega_{v_{n-1}} & 0 \\
B_{v_n}^\transpose & 0 & \cdots & \cdots & 0 & \Omega_{v_n}
\end{matrix} \right)
\end{equation}
for a $\graphgenus\times\graphgenus$ matrix $\Omega_G$, and $\graphgenus \times \genusfunction(v_k)$ matrices $B_{v_k}$, $k=1, \dots, n$.

\smallskip

The matrices $B_{v_k}$, $k = 1,\dots, n$ admit the following description. For each vertex $v \in V$ and each $j\in [\graphgenus]$, consider the divisor $D_{j,v}$ defined in \eqref{eq:DivisorFromCycle}. Let $\gamma_{j,v}$ be the unique (homology class of) one-chain in $C_v$,  which has boundary $D_{j,v}$ and verifies 
\[\langle a_i , \gamma_{j,v} \rangle = \langle b_i , \gamma_{j,v} \rangle  =0\]
for all cycles $a_i$, $b_i$ in the fixed symplectic basis of $C_v$.  Denote by $\omega_{i,k}$, $i\in I_k$, the collection of holomorphic forms dual to the basis $a_i, b_i$, $i\in I_k$. Recall that $\omega_{i,k} = \omega_{i,0}\rest{C_{v_k}}$ (see Section~\ref{ss:AsymptoticsHolomorphicForms}).

\begin{prop} \label{prop:LPM2}  For each $k=1, \dots, n$, write  $B_{v_k} = \Bigl (B_{v_k}(i,j) \Bigr)_{i\in[\graphgenus], \,  j \in I_k}$. Then
\[
B_{v_k}(i,j)  =\int_{\gamma_{i,v_k}} \omega_{j,k}.
\]  
\end{prop}

\begin{proof} The coordinate $B_{v_k}(i,j)$ is monodromy invariant and coincides with the limit
\[
B_{v_k}(i,j)  = \lim_{t\to 0} \int_{b_{i,t}} \omega_{j,t}.
\]
Since for $i\in I_k$, the limit form $\omega_{i,0}$ has restriction zero on $C_u$ for all the vertices $u \neq v_k$,
\[
B_{v_k}(i,j)   = \int_{\gamma_{i,v_k}}\omega_{j,0}\rest{C_{v_k}}
\]
which gives the result.
\end{proof}

In Section~\ref{ss:DecsrApImPM}, we continue to investigate the limit period matrix  $\Lambda_0$ and give a precise description of its imaginary part. In the next Section~\ref{ss:Ap}, we first collect necessary notions and definitions.

\subsection{The approximative period matrix} \label{ss:Ap} In this section, we introduce a certain \emph{approximative matrix} $\whAp$, which provides an approximation for the imaginary part of the period matrix $\Im(\bipm_t)$, $t \in B^\ast$, and turns out to be quite useful in the proof of our main results.

\smallskip

Let $(S, r_1, \dots r_\nmark)$ be a marked stable Riemann surface of genus $g$. Denote by $G = (V,E)$ the underlying graph of genus $\graphgenus$, and let $\mgr$ be a metric graph arising from an edge length function $l \colon E \to (0, + \infty)$.  Let $p,q,x,y$ be four distinct, non-nodal points on $S$. Suppose that for each node $p^e$ associated to an edge $e \in E$, we have fixed local holomorphic parameters $z^e_u$, $z^e_v$ around the point $p^e$, viewed in the components $C_u$ and $C_v$, respectively.

\smallskip

Fix a symplectic basis $a_j$, $b_j$, $j \in I_{k}$ on each of the components $C_{v_k}$ corresponding to a vertex $v_k$ of $V$. Here, we have enumerated the vertices of $G$ as $v_1, \dots v_n$ and used the index sets
\[ I_k := I_{v_k} := \big \{h + \genusfunction(v_1) + \dots + \genusfunction(v_{k-1}) + 1, \dots,  h + \genusfunction(v_1) + \dots + \genusfunction(v_k) \big \}, \qquad k=1, \dots, n.
\]
Let $T= (V, E(T))$ be a spanning tree of the graph $G=(V,E)$. 
In our context, the above data will be encoded in the chosen admissible symplectic basis as defined in Section~\ref{ss:AdmissibleSymplecticBasis}, which we assume to be induced from a fixed spanning tree $T$ of $G$.

We denote by $\omega_i$, $i \in I_{k}$, the associated basis of the space of holomorphic differentials on $C_{v_k}$. 

\smallskip

Let $\gamma_1, \dots, \gamma_h$ be the basis of $H_1(G, \Z)$ obtained from $T$. That is, we enumerate the edges in $E \setminus E(T)$ as $e_1, \dots, e_h$ and define $\gamma_i$ as the unique cycle in $G$ created by adding $e_i$ to the spanning tree $T$.

\smallskip
To the above choices, we associate the following matrix $\whAp(S, \mgr; p,q,x,y) \in \R^{(g+1) \times(g +1)}$. First of all, we need to introduce the following relative quantities.

\smallskip
\begin{itemize}
\item For each graph cycle $\gamma_j$, $j\in [h]$, and Riemann surface component $C_v$, $v \in V$ of $S$, consider the degree zero divisor $D_{j,v} \in \Div^0(C_v)$
\[
D_{j,v} =  \sum_{e = vv'\in \E} \gamma_j(e) p^e_v.
\]
\smallskip

\item Let $\omega_{j,v}$ be the meromorphic one-form on $C_v$ associated to $D_{j,v}$. That is, $\omega_{j,v}$ has residue divisor $D_{j,v}$ and $\int_{a_i} \omega_{j,v} = 0$ for all cycles $a_i$ in the symplectic basis of $C_v$.
\smallskip

\item Let $\gamma_{j, v}$, $v \in V$, be the unique (homology class of a) one-chain on the component $C_v$ with boundary divisor $D_{r \tiret s, v}$ and
\[
\langle \gamma_{r \tiret s, v}, a_j \rangle = \langle \gamma_{r \tiret s, v}, b_j \rangle = 0
\]
for all cycles $a_j$, $b_j$ in our fixed symplectic basis on $C_v$.
\end{itemize}
For two points $r,s$ lying on the smooth part of  $S$...
\begin{itemize}
\item $\marking(r)$ and $\marking(s)$ denote the vertices in $G$ whose components $C_{\marking(r)}$ and $C_{\marking(s)}$ accommodate $r$ and $s$, respectively.
\smallskip

\item $\varrho_{r \tiret s}$ denotes the unique (oriented) path from $\marking(r)$ to $\marking(s)$ in the spanning tree $T$.
\smallskip

\item  $D_{r\tiret s,v}$, $v \in V$, is the degree zero divisor on the component $C_v$ induced by
\begin{equation} \label{eq:DivisorsFromPaths}
D_{r\tiret s,v} = \delta_{v,\marking(r)}  r - \delta_{v,\marking(s)} s + \sum_{e = vv'\in \E} \varrho_{r\tiret s}(e) p^e_v.
\end{equation}
\item $\omega_{r \tiret s, v}$, $v \in V$, is the meromorphic one-form on the Riemann surface component $C_v$ associated with $D_{r \tiret s, v}$. That is, $\omega_{r \tiret s, v}$ has residue divisor $D_{r \tiret s,v}$ and $\int_{a_j} \omega_{r \tiret s,v} = 0$ for all cycles $a_i$ in the fixed symplectic basis for $C_v$.
\smallskip

\item $\gamma_{r \tiret s, v}$, $v \in V$, is the unique (homology class of a) one-chain on the component $C_v$ with boundary divisor $D_{r \tiret s, v}$ and
\[
\langle \gamma_{r \tiret s, v}, a_j \rangle = \langle \gamma_{r \tiret s, v}, b_j \rangle = 0
\]
for all cycles $a_j$, $b_j$ in our fixed symplectic basis on $C_v$.
\end{itemize} 

\subsubsection{Definition} The \emph{approximative matrix} $\whAp \in \R^{(g+1) \times (g+1)}$ is defined via the following decomposition 
\begin{equation} \label{eq:DefinewhAp}
\whAp (S, \mgr; p,q,x,y) := \whApc + (2\pi)^{-1} \widehat \rmM_l,
\end{equation}
for the following \emph{complex part} $\whApc \in \R^{(g+1) \times (g+1)}$ and \emph{graph part} $\widehat\rmM_l \in \R^{(g+1) \times (g+1)}$.

\subsubsection{The graph part} The graph part $\widehat \rmM_l \in \R^{(g+1) \times (g+1)}$ of $\whAp$ is given by the decomposition
\[
\widehat \rmM_l :=  \begin{pmatrix} \~{\rmM}_l  & \~{\rmW}_l \\
 \~{\rmZ}_l & \rho_l \end{pmatrix}, \quad \~{\rmM}_l = \begin{pmatrix} {\rmM}_l  & 0 \\ 0 & 0 \end{pmatrix}\in \R^{g \times g}, \quad  \~{\rmW}_l = \begin{pmatrix} {\rmW}_{l}   \\ 0  \end{pmatrix} \in \R^{g }, \quad \~{\rmZ}_{l}^\transpose = \begin{pmatrix} {\rmW}_{l}   \\ 0  \end{pmatrix}^\transpose  \in \R^{ g},
 \]
and $\rho_l \in \R$. In short, the parts of $\widehat \rmM_l$ describe the pairings of paths and one-forms in the metric graph $\mgr$ associated to $S$ and the points $p,q,x,y$. 
\smallskip

\begin{itemize} \setlength\itemsep{0.2 mm}
\item $ \~{\rmM}_{l} \in \R^{g \times g} $ is the graph period matrix $\rmM_l \in \R^{\graphgenus \times \graphgenus}$ (see \eqref{eq:M_l}) filled up with zeros.
\smallskip

\item the column vector $\~{\rmW}_{l} \in \R^{g}$ and the row vector  $\~\rmZ_l \in \R^{g}$ are the vectors
\begin{align} \label{eq:DefGraphVectors}
&\rmW_l := P_l(\varrho_{p \tiret q}) \in \R^{\graphgenus},  &\rmZ_l := P_l (\varrho_{x-y}) \in \R^{\graphgenus}
\end{align}
filled up with zeros to dimension $g$. Here, $P_l$ is the $\graphgenus \times |E|$ matrix from \eqref{eq:MatrixP}), that is,
\[
\~{\rmW}_{l}(i ) = {\rmW_l}(i) = \innone{l}{\gamma_i, \varrho_{p \tiret q}}, \qquad \~{\rmZ}_{l}(i ) = {\rmZ_l}(i) = \innone{l}{\gamma_i, \varrho_{x \tiret y}}, \qquad i \in [h].
\]
\smallskip

\item $\rho_l \in \R$ is the scalar $\rho_l = \innone{l}{\varrho_{p\tiret q}, \varrho_{x \tiret y}}$.
\end{itemize}

\subsubsection{The complex part} The complex part $\whApc \in \R^{(g+1) \times (g+1)}$ of $\whAp$ is given by
\begin{equation} \label{eq:ComplexPartApPM}
\whApc :=  \begin{pmatrix} \Apc  & \mathcal{W} \\
\mathcal{Z} & \mathcal{R} \end{pmatrix}, \qquad \Apc   \in \R^{g \times g}, \quad  \mathcal{W} \in \R^{g \times 1}, \quad \mathcal{Z}  \in \R^{1 \times g}, \quad \mathcal{R} \in \R.  
\end{equation}
 The entries of $\whApc$ are (regularized) pairings between certain paths and meromorphic one-forms on the components $C_v$, $v \in V$, of $S$. Note in particular that the regularizations are well-defined, since we have fixed local coordinates around the nodal points $p^e$, $e \in E$, of $S$.

\begin{itemize} \setlength\itemsep{0.2 mm}
\item The matrix $\Apc \in \R^{g \times g}$ has the block decomposition
\[ \Apc :=  \left(
 \begin{matrix}
\Apc_G & \Im B_{v_1} & \Im B_{v_2} & \cdots & \Im B_{v_n}\\
\Im B_{v_1}^\transpose &\Im  \Omega_{v_1} & 0 & \cdots & 0\\
\Im B_{v_2}^\transpose & 0 & \Im \Omega_{v_2}  & \cdots & 0\\ 
\vdots & \vdots  & \vdots & \ddots &0 \\
\Im B_{v_n}^\transpose & 0 & \cdots & \cdots & \Im\Omega_{v_n}
\end{matrix} \right).
\]
The upper left matrix $\Apc_G \in \R^{\graphgenus \times \graphgenus}$ is given by
\[
\Apc_G(i,j) := \sum_{v\in V} \hp{C_v}{\gamma_{i,v}, \omega_{j,v }}', \qquad i,j \in [\graphgenus].
\]
For any $k \in [\graphgenus]$, the matrices $\Im \Omega_{v_k} \in \R^{I_k \times I_k} $ and $\Im(B_{v_k}) \in \R^{\graphgenus \times I_k}$ are the imaginary parts of, respectively, the period matrix $\Omega_{v_k} $ of  $C_{v_k}$ (for the fixed symplectic basis of $H_1(C_{v_k}, \Z)$), and the matrix $B_{v_k} \in \C^{\graphgenus \times I_k}$ given by
\begin{equation} \label{eq:DefApc1}
B_{v_k}(i,j)  =\int_{\gamma_{i,v_k}} \omega_{j,k}, \qquad i \in [\graphgenus], \, j \in I_k.
\end{equation}
\smallskip

\item The column vector $\mathcal{W} \in \R^{g}$ and the row vector  $\mathcal{Z} \in \R^{g}$ are decomposed as
\begin{align*}
&\mathcal{W} =  \begin{pmatrix} \mathcal{W}^\grup \\ \mathcal{W}^\surfup  \end{pmatrix} \in \R^{\graphgenus} \times \R^{g - \graphgenus}, &\mathcal{Z}^\transpose = \begin{pmatrix} \mathcal{Z}^\grup \\ \mathcal{Z}^\surfup  \end{pmatrix} \in \R^{\graphgenus} \times \R^{g - \graphgenus}.
\end{align*}
The entries of the \emph{graph parts} $ \mathcal{W}^\grup \in \R^\graphgenus$ and $\mathcal{Z}^\grup \in \R^\graphgenus$ are given by
 \begin{align*}
 &\mathcal{W}^\grup(i) = \sum_{v \in V} \hp{C_v}{\gamma_{i,v}, \omega_{p \tiret q, v}}', &\mathcal{Z}^\grup(i) = \sum_{v \in V} \hp{C_v}{\gamma_{x \tiret y,v}, \omega_{i,v}}', 
\end{align*}
for $i \in [\graphgenus]$. The \emph{surface parts} $\mathcal{W}^\surfup\in \R^{g - \graphgenus}$ and $\mathcal{Z}^\surfup \in \R^{g - \graphgenus}$ are given by the pairings

\begin{align*}
 & \mathcal{W}^\surfup(i) =  \Im \Big ( \int_{b_i} \omega_{p \tiret q, v_k} \Big ), &\mathcal{Z}^\surfup(i) =  \Im \Big ( \int_{\gamma_{x \tiret y,v_k}} \omega_i \Big ),
\end{align*}
on the Riemann surface components $C_{v_k}$ for every $k \in [n]$ and $i \in I_k$.
\smallskip

\item Finally, the scalar $\mathcal{R} \in \R$ is defined by
\[
\mathcal{R} = \sum_{v \in V} \hp{C_v}{\gamma_{x \tiret y,v}, \omega_{p \tiret q, v}}' .
\]
\end{itemize}

\begin{remark}
As is clear from the definition, the edge length function $l \colon E \to (0, + \infty)$ plays no role in the definition of the complex part  $\widehat \Apc$. To any stable Riemann surface $S$ and four distinct, smooth points $p,q,x,y$ on $S$, we can associate the matrix in \eqref{eq:ComplexPartApPM}, which we denote by $\widehat \Apc(S_0; p,q,x,y)$.
\end{remark}

\subsubsection{Alternative splitting} In the context of the height pairing, we will also decompose $\whAp$ as 
\[
\whAp =  \begin{pmatrix} \Ap & \Ap_W \\
 \Ap_Z & \Ap_\rho
   \end{pmatrix} :=  \begin{pmatrix} \Apc +  (2\pi)^{-1} \~{\rmM_l}  & \mathcal{Z} + (2\pi)^{-1} \~{\rmZ_l} \\
 \mathcal{W} + (2\pi)^{-1} \~{\rmW_l} &   \mathcal{R} + (2\pi)^{-1} \rho_l
   \end{pmatrix} 
\]
for the respective matrices $\Ap \in \R^{g \times g}$, $\Ap_W \in \R^{g \times 1}$, $\Ap_Z \in \R^{1 \times g}$ and $\Ap_\varrho \in \R$.

\subsection{Description and approximation of imaginary parts} \label{ss:DecsrApImPM}
In this section, we use the approximative matrix $\whAp$ to describe the imaginary part $\Im(\ELambda_0)$ of the limit period matrix and approximate the imaginary part $\Im(\bipm_t)$ of the period matrix.

\smallskip

Recall that we have fixed an adapted system of parameters for $B$ and $\rsf$. In particular, for each edge $e=\{u,v\} \in E$, we fix local coordinates around the point $p^e$ on $\rsf$ given by $z^e_u$ and $z_v^e$ with $z^e_v z^e_u = z_e$. On the central fiber $S_0$,  $z^e_u$ and $z^e_v$ are local parameters around the point $p^e_u$ and $p^e_v$ of $C_u$ and $C_v$, respectively. We will use these parameters on stable Riemann surfaces when defining regularized integrals and height pairings. The resulting regularizations only depend on the choice of parameters $z_e$ and not on $z^e_u$ and $z^e_v$.

\smallskip

We obtain the following description of the limit period matrix.

\begin{thm} \label{thm:ImPM}
The imaginary part of the limit period matrix coincides with
\[
\Im(\ELambda_0) = \whApc(S; p_0,q_0,x_0,y_0) = \whApc= \begin{pmatrix} \Apc & \mathcal{Z} \\ \mathcal{W} & \mathcal{R} \end{pmatrix} 
\]
that is, it is equal to the matrix \eqref{eq:ComplexPartApPM} for the stable Riemann surface $S$ with points $p_0,q_0,x_0, y_0$.
\end{thm}

As an immediate consequence, we obtain the following approximation result for $\Im(\bipm_t)$.

\smallskip

For any point $t \in B$, let $t_0 \in B$ be its projection to the stratum $D_E$. That is, $t_0$ has coordinates
\[
z_e(t_0) := z_e(t), \quad e \in [N] \setminus E, \qquad \qquad z_e(t_0) := 0, \quad e \in E.
\]
To $t \in B^\ast$ we associate the metric graph $\mgr_t$ given by $G  =(V,E)$ with the edge length function
$\ell_t(e) := - \log|z_e(t)|$, $e \in E$.  Recall as well that for every index set  $I_k$, $k \in [n]$, the cycles $a_{i,t}$, $b_{i, t}$ for $i \in I_k$ are monodromy invariant and have been spread out to all fibers $\rsf_t$, $t \in B$. For each base point $t \in D_E$, the cycles $a_{i,t}$, $b_{i, t}$ with $i \in I_k$ give a symplectic basis on the component $C_{v_k, t}$ of  the stable Riemann surface $\rsf_t$. The approximative matrices in the following statements are taken with respect to this symplectic basis.

\begin{thm}\label{thm:ImaginaryAp}
Notations as above, for $t \in B^\ast$, 
\[
\Im(\bipm_t) = \widehat \Ap(\rsf_{t_0}, \mgr_t; p_{t_ 0}, q_{t_0}, x_{t_0}, y_{t_0}) + R(t)= \widehat \Apc + (2\pi)^{-1} \widehat  \rmM_{\ell_t} + R(t),
\]
where the error $R(t) \in \R^{(g+1) \times (g+1)}$ goes to zero uniformly if $z_e(t) \to 0$ for all $e \in E$. In particular,
\[
\Im(\bipm_t)  = \whAp(\rsf_{t_0}, \mgr_t; p_{t_0}, q_{t_0}, x_{t_0}, y_{t_0}) + o(1)
\]
if $t \in B^\ast$ converges to $0$ in $B$.
\end{thm}

\begin{proof}[Proof of Theorem~\ref{thm:ImPM}] 
Taken into account that, $\ExtM_{\ell_t} = \sum_{e} \ell_t(e) \ExtM_e$ for all $t  \in B^\ast$ in the notation of Theorem~\ref{thm:nilpotentorbit2}, we have to prove that
\[\lim_{t \to 0} \Im(\bipm_t) - (2\pi)^{-1} \ExtM_{\ell_t} = \whApc(S_0; p_0, q_0, x_0, y_0).
\] Applying Proposition~\ref{prop:LPM1} and Proposition~\ref{prop:LPM1}, it suffices to show that
\begin{align*}
&\Im(\Omega_t) - (2\pi)^{-1} \~\rmM_{\ell_t} \to \Apc_G,	&&\Im(\rmW_t) - (2\pi)^{-1} \~\rmW_{\ell_t} \to \mathcal{W}, \\
&\Im(\rmZ_t) - (2\pi)^{-1} \~\rmZ_{\ell_t} \to \mathcal{Z},	&&\Im(\rho_t) - (2\pi)^{-1} \rho_{\ell_t} \to \mathcal{R},
\end{align*}
as $t \in B^\ast$ converges to $0$ in $B$. In the following, we prove the first limiting behavior. The remaining claims can be obtained analogously.

Fix $i,j \in [h]$. For any point $t \in B^\ast$ and $\tilde t \in \widetilde{B^\ast}$ lying above $t$, the $(i,j)$-th entries satisfy
\begin{equation} \label{eq:ProofPeriodMatrixEntries} 
 \Im(\Omega_t) (i,j) - (2\pi)^{-1}  \~\rmM_{\ell_t} (i,j) = \Im(\int_{b_{i, \tilde t}} \omega_j) -\sum_{e} \Im(\zeta_e(\tilde t)) \innone{e}{\gamma_i, \gamma_j}. 
\end{equation}
Based on the decompositions of the Riemann surface $\rsf_t$ obtained by adapted coordinates (see \eqref{eq:AdaptedCoordinates}), we decompose the difference \eqref{eq:ProofPeriodMatrixEntries} into terms whose asymptotics we can describe.

\smallskip

Consider the decomposition of the Riemann surface $\rsf_t$ into the regions $B_{e,t}$, $e\in E$, and the Riemann surfaces with boundary 
\[Z_{v, t}\coloneqq Y_{v, t} \cup \bigsqcup_{e \sim v} A^e_{v,t}, \qquad v\in V.\]

Recall that for $e$ with extremities $u$ and $v$, we set (see Section~\ref{sec:hybrid_log_map_einf})
  \begin{align*}
   A^e_{u,t} &: =\Bigl\{ (z_u^e, z^e_v) \,\bigl|\, z_u^e z^e_v =z_e(t),\, |z_u^e|, |z_v^e| \leq 1 \, \textrm{and} \,  |z_u^e| \geq \varrho_e \Bigr\} \\
A^e_{v,t} &: =\Bigl\{ (z_u^e, z^e_v) \,\bigl|\, z_u^e z^e_v =z_e(t),\, |z_u^e|, |z_v^e| \leq 1 \, \textrm{and} \,  |z_v^e| \geq \varrho_e \Bigr\},
\end{align*}
 where $\varrho_e := \varrho_e(t) = 1/ \ell_t(e)$, and that $B_{e,t} := W_{e, t} \setminus ( A^e_{u,t} \sqcup  A^e_{v,t} )$.  We write 
 \[
 b_{i, \tilde t} = \sum_{v\in V} b^{v}_{i, \tilde t} + \sum_{e\in E} b^e_{i, \tilde t}
 \]
 where $b^v_{i, \tilde t}$ is the part of $b_{i, \tilde t}$ lying on $Z_{v,t}$, and $b^e_{i, \tilde t}$ is the part lying in $B_{e,t}$. 

\medskip

We can decompose $\Im(\int_{b_{i, \tilde t}} \omega_j)$ as the sum over the above paths 
\[\Im(\int_{b_{i, \tilde t}} \omega_j) =\sum_{v\in V} \Im(\int_{b^{v}_{i, \tilde t}} \omega_j) +\sum_{e\in E}\Im(\int_{b^{e}_{i, \tilde t}} \omega_j). \]

In particular, the difference in \eqref{eq:ProofPeriodMatrixEntries} can be written as
\begin{align*}
\sum_{v\in V} \Im(\int_{b^{v}_{i, \tilde t}} \omega_j) + \sum_{e \sim v} \frac{1}{2\pi} \log \varrho_e \innone{e}{\gamma_i, \gamma_j} 
+ \sum_{e\in E}  \Im(\int_{b^{e}_{i, \tilde t}} \omega_j) -  \frac{1}{2\pi}  \bigl ( \ell_t(e) + 2 \log \varrho_e \bigr)  \innone{e}{\gamma_i, \gamma_j}.  \end{align*}
In order to conclude, observe that for each vertex $v \in V$,
\[\lim_{t\to 0} \, \Im(\int_{b^{v}_{i, \tilde t}} \omega_j) + \sum_{e \sim v} \frac{1}{2\pi} \log \varrho_e  \innone{e}{\gamma_i, \gamma_j} = \hp{C_v}{ \gamma_{i, v} ,  \omega_{j, v}  }',\]
and for each edge $e \in E$,
\[\lim_{t\to 0}   \, \Im(\int_{b^{e}_{i, \tilde t}} \omega_j) -\frac{1}{2\pi} \bigl ( \ell_t(e) + 2 \log \varrho_e \bigr)  \innone{e}{\gamma_i, \gamma_j} = 0,\]
since $\omega_{j,0}\rest{C_v} = \omega_{j, v}$ and $b^v_{i, \tilde t}$ converges to $\gamma_{i, v}$. Moreover, the residue of the meromorphic form $\omega_{j,v}$ at $p^e_v$ coincides with $\gamma_j(e)$, so the terms calculate precisely the regularized integration. 
\end{proof}


\section{Tameness of the height pairing: Proof of Theorem~\ref{thm:tameness-height-pairing1}} \label{sec:ProofHeightPairing}
In this section we prove the stratumwise tameness of the height pairing, as stated in Theorem~\ref{thm:tameness-height-pairing1}. Recall that, given a smooth Riemann surface $\rsf_t$ with four markings $p_t,q_t,x_t,y_t$, we are interested in the limit of the height pairing $\hp{\rsf_t}{p_t - q_t, x_t-y_t}$, as $\rsf_t$ degenerates to a hybrid curve $\curve$.

\smallskip

Our approach can be summarized as follows:

\medskip

\noindent -  From the preceding sections, see Lemma~\ref{lem:HeightPairingByPeriodMap} and Theorem~\ref{thm:ImaginaryAp}, we have an expression of the height pairing using the imaginary part $\Im(\bipm_t)$ of a Hodge theoretic period map and also the  asymptotics of $\Im(\bipm_t)$ in terms of the approximative matrix $\widehat \Ap$.
\smallskip

\noindent - In the present section, we use the shape of the approximative matrix $\whAp$ in order to separate the asymptotic of the height pairing $\hp{\rsf_t}{p_t - q_t, x_t - y_t}$ into two parts: a metric graph term and a complex term originating from Riemann surfaces. We state and prove the important Separation Theorem~\ref{thm:SeparateComplexGraph} (see Sections~\ref{ss:SeparationThmTameHeightPairing} through~\ref{ss:ProofEasySeparation}). The formalization of this is provided by an auxiliary notion of height pairing on metrized complexes (see Section~\ref{ss:AuxiliaryMetrizedComplexPairing}). The formulation of the Separation Theorem~\ref{thm:SeparateComplexGraph}  is slightly elaborate, however, we will need this version later in the proof of the uniform tameness of the height pairing.  Our separation theorem refines and expands the main result of~\cite{ABBF}, by providing the missing complex part in~\cite{ABBF}, dealing as well with the more delicate situation where the points go to the same limit, including the singular points.

\smallskip

\noindent - Using the topological connections between the moduli spaces of  hybrid curves and tropical curves, we then conclude by applying the results from the second part of our paper.

\smallskip

Using the above approach, we prove Theorem~\ref{thm:tameness-height-pairing1} in Section~\ref{ss:ProofTamenessHP1}.

\subsection{Height pairing on metrized complexes} \label{ss:AuxiliaryMetrizedComplexPairing} In what follows, we introduce a  \emph{height pairing} between degree zero divisors on a \emph{metrized complex}. This provides the formal framework to split the asymptotics of height pairings into a graph and complex part.

\smallskip

Let $\mc$ be a metrized complex, given by a stable Riemann surface $S$ (possibly with markings) and an edge length function $l \colon E \to (0, + \infty)$ on the edges of the dual graph $G = (V,E)$. The underlying metric graph is denoted by $\mgr$. Consider the natural contraction map $\mc \to \mgr$, which contracts each component $C_v$, $v \in V$, of $S$ to the vertex $v \in \mgr$. Using our previous convention, the image of $x \in \mc$ under this map is denoted by $\bar x \in \mgr$. 

\smallskip

Given two points $p,q \in \mc$, we associate the following collection of degree zero divisors  $D^\cmc_{p \tiret q,v}$, $v \in V$, on the components $C_v$ of $S$. Let $f \colon \mgr \to \R$  be a solution to the Poisson equation
\begin{equation} \label{eq:PoissonMCPairing}
\Delta_\mgr f = \delta_{\bar p} - \delta_{\bar q}
\end{equation}
on the metric graph $\mgr$ (note that $f$ is unique up to an additive constant). For a Riemann surface component $C_v$, $v \in V$, of $S$, define the degree zero divisor $D^\cmc_{p \tiret q,v} \in \Div^0(C_v)$ as
\begin{equation} \label{eq:ComplexDivisorsMetrizedComplex}
D^\cmc_{p \tiret q,v} = \delta_{v, \bar p}\, p -  \delta_{v, \bar q} \, q - \sum_{e \sim v} \slp_e f(v) \cdot p^e_v.
\end{equation}
The divisors $D^\cmc_{p \tiret q, v}$, $v \in V$, only depend on the points $p,q \in \mc$ (and not on the choice of $f$).

\smallskip

Fix further a local parameter $z^e_u$ on $C_u$ around each attachment point $p^e_u$, $e \sim u$, on each component $C_u$, $u \in V$ of $S$. For four points $p,q,x,y$ with $\{p, q \} \cap \{x,y\} = \varnothing$, the quantity
\begin{equation}
\hp{\mc}{p - q, x - y} := \hp{\mgr}{\bar p - \bar q, \bar x - \bar y} + \sum_{v \in V} \hp{C_v}{D^\cmc_{p \tiret q, v} , D^\cmc_{x \tiret y, v}}'
\end{equation}
is called the \emph{metrized complex height pairing} between the divisors $p - q$ and $x - y$ on $\mc$ (with respect to the local parameters $z^e_u$, $e \sim u$, $u \in V$).
 
 The following proposition relates the metrized complex height pairing to the hybrid height pairing on hybrid curves (see Section~\ref{ss:HeightPairingHybridCurves}).
  
\begin{prop} Let  $\curve$ be the hybrid curve of rank one associated to $\mc$ and let $\mccan$ be the canonical representative of $\curve$ with normalized edge length. Let $L=\sum_{e\in E} l(e)$ be the total length of the metric graph $\mgr$. We have
\[\hp{\mc}{p - q, x - y} = L\,\lhp{\curve,1}{p-q, x-y} + \lhp{\curve, \smallcc}{p-q, x-y}\]
where $\lhp{\curve,1}{p-q, x-y}$ and $\lhp{\curve,\smallcc}{p-q, x-y}$ denote the components of the hybrid height pairing $\lhp{\curve}{p-q, x-y}$ between the two divisors $p-q$ and $x-y$ on $\curve$.
\end{prop}
\begin{proof} The proof follows by observing that \[\hp{\mgr}{\bar p - \bar q, \bar x - \bar y} = L \hp{\Gamma^1}{\bar p - \bar q, \bar x - \bar y} = L\,\lhp{\curve,1}{p-q, x-y}.\] 
\end{proof}

\subsection{Separation Theorem} \label{ss:SeparationThmTameHeightPairing}

Let $\rsf / B$ be a versal family of a stable Riemann surface $S_0$ (possibly carrying markings $r_1, \dots r_\nmark$) with dual graph $G = (V,E)$, together with a fixed choice of adapted coordinates. Consider the associated family of hybrid curves $\rsf^\hyb/ B^\hyb$.

\smallskip

Let $\shy$ be a fixed point in the hybrid boundary $\partial_\infty B^\hyb$ lying above the point $0$ in $B$. That is, $\shy$ belongs to the hybrid stratum $D_\pi^\hyb$ of an ordered partition $\pi = (\pi_1, \dots, \pi_r)$ on $E$, and has the form $\shy = (l,0)$ for the complex coordinate $s = 0$ and simplicial coordinates $l = (l_e)_{e \in E}$.  

\smallskip

For a base point $t \in B^\ast$, let $t_0$ be its projection to the stratum $D_E$, given in adapted coordinates by
\[
 	z_e(t_0) := z_e(t), \qquad e \in [N] \setminus E, \qquad \qquad z_e(t_0) := 0, \qquad e \in E.
\]
The point $t_0$ belongs to the open stratum $\inn D_E$, which contains also $0$. In particular, the dual graph of the stable Riemann surface $\rsf_{t_0}$ coincides with $G$. Equipping $G$ with the edge lengths 
\[
\ell (e)  = \ell_t(e) := - \log|z_e(t)|, \quad e \in E,
\]
we obtain a metric graph $\mgr_t$. Together, $\rsf_{t_0}$ and $\mgr_t$ define a metrized complex $\mc_t$, $t \in B^\ast$.

\smallskip 

Recall from Section~\ref{sec:hybrid_log_map_einf} that using the adapted coordinates, we decompose every smooth fiber $\rsf_t$, $t \in B^\ast$, as
\[
\rsf_t = \bigsqcup_{v \in V} Y_{v,t} \sqcup \bigsqcup_{e = uv \in E} A^e_{u,t} \sqcup A^e_{v,t} \sqcup B_{e,t}.
\]
This allows to view any point $y$ on $\rsf_t \setminus \bigcup_{e \in E} B_{e,t}$ as a point on $\rsf_{t_0}$. Indeed, it suffices to identify, for any vertex $v$, the part $Y_{v,t_0}$ with the corresponding subset of $C_{v,t_0}$, and for each edge $e = uv$, the region $A^e_{u,t}$ with a disc $\{\varrho_e(t) \le |z^e_u| \le 1\}$ on $C_{u, t_0}$ (and similar for $A^e_{v,t}$).

\smallskip

Consider now a base point $t \in B^\ast$ close to $0$ in $B$. For four pairwise distinct points $p,q,x,y$ on $\rsf_t \setminus \bigcup_{e \in E} B_{e,t}$, we can hence consider the approximative matrix
\[
\widehat \Ap (\rsf_{t_0},\mgr_t; p,q,x,y)
\]
introduced in Section~\ref{ss:Ap}. For this, we fix a symplectic basis relative to a layered spanning tree $T$ for $(G, \pi)$, and a symplectic basis $a_j$, $b_j$ on each component of $S_0$. We then spread them out to a symplectic basis $a_{j, s'}, b_{j, s'}$ for the components of the fibers $\rsf_{s'}$ close to $S_0$; the local coordinates $z^e_u$ in the regularization are induced from the adapted coordinates on $\rsf / B$.
 
\smallskip

The following result provides the aforementioned separation of the height pairing.

\begin{thm}[Separation Theorem] \label{thm:SeparateComplexGraph}
For $t \in B^\ast$ close to $s$ in $B$, consider a region $R_t$ contained in
\begin{equation} \label{eq:SepThAssumption}
R_t \subseteq \Bigl \{\,(p,q,x,y) \in \rsf_t^4 \, \st \, p,q,x,y \in \rsf_t \setminus  \bigcup_{e \in E} B_{e,t} \,\Bigr\}.
\end{equation}
For $(p,q,x,y) \in R_t$, consider a real matrix $(g+1) \times (g+1)$ matrix
\[
\widehat \Ap' (p,q,x,y) = \begin{pmatrix} \Ap'  (p,q,x,y)  &  \mathcal{A}_\rmZ'  (p,q,x,y) \\  \Ap_\rmW'  (p,q,x,y) &  \Ap_\varrho'  (p,q,x,y)\end{pmatrix}
\]
which is close to the approximative matrix $\widehat \Ap$, in the sense that
\begin{equation} \label{eq:SepTh1}
 \widehat \Ap'(p,q,x,y) =  \widehat \Ap(\rsf_{t_0}, \mgr_t; p,q,x,y) + o(1),
\end{equation}
where the $o(1)$ term goes to zero uniformly for $(p,q,x,y) \in R_t$ as $t \to \shy$ in $B^\hyb$.

Then, the following asymptotics holds
\begin{equation*} \label{eq:SepTh2}
 2 \pi \Big( \Ap_\varrho'(p,q,x,y) -  \Ap_\rmW'(p,q,x,y)   \Ap'(p,q,x,y)^{-1}   \Ap_\rmZ'(p,q,x,y) \Big)  = \hp{\mc_t}{p - q, x - y} + o(1),
\end{equation*}
where the $o(1)$ term goes to zero uniformly for $(p,q,x,y) \in R_t$ as $t \to \shy$ in $B^\hyb$.

\end{thm}
Note also that the resulting expression in Theorem~\ref{thm:SeparateComplexGraph} is equal to
\begin{equation} \label{eq:RecallMCPairing}
\hp{\mc_t}{p - q, x - y} = \hp{\mgr_t}{\bar p - \bar q, \bar x - \bar y} + \sum_{v \in V} \hp{C_{v, t_0}}{D^\cmc_{p \tiret q, v} , D^\cmc_{x \tiret y, v}}'
\end{equation}
which justifies the name given to the theorem. 

\begin{remark}
Recall that, in the definition of the regions $A^e_{u,t}$ and $B_{e,t}$ on the smooth fiber $\rsf_t$, $t \in B^\ast$, we used the radius function $\varrho_e(t) := -1/\log|z_e(t)|$. As the proof shows, the statement remains valid for other functions $\varrho_e'(t)$ such that $\lim_{t \to 0} \varrho_e'(t) / \ell^{\operatorname{ess}}_t(e) = 0$ for all edges $e \in E$. Here $ \ell^{\operatorname{{\tiny ess}}}_t(e)$ is the sum of the edge lengths of the longest path in $\mgr_t$, which does not repeat any edges, contains $e$ and whose intermediate vertices are  all of degree $2$.
\end{remark}

\subsection{The inverse of the approximative matrix} The proof of Theorem~\ref{thm:SeparateComplexGraph} requires a precise description of the inverse of the approximative matrix close to a hybrid point $\shy$. This information is provided in the following and extends the results in \cite[Section 9]{AN}.

\smallskip

Notations as in the previous subsection, consider a base point $t \in B^\ast$. Let
\[
\Ap(t) := \Ap(\rsf_{t_0}, \mgr_t), \qquad t \in B^\ast,
\]
be the approximative matrix defined in Section~\ref{ss:Ap} (note that $\Ap$ is actually defined independently of the points $p,q,x,y$, justifying our notation). 

\smallskip

In the following, it will be convenient to introduce the new edge length function
\begin{equation} \label{eq:piEdgeLengths}
\fell_t(e) := - \frac{1}{2 \pi} \log|z_e(t)|, \qquad e \in E,
\end{equation}
on $G$ and write $\Ap(t)$ in the form 
\[
\Ap(t) = \begin{pmatrix} \Apc_{G,t} + \rmM_{\fell_t} & \Im(B_t) \\  \Im(B_t)^\transpose & \Im(\Omega_{S,t})
\end{pmatrix}.
\]
where $\rmM_{\fell_t}$  is the graph period matrix~\eqref{eq:M_l} for $\fell_t$.  Note that $2 \pi \fell_t(e) = \ell_t(e)$ for the edge lengths $\ell_t(e)$, $e \in E$, considered in Theorem~\ref{thm:SeparateComplexGraph}. However, the choice of this normalization allows to considerably streamline the calculations.

\smallskip

Moreover, $B_t \in \C^{\graphgenus \times (g - \graphgenus)}$ is given by
\[
B_t  = \Big ( B_{v_1, t}, \dots, B_{v_n,t} \Big),
\]
where $B_{v_k,t} \in \C^{\graphgenus \times I_k}$ is the matrix from \eqref{eq:DefApc1} for the component $C_{v_k, t_0}$ of $\rsf_{t_0}$, and
\[
\Omega_{S,t}  =\left(
 \begin{matrix}
\Im ( \Omega_{v_1, t} ) & 0 & \cdots & 0\\
 0 & \Im( \Omega_{v_2, t})  & \cdots & 0\\ 
\cdots & \vdots & \ddots &0 \\
 0 & \cdots & \cdots & \Im(\Omega_{v_n, t})
\end{matrix} \right),
\]
where $\Omega_{v_k,t} \in \C^{I_k \times I_k}$ are the period matrices of the components $C_{v_k, t_0}$ of $\rsf_{t_0}$. 

\smallskip

In the following, in order to simplify the presentation, we often omit the dependence on $t$ in the notation and simply write $\Apc_G$, $M_\fell$, $\Omega_S$, etc.

\smallskip
Let $\shy =(l,0) \in B^{\hyb}$ with complex coordinate $s = 0$, lying in the hybrid stratum $D_\pi^\hyb$ of an ordered partition $\pi=(\pi_1, \dots, \pi_r)$ on $E$. Consider the growth functions $\fEll_j \colon B^\ast \to (0, + \infty)$ given by
\[
	\fEll_j(t) := \sum_{e \in \pi_j} \fell_t(e) = \sum_{e \in \pi_j} -\frac{1}{2\pi}\log\abs{z_e(t)} , \qquad j=1,\dots, r,
\] 
and $\fEll_{r+1} \equiv 1$ on $B^\ast$. First of all, recall from \cite[Theorem 9.4]{AN} that the inverse has the block structure
\begin{equation} \label{eq:InvPerOld}
\mathcal{A}(t)^{-1} = \Big (\fEll_{\min\{j,k\}}^{-1} \bigl(C_{jk} +o(1) \bigr) \Big )_{j,k \in [r+1]} , \qquad t \to \shy,
\end{equation}
for some matrices $C_{jk} \in \R^{h_\pi^j \times h_\pi^k}$, where $h_\pi^j$ is the genus of the $j$-th graded minor $\gr_\pi^j(G)$ of $G$, for $j\in [r]$, and $h_{r+1} := g-h$.  The matrices on the diagonal are explicitly given by
\[
	C_{r+1, r+1} = \Im(\Omega_{S, 0}) =\left(
 \begin{matrix}
\Im ( \Omega_{v_1, 0} ) & 0 & \cdots & 0\\
 0 & \Im( \Omega_{v_2, 0})  & \cdots & 0\\ 
\cdots & \vdots & \ddots &0 \\
 0 & \cdots & \cdots & \Im(\Omega_{v_n, 0})
\end{matrix} \right),
\]
where $\Omega_{v_k, 0}$ is the period matrix of the component $C_{v_k}$ of $S_0$, and
\[
	C_{jj} = \rmM_{j,l}^\pi, \qquad j=1, \dots, r,
\]
where the matrix $\rmM_{j,l}^\pi \in \R^{h_\pi^j \times h_\pi^j}$ is the graph period matrix \eqref{sec:explicit_extension_form} for the $j$-th graded minor $\gr_\pi^j(G)$ of $G$, equipped with the edge lengths $l(e)$, $e \in E$, given by the hybrid point $\shy = (l,0)$.

\smallskip

The description of the height pairing requires also information on the off-diagonal terms. The following proposition suffices in the proof of Theorem~\ref{thm:SeparateComplexGraph}.

\begin{prop} \label{prop:InvLem2.0}
Notations as above, in the regime when $t \to \shy$ in $B^\hyb$, the inverse of $\Ap (t)$ is given by
\[
\Ap(t)^{-1} = 
 \left( \begin{matrix}  \rmM_\fell^{-1} + \rmM_\fell^{-1} \big (\Im(B)  \Im(\Omega_S)^{-1} \Im(B)^\transpose -  \Apc_G \big ) \rmM_\fell^{-1}  & - \rmM_\fell^{-1} \Im(B) \Im(\Omega_S)^{-1}   \\ - \Im(\Omega_S)^{-1} \Im(B)^\transpose \rmM_\fell^{-1}  & \Im(\Omega_S)^{-1}  \end{matrix} \right ) + R
\]
where there error term $R \in \R^{g\times g}$ has the following order
\[
R(t) =  \left( \begin{matrix} \rmM_\fell^{-1} o(1) \rmM_\fell^{-1} &  \rmM_\fell^{-1}o(1) \\ o(1) \rmM_\fell^{-1} & o(1) \end{matrix} \right ), \qquad \text{as } t\to \shy.
\]
\end{prop}
\begin{proof}
Recall that the inverse of a symmetric $2 \times 2$- block matrix is expressed in terms of the Schur complement, that is,
\[
\left( \begin{matrix} A & B \\ B^\transpose & D \end{matrix} \right )^{-1} = \left( \begin{matrix}  A^{-1} + A^{-1} B \zeta^{-1} B^\transpose A^{-1} & - A^{-1} B^\transpose \zeta^{-1}   \\ - \zeta^{-1} B^\transpose A^{-1}  & \zeta^{-1} \end{matrix} \right )
\]
where $\zeta =   D - B^\transpose A^{-1} B$. To obtain the proof, we need to apply the above to the matrix
\[
\Ap(t) =  \left(
 \begin{matrix}
\Apc_G + \rmM_\fell & \Im(B_{v_1})  &\Im(B_{v_2}) & \cdots & \Im(B_{v_n})\\
\Im(B_{v_1})^\transpose & \Im(\Omega_{v_1}) & 0 & \cdots & 0\\
\Im(B_{v_2})^\transpose & 0 & \Im(\Omega_{v_2})  & \ddots & 0\\ 
\vdots & \vdots &  \ddots & \ddots &0 \\
\Im(B_{v_n})^\transpose & 0 & \cdots & \cdots & \Im(\Omega_{v_n})
\end{matrix} \right) = \left (   \begin{matrix} 
\Apc_G + \rmM_\fell & \Im(B)\\
\Im(B)^\transpose & \Im(\Omega_S).
\end{matrix} \right)
\]
We start by deriving the asymptotics of the Schur complement matrix $\zeta$. Applying \cite[Lemma 9.3]{AN}, we conclude that
\begin{equation} \label{eq:GraphInverse}
\rmM_\fell^{-1} = \Big (\fEll_{\min\{j,k\}}^{-1} \bigl(C_{jk} +o(1)\bigr) \Big )_{j,k \in[r+1]}
\end{equation}
for some matrices $C_{jk} \in \R^{h_\pi^j \times h_\pi^k}$ as $t \to \shy$. In particular, $\rmM_{\fell}^{-1} = o(1)$ as $t \to \shy = (l,0)$, which in turn implies that $\zeta = \Im(\Omega_S) + o(1)$ and
\[
\zeta^{-1} = \Im(\Omega_S)^{-1} + o(1).
\]
Hence we have proven the claimed asymptotics for the lower right block of the inverse $\Ap(t)^{-1}$. To describe the other terms, note first that
\begin{align*}
A^{-1} &= (\Apc_G + \rmM_\fell)^{-1} = \rmM_\fell^{-1} (\operatorname{Id} + \Apc_G \rmM_\fell^{-1})^{-1} \\
&= \rmM_\fell^{-1} - \rmM_\fell^{-1} \Apc_G \rmM_\fell^{-1} + \rmM_\fell^{-1} o(1) \rmM_\fell^{-1}.
\end{align*}
Taking into account the formula involving the Schur complement, we arrive at the result.
\end{proof}

\subsection{An auxiliary proposition}

Recall that in Theorem~\ref{thm:SeparateComplexGraph}, we want to describe the left hand side in \eqref{eq:SepTh2}. The assumption \eqref{eq:SepTh1} suggests that, in the limit $t \to \shy$, we can replace $\widehat \Ap'$ by the approximative  matrix $\widehat \Ap$. The next proposition ensures that this is indeed the case.

\smallskip

Recall that the vectors $\Ap_\rmZ, \Ap_\rmW \in \R^g$ are defined as the sums (see Section~\ref{ss:Ap})
\[
\Ap_\rmZ = \mathcal{Z} + \widetilde \rmZ_\fell = \begin{pmatrix} \mathcal{Z}^\grup \\ \mathcal{Z}^\surfup \end{pmatrix} + \begin{pmatrix} \rmZ_\fell \\0 \end{pmatrix} \in \R^{\graphgenus} \times \R^g, \quad \Ap_\rmW = \mathcal{W} + \widetilde \rmW_\fell = \begin{pmatrix} \mathcal{W}^\grup \\ \mathcal{W}^\surfup \end{pmatrix} + \begin{pmatrix} \rmW_\fell \\0 \end{pmatrix} \in \R^{\graphgenus} \times \R^g,
\]
for the vectors $\mathcal{Z}, \mathcal{W} \in \R^g$, and $\mathcal{Z}^\grup, \mathcal{W}^\grup, \rmZ_\fell, \rmW_\fell \in \R^\graphgenus$, and $\mathcal{Z}^\surfup, \mathcal{W}^\surfup \in \R^{g - \graphgenus}$. More precisely, the graph vectors
\[
\rmZ_\fell =\rmZ_{\fell_t}(x,y) \in \R^g, \qquad \qquad \rmW_{\fell} = \rmW_{\fell_t}(p,q) \in \R^g,
\]
are defined by \eqref{eq:DefGraphVectors}, using the edge length function $\fell = \fell_t$ from \eqref{eq:piEdgeLengths}.

\begin{prop} \label{prop:AuxPropTHP} Assumptions as in Theorem~\ref{thm:SeparateComplexGraph}, as $t \in B^\ast$ converges to $\shy = (l,0)$,
\[
\Ap_\rho' - \Ap_\rmW'  (\Ap')^{-1} \Ap_\rmZ' = \Ap_\rho + \Ap_\rmW \Ap^{-1} \Ap_\rmZ + o(1) 
\]
where the $o(1)$-term goes to zero uniformly for $(p,q,x,y) \in R_t$ as $t \to \shy$ in $B^\hyb$.
\end{prop}
\begin{proof}
We begin by estimating the size of the vectors $\Ap_\rmW, \Ap_\rmZ \in \R^g$. Recall that we decomposed $[g] = \{1, \dots, g\}$ into the index sets $I_\pi^j$, $j = 1, \dots, r+1$, according to the decomposition $g = \sum_{j=1}^{r+1} h_\pi^j$ (as above, $h_\pi^j$ is the genus of the $j$-th graded minor, $j\in [r]$, and $h_\pi^{r+1}=g-h$). Since our basis $\gamma_1, \dots, \gamma_h$ of $H_1(G, \Z)$ is admissible for the ordered partition $\pi$, it is easily verified that,  as $t \to \shy$ in $B^\hyb$,
\begin{align} \label{eq:VectorAs1} 
\rmZ_\fell \rest{J_\pi^j}  = O(\fEll_j(t)) \text{ and } \rmW_\fell \rest{J_\pi^j}  = O(\fEll_j(t)), \qquad j=1, \dots, r,
\end{align}
uniformly for $(p,q,x,y)$ in $R_t$. With a little effort one can show that for $j=1, \dots, r$, 
\begin{align} \label{eq:VectorAs2} 
\mathcal{Z}^\grup \rest{J_\pi^j} &= O \Big (\max_{e \in \pi_j \cup \dots \pi_r} |\log \varrho_e(t)| \Big)   = o(1) \fEll_j(t), \\ \label{eq:VectorAs3} 
\mathcal{W}^\grup \rest{J_\pi^j}  &= O \Big (\max_{e \in \pi_j \cup \dots \pi_r} |\log \varrho_e(t)| \Big)  = o(1) \fEll_j(t),
\end{align}
uniformly for $(p,q,x,y) \in R_t$ as $t \to \shy$. Indeed, the assumption \eqref{eq:SepThAssumption} on the region $R_t \subseteq \rsf_t^4$  ensures that the points $p,q,x,y$, seen on the singular Riemann surface $\rsf_{t_0}$, stay outside the circles $|z^e_u| \le \varrho_e(t)$ around the nodes $p^e_v$, which gives the logarithmic estimate. Moreover, it is clear that
\begin{equation} \label{eq:VectorAs4} 
\mathcal{Z}^\surfup = O(1) \qquad \text{and}  \qquad \mathcal{W}^\surfup = O(1)
\end{equation}
uniformly. Altogether, the entries of $\Ap_\rmZ$ and $\Ap_\rmW$ have orders
\begin{equation} \label{eq:VectorAs5} 
\Ap_\rmZ\rest{J_\pi^j} = O(\fEll_j(t))  \qquad \text{ and }  \qquad \Ap_\rmW\rest{J_\pi^j}  = O(\fEll_j(t)), \qquad j=1, \dots, r+1,
\end{equation}
uniformly for $(p,q,x,y)$ in $R_t$ as $t \to \shy$ in $B^\hyb$.  Taking into account the asymptotics of $\Ap^{-1}$ given by \eqref{eq:InvPerOld}, it follows that
\[
\Ap_\rmZ \Ap^{-1}   = O(1) \qquad \text{and} \qquad   \Ap^{-1} \Ap_\rmW  = O(1).
\]
Combining the above, we readily conclude that
\[
\Ap_\rmZ'  (\Ap')^{-1} \Ap_\rmW'= (\Ap_\rmZ + o(1))  (\Ap + o(1))^{-1} (\Ap_\rmW+ o(1)) =\Ap_\rmZ  \Ap^{-1} \Ap_\rmW + o(1),
\]
and the claim follows upon recalling that $\Ap_\rho' = \Ap_\rho + o(1)$ uniformly for $(p,q,x,y) \in R_t$.
\end{proof}

\subsection{Proof of Theorem~\ref{thm:SeparateComplexGraph}} \label{ss:ProofEasySeparation}
By the preceding discussion, in order to prove Theorem~\ref{thm:SeparateComplexGraph}, it suffices to understand the behavior of the quantity
\[
\Ap_\rho (p,y,x,y) - \Ap_W (p,y,x,y)  \Ap (p,y,x,y)^{-1} \Ap_Z (p,y,x,y) \]
as $t$ converges to $\shy = (l,0)$ in $B^\hyb$ and $(p,q,x,y)$ belongs to $R_t$. To analyze this expression, we order the appearing terms as
\begin{equation} \label{eq:AnalyzableTerms}
	\Ap_\rho - \Ap_\rmW  \Ap^{-1} \Ap_\rmZ = A_1 + A_2 + A_3 + A_4,
\end{equation}
with the following four parts (see Section~\ref{ss:Ap} for the definition of the appearing matrices)
\begin{align*}
A_1 &= \Ap_\rho -  \Ap_\rmZ \Ap^{-1} \Ap_\rmW ,   &&A_2= \varrho_\fell - \widetilde \rmZ_\fell \Ap^{-1} \widetilde \rmW_\fell \\
A_3& = - \widetilde \rmZ_\fell \Ap^{-1} \Ap_\rmW     &&A_4= - \Ap_\rmZ \Ap^{-1} \widetilde \rmW_\fell.
\end{align*}

In order to abbreviate, here and in the following, we often omit the dependence on $p,q,x,y$ and $t$ in the notation.

\smallskip

It turns out that all $A_i$'s can be expressed in terms of (regularized) height pairings, either on graphs or on Riemann surfaces, as we explain next.

\medskip

Given $t \in B^\ast$, let $\projhar\colon C^1(G, \R) \to \Omega^1(G)$ be the orthogonal projection with respect to the scalar product $\innone{\fell_t}{\cdot\,, \cdot}$ given by the edge lengths $\fell_t(e)$, $e \in E$ (see \eqref{eq:piEdgeLengths} and Section~\ref{ss:ProjectionHarmonicOneForms}). Consider the harmonic one-forms $\chi\indbis{x\tiret y}, \chi\indbis{p\tiret q} \in \Omega^1(G, \R)$ given by projecting the paths $\varrho_{p \tiret q}$ and $\varrho_{x\tiret y}$ to the space of harmonic one-forms, that is,
\begin{align*}
	&\chi_{x\tiret y} := \projhar (\varrho_{x \tiret y}),  &\chi_{p\tiret q} := \projhar(\varrho_{p\tiret q}).
\end{align*}

 Note that, the paths $\varrho_{p \tiret q}$ and $\varrho_{x\tiret y}$ in the spanning tree $T$ of $G$ are naturally viewed as one-forms on $G$. 
 
 \smallskip
 
The harmonic one-forms $\chi_{p\tiret q}$ and $\chi_{x\tiret y}$ also induce the following divisors
\[
D_{x\tiret y, v}^\perp := \sum_{e = vu\in \E} \chi_{x\tiret y}(e)\,  p^e_v, \qquad D_{p\tiret q,v}^\perp = \sum_{e = vu\in \E} \chi_{p\tiret q}(e) \, p^e_v,
\]
on the components $C_{v, t_0}$, $v \in V$, of the stable Riemann surface $\rsf_{t_0}$. Since $\chi_{p\tiret q}$ and $\chi_{x\tiret y}$ are harmonic one-forms, all divisors $D_{x\tiret y, v}^\perp$ and $D_{p\tiret q, v}^\perp$ are of degree zero. 

The relevance of this construction is that the divisors $D^\cmc_{p\tiret q, v}$, $D^\cmc_{x\tiret y, v}$ on $C_{v, t_0}$ in the height pairing on $\mc_t$ (see \eqref{eq:RecallMCPairing} and \eqref{eq:ComplexDivisorsMetrizedComplex}) can be recovered from the formula
\begin{equation}
D_{x\tiret y,v}^\cmc = D_{x\tiret y,v} - D_{x\tiret y, v}^\perp, \qquad D_{p\tiret q,v}^\cmc = D_{p\tiret q,v} - D_{p\tiret q, v}^\perp,
\end{equation}
where $D_{x\tiret y,v}$ and $D_{p\tiret q,v}$ are the divisors on $C_{v, t_0}$ associated to the paths $\varrho_{x \tiret y}$ and $\varrho_{p \tiret q}$ in \eqref{eq:DivisorsFromPaths}. 
(By the proof of \eqref{eq:HeightPairingAlternative}, the slopes of any solution $f$ to \eqref{eq:PoissonMCPairing} are given by $df = \varrho_{p\tiret q} - \chi_{p \tiret q}$.)

\begin{prop} \label{prop:terms}
Assumptions as in Theorem~\ref{thm:SeparateComplexGraph}, the following asymptotics hold uniformly for $(p,q,x,y) \in R_t$ as $t \in B^\ast$ converges to $\shy = (x,0)$ in $B^\hyb$.
\begin{itemize}
\item [$(i)$] The first term $A_1$ gives the sum of regularized pairings
\[
A_1 = (2 \pi)^{-1} \sum_{v\in V} \hp{C_{v, t_0}}{ D_{x\tiret y, v}, D_{p\tiret q, v}}' + o(1).
\]
\item [$(ii)$] The following equality holds (here, $\fell_t$ is the edge length function from \eqref{eq:piEdgeLengths} and $\mgr_t$ is the metric graph with edge length function $\ell_t = 2 \pi \fell_t$, see Theorem~\ref{thm:SeparateComplexGraph})
\[
\rho_{\fell_t} - \rmZ_{\fell_t}\rmM_{\fell_t}^{-1} \rmW_{\fell_t}= (2 \pi)^{-1} \hp{\mgr_t}{\bar x -\bar y, \bar p - \bar q}.
\] 
Moreover, the second term $A_2$ gives the sum of pairings
\[
A_2 = (2 \pi)^{-1} \hp{\mgr_t}{\bar x -\bar y, \bar p - \bar q} + (2 \pi)^{-1} \sum_{v \in V}  \hp{C_{v, t_0}}{D_{x\tiret y, v}^\perp , D_{p\tiret q, v}^\perp }' + o(1).
\]
\item[$(iii)$] The remaining terms $A_3$ and $A_4$ give the sums of regularized pairings
\begin{align*}
A_3 &=  - (2 \pi)^{-1} \sum_{v \in V} \hp{C_{v, t_0}}{D_{x\tiret y, v}^\perp, D_{p\tiret q, v}}' + o(1), \\
A_4 &= - (2 \pi)^{-1} \sum_{v \in V} \hp{C_{v, t_0}}{D_{x\tiret y, v}, D_{p\tiret q, v}^\perp}' + o(1).
\end{align*}

\end{itemize}
\end{prop}
\begin{proof}
$(i)$ Combining Proposition~\ref{prop:InvLem2.0} with the estimates \eqref{eq:VectorAs1} --\eqref{eq:VectorAs5}, one verifies that
\[
A_1 = \mathcal{R} - \mathcal{Z}^\surfup \Im(\Omega_S)^{-1} \mathcal{W}^\surfup + o(1)
\]
uniformly for $(p,q,x,y) \in R_t$ as $t \to \shy$ in $B^\hyb$. Hence the claim follows from Lemma~\ref{lem:RegularizedHeightFormula}, after considering the shape of the approximative matrix $\widehat \Ap(\rsf_{t_0}, \mgr_t; p,q,x,y)$ (see Section~\ref{ss:Ap}).

\smallskip

$(ii)$ By Proposition~\ref{prop:OrthoProjection}, the harmonic one-forms $\chi_{p\tiret q}$ and $\chi_{x\tiret y}$ on $G$ are given by
\[
\chi_{x\tiret y} = \sum_{j=1}^h (\rmM_\fell^{-1} \rmZ_\fell)_j \, \gamma_j, \qquad \chi_{p\tiret q} = \sum_{j=1}^h (\rmM_\fell^{-1} \rmW_\fell)_j \,  \gamma_j.
\]
The first equality is then an immediate consequence of the expression \eqref{eq:HeightPairingAlternative} for the height pairing on metric graphs. By Proposition~\ref{prop:InvLem2.0},
\begin{align*}
\widetilde \rmZ_\fell \Ap^{-1} \widetilde \rmW_\fell = \rmZ_\fell \rmM_\fell^{-1}\Big (\rmM_\fell + (\Im(B) \Im(\Omega_S)^{-1} \Im(B)^\transpose  - \Apc_G)   + o(1) \Big) \rmM_\fell^{-1} \rmW_\fell.
\end{align*}
Taking into account the shape of $\rmW_\fell$, $\rmZ_\fell$ (see \eqref{eq:VectorAs1}) and the inverse $\rmM_\fell^{-1}$ (see \eqref{eq:GraphInverse} for the asymptotics), it follows that
\[
\rmZ_\fell \rmM_\fell^{-1}  o(1) \rmM_\fell^{-1} \rmW_\fell  = o(1), \qquad t \to \shy,
\]
uniformly for $(p,q,x,y) \in R_t$. Hence the proof of (ii) is complete upon observing that 
\[
\rmZ_\fell \rmM_\fell^{-1} \Big (\Apc_G -  \Im(B) \Im(\Omega_S)^{-1}  \Im(B)^\transpose \Big )  \rmM_\fell^{-1} \rmW_\fell = (2 \pi)^{-1}\sum_{v \in V}  \hp{C_{v, t_0}}{D_{x\tiret y, v}^\perp, D_{p\tiret q, v}^\perp}',
\]
by the formula in Lemma~\ref{lem:RegularizedHeightFormula} and the shape of $\whAp(\rsf_{t_0}, \mgr_t; p,q,x,y)$ (see Section~\ref{ss:Ap}).

\smallskip

(iii) To prove the remaining statements, we conclude from Proposition~\ref{prop:InvLem2.0} and the estimates \eqref{eq:VectorAs5}--\eqref{eq:VectorAs5} together with \eqref{eq:GraphInverse} that, as $t \to \shy$, we have the uniform asymptotic
\begin{align*}
- A_3 &= \widetilde \rmZ_\fell  \Ap^{-1} \mathcal{W} \\
&= \rmZ_\fell  \Big(  \begin{matrix} \rmM_\fell^{-1} + \rmM_\fell^{-1} O(1) \rmM_\fell^{-1}  &  - \rmM_\fell^{-1} \Im(B) \Im(\Omega_S)^{-1} + \rmM_\fell^{-1} o(1) \end{matrix} \Big )  \mathcal{W}  \\
&= \rmZ_\fell \rmM_\fell^{-1} \mathcal{W}^\grup -  \rmZ_\fell \rmM_\fell^{-1} \Im(B) \Im(\Omega_S)^{-1} \mathcal{W}^\surfup + o(1).
\end{align*} 
By the shape of the approximative matrix and Lemma~\ref{lem:RegularizedHeightFormula}, we infer that
\[
 \rmZ_\fell \rmM_\fell^{-1} \mathcal{W}^\grup -  \rmZ_\fell \rmM_\fell^{-1} \Im(B) \Im(\Omega_S)^{-1}  \mathcal{W}\surfup = (2 \pi)^{-1} \sum_{v \in V} \hp{C_{v, t_0}}{D_{x\tiret y, v}^\perp, D_{p\tiret q, v}}'
\]
and hence we have described the asymptotic of $A_3$. The same reasoning applies to the remaining term $A_4$. 
\end{proof}
Theorem~\ref{thm:SeparateComplexGraph} now follows from Proposition~\ref{prop:AuxPropTHP}, Equation \eqref{eq:AnalyzableTerms} and Proposition~\ref{prop:terms}.
\qed

\subsection{Proof of Theorem~\ref{thm:tameness-height-pairing1}} \label{ss:ProofTamenessHP1}
Let $\rsf \to B$ be a versal family of the stable marked Riemann surface $(S_0, p_0, q_0, x_0, y_0)$. For every hybrid boundary stratum $D_\pi^\hyb$, and every point $\shy$ in the closure $\bar D_\pi^\hyb$ of $D_\pi^\hyb$, we have to verify the convergence \eqref{eq:StratumwiseTamenessHeightPairingSpelldOut}, as $t \in B^\ast$ converges tamely to $\shy$ in $B^\hyb$. Since all our considerations are local, it suffices to treat the case where $\shy = (l,0)$ has complex coordinate $s = 0$.

\medskip

Suppose first that $D_\pi^\hyb$ is associated to an ordered partition $\pi$ on the full edge set $E$ of the dual graph $S_0$, and $\shy$ is a point of the form $\shy = (l,0)$ in $D_\pi^\hyb$ (instead of the closure). Using the period map of the associated biextension mixed Hodge structure (see Section~\ref{sec:period}), the height pairing an be written as (see Lemma~\ref{lem:HeightPairingByPeriodMap})
\[
\hp{\rsf_t}{x_t - y_t , p_t-q_t} = 2\pi \Big (  \Im(\rho_t) - \Im(\rmW_t) \Im(\Omega_t)^{-1} \Im(\rmZ_t) \Big ).
\]
For a point $t \in B^\ast$, the imaginary part of the period matrix can be uniformly approximated by $\mathcal{A}(\rsf_{t_0}, \mgr_t; p_{t_0}, q_{t_0}, x_{t_0}, y_{t_0})$ (see Theorem~\ref{thm:ImaginaryAp}). Applying Theorem~\ref{thm:SeparateComplexGraph} (to the singleton region $R_t = \{(p_t, q_t, x_t, y_t)\}$), we conclude that
\begin{align*}
\hp{\rsf_t}{x_t - y_t , p_t-q_t}  = \hp{\mgr_t}{\bar x_0 -\bar y_0, \bar p_0 + \bar q_0} + \sum_{v\in V} \hp{C_{v, t}}{D^\cmc_{x_t \tiret y_t , v} , D^\cmc_{p_t \tiret q_t, v}}' + o(1),
\end{align*}
if $t \in B^\ast$ converges to $\shy$. The convergence \eqref{eq:StratumwiseTamenessHeightPairingSpelldOut} now follows from our results on metric graphs (see Theorem~\ref{thm:HeightPairingGraphs} and Theorem~\ref{thm:MainGraphDiscrete}), which can be applied since the metric graph $\mgr_t$ converges tamely to the underlying tropical curve $\curve^\trop_\shy$ of the hybrid curve $\rsf^\hyb_\shy$ in $\mgtrop{\grind{G}}$.

Namely, Theorem~\ref{thm:HeightPairingGraphs} implies that the term $\hp{\mgr_t}{\bar x_0 -\bar y_0, \bar p_0 + \bar q_0}$ can be split further into contributions from the graded minors. Theorem~\ref{thm:MainGraphDiscrete} implies that the coefficients of $D^\cmc_{x_t \tiret y_t, v}$ and $D^\cmc_{p_t\tiret q_t, v}$ on the metrized complex $\mc_t$ (recall that they are induced by $\mgr_t$ and hence depend on $t$) converge to the coefficients of the corresponding divisors $D^\smallcc_{p_{t} \tiret q_{t} ,v}$ and   $D^\smallcc_{x_t \tiret y_t ,v}$ appearing in the definition of the hybrid height pairing (see Section~\ref{ss:HeightPairingHybridCurves}). We conclude that the regularized height pairings $\hp{C_{v,t}}{D^\cmc_{x_t\tiret y_t, v} , D^\cmc_{p_t\tiret q_t, v}}'$ can be as well replaced by $\hp{C_{v,t}}{D^\smallcc_{x_t\tiret y_t, v} , D^\smallcc_{p_t\tiret q_t, v}}'$.

\medskip

 Suppose, more generally, that $\shy = (l,0)$ lies in the closure $\bar D_\pi^\hyb$ of an ordered partition $\pi = (\pi_1, \dots, \pi_r)$ of some subset $F \subseteq E$. 
 The proof of \eqref{eq:StratumwiseTamenessHeightPairingSpelldOut}  in this case uses the special case treated above and an appropriate induction. We only indicate the main steps here. Let $G_F=\rquot{G}{(E\setminus F)}$ be the graph obtained from $G = (V,E)$ by contracting all edges in $E \setminus F$. Using the ordered partition $\pi\rest F$ obtained from $\pi$ by restricting to $F$, we get a layered graph $(G_F, \pi)$. Consider also the graph $G-F = (V, E \setminus F)$, obtained from $G = (V,E)$ by deleting all edges in $F$. For each connected component $K$ of $G-F$, we get a stable marked Riemann surface denoted by $S_K$, by glueing the marked components $C_v$, $v \in K$, along the nodes $p^e_v$ for $v,e \in K$, and by marking with $p^e_v$ for $v\in K$ and $e\in F$.  Consider the corresponding deformation family $\rsf_K / B_K$.
 
We proceed as in the proof of the theorem previously exposed. By Theorem \ref{thm:HeightPairingGraphs}, we can approximate the height pairing $\hp{\rsf_t}{p_t - q_t, x_t - y_t}$ on a smooth fiber $\rsf_t$ close to $\rsf_\shy^\hyb$, by a graph part involving the metrization $\mgr_t$ and a complex part on the surfaces $C_{v,t_0}$. Applying Lemma~\ref{lem:JfunctionTame} to the ordered partition $(\pi_1, \dots, \pi_r, \pi_\fin = E \setminus F)$ of $E$, we can split the pairing on $\mgr_t$ into contributions from the graded minors $\Gamma_{1, t}, \dots \Gamma_{r,t}$ of $G_F$, and a contribution from the metrization $\mgr_{E \setminus F, t}$ of $G - F$. However, the case previously established of \eqref{eq:StratumwiseTamenessHeightPairingSpelldOut}, applied to the family $\rsf^\hyb_K / B^\hyb_K$,  implies that the contributions from $\mgr_{E \setminus F, t}$ and the contributions from the surfaces $C_{v, t_0}$, $v \in K$, together give regularized height pairings on certain fibers $\rsf_{K, t_K}$ in the families $\rsf_K / B_K$. Finally, taking a closer look at the definition of the hybrid log map $\loghyb_\pi$, the contributions from the minors $\Gamma_{1, t}, \dots, \Gamma_{r,t}$ and from the fibers $\rsf_{K, t_K}$ are precisely the terms in the sum from \eqref{eq:StratumwiseTamenessHeightPairingSpelldOut}, and the theorem follows.


\section{Uniform tameness of the height pairing: Proof of Theorem~\ref{thm:FinalTamenessHeightPairing}} \label{sec:HybridJFctAsymptotics}
In this section, we prove the uniform tameness of the height pairing as stated in Theorem~\ref{thm:FinalTamenessHeightPairing}.

\smallskip

In the preceding section,  we have already obtained the asymptotics of the height pairing $\hp{\rsf_t}{p_t - q_t, x_t -  y_t}$, where the smooth Riemann surface $\rsf_t$ degenerates to a hybrid curve, and $p_t$, $q_t$, $x_t$, $y_t$ are sections remaining entirely in the smooth parts of the respective curves (including for the limit stable curve underlying the hybrid curve). Proving the uniform tameness of the height pairing, however, corresponds to allowing that $p,q,x,y$ are arbitrary points on $\rsf_t$ (instead of smooth sections staying away from the appearing nodes), and that their limits are arbitrary points on the hybrid curve. This generalization is crucial in obtaining the asymptotics of the Arakelov Green function (and solutions to more general Poisson equations). 

\smallskip

The key auxiliary result in this context is Theorem~\ref{thm:HPNoSections} in Section~\ref{ss:UniformApproximationHeightPairing}.
This theorem splits the asymptotics of $\hp{\rsf_t}{p - q, x -  y}$ into two parts: a graph term originating from metric graphs and a complex term originating from Riemann surfaces. The separation is formalized by the notion of height pairing on metrized complexes from  Section~\ref{ss:AuxiliaryMetrizedComplexPairing} and allows to apply the results on metric graphs from Part~\ref{part:TropicalCurves}.  The formulation and proof of Theorem~\ref{thm:HPNoSections} occupy Sections~\ref{ss:UniformApproximationHeightPairing}--\ref{ss:UniformApproximationHeightPairingRevisited}. In Section~\ref{thm:FinalTamenessHeightPairing}, we finally prove Theorem~\ref{thm:FinalTamenessHeightPairing} by combining Theorem~\ref{lem:CrucialLemma} with the previous results on metric graphs.

\subsection{Uniform approximation of the height pairing} \label{ss:UniformApproximationHeightPairing}
Let $(S_0, q_1, \dots, q_\nmark)$ be a stable marked Riemann surface of genus $g$ with dual graph $G = (V,E)$. Consider a versal family $\rsf \to B$ for $S_0$ and let $\rsf^\hyb \to B^\hyb$ be the associated family of hybrid curves.

\smallskip

We {fix an adapted system of coordinates} on the family $\rsf \to B$. Recall in particular that the coordinates of $B \cong \Delta^N$, $N = 3g-3+\nmark$ are naturally decomposed into $\underline z = \underline z_E\times \underline z_{E^c}$, where $E^c = [N] \setminus E$. Moreover, the fibers $\rsf_t$, $t \in B$, can be decomposed as (see \eqref{eq:AdaptedCoordinates} for details)
\begin{equation} \label{eq:AdaptedDecompositionFinals}
 \rsf_{t}  = \bigsqcup_{v\in V} Y_{v,t} \sqcup \bigsqcup_{e\in E} W_{e, t}.
\end{equation}

We are interested in \emph{approximating the height pairing} $\hp{\rsf_t}{p - q, x - y}$ on a smooth fiber $\rsf_t$, $t \in B^\ast$, close to a hybrid curve $\curve$ with underlying stable Riemann surface $S_0$. Hence we fix a hybrid base point $\shy = (l,0)$ with complex coordinate $s = 0 \in B$ and consider the hybrid curve $\curve = \rsf_\shy^\hyb$.  The point $\shy$ lies in the hybrid stratum $D_\pi^\hyb$ of some ordered partition $\pi = (\pi_1, \dots, \pi_r) \in \Pifs(E)$ on the edge set $E$. Consider the associated stratumwise log map in the  commutative diagram
\[
\begin{tikzcd}
 \rsf^\ast \arrow[d]\arrow[r, "\loghybdiagpi{\pi}"] & \arrow[d] \rsf_\pi^\hyb\\
B^\ast\arrow[r, "\loghybdiagpi{\pi}"] &D_\pi^\hyb.
\end{tikzcd}
\]
On the base $B^\hyb$, the log map sends every point $t \in B^\ast$, representing a Riemann surface $\rsf_t$, to a point $\thy = \loghyb_\pi(t)$ in $D_\pi^\hyb$, representing a hybrid curve $\rsf^\hyb_\thy$ of type $(G, \pi)$. The log map on the family is given by an additional map $\psi^t_\thy \colon \rsf_t \to \rsf^\hyb_\thy$ between the fibers (see Section~\ref{sec:hybrid_log_map_rs}).

\medskip

We will approximate the height pairing on smooth fibers $\rsf_t$, $t \in B^\ast$, by the height pairings on the following metrized complexes. For $t \in B^\ast$, let $\mc_t$ be the metrized complex which is conformally equivalent to $\rsf_\thy^\hyb$, $\thy = \loghyb_\pi(t)$, and has edge lengths $\ell_t(e) = -\log|z_e(t)|$, $e \in E$. Equivalently, the metrized complex $\mc_t$, $t \in B^\ast$, is given by
\begin{itemize}
\item the metric graph $\mgr_t$, obtained from the dual graph $G = (V,E)$ of $S_0$ and the above edge lengths $\ell_t(e) = -\log|z_e(t)|$, $e \in E$,
\item and the Riemann surface components $C_{v,t}$, $v \in V$, associated to $\rsf_t$ by the adapted coordinates (see Section~\ref{sec:AdaptedCoordinatesFamily}).
\end{itemize}

The homotecies on intervals define a natural homeomorphism $\rsf^\hyb_\thy \cong \mc_t$. Here we slightly abuse the notation and use the letter $\rsf^\hyb_\thy$ also for the canonical metrized complex $\rsf^\hyb_\thy = \mccan_\thy$ representing the hybrid curve $\rsf^\hyb_\thy$ (the one having layerwise normalized edge lengths). The log map $\psi^t_\thy \colon \rsf_t \to \rsf^\hyb_\thy$ and the contraction from $\mc_t$ to $\mgr_t$ induce a sequence of maps
\begin{equation} \label{eq:MapSurfaceUnscaledGraph}
\begin{tikzcd}
\rsf_{t}  \arrow[r, "\varphi_t"] &  \mc_t \arrow[r, "\kappa_t"]  & \mgr_t 
\end{tikzcd} 
\end{equation}
for every $t \in B^\ast$. 
Abbreviating the notation, we denote the image of a point $y \in \rsf_t$ in $\mgr_t$ by $\bar y := (\kappa_t \circ \varphi_t)(y)$.  By an abuse of the notation, we use the same letter $\bar y$ for the image $ \varphi_t(y)$ of $y \in \rsf_t$ in the metrized complex $\mc_t$. The distinction between $\bar y$ in the metric graph and in the metrized complex will be clear from the context.

\medskip

The adapted coordinates on the family induce local coordinates $z^e_u$ around all attachment points $p^e_u(t)$ of the metrized complex $\mc_t$, $t \in B^\ast$. Using those parameters, the height pairing 
\begin{equation} \label{eq:RepeatMetrizedHeightPairingDefinition}
\hp{\mc_{t}}{\bar p - \bar q, \bar x - \bar y} = \hp{\mgr_{t}}{\bar p - \bar q, \bar x - \bar y} + \sum_{v\in V} \hp{C_{v,t}}{D_{\bar p \tiret \bar q,v}^\cmc, D_{\bar x \tiret \bar y ,v}^\cmc}'
\end{equation}
is defined for $p,q,x,y$ with $\{p, q \} \cap \{x,y\} = \varnothing$ on a smooth fiber $\rsf_t$, $t \in B^\ast$ (see Section~\ref{ss:AuxiliaryMetrizedComplexPairing}).

\smallskip

The next theorem approximates the height pairing on $\rsf_t$ in terms of the height pairing on the metrized complex $\mc_t$. Recall that $\varepsilon_\pi(p,q,x,y)$ denotes the correction term from \eqref{eq:CorrectionFourPoints}.

\begin{thm}[Uniform approximation of the height pairing] \label{thm:HPNoSections}
Suppose that $t \in B^\ast$ converges in $B^\hyb$ to the fixed hybrid point $\shy = (l,0) \in D_\pi^\hyb$ lying above $0\in B$. Then for $p,q,x,y \in \rsf_t$,
\[
\hp{\rsf_{t}}{ p - q, x  - y} =  \hp{\mc_{t}}{\bar p - \bar q, \bar x - \bar y} + \varepsilon_\pi(p,q,x,y) +  o(1)
\]
where the $o(1)$ term goes to zero uniformly for $p,q,x,y \in \rsf_t $. 
\end{thm}

\begin{remark} Note that the metrized complex $\mc_t$, the map $\varphi_t \colon \rsf_t \to \mc_t$ and the correction term $\varepsilon_\pi(p,q,x,y)$ are independent of the ordered partition  $\pi$ and defined solely in terms of the point $t \in B^\ast$. The uniform asymptotics in \eqref{lem:CrucialLemma} hence actually holds if $t \in B^\ast$ converges to $s = 0$ in  $B$ (as follows from the compactness properties of $B^\hyb$).
\end{remark}

As explained below, one can reduce the proof of Theorem~\ref{thm:HPNoSections} to the special case where the points $p$ and $x$ lie on well-behaved sections.

\smallskip

Namely, suppose that $p_t$ and $x_t$, $t \in B$, are \emph{adapted sections} of the family $\rsf \to B$. Hereby we mean that there is a vertex $v \in V$ such that $p_t$ lies in the region $Y_{v,t}$ for all $t \in B$ and $p_t$ only depends on the corresponding $3 \genusfunction(v) - 3 + \deg(v)$ coordinates in $B$, and similarly for $x_t$. Moreover, we require that $p_t \neq x_t$ for all $t \in B$. In particular, these sections naturally extend to two disjoint sections $p_\thy$ and $q_\thy$, $\thy \in B^\hyb$, of the hybrid family $\rsf^\hyb \to B^\hyb$.

(If the stable Riemann surface $S_0$ carries marked points $q_1, \dots, q_\nmark$, then any two different markings $q_i \neq q_j$ produce adapted sections.)  

\smallskip

For convenience, we restate Theorem~\ref{thm:HPNoSections} in the presence of two sections. Recall that $\varepsilon_\pi(q,y)$ denotes the correction term defined in \eqref{eq:CorrectionTwoPoints}.

\begin{thm}[Uniform approximation of the height pairing: special case with sections] \label{lem:CrucialLemma}
Let $x_t$ and $p_t$, $t \in B$, be two adapted sections of $\rsf \to B$. Suppose that $t \in B^\ast$ converges in $B^\hyb$ to the fixed hybrid point $\shy = (l,0) \in D_\pi^\hyb$ lying above $0\in B$.  Then for $q,y \in \rsf_t$,
\begin{equation} \label{eq:CrucialAsymptotics}
\hp{\rsf_{t}}{ p_t - q, x_t - y} =  \hp{\mc_{t}}{\bar p_t -  \bar q, \bar x_t - \bar y} + \varepsilon_\pi(q,y) + o(1)
\end{equation}
where the $o(1)$ term goes to zero uniformly for $q, y \in \rsf_t $. 
\end{thm}
As mentioned above, Theorem~\ref{thm:HPNoSections} follows directly from the special case in Theorem~\ref{lem:CrucialLemma}.
\begin{proof}[Proof of Theorem~\ref{thm:HPNoSections} from Theorem~\ref{lem:CrucialLemma}]
Close to $s=0$ in $B$, we can choose two auxiliary sections $p_t$ and $x_t$, $t \in B$, of the family $\rsf / B$ with the above adaptedness properties. The height pairings on $\rsf_t$ and $\mc_t$ are both bilinear. Hence it suffices to use the equality
 \[
 \hp{\rsf_{t}}{ p - q, x - y} = \hp{\rsf_{t}}{ p - p_t, x - x_t} + \hp{\rsf_{t}}{ p_t - q, x  - y} + \hp{\rsf_{t}}{ p_t - q, x  - x_t} + \hp{\rsf_{t}}{ p_t - q, x_t  - y}
 \]
and apply Theorem~\ref{lem:CrucialLemma}.
\end{proof}
It remains to prove Theorem~\ref{lem:CrucialLemma}. The idea of our method is to treat the points $q$ and $y$ as markings of $\rsf_t$, modifying if necessary the limit marked stable Riemann surface, and then apply the results from the preceding sections. On a larger scale, the uniformity of height pairing asymptotics is a consequence of the \emph{uniformity of the estimates in the nilpotent orbit theorem}.

\smallskip

In what follows, we work with deformation families of certain auxiliary marked Riemann surfaces.  Depending on the position of the points $q,y$ on $\rsf_t$ in the decomposition \eqref{eq:AdaptedDecompositionFinals}, we use different marked Riemann surfaces. We distinguish the following four cases (see Figure~\ref{fig:ProofPicCases}):

\begin{figure}[!t]
\centering
   \scalebox{.4}{\input{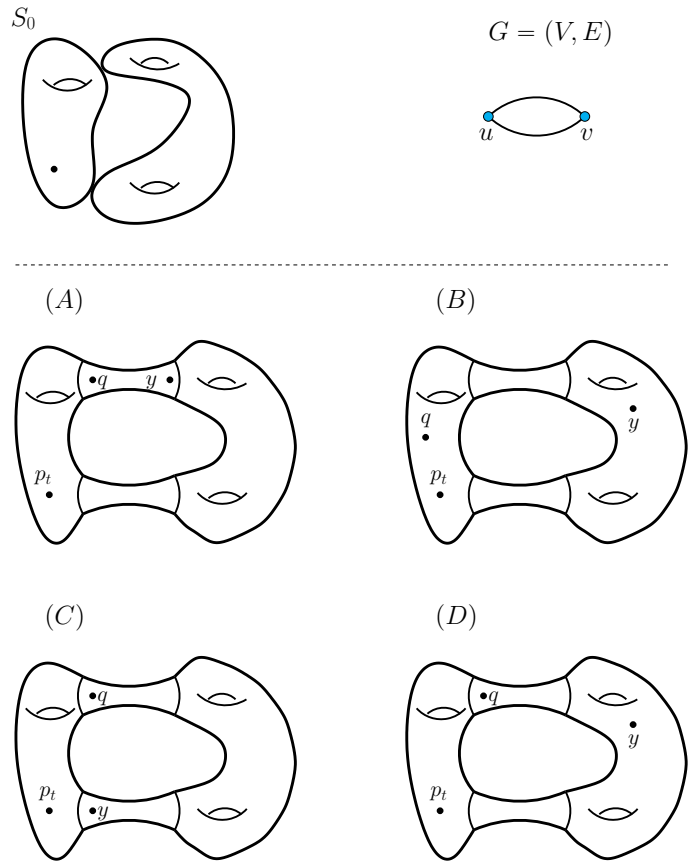}}
\caption{The stable Riemann surface $S_0$ and the four different positions for points $q,y$ on a nearby fiber $\rsf_t$, decomposed according to the adapted coordinates.}
\label{fig:ProofPicCases}
\end{figure}
\begin{itemize}
\setlength\itemsep{0.2 mm}
\item[(A)] $q,y$ belong to the same cylinder $W_{e,t}$ for some edge $e \in E$,
\item[(B)] $q,y,$ both belong to $\bigcup_{v \in V} Y_{v,t}$,
\item[(C)] $q \in W_{e,t}$ and $y \in W_{e',t}$ for two different cylinders $W_{e,t}$ and $W_{e',t}$, $e \neq e'$,
\item[(D)] $q \in \bigcup_{v \in V} Y_{v,t}$ and $y \in W_{e,t}$ for some $e \in E$ or vice versa.
\end{itemize}
In the next sections, we obtain the uniform asymptotic \eqref{lem:CrucialLemma} in each of these cases.

\subsection{Proof of Theorem~\ref{lem:CrucialLemma} in Case (A): both points in the same cylinder} Suppose that $q,y \in W_{e,t}$  for the same cylinder $W_{e,t} \subseteq \rsf_t$,
\[
W_{e,t} = \big  \{ (t, z^e_u, z^e_v) \, \st \, z^e_u z^e_v = z_e(t) \big  \} 
\]
associated to some edge $e = uv$. By symmetry, it suffices to prove the asymptotic \eqref{lem:CrucialLemma} in the case that $|z^e_u(q)| > |z^e_u(y)|$.
\subsubsection{The auxiliary stable Riemann surface $\widetilde S$ and its deformation space}
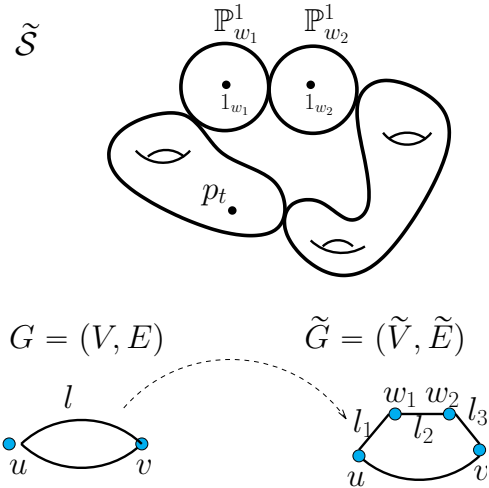
\begin{figure}[!h]
\centering
 \scalebox{.5}{\input{ProofPicSphereChain.tikz}}
\caption{The auxiliary stable Riemann surface $\~S$ and its dual graph $\~G = (\~V, \~E)$, together with the dual graph $G = (V,E)$ of $S_0$.}
\label{fig:ProofPicSphereChain}
\end{figure}

As outlined above, we utilize the deformation space of an \emph{auxiliary Riemann surface} with three marked points. Recall that $S_0$ denotes the stable marked Riemann surface underlying our fixed hybrid curve $\curve$. Consider the marked stable Riemann surface $\widetilde{S}$ obtained from $S_0$ by replacing the node $p_e$ between $C_u$ and $C_v$ with a chain of two Riemann spheres $\mathbb{P}_{w_1}$, $\mathbb{P}_{w_2}$, each of them accommodating a marked point denoted by  $1_{w_1}$ and $1_{w_2}$, respectively (see Figure~\ref{fig:ProofPicSphereChain}). Compared to the graph $G = (V,E)$ of $S_0$, the underlying graph $\widetilde G = (\widetilde V, \widetilde E)$ of $\widetilde S$ has two new vertices $\widetilde V = V \cup \{w_1, w_2\}$ (where $w_i$ represents $\mathbb{P}_{w_i}$ , $i=1,2$), and the edge $e=uv$ has been subdivided into the three edges $e_1 = u w_1$, $e_2 = w_1 w_2$ and $e_3 = w_2 v$.

Moreover, we fix global coordinates $\zeta_{w_i}$ on the spheres $\mathbb{P}_{w_i}$, $i=1,2$ such that
\begin{itemize}
\item the marking $1_{w_i}$ on  $\mathbb{P}_{w_i}$ corresponds to  $\zeta_{w_i}(1_{w_i}) = 1$,
\item the node $p_{e_2} = \mathbb{P}_{w_1} \cap \mathbb{P}_{w_2}$ corresponds to $\zeta_{w_1} = \infty$ and $\zeta_{w_2} = 0$,
\item the node $p_{e_1} = C_u \cap \mathbb{P}_{w_1}$ corresponds to $\zeta_{w_1} = 0$ and $z^e_u = 0$,
\item and the node $p_{e_3} = C_v \cap \mathbb{P}_{w_1}$ corresponds to $\zeta_{w_2} = \infty$ and $z^e_v = 0$
\end{itemize}
Next we compare the deformation families $\rsf \to B$ and $\~ \rsf \to \~ B$ for $S_0$ and $\~ S$. Recall that we have fixed an adapted coordinate system for the family $\rsf \to B$. This means that $\rsf \to B$ together with its adapted coordinates is obtained by plumbing construction in Section~\ref{sec:plumbing}.  Recall that, for every marked component $C_w$, $w \in V$ of $S_0$, we take a deformation family $\rsf_w \to B_w$, then fix local coordinates $z^{e'}_w$, $e'\in E$, around the markings corresponding to nodes in $S_0$, cut out the unit discs $\abs{z^{e'}_w} < 1$ from the $C_w$'s to get the surfaces with boundaries $Y_{w,t}$, and then plug the cylinders $W_{e',t}$ instead.

\smallskip

Based on the plumbing construction of $\rsf \to B$, we can provide a plumbing construction for the deformation family $\~\rsf \to \~B$ for $\~S$ such that \emph{the adapted coordinate systems are compatible}. Note first that, for the base spaces $B$ and $\~B$,
\begin{equation} \label{eq:ComparisonDeformationBases}
B = \prod_{v \in V} B_v \times \prod_{e' \in E} \Delta_{e'}, \qquad \~B = \prod_{v \in V} B_v \times \prod_{e' \in E \setminus\{e\}} \Delta_{e'} \times \Delta_{e_1} \times \Delta_{e_2} \times \Delta_{e_3},
\end{equation}
since for any original vertex $w \in V$, the associated marked component $C_w$ in $S_0$ and $\~S$ is the same, and the additional marked spheres in $\~ S$ do not support a parameter.  To construct the family $\~\rsf \to \~B$, we now start with the \emph{same deformation families} $\rsf_v \to B_v$ as for $\rsf \to B$, choose \emph{the same local coordinates} $z^{e'}_w$ around markings representing nodes, and cut out the unit discs $\abs{z^{e'}_w} < 1$. To find adapted coordinates for the sphere $\mathbb{P}_{w_i}$, $i=1,2$ we use the global coordinate $\zeta_{w_i}$ and cut out the discs $\abs{\zeta_{w_i}} < 1$ and $\abs{\zeta_{w_i}^{-1}} < 1$. Finally, we  plug cylinders $W_{e', \tilde{t}}$, $e' \in \~E$ to obtain the family $\~{\rsf} \to \~B$. The sections $p_t$ and $x_t$ of $\rsf \to B$ extend naturally to sections of $\~\rsf \to \~B$. By an abuse of notation, we denote the latter again by $p_t$ and $x_t$ (this is justified, since the image is independent of the coordinates of $\~B$ which are not already in $B$).

\medskip

Fix now a point $t \in B^\ast$ and suppose that  $q,y$ belong to the cylinder $ W_{e,t}$ with $|z^e_u(q)| > |z^e_u(y)|$. We can associate with this data a point $\tilde t \in \~B$, given in the decomposition \eqref{eq:ComparisonDeformationBases} by
\begin{equation} \label{eq:MapDeformationSpaces}
(t, q,y) \,  \mapsto \,  \tilde t = \Big ( (\tilde t_w = t_w)_{w \in V}, \,  (\tilde t_{e'} = t_{e'})_{\substack{e' \in E\\ e' \neq e}}, \,  \tilde t_{e_1} = z^e_u(q), \,  \tilde{t}_{e_2} = \frac{z^e_u(y)}{z^e_u(q)},  \, \tilde{t}_{e_3} = z^e_v(y)  \Big ).
\end{equation}
We can view the data $(\rsf_t, p_t, q, y)$ as a marked Riemann surface. As we show in the following lemma, this will be isomorphic to the marked Riemann surface $(\~\rsf_{\tilde{t}}, p_t, 1_{w_1}, 1_{w_2})$. Moreover, we get the following equality.

\begin{lem} \label{lem:biholomirphism-marking}
Notations as above, let $t \in B^\ast$. Suppose that $q,y \in W_{e,t}$ for an edge $e$ and $|z^e_u(q)| > |z^e_u(y)|$. Then, the marked Riemann surfaces $(\rsf_t, p_t, q, y)$ and $(\~\rsf_{\tilde{t}}, p_t, 1_{w_1}, 1_{w_2})$ are isomorphic, and we have
\[
\hp{\rsf_t}{p_t - q, y - x_t} = \hp{\widetilde{\rsf}_{\tilde t}}{p_t - 1_{w_1}, 1_{w_2} - x_t} 
\]
\end{lem}
\begin{proof} It suffices to construct a biholomorphic map $\phi \colon \rsf_t \to \~\rsf_{\tilde t}$ respecting the marking. The adapted systems of coordinates give decompositions
\[
\rsf_t = \bigsqcup_{w \in V} Y_{w,t} \sqcup \bigsqcup_{e' \in E} W_{e',t}, \qquad \rsf_{\tilde {t}} = \bigsqcup_{w \in V} Y_{v,{\tilde{t}}} \sqcup \bigsqcup_{e' \in \~E} W_{e',\tilde {t}},
\]
By construction, $Y_{w,t} = Y_{w, \tilde{t}}$ for $w \in V$ and $W_{e', \tilde{t}} = W_{e', t}$ for $e' \neq e$ in $E$. Setting $\phi$ as the identity on these parts, it remains to construct a map from $W_{e, t}$ to $W_{e_1, \tilde t} \cup W_{e_2, \tilde t} \cup W_{e_3, \tilde t}$. We split $W_{e,t}$ as $W_{e,t} = W_{1,t} \cup W_{2,t} \cup W_{3,t}$ with
\begin{align*}
W_{1,t} &= \bigl\{ (t, z^e_u, z^e_v) \in W_{e,t} \, \st \, |z^e_u| \ge |z^e_u(q)|  \bigr\}, && W_{3,t} = \bigl\{ (t, z^e_u, z^e_v) \in W_{e,t} \,\st \,|z^e_u(y)| \ge |z^e_u|   \bigr\} \\
W_{2,t} &= \bigl\{ (t, z^e_u, z^e_v) \in W_{e,t} \, \st \,|z^e_u(q)| \ge |z^e_u| \ge |z^e_u(y)|  \bigr\}.
\end{align*}
We send $W_{i, t}$ to $W_{e_i, \tilde{t}}$ by setting
\begin{align*}
\phi(t, z^e_u, z^e_v) &= \Big (\tilde {t}, z^e_u, \frac{z^e_u(q)}{z^e_u} \Big)  \text{ on } W_{1,t} &&
\phi(t, z^e_u, z^e_v) = \Big  (\tilde {t}, \frac{z^e_u}{z^e_u(q)}, \frac{z^e_u(y)}{z^e_u} \Big )  \text{ on } W_{2,t} \\
\phi(t, z^e_u, z^e_v) &=\Big  (\tilde {t},  \frac{z^e_v(y)}{z^e_v}, z^e_v \Big )  \text{ on } W_{3,t}.
\end{align*}
As is easily checked, this map is biholomorphic and sends the points $p_t, q,y,x_t$ to points $p_t, 1_{w_1}, 1_{w_2}, x_t$, respectively.
\end{proof} 

\begin{figure}[!t]
\centering
 \scalebox{.5}{\input{ProofPicSubcases.tikz}}
\caption{The decomposition of the cylinder $W_{e,t}$ and the position of the points $q,y$ on $W_{e,t}$ for the six subcases (a)--(f). The term $\varepsilon(q,y)$ is non-zero precisely for the cases (e) and (f).}
\label{fig:ProofPicSubcases}
\end{figure}
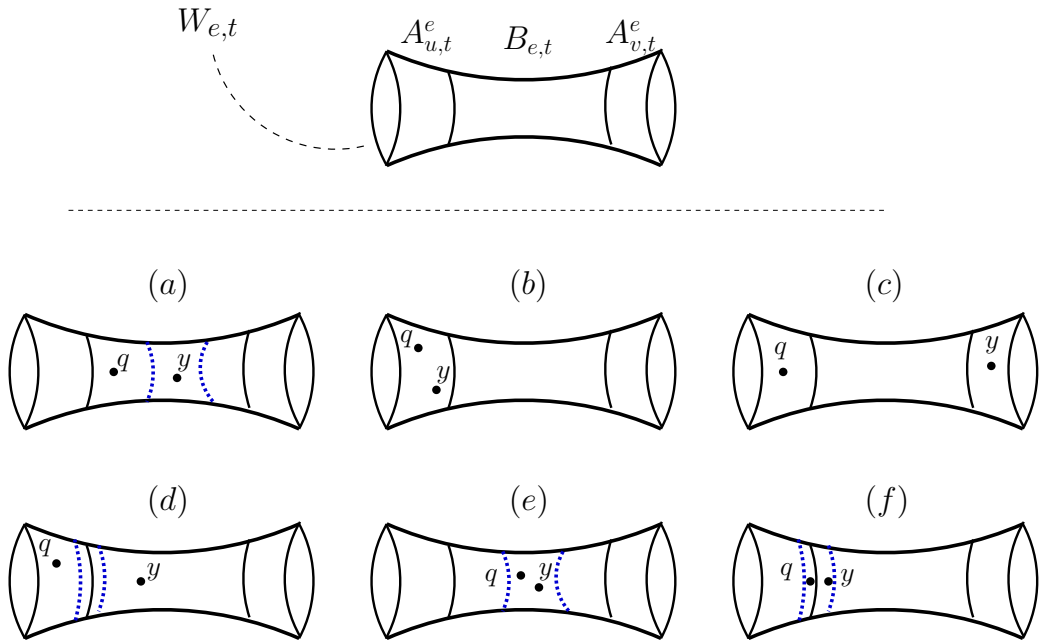

In order to prove the uniform asymptotics \eqref{lem:CrucialLemma}, we split the analysis in six subcases, according to the relative position of $q,y$ in the cylinder $W_{e,t}$ (see Figure~\ref{fig:ProofPicSubcases}). Recall that, in addition, we always assume that $|z^e_u(q)| > |z^e_u(y)|$. By the symmetry of the appearing expressions in $z^e_u$ and $z^e_v$, one deduces the asymptotics \eqref{lem:CrucialLemma} for all points $q,y \in W_{e,t}$ from the following six cases.

\begin{itemize}
\setlength\itemsep{0.2 mm}
\item[(a)] $q,y \in B_{e,t}$ and $\frac{|z^e_u(y)|}{|z^e_u(q)|} < \varrho_e(t)$
\item[(b)] $q, y \in A^e_{u,t}$
\item[(c)] $q \in A^e_{u,t}$ and $y \in A^e_{v,t}$
\item[(d)] $q \in A^e_{u,t}, y \in B_{e,t}$ and either $|z^e_u(q)| > {\varrho_e(t)}^{1/2}$ or ${\varrho_e(t)}^{3/2} > |z^e_u(y)|$
\item[(e)] $q,y \in B_{e,t}$ and $\frac{|z^e_u(y)|}{|z^e_u(q)|} > \varrho_e(t)$
\item[(f)] $q \in A^e_{u,t}, y \in B_{e,t}$ and ${\varrho_e(t)}^{1/2} > |z^e_u(q)| >  |z^e_u(y)| > {\varrho_e(t)}^{3/2}$. 
\end{itemize}

\subsubsection{Interlude: the proof of Theorem~\ref{thm:tameness-height-pairing1} in a nutshell} \label{ss:nutshell}
Before proceeding further with the proof of the uniform tameness of the height pairing, we briefly recall the method of proof of the (non-uniform) tameness of the height pairing (see Theorem~\ref{thm:tameness-height-pairing1}).

Recall that we considered a deformation family $\rsf \to B$ of a stable Riemann surface $S_0$ with four markings $p_0, q_0, x_0, y_0$ and dual graph $G= (V,E)$, equipped with adapted coordinates. The marked points produce four holomorphic, disjoint sections $p_t$, $q_t$, $x_t$ and $y_t$, $t \in B$, of the family $\rsf \to B$. We are interested in describing the height pairing $\hp{\rsf_t}{p_t - q_t, x_t - y_t }$ close to a hybrid point $\shy = (l, 0)$, lying over $0 \in B$.

\begin{itemize}
\item [(i)] In Section~\ref{sec:period}, we have studied the variation of mixed Hodge structures associated with the family of marked Riemann surfaces $\rsf \to B$.  The imaginary part $\Im(\widehat \Omega_t)$, $t \in B^\ast$, of the associated period matrix is defined over the original base $B$ (instead of the universal cover). By Lemma~\ref{lem:HeightPairingByPeriodMap}, the height pairing can be computed from $\Im(\widehat\Omega_t)$ as 
\[
	\hp{\rsf_t}{ p_t - q_t, x_t - y_t} = F \big ( \Im(\bipm_{t})  \big ),
\]
where $F$ is the function sending a real $(g+1)\times(g+1)$ matrix $\widehat{\mathcal{B}}$ of the form
\[
\widehat{\mathcal{B}}  = \begin{pmatrix} \mathcal{B} & \mathcal{B}_\rmW \\ \mathcal{B}_\rmZ & \mathcal{B}_\varrho  \end{pmatrix} , \qquad \qquad \mathcal{B} \in \R^{g \times g}, \mathcal{B}_\rmW \in \R^{g \times 1}, \mathcal{B}_\rmZ \in \R^{1 \times g}, \mathcal{B}_\varrho \in \R
\]
with an invertible submatrix $\mathcal{B} \in \R^{g \times g}$ to
\[
F (\widehat{\mathcal{B}}) :=  2 \pi \Big (  \mathcal{B}_\varrho - \mathcal{B}_\rmZ  \mathcal{B}^{-1} \mathcal{B}_\rmW \Big ).
\]

\item [(ii)] The nilpotent orbit theorem provides an approximation (see Theorem~\ref{thm:ImaginaryAp})
\[
\Im(\bipm_{t}) = \widehat \Ap(\rsf_{t_0}, \mgr_t; p_{t_0}, q_{t_0}, x_{t_0}, y_{t_0}) + o(1), \qquad t \in B^\ast,
\]
where the error term goes to zero uniformly if the coordinates $z_e(t) \to 0$ for all $e \in E$.

The entries of the approximative matrix $\widehat \Ap$ are (regularized) pairings between paths and one-forms on the metric graph $\mgr_t$ (the metric realization of $G$ equipped with lengths $\ell_t(e) = -\log|z_e(t)|$) and the stable Riemann surface $\rsf_{t_0}$ (here $t_0\in B$ is obtained from $t$ by setting $z_e(t_0) = 0$ for all $e \in E$).

\item [(iii)] As follows from Theorem~\ref{thm:SeparateComplexGraph}, the approximation on the height pairing is asymptotically equal to
\[
F \Big ( \widehat \Ap(\rsf_{t_0}, \mgr_t; p_{t_0}, q_{t_0}, x_{t_0}, y_{t_0}) + o(1) \Big ) = \hp{\mc_t}{x_{t_0} - y_{t_0}, p_{t_0} - q_{t_0}} + o(1), \qquad \text{ as $t \to \shy$},
\]
where $\hp{\mc_t}{\cdot, \cdot} $ is the height pairing on the metrized complex $\mc_t$ associated to $\rsf_{t_0}$ and $\mgr_t$. Altogether, we conclude that
\begin{align} \label{eq:SchemeHeightPairingProof} \begin{split}
\hp{\rsf_t}{x_t - y_t, p_t - q_t} &= F \big ( \Im(\bipm_{t}) \big ) = F \big ( \widehat \Ap(\rsf_{t_0}, \mgr_t; p_{t_0}, q_{t_0}, x_{t_0}, y_{t_0}) + o(1) \big )  \\ &= \hp{\mc_t}{x_{t_0} - y_{t_0}, p_{t_0} - q_{t_0}} + o(1) \end{split}
\end{align}
\item [(iv)]  In Section~\ref{ss:ProofTamenessHP1}, we finally derive Theorem~\ref{thm:SeparateComplexGraph} by applying our results on degenerations of metric graphs to study the height pairing on the metrized complex $\mc_t$.
\end{itemize}

\smallskip 

After these preparations, we give the proof of the asymptotic \eqref{eq:CrucialAsymptotics} in the above six cases. Vaguely speaking, we apply the steps (i)--(iii) for the deformation family $\~{\rsf} \to \~{B}$. However, some modifications are needed in the respective cases. We begin with the cases (a)--(d). Note that these are precisely the ones where $\varepsilon_\pi (q, y) \equiv 0$

\subsubsection{Case \emph{(a)}} In this case, we can directly apply the above procedure in the deformation family $\~\rsf \to \~B$. More precisely, if $t \in B^\ast$ converges to the complex part $s = 0$ of the point $\shy$  in the original disc $B$, and $q, y \in W_{e,t}$ are points in the subregion (a), then the points $\tilde t$ in $\~B$ satisfy that $z_{e'}(\tilde t) \to 0$ for all $e' \in \~E$. By definition of this subregion and the properties of $\varrho_e(t)$, this convergence is uniform in $q,y$. By Theorem~\ref{thm:ImaginaryAp}, the imaginary part $\Im(\bipm_{\~t})$ of the period matrix on $\~B$ is equal to $ \Im(\bipm_{\~t}) = \widehat \Ap + o(1)$, where $\widehat \Ap$ the respective approximate matrix and the $o(1)$ term goes to zero uniformly for $q,y$ in the subregion (a). Applying the steps (i)--(iii) in Section~\ref{ss:nutshell} in the family $\~{\rsf} \to \~{B}$, we arrive at the approximation
\[
\hp{\rsf_t}{p_t - q, y - x_t} = \hp{\widetilde{\rsf}_{\tilde t}}{p_t - 1_{w_1}, 1_{w_2} - x_t}  = \hp{\~\mc_{\tilde t}}{p_{t} - 1_{w_1}, 1_{w_2} - {x_{t}}} + o(1)
\]
with a uniform error term. The metrized complex $\~\mc_{\tilde t}$ here is constructed from $\mc_t$ by adding two Riemann spheres $\mathbb{P}_{w_1}$, $\mathbb{P}_{w_2}$  with marked points $1_{w_1}$, $1_{w_2}$ at the positions $-\log|z^e_u(q)|$ and $-\log|z^e_u(y)|$ on the edge $e$. The regularization on these Riemann spheres is in global coordinates $\zeta_{w_i}$, $i=1,2$, such that the markings $1_{w_i}$ have coordinates equal to one, $\zeta_{w_i}(1_{w_i}) = 1$.

However, by Proposition~\ref{prop:explicit_regularization_sphere}, the regularized pairings coming from the sphere $\mathbb{P}_{w_i}$, $i=1,2$, are zero. Hence we recover precisely  the terms  in \eqref{eq:CrucialAsymptotics}, that is,
\[
\hp{\~\mc_{\~t}}{p_{t} - 1_{w_1}, 1_{w_2} - {x_{t}}} = \hp{\mc_{t}}{\bar p_t - \bar q , \bar y  -  \bar x_t }.
\]

\subsubsection{Case \emph{(b)}}
Suppose that both $q$ and $y$ belong to $A^e_{u,t}$. Note that in this case, the corresponding points $\tilde t$ do not necessarily satisfy $\tilde t_{e_1} \to 0 $ and $\tilde t_{e_2} \to 0 $, that is, they might not approach the stratum
\[
	D_{\~E} = \big \{ \tilde s \in \~B \, \st \,  \text{$\tilde s_{e'} = 0$ for all $e' \in \~E$} \big \}
\]
in $\~B$ if $t \to \shy$. In particular, we cannot proceed as above and apply Theorem~\ref{thm:ImaginaryAp} directly in the deformation space $\~B$. However, note that $\tilde t_{e'} \to 0$ for all $e' \in E \cup \{e_3\}$. By the properties of $\varrho_{e'}(t)$, this convergence is moreover uniform as $t \to 0$ in $B$ and $q,y$ belong to the subregion (b). Hence, we can apply the nilpotent orbit theorem for the projection to the closed stratum
\[
	D_{E \cup \{e_3\}} =  \big \{ \tilde s \in \~B \, \st \,  \text{$\tilde s_{e} = 0$ for all $e \in E \cup \{e_3\}$} \big \}
\]
in $\~B$. More formally, we fix the variables $z_{e_1 }$ and $z_{e_2}$ as $z_{e_1} :=z^e_u(q)$ and $z_{e_2} := z^e_u(y) / z^e_u(q)$, which gives a variation of Hodge structures over a disc of lower dimension. Applying the nilpotent orbit theorem to this new system, and after some calculations that we omit, we end up with the uniform approximation
\[
\Im(\bipm_{\~t}) = \widehat \Ap \big (\rsf_{t_0}, \mgr_{\tilde t} \, ; p_{t} , q, y, x_{t},  \big ) + o(1).
\]
Here, $\mgr_{\tilde t}$ is the metric graph obtained from $\mgr_t$ by shortening the length of the edge $e$ to $\tilde \ell(e) = -\log\abs{z^e_u(q)}$ and $\rsf_{t_0}$ is the fiber in $\rsf \to B$ over the point $t_0$, obtained from $t$ by setting $z_e(t_0) = 0$ for all $e \in E$. The approximative period matrix $\widehat \Ap \big (  \rsf_{t_0}, \mgr_{\tilde t} \, ; p_{t_0} , q, y, x_{t_0}, \big )$ is obtained by viewing the points $p_{t} , q, y, x_{t}$ on the stable surface $\rsf_{t_0}$, and combining (regularized) pairings between paths and one-forms on $\rsf_{t_0}$ and $\mgr_{\tilde t}$. The regularization on the component $C_{t_0, u}$ of $\rsf_{t_0}$ at the node $p^{e}_u$ is taken in the local coordinate $\tilde z^e_u = z^e_u / z^e_u(q)$.  Due to the uniformity of the constants in the nilpotent orbit theorem (e.g., \cite{CKS}), the $o(1)$ term goes to zero uniformly as $t \to 0$ in $B$ and $q,y$ belong to the subregion (b).

\smallskip

The following lemma allows to replace the metric graph $\mgr_{\tilde t}$ by the original metric graph $\mgr_t$.

\smallskip

\begin{lem} \label{lem:Magic}
Let $S$ be a stable Riemann surface with dual graph $G = (V,E)$, together with fixed local parameters $z^e_u$, $e \sim u$, around the nodal points $p^e_u$ on the components $C_u$, $u \in V$, of $S$. Let $l\colon E \to (0, + \infty)$ be an edge length function and $\mgr$ the metric graph associated to $(G,l)$.

Consider a path $\gamma$ on $S$ connecting two smooth points $x, y \in S$, and a meromorphic one-form $\omega$ on  $S$ with real-valued residue divisor. Denote by $\gamma_v$ and $\omega_v$ the restrictions of $\gamma$ and $\omega$ to the components $C_v$, $v \in V$, of $S$. Let $\bar \gamma$ be the induced path on $G$ and $\bar \omega$ the induced one-form on the graph $G$ (obtained by taking the residues of $\omega$ along nodes).

Suppose that, instead of $z^e_u$, the regularization at some node $p^e_u$, $e\in E$, is taken in the new local coordinate $\tilde z^e_u = a z^ e_u$ for some $a \in \C^*$. Denote by $\~\mgr$ the metric graph obtained from $\mgr$ by changing the length of $e$ to $\tilde l(e) = l(e) - (2\pi)^{-1} \log|a|$. Then, we have 

\[
\sum_{v \in V}  \hp{C_v}{\gamma_v, \omega_v}''+  \hp{\~\mgr}{\bar \gamma, \bar \omega}  = \sum_{v \in V}  \hp{C_v}{\gamma_v, \omega_v}' +  \hp{\mgr}{\bar \gamma, \bar \omega},
\]
where $\hp{C_v}{\gamma_v, \omega_v}''$, $v\in V$, denotes the regularized height pairing using the new parameter $\~z^e_u$ instead of $z^e_u$ at $p^e_u$ (the other coordinates remain unchanged).
\end{lem}
\begin{proof} Clearly, the change in the regularized pairing on the component $C_u$ is given by
\[
\hp{C_u}{\gamma_u, \omega_u}'' = \hp{C_u}{\gamma_u, \omega_u}' +(2\pi)^{-1}  \log|a| \bar \gamma (e) \bar \omega (e).
\]
The claim immediately follows.
\end{proof}
Since all entries of the approximative period matrices have the form in Lemma~\ref{lem:Magic}, we infer the equality
\[
\widehat \Ap \big (\rsf_{t_0}, \mgr_{\~t} \, ; p_{t} ,q,y, x_{t}  \big ) = \widehat \Ap  \big (\rsf_{t_0}, \mgr_{t} ; p_{t} , q, y, x_{t} \big ),
\]
where now the approximative period is computed by regularization in the original local coordinate $z^e_u$ on $\rsf_{t_0}$ and using the original metric graph $\mgr_t$. Hence Theorem~\ref{thm:SeparateComplexGraph} (see also \eqref{eq:SchemeHeightPairingProof}) implies that
\[
\hp{\rsf_t}{p_t - q, y - x_t} = F \Big ( \widehat \Ap \bigl ( \rsf_{t_0}, \mgr_{t}; p_{t} , q, y, x_{t} \bigr) + o(1)  \Big )  = \hp{\mc_t}{\bar p_t - \bar q, \bar y - \bar x_t} + o(1).
\]
This finishes the proof of Theorem~\ref{lem:CrucialLemma} in the case (b).

\subsubsection{Case \emph{(c)}}

The case (c) can be treated analogous to (b), by projecting to the stratum
\[
	\~D_{E \cup \{e_2\} } = \bigl \{ \tilde s \in \~B \, \st \,  \text{$\tilde s_{e'} = 0$ for all $e' \in E$ and $\tilde s_{e_2} = 0$} \bigr \}
\]
in $\~B$ and again applying Lemma~\ref{lem:Magic}. 

\subsubsection{Case \emph{(d)}}
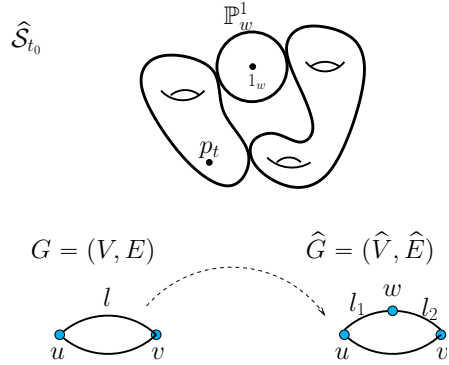
\begin{figure}[!t]
\centering
 \scalebox{.4}{\input{ProofPicAuxiliaryCase.tikz}}
\caption{The stable Riemann surface $\widehat \rsf_{t_0}$ and its dual graph $\widehat G = (\widehat V, \widehat E)$, together with the dual graph $G = (V,E)$ of $S_0$.}
\label{fig:ProofPicAuxiliaryCase}
\end{figure}

Suppose that $q \in A^e_{u,t}, y \in B_{e,t}$, and either $|z^e_u(q)| > {\varrho_e(t)}^{1/2}$ or ${\varrho_e(t)}^{3/2} > |z^e_u(y)|$. By the properties of $\varrho_e(t)$,  we conclude that
\[
\tilde t_{e_2} = \frac{z^e_u(y)}{z^e_u(q)} \to 0, \qquad \tilde t_{e_3} = z^e_v(y) \to 0,
\]
uniformly as $t \to s$ in $B$ and $q,y$ belong to the subregion (d) of $\rsf_t$. In particular, we can apply the nilpotent orbit theorem to the projection to the stratum
\[
	\~D_{E \cup \{e_2, e_3\} } = \bigl \{ \tilde s \in \~B \, \st \,  \text{$\tilde s_{e'} = 0$ for all $e' \in E \cup \{e_2, e_3\}$} \bigr \}
\]
in $\~B$. As is readily verified (using again Lemma~\ref{lem:Magic}), we get the following approximation
\[
\Im(\bipm_{\~t})  = \mathcal{A} \big (\widehat{\rsf}_{t}, \widehat{\mgr}_{\tilde t} \, ; p_{t} , q, 1_w, x_{t} \big ) + o(1),
\]
where the $o(1)$ term goes to zero uniformly. Here, the stable Riemann surface $\widehat{\rsf}_{t_0}$ is obtained from $\rsf_{t_0}$ by adding a Riemann sphere $\mathbb{P}_w$ with a marked point $1_w$ at the node $p_e$. The underlying graph $\widehat G$ has a new vertex $w$, and the edge $e=uv$ is subdivided into two new edges $e_1 = uw$ and  $e_2= wv$ (see Figure~\ref{fig:ProofPicAuxiliaryCase}). The metric graph $\widehat \mgr_{\tilde t}$ is obtained from $\mgr_t$ by setting $\ell_{e_1} = -\log|z^e_u(y)|$ and $\ell_{e_2} = -\log|z^e_v(y)|$. The regularization on the components of $S_0$ is taken in the original local coordinates. On the Riemann sphere $\mathbb{P}_w$, the regularization is taken in  a global coordinate $\zeta$ such that $\zeta(1_w) = 1$, and $\zeta(p^{e_1}_w) =0, \zeta(p^{e_2}_w)=\infty$. Moreover, $q$ is seen as a point on the component of $\widehat{\rsf}_{t_0}$ corresponding to the vertex $u$.

By Theorem~\ref{thm:SeparateComplexGraph} (see \eqref{eq:SchemeHeightPairingProof}), we get
\[
\hp{\rsf_t}{p_t - q, y - x_t} = \hp{\widehat{\mc}_{\~t}}{ p_{t} - q,  1_w -  x_{t} } + o(1)
\]
for the metrized complex $\widehat{\mc}_{\~t}$ obtained from the stable Riemann surface $\widehat{\rsf}_{t_0}$ and the metric graph $\widehat{\mgr}_{\tilde t} $ (strictly speaking, Theorem~\ref{thm:SeparateComplexGraph} is applied in the deformation family $\widehat \rsf \to \widehat B$ of the marked surface $\widehat S_0 = \widehat \rsf_0$). However, $\mgr_t$ is isomorphic to $\widehat{\mgr}_{\tilde t}$ as a metric graph, and the term in the regularized height pairing coming from the sphere $\mathbb{P}_w$ vanishes by Proposition~\ref{prop:explicit_regularization_sphere}. Hence we infer that
\[
\hp{\widehat{\mc}_{\tilde t}}{ p_{t} - q,  1_w -  x_{t} } = \hp{\mc_{t}}{\bar p_t - \bar q, \bar x_t - \bar y} ,
\] 
which finishes the proof in this case.

\subsubsection{Cases \emph{(e)} and \emph{(f)}} 
It remains to describe the asymptotic in the cases (e) and (f). These are precisely the ones where the additional term $\varepsilon_\pi(q,y)$ is non-zero.

\smallskip

If $t$ converges to $s=0$ in $B$ and $q,y$ belong to one of the subregions (e) or (f), then
\[
\tilde t_{e'} \to 0 \text{ for all $e' \in E \cup \{e_1, e_3\}$}
\]
uniformly. Moreover, we have the uniform lower estimate
\[
|\tilde t_{e_2}| \ge \varrho_e(t).
\]
Following the previous approach, we can use the nilpotent orbit theorem for the projection to the stratum
\[
	\~D_{E \cup \{e_1, e_3\} } = \bigl \{ \tilde s \in \~B \,\st \,  \text{ for all $e' \in E \cup \{e_1, e_3\}$} \bigr \}
\]
in the base $\~B$. After applying Lemma~\ref{lem:Magic} once again, we get a uniform approximation
\[
\Im(\bipm_{\~t})  = \widehat \Ap ( \widehat{\rsf}_{t_0}, \widehat{\mgr}_{\tilde t} \, ; p_{t}, q_w, 1_w, x_{t}) + o(1).
\]
Here the approximative matrix $\widehat \Ap ( \widehat{\rsf}_{t_0}, \widehat{\mgr}_{\tilde t} \, ; p_{t}, q_w, 1_w, x_{t})$ is computed from the same stable Riemann surface $\widehat S_{t_0}$ and the metric graph $\widehat{\mgr}_{\tilde t}$ as in the previous case (d). Recall that $\widehat S_{t_0}$ is obtained by adding to $\rsf_{t_0}$ a  Riemann sphere $\mathbb{P}_w$ with a marked point $1_w$ at the node $p_e$ (see Figure~\ref{fig:ProofPicAuxiliaryCase}). As in case (d), the metric graph $\widehat{\mgr}_{\tilde t}$ is the metric realization of the dual graph $\widehat{G} = (\widehat{V}, \widehat{E})$ with lengths $\ell_{e_1'} = -\log\abs{z^e_u(y)}$ and $\ell_{e_2'} = -\log\abs{z^e_v(y)}$. The regularization on $\mathbb{P}_w$ is in the global coordinate $\zeta$ with $\zeta(1_w) = 1$ and such that $\zeta(p^{e_1}_w) = 0$ and $\zeta(p^{e_2}_w) = +\infty$. Finally, $q_w$ is the point on the sphere $\mathbb{P}_w$ with coordinates $\zeta(q_w) = \frac{z^e_u(y)}{z^e_u(y)}$.

\smallskip

Applying Theorem~\ref{thm:SeparateComplexGraph} (see also \eqref{eq:SchemeHeightPairingProof}), we hence infer that
\[
\hp{\rsf_t}{p_t - q, y - x_t} = \hp{\widehat{\mc}_{\tilde t}}{ p_t - q_w,  1_w -  x_t } + o(1),
\]
where $o(1)$ term goes to zero uniformly, if $t$ converges to $s$ in $B$ and $q,y$ satisfy (e) or (f). Recall that we can write
\begin{equation} \label{eq:CC0}
\hp{\widehat{\mc}_{\tilde t}}{ p_t - q_w,  1_w -  x_t } = \hp{\widehat{\mgr}_{\tilde t}}{\vartheta, \theta} + \sum_{a \in \widehat{V}}  \hp{C_{a, t}}{ D^\cmc_{p_t \tiret q_w, a}, D^\cmc_{1_w - x_t ,a }}',
\end{equation}
where $\vartheta$, $\theta$ denote the exact one forms on $\widehat{\mgr}_{\tilde t}$ with
\begin{align*}
&\partial \vartheta = \bar p_{t} - w, & \partial \theta = w - \bar x_{t},
\end{align*}
where $ \bar p_{t} $, $\bar x_{t}$ are the vertices of $\widehat G$ accommodating the sections $p_t$ and $x_t$. Moreover, the divisors $D^\cmc_{p_{t} \tiret q_w, a}$ and $D^\cmc_{1_w - x_{t},a }$, $a \in \widehat{V}$, are obtained by modifying $p_{t} - q_w$ and $1_w -x_{t}$ with the coefficients of $\theta$ and $\vartheta$, respectively. In fact, we have
\[
D^\cmc_{1_w \tiret x_t, a} =  D_{\bar y \tiret \bar x_t,a}^\cmc  
\]
for all vertices $a \neq w$ of $\widehat G$. By Proposition~\ref{prop:explicit_regularization_sphere}, the term originating from the additional Riemann sphere $\mathbb{P}_w$ is given by
\begin{equation} \label{eq:CC1}
\hp{\mathbb{P}_w}{ D^\cmc_{p_t \tiret q_w, w}, D^\cmc_{1_w \tiret x_t,w }}' =  \log  \Big | \frac{z^e_u(y)}{z^e_u(q)} - 1 \Big  |  - \log \Big | \frac{z^e_u(y)}{z^e_u(q)} \Big | \cdot \theta(e) = \varepsilon_\pi(q,y) -  \log \Big | \frac{z^e_u(y)}{z^e_u(q)} \Big | \cdot \theta(e).
\end{equation}
From this point on, on we distinguish the cases (e) and (f).

\smallskip
We begin with case (f). Denote by $\vartheta'$ the unique exact one-form on $\widehat{\mgr}_{\tilde t}$ with
\[
\partial \vartheta' = \bar p_{t} - u.
\]
(recall that $e = uv$ is the edge connecting the vertices $u$ and $v$ in the original graph $G$).  Then, we have
\[
\vartheta - \vartheta' = \chi_{e_1} + \gamma,
\]
where $\chi_{e_1}$ is the one-form on $\widehat{\mgr}_{\tilde t}$ with $\chi_{e_1}(e') = \delta_{e_1, e'}$ for all edges $e' \in \widehat{E}$, and $\gamma$ is a cycle on $\widehat G$. Moreover, it turns out that
\begin{equation} \label{eq:ModificationSlopesZero}
	\sup_{e' \in \widehat E} |\gamma(e')| \to 0
\end{equation}
uniformly as $t$ converges to $s$ in $B$ and $q,y$ satisfying (f). This can be deduced from Proposition~\ref{prop:OrthoProjection}, after taking into account that $|z^e_u(y)| \ge \varrho_e(t)^{3/2}$ and $\lim_{t \to s} \log\varrho_e(t) / \log|z_e(t)| = 0$ by assumption.

\smallskip

Since $\theta$ is an exact one-form on $\widehat{\mgr}_{\tilde t}$, we have
\[
\hp{\widehat{\mgr}_{\tilde t}}{\vartheta, \theta} = \hp{\widehat{\mgr}_{\tilde t}}{\vartheta', \theta} + \log|z^e_u(y)| \theta(e) = \hp{\mgr_{t}}{\bar p_t - \bar q, \bar x_t - \bar y}   + \log|z^e_u(y)| \theta(e).
\]
Moreover, the asymptotic \eqref{eq:ModificationSlopesZero} implies that
\begin{align*}
D^\cmc_{p_{t} \tiret q_w, u} &=  D_{\bar p_t \tiret \bar q,u}^\cmc  + (q - p^e_u)  + o(1) \\
D^\cmc_{p_{t} \tiret q_w, a} &=  D_{\bar p_t \tiret \bar q,a}^\cmc + o(1), \qquad \text{ for all vertices $a \neq u $ of $\widehat G$}.
\end{align*}
Since $|z^e_u(q)| \to 0$ uniformly in the subregion (f), we get that
\begin{equation}  \label{eq:CC2}
\hp{C_{u,t_0}}{ D_{\bar y \tiret \bar x_t,u}^\cmc  , q - p^e_u}'  = - \log|z^e_u(a)| \theta(e) + o(1).
\end{equation}
Altogether, we conclude that
\begin{equation} \label{eq:CC3}
\sum_{a \in {V}}  \hp{C_{a, t_0}}{ D^\cmc_{p_t \tiret q_w, a}, D^\cmc_{1_w - x_t ,a }}' = \sum_{a \in {V}}  \hp{C_{a, t_0}}{ D^\cmc_{p_t \tiret q_w, a}, D_{\bar p_t \tiret \bar q,a}^\cmc}' - \log|z^e_u(a)| \theta(e) + o(1),
\end{equation}
where the $o(1)$ term goes to zero uniformly for $q,y \in W_{e,t}$ satisfying (f), as $t$ converges to $s$ in $B$. Combining \eqref{eq:CC1}, \eqref{eq:CC2} and \eqref{eq:CC3}, we arrive at the claimed asymptotic \eqref{eq:CrucialAsymptotics}.

\smallskip

It remains to treat the case (e). We follow the approach in case (f). However, in order to be formally correct, we have to consider the dual graph $\~G = (\~V, \~E)$ associated to the stable marked Riemann surface $\~\rsf$. We consider the metric graph $\~\mgr_{\tilde t}$ obtained from $\~G$ by setting $\tilde \ell(f) = -\log|\tilde z_{e'}(\tilde t)|$ for $e'\in \~E$. Let $\vartheta'$ be the unique exact one-form on $\~\mgr_{\tilde t}$ such that
\[
	\partial \vartheta' = \bar p_{t} - w_1
\]
in $\~\mgr_{\tilde t}$. Comparing to the original one-form $\theta$ (seen here as an exact one-form on $\~\mgr_{\tilde t}$), we get that
\[
\vartheta - \vartheta' = \chi_{e_2} + \gamma
\]
where $\chi_{e_2}$ is the one-form on $\~{\mgr}_{\tilde t}$ with $\chi_{e_2}(e') = \delta_{e_2, e'}$ for all edges $e' \in \~{E}$, and $\gamma$ is a cycle on $\~G$. Similar to case (f),
\begin{equation} \label{eq:ModificationSlopesZero}
	\sup_{e' \in \~E} |\gamma(e')| \to 0
\end{equation}
uniformly as $t$ converges to $s$ in $B$ and $q,y$ satisfy (e). This follows again from Proposition~\ref{prop:OrthoProjection}, after taking into account that $|z^e_u(y)| / |z^e_u(q) | \ge \varrho_e(t)$ and $\lim_{t \to s} \log\varrho_e(t) / \log|z_e(t)| = 0$ by assumption. Since $\theta$ is an exact one-form on $\~{\mgr}_{\tilde t }$, it follows that
\[
\hp{\~{\mgr}_{\tilde t}}{\vartheta, \theta} = \hp{\~\mgr_{\tilde t}}{\vartheta', \theta} + \log\Big| \frac{z^e_u(y)}{z^e_u(y)} \Big| \theta(e) = \hp{\mgr_{t}}{\bar p_t - \bar q, \bar x_t - \bar y}   + \log\Big |\frac{z^e_u(y)}{z^e_u(y)} \Big| \theta(e) .
\]
Following the steps in step (f), we conclude that
\[
\sum_{a \in V}  \hp{C_{a, t_0}}{ D^\cmc_{p_{t} \tiret q_w, a}, D^\cmc_{1_w - x_{t},a }}' = \sum_{a \in {V}}  \hp{C_{a, t_0}}{ D^\cmc_{p_{t} \tiret q_w, a}, D_{\bar p_t \tiret \bar q,a}^\cmc}' + o(1).
\] 
Taking into account \eqref {eq:CC0} and \eqref{eq:CC1}, we arrive at the desired conclusion, which finishes the proof of the case (A) in Theorem~\ref{lem:CrucialLemma}. \qed

\subsection{Proof of Theorem~\ref{lem:CrucialLemma} in Case (B): both points in the smooth part}
In order to treat the second case (B), let us start with some general considerations. Let  $\mgg{\combind{g,4}}$ be the moduli space of Riemann surfaces $S$ with four markings, denoted by $p,q, x,y$.
Consider the following map between different moduli spaces:
\begin{equation} \label{eq:moduli-diagram-uniformity}
\begin{tikzcd}
\arrow[r, "\proj_{4,1}"]\mgbarg{\combind{g,4}} & \mgbarg{\combind{g,1}}
\end{tikzcd}
\end{equation}
where $\proj_{4,1}$ is the forgetful map which forgets the last three points $q,x,y$ (and contracts the unstable Riemann spheres). Recall we had fixed a stable Riemann surface $S_0$ with one marking $p_0$ and with dual graph $G = (V, E, \genusfunction, \marking)$. Denote by $s_0$ the corresponding point in $\mgbarg{\combind{g,1}}$. Since the map $\proj_{4,1}$ is proper, the fiber  $K:=\proj_{4,1}^{-1}(s_0)$ is a compact subset of $\mgbarg{\combind{g,4}}$.

\smallskip

 The points in $K$ are stratified into their combinatorial types. These are precisely the stable graphs with four markings which by forgetting the last three markings  yield the same stable graph $G$.

\smallskip

As above, denote by $\rsf \to B$ the versal family over versal deformation space $B$ of $S_0$. Shrinking $B$ if necessary, for each point in $K$ we choose an \'etale chart consisting of the versal deformation space of the corresponding curve with four markings, which lies also above $B$.

\smallskip

We cover $K$ by \'etale charts given by the versal deformation spaces of the corresponding stable Riemann surfaces  with four markings. Since $K$ is compact, it is covered by finitely many of these open sets $B_j$, $j\in J$, for a finite set $J$. Shrinking $B$ if necessary, we moreover suppose that the following statement holds.

   \begin{prop} For each $t$ close enough to $0$ in $B$, the fiber $\proj_4^{-1}(t)$ is entirely covered by the open sets $B_j$, $j\in J$.
   \end{prop}

Assume that each $B_j$ is the versal deformation space of a stable Riemann surface $S_j$ with four markings $p_j, q_j, x_j, y_j$,  $j\in J$, and denote by $G_j=(V_j, E_j)$ the corresponding stable graph with four markings.

 \smallskip 
 
 Consider now the surface $\rsf_t$, $t \in B^\ast$, together with the marking $p_t$ in the original family $\rsf \to B$. Let $q,y$ be two points satisfying the assumptions (B). Considering also the section $x_t$, we get a Riemann surface with four markings. In particular, we can view the data $(\rsf_t, p_t, q,x_t,y)$ as a point $\tilde t (q,y) $ in the deformation space $B_j$ of one of the marked Riemann surfaces $S_j$. Using the continuity of the forgetful map $\proj_{4,1} \colon \mgbarg{\combind{g,4}} \to \mgbarg{\combind{g,1}}$ and the assumption (B), we can moreover suppose the following: the dual graph $G_j$ is either equal to $G$, that is, $G_j = G$, or it is obtained from $G$ by adding finitely many trees at the vertices of $G$. In both cases, we can identify $E$ with a subset $E \subseteq E_j$ of the edge set $E_j$ of $G_j$.
 
 \smallskip
 Suppose now that $q, y \in \rsf_t$ satisfy (B) and $\tilde t =\tilde t(q,y)$ belongs to the deformation space $B_j$ of some $S_j$ with the above properties. Consider the stratum
 \[
 D_{E} = \bigl\{\tilde s \in B_j\, \st \, s_{e'} = 0 \text{ for all $e' \in E \subseteq E_j$} \bigr\}
 \]
 in the deformation base $B_j$. As $t$ converges to $s$ in $B$, we have that
 \[
 \tilde z_{e'}(\tilde t) = \tilde z_{e'}(t) \to 0
 \] 
for all $e' \in E \subseteq E_j$. In particular, we can apply the nilpotent orbit theorem to the projection onto the stratum $D_E$ in $B_j$. As one readily computes, the resulting approximative period matrix is given by
 \[
 \Im(\widehat \Omega_{\tilde t}) = \mathcal{A}( \rsf_{t_0}, \mgr_t; p_{t}, q , x_t , y)  + o(1),
\]
where the error term goes to zero as $t$ converges to $s$ in $B$, uniformly for all points $q,y \in \rsf_t$ such that $\tilde t = \tilde t(q,y)$ belongs to $B_j$.  Applying the results of Section~\ref{sec:ProofHeightPairing}, we arrive at the conclusion.

\subsection{Proof of Theorem~\ref{lem:CrucialLemma} in Case (C) and Case (D)}
It remains to establish Theorem~\ref{lem:CrucialLemma} in the above cases (C) and (D). Consider first the case (C), that is, suppose that $q \in W_{e,t}$ and  $y \in W_{e',t}$ for two cylinders $W_{e,t}$ and  $W_{e',t}$ associated to edges $e \neq e'$ in $G$. We outline how to obtain the proof, using similar steps as in case (A).

\smallskip

\begin{figure}[!h]
\centering
 \scalebox{.5}{\input{ProofPicTwoCylinders.tikz}}
\caption{$S_0$ and its dual graph $G$, together with the auxiliary stable Riemann surface $\~S$, used in the analysis of case (C), and its dual graph $\~G = (\~V, \~E)$.}
\label{fig:ProofPicTwoCylinders}
\end{figure}
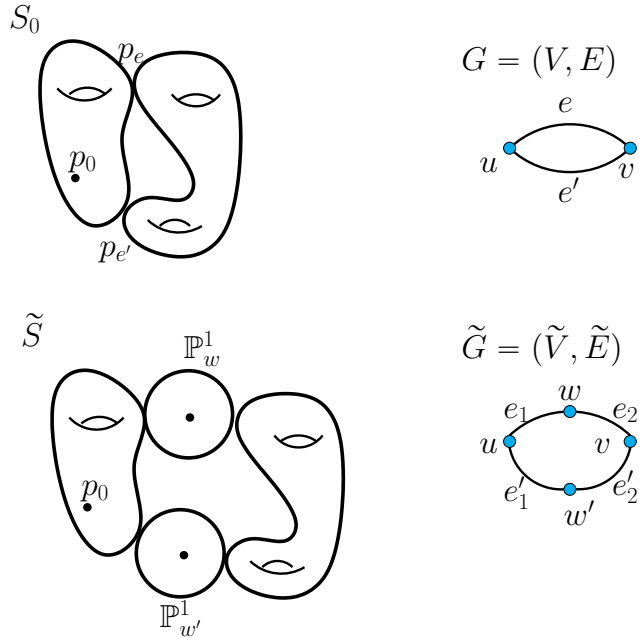

Again we work with a deformation space $\~\rsf \to \~B$ of an auxiliary marked Riemann surface $\~S$. In this case, $\~S$ is obtained by replacing the nodes $p_e$ and $p_{e'}$ in the stable Riemann surface $S_0$ with two Riemann spheres $\mathbb{P}_w$ and $\mathbb{P}_{w'}$, each of them accommodating one marked point (see Figure~\ref{fig:ProofPicTwoCylinders}). Compared to the dual graph $G$ of $S_0$, the dual graph $\~G$ has two new vertices $w$ and $w'$, and the edges $e = uv$ and $e'= u'v'$ are subdivided into new edges $e_1 = u w$, $e_2 = w v$ and $e_1' = u' w'$, $e_2' = w' v'$, respectively. Similar as for case (A), we consider a deformation family $\~\rsf \to \~B$ with adapted coordinates compatible to $\rsf \to B$. Analogous to \eqref{eq:MapDeformationSpaces}, we can interpret $\rsf_t$ with the additional markings $q,y$ as a point $\tilde t \in \~B$, whose new coordinates $\underline{\tilde z}$ are given by 
\[
	\tilde z_{e_1}(\tilde t) = z^e_u(q), \quad \tilde z_{e_2}(\tilde t) = z^e_v(q), \quad \tilde z_{e_1'}(\tilde t) = z^{e'}_{u'}(q),  \quad \tilde z_{e_2'}(\tilde t) = z^{e'}_{v'}(q).
\]

After these preparations, the proof follows by using the same steps as in case (A).  Note that the main difficulty there was due to the fact that the relative position of the points $q,y$ on the cylinder $W_{e,t}$ matters as well, making the term $\varepsilon(q,y)$ appear. Since this does not play a role for case (C), the analysis becomes in fact simpler, and hence we omit the details.

\smallskip

Clearly, the remaining case (D) can be treated as a combination of the above, e.g. by first adding a marked Riemann sphere in the middle of the edge $e$ and then performing the steps from case (B). Hence we omit again the details.


\subsection{Uniform approximation of the height pairing: revisited} \label{ss:UniformApproximationHeightPairingRevisited}
We resume the discussion from Section~\ref{ss:UniformApproximationHeightPairing} and, notations as there, consider a family of hybrid curves $\rsf^\hyb \to B^\hyb$ associated to a marked Riemann surface $(S_0, q_1, \dots, q_n)$ with dual graph $G = (V,E)$.

\smallskip

Recall that Theorem~\ref{thm:HPNoSections} describes the height pairing on smooth fibers $\rsf_t$ close to a hybrid curve $\curve = \rsf^\hyb_\shy$ with underlying stable Riemann surface $S_0$. That is, the corresponding hybrid point $\shy = (l,0)$ has complex coordinate $s = 0 \in B$. However, since stratumwise tameness concerns all boundary strata in $B^\hyb$, we need to consider general hybrid fibers in $\rsf^\hyb$ and also vary the ordered partition $\pi$ in the definition of the metrized complex $\mc_t$. The resulting refinement of Theorem~\ref{lem:CrucialLemma} is collected here.

\smallskip

Fix a hybrid base point $\shy = (l,s)$ in the boundary $\partial_\infty B^\hyb$ and consider the hybrid curve $ \rsf_\shy^\hyb$. Suppose that $\shy$ belongs to the closure $\bar D_\pi^\hyb$ of the hybrid stratum $D_\pi^\hyb$, associated to an ordered partition $\pi \in \Pi(E_\pi)$ on a subset $E_\pi$ of $E$.  Let $G_\pi = (V_\pi, E_\pi)$ be the graph obtained from $G = (V,E)$ by contracting the edges $e \in E \setminus E_\pi$. Note that $G_\pi$ coincides with the dual graph of all fibers $\rsf_t$, $t \in \inn D_{E_\pi}$. Moreover, $G_\pi$ is either equal to or a contraction of the underlying graph of $\rsf_\shy^\hyb$.

\smallskip

Following the steps in Section~\ref{ss:UniformApproximationHeightPairing}, we consider the stratumwise log maps $\loghyb_\pi \colon B^\ast \to D_\pi^\hyb$ and $\loghyb_\pi \colon \rsf^\ast \to \rsf_\pi^\hyb$. Similarly, to a point $t \in B^\ast$, we associate the metrized complex ${\mc}^\pi_t$ which is conformally equivalent to $\rsf_\thy^\hyb$ and has lengths $\ell_t(e) := - \log|z_e(t)|$, $e \in E_\pi$.
Again, the hybrid log map $\loghyb_\pi \colon \rsf^\ast \to \rsf_\pi^\hyb$  leads to a map $\varphi_t \colon \rsf_t \to {\mc}^\pi_t$ for every $t \in B^\ast$. We denote by $\bar y = \varphi_t(y)$ the image of $y \in \rsf_t$ in $\mc_t$.

\smallskip

Allowing the complex part of $\shy$ to be any point $s$ on the divisor $D \subset B$, we obtain the following slight generalization of Theorem~\ref{lem:CrucialLemma}. 
\begin{thm} \label{lem:AuxiliarLemmaClosedStrataJFunc}
Let $\shy = (l,s)$ be a hybrid point in the closure $\bar D_\pi^\hyb$ of some boundary stratum $D_\pi^\hyb$ in $B^\hyb$. Suppose that $t \in B^\ast$ converges tamely to $\shy$. Then for $p,q,x,y \in \rsf_t$,
\begin{equation} \label{eq:CrucialAsymptotics}
\hp{\rsf_{t}}{p - q, x - y} =  \hp{{\mc}^\pi_t}{\bar p-  \bar q, \bar x - \bar y} + \varepsilon_{\pi}(p,q,x,y) + o(1),
\end{equation}
where the $o(1)$ term goes to zero uniformly for $q, y,x,y \in \rsf_t$.
\end{thm}
Recall that the appearing correction term $\varepsilon_{\pi}(p,q,x,y)$ is defined in \eqref{eq:CorrectionFourPoints}.


\subsection{Proof of Theorem~\ref{thm:FinalTamenessHeightPairing}} \label{ss:FinalTamenessHeightPairing}
We begin by recalling the framework around Theorem~\ref{thm:FinalTamenessHeightPairing}. Consider a stable  (possibly marked) Riemann surface $S_0$ having dual graph $G = (V,E)$. Let $\rsf \to B$ be a deformation family for $S_0$, equipped with a fixed choice of adapted coordinates, and $\rsf^\hyb \to B^\hyb$ the associated family of hybrid curves.

We have to verify the following property: for any boundary stratum $D_\pi^\hyb$ in $B^\hyb$ and any point $\shy$ in the closure $\bar D_\pi^\hyb$ of $D_\pi^\hyb$, we have
\begin{equation*} \label{eq:ToProveUniformTameHP}
 \hp{\rsf_t}{p \tiret q , x \tiret y } =  \sum_{j \in \{1, \dots, r, \smallcc\}} L_j(t)  \lhp{\rsf^\hyb_\thy, j}{\Lognoind(p) \tiret \Lognoind(q) , \Lognoind(
 x) \tiret  \Lognoind(y)} + \varepsilon_{\pi}(p,q,x,y) + o(1), \, \, \thy = \Lognoind (t),
\end{equation*}
uniformly for $p,q,x,y \in \rsf_t$, if $t \in B^\ast$ converges tamely to $\shy$ in $B^\hyb$. Here, we set $L_j = L_j(t) := \sum_{e \in \pi_j} -\log|z_e(t)|$ for $j=1, \dots, r$ (by convention $L_\smallcc \equiv 1$); we abbreviate by $\Lognoind := \loghyb_\pi$ the stratumwise hybrid log maps on $B^\ast$ and $\rsf^\ast$; and $\lhp{\rsf^\hyb_\thy, j}{\cdot , \cdot}$ for $j \in \{ 1, \dots, r, \smallcc\}$ is the $j$-th component of the hybrid height pairing on the hybrid curve $\rsf^\hyb_\thy$.

Applying Theorem~\ref{thm:HPNoSections} we get a uniform approximation (see also \eqref{eq:RepeatMetrizedHeightPairingDefinition})
\[
\hp{\rsf_{t}}{ p - q, x  - y} = \hp{\mgr_{t}}{\bar p - \bar q, \bar x - \bar y} + \sum_{v\in V} \hp{C_{v,t}}{D_{\bar p \tiret \bar q,v}^\cmc, D_{\bar x \tiret \bar y ,v}^\cmc}' + \varepsilon_\pi(p,q,x,y) +  o(1)
\]
as $t \in B^\ast$ converges to $\shy$, where $\mc_t$ is the metrized complex with underlying metric graph $\mgr_t$ defined in Section~\ref{ss:UniformApproximationHeightPairing}. If $t \in B^\ast$ converges tamely to $\shy$ in $B^\hyb$, then the metric graphs $\mgr_t$ converge tamely to the underlying tropical curve $\curve^\trop$ of the hybrid curve $\rsf_\shy^\hyb$ (in the sense that this convergence holds for the associated points in $\mgtropcombin{\grind{G}}$). Applying Lemma~\ref{lem:JfunctionTame} and taking into account the connection between height pairings and $j$-functions (see \eqref{eq:MGHPvsJ}, Theorem~\ref{thm:TropicalHPvsJFunction}, and also Remark~\ref{rem:TropicalVSHybridHeightPairing}), it follows that
\[
\hp{\mgr_{t}}{\bar p - \bar q, \bar x - \bar y} = \sum_{j  = 1}^r L_j(t) \,  \lhp{\rsf^\hyb_\thy, j}{\Lognoind(p) \tiret \Lognoind(q) , \Lognoind(
 x) \tiret  \Lognoind(y)} + o(1)
\]
where $\curve^\trop_\thy$ is the underlying tropical curve of the hybrid curve $\rsf_\thy^\hyb$, $\thy = \Lognoind(t)$, and the $o(1)$ term goes to zero uniformly for $p,q,x,y \in \rsf_t$. 

It remains to prove that, as $t \in B^\ast$ converges tamely to $\shy$, 
\begin{equation} \label{eq:ChangeDivisors}
\sum_{v \in V} \hp{C_{v,t}}{D_{\bar p \tiret \bar q,v}^\cmc, D_{\bar y - \bar x ,v}^\cmc}'  =  \lhp{\rsf^\hyb_\thy, \smallcc}{\Lognoind(p) - \Lognoind(q), \Lognoind(y) - \Lognoind(x)} + o(1),
\end{equation}
where the $o(1)$ term goes to zero uniformly. Recall that in the definition of the height pairing on $\rsf_\thy^\hyb$ (see Definition~\ref{def:HybridHeightPairing}), we associated to every degree zero divisor $D$ a collection of divisors $D^\smallcc_v$, $v \in V$, on the Riemann surface components  of $\rsf_\thy^\hyb$. Hence we get a collection of divisors $D^\smallcc_{\Lognoind(p)  \tiret   \Lognoind(q), v}$ and $D^\smallcc_{\Lognoind(y) \tiret \Lognoind(x),v }$ on the components $C_{v,t}$, $v \in V$ of $\rsf_\thy^\hyb$. We have
\begin{align*}
&D_{\bar p \tiret \bar q,v}^\cmc =  D^\smallcc_{\Lognoind(p)  \tiret  \Lognoind(q), v} + R_{ v}, &D_{\bar y \tiret \bar x,v}^\cmc =  D^\smallcc_{\Lognoind(y)  \tiret  \Lognoind(x),v } + R_{ v}',
\end{align*}
where the degree zero divisors $R_{v}$ and $R'_{v} $ on $C_{v,t}$ are supported on the attachement points $p^e_v$, $e \sim v$, and have small coefficients. More precisely,
\begin{align*}
&R_{ v} = \sum_{e \in v} a_e \delta_{p^e_v}, &R'_{ v} = \sum_{e \in v} a_e' \delta_{p^e_v}
\end{align*}
and if the edge $e$ belongs to the $j$-th layer $\pi_j$, then 
\[
 \max\{|a_e|, |a_e'|\} \le C \, \frac{L_{j+1}(t)}{L_j(t)}
\]
as $t \in B^\ast$ converges tamely to $\shy$,  where the constant $C$ is uniform for $p,q,x, y \in \rsf^\hyb_t$. Indeed, this follows from Lemma~\ref{lem:JfunctionTame}, since the metric graph $\mgr_t$ converges tamely to the underlying tropical curve $ \curve^\trop$ of $ \rsf^\hyb_\shy$ in $\mgtropcombin{\grind{G}}$. This in turn implies that
\begin{equation} \label{eq:ChangeDivisors2}
\hp{C_{v,t}}{ D_{\bar p \tiret \bar q,v}^\cmc  , D_{\bar y \tiret \bar x,v}^\cmc  }'  = 
\hp{C_{v,t}}{ D^\smallcc_{\Lognoind(p)  \tiret  \Lognoind(q), v}  ,  D^\smallcc_{\Lognoind(x)  \tiret  \Lognoind(y), v}  }'  + o(1),
\end{equation}
where the $o(1)$ term goes to zero uniformly for $p,q,x, y \in \rsf_t$, as $t \to \shy$ tamely. Note that here we implicitly used properties of the hybrid log map $\Lognoind$ in order to ensure that
\begin{align*}
\hp{C_{v,t}}{R_v, D_{\bar y \tiret \bar x  ,v}^\cmc}' = o(1)
\end{align*}
uniformly (and similar for the other appearing terms). Namely, the image points $\bar y=\varphi_t(y)$ and $\bar x=\varphi_t(x)$ stay outside the circle $\{|z^e_u| \ge \varrho_e(t)\}$  around every attachement point $p^e_v$, $e \sim v$ (here $z^e_v$ is the adapted local coordinate). This causes a potential pole arising from $x$ or $y$ coming close to $p^e_v$ to be controlled by $|a_e \log\varrho_e(t)|$, which goes to goes to zero uniformly if $t \in B^\ast$ converges tamely to $\shy$ by the above estimate and the properties of the tame topology. Summing \eqref{eq:ChangeDivisors2} over all vertices $v \in V$, we have established \eqref{eq:ChangeDivisors}.

\qed


\section{Tameness of the Arakelov Green function: Proofs of Theorems~\ref{thm:MainGeneralMeasures},~\ref{thm:WeakConvergence},~\ref{thm:WeakTamenessArakelovGreenFunction},~\ref{thm:MainArakelovDetails}, and~\ref{thm:MainArakelovDetails2}}
\label{sec:Green_functions_asymptotics}

In this section we prove our main results on asymptotic behavior of the Arakelov Green function and the solutions of Poisson equations on degenerating families of Riemann surfaces.

We begin with the proofs of our asymptotic results on the Arakelov Green function, Theorem~\ref{thm:MainArakelovDetails} (Section~\ref{ss:ProofMainArakelovDetails}) and Theorem~\ref{thm:MainArakelovDetails2} (Section~\ref{ss:ProofMainArakelovDetails2}). Theorem~\ref{thm:WeakTamenessArakelovGreenFunction} is then proved in Section~\ref{ss:ProofMainGreenFunction}. Finally, Section~\ref{ss:ProofMainGeneralMeasures} contains the proofs of Theorem~\ref{thm:MainGeneralMeasures} and Theorem~\ref{thm:WeakConvergence}.    
\subsection{Proof of Theorem~\ref{thm:MainArakelovDetails}} \label{ss:ProofMainArakelovDetails}
In the following, we give the proof of Theorem~\ref{thm:MainArakelovDetails}. We use the notation introduced in the respective Section~\ref{ss:MainArakelovDetails}.

The starting point of the proof is the \emph{integral representation of the Green function} in terms of the height pairing and the canonical measure:
\begin{equation} \label{eq:GreenIntegral}
\gri{S}(p,y) = \int_S \hp{S}{p-q, y-x} d\mu(q) - \int_S \int_S \hp{S}{p-q, \eta-x} d\mu(q) d \mu(\eta).
\end{equation}
for two distinct points $p, y$ on a smooth Riemann surface $S$. This expression is the precise analogue of the formula for Green functions on metric graphs (see \eqref{eq:GreenByJFunction} and also \eqref{eq:MGHPvsJ}).

\subsubsection{Proof of Theorem~\ref{thm:MainArakelovDetails}(i)} \label{ss:ProofClaimA} 
Recall that the marking $p_0$ of the stable marked Riemann surface $S_0$ leads to a section $p_t$, $t \in B$, of the family $\rsf \to B$. In a small neighborhood $U$ of $s = 0$ in $B$, we can fix another \emph{auxiliary section} $x_t$, $t \in U$, of $\rsf \to B$, such $x_t$ belongs to the set $Y_{w,t}$ for some fixed vertex $w$ and all $t$ (see the decomposition~\eqref{eq:AdaptedDecompositionFinals}), and $x_t$ only depends on the adapted coordinates associated to the component $C_w$ of $S_0$, and $p_t \neq x_t$ for all $t$.

\smallskip

Consider the smooth marked Riemann surface $(\rsf_t, p_t)$ for $t \in B^\ast$. Combining the representation \eqref{eq:GreenIntegral} with the uniform tameness of the height pairing (see Theorem~\ref{thm:FinalTamenessHeightPairing}),
\begin{equation} \label{eq:GFFourTerms}
\gri{\rsf_t}(p_t,y) = h_1(y) + h_2(y) + h_3(y) + h_4(y),
\end{equation}
for any point $y \neq p_t$ on the smooth Riemann surface $\rsf_t$, $t \in B^\ast$, where
\begin{align*}
h_1(y) &= \sum_{j=1}^r L_j(t) \int_{\rsf_t}  \lhp{\rsf^\hyb_\thy, j}{\Lognoind(p_t) - \Lognoind(q) ,  \Lognoind(y) - \Lognoind(
 x_t)}  d\mu_t(q) \\
 & \qquad \qquad - L_j(t)\int_{\rsf_t} \int_{\rsf_t} \lhp{\rsf^\hyb_\thy, j}{\Lognoind(p_t) - \Lognoind(q) ,  \Lognoind(\eta) -  \Lognoind(
 x_t)}   d\mu_t(q) d \mu_t(\eta) \\
h_2(y) &= \int_{\rsf_t}   \lhp{\rsf^\hyb_\thy, \smallcc}{\Lognoind(p_t) - \Lognoind(q) ,  \Lognoind(y) -  \Lognoind(
 x_t)} d\mu_t (q)  \\
 & \qquad \qquad - \int_{\rsf_t} \int_{\rsf_t}  \lhp{\rsf^\hyb_\thy, \smallcc}{\Lognoind(p_t) - \Lognoind(q) ,  \Lognoind(\eta) -  \Lognoind(
 x_t)}  d\mu_t(q) d \mu_t(\eta) \\
h_3(y) &= - \int_{\rsf_t} \varepsilon_\pi(q,y)  \, d\mu_t(q) + \int_{\rsf_t} \int_{\rsf_t} \varepsilon_\pi(q,\eta) \, d\mu_t(q) d \mu_t(\eta) \\
h_4(y) &= \int_{\rsf_t} o(1) \, d\mu_t(q) - \int_{\rsf_t} \int_{\rsf_t} o(1) \, d\mu_t(q) d \mu_t(\eta).
\end{align*}
Here, $\varepsilon_\pi(q,y)$ is the correction term defined in \eqref{eq:CorrectionTwoPoints}, $\lhp{\rsf^\hyb_\thy, j}{\cdot , \cdot}$ for $j = 1, \dots, r, \smallcc$ is the $j$-th component of the height pairing on the hybrid curve $\rsf^\hyb_\thy$, $\thy = \loghyb_\pi(t)$, and  $\Lognoind := \loghyb_\pi$ abbreviates the hybrid log map from $\rsf_t$ to $\rsf^\hyb_\thy$. Moreover, $o(1)$ denotes a term depending on $q$ and $y$ on $\rsf_t$, which goes to zero uniformly as $t \in B^\ast$ converges tamely to $\shy$.

\smallskip

In the following, we analyze the four terms $h_i(y)$ in \eqref{eq:GFFourTerms} separately. As the next two propositions show, $h_1(y)$ and $h_2(y)$ give the contributions stemming from graph minors and the complex part, respectively.

Let the measure $\~\mu_t$ be the push-out of $\mu_t$ from $\rsf_t$ to the hybrid curve $\curve_\thy$ via the log map $\Lognoind \colon \rsf_t \to \rsf^\hyb_\thy$, and denote by $\~\mu_t^\trop$ the further push-out of this measure to the underlying tropical curve $\curve^\trop_\thy$ via the contraction map $\rsf^\hyb_\thy \to \curve^\trop_\thy$. Here, and in the following, when writing $\rsf^\hyb_\thy$, we mean the fiber over $\thy = \loghyb_\pi(t)$ in the family $\rsf^\hyb_\thy$, that we further identify with  the canonical metrized complex  in its conformal equivalence with normalized edge lengths in each layer.  Moreover, the measure $\widetilde \mu_t$ is seen as a (non-layered) measure on the topological space $\rsf^\hyb_\thy$. 

Moreover, adapting our previous convention, we denote by $\bar q$ the image of $q \in  \rsf_t$ under the composition of the log map $\Lognoind \colon \rsf_t \to \rsf^\hyb_\thy$  with  the contraction $\rsf^\hyb_\thy \to \curve^\trop_\thy$. 

\begin{prop} \label{prop:h1}
Consider the smooth Riemann surface $\rsf_t$ for $t \in B^\ast$. Then 
\begin{align*}
&h_1(y) = \sum_{j=1}^r L_j(t) \grihat{t,j}(p_t, y)  \qquad  \text{for all $y \neq p_t$ on $\rsf_t$}.
\end{align*}
\end{prop}
\begin{proof}
By Proposition~\ref{prop:BasicPropertiesHybridGreenFunctions}, Green functions on the hybrid curve $\rsf^\hyb_\thy$, $\thy = \loghyb_\pi(t)$, and its underlying tropical curve $\curve^\trop_\thy$ are closely connected. In particular, $\grihat{t,1}(p_t, \cdot), \dots, \grihat{t,r}(p_t, \cdot)$ are precisely the pullbacks of the components of the Green function on $\curve^\trop_\thy$ for the measure $\~{\mu_t}^\trop$. The claim now follows from the integral representation in Lemma~\ref{lem:SolutionFormulaJFunctionTropical} using the connection between tropical and hybrid height pairings, Remark~\ref{rem:TropicalVSHybridHeightPairing}, and between the tropical height pairings and $j$-functions, Theorem~\ref{thm:TropicalHPvsJFunction}.
\end{proof}
\begin{prop} \label{prop:h2}
Consider the smooth Riemann surface $\rsf_t$ for $t \in B^\ast$. Then
\[
h_2(y) = \grihat{t,\smallcc}(p_t,  y) \qquad  \text{for all $y \neq p_t$ on $\rsf_t$}.
\]
\end{prop}
\begin{proof}
We have to prove that
\[
y \mapsto \int_{\rsf^\hyb_\thy} \lhp{\smallcc}{\Lognoind(p_t) - q, y - \Lognoind(x_t)} d\widetilde{\mu}_t(q) -  \int_{\rsf^\hyb_\thy}  \int_{\rsf^\hyb_\thy} \lhp{ \smallcc}{\Lognoind(p_t) - q, y' - \Lognoind(x_t)} d\widetilde{\mu}_t(q)  d\widetilde{\mu}_t(y') 
 \]
is precisely the complex part $y \in \rsf^\hyb_\thy \mapsto \widetilde{g}_{t,\smallcc}(p_t, y) $ of the Green function $\widetilde{\lgr}(p_t, \cdot)$ for the measure $\widetilde{\mu}_t$ on the hybrid curve $\rsf^\hyb_\thy$.

Subtracting the double integral ensures that the integral over $\rsf_\thy^\hyb$ with respect to $\widetilde{\mu}_t$ vanishes. Hence, it suffices to prove that the function $f_\smallcc \colon \rsf_\thy^\hyb \to \C$ given by
\[
	f_\smallcc (y) := \int_{\rsf^\hyb_\thy} \lhp{\smallcc}{\Lognoind(p_t) - q, y - \Lognoind(x_t)} d\widetilde{\mu}_t(q) , \qquad y \in \rsf_\thy^\hyb,
\]
has the following two properties: 

\smallskip

(i) It satisfies the respective Poisson equation, that is,
\[
\int_{\pi_\smallcc} f_\smallcc \, \Delta_\smallcc h = h(\Lognoind(p_t))  -  \int_{\pi_\smallcc} h \, d \widetilde{\mu}_t + \sum_{j=1}^r   h \Big( \divind{j}{\smallcc}(f_j) \Big)
\]
for any smooth function $h \colon \pi_\smallcc \to \R$ on the disjoint union of components $\pi_\smallcc = \bigsqcup_v C_{v,t}$. Here, $\divind{j}{\smallcc}(f_j)$ is the divisor defined in \eqref{eq:DefineDivisorsHybridLaplacian} and obtained from the $j$-th part $f_j$ of the solution $\lf = (f_1, \dots f_r)$ of the tropical Poisson equation
 \begin{align*}
\Deltatrop \lf = \bm{\delta}_{\bar p_t} - \widetilde{\mu}_t,  \qquad \qquad
\lf \text{ is harmonically arranged},
\qquad \qquad \int_{\curve^\trop_\thy} \lf \, d \bm{\delta}_{\bar x_t} = 0.
\end{align*}
on the underlying tropical curve $\curve^\trop_\thy$ of $\rsf_\thy^\hyb$.

\smallskip

(ii) On each edge $e$ of $\rsf_\thy^\hyb$, the function $y \mapsto f_\smallcc(y)$ is linear, the values $f_\smallcc(p^e_v)$ and $f_\smallcc(p^e_u)$ in the endpoints coincide with the regularized values on the components $C_{u,t}$ and $C_{v,t}$, and $f_\smallcc$ is harmonically arranged. 

\smallskip

By Theorem~\ref{thm:HybridHPvsJFunction}, for any fixed $q \in \rsf_\thy^\hyb$ the function $y \in \rsf_\thy^\hyb \mapsto  \lhp{\smallcc}{\Lognoind(p_t) - q, y - \Lognoind(x_t)}$ is precisely the complex part of the hybrid $j$-function. In particular, we have the equality
\[
\int_{\pi_\smallcc} \lhp{\smallcc}{\Lognoind(p_t) - q, \cdot - \Lognoind (x_t)}  \Delta_\smallcc h = h(\Lognoind(p_t)) + \sum_{k=1} h ( \divind{j}{\smallcc}( \zeta_{q,j}) ) - h (q)
\]
for every smooth function $h \colon \pi_\smallcc \to \R$ (by convention, $h(q) := 0$ if $q \in \rsf_\thy^\hyb \setminus \pi_\smallcc$). Here $\divind{j}{\smallcc}( \zeta_{q,j})$ is the divisor induced by the $j$-th part $\zeta_{q,j}$ of the tropical solution $\bm{\zeta}_q = (\zeta_{q,1}, \dots, \zeta_{q,r})$ of the tropical Poisson equation
 \begin{align*}
\Deltatrop \bm{\zeta}_q = \bm{\delta}_{\bar p_t} - \bm{\delta}_{\bar q}, \qquad \qquad  
\bm{\zeta} \text{ is harmonically arranged}, \qquad \qquad
\int_{\curve^\trop_\thy} \bm{\zeta}_q \, d \bm{\delta}_{\bar x_t} = 0
\end{align*}
on the underlying tropical curve $\curve^\trop_\thy$ of $\rsf_\thy^\hyb$ (see \eqref{eq:DefineDivisorsHybridLaplacian}). This clearly implies that
\begin{align*}
\int_{\pi_\smallcc}  f_\smallcc \, \Delta_\smallcc h &= \int_{\rsf^\hyb_\thy} \Big (\int_{\pi_\smallcc} \lhp{\smallcc}{\Lognoind(p_t) - q, \cdot - \Lognoind(x_t)}  \Delta_\smallcc h  \Big) \, d \widetilde{\mu}_t \\
&= h(\Lognoind(p_t)) - \int_{\pi_\smallcc} h \, d \widetilde{\mu}_t  + \sum_{k=1} h( \divind{j}{\smallcc}(\zeta_{\mu,j}) ),
\end{align*}
where $\divind{j}{\smallcc}(\zeta_{\mu,j})$ is the divisor obtained as the integral 
\[
\divind{j}{\smallcc}(\zeta_{\mu,j}) := \int_{\rsf^\hyb_\thy}  \divind{j}{\smallcc}( \zeta_{q,j}) \, d \widetilde{\mu}_t(q). 
\]
in the sense of pointwise integration of the coefficient). As follows from Lemma~\ref{lem:SolutionFormulaJFunctionTropical},
\begin{align*}
\divind{j}{\smallcc}(\zeta_{\mu,j}) =  
- \int_{\rsf^\hyb_\thy}  \sum_{\substack{e \in \pi_j} } \sum_{\substack{v \in e}} \slp_e \zeta_{q,j} (v) d\~{\mu}_t(q) \, \delta_{p^e_v} =  -\sum_{\substack{e \in \pi_j} } \sum_{\substack{v \in e}} \slp_e f_j (v) \, \delta_{p^e_v} = \divind{j}{\smallcc}(f_j)
\end{align*}
and hence (i) is proven. To verify the properties in (ii), notice that for fixed $q \in \rsf_\thy^\hyb$ the function $y \in \rsf_\thy^\hyb \mapsto  \lhp{\smallcc}{\Lognoind(p_t) - q, y - \Lognoind(x_t)}$ is linear on edges and harmonically arranged (see again Theorem~\ref{thm:HybridHPvsJFunction}). These properties are preserved by integrating in $q$ with respect to $\widetilde{\mu}_t$ and hence carry over to $ f_\smallcc(y)$. The proof is complete.

\end{proof}

Finally, it turns out that the remaining two terms $h_3(y)$ and $h_4(y)$ vanish in the limit $t \to \shy$.   In order to prove this, we need the following basic bounds on the canonical measure. For an edge $e \in E$, consider the cylinder $W_{e,t} \subseteq \rsf_t$ (see ) and its subregion
\[
\~{W}_{e,t} = \Big \{ (\underline{z}, z^e_u, z^e_v) \in W_{e,t} \, \st \, |z^e_u| \le |\log |z_e(t)||^{-1/2} \text{ and }  |z^e_v| \le |\log|z_e(t)||^{-1/2} \Big \}.
\]
Recall moreover that, in Section~\ref{sec:hybrid_log_map_einf}, the cylinder $W_{e,t}$ was decomposed into
\[
W_{e,t} = A^e_{u,t} \sqcup B_{e,t} \sqcup A^e_{v,t}.
\]
Note that $\~{W}_{e,t}$ intersects all three sets in this decomposition.

\begin{lem} \label{lem:CanonicalMeasureBound} As $t \in B^\ast$  converges to $s=0$ in $B$, the canonical measure $\mu_t$ on the smooth fiber $\rsf_t$, $t \in B^\ast$, satisfies the following bounds.

\smallskip

(a) In the adapted coordinate $z^e_u$ on the cylinder $W_{e,t}$, 
\[
\mu_t \rest{W_{e,t}} = f_t \, \, dz^e_u  \wedge d \overline{z^e_u}
\]
for a continuous function $f_t \colon W_{e,t} \to \R$ such that
 \[
 	|f_t(q)| \le C \frac{1}{|\log|z_e(t)||} \frac{1}{|z^e_u(q)|^2}, \qquad q \in \widetilde{W}_{e,t},
 \]
and
\[
 	|f_t(q)| \le C \Big (1 + \frac{1}{|\log|z_e(t)||} \frac{1}{|z^e_u(q)|^2} \Big )  , \qquad q \in A^e_{u,t},
\]
for a constant $C >0$ independent of $t$.

\smallskip

(b) Consider a fixed point $p$ in the smooth part of the fiber $\rsf_0$ over $s = 0$. Let $U$ be a small open neighborhood of $p$ in the family $\rsf \to B$, parametrized by $z_1, \dots, z_N, z$, where  $\underline z = (z_1, \dots, z_N)$ are the adapted coordinates from the base $B= \Delta^N$ and $z$ represents an additional parameter on the fibers of $\rsf \to B$. Then on the subset $U_t = U \cap \rsf_t$ of $\rsf_t$, $t \in B^\ast$, 
\[
\mu_t \rest{U_t} = f_t \, \, dz \wedge d \overline z
\]
for a continuous function $f_t \colon U_t \to \R$ such that
 \[
 	|f_t(q)| \le C, \qquad q \in U_t,
\]
for a constant $C >0$ independent of $t$.
\end{lem}
\begin{proof} The claims are a straightforward consequence of the discussion in \cite[Section 9]{AN}.
\end{proof}

\begin{prop} \label{prop:h3h4}
Notations as above,
\begin{align*}
&\sup_{y \in \rsf_t}|h_3(y)| \to 0, &\sup_{y \in \rsf_t} |h_4(y)| \to 0,
\end{align*}
as $t \in B^\ast$ converges to the fixed hybrid point $\shy$ in $B^\hyb$.
\end{prop}
\begin{proof} Since the canonical measure $\mu_t$ on $\rsf_t$ is a positive measure of total mass one, $h_4(y)$ clearly goes to zero uniformly as $t \in B^\ast$ converges to $\shy$ in $B^\hyb$. In order to prove that $h_3(y)$ goes to zero uniformly, it suffices to verify that
\begin{equation} \label{eq:IntegralEpsilons}
\sup_{y \in \rsf_t} \, \int_{\rsf_t} |\varepsilon_\pi(q,y)| \, d \mu_t (q) \to 0
\end{equation}
as $t \in B^\ast$ converges to $s = 0$ in $B$. Taking into account Lemma~\ref{lem:CrucialLemma} and the symmetry of all appearing terms, this reduces to proving that
\begin{equation} \label{eq:ProofFunctionh3}
\sup_{y \in W_{e,t}} \, \int_{\substack{\{ q \in W_{e,t} \colon |z^e_u(q)| > |z^e_u(y)|\} }} |\varepsilon_\pi(q,y)| \, d \mu_t (q) \to 0
\end{equation}
for every cylinder $W_{e,t}$, $e \in E$. Note that if $y \in W_{e,t}$, and $q \in W_{e,t}$ is a point with $\varepsilon_e(q,y) \neq 0$ whose adapted coordinates satisfy $|z^e_u(q)| > |z^e_u(y)|$, then we are in one of the cases (i), (iii) or (v) from Lemma~\ref{lem:CrucialLemma}. In particular, the condition $\lim_{t \to \shy} \log\varrho_e(t)^2 / \log|z_e(t)| =0$ implies that $q$ belongs to the region $\widetilde{W}_{e,t}$ from Lemma~\ref{lem:CanonicalMeasureBound} and $|z^e_u(y)| / |z^e_u(q)| > \varrho_e(t)$.

\smallskip

Introducing the complex coordinate $\zeta := z^e_u(y) / z^e_u(q) = z^e_v(q) / z^e_u(y)$ and using Lemma~\ref{lem:CanonicalMeasureBound}, the supremum in \eqref{eq:ProofFunctionh3} is controlled by
\[
 \frac{1}{|\log |z_e(t)||} \int_{\{ \varrho_e(t) \le |\zeta| \le 1\} } \frac{1}{|\zeta|^2} \big | \log|\zeta - 1| \big |  \, \big( \frac{i}{2} d \zeta \wedge \overline{\zeta} \big ).
\]
Passing to polar coordinates $r = |\zeta|$ and $\theta = \arg(\zeta)$, we easily conclude that
\[
\int_{\{ \varrho_e(t) \le |\zeta| \le 1/2 \} } \frac{1}{|\zeta|^2} \big | \log|\zeta - 1| \big |  \, \big ( \frac{i}{2} d \zeta \wedge \overline{\zeta} \big )  \lesssim \int_{\varrho_e(t)}^{1/2} \frac{1}{r} \, dr \lesssim |\log \varrho_e(t)|.
\]
Moreover, we have the simple estimate
\[
\int_{\{ \frac{1}{2} \le |\zeta| \le 1\} } \frac{1}{|\zeta|^2} \big | \log|\zeta - 1| \big |  \, \Big ( \frac{i}{2} d \zeta \wedge \overline{\zeta} \Big )  \le \int_{\{ \frac{1}{2} \le |\zeta| \le 1\} } \big | \log|\zeta - 1| \big |  \, \Big ( 2 i d \zeta \wedge \overline{\zeta} \Big) < \infty.
\]
Since $\lim_{t \to \shy} \log\varrho_e(t) / \log|z_e(t)| =0$ by assumption, we have proven the claim on $h_3(y)$.
\end{proof}

Combining the expression \eqref{eq:GFFourTerms} with Proposition~\ref{prop:h1}, Proposition~\ref{prop:h2}, and Proposition~\ref{prop:h3h4}, we conclude with the proof of the statement in Theorem~\ref{thm:MainArakelovDetails}(i).

\qed

\subsubsection{Proof of Theorem~\ref{thm:MainArakelovDetails}(ii)} \label{ss:ClaimB} Consider the hybrid curve $\curve$ represented by the point $\shy = (l,0)$ in $B^\hyb$. Let $\mccan$ be the associated metrized complex with normalized edge lengths in each layer. Recall that we are considering the functions $\gritilde{t, 1}(p_\shy, \cdot), \dots, \gritilde{t, r}(p_\shy, \cdot), \gritilde{t, \smallcc}(p_\shy, \cdot)$, $t \in \loghyb_\pi^{-1}(\shy)$, and functions $\gri{1}(p_\shy, \cdot), \dots, \gri{r}(p_\shy, \cdot),  \gri{\smallcc}(p_\shy, \cdot)$, components of the hybrid Green functions associated to $\~\mu_t$ and $\mu_\shy$ on $\curve$. We regard all these as functions defined on the space $\mccan \setminus \{p_\shy\}$ with values in $\R$. As before, we view the pushout measures $\~{\mu}_t$ and the canonical measures $\mu_\shy $ of $\curve$ as  measures defined on the space $\mccan$. We want to prove that each $\gritilde{t, j}(p_\shy, \cdot)$ converges to $\gri{j}(p_\shy, \cdot)$, in the sense explained in Theorem~\ref{thm:MainArakelovDetails}(ii). 

\smallskip

We begin with the convergence of the coefficients at the logarithmic poles.

\begin{lem}  \label{lem:ConvergenceLogCoeff}
For every point $p \in \Sigma$ and every attachment point $p^e_u$ of an edge $e = uv$, the coefficients in \eqref{eq:PoleCoeff1} and \eqref{eq:PoleCoeff2} satisfy
\begin{equation} \label{eq:ConvergenceLogCoeff}
\~c_t(p, p^e_u) \to c_\shy(p, p^e_u), \qquad \qquad \text{as $t \in \loghyb_\pi^{-1}(\shy) \to \shy$}.
\end{equation}

\end{lem}
\begin{proof}
By the integral representation in Lemma~\ref{lem:SolutionFormulaJFunctionTropical}, the coefficients at $p^e_u$ for an edge $e = uv$ in the $j$-th layer $\pi_j$ have an integral representation as well. More precisely,
\[
\~{c}_t (p, p^e_u) = \int_\mccan \slp_e \zeta_{q,j} (u) \, d \~{\mu}_t (q), \qquad c_\shy(p, p^e_u) =  \int_\mccan \slp_e \zeta_{q,j} (u) \, d \mu_\shy (q).
\]
Here $\zeta_{q,j} \colon \Gamma \to \R$ is the function on the underlying normalized metric graph $\Gamma$ of $\mccan$, which maps $y \in \Gamma$ to the $j$-th component of the height pairing $\zeta_{q,j}(y) := \lhp{j}{\bar p_\shy - \bar q, y - \bar x}$ on the underlying tropical curve $\curve^\trop$ of $\curve$. The convergence then follows from the weak convergence of the measures $\~{\mu}_t$ to $\mu_\shy$. Note that, in general, weak convergence of measures on a metric graph does not imply the convergence of slopes for the corresponding Poisson equations, cf. \cite[Proposition 3.11]{BR10}. However, the additional bounds in Lemma~\ref{lem:CanonicalMeasureBound} imply \eqref{eq:ConvergenceLogCoeff}.
\end{proof}
By (the proof of) Proposition~\ref{prop:h2}, the complex parts of Green functions can be written as
\[
\gritilde{t, \smallcc} (p_\shy, y) = f_t(y) - \int_{ \mccan} f_t(\eta) \, d\widetilde{\mu}_t(\eta), \qquad \gri{\smallcc} (p_\shy, y) = f (y) - \int_{ \mccan} f_t(\eta) \, d \mu_\shy (\eta),
\] 
for $y \in \Sigma \setminus \{p_\shy\}$, using the auxiliary functions $f_t \colon \mccan \to \R$, $t \in \loghyb^{-1}_\pi(\shy)$, and $f \colon  \mccan \to \R$ given by
\begin{equation} \label{eq:DefAuxFcts}
f_t (y) = \int_{ \mccan } \lhp{\smallcc}{p_\shy - q, y - x} d\widetilde{\mu}_t(q) , \qquad \qquad f (y) = \int_{\mccan} \lhp{\smallcc}{p_\shy - q, y - x} d{\mu}_\shy(q),
\end{equation}
for $y \in  \mccan$. Here $x \neq p_\shy$ is a fixed auxiliary point in the middle of a component $C_u$ of $\curve$, and $\lhp{\smallcc}{\cdot, \, \cdot }$ is the last component of the hybrid height pairing on $\curve$.

\smallskip

Recall that we denote by $\Delta_{\delta}(p^{e}_u)$ the open disc of radius $\delta$ around $p^e_u$, with respect to the local coordinate $z^e_u$, and by $\inn\Delta_\delta(p^e_u)$ the corresponding punctured disc (see \eqref{eq:OpenPunctures}).

\smallskip

The following convergence holds.

\begin{lem} \label{lem:EssentialConvergence} Suppose that $t \in \loghyb_\pi^{-1}(\shy)$ converges to $\shy$. Then:
\begin{itemize}
\item [(i)] The functions $f_t$ converge pointwise to $f$. That is, $f_t(y) \to f (y)$ for all $y \in \mccan \setminus \{ p_\shy\}$.
\item [(ii)] The convergence is uniform on every subset $A \subseteq \mccan \setminus \{p_\shy\}$ which does not intersect the punctured discs $\inn \Delta_\delta(p^e_u) \subseteq C_u$ around the attachment points $p^e_u$, $u \in V$, $e \sim u$ for some $\delta > 0$.
\item [(iii)] The regularizations satisfy that $f_t^\reg(y) \to  f^\reg(y)$ uniformly for all $y \in \mccan \setminus \{ p_\shy\}$.
\end{itemize}
\end{lem}
Since the functions $f_t$ and $f$ have an integral representation, Lemma~\ref{lem:EssentialConvergence} essentially follows from the weak convergence of measures $\~{\mu}_t \to \mu_\shy$ on $\Sigma$. However, the integrands in \eqref{eq:DefAuxFcts} are not continuous on $\Sigma$ (they have logarithmic poles), and hence we need to additionally use the bounds of Lemma~\ref{lem:CanonicalMeasureBound} on the canonical measure.

\begin{proof}
As a preparatory step, observe that, for given $\varepsilon >0$ we can find a radius $r > 0$ with
\begin{equation} \label{eq:BasicIntBound}
\int_{0 \le |z| \le 2r} (1 + |\log|z||) i dz \wedge d \bar z \le \varepsilon, \quad \int_{\varrho_e(t) \le |z| \le 2r} (1 + |\log|z||) \Big (1 + \frac{\varrho_e(t)}{|z|^2} \Big) i dz \wedge d\bar z \le \varepsilon
\end{equation}
for all $t \in B^\ast$ sufficiently  close to $\shy$ in $B^\hyb$. Here, $\varrho_e(t)$ is the auxiliary function from \eqref{eq:Defrhoe} for some edge $e \in E$.

\smallskip

In order to prove (i), note that for fixed $y \in \Sigma$, the function $q \in \Sigma \mapsto  \lhp{\smallcc}{p_\shy - q, y - x}$ is continuous on $\Sigma$, except for logarithmic poles at $y$, $x$ and the attachment points $p^e_u$. For given $\varepsilon >0$, using \eqref{eq:BasicIntBound} and Lemma~\ref{lem:CanonicalMeasureBound}, we can find small neighborhoods $U_x$ and $U_y$ of $x$ and $y$, and punctured discs $\inn \Delta_\delta (p^e_u)$ around the points $p^e_u$ such that on $U := U_x \cup U_y \cup \bigcup \inn \Delta_\delta (p^e_u)$,
\[
\int_{U} | \lhp{\smallcc}{p_\shy - q, y - x}| \, d|\widetilde{\mu}_t|(q) \le \varepsilon,  \qquad  \int_{U} | \lhp{\smallcc}{p_\shy - q, y - x}| \, d|{\mu}_\shy|(q) \le \varepsilon
\]
uniformly for all $t \in B^\ast$ close to $\shy$ in $B^\hyb$. The first claim (i) now follows from the continuity of the map $q  \mapsto  \lhp{\smallcc}{p_\shy - q, y - x}$ on the complement $\Sigma \setminus U $ and the weak convergence  $\~\mu_t \to \mu_\shy$.

\smallskip

The claim in (ii) is deduced by standard arguments and we omit the details.

\smallskip

In order to prove the third claim, taking into account (ii), it suffices to prove the following statement: for given $\varepsilon >0$ and an attachment point $p^e_u$, we can find a punctured disc $\inn \Delta_\delta(p^e_u)$ such that
\begin{equation} \label{eq:ToDo}
\Big | f^\reg_t(y_1) - f_t^\reg(y_0) \Big| = \Big| \int_\Sigma \hp{C_u}{D^\smallcc_{p_\shy -q, u}, y_1 - y_0} \, d\~{\mu}_t \Big | \le \varepsilon
\end{equation}
for all $y_0, y_1 \in \inn \Delta_\delta(p^e_u)$ and all $t \in B^\ast$ close to $\shy$ in $B^\hyb$. Fix first a preliminary radius $r >0$ satisfying \eqref{eq:BasicIntBound} and consider the punctured disc $\inn \Delta_{2r} (p^e_u)$. Suppose that $y_1, y_0$ belong to a punctured disc $\inn \Delta_{\delta} (p^e_u)$ with radius $\delta \ll r$. Let $w \in \inn \Delta_{2r} (p^e_u)$ be an auxiliary point with adapted coordinates $|z^e_u(w)| = r$.  For every point $q \in \Sigma$, we decompose the degree zero divisor $D^\smallcc_{p_\shy - q, u} $ on the component $C_u$ of $p^e_u$ as
\[
D^\smallcc_{p_\shy - q, u} = D_q + a(q) (p^e_u -w) := (D^\smallcc_{p_\shy - q, u} - a(q) (p^e_u -w)) +  a(q) (p^e_u -w),
\]
where $a(q) := D^\smallcc_{p_\shy - q, u}(p^e_u)$ is the coefficient at the point $p^e_u$. If $\delta$ is sufficiently small, then
\[
|a(q)|  \big |\hp{C_u}{p^e_u - w, y_1 - y_0} - \log|z^e_u(y_1)| + \log|z^e_u(y_0)| \big | \le \varepsilon.
\]
for all $y_1, y_0 \in \inn \Delta_{\delta} (p^e_u)$. Since moreover the coefficient at the pole $p^e_u$ of $f$ is equal to $\~{c}_t(p, p^e_u) = \int_\Sigma a(q) \, d\~{\mu}_t(q)$, it follows that
\[
\Big | f^\reg_t(y_1) - f_t^\reg(y_0) \Big| \le  \Big | \int_\Sigma \hp{C_u}{D_q, y_1 - y_0} \, d \~{\mu}(q) \Big | + \varepsilon
\]
for all $y_1, y_0 \in \inn \Delta_\delta(p^e_u)$. However, upon choosing $\delta$ sufficiently small, $|\hp{C_u}{D_q, y_1 - y_0} | \le \varepsilon$ uniformly for $q \in \Sigma \setminus \inn \Delta_{r}(p^e_u)$ (since the divisor $D_q$ is supported bounded away from $\inn \Delta_\delta(p^e_u)$). Moreover, the assumption \eqref{eq:BasicIntBound} on the preliminary radius $r>0$ implies that
\[
\int_{\inn \Delta_{r}(p^e_u)} | \hp{C_u}{D_q, y_1 - y_0}  \, d |\~{\mu}_t|(q) \le \varepsilon 
\]
uniformly for all $y_1, y_0 \in \inn \Delta_\delta(p^e_u)$, if $t \in B^\ast$ is sufficiently close to $\shy$. Combining all these statements, we have proven \eqref{eq:ToDo}.
\end{proof}

After these preparations, we finish the proof of Theorem~\ref{thm:MainArakelovDetails}.
\begin{proof}[Proof of Theorem~\ref{thm:MainArakelovDetails}(ii)] The convergence of the coefficients at the logarithmic poles was proved in Lemma~\ref{lem:ConvergenceLogCoeff}. The convergence of the graph pieces $\gritilde{t, j}(p_t, y)$, $j \in [r]$ follows from Theorem~\ref{thm:GreenFunctionConvergence}.  Indeed the metric graphs $\mgr_t$, obtained by equipping $G$ with edge lengths $\ell_t(e) := - \log|z_e(t)|$, converge tamely in $\mgtropcombin{\grind{G}}$ to the underlying tropical curve $\curve^\trop$ of $\curve$. Moreover, the pushouts of the measures $\~{\mu}_t$ from $\curve$ to $\mgr_t$ converge weakly to the pushout of the measure $\mu_\thy$ from $\curve$ to $\curve^\trop$.

It remains to prove the uniform convergence of the regularized complex parts. Note that $\gritildereg{t, \smallcc}(p_\shy, y) = f_t^\reg(y) - I_t$ and $\grireg{\smallcc}(p_t, \cdot) = f^\reg -  I$ on $\mccan \setminus \{ p_\shy\}$, where $I_t$ and $I$ are the integrals
\begin{align*}
&I_t := \int_\mccan  f_t(\eta) \, d \~{\mu}_t(\eta) & I := \int_\mccan  f (\eta) \, d \mu_\shy(\eta).
\end{align*}

By Lemma~\ref{lem:EssentialConvergence}(iii), it suffices to prove that $I_t \to I$.  The functions $f$ and $f_t$ are continuous on $\Sigma$, except for logarithmic poles at the attachment points $p^e_u$ and at the point $p_\shy$, and the respective coefficients remain bounded. Using the bounds from Lemma~\ref{lem:CanonicalMeasureBound}, we deduce that for every $\varepsilon >0$, one can choose punctured discs $\inn \Delta_\delta(p^e_u)$ around the points $p^e_u$ and a small open neighborhood $U$ around $p_\shy$ such that
\[
\int_{U \cup \bigcup \inn \Delta_\delta(p^e_u)} |f_t(\eta)|  \, d |\~{\mu}_t |(\eta) + \int_{U \cup \bigcup \inn \Delta_\delta(p^e_u)} |f (\eta)|  \, d |{\mu}_\shy |(\eta)  \le \varepsilon
\]
uniformly for $t \in \loghyb_\pi^{-1}(\shy)$ close to $\shy$ in $B^\hyb$. On the complement $K := \mccan \setminus (U \cup \bigcup \inn \Delta_\delta(p^e_u))$, the functions $f_t$ converge uniformly to $f$ by Lemma~\ref{lem:EssentialConvergence}(ii). Hence, we have $I_t \to I$ by the weak convergence of the measures $\~{\mu}_t$ to $\mu_\shy$.
\end{proof}

\subsection{Proof of Theorem~\ref{thm:MainArakelovDetails2}}  \label{ss:ProofMainArakelovDetails2}
In the following, we give the proof of Theorem~\ref{thm:MainArakelovDetails2}. We use the notation introduced in the respective Section~\ref{ss:MainArakelovDetails2}.

(i) By the connection between the height pairing and the Poisson equation (see \eqref{eq:HPJFuncSmooth}) on the smooth Riemann surface $\rsf_t$, $t \in B^\ast$,
\begin{equation} \label{eq:TwoVariablesProof1}
\gri{t}(p,y) - \gri{t}(p',y) =  \hp{\rsf_t}{p - p', y- x} - \int_{\rsf_t} \hp{\rsf_t}{p - p', \eta - x} d \mu_t(\eta)
\end{equation}
for $p,p', x,y \in \rsf_t$. Hence Theorem~\ref{thm:MainArakelovDetails2}(i) reduces to the case where $p$ lies on a smooth section (Theorem~\ref{thm:MainArakelovDetails}) and the convergence of the height pairing (Theorem~\ref{thm:FinalTamenessHeightPairing}).

\smallskip

More precisely, fix two disjoint sections $p_t$ and $x_t$, $t \in B$, of $\rsf \to B$ with the adaptedness properties from Theorem~\ref{lem:CrucialLemma}. By the above equality and Theorem~\ref{thm:FinalTamenessHeightPairing},
\begin{align*}
\gri{t}(p,y) - \gri{t}(p_t,y) &= \sum_{j} L_j(t) \Big( \lhp{j}{\Lognoind(p) - \Lognoind(p_t), \Lognoind(y) - \Lognoind(x_t)} \\& \qquad \qquad  - \int_{\rsf_t} \lhp{j}{\Lognoind(p) - \Lognoind(p_t), \Lognoind(\eta) - \Lognoind(x_t)} \, d \mu(\eta)\Big ) \\
& \qquad \qquad +  \varepsilon_\pi(p,y) + \int_{\rsf_t} \varepsilon_\pi(p,\eta) d \mu(\eta) + o(1)
\end{align*}
uniformly. The integral over $\varepsilon_\pi$ vanishes in the limit by \eqref{eq:IntegralEpsilons}. By the connection between the hybrid height pairing and Poisson equation (Theorem~\ref{thm:HybridHPvsJFunction}), we have the analogue of \eqref{eq:TwoVariablesProof1},
\begin{align*}
\grihat{t, j}(p, y) - \grihat{t, j}(p, y) &=  \lhp{j}{\Lognoind(p) - \Lognoind(p_t), \Lognoind(y) - \Lognoind(x_t)} \\ 
&\qquad \qquad - \int_{\rsf_t} \lhp{j}{\Lognoind(p) - \Lognoind(p_t), \Lognoind(\eta) - \Lognoind(x_t)} \, d \mu_t(\eta).
\end{align*}
Applying Theorem~\ref{thm:MainArakelovDetails}(i), the proof of Theorem~\ref{thm:MainArakelovDetails2}(i) is complete.

\smallskip

(ii) The convergence $\~c_t(p, p^e_u) \to c_\shy (p, p^e_u)$ was established already in \ref{lem:ConvergenceLogCoeff}. The uniform convergence of the graph functions $\gritilde{t, j}(\cdot, \cdot)$ to $\gri{\shy, j}(\cdot, \cdot)$ for $j \in [r]$ follows from Theorem~\ref{thm:GreenFunctionConvergence} (see the proof of Theorem~\ref{thm:MainArakelovDetails}(ii)). In order to treat the complex parts, fix two points $p_0 \neq x_0$ on the smooth part of $\rsf^\hyb_\shy$. By similar steps as in (i),
\[
\grireg{\shy, \smallcc} (p, y)  - \gritildereg{t, \smallcc} (p, y)  = \grireg{\shy, \smallcc} (p_0, y)  - \gritildereg{t, \smallcc} (p_0, y) + f^\reg(p) - f^\reg_t(p)
\]
for general points $p,y$ on the fiber $\rsf^\hyb_\thy$. Here, $f_t^\reg(p)$ and $f_t^\reg(p)$ are the regularizations of
\begin{align*}
&f(p) = \int_{\rsf^\hyb_\shy}  \lhp{\smallcc}{x_0 - \eta, p - p_0 } \, d\mu_\shy(\eta), & f_t(p) = \int_{\rsf^\hyb_\shy}  \lhp{\smallcc}{x_0 - \eta, p - p_0 } \, d\~{\mu}_t (\eta).
\end{align*}
The claim then follows from Theorem~\ref{thm:MainArakelovDetails}(i) and Lemma~\ref{lem:EssentialConvergence}(iii).

\qed

\subsection{Proof of Theorem~\ref{thm:WeakTamenessArakelovGreenFunction}}  \label{ss:ProofMainGreenFunction} We begin by recalling the setting of Theorem~\ref{thm:WeakTamenessArakelovGreenFunction}. Recall that we consider a stable marked Riemann surface $S_0$ with $\nmark \ge 1$ markings $p_1, \dots, p_\nmark$, and an associated versal deformation family $\rsf / B$ equipped with adapted coordinates. Let $\rsf^\hyb / B^\hyb$ be the corresponding family of hybrid curves.

\smallskip

We fix one of the markings $p := r_i$ of $S_0$ and consider the associated section $p_\thy$, $\thy \in B^\hyb$, of the hybrid family $\rsf^\hyb \to B^\hyb$. For each base point $\thy \in B^\hyb$, introduce the hybrid function
\[
\lf_\thy := \lgri{\thy}(p_\thy, \cdot),
\]
on the hybrid curve $\rsf^\hyb_\thy$. That is, $\lf_\thy$ is the hybrid function on the hybrid curve $\rsf^\hyb_\thy$ solving the hybrid Poisson equation \eqref{eq:green_general_measure_hybrid} for the point $p_\thy$ and the canonical measure $\lmu_\thy$ on $\rsf^\hyb_\thy$. In the following we prove that the family $\mathscr{F} = (\lf_\thy)_{\thy \in B^\hyb}$ is weakly tame in the sense of Section~\ref{ss:HybridWeakTame}. 

\smallskip

Recall first that, in the preceding section, we have proven that
\begin{equation} \label{eq:StartWeakTamenessArakelov}
\sup_{q \in \rsf_t} \Big | f_t(q) - \lf_{t, \pi}^\ast(q) \Big | \to 0,
\end{equation}
whenever $t \in B^\ast$ converges tamely to a point $\shy$ of the form $\shy = (0, l)$ lying in a hybrid stratum $D_\pi^\hyb$. Here, for any point $t \in B^\ast$ and hybrid stratum $D_\pi^\hyb$, we denote by $\lf_{t, \pi}$ the hybrid solution of the Poisson equation \eqref{eq:green_general_measure_hybrid} on the hybrid curve $\rsf^\hyb_\thy$, $\thy := \loghyb_\pi(t)$, for the point $p = \loghyb_\pi(p_t)$ and the measure on $\rsf^\hyb_\thy$ obtained by pushing out the canonical measure $\mu_t$ on $\rsf_t$ via the log map $\loghyb_\pi \colon \rsf_t \to \rsf^\hyb_\pi$. Moreover, $\lf_{t, \pi}^\ast$ is the pullback of the hybrid function $\lf_{t, \pi}$ from the hybrid curve $\rsf^\hyb_\thy$ to the smooth Riemann surface $\rsf_t$.

\smallskip 

On the other hand, we can clearly carry out similar considerations locally around general hybrid points $\shy = (l,s)$. Using Theorem~\ref{thm:FinalTamenessHeightPairing} and the same steps of Section~\ref{ss:ProofClaimA}, one readily verifies that \eqref{eq:StartWeakTamenessArakelov} holds true whenever $t \in B^\ast$ converges tamely to a hybrid point $\shy = (l,s)$, which belongs to the closure $\bar D_\pi^\hyb$ of the stratum $D_\pi^\hyb$ (here, $\pi = (\pi_1, \dots, \pi_r)$ is a general ordered partition of full sedentarity on some subset $F \subseteq E$).

\smallskip

We will associate to every point in the base $\thy \in B^\hyb$ a hybrid function $\layh_\thy$ on the hybrid curve $\rsf^\hyb_\shy$, $\shy = \loghyb(\thy)$, with the properties in Section~\ref{ss:HybridWeakTame}. Recall that the domain of definition $\mathscr{U}$ of the hybrid log map $\loghyb$ is the disjoint union $\mathscr U = \bigsqcup \inn R_\pi$ of the regions $\inn R_\pi := \loghyb^{-1}(D_\pi^\hyb)$. 

In the following, we define the functions $\layh_\thy$ for $\thy$ in a fixed region $\inn R_\pi$. For a base point $t$ in $D_\pi^\hyb = \inn R_\pi \setminus B^\ast$, we simply set $\layh_\thy := \lf_\thy$. Bearing in mind \eqref{eq:StartWeakTamenessArakelov}, we would like to ensure that
\begin{equation} \label{eq:MiddleWeakTamenessArakelov}
\layh_t^\ast = \lf_{t, \pi}^\ast, \qquad \text{ for all } t \in B^\ast \cap \inn R_\pi.
\end{equation}
Given $t \in  B^\ast \cap \inn R_\pi$, consider the two points $\shy := \loghyb(t)$ and $\shy' := \loghyb_\pi(t)$ in the hybrid stratum $D_\pi^\hyb$. The underlying complex points $s, s' \in B$ of $\shy = (l, s)$ and $\shy' = (l',s')$, respectively, satisfy
\begin{align*}
s_e = s'_e = 0 \text{ for $e \in F$},  \qquad \arg(s_e) = \arg(s'_e) \text{ for $e \in E \setminus F$}, \qquad s_e = s'_e \text{ for $e \in [N] \setminus E$},
\end{align*}
in the adapted coordinates on the polydisc $B = \Delta^N$. In particular, applying the construction in Section~\ref{sec:LogSurfacesHybridCurves}, we get a natural homeomorphism $\psi^\shy_{\shy’} \colon \rsf^\hyb_\shy \to \rsf^\hyb_{\shy’}$. 

For $t \in \inn R_\pi \cap B^\ast$, we define the hybrid function $\layh_t$ on $\rsf^\hyb_\shy$ as the composition $ \layh_t := \lf_{t, \pi} \circ \psi^{\shy}_{\shy'}$  (by composing all components $h_{\thy, j}$, $j \in \{1, \dots, r, \smallcc\}$, viewed as functions on the metrized complex $\rsf^\hyb_\shy$, with the map $\psi^\shy_{\shy’}$). Since $\psi^\shy_{\shy'} \circ \loghyb = \loghyb_\pi$ on $\rsf_t$, the property \eqref{eq:MiddleWeakTamenessArakelov} holds. 

It remains to prove that the functions $h_\thy$, $\thy \in B^\hyb$ satisfy the conditions in Section~\ref{ss:HybridWeakTame}. The continuity of the map $\thy \mapsto \layh_\thy$ on each region $\inn R_\pi$ follows easily by the same steps as in Section~\ref{ss:ClaimB}. Since there are only finitely many regions $\inn R_\pi$, it remains to prove \eqref{eq:Differencev} in case that $t \in B^\ast \cap \inn R_\pi$  converges tamely to a point $\shy \in \partial _\infty B^\hyb$. However, by the continuity of $\loghyb$, the limit point $\shy$ belongs to the closure $\bar D_\pi^\hyb$ of the stratum $D_\pi^\hyb$ and hence \eqref{eq:Differencev} follows from \eqref{eq:StartWeakTamenessArakelov} and \eqref{eq:MiddleWeakTamenessArakelov}.

\qed


\subsection{Proofs of Theorem~\ref{thm:MainGeneralMeasures} and Theorem~\ref{thm:WeakConvergence}} \label{ss:ProofMainGeneralMeasures} The proofs of Theorem~\ref{thm:MainGeneralMeasures} and Theorem~\ref{thm:WeakConvergence} are analogous to the proofs of Theorem~\ref{thm:MainArakelovDetails} and Theorem~\ref{thm:WeakTamenessArakelovGreenFunction}. 

\smallskip

Namely, note that these proofs were based only on the integral representation \eqref{eq:GreenIntegral}, the weak continuity of the family of canonical measures $\lmu^{\can}_\thy$, $\thy \in B^\hyb$,  on $\rsf^\hyb \to B^\hyb$ and the bounds in Lemma ~\ref{lem:CanonicalMeasureBound}. They did not use other specific properties of the canonical measure. 

\smallskip

For two families of measures $\lmu_\thy$ and $\lnu_\thy$, $\thy \in B^\hyb$, as in Section~\ref{ss:SurfacesPoissonGeneralLimit}, the solutions $f_t$ of the respective Poisson equations on smooth fibers $\rsf_t$, $t \in B^\ast$, have the integral representation
\begin{equation}
f_t(y) = \int_{\rsf_t} \hp{\rsf_t}{p-q, x-y} d\mu_t(q) - \int_{\rsf_t} \int_{\rsf_t} \hp{\rsf_t}{p-q, x- \eta} d\mu_t(q) d \nu_t(\eta), \qquad y \in  \rsf_t,
\end{equation}
where $p, x$ are two distinct auxiliary points on the smooth Riemann surface $\rsf_t$. Moreover, under the assumptions of Theorem~\ref{thm:MainGeneralMeasures} and Theorem~\ref{thm:WeakConvergence}, the families of measures $\lmu_\thy$ and $\lnu_\thy$, $\thy \in B^\hyb$, are weakly continuous and satisfy the properties (1)-(2)-(3) in Section~\ref{ss:SurfacesPoissonGeneralLimit}, which are the precise analogue of the bounds in Lemma ~\ref{lem:CanonicalMeasureBound}. Hence the steps in Section \ref{ss:ProofMainArakelovDetails} and Section \ref{ss:ProofMainGreenFunction} carry over to this setting and give the proofs of Theorem~\ref{thm:MainGeneralMeasures} and Theorem~\ref{thm:WeakConvergence}.

\qed

\addtocontents{toc}{\SkipTocEntry}
\section*{Acknowledgments}%

We thank Raman Sanyal for pointing to us references relative to the polytope constructions in Section~\ref{ss:Preliminariespermutohedra} and Gregory Berkolaiko for hints with respect to the literature in Remark~\ref{rem:SpectralTheory}.

It is a pleasure to thank Spencer Bloch, Jos\'e Burgos Gil and Javier Fres\'an for previous collaboration and discussions on related topics. Special thanks to Hern\'an Iriarte for ongoing collaboration on higher rank non-Archimedean geometry. We also thank Aleksey Kostenko for useful discussions and helpful remarks in context with the spectral theory of weighted metric graphs.

The authors are grateful to the hospitality of the mathematics department CMLS at \'Ecole Polytechnique, where the majority of this research was carried out.  
\smallskip

O.A. thanks CNRS and ANR project ANR-18-CE40-0009. 

N.N. acknowledges financial support by the Austrian Science Fund (FWF) under Grant No. J 4497.


\bibliographystyle{alpha}
\bibliography{bibliography}

\end{document}

%% file: example3-2.tikz
{\pgfkeys{/pgf/fpu/.try=false}%
\ifx\XFigwidth\undefined\dimen1=0pt\else\dimen1\XFigwidth\fi
\divide\dimen1 by 5875
\ifx\XFigheight\undefined\dimen3=0pt\else\dimen3\XFigheight\fi
\divide\dimen3 by 3626
\ifdim\dimen1=0pt\ifdim\dimen3=0pt\dimen1=3946sp\dimen3\dimen1
  \else\dimen1\dimen3\fi\else\ifdim\dimen3=0pt\dimen3\dimen1\fi\fi
\tikzpicture[x=+\dimen1, y=+\dimen3]
{\ifx\XFigu\undefined\catcode`\@11
\def\temp{\alloc@1\dimen\dimendef\insc@unt}\temp\XFigu\catcode`\@12\fi}
\XFigu3946sp
\ifdim\XFigu<0pt\XFigu-\XFigu\fi
\pgfdeclarearrow{
  name = xfiga0,
  parameters = {
    \the\pgfarrowlinewidth \the\pgfarrowlength \the\pgfarrowwidth},
  defaults = {
	  line width=+7.5\XFigu, length=+120\XFigu, width=+60\XFigu},
  setup code = {
    \dimen7 2.15\pgfarrowlength\pgfmathveclen{\the\dimen7}{\the\pgfarrowwidth}
    \dimen7 2\pgfarrowwidth\pgfmathdivide{\pgfmathresult}{\the\dimen7}
    \dimen7 \pgfmathresult\pgfarrowlinewidth
    \pgfarrowssettipend{+\dimen7}
    \pgfarrowssetbackend{+-\pgfarrowlength}
    \dimen9 -0.5\pgfarrowlinewidth
    \pgfarrowssetvisualbackend{+\dimen9}
    \pgfarrowssetlineend{+-0.5\pgfarrowlinewidth}
    \pgfarrowshullpoint{+\dimen7}{+0pt}
    \pgfarrowsupperhullpoint{+-\pgfarrowlength}{+0.5\pgfarrowwidth}
    \pgfarrowssavethe\pgfarrowlinewidth
    \pgfarrowssavethe\pgfarrowlength
    \pgfarrowssavethe\pgfarrowwidth
  },
  drawing code = {\pgfsetdash{}{+0pt}
    \ifdim\pgfarrowlinewidth=\pgflinewidth\else\pgfsetlinewidth{+\pgfarrowlinewidth}\fi
    \pgfpathmoveto{\pgfqpoint{-\pgfarrowlength}{0.5\pgfarrowwidth}}
    \pgfpathlineto{\pgfqpoint{0pt}{0pt}}
    \pgfpathlineto{\pgfqpoint{-\pgfarrowlength}{-0.5\pgfarrowwidth}}
    \pgfusepathqstroke
  }
}
\definecolor{blue3}{rgb}{0,0,0.82}
\definecolor{red3}{rgb}{0.82,0,0}
\clip(6747,-4803) rectangle (12622,-1177);
\tikzset{inner sep=+0pt, outer sep=+0pt}
\pgfsetfillcolor{blue3}
\pgftext[base,left,at=\pgfqpointxy{9525}{-3375}] {\fontsize{27}{32.4}\normalfont $\pi_2$}
\pgfsetbeveljoin
\pgfsetlinewidth{+45\XFigu}
\pgfsetdash{}{+0pt}
\pgfsetstrokecolor{black}
\draw (7289,-1327)--(7316,-1342)--(7344,-1358)--(7373,-1375)--(7402,-1393)--(7431,-1412)
  --(7461,-1432)--(7492,-1452)--(7523,-1474)--(7555,-1497)--(7587,-1521)--(7620,-1545)
  --(7653,-1571)--(7686,-1597)--(7719,-1624)--(7752,-1652)--(7785,-1680)--(7817,-1708)
  --(7849,-1737)--(7881,-1767)--(7912,-1796)--(7942,-1825)--(7971,-1855)--(7999,-1884)
  --(8026,-1914)--(8052,-1943)--(8077,-1971)--(8101,-2000)--(8123,-2028)--(8144,-2055)
  --(8164,-2083)--(8183,-2109)--(8200,-2136)--(8217,-2162)--(8232,-2188)--(8247,-2217)
  --(8262,-2246)--(8275,-2275)--(8286,-2304)--(8297,-2334)--(8306,-2363)--(8315,-2393)
  --(8322,-2424)--(8328,-2455)--(8333,-2486)--(8336,-2517)--(8339,-2549)--(8341,-2580)
  --(8342,-2612)--(8342,-2644)--(8341,-2676)--(8339,-2708)--(8337,-2740)--(8333,-2771)
  --(8330,-2802)--(8325,-2833)--(8321,-2864)--(8316,-2894)--(8310,-2924)--(8304,-2953)
  --(8298,-2982)--(8292,-3011)--(8286,-3039)--(8280,-3068)--(8274,-3097)--(8267,-3127)
  --(8261,-3159)--(8254,-3190)--(8248,-3223)--(8241,-3256)--(8235,-3289)--(8229,-3323)
  --(8223,-3358)--(8217,-3393)--(8211,-3429)--(8206,-3465)--(8201,-3501)--(8197,-3537)
  --(8193,-3572)--(8190,-3608)--(8187,-3643)--(8185,-3677)--(8183,-3711)--(8183,-3744)
  --(8183,-3776)--(8183,-3806)--(8185,-3836)--(8187,-3864)--(8190,-3892)--(8194,-3918)
  --(8199,-3943)--(8204,-3967)--(8211,-3991)--(8220,-4019)--(8231,-4047)--(8243,-4074)
  --(8257,-4100)--(8272,-4125)--(8288,-4150)--(8305,-4175)--(8322,-4199)--(8341,-4222)
  --(8359,-4245)--(8378,-4268)--(8397,-4289)--(8415,-4310)--(8432,-4330)--(8449,-4349)
  --(8464,-4368)--(8478,-4386)--(8491,-4403)--(8502,-4419)--(8511,-4435)--(8519,-4451)
  --(8525,-4467)--(8529,-4483)--(8532,-4499)--(8532,-4515)--(8531,-4531)--(8529,-4548)
  --(8525,-4564)--(8519,-4581)--(8511,-4598)--(8503,-4614)--(8492,-4630)--(8481,-4646)
  --(8468,-4661)--(8455,-4675)--(8441,-4688)--(8426,-4701)--(8410,-4712)--(8394,-4722)
  --(8379,-4731)--(8362,-4738)--(8346,-4745)--(8330,-4750)--(8314,-4754)--(8296,-4757)
  --(8278,-4759)--(8260,-4759)--(8241,-4758)--(8222,-4756)--(8203,-4752)--(8183,-4747)
  --(8163,-4741)--(8144,-4733)--(8124,-4725)--(8105,-4715)--(8087,-4704)--(8069,-4693)
  --(8052,-4681)--(8036,-4669)--(8021,-4656)--(8007,-4643)--(7995,-4629)--(7983,-4616)
  --(7972,-4602)--(7964,-4592)--(7957,-4581)--(7950,-4570)--(7943,-4559)--(7936,-4546)
  --(7929,-4533)--(7922,-4518)--(7915,-4502)--(7907,-4484)--(7899,-4465)--(7891,-4445)
  --(7882,-4422)--(7872,-4398)--(7862,-4371)--(7851,-4343)--(7839,-4313)--(7827,-4281)
  --(7814,-4247)--(7800,-4210)--(7785,-4172)--(7770,-4132)--(7753,-4091)--(7736,-4047)
  --(7717,-4001)--(7698,-3952)--(7678,-3902)--(7656,-3849)--(7633,-3793)--(7618,-3758)
  --(7603,-3721)--(7587,-3684)--(7570,-3644)--(7553,-3604)--(7536,-3563)--(7518,-3520)
  --(7499,-3476)--(7480,-3430)--(7460,-3384)--(7440,-3336)--(7419,-3287)--(7398,-3237)
  --(7377,-3187)--(7355,-3135)--(7333,-3082)--(7311,-3029)--(7289,-2975)--(7266,-2921)
  --(7244,-2866)--(7221,-2811)--(7199,-2755)--(7176,-2700)--(7154,-2645)--(7132,-2590)
  --(7111,-2535)--(7090,-2481)--(7069,-2428)--(7049,-2375)--(7029,-2323)--(7010,-2272)
  --(6992,-2221)--(6974,-2172)--(6957,-2124)--(6941,-2078)--(6925,-2032)--(6910,-1988)
  --(6896,-1946)--(6883,-1904)--(6871,-1864)--(6860,-1826)--(6849,-1789)--(6839,-1753)
  --(6831,-1719)--(6823,-1686)--(6816,-1654)--(6809,-1621)--(6803,-1590)--(6799,-1559)
  --(6795,-1531)--(6793,-1503)--(6791,-1477)--(6791,-1452)--(6792,-1429)--(6793,-1406)
  --(6796,-1385)--(6800,-1365)--(6804,-1347)--(6810,-1329)--(6817,-1313)--(6825,-1299)
  --(6833,-1285)--(6843,-1273)--(6854,-1262)--(6865,-1253)--(6878,-1244)--(6891,-1237)
  --(6905,-1232)--(6920,-1227)--(6935,-1224)--(6951,-1222)--(6968,-1221)--(6986,-1221)
  --(7004,-1222)--(7022,-1224)--(7041,-1228)--(7060,-1232)--(7080,-1237)--(7099,-1243)
  --(7120,-1249)--(7140,-1257)--(7161,-1265)--(7181,-1274)--(7202,-1283)--(7224,-1293)
  --(7245,-1304)--(7267,-1315)--cycle;
\pgfsetlinewidth{+30\XFigu}
\pgfsetdash{}{+0pt}
\draw (7688,-2930) arc[start angle=+158.4, end angle=+40.3, radius=+208.6];
\draw (7040,-1870) arc[start angle=+-141.6, end angle=+-38.4, radius=+487.8];
\pgfsetlinewidth{+15\XFigu}
\pgfsetdash{}{+0pt}
\draw (8325,-4425)--(8475,-4575);
\draw (8475,-4425)--(8325,-4575);
\draw (7875,-2100)--(8025,-2250);
\draw (8025,-2100)--(7875,-2250);
\draw (10125,-2025)--(10275,-2175);
\draw (10275,-2025)--(10125,-2175);
\draw (10650,-4425)--(10800,-4575);
\draw (10800,-4425)--(10650,-4575);
\pgfsetlinewidth{+45\XFigu}
\pgfsetdash{}{+0pt}
\pgfsetstrokecolor{red3}
\draw (7924,-2165) arc[start angle=+111.82, end angle=+71.12, radius=+3304];
\pgfsetlinewidth{+30\XFigu}
\pgfsetdash{}{+0pt}
\pgfsetstrokecolor{black}
\draw (10752,-2047) arc[start angle=+122.6, end angle=+57.4, radius=+382.7];
\draw (10575,-2047) arc[start angle=+-108.45, end angle=+-71.55, radius=+1210.2];
\draw (11812,-2518) arc[start angle=+150.3, end angle=+29.7, radius=+237.8];
\draw (11695,-2400) arc[start angle=+-143.6, end angle=+-46.8, radius=+434.9];
\draw (10575,-4285) arc[start angle=+119.77, end angle=+60.23, radius=+890.3];
\draw (10340,-4227) arc[start angle=+-109.80, end angle=+-70.20, radius=+2000.5];
\draw (7216,-1988) arc[start angle=+110.38, end angle=+82.30, radius=+1101.4];
\pgfsetlinewidth{+45\XFigu}
\pgfsetdash{}{+0pt}
\pgfsetstrokecolor{blue3}
\draw (8395,-4521) arc[start angle=+-112.65, end angle=+-67.35, radius=+3060.5];
\draw (11930,-2872) arc[start angle=+36.94, end angle=+-36.94, radius=+1175.5];
\pgfsetarrows{[line width=15\XFigu, length=180\XFigu]}
\pgfsetarrowsend{xfiga0}
\pgfsetlinewidth{+15\XFigu}
\pgfsetdash{{+90\XFigu}{+90\XFigu}}{++0pt}
\pgfsetstrokecolor{black}
\draw (9900,-3450) arc[start angle=+117.08, end angle=+195.66, radius=+1034.2];
\draw (9975,-3450) arc[start angle=+114.06, end angle=+57.77, radius=+2231.9];
\pgfsetlinewidth{+45\XFigu}
\pgfsetdash{}{+0pt}
\draw  (10988,-4285) ellipse [x radius=+1355,y radius=+471];
\draw (12313,-2940)--(12289,-2956)--(12264,-2971)--(12240,-2984)--(12214,-2997)--(12189,-3009)
  --(12162,-3019)--(12136,-3029)--(12109,-3037)--(12082,-3044)--(12055,-3050)
  --(12028,-3054)--(12000,-3057)--(11973,-3059)--(11947,-3059)--(11921,-3058)
  --(11895,-3056)--(11870,-3052)--(11846,-3047)--(11823,-3041)--(11801,-3033)
  --(11780,-3024)--(11760,-3014)--(11742,-3003)--(11724,-2991)--(11708,-2977)
  --(11693,-2963)--(11679,-2947)--(11666,-2930)--(11653,-2911)--(11641,-2890)
  --(11629,-2868)--(11619,-2844)--(11608,-2820)--(11598,-2793)--(11588,-2766)
  --(11578,-2738)--(11568,-2709)--(11557,-2680)--(11546,-2650)--(11534,-2620)
  --(11522,-2591)--(11508,-2561)--(11494,-2533)--(11478,-2506)--(11462,-2479)
  --(11444,-2454)--(11424,-2431)--(11404,-2410)--(11381,-2390)--(11357,-2372)
  --(11332,-2356)--(11304,-2343)--(11275,-2331)--(11243,-2322)--(11216,-2316)
  --(11188,-2312)--(11159,-2308)--(11128,-2306)--(11095,-2304)--(11061,-2304)
  --(11026,-2304)--(10988,-2306)--(10950,-2307)--(10910,-2310)--(10868,-2313)
  --(10826,-2317)--(10783,-2321)--(10739,-2326)--(10694,-2330)--(10649,-2335)
  --(10604,-2340)--(10559,-2345)--(10514,-2350)--(10469,-2354)--(10425,-2359)
  --(10382,-2363)--(10340,-2366)--(10299,-2369)--(10259,-2371)--(10220,-2373)
  --(10183,-2373)--(10147,-2373)--(10113,-2372)--(10080,-2370)--(10049,-2367)
  --(10020,-2363)--(9992,-2358)--(9966,-2351)--(9946,-2345)--(9928,-2339)--(9911,-2331)
  --(9894,-2323)--(9879,-2314)--(9864,-2305)--(9850,-2294)--(9838,-2283)--(9826,-2271)
  --(9815,-2259)--(9806,-2246)--(9797,-2233)--(9790,-2219)--(9784,-2204)--(9779,-2189)
  --(9775,-2174)--(9773,-2158)--(9772,-2142)--(9773,-2126)--(9774,-2110)--(9778,-2093)
  --(9782,-2076)--(9788,-2060)--(9796,-2043)--(9805,-2027)--(9815,-2011)--(9827,-1995)
  --(9840,-1979)--(9854,-1964)--(9870,-1949)--(9888,-1934)--(9906,-1920)--(9926,-1907)
  --(9948,-1894)--(9970,-1882)--(9994,-1870)--(10019,-1859)--(10045,-1848)--(10073,-1839)
  --(10102,-1830)--(10132,-1822)--(10164,-1814)--(10196,-1808)--(10231,-1802)
  --(10265,-1797)--(10300,-1793)--(10337,-1789)--(10376,-1786)--(10416,-1784)
  --(10457,-1783)--(10499,-1782)--(10543,-1782)--(10589,-1782)--(10635,-1784)
  --(10683,-1786)--(10732,-1788)--(10783,-1792)--(10834,-1796)--(10886,-1800)
  --(10939,-1806)--(10993,-1812)--(11048,-1818)--(11102,-1825)--(11158,-1833)
  --(11213,-1842)--(11269,-1851)--(11325,-1860)--(11380,-1871)--(11435,-1881)
  --(11490,-1892)--(11544,-1904)--(11597,-1916)--(11650,-1928)--(11702,-1941)
  --(11752,-1954)--(11801,-1967)--(11849,-1981)--(11896,-1995)--(11941,-2009)
  --(11985,-2023)--(12027,-2038)--(12068,-2052)--(12107,-2067)--(12145,-2082)
  --(12180,-2097)--(12214,-2113)--(12247,-2128)--(12277,-2143)--(12306,-2159)
  --(12333,-2175)--(12364,-2194)--(12393,-2214)--(12419,-2234)--(12443,-2254)
  --(12465,-2274)--(12485,-2295)--(12503,-2317)--(12519,-2339)--(12533,-2361)
  --(12545,-2383)--(12555,-2406)--(12564,-2429)--(12570,-2452)--(12575,-2476)
  --(12577,-2499)--(12578,-2523)--(12578,-2546)--(12576,-2570)--(12572,-2593)
  --(12566,-2616)--(12560,-2639)--(12551,-2661)--(12542,-2683)--(12531,-2705)
  --(12520,-2726)--(12507,-2746)--(12494,-2766)--(12480,-2786)--(12465,-2804)
  --(12449,-2822)--(12433,-2839)--(12417,-2856)--(12400,-2872)--(12383,-2887)
  --(12366,-2901)--(12348,-2915)--(12331,-2928)--cycle;
\pgfsetfillcolor{red3}
\pgftext[base,left,at=\pgfqpointxy{8625}{-1800}] {\fontsize{27}{32.4}\normalfont $\pi_1$}
\pgfsetarrowsend{}
\pgfsetlinewidth{+30\XFigu}
\pgfsetdash{}{+0pt}
\draw (7452,-2872) arc[start angle=+-110.67, end angle=+-60.52, radius=+906.5];
\endtikzpicture}%

%% file: example2-2.tikz
{\pgfkeys{/pgf/fpu/.try=false}%
\ifx\XFigwidth\undefined\dimen1=0pt\else\dimen1\XFigwidth\fi
\divide\dimen1 by 14612
\ifx\XFigheight\undefined\dimen3=0pt\else\dimen3\XFigheight\fi
\divide\dimen3 by 5287
\ifdim\dimen1=0pt\ifdim\dimen3=0pt\dimen1=3946sp\dimen3\dimen1
  \else\dimen1\dimen3\fi\else\ifdim\dimen3=0pt\dimen3\dimen1\fi\fi
\tikzpicture[x=+\dimen1, y=+\dimen3]
{\ifx\XFigu\undefined\catcode`\@11
\def\temp{\alloc@1\dimen\dimendef\insc@unt}\temp\XFigu\catcode`\@12\fi}
\XFigu3946sp
\ifdim\XFigu<0pt\XFigu-\XFigu\fi
\catcode`\@11
\pgfutil@ifundefined{pgf@pattern@name@xfigp0}{
\pgfdeclarepatternformonly{xfigp0}
{\pgfqpoint{-1bp}{-1bp}}{\pgfqpoint{9bp}{5bp}}{\pgfqpoint{8bp}{4bp}}
{	\pgfsetdash{}{0pt}\pgfsetlinewidth{0.45bp}
	\pgfpathqmoveto{-1bp}{4.5bp}\pgfpathqlineto{9bp}{-0.5bp}
	\pgfusepathqstroke
}
}{}
\catcode`\@12
\pgfdeclarearrow{
  name = xfiga0,
  parameters = {
    \the\pgfarrowlinewidth \the\pgfarrowlength \the\pgfarrowwidth},
  defaults = {
	  line width=+7.5\XFigu, length=+120\XFigu, width=+60\XFigu},
  setup code = {
    \dimen7 2.15\pgfarrowlength\pgfmathveclen{\the\dimen7}{\the\pgfarrowwidth}
    \dimen7 2\pgfarrowwidth\pgfmathdivide{\pgfmathresult}{\the\dimen7}
    \dimen7 \pgfmathresult\pgfarrowlinewidth
    \pgfarrowssettipend{+\dimen7}
    \pgfarrowssetbackend{+-\pgfarrowlength}
    \dimen9 -0.5\pgfarrowlinewidth
    \pgfarrowssetvisualbackend{+\dimen9}
    \pgfarrowssetlineend{+-0.5\pgfarrowlinewidth}
    \pgfarrowshullpoint{+\dimen7}{+0pt}
    \pgfarrowsupperhullpoint{+-\pgfarrowlength}{+0.5\pgfarrowwidth}
    \pgfarrowssavethe\pgfarrowlinewidth
    \pgfarrowssavethe\pgfarrowlength
    \pgfarrowssavethe\pgfarrowwidth
  },
  drawing code = {\pgfsetdash{}{+0pt}
    \ifdim\pgfarrowlinewidth=\pgflinewidth\else\pgfsetlinewidth{+\pgfarrowlinewidth}\fi
    \pgfpathmoveto{\pgfqpoint{-\pgfarrowlength}{0.5\pgfarrowwidth}}
    \pgfpathlineto{\pgfqpoint{0pt}{0pt}}
    \pgfpathlineto{\pgfqpoint{-\pgfarrowlength}{-0.5\pgfarrowwidth}}
    \pgfusepathqstroke
  }
}
\definecolor{red3}{rgb}{0.82,0,0}
\clip(1018,-11714) rectangle (15630,-6427);
\tikzset{inner sep=+0pt, outer sep=+0pt}
\pgfsetfillcolor{black}
\pgftext[base,left,at=\pgfqpointxy{10875}{-6975}] {\fontsize{28}{33.6}\normalfont $\Psi$}
\pgfsetlinewidth{+30\XFigu}
\pgfsetstrokecolor{black}
\draw (5878,-7405) arc[start angle=+135.0, end angle=+45.0, radius=+300.5];
\draw (6425,-7695) arc[start angle=+-138.4, end angle=+-41.6, radius=+511.8];
\draw (6595,-7780) arc[start angle=+130.6, end angle=+49.4, radius=+391.8];
\draw (8737,-10940) arc[start angle=+-123.59, end angle=+-26.44, radius=+714.4];
\draw (8999,-10994) arc[start angle=+146.2, end angle=+64.1, radius=+543.5];
\draw (1619,-7323) arc[start angle=+-163.3, end angle=+-66.5, radius=+424.5];
\draw (1717,-7446) arc[start angle=+105.9, end angle=+24.3, radius=+323];
\pgfsetlinewidth{+45\XFigu}
\pgfsetdash{}{+0pt}
\pgfsetstrokecolor{red3}
\draw (2775,-6900) arc[start angle=+-102.620, end angle=+-81.075, radius=+6222.7];
\pgfsetstrokecolor{black}
\draw (8025,-10275) arc[start angle=+-109.82, end angle=+-1.61, radius=+1232.6];
\pgfsetstrokecolor{red3}
\draw (5100,-7275) arc[start angle=+82.847, end angle=+97.153, radius=+9637.5];
\pgfsetstrokecolor{black}
\draw (2700,-7275) arc[start angle=+125.00, end angle=+174.87, radius=+1952.9];
\pgfsetlinewidth{+30\XFigu}
\pgfsetdash{{+120\XFigu}{+120\XFigu}}{++0pt}
\draw (2775,-6900) arc[start angle=+108.4, end angle=+228.9, radius=+220.2];
\draw (2700,-7275) arc[start angle=+-37.9, end angle=+15.3, radius=+427.6];
\draw (5100,-7275) arc[start angle=+-90.0, end angle=+90.0, radius=+150];
\draw (5100,-6975) arc[start angle=+126.9, end angle=+233.1, radius=+187.5];
\draw (1875,-8700) arc[start angle=+116.6, end angle=+206.6, radius=+335.4];
\draw (1725,-9150) arc[start angle=+-63.4, end angle=+26.6, radius=+335.4];
\draw (7875,-10500) arc[start angle=+-97.1, end angle=+29.7, radius=+151.2];
\draw (8025,-10275) arc[start angle=+82.9, end angle=+209.7, radius=+151.2];
\draw (9675,-9150) arc[start angle=+-71.6, end angle=+18.4, radius=+118.6];
\draw (9750,-9000) arc[start angle=+36.9, end angle=+270.0, radius=+93.8];
\draw (7725,-8025) arc[start angle=+187.1, end angle=+60.3, radius=+151.2];
\draw (7725,-8025) arc[start angle=+-108.4, end angle=+18.4, radius=+118.6];
\pgfsetlinewidth{+45\XFigu}
\pgfsetdash{}{+0pt}
\pgfsetstrokecolor{red3}
\draw (13800,-9300) arc[start angle=+173.45, end angle=+300.32, radius=+1151.3];
\pgfsetarrows{[line width=15\XFigu, length=180\XFigu]}
\pgfsetarrowsend{xfiga0}
\pgfsetlinewidth{+30\XFigu}
\pgfsetdash{}{+0pt}
\pgfsetstrokecolor{black}
\draw (8475,-7725) arc[start angle=+122.19, end angle=+33.41, radius=+4058.8];
\pgfsetlinewidth{+45\XFigu}
\pgfsetdash{}{+0pt}
\pgfsetfillpattern{xfigp0}{black}
\draw[pattern,preaction={fill=black}]  (13800,-9300) circle [radius=+75];
\draw[pattern,preaction={fill=black}]  (14700,-9300) circle [radius=+75];
\draw[pattern,preaction={fill=black}]  (15525,-10425) circle [radius=+75];
\pgfsetstrokecolor{red3}
\pgfsetarrowsend{}
\draw (13875,-9300)--(14625,-9300);
\draw (14775,-9375)--(15450,-10350);
\pgfsetbeveljoin
\pgfsetstrokecolor{black}
\draw (5100,-6975)--(5102,-6974)--(5108,-6973)--(5118,-6971)--(5132,-6968)--(5151,-6964)
  --(5173,-6959)--(5199,-6954)--(5226,-6948)--(5254,-6942)--(5283,-6937)--(5311,-6932)
  --(5338,-6927)--(5365,-6922)--(5392,-6918)--(5419,-6914)--(5447,-6910)--(5476,-6907)
  --(5506,-6903)--(5538,-6900)--(5558,-6898)--(5580,-6896)--(5602,-6894)--(5626,-6892)
  --(5650,-6890)--(5676,-6888)--(5702,-6886)--(5730,-6884)--(5759,-6883)--(5789,-6882)
  --(5819,-6881)--(5851,-6880)--(5884,-6880)--(5918,-6880)--(5952,-6881)--(5987,-6882)
  --(6022,-6883)--(6058,-6885)--(6094,-6887)--(6131,-6889)--(6167,-6893)--(6204,-6896)
  --(6241,-6900)--(6277,-6905)--(6313,-6910)--(6349,-6916)--(6385,-6922)--(6421,-6929)
  --(6456,-6937)--(6492,-6945)--(6527,-6953)--(6563,-6963)--(6596,-6972)--(6630,-6982)
  --(6664,-6993)--(6698,-7004)--(6733,-7016)--(6769,-7029)--(6804,-7043)--(6841,-7057)
  --(6877,-7072)--(6914,-7087)--(6952,-7103)--(6989,-7120)--(7027,-7137)--(7064,-7154)
  --(7102,-7172)--(7139,-7190)--(7176,-7209)--(7212,-7227)--(7248,-7246)--(7283,-7265)
  --(7317,-7283)--(7350,-7302)--(7383,-7320)--(7414,-7339)--(7444,-7356)--(7473,-7374)
  --(7501,-7391)--(7528,-7408)--(7553,-7424)--(7577,-7440)--(7600,-7456)--(7622,-7471)
  --(7643,-7486)--(7663,-7500)--(7694,-7524)--(7722,-7546)--(7748,-7568)--(7771,-7590)
  --(7793,-7612)--(7813,-7635)--(7831,-7658)--(7849,-7681)--(7865,-7706)--(7881,-7731)
  --(7895,-7757)--(7908,-7782)--(7920,-7806)--(7930,-7827)--(7938,-7845)--(7943,-7858)
  --(7947,-7867)--(7949,-7873)--(7950,-7875);
\draw (5100,-7275)--(5102,-7276)--(5106,-7279)--(5114,-7284)--(5127,-7292)--(5144,-7302)
  --(5166,-7316)--(5193,-7333)--(5225,-7353)--(5260,-7374)--(5298,-7398)--(5339,-7423)
  --(5380,-7448)--(5422,-7474)--(5464,-7499)--(5505,-7524)--(5545,-7548)--(5583,-7571)
  --(5620,-7593)--(5655,-7614)--(5689,-7634)--(5722,-7652)--(5753,-7670)--(5783,-7687)
  --(5812,-7703)--(5840,-7718)--(5868,-7733)--(5895,-7747)--(5923,-7761)--(5950,-7775)
  --(5978,-7789)--(6005,-7802)--(6034,-7816)--(6063,-7829)--(6092,-7843)--(6122,-7856)
  --(6152,-7870)--(6183,-7883)--(6215,-7896)--(6247,-7909)--(6280,-7922)--(6313,-7935)
  --(6347,-7948)--(6380,-7960)--(6414,-7972)--(6448,-7983)--(6482,-7994)--(6516,-8005)
  --(6549,-8015)--(6582,-8025)--(6615,-8034)--(6647,-8042)--(6679,-8050)--(6710,-8057)
  --(6741,-8063)--(6771,-8069)--(6801,-8075)--(6830,-8080)--(6859,-8084)--(6888,-8088)
  --(6916,-8091)--(6945,-8093)--(6974,-8095)--(7003,-8097)--(7033,-8098)--(7064,-8098)
  --(7095,-8098)--(7129,-8097)--(7163,-8096)--(7199,-8093)--(7237,-8091)--(7277,-8088)
  --(7318,-8084)--(7361,-8080)--(7406,-8075)--(7451,-8070)--(7496,-8065)--(7541,-8060)
  --(7585,-8054)--(7626,-8049)--(7665,-8044)--(7699,-8039)--(7728,-8035)--(7752,-8032)
  --(7771,-8029)--(7785,-8027)--(7793,-8026)--(7798,-8025)--(7800,-8025);
\draw (7872,-10525)--(7871,-10528)--(7870,-10534)--(7869,-10545)--(7866,-10561)--(7862,-10581)
  --(7859,-10605)--(7855,-10632)--(7851,-10660)--(7847,-10688)--(7844,-10715)
  --(7842,-10742)--(7841,-10767)--(7841,-10791)--(7842,-10815)--(7844,-10840)
  --(7846,-10859)--(7849,-10878)--(7852,-10898)--(7854,-10920)--(7857,-10941)
  --(7859,-10964)--(7861,-10987)--(7863,-11011)--(7865,-11035)--(7867,-11060)
  --(7870,-11085)--(7872,-11109)--(7875,-11134)--(7879,-11158)--(7884,-11182)
  --(7890,-11205)--(7897,-11228)--(7906,-11250)--(7917,-11271)--(7929,-11291)
  --(7943,-11311)--(7960,-11330)--(7976,-11345)--(7993,-11361)--(8012,-11375)
  --(8033,-11390)--(8055,-11404)--(8078,-11419)--(8102,-11433)--(8128,-11447)
  --(8155,-11462)--(8182,-11476)--(8210,-11490)--(8238,-11504)--(8267,-11518)
  --(8296,-11532)--(8325,-11545)--(8354,-11558)--(8383,-11571)--(8412,-11583)
  --(8440,-11595)--(8467,-11606)--(8495,-11616)--(8521,-11626)--(8548,-11634)
  --(8573,-11642)--(8599,-11649)--(8624,-11655)--(8651,-11660)--(8678,-11664)
  --(8705,-11667)--(8733,-11669)--(8760,-11670)--(8787,-11670)--(8814,-11669)
  --(8841,-11668)--(8868,-11667)--(8895,-11665)--(8922,-11662)--(8948,-11660)
  --(8975,-11657)--(9002,-11654)--(9029,-11650)--(9056,-11646)--(9083,-11642)
  --(9109,-11638)--(9136,-11634)--(9163,-11629)--(9189,-11624)--(9215,-11618)
  --(9242,-11611)--(9268,-11604)--(9292,-11597)--(9316,-11588)--(9340,-11579)
  --(9364,-11570)--(9388,-11560)--(9413,-11550)--(9437,-11539)--(9461,-11528)
  --(9486,-11516)--(9510,-11504)--(9535,-11492)--(9559,-11480)--(9584,-11468)
  --(9608,-11455)--(9633,-11442)--(9657,-11430)--(9680,-11416)--(9704,-11403)
  --(9727,-11390)--(9750,-11376)--(9773,-11362)--(9795,-11347)--(9817,-11332)
  --(9838,-11317)--(9859,-11301)--(9879,-11285)--(9899,-11268)--(9918,-11250)
  --(9937,-11232)--(9955,-11213)--(9973,-11194)--(9991,-11174)--(10009,-11153)
  --(10026,-11133)--(10044,-11112)--(10061,-11090)--(10078,-11069)--(10095,-11047)
  --(10112,-11025)--(10128,-11003)--(10145,-10981)--(10161,-10958)--(10177,-10936)
  --(10193,-10914)--(10208,-10892)--(10223,-10870)--(10238,-10847)--(10252,-10825)
  --(10265,-10803)--(10278,-10781)--(10290,-10758)--(10302,-10736)--(10314,-10711)
  --(10325,-10686)--(10335,-10661)--(10344,-10636)--(10353,-10610)--(10362,-10584)
  --(10370,-10557)--(10377,-10531)--(10384,-10504)--(10391,-10477)--(10398,-10450)
  --(10405,-10422)--(10411,-10395)--(10417,-10368)--(10423,-10341)--(10429,-10314)
  --(10434,-10288)--(10439,-10261)--(10444,-10235)--(10449,-10210)--(10453,-10184)
  --(10456,-10159)--(10460,-10134)--(10462,-10110)--(10464,-10084)--(10465,-10058)
  --(10465,-10032)--(10465,-10007)--(10464,-9982)--(10462,-9957)--(10461,-9932)
  --(10458,-9908)--(10456,-9883)--(10453,-9859)--(10450,-9835)--(10447,-9811)
  --(10443,-9787)--(10440,-9763)--(10436,-9739)--(10432,-9715)--(10428,-9691)
  --(10423,-9667)--(10419,-9643)--(10414,-9619)--(10408,-9596)--(10402,-9572)
  --(10395,-9549)--(10388,-9525)--(10380,-9502)--(10372,-9478)--(10364,-9454)
  --(10356,-9429)--(10347,-9404)--(10338,-9379)--(10329,-9354)--(10320,-9329)
  --(10311,-9304)--(10301,-9279)--(10292,-9254)--(10281,-9230)--(10271,-9207)
  --(10260,-9185)--(10249,-9164)--(10237,-9143)--(10224,-9125)--(10211,-9107)
  --(10196,-9091)--(10181,-9077)--(10163,-9063)--(10143,-9051)--(10121,-9041)
  --(10096,-9032)--(10068,-9024)--(10037,-9017)--(10004,-9010)--(9969,-9005)
  --(9932,-9000)--(9895,-8996)--(9858,-8992)--(9824,-8989)--(9793,-8986)--(9766,-8984)
  --(9744,-8982)--(9728,-8981)--(9718,-8981)--(9712,-8980)--(9709,-8980);
\pgfsetstrokecolor{red3}
\draw (7950,-7875)--(7952,-7876)--(7956,-7880)--(7963,-7886)--(7974,-7895)--(7990,-7908)
  --(8010,-7924)--(8036,-7945)--(8065,-7969)--(8099,-7996)--(8136,-8026)--(8175,-8057)
  --(8216,-8090)--(8258,-8124)--(8300,-8158)--(8342,-8191)--(8383,-8224)--(8423,-8255)
  --(8462,-8285)--(8499,-8314)--(8535,-8341)--(8569,-8367)--(8601,-8391)--(8632,-8414)
  --(8662,-8436)--(8691,-8456)--(8719,-8476)--(8746,-8494)--(8772,-8512)--(8798,-8530)
  --(8824,-8546)--(8850,-8563)--(8878,-8580)--(8905,-8596)--(8933,-8613)--(8961,-8629)
  --(8991,-8645)--(9021,-8661)--(9052,-8678)--(9084,-8695)--(9118,-8712)--(9154,-8729)
  --(9191,-8748)--(9230,-8766)--(9271,-8786)--(9314,-8806)--(9358,-8826)--(9403,-8846)
  --(9448,-8867)--(9492,-8887)--(9536,-8906)--(9577,-8924)--(9615,-8941)--(9649,-8956)
  --(9678,-8969)--(9703,-8979)--(9721,-8988)--(9735,-8993)--(9743,-8997)--(9748,-8999)
  --(9750,-9000);
\draw (7800,-8025)--(7802,-8026)--(7805,-8029)--(7813,-8034)--(7824,-8042)--(7839,-8053)
  --(7859,-8067)--(7885,-8085)--(7914,-8106)--(7948,-8130)--(7986,-8156)--(8027,-8185)
  --(8070,-8215)--(8115,-8246)--(8160,-8278)--(8206,-8310)--(8251,-8341)--(8295,-8372)
  --(8338,-8401)--(8380,-8430)--(8420,-8457)--(8458,-8483)--(8494,-8508)--(8529,-8531)
  --(8562,-8553)--(8594,-8574)--(8624,-8594)--(8653,-8613)--(8682,-8631)--(8709,-8649)
  --(8735,-8665)--(8761,-8681)--(8787,-8697)--(8813,-8713)--(8841,-8730)--(8869,-8746)
  --(8898,-8763)--(8927,-8779)--(8956,-8795)--(8986,-8811)--(9016,-8828)--(9048,-8845)
  --(9081,-8862)--(9115,-8879)--(9151,-8898)--(9188,-8916)--(9227,-8936)--(9268,-8956)
  --(9309,-8976)--(9351,-8996)--(9393,-9017)--(9435,-9037)--(9476,-9056)--(9514,-9074)
  --(9550,-9091)--(9581,-9106)--(9609,-9119)--(9631,-9129)--(9648,-9138)--(9661,-9143)
  --(9669,-9147)--(9673,-9149)--(9675,-9150);
\draw (1725,-9150)--(1727,-9150)--(1732,-9151)--(1740,-9151)--(1754,-9153)--(1772,-9155)
  --(1797,-9157)--(1826,-9160)--(1860,-9163)--(1898,-9167)--(1939,-9172)--(1983,-9176)
  --(2027,-9181)--(2072,-9186)--(2117,-9191)--(2161,-9196)--(2204,-9202)--(2245,-9207)
  --(2284,-9212)--(2321,-9217)--(2357,-9223)--(2391,-9228)--(2423,-9233)--(2454,-9239)
  --(2484,-9244)--(2514,-9250)--(2542,-9256)--(2570,-9262)--(2597,-9268)--(2625,-9275)
  --(2653,-9282)--(2680,-9289)--(2708,-9297)--(2737,-9305)--(2765,-9313)--(2795,-9322)
  --(2825,-9331)--(2855,-9341)--(2886,-9351)--(2918,-9361)--(2950,-9371)--(2983,-9382)
  --(3016,-9392)--(3050,-9403)--(3084,-9414)--(3118,-9425)--(3152,-9436)--(3186,-9447)
  --(3221,-9458)--(3255,-9469)--(3289,-9479)--(3323,-9490)--(3356,-9500)--(3390,-9510)
  --(3423,-9519)--(3456,-9528)--(3489,-9537)--(3522,-9546)--(3554,-9554)--(3588,-9563)
  --(3617,-9569)--(3647,-9576)--(3677,-9583)--(3708,-9590)--(3740,-9596)--(3772,-9603)
  --(3805,-9610)--(3839,-9616)--(3873,-9623)--(3908,-9630)--(3944,-9637)--(3980,-9643)
  --(4017,-9650)--(4055,-9657)--(4093,-9664)--(4131,-9671)--(4169,-9678)--(4208,-9686)
  --(4247,-9693)--(4285,-9700)--(4324,-9708)--(4362,-9715)--(4400,-9723)--(4437,-9730)
  --(4474,-9738)--(4511,-9746)--(4547,-9754)--(4583,-9762)--(4618,-9770)--(4652,-9778)
  --(4687,-9786)--(4720,-9795)--(4754,-9804)--(4788,-9813)--(4821,-9822)--(4854,-9831)
  --(4888,-9841)--(4922,-9851)--(4956,-9862)--(4991,-9873)--(5026,-9884)--(5061,-9896)
  --(5097,-9908)--(5133,-9920)--(5170,-9933)--(5207,-9946)--(5244,-9959)--(5282,-9972)
  --(5319,-9985)--(5357,-9999)--(5395,-10013)--(5432,-10026)--(5470,-10040)--(5507,-10053)
  --(5544,-10066)--(5580,-10079)--(5616,-10092)--(5651,-10105)--(5686,-10117)
  --(5721,-10129)--(5754,-10141)--(5787,-10152)--(5820,-10163)--(5852,-10174)
  --(5883,-10184)--(5914,-10194)--(5945,-10203)--(5975,-10213)--(6009,-10223)
  --(6044,-10232)--(6078,-10242)--(6113,-10251)--(6148,-10260)--(6183,-10269)
  --(6219,-10278)--(6255,-10286)--(6291,-10294)--(6328,-10302)--(6364,-10309)
  --(6401,-10316)--(6438,-10323)--(6475,-10330)--(6511,-10336)--(6547,-10342)
  --(6583,-10348)--(6618,-10353)--(6652,-10358)--(6686,-10362)--(6719,-10366)
  --(6751,-10370)--(6782,-10373)--(6812,-10376)--(6841,-10379)--(6869,-10381)
  --(6897,-10383)--(6924,-10385)--(6950,-10386)--(6975,-10388)--(7006,-10389)
  --(7037,-10390)--(7067,-10390)--(7097,-10391)--(7127,-10391)--(7157,-10391)
  --(7187,-10390)--(7216,-10390)--(7245,-10389)--(7274,-10388)--(7302,-10387)
  --(7329,-10386)--(7356,-10385)--(7381,-10383)--(7406,-10382)--(7429,-10381)
  --(7451,-10380)--(7472,-10378)--(7492,-10378)--(7511,-10377)--(7528,-10376)
  --(7545,-10375)--(7560,-10375)--(7575,-10375)--(7599,-10375)--(7621,-10376)
  --(7642,-10378)--(7661,-10380)--(7680,-10383)--(7697,-10387)--(7713,-10391)
  --(7727,-10395)--(7740,-10400)--(7752,-10405)--(7762,-10410)--(7772,-10415)
  --(7780,-10420)--(7788,-10425)--(7798,-10432)--(7807,-10439)--(7817,-10447)
  --(7828,-10457)--(7840,-10467)--(7852,-10479)--(7863,-10489)--(7871,-10496)
  --(7874,-10499)--(7875,-10500);
\draw (1875,-8700)--(1876,-8701)--(1880,-8702)--(1886,-8704)--(1896,-8708)--(1910,-8714)
  --(1928,-8722)--(1952,-8731)--(1980,-8743)--(2013,-8756)--(2051,-8772)--(2093,-8789)
  --(2139,-8807)--(2188,-8827)--(2240,-8848)--(2294,-8869)--(2349,-8891)--(2404,-8913)
  --(2460,-8935)--(2515,-8957)--(2570,-8979)--(2624,-9000)--(2677,-9020)--(2728,-9040)
  --(2779,-9059)--(2827,-9077)--(2874,-9094)--(2920,-9111)--(2965,-9128)--(3008,-9143)
  --(3051,-9158)--(3092,-9172)--(3133,-9186)--(3173,-9200)--(3213,-9213)--(3253,-9226)
  --(3293,-9238)--(3332,-9251)--(3372,-9263)--(3413,-9275)--(3449,-9286)--(3487,-9297)
  --(3524,-9308)--(3563,-9318)--(3602,-9329)--(3641,-9340)--(3681,-9351)--(3722,-9362)
  --(3764,-9373)--(3806,-9384)--(3849,-9395)--(3893,-9407)--(3937,-9418)--(3982,-9429)
  --(4028,-9440)--(4074,-9452)--(4121,-9463)--(4168,-9474)--(4216,-9486)--(4264,-9497)
  --(4312,-9509)--(4360,-9520)--(4409,-9531)--(4458,-9542)--(4506,-9553)--(4555,-9564)
  --(4603,-9575)--(4651,-9586)--(4699,-9596)--(4746,-9607)--(4793,-9617)--(4840,-9627)
  --(4886,-9637)--(4932,-9647)--(4977,-9657)--(5022,-9666)--(5066,-9675)--(5109,-9685)
  --(5152,-9694)--(5195,-9703)--(5237,-9711)--(5279,-9720)--(5321,-9729)--(5363,-9738)
  --(5408,-9747)--(5453,-9756)--(5499,-9766)--(5545,-9775)--(5591,-9785)--(5637,-9794)
  --(5683,-9804)--(5730,-9813)--(5777,-9823)--(5824,-9832)--(5872,-9842)--(5920,-9852)
  --(5967,-9861)--(6015,-9871)--(6063,-9880)--(6111,-9890)--(6159,-9899)--(6206,-9909)
  --(6253,-9918)--(6299,-9927)--(6345,-9936)--(6391,-9945)--(6435,-9954)--(6479,-9962)
  --(6522,-9970)--(6564,-9978)--(6605,-9986)--(6645,-9993)--(6684,-10001)--(6722,-10007)
  --(6758,-10014)--(6794,-10021)--(6828,-10027)--(6861,-10033)--(6893,-10038)
  --(6924,-10044)--(6954,-10049)--(6983,-10053)--(7010,-10058)--(7038,-10063)
  --(7081,-10069)--(7122,-10076)--(7162,-10082)--(7201,-10087)--(7238,-10092)
  --(7274,-10097)--(7309,-10102)--(7342,-10106)--(7375,-10110)--(7406,-10113)
  --(7435,-10117)--(7464,-10120)--(7490,-10123)--(7516,-10125)--(7540,-10128)
  --(7562,-10131)--(7583,-10133)--(7603,-10135)--(7621,-10138)--(7639,-10140)
  --(7655,-10142)--(7671,-10145)--(7686,-10147)--(7700,-10150)--(7719,-10154)
  --(7738,-10159)--(7756,-10164)--(7776,-10170)--(7796,-10177)--(7817,-10184)
  --(7840,-10193)--(7864,-10203)--(7890,-10214)--(7916,-10225)--(7942,-10237)
  --(7966,-10248)--(7987,-10257)--(8004,-10265)--(8015,-10270)--(8022,-10274)
  --(8025,-10275);
\pgfsetstrokecolor{black}
\draw (2784,-6926)--(2782,-6923)--(2778,-6917)--(2770,-6907)--(2758,-6891)--(2742,-6870)
  --(2723,-6845)--(2701,-6817)--(2677,-6787)--(2652,-6756)--(2627,-6725)--(2603,-6696)
  --(2580,-6669)--(2558,-6644)--(2537,-6621)--(2517,-6601)--(2498,-6583)--(2480,-6566)
  --(2462,-6552)--(2444,-6539)--(2427,-6527)--(2409,-6517)--(2394,-6509)--(2379,-6501)
  --(2363,-6495)--(2347,-6489)--(2330,-6483)--(2312,-6479)--(2294,-6476)--(2276,-6473)
  --(2257,-6472)--(2237,-6471)--(2217,-6472)--(2196,-6474)--(2175,-6478)--(2154,-6482)
  --(2133,-6488)--(2111,-6496)--(2090,-6504)--(2068,-6514)--(2046,-6525)--(2025,-6537)
  --(2004,-6551)--(1983,-6566)--(1961,-6582)--(1940,-6599)--(1919,-6618)--(1898,-6638)
  --(1882,-6654)--(1865,-6671)--(1848,-6689)--(1831,-6709)--(1814,-6729)--(1796,-6751)
  --(1778,-6773)--(1760,-6797)--(1741,-6822)--(1722,-6849)--(1703,-6876)--(1683,-6905)
  --(1664,-6935)--(1644,-6966)--(1624,-6997)--(1604,-7030)--(1585,-7064)--(1565,-7098)
  --(1546,-7133)--(1526,-7168)--(1507,-7204)--(1489,-7240)--(1471,-7277)--(1453,-7313)
  --(1435,-7350)--(1418,-7387)--(1402,-7424)--(1386,-7461)--(1371,-7498)--(1356,-7534)
  --(1341,-7571)--(1327,-7608)--(1314,-7645)--(1301,-7683)--(1289,-7718)--(1277,-7754)
  --(1265,-7790)--(1254,-7827)--(1243,-7865)--(1232,-7903)--(1221,-7942)--(1210,-7981)
  --(1200,-8021)--(1190,-8061)--(1180,-8102)--(1171,-8143)--(1161,-8185)--(1153,-8227)
  --(1144,-8269)--(1136,-8310)--(1128,-8352)--(1120,-8393)--(1113,-8434)--(1107,-8475)
  --(1101,-8515)--(1095,-8554)--(1090,-8592)--(1085,-8630)--(1080,-8666)--(1077,-8702)
  --(1073,-8736)--(1070,-8769)--(1068,-8801)--(1066,-8832)--(1064,-8861)--(1063,-8890)
  --(1062,-8917)--(1062,-8943)--(1062,-8968)--(1062,-8992)--(1063,-9021)--(1065,-9049)
  --(1067,-9076)--(1071,-9102)--(1074,-9126)--(1079,-9149)--(1084,-9171)--(1090,-9192)
  --(1096,-9212)--(1104,-9230)--(1112,-9247)--(1120,-9263)--(1129,-9278)--(1139,-9292)
  --(1149,-9304)--(1160,-9314)--(1171,-9324)--(1182,-9332)--(1194,-9339)--(1206,-9345)
  --(1218,-9349)--(1230,-9353)--(1242,-9355)--(1255,-9357)--(1267,-9357)--(1280,-9357)
  --(1293,-9356)--(1306,-9354)--(1323,-9351)--(1340,-9347)--(1358,-9341)--(1378,-9334)
  --(1398,-9326)--(1420,-9315)--(1444,-9303)--(1470,-9289)--(1497,-9274)--(1527,-9256)
  --(1558,-9237)--(1591,-9217)--(1623,-9197)--(1655,-9176)--(1685,-9157)--(1711,-9140)
  --(1732,-9126)--(1748,-9115)--(1759,-9108)--(1765,-9104)--(1768,-9102);
\pgfsetfillcolor{red3}
\pgftext[base,left,at=\pgfqpointxy{14100}{-9075}] {\fontsize{24}{28.8}\normalfont $\ell_1$}
\pgftext[base,left,at=\pgfqpointxy{15300}{-9825}] {\fontsize{24}{28.8}\normalfont $\ell_2$}
\pgftext[base,left,at=\pgfqpointxy{13725}{-10500}] {\fontsize{24}{28.8}\normalfont $\ell_3$}
\pgftext[base,left,at=\pgfqpointxy{3600}{-6975}] {\fontsize{24}{28.8}\normalfont $\ell_1$}
\pgftext[base,left,at=\pgfqpointxy{4500}{-10275}] {\fontsize{24}{28.8}\normalfont $\ell_3$}
\pgftext[base,left,at=\pgfqpointxy{9075}{-8550}] {\fontsize{24}{28.8}\normalfont $\ell_2$}
\pgfsetlinewidth{+30\XFigu}
\pgfsetdash{}{+0pt}
\draw (5708,-7320) arc[start angle=+-134.22, end angle=+-58.47, radius=+626.9];
\endtikzpicture}%

%% file: example8.tikz
{\pgfkeys{/pgf/fpu/.try=false}%
\ifx\XFigwidth\undefined\dimen1=0pt\else\dimen1\XFigwidth\fi
\divide\dimen1 by 2279
\ifx\XFigheight\undefined\dimen3=0pt\else\dimen3\XFigheight\fi
\divide\dimen3 by 2350
\ifdim\dimen1=0pt\ifdim\dimen3=0pt\dimen1=3946sp\dimen3\dimen1
  \else\dimen1\dimen3\fi\else\ifdim\dimen3=0pt\dimen3\dimen1\fi\fi
\tikzpicture[x=+\dimen1, y=+\dimen3]
{\ifx\XFigu\undefined\catcode`\@11
\def\temp{\alloc@1\dimen\dimendef\insc@unt}\temp\XFigu\catcode`\@12\fi}
\XFigu3946sp
\ifdim\XFigu<0pt\XFigu-\XFigu\fi
\catcode`\@11
\pgfutil@ifundefined{pgf@pattern@name@xfigp0}{
\pgfdeclarepatternformonly{xfigp0}
{\pgfqpoint{-1bp}{-1bp}}{\pgfqpoint{9bp}{5bp}}{\pgfqpoint{8bp}{4bp}}
{	\pgfsetdash{}{0pt}\pgfsetlinewidth{0.45bp}
	\pgfpathqmoveto{-1bp}{4.5bp}\pgfpathqlineto{9bp}{-0.5bp}
	\pgfusepathqstroke
}
}{}
\catcode`\@12
\pgfdeclarearrow{
  name = xfiga0,
  parameters = {
    \the\pgfarrowlinewidth \the\pgfarrowlength \the\pgfarrowwidth},
  defaults = {
	  line width=+7.5\XFigu, length=+120\XFigu, width=+60\XFigu},
  setup code = {
    \dimen7 2.15\pgfarrowlength\pgfmathveclen{\the\dimen7}{\the\pgfarrowwidth}
    \dimen7 2\pgfarrowwidth\pgfmathdivide{\pgfmathresult}{\the\dimen7}
    \dimen7 \pgfmathresult\pgfarrowlinewidth
    \pgfarrowssettipend{+\dimen7}
    \pgfarrowssetbackend{+-\pgfarrowlength}
    \dimen9 -0.5\pgfarrowlinewidth
    \pgfarrowssetvisualbackend{+\dimen9}
    \pgfarrowssetlineend{+-0.5\pgfarrowlinewidth}
    \pgfarrowshullpoint{+\dimen7}{+0pt}
    \pgfarrowsupperhullpoint{+-\pgfarrowlength}{+0.5\pgfarrowwidth}
    \pgfarrowssavethe\pgfarrowlinewidth
    \pgfarrowssavethe\pgfarrowlength
    \pgfarrowssavethe\pgfarrowwidth
  },
  drawing code = {\pgfsetdash{}{+0pt}
    \ifdim\pgfarrowlinewidth=\pgflinewidth\else\pgfsetlinewidth{+\pgfarrowlinewidth}\fi
    \pgfpathmoveto{\pgfqpoint{-\pgfarrowlength}{0.5\pgfarrowwidth}}
    \pgfpathlineto{\pgfqpoint{0pt}{0pt}}
    \pgfpathlineto{\pgfqpoint{-\pgfarrowlength}{-0.5\pgfarrowwidth}}
    \pgfusepathqstroke
  }
}
\clip(2310,-5911) rectangle (4589,-3561);
\tikzset{inner sep=+0pt, outer sep=+0pt}
\pgfsetfillcolor{black}
\pgftext[base,left,at=\pgfqpointxy{4200}{-3900}] {\fontsize{24}{28.8}\normalfont $\ell_2$}
\pgfsetlinewidth{+45\XFigu}
\pgfsetstrokecolor{black}
\draw (3750,-5700) arc[start angle=+237.53, end angle=+122.47, radius=+977.9];
\pgfsetarrows{[line width=7.5\XFigu, length=180\XFigu]}
\pgfsetarrowsend{xfiga0}
\pgfsetlinewidth{+7.5\XFigu}
\pgfsetdash{{+60\XFigu}{+60\XFigu}}{++0pt}
\draw (2550,-4650) arc[start angle=+172.9, end angle=+299.7, radius=+453.5];
\draw (4575,-5550) arc[start angle=+-12.1, end angle=+114.8, radius=+268.5];
\draw (4524,-3988) arc[start angle=+0.86, end angle=+-83.59, radius=+805];
\pgfsetlinewidth{+45\XFigu}
\pgfsetdash{}{+0pt}
\pgfsetfillpattern{xfigp0}{black}
\draw[pattern,preaction={fill=black}]  (3750,-5700) circle [radius=+84];
\draw[pattern,preaction={fill=black}]  (3750,-4050) circle [radius=+84];
\pgfsetarrowsend{}
\draw (3750,-4050)--(3750,-5700);
\pgfsetfillcolor{black}
\pgftext[base,left,at=\pgfqpointxy{2325}{-4575}] {\fontsize{24}{28.8}\normalfont $\ell_1$}
\pgftext[base,left,at=\pgfqpointxy{4275}{-5775}] {\fontsize{24}{28.8}\normalfont $\ell_3$}
\draw (3750,-4050) arc[start angle=+57.53, end angle=+-57.53, radius=+977.9];
\endtikzpicture}%

%% file: example-adcoord.tikz
{\pgfkeys{/pgf/fpu/.try=false}%
\ifx\XFigwidth\undefined\dimen1=0pt\else\dimen1\XFigwidth\fi
\divide\dimen1 by 25347
\ifx\XFigheight\undefined\dimen3=0pt\else\dimen3\XFigheight\fi
\divide\dimen3 by 5518
\ifdim\dimen1=0pt\ifdim\dimen3=0pt\dimen1=3946sp\dimen3\dimen1
  \else\dimen1\dimen3\fi\else\ifdim\dimen3=0pt\dimen3\dimen1\fi\fi
\tikzpicture[x=+\dimen1, y=+\dimen3]
{\ifx\XFigu\undefined\catcode`\@11
\def\temp{\alloc@1\dimen\dimendef\insc@unt}\temp\XFigu\catcode`\@12\fi}
\XFigu3946sp
\ifdim\XFigu<0pt\XFigu-\XFigu\fi
\clip(-165,-7363) rectangle (25182,-1845);
\tikzset{inner sep=+0pt, outer sep=+0pt}
\pgfsetfillcolor{black}
\pgftext[base,left,at=\pgfqpointxy{10350}{-6600}] {\fontsize{44}{52.8}\normalfont $G=(V,E)$}
\pgfsetlinewidth{+45\XFigu}
\pgfsetstrokecolor{black}
\draw (7200,-2475) arc[start angle=+-104.04, end angle=+-80.05, radius=+5056.2];
\pgfsetlinewidth{+30\XFigu}
\pgfsetdash{}{+0pt}
\draw (13500,-4275) arc[start angle=+150.26, end angle=+209.74, radius=+604.7];
\draw (13500,-2475) arc[start angle=+150.26, end angle=+209.74, radius=+604.7];
\pgfsetdash{}{+0pt}
\draw (12610,-2970) arc[start angle=+-138.5, end angle=+-41.5, radius=+508.3];
\pgftext[base,left,at=\pgfqpointxy{7725}{-2475}] {\fontsize{44}{52.8}\normalfont $W_{e,s}$}
\pgftext[base,left,at=\pgfqpointxy{10650}{-2550}] {\fontsize{44}{52.8}\normalfont $Y_{v,s}$}
\pgftext[base,left,at=\pgfqpointxy{14325}{-2475}] {\fontsize{44}{52.8}\normalfont $W_{e,\cdot}$}
\pgftext[base,left,at=\pgfqpointxy{18525}{-2550}] {\fontsize{44}{52.8}\normalfont $Y_{v,\cdot}$}
\pgftext[base,left,at=\pgfqpointxy{13500}{-7159}] {\fontsize{44}{52.8}\normalfont $u$}
\pgftext[base,left,at=\pgfqpointxy{15934}{-7159}] {\fontsize{44}{52.8}\normalfont $v$}
\pgftext[base,left,at=\pgfqpointxy{21450}{-2400}] {\fontsize{44}{52.8}\normalfont $S_{0}$}
\pgftext[base,left,at=\pgfqpointxy{7425}{-4200}] {\fontsize{44}{52.8}\normalfont $W_{e',s}$}
\pgftext[base,left,at=\pgfqpointxy{11775}{-2400}] {\fontsize{44}{52.8}\normalfont $Y_{u,\cdot}$}
\pgftext[base,left,at=\pgfqpointxy{13800}{-4200}] {\fontsize{44}{52.8}\normalfont $W_{e',\cdot}$}
\pgftext[base,left,at=\pgfqpointxy{5475}{-2475}] {\fontsize{44}{52.8}\normalfont $Y_{u,s}$}
\pgftext[base,left,at=\pgfqpointxy{1800}{-4200}] {\fontsize{44}{52.8}\normalfont $W_{e',t}$}
\pgftext[base,left,at=\pgfqpointxy{-150}{-2475}] {\fontsize{44}{52.8}\normalfont $Y_{u,t}$}
\pgfsetlinewidth{+45\XFigu}
\draw (7200,-3075) arc[start angle=+98.340, end angle=+81.660, radius=+7239.1];
\pgfsetdash{}{+0pt}
\draw (13500,-3075) arc[start angle=+138.89, end angle=+221.11, radius=+912.5];
\pgfsetlinewidth{+30\XFigu}
\pgfsetdash{}{+0pt}
\draw (12778,-3055) arc[start angle=+130.6, end angle=+49.4, radius=+390.7];
\pgfsetlinewidth{+45\XFigu}
\draw (13500,-2475) arc[start angle=+-102.842, end angle=+-79.450, radius=+9251.1];
\draw (13500,-3075) arc[start angle=+104.683, end angle=+75.317, radius=+7249.2];
\draw (13500,-4275) arc[start angle=+-108.662, end angle=+-71.338, radius=+5976.8];
\draw (13500,-4875) arc[start angle=+99.000, end angle=+81.000, radius=+12225.5];
\pgfsetlinewidth{+30\XFigu}
\pgfsetdash{}{+0pt}
\draw (17218,-2510) arc[start angle=+29.87, end angle=+-25.66, radius=+643.3];
\draw (17279,-4309) arc[start angle=+30.13, end angle=+-25.73, radius=+639.8];
\pgfsetdash{}{+0pt}
\draw (17817,-4646) arc[start angle=+133.2, end angle=+51.1, radius=+389.5];
\draw (17644,-4566) arc[start angle=+-136.3, end angle=+-40.0, radius=+513.9];
\draw (17576,-2769) arc[start angle=+-132.30, end angle=+-56.04, radius=+623.5];
\draw (17750,-2848) arc[start angle=+137.6, end angle=+46.4, radius=+296.3];
\draw (9942,-4646) arc[start angle=+133.2, end angle=+51.1, radius=+389.5];
\pgfsetdash{}{+0pt}
\draw (9404,-4309) arc[start angle=+30.13, end angle=+-25.73, radius=+639.8];
\pgfsetdash{}{+0pt}
\draw (9769,-4566) arc[start angle=+-136.3, end angle=+-40.0, radius=+513.9];
\draw (3430,-4620) arc[start angle=+-138.4, end angle=+-41.6, radius=+511.8];
\draw (6310,-2970) arc[start angle=+-138.5, end angle=+-41.5, radius=+508.3];
\pgfsetlinewidth{+45\XFigu}
\pgfsetdash{}{+0pt}
\draw (7200,-3075) arc[start angle=+138.89, end angle=+221.11, radius=+912.5];
\pgfsetlinewidth{+30\XFigu}
\pgfsetdash{}{+0pt}
\draw (9875,-2848) arc[start angle=+137.6, end angle=+46.4, radius=+296.3];
\pgfsetdash{}{+0pt}
\draw (7200,-4275) arc[start angle=+150.26, end angle=+209.74, radius=+604.7];
\draw (3075,-4350) arc[start angle=+28.07, end angle=+-28.07, radius=+637.5];
\pgfsetlinewidth{+45\XFigu}
\pgfsetdash{}{+0pt}
\draw (1650,-4950) arc[start angle=+113.84, end angle=+66.16, radius=+1762.9];
\draw (1650,-4350) arc[start angle=+-113.84, end angle=+-66.16, radius=+1762.9];
\pgfsetlinewidth{+30\XFigu}
\pgfsetdash{}{+0pt}
\draw (3600,-4705) arc[start angle=+130.6, end angle=+49.4, radius=+391.8];
\draw (3425,-2820) arc[start angle=+-134.22, end angle=+-58.47, radius=+626.9];
\draw (3595,-2905) arc[start angle=+135.0, end angle=+45.0, radius=+300.5];
\pgfsetdash{}{+0pt}
\draw (1650,-4350) arc[start angle=+150.26, end angle=+209.74, radius=+604.7];
\pgfsetlinewidth{+45\XFigu}
\pgfsetdash{}{+0pt}
\draw (1650,-3150) arc[start angle=+138.89, end angle=+221.11, radius=+912.5];
\pgfsetlinewidth{+30\XFigu}
\pgfsetdash{}{+0pt}
\draw (928,-3130) arc[start angle=+130.6, end angle=+49.4, radius=+390.7];
\draw (760,-3045) arc[start angle=+-138.5, end angle=+-41.5, radius=+508.3];
\pgfsetdash{}{+0pt}
\draw (3075,-2550) arc[start angle=+28.07, end angle=+-28.07, radius=+637.5];
\pgfsetlinewidth{+45\XFigu}
\pgfsetdash{}{+0pt}
\draw (1650,-3150) arc[start angle=+113.84, end angle=+66.16, radius=+1762.9];
\draw (1650,-2550) arc[start angle=+-113.84, end angle=+-66.16, radius=+1762.9];
\pgfsetlinewidth{+30\XFigu}
\pgfsetdash{}{+0pt}
\draw (1650,-2550) arc[start angle=+150.26, end angle=+209.74, radius=+604.7];
\pgfsetdash{}{+0pt}
\draw (6478,-3055) arc[start angle=+130.6, end angle=+49.4, radius=+390.7];
\pgfsetdash{}{+0pt}
\draw (7200,-2475) arc[start angle=+150.26, end angle=+209.74, radius=+604.7];
\draw (9343,-2510) arc[start angle=+29.87, end angle=+-25.66, radius=+643.3];
\pgfsetdash{}{+0pt}
\draw (9701,-2769) arc[start angle=+-132.30, end angle=+-56.04, radius=+623.5];
\pgfsetlinewidth{+45\XFigu}
\draw (23078,-4182) arc[start angle=+34.60, end angle=+-24.66, radius=+607];
\pgfsetlinewidth{+30\XFigu}
\draw (22252,-3030) arc[start angle=+135.9, end angle=+54.2, radius=+388.9];
\draw (23817,-4646) arc[start angle=+133.2, end angle=+51.1, radius=+389.5];
\draw (23644,-4566) arc[start angle=+-136.3, end angle=+-40.0, radius=+513.9];
\draw (23576,-2769) arc[start angle=+-132.30, end angle=+-56.04, radius=+623.5];
\draw (23750,-2848) arc[start angle=+137.6, end angle=+46.4, radius=+296.3];
\pgfsetlinewidth{+45\XFigu}
\pgfsetdash{}{+0pt}
\draw (22972,-2987) arc[start angle=+143.96, end angle=+226.18, radius=+912.3];
\pgfsetlinewidth{+30\XFigu}
\pgfsetdash{}{+0pt}
\draw (22076,-2960) arc[start angle=+-133.2, end angle=+-36.7, radius=+511.7];
\draw (13732,-6580) arc[start angle=+133.64, end angle=+46.36, radius=+1679.5];
\draw (13732,-6580) arc[start angle=+-133.55, end angle=+-46.45, radius=+1682.1];
\pgfsetlinewidth{+7.5\XFigu}
\pgfsetfillcolor{cyan}
\filldraw  (16050,-6580) circle [radius=+116];
\filldraw  (13732,-6580) circle [radius=+116];
\pgfsetbeveljoin
\pgfsetlinewidth{+45\XFigu}
\pgfsetdash{}{+0pt}
\draw (13500,-2475)--(13497,-2473)--(13490,-2470)--(13479,-2464)--(13461,-2454)--(13438,-2442)
  --(13410,-2428)--(13378,-2411)--(13344,-2394)--(13309,-2376)--(13273,-2359)
  --(13239,-2343)--(13207,-2328)--(13176,-2315)--(13147,-2303)--(13121,-2293)
  --(13096,-2285)--(13072,-2278)--(13050,-2272)--(13029,-2268)--(13008,-2265)
  --(12988,-2263)--(12970,-2262)--(12953,-2261)--(12936,-2262)--(12919,-2263)
  --(12901,-2266)--(12883,-2269)--(12866,-2274)--(12848,-2279)--(12830,-2286)
  --(12812,-2294)--(12794,-2303)--(12776,-2314)--(12759,-2326)--(12741,-2339)
  --(12724,-2353)--(12708,-2369)--(12692,-2386)--(12677,-2404)--(12662,-2423)
  --(12648,-2443)--(12634,-2464)--(12621,-2487)--(12609,-2510)--(12597,-2535)
  --(12586,-2560)--(12575,-2588)--(12567,-2609)--(12559,-2632)--(12552,-2655)
  --(12544,-2680)--(12537,-2706)--(12530,-2733)--(12523,-2761)--(12517,-2791)
  --(12511,-2822)--(12504,-2853)--(12499,-2887)--(12493,-2921)--(12488,-2956)
  --(12483,-2992)--(12479,-3029)--(12475,-3067)--(12471,-3106)--(12467,-3146)
  --(12465,-3185)--(12462,-3226)--(12460,-3266)--(12459,-3307)--(12457,-3348)
  --(12457,-3388)--(12457,-3429)--(12457,-3470)--(12458,-3510)--(12459,-3550)
  --(12460,-3590)--(12462,-3630)--(12465,-3669)--(12468,-3709)--(12471,-3748)
  --(12475,-3788)--(12479,-3825)--(12483,-3862)--(12488,-3900)--(12494,-3939)
  --(12499,-3977)--(12506,-4017)--(12512,-4056)--(12519,-4096)--(12527,-4137)
  --(12535,-4178)--(12543,-4219)--(12552,-4260)--(12561,-4302)--(12571,-4343)
  --(12581,-4385)--(12591,-4426)--(12601,-4467)--(12612,-4508)--(12623,-4548)
  --(12634,-4588)--(12645,-4627)--(12657,-4664)--(12668,-4701)--(12680,-4737)
  --(12691,-4772)--(12703,-4806)--(12714,-4839)--(12725,-4870)--(12737,-4900)
  --(12748,-4928)--(12759,-4956)--(12770,-4982)--(12781,-5007)--(12791,-5031)
  --(12802,-5053)--(12813,-5075)--(12826,-5101)--(12839,-5126)--(12853,-5149)
  --(12867,-5171)--(12880,-5191)--(12894,-5210)--(12908,-5228)--(12922,-5245)
  --(12937,-5260)--(12951,-5273)--(12966,-5286)--(12980,-5296)--(12995,-5306)
  --(13009,-5314)--(13023,-5321)--(13038,-5326)--(13052,-5330)--(13065,-5332)
  --(13079,-5334)--(13092,-5334)--(13105,-5333)--(13118,-5331)--(13130,-5328)
  --(13142,-5324)--(13154,-5319)--(13165,-5313)--(13176,-5307)--(13188,-5300)
  --(13202,-5290)--(13215,-5279)--(13230,-5266)--(13244,-5251)--(13259,-5235)
  --(13275,-5216)--(13291,-5195)--(13309,-5172)--(13327,-5145)--(13347,-5117)
  --(13367,-5087)--(13388,-5055)--(13409,-5022)--(13429,-4990)--(13448,-4961)
  --(13464,-4934)--(13478,-4912)--(13488,-4895)--(13494,-4884)--(13498,-4878)
  --(13500,-4875);
\draw (17240,-3109)--(17243,-3110)--(17250,-3114)--(17262,-3119)--(17278,-3127)--(17300,-3137)
  --(17324,-3148)--(17350,-3160)--(17376,-3172)--(17402,-3184)--(17427,-3196)
  --(17450,-3206)--(17471,-3216)--(17492,-3226)--(17511,-3235)--(17531,-3244)
  --(17549,-3253)--(17569,-3261)--(17586,-3269)--(17604,-3278)--(17623,-3286)
  --(17642,-3296)--(17662,-3305)--(17682,-3315)--(17703,-3325)--(17724,-3336)
  --(17745,-3347)--(17766,-3358)--(17787,-3370)--(17808,-3381)--(17827,-3393)
  --(17846,-3404)--(17864,-3416)--(17881,-3428)--(17897,-3439)--(17911,-3450)
  --(17925,-3462)--(17938,-3473)--(17952,-3486)--(17964,-3500)--(17976,-3514)
  --(17987,-3528)--(17997,-3543)--(18006,-3559)--(18015,-3576)--(18023,-3593)
  --(18029,-3610)--(18035,-3628)--(18040,-3646)--(18044,-3664)--(18047,-3682)
  --(18049,-3699)--(18050,-3717)--(18050,-3734)--(18050,-3751)--(18049,-3768)
  --(18048,-3784)--(18046,-3800)--(18043,-3816)--(18040,-3833)--(18035,-3850)
  --(18030,-3868)--(18024,-3886)--(18017,-3905)--(18009,-3923)--(18001,-3942)
  --(17991,-3960)--(17980,-3978)--(17968,-3996)--(17955,-4013)--(17941,-4029)
  --(17927,-4045)--(17912,-4060)--(17896,-4074)--(17879,-4088)--(17861,-4101)
  --(17845,-4111)--(17828,-4121)--(17810,-4131)--(17790,-4141)--(17769,-4151)
  --(17746,-4161)--(17720,-4171)--(17692,-4182)--(17662,-4192)--(17629,-4204)
  --(17593,-4215)--(17556,-4227)--(17517,-4239)--(17478,-4251)--(17440,-4263)
  --(17403,-4274)--(17369,-4283)--(17341,-4292)--(17317,-4298)--(17300,-4303)
  --(17288,-4306)--(17282,-4308)--(17279,-4309);
\draw (17302,-4908)--(17304,-4910)--(17308,-4915)--(17316,-4923)--(17326,-4935)--(17340,-4950)
  --(17357,-4969)--(17376,-4988)--(17395,-5009)--(17415,-5029)--(17435,-5048)
  --(17455,-5067)--(17474,-5083)--(17493,-5098)--(17513,-5112)--(17534,-5126)
  --(17550,-5136)--(17567,-5146)--(17585,-5156)--(17603,-5167)--(17622,-5179)
  --(17641,-5191)--(17661,-5203)--(17681,-5216)--(17701,-5230)--(17722,-5243)
  --(17743,-5256)--(17764,-5269)--(17785,-5282)--(17807,-5294)--(17828,-5305)
  --(17850,-5314)--(17873,-5323)--(17896,-5329)--(17919,-5334)--(17943,-5337)
  --(17967,-5337)--(17992,-5336)--(18014,-5333)--(18036,-5329)--(18060,-5323)
  --(18084,-5316)--(18109,-5307)--(18134,-5298)--(18161,-5287)--(18188,-5276)
  --(18215,-5264)--(18243,-5251)--(18272,-5238)--(18300,-5224)--(18329,-5210)
  --(18357,-5195)--(18386,-5180)--(18414,-5166)--(18442,-5151)--(18469,-5136)
  --(18495,-5121)--(18521,-5106)--(18546,-5091)--(18570,-5076)--(18593,-5061)
  --(18615,-5045)--(18636,-5029)--(18656,-5013)--(18677,-4995)--(18696,-4976)
  --(18715,-4956)--(18733,-4936)--(18751,-4915)--(18768,-4894)--(18784,-4873)
  --(18800,-4851)--(18815,-4829)--(18830,-4806)--(18845,-4784)--(18859,-4761)
  --(18873,-4738)--(18887,-4715)--(18901,-4691)--(18915,-4668)--(18928,-4645)
  --(18941,-4621)--(18954,-4597)--(18967,-4573)--(18979,-4549)--(18990,-4525)
  --(19002,-4500)--(19012,-4475)--(19021,-4451)--(19029,-4427)--(19037,-4403)
  --(19045,-4378)--(19052,-4353)--(19058,-4327)--(19065,-4301)--(19071,-4275)
  --(19077,-4249)--(19083,-4222)--(19089,-4195)--(19094,-4168)--(19100,-4141)
  --(19105,-4114)--(19110,-4087)--(19115,-4060)--(19119,-4033)--(19123,-4006)
  --(19127,-3979)--(19130,-3952)--(19133,-3926)--(19135,-3899)--(19137,-3873)
  --(19138,-3847)--(19138,-3821)--(19138,-3795)--(19137,-3769)--(19135,-3743)
  --(19132,-3717)--(19128,-3691)--(19124,-3665)--(19119,-3638)--(19114,-3612)
  --(19109,-3585)--(19103,-3559)--(19096,-3532)--(19090,-3506)--(19083,-3479)
  --(19076,-3452)--(19069,-3426)--(19061,-3399)--(19054,-3373)--(19046,-3346)
  --(19038,-3321)--(19030,-3295)--(19022,-3269)--(19013,-3244)--(19004,-3220)
  --(18995,-3195)--(18986,-3171)--(18976,-3148)--(18965,-3125)--(18953,-3101)
  --(18940,-3077)--(18926,-3053)--(18912,-3030)--(18897,-3007)--(18882,-2984)
  --(18866,-2961)--(18850,-2939)--(18833,-2917)--(18816,-2894)--(18799,-2872)
  --(18781,-2850)--(18764,-2829)--(18746,-2807)--(18729,-2786)--(18711,-2765)
  --(18694,-2744)--(18676,-2724)--(18659,-2704)--(18641,-2684)--(18624,-2665)
  --(18606,-2647)--(18589,-2629)--(18571,-2612)--(18551,-2594)--(18532,-2578)
  --(18512,-2562)--(18492,-2546)--(18471,-2531)--(18451,-2517)--(18431,-2504)
  --(18410,-2490)--(18389,-2477)--(18369,-2464)--(18348,-2452)--(18327,-2440)
  --(18306,-2428)--(18285,-2416)--(18264,-2404)--(18243,-2392)--(18221,-2381)
  --(18200,-2369)--(18178,-2359)--(18156,-2348)--(18134,-2338)--(18112,-2328)
  --(18089,-2319)--(18066,-2310)--(18043,-2302)--(18019,-2293)--(17995,-2285)
  --(17970,-2277)--(17945,-2269)--(17920,-2260)--(17894,-2252)--(17869,-2244)
  --(17843,-2236)--(17817,-2228)--(17792,-2221)--(17767,-2214)--(17742,-2208)
  --(17717,-2203)--(17693,-2199)--(17670,-2196)--(17647,-2194)--(17625,-2194)
  --(17604,-2196)--(17583,-2199)--(17561,-2205)--(17539,-2213)--(17517,-2224)
  --(17495,-2238)--(17471,-2255)--(17447,-2275)--(17422,-2297)--(17396,-2322)
  --(17370,-2347)--(17344,-2374)--(17318,-2400)--(17295,-2425)--(17274,-2448)
  --(17256,-2468)--(17241,-2484)--(17231,-2496)--(17224,-2504)--(17220,-2508)
  --(17218,-2510);
\draw (9365,-3109)--(9368,-3110)--(9375,-3114)--(9387,-3119)--(9403,-3127)--(9425,-3137)
  --(9449,-3148)--(9475,-3160)--(9501,-3172)--(9527,-3184)--(9552,-3196)--(9575,-3206)
  --(9596,-3216)--(9617,-3226)--(9636,-3235)--(9656,-3244)--(9674,-3253)--(9694,-3261)
  --(9711,-3269)--(9729,-3278)--(9748,-3286)--(9767,-3296)--(9787,-3305)--(9807,-3315)
  --(9828,-3325)--(9849,-3336)--(9870,-3347)--(9891,-3358)--(9912,-3370)--(9933,-3381)
  --(9952,-3393)--(9971,-3404)--(9989,-3416)--(10006,-3428)--(10022,-3439)--(10036,-3450)
  --(10050,-3462)--(10063,-3473)--(10077,-3486)--(10089,-3500)--(10101,-3514)
  --(10112,-3528)--(10122,-3543)--(10131,-3559)--(10140,-3576)--(10148,-3593)
  --(10154,-3610)--(10160,-3628)--(10165,-3646)--(10169,-3664)--(10172,-3682)
  --(10174,-3699)--(10175,-3717)--(10175,-3734)--(10175,-3751)--(10174,-3768)
  --(10173,-3784)--(10171,-3800)--(10168,-3816)--(10165,-3833)--(10160,-3850)
  --(10155,-3868)--(10149,-3886)--(10142,-3905)--(10134,-3923)--(10126,-3942)
  --(10116,-3960)--(10105,-3978)--(10093,-3996)--(10080,-4013)--(10066,-4029)
  --(10052,-4045)--(10037,-4060)--(10021,-4074)--(10004,-4088)--(9986,-4101)
  --(9970,-4111)--(9953,-4121)--(9935,-4131)--(9915,-4141)--(9894,-4151)--(9871,-4161)
  --(9845,-4171)--(9817,-4182)--(9787,-4192)--(9754,-4204)--(9718,-4215)--(9681,-4227)
  --(9642,-4239)--(9603,-4251)--(9565,-4263)--(9528,-4274)--(9494,-4283)--(9466,-4292)
  --(9442,-4298)--(9425,-4303)--(9413,-4306)--(9407,-4308)--(9404,-4309);
\draw (3075,-3150)--(3078,-3152)--(3085,-3155)--(3096,-3161)--(3113,-3169)--(3134,-3180)
  --(3158,-3192)--(3184,-3204)--(3210,-3218)--(3236,-3230)--(3260,-3242)--(3283,-3254)
  --(3304,-3264)--(3324,-3275)--(3344,-3284)--(3362,-3294)--(3381,-3303)--(3400,-3313)
  --(3417,-3321)--(3435,-3330)--(3454,-3339)--(3473,-3349)--(3492,-3359)--(3512,-3370)
  --(3533,-3381)--(3553,-3392)--(3574,-3404)--(3595,-3416)--(3615,-3428)--(3635,-3440)
  --(3654,-3453)--(3673,-3465)--(3690,-3477)--(3707,-3490)--(3722,-3502)--(3737,-3514)
  --(3750,-3526)--(3763,-3538)--(3775,-3551)--(3787,-3565)--(3799,-3579)--(3809,-3594)
  --(3819,-3610)--(3827,-3626)--(3835,-3643)--(3843,-3661)--(3849,-3678)--(3854,-3696)
  --(3858,-3715)--(3861,-3733)--(3863,-3751)--(3865,-3768)--(3865,-3786)--(3865,-3803)
  --(3864,-3820)--(3863,-3838)--(3860,-3855)--(3857,-3873)--(3853,-3891)--(3849,-3910)
  --(3843,-3929)--(3836,-3949)--(3828,-3969)--(3818,-3989)--(3808,-4009)--(3796,-4029)
  --(3783,-4048)--(3769,-4067)--(3753,-4085)--(3737,-4103)--(3720,-4119)--(3702,-4134)
  --(3683,-4149)--(3663,-4163)--(3646,-4172)--(3629,-4182)--(3611,-4191)--(3591,-4200)
  --(3569,-4209)--(3545,-4219)--(3520,-4228)--(3491,-4237)--(3461,-4247)--(3427,-4257)
  --(3392,-4268)--(3354,-4278)--(3315,-4289)--(3276,-4299)--(3237,-4310)--(3200,-4319)
  --(3166,-4328)--(3137,-4335)--(3113,-4341)--(3096,-4345)--(3084,-4348)--(3078,-4349)
  --(3075,-4350);
\draw (3075,-4950)--(3077,-4952)--(3081,-4957)--(3088,-4966)--(3099,-4978)--(3112,-4994)
  --(3128,-5012)--(3146,-5033)--(3165,-5054)--(3185,-5074)--(3204,-5094)--(3223,-5113)
  --(3242,-5130)--(3260,-5146)--(3280,-5161)--(3300,-5175)--(3316,-5185)--(3332,-5196)
  --(3350,-5207)--(3368,-5218)--(3386,-5230)--(3405,-5243)--(3424,-5256)--(3444,-5270)
  --(3463,-5284)--(3484,-5298)--(3504,-5312)--(3524,-5326)--(3545,-5339)--(3566,-5352)
  --(3588,-5363)--(3609,-5374)--(3631,-5383)--(3654,-5390)--(3677,-5396)--(3701,-5399)
  --(3725,-5401)--(3750,-5400)--(3772,-5398)--(3795,-5394)--(3818,-5389)--(3843,-5383)
  --(3868,-5375)--(3894,-5367)--(3920,-5357)--(3948,-5347)--(3976,-5335)--(4004,-5324)
  --(4033,-5311)--(4062,-5298)--(4091,-5285)--(4120,-5272)--(4149,-5258)--(4178,-5244)
  --(4206,-5230)--(4233,-5216)--(4260,-5202)--(4287,-5188)--(4312,-5174)--(4337,-5160)
  --(4360,-5145)--(4383,-5131)--(4404,-5116)--(4425,-5100)--(4446,-5082)--(4467,-5064)
  --(4486,-5045)--(4505,-5026)--(4523,-5006)--(4541,-4985)--(4558,-4964)--(4574,-4943)
  --(4591,-4921)--(4606,-4899)--(4622,-4877)--(4637,-4855)--(4652,-4833)--(4667,-4810)
  --(4682,-4787)--(4696,-4765)--(4710,-4742)--(4724,-4718)--(4738,-4695)--(4751,-4672)
  --(4764,-4648)--(4777,-4624)--(4789,-4600)--(4800,-4575)--(4810,-4552)--(4819,-4528)
  --(4828,-4504)--(4836,-4479)--(4844,-4454)--(4852,-4429)--(4859,-4403)--(4866,-4377)
  --(4873,-4351)--(4880,-4325)--(4887,-4298)--(4893,-4271)--(4900,-4245)--(4906,-4218)
  --(4912,-4191)--(4917,-4164)--(4923,-4137)--(4928,-4110)--(4932,-4084)--(4937,-4057)
  --(4940,-4031)--(4944,-4004)--(4946,-3978)--(4948,-3952)--(4950,-3926)--(4950,-3900)
  --(4950,-3874)--(4948,-3848)--(4946,-3822)--(4944,-3796)--(4940,-3769)--(4937,-3743)
  --(4932,-3716)--(4928,-3690)--(4923,-3663)--(4917,-3636)--(4912,-3609)--(4906,-3582)
  --(4900,-3555)--(4893,-3529)--(4887,-3502)--(4880,-3475)--(4873,-3449)--(4866,-3423)
  --(4859,-3397)--(4852,-3371)--(4844,-3346)--(4836,-3321)--(4828,-3296)--(4819,-3272)
  --(4810,-3248)--(4800,-3225)--(4789,-3200)--(4777,-3176)--(4764,-3152)--(4750,-3128)
  --(4736,-3105)--(4722,-3082)--(4707,-3058)--(4691,-3035)--(4676,-3013)--(4659,-2990)
  --(4643,-2967)--(4627,-2945)--(4610,-2923)--(4593,-2901)--(4576,-2879)--(4560,-2857)
  --(4543,-2836)--(4526,-2815)--(4509,-2794)--(4493,-2774)--(4476,-2755)--(4459,-2736)
  --(4442,-2718)--(4425,-2700)--(4406,-2682)--(4387,-2664)--(4368,-2648)--(4348,-2632)
  --(4328,-2616)--(4308,-2601)--(4288,-2587)--(4268,-2573)--(4248,-2559)--(4227,-2546)
  --(4207,-2532)--(4187,-2519)--(4166,-2507)--(4145,-2494)--(4125,-2481)--(4104,-2469)
  --(4083,-2457)--(4062,-2445)--(4040,-2433)--(4019,-2422)--(3997,-2411)--(3975,-2400)
  --(3953,-2390)--(3930,-2380)--(3907,-2371)--(3883,-2362)--(3859,-2353)--(3835,-2343)
  --(3810,-2334)--(3785,-2325)--(3760,-2315)--(3734,-2306)--(3709,-2297)--(3684,-2288)
  --(3658,-2280)--(3633,-2272)--(3608,-2265)--(3584,-2259)--(3560,-2254)--(3537,-2250)
  --(3514,-2248)--(3492,-2247)--(3471,-2248)--(3450,-2250)--(3428,-2255)--(3406,-2262)
  --(3383,-2273)--(3361,-2286)--(3337,-2302)--(3312,-2321)--(3286,-2343)--(3259,-2367)
  --(3232,-2392)--(3205,-2418)--(3179,-2443)--(3155,-2467)--(3133,-2490)--(3114,-2509)
  --(3099,-2524)--(3088,-2536)--(3081,-2544)--(3077,-2548)--(3075,-2550);
\draw (1650,-2550)--(1647,-2548)--(1640,-2545)--(1629,-2539)--(1611,-2529)--(1588,-2517)
  --(1560,-2503)--(1528,-2486)--(1494,-2469)--(1459,-2451)--(1423,-2434)--(1389,-2418)
  --(1357,-2403)--(1326,-2390)--(1297,-2378)--(1271,-2368)--(1246,-2360)--(1222,-2353)
  --(1200,-2347)--(1179,-2343)--(1158,-2340)--(1138,-2338)--(1120,-2337)--(1103,-2336)
  --(1086,-2337)--(1069,-2338)--(1051,-2341)--(1033,-2344)--(1016,-2349)--(998,-2354)
  --(980,-2361)--(962,-2369)--(944,-2378)--(926,-2389)--(909,-2401)--(891,-2414)
  --(874,-2428)--(858,-2444)--(842,-2461)--(827,-2479)--(812,-2498)--(798,-2518)
  --(784,-2539)--(771,-2562)--(759,-2585)--(747,-2610)--(736,-2635)--(725,-2663)
  --(717,-2684)--(709,-2707)--(702,-2730)--(694,-2755)--(687,-2781)--(680,-2808)
  --(673,-2836)--(667,-2866)--(661,-2897)--(654,-2928)--(649,-2962)--(643,-2996)
  --(638,-3031)--(633,-3067)--(629,-3104)--(625,-3142)--(621,-3181)--(617,-3221)
  --(615,-3260)--(612,-3301)--(610,-3341)--(609,-3382)--(607,-3423)--(607,-3463)
  --(607,-3504)--(607,-3545)--(608,-3585)--(609,-3625)--(610,-3665)--(612,-3705)
  --(615,-3744)--(618,-3784)--(621,-3823)--(625,-3863)--(629,-3900)--(633,-3937)
  --(638,-3975)--(644,-4014)--(649,-4052)--(656,-4092)--(662,-4131)--(669,-4171)
  --(677,-4212)--(685,-4253)--(693,-4294)--(702,-4335)--(711,-4377)--(721,-4418)
  --(731,-4460)--(741,-4501)--(751,-4542)--(762,-4583)--(773,-4623)--(784,-4663)
  --(795,-4702)--(807,-4739)--(818,-4776)--(830,-4812)--(841,-4847)--(853,-4881)
  --(864,-4914)--(875,-4945)--(887,-4975)--(898,-5003)--(909,-5031)--(920,-5057)
  --(931,-5082)--(941,-5106)--(952,-5128)--(963,-5150)--(976,-5176)--(989,-5201)
  --(1003,-5224)--(1017,-5246)--(1030,-5266)--(1044,-5285)--(1058,-5303)--(1072,-5320)
  --(1087,-5335)--(1101,-5348)--(1116,-5361)--(1130,-5371)--(1145,-5381)--(1159,-5389)
  --(1173,-5396)--(1188,-5401)--(1202,-5405)--(1215,-5407)--(1229,-5409)--(1242,-5409)
  --(1255,-5408)--(1268,-5406)--(1280,-5403)--(1292,-5399)--(1304,-5394)--(1315,-5388)
  --(1326,-5382)--(1338,-5375)--(1352,-5365)--(1365,-5354)--(1380,-5341)--(1394,-5326)
  --(1409,-5310)--(1425,-5291)--(1441,-5270)--(1459,-5247)--(1477,-5220)--(1497,-5192)
  --(1517,-5162)--(1538,-5130)--(1559,-5097)--(1579,-5065)--(1598,-5036)--(1614,-5009)
  --(1628,-4987)--(1638,-4970)--(1644,-4959)--(1648,-4953)--(1650,-4950);
\draw (7200,-2475)--(7197,-2473)--(7190,-2470)--(7179,-2464)--(7161,-2454)--(7138,-2442)
  --(7110,-2428)--(7078,-2411)--(7044,-2394)--(7009,-2376)--(6973,-2359)--(6939,-2343)
  --(6907,-2328)--(6876,-2315)--(6847,-2303)--(6821,-2293)--(6796,-2285)--(6772,-2278)
  --(6750,-2272)--(6729,-2268)--(6708,-2265)--(6688,-2263)--(6670,-2262)--(6653,-2261)
  --(6636,-2262)--(6619,-2263)--(6601,-2266)--(6583,-2269)--(6566,-2274)--(6548,-2279)
  --(6530,-2286)--(6512,-2294)--(6494,-2303)--(6476,-2314)--(6459,-2326)--(6441,-2339)
  --(6424,-2353)--(6408,-2369)--(6392,-2386)--(6377,-2404)--(6362,-2423)--(6348,-2443)
  --(6334,-2464)--(6321,-2487)--(6309,-2510)--(6297,-2535)--(6286,-2560)--(6275,-2588)
  --(6267,-2609)--(6259,-2632)--(6252,-2655)--(6244,-2680)--(6237,-2706)--(6230,-2733)
  --(6223,-2761)--(6217,-2791)--(6211,-2822)--(6204,-2853)--(6199,-2887)--(6193,-2921)
  --(6188,-2956)--(6183,-2992)--(6179,-3029)--(6175,-3067)--(6171,-3106)--(6167,-3146)
  --(6165,-3185)--(6162,-3226)--(6160,-3266)--(6159,-3307)--(6157,-3348)--(6157,-3388)
  --(6157,-3429)--(6157,-3470)--(6158,-3510)--(6159,-3550)--(6160,-3590)--(6162,-3630)
  --(6165,-3669)--(6168,-3709)--(6171,-3748)--(6175,-3788)--(6179,-3825)--(6183,-3862)
  --(6188,-3900)--(6194,-3939)--(6199,-3977)--(6206,-4017)--(6212,-4056)--(6219,-4096)
  --(6227,-4137)--(6235,-4178)--(6243,-4219)--(6252,-4260)--(6261,-4302)--(6271,-4343)
  --(6281,-4385)--(6291,-4426)--(6301,-4467)--(6312,-4508)--(6323,-4548)--(6334,-4588)
  --(6345,-4627)--(6357,-4664)--(6368,-4701)--(6380,-4737)--(6391,-4772)--(6403,-4806)
  --(6414,-4839)--(6425,-4870)--(6437,-4900)--(6448,-4928)--(6459,-4956)--(6470,-4982)
  --(6481,-5007)--(6491,-5031)--(6502,-5053)--(6513,-5075)--(6526,-5101)--(6539,-5126)
  --(6553,-5149)--(6567,-5171)--(6580,-5191)--(6594,-5210)--(6608,-5228)--(6622,-5245)
  --(6637,-5260)--(6651,-5273)--(6666,-5286)--(6680,-5296)--(6695,-5306)--(6709,-5314)
  --(6723,-5321)--(6738,-5326)--(6752,-5330)--(6765,-5332)--(6779,-5334)--(6792,-5334)
  --(6805,-5333)--(6818,-5331)--(6830,-5328)--(6842,-5324)--(6854,-5319)--(6865,-5313)
  --(6876,-5307)--(6888,-5300)--(6902,-5290)--(6915,-5279)--(6930,-5266)--(6944,-5251)
  --(6959,-5235)--(6975,-5216)--(6991,-5195)--(7009,-5172)--(7027,-5145)--(7047,-5117)
  --(7067,-5087)--(7088,-5055)--(7109,-5022)--(7129,-4990)--(7148,-4961)--(7164,-4934)
  --(7178,-4912)--(7188,-4895)--(7194,-4884)--(7198,-4878)--(7200,-4875);
\draw (9427,-4908)--(9429,-4910)--(9433,-4915)--(9441,-4923)--(9451,-4935)--(9465,-4950)
  --(9482,-4969)--(9501,-4988)--(9520,-5009)--(9540,-5029)--(9560,-5048)--(9580,-5067)
  --(9599,-5083)--(9618,-5098)--(9638,-5112)--(9659,-5126)--(9675,-5136)--(9692,-5146)
  --(9710,-5156)--(9728,-5167)--(9747,-5179)--(9766,-5191)--(9786,-5203)--(9806,-5216)
  --(9826,-5230)--(9847,-5243)--(9868,-5256)--(9889,-5269)--(9910,-5282)--(9932,-5294)
  --(9953,-5305)--(9975,-5314)--(9998,-5323)--(10021,-5329)--(10044,-5334)--(10068,-5337)
  --(10092,-5337)--(10117,-5336)--(10139,-5333)--(10161,-5329)--(10185,-5323)
  --(10209,-5316)--(10234,-5307)--(10259,-5298)--(10286,-5287)--(10313,-5276)
  --(10340,-5264)--(10368,-5251)--(10397,-5238)--(10425,-5224)--(10454,-5210)
  --(10482,-5195)--(10511,-5180)--(10539,-5166)--(10567,-5151)--(10594,-5136)
  --(10620,-5121)--(10646,-5106)--(10671,-5091)--(10695,-5076)--(10718,-5061)
  --(10740,-5045)--(10761,-5029)--(10781,-5013)--(10802,-4995)--(10821,-4976)
  --(10840,-4956)--(10858,-4936)--(10876,-4915)--(10893,-4894)--(10909,-4873)
  --(10925,-4851)--(10940,-4829)--(10955,-4806)--(10970,-4784)--(10984,-4761)
  --(10998,-4738)--(11012,-4715)--(11026,-4691)--(11040,-4668)--(11053,-4645)
  --(11066,-4621)--(11079,-4597)--(11092,-4573)--(11104,-4549)--(11115,-4525)
  --(11127,-4500)--(11137,-4475)--(11146,-4451)--(11154,-4427)--(11162,-4403)
  --(11170,-4378)--(11177,-4353)--(11183,-4327)--(11190,-4301)--(11196,-4275)
  --(11202,-4249)--(11208,-4222)--(11214,-4195)--(11219,-4168)--(11225,-4141)
  --(11230,-4114)--(11235,-4087)--(11240,-4060)--(11244,-4033)--(11248,-4006)
  --(11252,-3979)--(11255,-3952)--(11258,-3926)--(11260,-3899)--(11262,-3873)
  --(11263,-3847)--(11263,-3821)--(11263,-3795)--(11262,-3769)--(11260,-3743)
  --(11257,-3717)--(11253,-3691)--(11249,-3665)--(11244,-3638)--(11239,-3612)
  --(11234,-3585)--(11228,-3559)--(11221,-3532)--(11215,-3506)--(11208,-3479)
  --(11201,-3452)--(11194,-3426)--(11186,-3399)--(11179,-3373)--(11171,-3346)
  --(11163,-3321)--(11155,-3295)--(11147,-3269)--(11138,-3244)--(11129,-3220)
  --(11120,-3195)--(11111,-3171)--(11101,-3148)--(11090,-3125)--(11078,-3101)
  --(11065,-3077)--(11051,-3053)--(11037,-3030)--(11022,-3007)--(11007,-2984)
  --(10991,-2961)--(10975,-2939)--(10958,-2917)--(10941,-2894)--(10924,-2872)
  --(10906,-2850)--(10889,-2829)--(10871,-2807)--(10854,-2786)--(10836,-2765)
  --(10819,-2744)--(10801,-2724)--(10784,-2704)--(10766,-2684)--(10749,-2665)
  --(10731,-2647)--(10714,-2629)--(10696,-2612)--(10676,-2594)--(10657,-2578)
  --(10637,-2562)--(10617,-2546)--(10596,-2531)--(10576,-2517)--(10556,-2504)
  --(10535,-2490)--(10514,-2477)--(10494,-2464)--(10473,-2452)--(10452,-2440)
  --(10431,-2428)--(10410,-2416)--(10389,-2404)--(10368,-2392)--(10346,-2381)
  --(10325,-2369)--(10303,-2359)--(10281,-2348)--(10259,-2338)--(10237,-2328)
  --(10214,-2319)--(10191,-2310)--(10168,-2302)--(10144,-2293)--(10120,-2285)
  --(10095,-2277)--(10070,-2269)--(10045,-2260)--(10019,-2252)--(9994,-2244)
  --(9968,-2236)--(9942,-2228)--(9917,-2221)--(9892,-2214)--(9867,-2208)--(9842,-2203)
  --(9818,-2199)--(9795,-2196)--(9772,-2194)--(9750,-2194)--(9729,-2196)--(9708,-2199)
  --(9686,-2205)--(9664,-2213)--(9642,-2224)--(9620,-2238)--(9596,-2255)--(9572,-2275)
  --(9547,-2297)--(9521,-2322)--(9495,-2347)--(9469,-2374)--(9443,-2400)--(9420,-2425)
  --(9399,-2448)--(9381,-2468)--(9366,-2484)--(9356,-2496)--(9349,-2504)--(9345,-2508)
  --(9343,-2510);
\pgfsetdash{}{+0pt}
\draw (22920,-2389)--(22921,-2389)--(22925,-2391)--(22934,-2395)--(22949,-2402)--(22966,-2411)
  --(22984,-2420)--(23001,-2430)--(23015,-2439)--(23028,-2448)--(23038,-2457)
  --(23046,-2467)--(23054,-2477)--(23059,-2487)--(23064,-2498)--(23068,-2510)
  --(23072,-2523)--(23076,-2537)--(23078,-2552)--(23080,-2568)--(23081,-2585)
  --(23081,-2602)--(23081,-2620)--(23080,-2637)--(23078,-2655)--(23076,-2672)
  --(23073,-2690)--(23069,-2705)--(23066,-2721)--(23061,-2739)--(23055,-2757)
  --(23049,-2778)--(23041,-2801)--(23032,-2826)--(23022,-2854)--(23012,-2882)
  --(23001,-2910)--(22991,-2936)--(22983,-2958)--(22977,-2973)--(22974,-2983)
  --(22972,-2986)--(22972,-2987);
\pgfsetdash{}{+0pt}
\draw (23302,-4908)--(23304,-4910)--(23308,-4915)--(23316,-4923)--(23326,-4935)--(23340,-4950)
  --(23357,-4969)--(23376,-4988)--(23395,-5009)--(23415,-5029)--(23435,-5048)
  --(23455,-5067)--(23474,-5083)--(23493,-5098)--(23513,-5112)--(23534,-5126)
  --(23550,-5136)--(23567,-5146)--(23585,-5156)--(23603,-5167)--(23622,-5179)
  --(23641,-5191)--(23661,-5203)--(23681,-5216)--(23701,-5230)--(23722,-5243)
  --(23743,-5256)--(23764,-5269)--(23785,-5282)--(23807,-5294)--(23828,-5305)
  --(23850,-5314)--(23873,-5323)--(23896,-5329)--(23919,-5334)--(23943,-5337)
  --(23967,-5337)--(23992,-5336)--(24014,-5333)--(24036,-5329)--(24060,-5323)
  --(24084,-5316)--(24109,-5307)--(24134,-5298)--(24161,-5287)--(24188,-5276)
  --(24215,-5264)--(24243,-5251)--(24272,-5238)--(24300,-5224)--(24329,-5210)
  --(24357,-5195)--(24386,-5180)--(24414,-5166)--(24442,-5151)--(24469,-5136)
  --(24495,-5121)--(24521,-5106)--(24546,-5091)--(24570,-5076)--(24593,-5061)
  --(24615,-5045)--(24636,-5029)--(24656,-5013)--(24677,-4995)--(24696,-4976)
  --(24715,-4956)--(24733,-4936)--(24751,-4915)--(24768,-4894)--(24784,-4873)
  --(24800,-4851)--(24815,-4829)--(24830,-4806)--(24845,-4784)--(24859,-4761)
  --(24873,-4738)--(24887,-4715)--(24901,-4691)--(24915,-4668)--(24928,-4645)
  --(24941,-4621)--(24954,-4597)--(24967,-4573)--(24979,-4549)--(24990,-4525)
  --(25002,-4500)--(25012,-4475)--(25021,-4451)--(25029,-4427)--(25037,-4403)
  --(25045,-4378)--(25052,-4353)--(25058,-4327)--(25065,-4301)--(25071,-4275)
  --(25077,-4249)--(25083,-4222)--(25089,-4195)--(25094,-4168)--(25100,-4141)
  --(25105,-4114)--(25110,-4087)--(25115,-4060)--(25119,-4033)--(25123,-4006)
  --(25127,-3979)--(25130,-3952)--(25133,-3926)--(25135,-3899)--(25137,-3873)
  --(25138,-3847)--(25138,-3821)--(25138,-3795)--(25137,-3769)--(25135,-3743)
  --(25132,-3717)--(25128,-3691)--(25124,-3665)--(25119,-3638)--(25114,-3612)
  --(25109,-3585)--(25103,-3559)--(25096,-3532)--(25090,-3506)--(25083,-3479)
  --(25076,-3452)--(25069,-3426)--(25061,-3399)--(25054,-3373)--(25046,-3346)
  --(25038,-3321)--(25030,-3295)--(25022,-3269)--(25013,-3244)--(25004,-3220)
  --(24995,-3195)--(24986,-3171)--(24976,-3148)--(24965,-3125)--(24953,-3101)
  --(24940,-3077)--(24926,-3053)--(24912,-3030)--(24897,-3007)--(24882,-2984)
  --(24866,-2961)--(24850,-2939)--(24833,-2917)--(24816,-2894)--(24799,-2872)
  --(24781,-2850)--(24764,-2829)--(24746,-2807)--(24729,-2786)--(24711,-2765)
  --(24694,-2744)--(24676,-2724)--(24659,-2704)--(24641,-2684)--(24624,-2665)
  --(24606,-2647)--(24589,-2629)--(24571,-2612)--(24551,-2594)--(24532,-2578)
  --(24512,-2562)--(24492,-2546)--(24471,-2531)--(24451,-2517)--(24431,-2504)
  --(24410,-2490)--(24389,-2477)--(24369,-2464)--(24348,-2452)--(24327,-2440)
  --(24306,-2428)--(24285,-2416)--(24264,-2404)--(24243,-2392)--(24221,-2381)
  --(24200,-2369)--(24178,-2359)--(24156,-2348)--(24134,-2338)--(24112,-2328)
  --(24089,-2319)--(24066,-2310)--(24043,-2302)--(24019,-2293)--(23995,-2285)
  --(23970,-2277)--(23945,-2269)--(23920,-2260)--(23894,-2252)--(23869,-2244)
  --(23843,-2236)--(23817,-2228)--(23792,-2221)--(23767,-2214)--(23742,-2208)
  --(23717,-2203)--(23693,-2199)--(23670,-2196)--(23647,-2194)--(23625,-2194)
  --(23604,-2196)--(23583,-2199)--(23561,-2205)--(23539,-2213)--(23517,-2224)
  --(23495,-2238)--(23471,-2255)--(23447,-2275)--(23422,-2297)--(23396,-2322)
  --(23370,-2347)--(23344,-2374)--(23318,-2400)--(23295,-2425)--(23274,-2448)
  --(23256,-2468)--(23241,-2484)--(23231,-2496)--(23224,-2504)--(23220,-2508)
  --(23218,-2510);
\pgfsetdash{}{+0pt}
\draw (23250,-2475)--(23247,-2478)--(23241,-2485)--(23231,-2496)--(23218,-2510)--(23204,-2527)
  --(23189,-2545)--(23175,-2562)--(23163,-2578)--(23153,-2593)--(23144,-2608)
  --(23136,-2622)--(23130,-2636)--(23125,-2650)--(23121,-2663)--(23117,-2677)
  --(23114,-2692)--(23112,-2707)--(23110,-2724)--(23108,-2741)--(23107,-2759)
  --(23107,-2777)--(23107,-2795)--(23108,-2813)--(23110,-2831)--(23112,-2849)
  --(23114,-2865)--(23117,-2882)--(23121,-2897)--(23125,-2913)--(23130,-2927)
  --(23135,-2942)--(23141,-2958)--(23148,-2974)--(23156,-2991)--(23166,-3009)
  --(23177,-3029)--(23189,-3050)--(23202,-3071)--(23215,-3093)--(23226,-3112)
  --(23236,-3128)--(23244,-3140)--(23248,-3147)--(23250,-3150);
\pgfsetdash{}{+0pt}
\draw (23240,-3109)--(23243,-3110)--(23250,-3114)--(23262,-3119)--(23278,-3127)--(23300,-3137)
  --(23324,-3148)--(23350,-3160)--(23376,-3172)--(23402,-3184)--(23427,-3196)
  --(23450,-3206)--(23471,-3216)--(23492,-3226)--(23511,-3235)--(23531,-3244)
  --(23549,-3253)--(23569,-3261)--(23586,-3269)--(23604,-3278)--(23623,-3286)
  --(23642,-3296)--(23662,-3305)--(23682,-3315)--(23703,-3325)--(23724,-3336)
  --(23745,-3347)--(23766,-3358)--(23787,-3370)--(23808,-3381)--(23827,-3393)
  --(23846,-3404)--(23864,-3416)--(23881,-3428)--(23897,-3439)--(23911,-3450)
  --(23925,-3462)--(23938,-3473)--(23952,-3486)--(23964,-3500)--(23976,-3514)
  --(23987,-3528)--(23997,-3543)--(24006,-3559)--(24015,-3576)--(24023,-3593)
  --(24029,-3610)--(24035,-3628)--(24040,-3646)--(24044,-3664)--(24047,-3682)
  --(24049,-3699)--(24050,-3717)--(24050,-3734)--(24050,-3751)--(24049,-3768)
  --(24048,-3784)--(24046,-3800)--(24043,-3816)--(24040,-3833)--(24035,-3850)
  --(24030,-3868)--(24024,-3886)--(24017,-3905)--(24009,-3923)--(24001,-3942)
  --(23991,-3960)--(23980,-3978)--(23968,-3996)--(23955,-4013)--(23941,-4029)
  --(23927,-4045)--(23912,-4060)--(23896,-4074)--(23879,-4088)--(23861,-4101)
  --(23845,-4111)--(23828,-4121)--(23810,-4131)--(23790,-4141)--(23769,-4151)
  --(23746,-4161)--(23720,-4171)--(23692,-4182)--(23662,-4192)--(23629,-4204)
  --(23593,-4215)--(23556,-4227)--(23517,-4239)--(23478,-4251)--(23440,-4263)
  --(23403,-4274)--(23369,-4283)--(23341,-4292)--(23317,-4298)--(23300,-4303)
  --(23288,-4306)--(23282,-4308)--(23279,-4309);
\pgfsetdash{}{+0pt}
\draw (23325,-4950)--(23322,-4947)--(23316,-4940)--(23306,-4929)--(23293,-4914)--(23279,-4896)
  --(23264,-4878)--(23250,-4860)--(23238,-4843)--(23228,-4827)--(23219,-4811)
  --(23211,-4795)--(23205,-4779)--(23200,-4763)--(23196,-4747)--(23192,-4731)
  --(23189,-4713)--(23187,-4695)--(23185,-4675)--(23183,-4655)--(23182,-4634)
  --(23182,-4613)--(23182,-4591)--(23183,-4570)--(23185,-4550)--(23187,-4530)
  --(23189,-4512)--(23192,-4494)--(23196,-4478)--(23200,-4463)--(23205,-4446)
  --(23211,-4430)--(23219,-4414)--(23228,-4398)--(23238,-4382)--(23250,-4365)
  --(23264,-4347)--(23279,-4329)--(23293,-4311)--(23306,-4296)--(23316,-4285)
  --(23322,-4278)--(23325,-4275);
\pgfsetdash{}{+0pt}
\draw (22920,-2389)--(22917,-2388)--(22910,-2385)--(22898,-2380)--(22880,-2372)--(22856,-2362)
  --(22826,-2350)--(22793,-2336)--(22758,-2322)--(22721,-2308)--(22684,-2294)
  --(22649,-2280)--(22615,-2269)--(22583,-2258)--(22554,-2249)--(22526,-2241)
  --(22501,-2235)--(22477,-2230)--(22454,-2226)--(22433,-2224)--(22412,-2222)
  --(22391,-2222)--(22374,-2223)--(22357,-2224)--(22340,-2226)--(22323,-2229)
  --(22305,-2233)--(22288,-2238)--(22271,-2244)--(22254,-2251)--(22236,-2259)
  --(22219,-2269)--(22202,-2279)--(22185,-2291)--(22169,-2305)--(22153,-2319)
  --(22137,-2335)--(22122,-2352)--(22108,-2370)--(22094,-2390)--(22081,-2410)
  --(22069,-2431)--(22057,-2454)--(22046,-2477)--(22036,-2502)--(22026,-2527)
  --(22017,-2554)--(22009,-2582)--(22003,-2604)--(21997,-2627)--(21992,-2651)
  --(21986,-2676)--(21981,-2703)--(21977,-2730)--(21973,-2759)--(21969,-2789)
  --(21965,-2820)--(21962,-2853)--(21959,-2886)--(21956,-2921)--(21954,-2956)
  --(21952,-2993)--(21951,-3030)--(21950,-3068)--(21950,-3107)--(21950,-3147)
  --(21951,-3187)--(21952,-3227)--(21953,-3268)--(21955,-3308)--(21958,-3349)
  --(21960,-3390)--(21964,-3430)--(21968,-3471)--(21972,-3511)--(21976,-3551)
  --(21982,-3590)--(21987,-3630)--(21993,-3669)--(21999,-3708)--(22006,-3747)
  --(22013,-3786)--(22021,-3823)--(22028,-3860)--(22036,-3897)--(22045,-3935)
  --(22054,-3973)--(22064,-4011)--(22074,-4050)--(22085,-4089)--(22096,-4129)
  --(22107,-4169)--(22119,-4209)--(22131,-4250)--(22144,-4291)--(22157,-4331)
  --(22171,-4372)--(22185,-4412)--(22199,-4452)--(22213,-4492)--(22227,-4531)
  --(22242,-4569)--(22257,-4607)--(22271,-4644)--(22286,-4679)--(22301,-4714)
  --(22315,-4748)--(22329,-4781)--(22344,-4812)--(22358,-4842)--(22372,-4871)
  --(22385,-4899)--(22399,-4925)--(22412,-4950)--(22425,-4974)--(22438,-4997)
  --(22450,-5018)--(22463,-5039)--(22478,-5064)--(22494,-5087)--(22509,-5109)
  --(22525,-5130)--(22540,-5149)--(22556,-5167)--(22571,-5183)--(22587,-5198)
  --(22602,-5212)--(22618,-5224)--(22633,-5235)--(22649,-5245)--(22664,-5253)
  --(22679,-5260)--(22694,-5265)--(22709,-5269)--(22723,-5272)--(22737,-5273)
  --(22751,-5273)--(22764,-5272)--(22777,-5270)--(22789,-5267)--(22801,-5263)
  --(22813,-5258)--(22824,-5252)--(22835,-5246)--(22845,-5238)--(22856,-5231)
  --(22869,-5219)--(22882,-5207)--(22895,-5193)--(22908,-5177)--(22921,-5159)
  --(22935,-5139)--(22950,-5117)--(22965,-5092)--(22981,-5064)--(22998,-5034)
  --(23016,-5002)--(23034,-4969)--(23052,-4935)--(23069,-4901)--(23085,-4870)
  --(23099,-4842)--(23111,-4819)--(23120,-4801)--(23125,-4790)--(23129,-4783)
  --(23130,-4780);
\pgfsetfillcolor{black}
\pgftext[base,left,at=\pgfqpointxy{19800}{-3825}] {\fontsize{36}{43.2}\normalfont $\bullet \,\,\bullet\,\, \bullet$}
\pgftext[base,left,at=\pgfqpointxy{1800}{-2475}] {\fontsize{44}{52.8}\normalfont $W_{e,t}$}
\pgftext[base,left,at=\pgfqpointxy{4500}{-2700}] {\fontsize{44}{52.8}\normalfont $Y_{v,t}$}
\pgfsetlinewidth{+30\XFigu}
\pgfsetdash{}{+0pt}
\draw (7200,-2475) arc[start angle=+28.07, end angle=+-28.07, radius=+637.5];
\draw (7200,-4275) arc[start angle=+29.74, end angle=+-29.74, radius=+604.7];
\draw (13500,-2475) arc[start angle=+28.07, end angle=+-28.07, radius=+637.5];
\draw (17218,-2510) arc[start angle=+153.60, end angle=+210.60, radius=+628.1];
\draw (9343,-2510) arc[start angle=+153.60, end angle=+210.60, radius=+628.1];
\draw (1650,-2550) arc[start angle=+28.07, end angle=+-28.07, radius=+637.5];
\draw (3065,-2546) arc[start angle=+153.60, end angle=+210.60, radius=+628.1];
\draw (3062,-4347) arc[start angle=+154.61, end angle=+209.79, radius=+647.2];
\draw (1650,-4350) arc[start angle=+29.74, end angle=+-29.74, radius=+604.7];
\draw (17279,-4309) arc[start angle=+154.61, end angle=+209.79, radius=+647.2];
\pgfsetlinewidth{+45\XFigu}
\draw (7200,-4275) arc[start angle=+-105.26, end angle=+-74.74, radius=+4275.7];
\pgfsetlinewidth{+30\XFigu}
\draw (9404,-4309) arc[start angle=+154.61, end angle=+209.79, radius=+647.2];
\pgfsetlinewidth{+45\XFigu}
\draw (7200,-4875) arc[start angle=+105.26, end angle=+74.74, radius=+4275.7];
\pgfsetlinewidth{+30\XFigu}
\draw (13500,-4275) arc[start angle=+29.74, end angle=+-29.74, radius=+604.7];
\endtikzpicture}%

%% file: example1.tikz
{\pgfkeys{/pgf/fpu/.try=false}%
\ifx\XFigwidth\undefined\dimen1=0pt\else\dimen1\XFigwidth\fi
\divide\dimen1 by 16097
\ifx\XFigheight\undefined\dimen3=0pt\else\dimen3\XFigheight\fi
\divide\dimen3 by 4670
\ifdim\dimen1=0pt\ifdim\dimen3=0pt\dimen1=3946sp\dimen3\dimen1
  \else\dimen1\dimen3\fi\else\ifdim\dimen3=0pt\dimen3\dimen1\fi\fi
\tikzpicture[x=+\dimen1, y=+\dimen3]
{\ifx\XFigu\undefined\catcode`\@11
\def\temp{\alloc@1\dimen\dimendef\insc@unt}\temp\XFigu\catcode`\@12\fi}
\XFigu3946sp
\ifdim\XFigu<0pt\XFigu-\XFigu\fi
\catcode`\@11
\pgfutil@ifundefined{pgf@pattern@name@xfigp0}{
\pgfdeclarepatternformonly{xfigp0}
{\pgfqpoint{-1bp}{-1bp}}{\pgfqpoint{9bp}{5bp}}{\pgfqpoint{8bp}{4bp}}
{	\pgfsetdash{}{0pt}\pgfsetlinewidth{0.45bp}
	\pgfpathqmoveto{-1bp}{4.5bp}\pgfpathqlineto{9bp}{-0.5bp}
	\pgfusepathqstroke
}
}{}
\catcode`\@12
\definecolor{red3}{rgb}{0.82,0,0}
\clip(285,-7883) rectangle (16382,-3213);
\tikzset{inner sep=+0pt, outer sep=+0pt}
\pgfsetfillcolor{black}
\pgftext[base,left,at=\pgfqpointxy{15000}{-6975}] {\fontsize{24}{28.8}\normalfont ${\mathit{gr}_{\pi}^{\mathfrak f}}(G)$}
\pgfsetlinewidth{+45\XFigu}
\pgfsetstrokecolor{red3}
\draw (2078,-7800)--(7274,-7800);
\pgfsetlinewidth{+7.5\XFigu}
\pgfsetstrokecolor{black}
\pgfsetfillpattern{xfigp0}{black}
\draw[pattern,preaction={fill=black}]  (5715,-6900) circle [radius=+75];
\pgfsetlinewidth{+60\XFigu}
\pgfsetstrokecolor{blue}
\draw (2078,-7800)--(3637,-6900);
\pgfsetlinewidth{+7.5\XFigu}
\pgfsetstrokecolor{black}
\draw[pattern,preaction={fill=black}]  (7274,-7800) circle [radius=+75];
\pgfsetlinewidth{+60\XFigu}
\pgfsetstrokecolor{blue}
\draw (5715,-6900)--(7274,-7800);
\pgfsetlinewidth{+7.5\XFigu}
\pgfsetstrokecolor{black}
\draw[pattern,preaction={fill=black}]  (3637,-6900) circle [radius=+75];
\pgfsetlinewidth{+60\XFigu}
\pgfsetstrokecolor{blue}
\draw (4676,-5100)--(4676,-3375);
\pgfsetlinewidth{+45\XFigu}
\pgfsetstrokecolor{black}
\draw (3637,-6900)--(4676,-5100);
\pgfsetstrokecolor{red3}
\draw (2078,-7800)--(4676,-3300);
\draw (4676,-3300)--(7274,-7800);
\pgfsetstrokecolor{black}
\draw (4676,-5100)--(5715,-6900);
\pgfsetlinewidth{+7.5\XFigu}
\draw[pattern,preaction={fill=black}]  (4676,-3300) circle [radius=+75];
\draw[pattern,preaction={fill=black}]  (4676,-5100) circle [radius=+75];
\pgfsetlinewidth{+45\XFigu}
\draw (5715,-6900)--(3637,-6900);
\pgfsetlinewidth{+7.5\XFigu}
\draw[pattern,preaction={fill=black}]  (14746,-5812) circle [radius=+37];
\draw[pattern,preaction={fill=black}]  (14236,-6095) circle [radius=+37];
\draw[pattern,preaction={fill=black}]  (15773,-5812) circle [radius=+37];
\draw[pattern,preaction={fill=black}]  (16337,-6126) circle [radius=+37];
\draw[pattern,preaction={fill=black}]  (15260,-4925) circle [radius=+37];
\draw[pattern,preaction={fill=black}]  (15260,-4509) circle [radius=+37];
\pgfsetlinewidth{+45\XFigu}
\draw (14746,-5812)--(15260,-4925);
\draw (15260,-4925)--(15773,-5812);
\draw (15773,-5812)--(14746,-5812);
\pgfsetstrokecolor{blue}
\draw (12180,-5535)--(12988,-6002);
\pgfsetlinewidth{+7.5\XFigu}
\pgfsetstrokecolor{black}
\draw[pattern,preaction={fill=black}]  (12180,-5535) circle [radius=+55];
\draw[pattern,preaction={fill=black}]  (12180,-4657) circle [radius=+55];
\pgfsetlinewidth{+45\XFigu}
\pgfsetstrokecolor{blue}
\draw (12180,-5535)--(12180,-4657);
\pgfsetlinewidth{+7.5\XFigu}
\pgfsetstrokecolor{black}
\draw[pattern,preaction={fill=black}]  (11420,-5974) circle [radius=+55];
\pgfsetlinewidth{+45\XFigu}
\pgfsetstrokecolor{blue}
\draw (11420,-5974)--(12180,-5535);
\pgfsetlinewidth{+7.5\XFigu}
\pgfsetstrokecolor{black}
\draw[pattern,preaction={fill=black}]  (12988,-6002) circle [radius=+55];
\pgfsetbeveljoin
\pgfsetlinewidth{+45\XFigu}
\pgfsetstrokecolor{red3}
\draw (9392,-5202)--(9394,-5202)--(9399,-5201)--(9409,-5199)--(9422,-5197)--(9441,-5194)
  --(9464,-5190)--(9490,-5186)--(9520,-5180)--(9551,-5174)--(9582,-5167)--(9613,-5160)
  --(9643,-5152)--(9671,-5144)--(9696,-5135)--(9718,-5125)--(9738,-5114)--(9754,-5102)
  --(9769,-5088)--(9781,-5073)--(9791,-5057)--(9799,-5040)--(9807,-5020)--(9814,-4999)
  --(9821,-4975)--(9827,-4949)--(9833,-4922)--(9838,-4893)--(9843,-4863)--(9848,-4833)
  --(9852,-4802)--(9857,-4771)--(9860,-4741)--(9863,-4712)--(9865,-4684)--(9866,-4658)
  --(9866,-4634)--(9865,-4613)--(9863,-4594)--(9859,-4578)--(9853,-4565)--(9845,-4555)
  --(9834,-4547)--(9820,-4542)--(9803,-4539)--(9784,-4539)--(9762,-4541)--(9737,-4544)
  --(9711,-4549)--(9683,-4556)--(9654,-4563)--(9625,-4571)--(9595,-4580)--(9566,-4589)
  --(9537,-4599)--(9510,-4609)--(9485,-4619)--(9462,-4631)--(9441,-4642)--(9422,-4655)
  --(9406,-4669)--(9392,-4684)--(9381,-4700)--(9372,-4717)--(9363,-4738)--(9356,-4761)
  --(9351,-4786)--(9346,-4815)--(9342,-4847)--(9338,-4881)--(9335,-4917)--(9333,-4955)
  --(9331,-4993)--(9330,-5030)--(9329,-5066)--(9328,-5098)--(9328,-5127)--(9327,-5152)
  --(9327,-5171)--(9327,-5185)--(9327,-5195)--(9327,-5200)--(9327,-5202);
\draw (9263,-5202)--(9262,-5200)--(9260,-5194)--(9257,-5185)--(9252,-5171)--(9245,-5152)
  --(9236,-5129)--(9225,-5102)--(9214,-5072)--(9201,-5040)--(9188,-5008)--(9175,-4977)
  --(9161,-4947)--(9148,-4920)--(9135,-4895)--(9122,-4873)--(9109,-4854)--(9096,-4838)
  --(9082,-4825)--(9068,-4814)--(9051,-4804)--(9033,-4797)--(9012,-4791)--(8990,-4786)
  --(8966,-4783)--(8941,-4780)--(8914,-4779)--(8887,-4778)--(8859,-4778)--(8831,-4778)
  --(8805,-4779)--(8779,-4780)--(8756,-4783)--(8735,-4786)--(8716,-4791)--(8701,-4797)
  --(8689,-4804)--(8680,-4814)--(8675,-4825)--(8672,-4838)--(8671,-4853)--(8672,-4870)
  --(8674,-4889)--(8678,-4911)--(8683,-4934)--(8688,-4958)--(8695,-4983)--(8702,-5008)
  --(8710,-5033)--(8718,-5058)--(8727,-5082)--(8736,-5105)--(8746,-5127)--(8756,-5146)
  --(8768,-5163)--(8780,-5178)--(8794,-5191)--(8809,-5202)--(8825,-5210)--(8842,-5216)
  --(8863,-5221)--(8886,-5225)--(8911,-5227)--(8940,-5228)--(8972,-5229)--(9006,-5228)
  --(9042,-5226)--(9080,-5224)--(9118,-5221)--(9155,-5218)--(9191,-5215)--(9223,-5212)
  --(9252,-5210)--(9277,-5207)--(9296,-5205)--(9310,-5204)--(9320,-5203)--(9325,-5202)
  --(9327,-5202);
\draw (9327,-5267)--(9328,-5269)--(9331,-5273)--(9337,-5280)--(9345,-5291)--(9357,-5306)
  --(9371,-5324)--(9388,-5347)--(9406,-5372)--(9427,-5400)--(9448,-5430)--(9469,-5460)
  --(9490,-5491)--(9509,-5521)--(9527,-5550)--(9542,-5578)--(9556,-5604)--(9567,-5629)
  --(9575,-5653)--(9581,-5676)--(9585,-5698)--(9586,-5720)--(9585,-5740)--(9582,-5761)
  --(9578,-5783)--(9572,-5806)--(9564,-5831)--(9555,-5856)--(9545,-5882)--(9534,-5910)
  --(9521,-5938)--(9508,-5967)--(9495,-5996)--(9481,-6025)--(9467,-6053)--(9452,-6081)
  --(9438,-6107)--(9424,-6132)--(9410,-6155)--(9397,-6176)--(9384,-6194)--(9372,-6210)
  --(9360,-6222)--(9349,-6231)--(9338,-6236)--(9327,-6238)--(9316,-6236)--(9305,-6231)
  --(9294,-6222)--(9282,-6210)--(9270,-6194)--(9257,-6176)--(9244,-6155)--(9230,-6132)
  --(9216,-6107)--(9202,-6081)--(9187,-6053)--(9173,-6025)--(9159,-5996)--(9146,-5967)
  --(9133,-5938)--(9120,-5910)--(9109,-5882)--(9099,-5856)--(9090,-5831)--(9082,-5806)
  --(9076,-5783)--(9072,-5761)--(9069,-5740)--(9068,-5720)--(9069,-5698)--(9073,-5676)
  --(9079,-5653)--(9087,-5629)--(9098,-5604)--(9112,-5578)--(9127,-5550)--(9145,-5521)
  --(9164,-5491)--(9185,-5460)--(9206,-5430)--(9227,-5400)--(9248,-5372)--(9266,-5347)
  --(9283,-5324)--(9297,-5306)--(9309,-5291)--(9317,-5280)--(9323,-5273)--(9326,-5269)
  --(9327,-5267);
\pgfsetlinewidth{+7.5\XFigu}
\pgfsetstrokecolor{black}
\draw[pattern,preaction={fill=black}]  (9327,-5202) circle [radius=+65];
\pgfsetdash{{+60\XFigu}{+60\XFigu}}{++0pt}
\draw (10676,-4425)--(10676,-6225);
\draw (13793,-4425)--(13793,-6225);
\pgfsetfillcolor{black}
\pgftext[base,left,at=\pgfqpointxy{300}{-3600}] {\fontsize{28}{33.6}\normalfont $(G, \pi=(({\color{red}\pi_1}, {\color{blue}\pi_2}), \pi_{\mathfrak f}))$}
\pgftext[base,left,at=\pgfqpointxy{9075}{-6975}] {\fontsize{24}{28.8}\normalfont ${\color{red}\mathit{gr}^1_{\pi}(G)}$}
\pgftext[base,left,at=\pgfqpointxy{11925}{-6975}] {\fontsize{24}{28.8}\normalfont ${\color{blue}\mathit{gr}_{\pi}^2(G)}$}
\pgfsetdash{}{+0pt}
\pgfsetfillpattern{xfigp0}{black}
\draw[pattern,preaction={fill=black}]  (2078,-7800) circle [radius=+75];
\endtikzpicture}%

%% file: example7.tikz
{\pgfkeys{/pgf/fpu/.try=false}%
\ifx\XFigwidth\undefined\dimen1=0pt\else\dimen1\XFigwidth\fi
\divide\dimen1 by 15394
\ifx\XFigheight\undefined\dimen3=0pt\else\dimen3\XFigheight\fi
\divide\dimen3 by 7659
\ifdim\dimen1=0pt\ifdim\dimen3=0pt\dimen1=3946sp\dimen3\dimen1
  \else\dimen1\dimen3\fi\else\ifdim\dimen3=0pt\dimen3\dimen1\fi\fi
\tikzpicture[x=+\dimen1, y=+\dimen3]
{\ifx\XFigu\undefined\catcode`\@11
\def\temp{\alloc@1\dimen\dimendef\insc@unt}\temp\XFigu\catcode`\@12\fi}
\XFigu3946sp
\ifdim\XFigu<0pt\XFigu-\XFigu\fi
\catcode`\@11
\pgfutil@ifundefined{pgf@pattern@name@xfigp0}{
\pgfdeclarepatternformonly{xfigp0}
{\pgfqpoint{-1bp}{-1bp}}{\pgfqpoint{9bp}{5bp}}{\pgfqpoint{8bp}{4bp}}
{	\pgfsetdash{}{0pt}\pgfsetlinewidth{0.45bp}
	\pgfpathqmoveto{-1bp}{4.5bp}\pgfpathqlineto{9bp}{-0.5bp}
	\pgfusepathqstroke
}
}{}
\catcode`\@12
\definecolor{green1}{rgb}{0,0.56,0}
\definecolor{pink4}{rgb}{1,0.88,0.88}
\clip(2696,-6045) rectangle (18090,1614);
\tikzset{inner sep=+0pt, outer sep=+0pt}
\pgfsetfillcolor{black}
\pgftext[base,left,at=\pgfqpointxy{16800}{-2250}] {\fontsize{24}{28.8}\normalfont $\mathbb C$}
\pgfsetlinewidth{+7.5\XFigu}
\pgfsetdash{}{+0pt}
\pgfsetstrokecolor{black}
\draw (15104,-2061) arc[start angle=+58.7, end angle=+121.3, radius=+154];
\draw (14944,-2061) arc[start angle=+-44.8, end angle=+-110.8, radius=+112.6];
\draw (14811,-2008) arc[start angle=+-129.5, end angle=+-50.5, radius=+62.9];
\draw (14731,-1981) arc[start angle=+155.2, end angle=+44.8, radius=+140.3];
\draw (14824,-2087) arc[start angle=+281.6, end angle=+160.9, radius=+81.2];
\draw (14838,-2021) arc[start angle=+151, end angle=+66, radius=+30.2];
\draw (14997,-1995) arc[start angle=+-142, end angle=+-38, radius=+33.5];
\draw (15011,-2008) arc[start angle=+159, end angle=+74, radius=+21.4];
\draw (15157,-2008) arc[start angle=+-126.2, end angle=+-37.9, radius=+67.4];
\draw (15184,-2021) arc[start angle=+196, end angle=+20, radius=+21];
\draw (15317,-2061) arc[start angle=+-54.6, end angle=+-125.4, radius=+184];
\draw (15184,-1915) arc[start angle=+121.6, end angle=+-36.9, radius=+100.5];
\draw (16106,-2353) arc[start angle=+40, end angle=+151, radius=+52.3];
\draw (15852,-2186) arc[start angle=+141.3, end angle=+306.4, radius=+122.3];
\draw (16346,-2384) arc[start angle=+-54.9, end angle=+-139.8, radius=+179.4];
\draw (15900,-2313) arc[start angle=+-149, end angle=+-43, radius=+45.3];
\draw (15924,-2329) arc[start angle=+176, end angle=+32, radius=+17.4];
\draw (15900,-2242) arc[start angle=+-125.2, end angle=+-63.3, radius=+92.6];
\draw (15940,-2249) arc[start angle=+130, end angle=+50, radius=+31.2];
\draw (15852,-2186) arc[start angle=+140.5, end angle=+43.5, radius=+154.2];
\draw (16210,-2337) arc[start angle=+-122.0, end angle=+-51.6, radius=+125.2];
\draw (16265,-2353) arc[start angle=+158, end angle=+40, radius=+28.4];
\draw (16249,-2210) arc[start angle=+-139.3, end angle=+-61.3, radius=+71.1];
\draw (16273,-2226) arc[start angle=+169, end angle=+11, radius=+20.4];
\draw (15224,-2291) arc[start angle=+131.0, end angle=+319.9, radius=+92.2];
\draw (15368,-2304) arc[start angle=+127.4, end angle=+61.5, radius=+154.9];
\draw (15355,-2420) arc[start angle=+154.3, end angle=+51.7, radius=+111.1];
\draw (15224,-2291) arc[start angle=+139.8, end angle=+51.5, radius=+282.9];
\draw (15407,-2291) arc[start angle=+-122.5, end angle=+-57.5, radius=+84.7];
\draw (15524,-2381) arc[start angle=+-134.2, end angle=+29.3, radius=+74.6];
\draw (15641,-2291) arc[start angle=+13.4, end angle=+52.0, radius=+70.1];
\draw (14186,-2324) arc[start angle=+184.2, end angle=+66.4, radius=+110.1];
\draw (14340,-2215) arc[start angle=+-104.9, end angle=+-31.8, radius=+150];
\draw (14506,-2149) arc[start angle=+150.6, end angle=+43.0, radius=+109.8];
\draw (14682,-2128) arc[start angle=+26.8, end angle=+-77.7, radius=+131];
\draw (14593,-2315) arc[start angle=+-66.1, end angle=+-95.4, radius=+530.3];
\draw (14186,-2324) arc[start angle=+-174, end angle=+-60, radius=+44.2];
\draw (14252,-2358) arc[start angle=+122.7, end angle=+57.3, radius=+71.3];
\draw (14252,-2292) arc[start angle=+165.1, end angle=+91.0, radius=+57.9];
\draw (14263,-2281) arc[start angle=+-119, end angle=+2, radius=+23];
\draw (14616,-2160) arc[start angle=+34.6, end angle=+-34.6, radius=+88.1];
\draw (14627,-2193) arc[start angle=+127, end angle=+233, radius=+27.5];
\draw  (14891,-2461) ellipse [x radius=+151,y radius=+150];
\draw (14747,-2487) arc[start angle=+-115.4, end angle=+-59.5, radius=+308.8];
\draw (14747,-2487) arc[start angle=+118.9, end angle=+66.2, radius=+325.9];
\draw (16330,-2145) arc[start angle=+61.2, end angle=+-53.5, radius=+142.2];
\draw (16187,-2162) arc[start angle=+138.1, end angle=+55.5, radius=+109];
\draw (16083,-2178) arc[start angle=+-115.9, end angle=+-46.6, radius=+92.4];
\pgfsetdash{{+60\XFigu}{+60\XFigu}}{++0pt}
\pgfsetfillcolor{pink4!5}
\filldraw (12525,675)--(14325,1425)--(17925,1425);
\filldraw (12525,675)--(16275,675)--(17850,1425);
\pgfsetfillcolor{green!5}
\filldraw (12525,-1500)--(14325,-750)--(17925,-750);
\filldraw (12525,-1500)--(16275,-1500)--(17850,-750);
\pgfsetlinewidth{+30\XFigu}
\pgfsetdash{}{+0pt}
\pgfsetstrokecolor{green1}
\draw (14796,-1376) arc[start angle=+192.6, end angle=+139.1, radius=+569.3];
\pgfsetlinewidth{+15\XFigu}
\pgfsetdash{}{+0pt}
\pgfsetstrokecolor{black}
\pgfsetfillpattern{xfigp0}{black}
\draw[pattern,preaction={fill=black}]  (14827,-1406) circle [radius=+44];
\draw[pattern,preaction={fill=black}]  (14237,-1314) circle [radius=+44];
\draw[pattern,preaction={fill=black}]  (15480,-1314) circle [radius=+44];
\draw[pattern,preaction={fill=black}]  (16288,-1314) circle [radius=+44];
\draw[pattern,preaction={fill=black}]  (14983,-879) circle [radius=+44];
\pgfsetlinewidth{+7.5\XFigu}
\pgfsetdash{{+60\XFigu}{+60\XFigu}}{++0pt}
\pgfsetfillcolor{cyan!5}
\filldraw (12540,-375)--(14340,375)--(17940,375);
\filldraw (12540,-375)--(16290,-375)--(17865,375);
\pgfsetlinewidth{+15\XFigu}
\pgfsetdash{}{+0pt}
\pgfsetfillpattern{xfigp0}{black}
\draw[pattern,preaction={fill=black}]  (14225,-250) circle [radius=+36];
\draw[pattern,preaction={fill=black}]  (15324,-250) circle [radius=+36];
\pgfsetlinewidth{+30\XFigu}
\pgfsetstrokecolor{blue}
\pgfsetdash{}{+0pt}
\draw (14225,-250)--(15274,-250);
\draw (15324,-250)--(16224,-250);
\pgfsetlinewidth{+15\XFigu}
\pgfsetdash{}{+0pt}
\pgfsetstrokecolor{black}
\draw[pattern,preaction={fill=black}]  (16224,-250) circle [radius=+36];
\pgfsetlinewidth{+30\XFigu}
\pgfsetdash{}{+0pt}
\pgfsetstrokecolor{blue}
\draw  (15824,81) circle [radius=+164];
\pgfsetlinewidth{+15\XFigu}
\pgfsetdash{}{+0pt}
\pgfsetstrokecolor{black}
\draw[pattern,preaction={fill=black}]  (15661,81) circle [radius=+26];
\draw[pattern,preaction={fill=black}]  (15674,100) circle [radius=+36];
\pgfsetlinewidth{+7.5\XFigu}
\pgfsetdash{{+60\XFigu}{+60\XFigu}}{++0pt}
\draw (12525,-2625)--(14325,-1875)--(17925,-1875);
\draw (12525,-2625)--(16275,-2625)--(17850,-1875);
\pgfsetlinewidth{+30\XFigu}
\pgfsetdash{}{+0pt}
\pgfsetstrokecolor{red}
\draw (14325,1125) arc[start angle=+-116.57, end angle=+-63.43, radius=+1425.5];
\pgfsetlinewidth{+15\XFigu}
\pgfsetdash{}{+0pt}
\pgfsetstrokecolor{black}
\draw[pattern,preaction={fill=black}]  (14325,1125) circle [radius=+53];
\pgfsetlinewidth{+30\XFigu}
\pgfsetdash{}{+0pt}
\pgfsetstrokecolor{red}
\draw (14325,1125) arc[start angle=+-140.70, end angle=+-46.46, radius=+820.4];
\draw (14325,1125) arc[start angle=+118.47, end angle=+61.53, radius=+1258.5];
\draw (14325,1125)--(15600,1125);
\draw  (15787,1125) circle [radius=+188];
\pgfsetlinewidth{+15\XFigu}
\pgfsetdash{}{+0pt}
\pgfsetstrokecolor{black}
\draw[pattern,preaction={fill=black}]  (15600,1125) circle [radius=+53];
\pgfsetlinewidth{+45\XFigu}
\pgfsetdash{}{+0pt}
\draw (13604,-4607) arc[start angle=+40.7, end angle=+149.5, radius=+262.9];
\draw (12367,-3794) arc[start angle=+141.07, end angle=+306.18, radius=+592.8];
\draw (14765,-4761) arc[start angle=+-54.93, end angle=+-140.18, radius=+864.7];
\draw (12598,-4413) arc[start angle=+-149.0, end angle=+-43.8, radius=+221];
\draw (12714,-4490) arc[start angle=+175.3, end angle=+32.3, radius=+84.1];
\draw (12598,-4065) arc[start angle=+-125.7, end angle=+-63.6, radius=+452.6];
\draw (12792,-4103) arc[start angle=+134.3, end angle=+45.7, radius=+139];
\draw (12367,-3794) arc[start angle=+140.23, end angle=+43.75, radius=+751.9];
\draw (14108,-4529) arc[start angle=+-121.89, end angle=+-51.69, radius=+606.2];
\draw (14379,-4607) arc[start angle=+159.2, end angle=+40.0, radius=+135.8];
\draw (14301,-3910) arc[start angle=+-139.7, end angle=+-61.1, radius=+340.9];
\draw (14418,-3988) arc[start angle=+169.5, end angle=+10.5, radius=+97.6];
\draw (13488,-3755) arc[start angle=+-116.0, end angle=+-46.7, radius=+448.3];
\draw (13992,-3678) arc[start angle=+137.4, end angle=+55.2, radius=+532.6];
\draw (14688,-3601) arc[start angle=+61.94, end angle=+-54.35, radius=+684.4];
\pgfsetdash{}{+0pt}
\draw (6353,-2311) arc[start angle=+-108.68, end angle=+-57.87, radius=+1094.8];
\draw (6957,-2805) arc[start angle=+57.71, end angle=+122.29, radius=+616.9];
\draw (6298,-2805) arc[start angle=+-43.8, end angle=+-111.2, radius=+457.1];
\draw (5748,-2586) arc[start angle=+-129.9, end angle=+-50.1, radius=+257];
\draw (5419,-2475) arc[start angle=+155.22, end angle=+44.70, radius=+577];
\draw (5803,-2915) arc[start angle=+281.3, end angle=+160.9, radius=+336.5];
\draw (5859,-2640) arc[start angle=+153.2, end angle=+63.2, radius=+122.1];
\draw (6518,-2530) arc[start angle=+-144.2, end angle=+-35.8, radius=+135.1];
\draw (6572,-2586) arc[start angle=+162.5, end angle=+71.5, radius=+86.5];
\draw (7177,-2586) arc[start angle=+-126.5, end angle=+-36.9, radius=+275.4];
\draw (7286,-2640) arc[start angle=+198.6, end angle=+17.6, radius=+86.8];
\draw (7836,-2805) arc[start angle=+-53.96, end angle=+-126.04, radius=+747];
\draw (7286,-2201) arc[start angle=+121.0, end angle=+-36.3, radius=+416.6];
\draw  (6702,-5394) circle [radius=+621];
\draw (6105,-5506) arc[start angle=+-115.14, end angle=+-59.49, radius=+1281.5];
\pgfsetlinewidth{+15\XFigu}
\pgfsetdash{{+90\XFigu}{+90\XFigu}}{++0pt}
\draw (6105,-5506) arc[start angle=+119.64, end angle=+65.73, radius=+1319.6];
\pgfsetlinewidth{+45\XFigu}
\pgfsetdash{}{+0pt}
\draw (2729,-4524) arc[start angle=+184.9, end angle=+66.2, radius=+454.4];
\draw (3365,-4069) arc[start angle=+-105.04, end angle=+-31.25, radius=+612.9];
\draw (4048,-3795) arc[start angle=+151.8, end angle=+42.3, radius=+449.4];
\draw (4776,-3705) arc[start angle=+26.6, end angle=+-77.2, radius=+543.2];
\draw (4411,-4478) arc[start angle=+-65.59, end angle=+-95.47, radius=+2144.9];
\draw (2729,-4524) arc[start angle=+-172.7, end angle=+-60.3, radius=+183.5];
\draw (3002,-4660) arc[start angle=+122.2, end angle=+57.8, radius=+298.4];
\draw (3002,-4387) arc[start angle=+166.7, end angle=+90.8, radius=+236.6];
\draw (3047,-4341) arc[start angle=+-119.6, end angle=+6.2, radius=+92];
\draw (4502,-3841) arc[start angle=+34.3, end angle=+-34.3, radius=+363.4];
\draw (4547,-3978) arc[start angle=+127.1, end angle=+232.9, radius=+113.5];
\pgfsetstrokecolor{green1}
\draw (6236,-2714) arc[start angle=+157.42, end angle=+211.81, radius=+2527.7];
\pgfsetstrokecolor{red}
\draw (4617,-3772)--(5675,-2776);
\pgfsetstrokecolor{blue}
\draw (6796,-2651) arc[start angle=+23.99, end angle=+-15.18, radius=+3631.9];
\pgfsetstrokecolor{red}
\draw (7792,-2651)--(8663,-4083);
\draw (4306,-4394)--(6236,-5328);
\pgfsetstrokecolor{blue}
\draw (4742,-4021) arc[start angle=+-112.18, end angle=+-86.22, radius=+5192.8];
\pgfsetstrokecolor{red}
\draw (7169,-5266)--(8539,-4394);
\pgfsetdash{}{+0pt}
\pgfsetstrokecolor{black}
\draw (8389,-4007) arc[start angle=+130.3, end angle=+319.7, radius=+381];
\draw (8980,-4061) arc[start angle=+127.47, end angle=+61.37, radius=+641.8];
\draw (8926,-4544) arc[start angle=+154.1, end angle=+51.9, radius=+460.7];
\draw (8389,-4007) arc[start angle=+140.06, end angle=+51.35, radius=+1158.6];
\draw (9141,-4007) arc[start angle=+-122.6, end angle=+-57.4, radius=+348.6];
\draw (9625,-4383) arc[start angle=+-134.1, end angle=+29.9, radius=+309.1];
\draw (10108,-4007) arc[start angle=+14.2, end angle=+53.0, radius=+291];
\pgfsetdash{}{+0pt}
\pgfsetstrokecolor{blue}
\draw (7169,-4394) arc[start angle=+106.05, end angle=+84.78, radius=+3560.8];
\draw (9846,-4332) arc[start angle=+-127.23, end angle=+-55.55, radius=+2179.8];
\pgfsetstrokecolor{red}
\draw (9908,-4021)--(12398,-4021);
\pgfsetfillcolor{black}
\pgftext[base,left,at=\pgfqpointxy{16800}{-975}] {\fontsize{24}{28.8}\normalfont $\infty_3$}
\pgftext[base,left,at=\pgfqpointxy{16800}{150}] {\fontsize{24}{28.8}\normalfont $\infty_2$}
\pgftext[base,left,at=\pgfqpointxy{16800}{1200}] {\fontsize{24}{28.8}\normalfont $\infty_1$}
\pgfsetlinewidth{+7.5\XFigu}
\pgfsetdash{}{+0pt}
\pgfsetstrokecolor{black}
\draw (14958,-1941) arc[start angle=+-108.2, end angle=+-58.7, radius=+271.6];
\endtikzpicture}%

%% file: example5.tikz
{\pgfkeys{/pgf/fpu/.try=false}%
\ifx\XFigwidth\undefined\dimen1=0pt\else\dimen1\XFigwidth\fi
\divide\dimen1 by 7009
\ifx\XFigheight\undefined\dimen3=0pt\else\dimen3\XFigheight\fi
\divide\dimen3 by 6409
\ifdim\dimen1=0pt\ifdim\dimen3=0pt\dimen1=3946sp\dimen3\dimen1
  \else\dimen1\dimen3\fi\else\ifdim\dimen3=0pt\dimen3\dimen1\fi\fi
\tikzpicture[x=+\dimen1, y=+\dimen3]
{\ifx\XFigu\undefined\catcode`\@11
\def\temp{\alloc@1\dimen\dimendef\insc@unt}\temp\XFigu\catcode`\@12\fi}
\XFigu3946sp
\ifdim\XFigu<0pt\XFigu-\XFigu\fi
\definecolor{blue3}{rgb}{0,0,0.82}
\definecolor{green1}{rgb}{0,0.56,0}
\definecolor{cyan3}{rgb}{0,0.82,0.82}
\definecolor{red3}{rgb}{0.82,0,0}
\clip(2656,-8744) rectangle (9665,-2335);
\tikzset{inner sep=+0pt, outer sep=+0pt}
\pgfsetlinewidth{+75\XFigu}
\pgfsetstrokecolor{red}
\pgfsetdash{}{+0pt}
\draw (8700,-2850)--(9300,-3600);
\pgfsetlinewidth{+45\XFigu}
\pgfsetstrokecolor{black}
\pgfsetdash{}{+0pt}
\draw (2700,-8700)--(2700,-4500);
\draw (7200,-8700)--(7200,-5400);
\draw (2700,-4500)--(6300,-4500);
\draw (7200,-8700)--(9600,-6300);
\draw (9600,-6300)--(9600,-5100);
\draw (7200,-5400)--(9000,-3900);
\draw (9600,-4500)--(9600,-6300);
\draw (6300,-4500)--(7200,-5400);
\draw (6300,-4500)--(8100,-3000);
\draw (5100,-2400)--(7800,-2400);
\pgfsetstrokecolor{green1}
\draw (9600,-4500)--(9300,-3600);
\pgfsetstrokecolor{red3}
\draw (9000,-3900)--(9600,-4500);
\draw (7800,-2400)--(8100,-3000);
\pgfsetstrokecolor{black}
\draw (2700,-4500)--(5100,-2400);
\pgfsetstrokecolor{green1}
\draw (7800,-2400)--(8700,-2850);
\pgfsetlinewidth{+7.5\XFigu}
\pgfsetstrokecolor{cyan3}
\pgfsetdash{{+15\XFigu}{+45\XFigu}}{+15\XFigu}
\draw (8100,-3000)--(7800,-2400);
\pgfsetlinewidth{+45\XFigu}
\pgfsetstrokecolor{red3}
\pgfsetdash{}{+0pt}
\draw (7800,-2400)--(8100,-3000);
\draw (9000,-3900)--(9600,-4500);
\draw (8700,-2850)--(9300,-3600);
\pgfsetstrokecolor{green1}
\draw (7800,-2400)--(8700,-2850);
\pgfsetstrokecolor{blue}
\draw (7800,-2400)--(8700,-2850);
\draw (8100,-3000)--(9000,-3900);
\draw (9300,-3600)--(9600,-4500);
\pgfsetstrokecolor{red3}
\pgfsetdash{}{+0pt}
\draw (9300,-3600)--(8700,-2850);
\pgfsetstrokecolor{green1}
\draw (8100,-3000)--(9000,-3900);
\pgfsetdash{}{+0pt}
\pgfsetstrokecolor{cyan}
\pgfsetfillcolor{cyan!65}
\filldraw (8100,-3000)--(7800,-2400)--(8700,-2850)--(9300,-3600)--(9600,-4500)--(9000,-3900)
  --(8100,-3000);
\pgfsetstrokecolor{blue3}
\draw (7800,-2400)--(8700,-2850);
\draw (9300,-3600)--(9600,-4500);
\draw (8100,-3000)--(9000,-3900);
\pgfsetstrokecolor{red3}
\draw (7800,-2400)--(8100,-3000);
\draw (8700,-2850)--(9300,-3600);
\draw (9000,-3900)--(9600,-4500);
\pgfsetstrokecolor{black}
\draw (9600,-4500)--(9600,-6300);
\draw (5100,-2400)--(7800,-2400);
\pgfsetlinewidth{+7.5\XFigu}
\pgfsetdash{{+60\XFigu}{+60\XFigu}}{++0pt}
\draw (6300,-2850)--(8700,-2850);
\draw (9300,-5400)--(9300,-3600);
\draw (5100,-2400)--(6300,-2850);
\draw (6300,-2850)--(6300,-5400);
\draw (6300,-5400)--(9300,-5400);
\draw (9300,-5400)--(9600,-6300);
\draw (2700,-8700)--(6300,-5400);
\pgfsetlinewidth{+75\XFigu}
\pgfsetstrokecolor{blue3}
\pgfsetdash{}{+0pt}
\draw (7800,-2400)--(8700,-2850);
\draw (9300,-3600)--(9600,-4500);
\draw (8100,-3000)--(9000,-3900);
\pgfsetstrokecolor{red}
\draw (7800,-2400)--(8100,-3000);
\draw (9000,-3900)--(9600,-4500);
\pgfsetlinewidth{+45\XFigu}
\pgfsetstrokecolor{black}
\pgfsetdash{}{+0pt}
\draw (2700,-8700)--(7200,-8700);
\endtikzpicture}%

%% file: example4.tikz
{\pgfkeys{/pgf/fpu/.try=false}%
\ifx\XFigwidth\undefined\dimen1=0pt\else\dimen1\XFigwidth\fi
\divide\dimen1 by 16137
\ifx\XFigheight\undefined\dimen3=0pt\else\dimen3\XFigheight\fi
\divide\dimen3 by 8061
\ifdim\dimen1=0pt\ifdim\dimen3=0pt\dimen1=3946sp\dimen3\dimen1
  \else\dimen1\dimen3\fi\else\ifdim\dimen3=0pt\dimen3\dimen1\fi\fi
\tikzpicture[x=+\dimen1, y=+\dimen3]
{\ifx\XFigu\undefined\catcode`\@11
\def\temp{\alloc@1\dimen\dimendef\insc@unt}\temp\XFigu\catcode`\@12\fi}
\XFigu3946sp
\ifdim\XFigu<0pt\XFigu-\XFigu\fi
\catcode`\@11
\pgfutil@ifundefined{pgf@pattern@name@xfigp0}{
\pgfdeclarepatternformonly{xfigp0}
{\pgfqpoint{-1bp}{-1bp}}{\pgfqpoint{9bp}{5bp}}{\pgfqpoint{8bp}{4bp}}
{	\pgfsetdash{}{0pt}\pgfsetlinewidth{0.45bp}
	\pgfpathqmoveto{-1bp}{4.5bp}\pgfpathqlineto{9bp}{-0.5bp}
	\pgfusepathqstroke
}
}{}
\pgfutil@ifundefined{pgf@pattern@name@xfigp4}{
\pgfdeclarepatternformonly{xfigp4}
{\pgfqpoint{-1bp}{-1bp}}{\pgfqpoint{9bp}{9bp}}{\pgfqpoint{8bp}{8bp}}
{	\pgfsetdash{}{0pt}\pgfsetlinewidth{0.45bp}
	\pgfpathqmoveto{-1bp}{-1bp}\pgfpathqlineto{9bp}{9bp}
	\pgfusepathqstroke
}
}{}
\catcode`\@12
\pgfdeclarearrow{
  name = xfiga0,
  parameters = {
    \the\pgfarrowlinewidth \the\pgfarrowlength \the\pgfarrowwidth},
  defaults = {
	  line width=+7.5\XFigu, length=+120\XFigu, width=+60\XFigu},
  setup code = {
    \dimen7 2.15\pgfarrowlength\pgfmathveclen{\the\dimen7}{\the\pgfarrowwidth}
    \dimen7 2\pgfarrowwidth\pgfmathdivide{\pgfmathresult}{\the\dimen7}
    \dimen7 \pgfmathresult\pgfarrowlinewidth
    \pgfarrowssettipend{+\dimen7}
    \pgfarrowssetbackend{+-\pgfarrowlength}
    \dimen9 -0.5\pgfarrowlinewidth
    \pgfarrowssetvisualbackend{+\dimen9}
    \pgfarrowssetlineend{+-0.5\pgfarrowlinewidth}
    \pgfarrowshullpoint{+\dimen7}{+0pt}
    \pgfarrowsupperhullpoint{+-\pgfarrowlength}{+0.5\pgfarrowwidth}
    \pgfarrowssavethe\pgfarrowlinewidth
    \pgfarrowssavethe\pgfarrowlength
    \pgfarrowssavethe\pgfarrowwidth
  },
  drawing code = {\pgfsetdash{}{+0pt}
    \ifdim\pgfarrowlinewidth=\pgflinewidth\else\pgfsetlinewidth{+\pgfarrowlinewidth}\fi
    \pgfpathmoveto{\pgfqpoint{-\pgfarrowlength}{0.5\pgfarrowwidth}}
    \pgfpathlineto{\pgfqpoint{0pt}{0pt}}
    \pgfpathlineto{\pgfqpoint{-\pgfarrowlength}{-0.5\pgfarrowwidth}}
    \pgfusepathqstroke
  }
}
\definecolor{blue3}{rgb}{0,0,0.82}
\definecolor{red3}{rgb}{0.82,0,0}
\definecolor{pink4}{rgb}{1,0.88,0.88}
\definecolor{xfigc32}{rgb}{0.000,0.000,0.000}
\clip(-315,-10486) rectangle (15822,-2425);
\tikzset{inner sep=+0pt, outer sep=+0pt}
\pgfsetfillcolor{black}
\pgftext[base,left,at=\pgfqpointxy{14400}{-3600}] {\fontsize{28}{33.6}\normalfont $\tilde{\mathscr{M}}_{G}^{^{\mathit{trop}}}$}
\pgfsetroundcap
\pgfsetlinewidth{+7.5\XFigu}
\pgfsetdash{{+60\XFigu}{+60\XFigu}}{++0pt}
\pgfsetstrokecolor{black}
\pgfsetfillpattern{xfigp0}{black}
\draw[pattern,preaction={fill=black}]  (150,-7950) circle [radius=+150];
\draw[pattern,preaction={fill=black}]  (2550,-7950) circle [radius=+150];
\pgfsetlinewidth{+60\XFigu}
\pgfsetbuttcap
\pgfsetdash{}{+0pt}
\draw (150,-7950)--(2550,-7950);
\pgfsetarrows{[line width=7.5\XFigu, width=90\XFigu, length=180\XFigu]}
\pgfsetarrowsend{xfiga0}
\pgfsetlinewidth{+30\XFigu}
\pgfsetdash{{+120\XFigu}{+120\XFigu}}{++0pt}
\draw (9900,-3600) arc[start angle=+119.74, end angle=+240.26, radius=+604.7];
\draw (12450,-3000) arc[start angle=+126.09, end angle=+180.78, radius=+1825.2];
\draw (14250,-4800) arc[start angle=+-27.90, end angle=+-98.97, radius=+1442.7];
\draw (14700,-8100) arc[start angle=+15.26, end angle=+127.87, radius=+855.1];
\draw (14850,-6300) arc[start angle=+45.78, end angle=+121.53, radius=+1106.2];
\pgfsetlinewidth{+90\XFigu}
\pgfsetdash{}{+0pt}
\draw[pattern,preaction={fill=black}]  (10327,-10131) circle [radius=+63];
\draw[pattern,preaction={fill=black}]  (11088,-10131) circle [radius=+63];
\draw[pattern,preaction={fill=black}]  (14500,-8250) circle [radius=+50];
\draw[pattern,preaction={fill=black}]  (15099,-8250) circle [radius=+50];
\draw[pattern,preaction={fill=black}]  (14515,-6496) circle [radius=+58];
\draw[pattern,preaction={fill=black}]  (15208,-6496) circle [radius=+58];
\draw[pattern,preaction={fill=black}]  (13627,-4581) circle [radius=+63];
\draw[pattern,preaction={fill=black}]  (14388,-4581) circle [radius=+63];
\pgfsetlinewidth{+60\XFigu}
\pgfsetdash{}{+0pt}
\pgfsetstrokecolor{red3}
\draw  (14706,-4581) circle [radius=+317];
\pgfsetlinewidth{+90\XFigu}
\pgfsetdash{}{+0pt}
\pgfsetstrokecolor{black}
\draw[pattern,preaction={fill=black}]  (11975,-2775) ellipse [x radius=+62,y radius=+63];
\pgfsetlinewidth{+60\XFigu}
\pgfsetdash{}{+0pt}
\pgfsetstrokecolor{red3}
\draw  (13037,-2775) ellipse [x radius=+312,y radius=+313];
\pgfsetlinewidth{+90\XFigu}
\pgfsetdash{}{+0pt}
\pgfsetstrokecolor{black}
\draw[pattern,preaction={fill=black}]  (9601,-3450) circle [radius=+50];
\draw[pattern,preaction={fill=black}]  (10200,-3450) circle [radius=+50];
\pgfsetlinewidth{+60\XFigu}
\pgfsetdash{}{+0pt}
\pgfsetstrokecolor{red3}
\draw  (10450,-3450) circle [radius=+250];
\pgfsetstrokecolor{blue3}
\draw  (15496,-6496) circle [radius=+288];
\pgfsetstrokecolor{black}
\draw  (15349,-8250) circle [radius=+250];
\draw  (11406,-10131) circle [radius=+317];
\pgfsetlinewidth{+90\XFigu}
\pgfsetdash{}{+0pt}
\draw[pattern,preaction={fill=black}]  (12724,-2775) ellipse [x radius=+62,y radius=+63];
\pgfsetlinewidth{+75\XFigu}
\pgfsetdash{}{+0pt}
\pgfsetfillpattern{xfigp4}{black}
\draw[pattern,preaction={fill=white}]  (13200,-6300) circle [radius=+75];
\draw[pattern,preaction={fill=white}]  (11700,-4800) circle [radius=+75];
\pgfsetlinewidth{+60\XFigu}
\pgfsetdash{}{+0pt}
\pgfsetarrowsend{}
\draw (10327,-10131)--(11088,-10131);
\draw (9601,-3450)--(10200,-3450);
\pgfsetstrokecolor{blue3}
\draw (11975,-2775)--(12724,-2775);
\pgfsetstrokecolor{red3}
\draw (13627,-4581)--(14388,-4581);
\draw (14515,-6496)--(15208,-6496);
\draw (14500,-8250)--(15099,-8250);
\pgfsetlinewidth{+75\XFigu}
\pgfsetstrokecolor{black}
\pgfsetdash{}{+0pt}
\draw (9000,-4800)--(11700,-4800);
\draw (11700,-4800)--(13200,-6300);
\pgfsetarrowsend{xfiga0}
\pgfsetlinewidth{+30\XFigu}
\pgfsetdash{{+120\XFigu}{+120\XFigu}}{++0pt}
\draw (10650,-10050)--(9900,-9300)--(9900,-8550);
\pgfsetlinewidth{+15\XFigu}
\pgfsetdash{{+90\XFigu}{+90\XFigu}}{++0pt}
\pgfsetarrowsend{}
\draw (9000,-9000)--(9000,-4800);
\draw (9000,-9000)--(13200,-9000);
\pgfsetlinewidth{+75\XFigu}
\pgfsetdash{}{+0pt}
\draw (13200,-6300)--(13200,-9000);
\pgfsetlinewidth{+30\XFigu}
\pgfsetdash{{+120\XFigu}{+120\XFigu}}{++0pt}
\pgfsetfillcolor{pink4!10}
\filldraw (9000,-4800)--(11700,-4800)--(13200,-6300)--(13200,-9000)--(9000,-9000)--(9000,-4800);
\pgfsetlinewidth{+45\XFigu}
\pgfsetdash{}{+0pt}
\draw (9000,-4800)--(11700,-4800);
\draw (13200,-6300)--(13200,-9000);
\draw (11700,-4800)--(13200,-6300);
\pgfsetfillcolor{black}
\pgftext[base,left,at=\pgfqpointxy{-300}{-7350}] {\fontsize{28}{33.6}\normalfont $G=(V,E)$}
\pgftext[base,left,at=\pgfqpointxy{4950}{-2850}] {\fontsize{28}{33.6}\normalfont ${\color{red}\pi_1}, {\color{blue}\pi_2}, \pi_{\mathfrak f}$}
\pgftext[base,left,at=\pgfqpointxy{11775}{-10200}] {\fontsize{28}{33.6}\normalfont $\tilde\mg^{^{\mathit{gr}}}_{G}$}
\pgfsetlinewidth{+60\XFigu}
\pgfsetdash{}{+0pt}
\draw  (3225,-7875) circle [radius=+679];
\endtikzpicture}%

%% file: example6.tikz
{\pgfkeys{/pgf/fpu/.try=false}%
\ifx\XFigwidth\undefined\dimen1=0pt\else\dimen1\XFigwidth\fi
\divide\dimen1 by 11055
\ifx\XFigheight\undefined\dimen3=0pt\else\dimen3\XFigheight\fi
\divide\dimen3 by 6904
\ifdim\dimen1=0pt\ifdim\dimen3=0pt\dimen1=3946sp\dimen3\dimen1
  \else\dimen1\dimen3\fi\else\ifdim\dimen3=0pt\dimen3\dimen1\fi\fi
\tikzpicture[x=+\dimen1, y=+\dimen3]
{\ifx\XFigu\undefined\catcode`\@11
\def\temp{\alloc@1\dimen\dimendef\insc@unt}\temp\XFigu\catcode`\@12\fi}
\XFigu3946sp
\ifdim\XFigu<0pt\XFigu-\XFigu\fi
\catcode`\@11
\pgfutil@ifundefined{pgf@pattern@name@xfigp0}{
\pgfdeclarepatternformonly{xfigp0}
{\pgfqpoint{-1bp}{-1bp}}{\pgfqpoint{9bp}{5bp}}{\pgfqpoint{8bp}{4bp}}
{	\pgfsetdash{}{0pt}\pgfsetlinewidth{0.45bp}
	\pgfpathqmoveto{-1bp}{4.5bp}\pgfpathqlineto{9bp}{-0.5bp}
	\pgfusepathqstroke
}
}{}
\catcode`\@12
\pgfdeclarearrow{
  name = xfiga0,
  parameters = {
    \the\pgfarrowlinewidth \the\pgfarrowlength \the\pgfarrowwidth},
  defaults = {
	  line width=+7.5\XFigu, length=+120\XFigu, width=+60\XFigu},
  setup code = {
    \dimen7 2.15\pgfarrowlength\pgfmathveclen{\the\dimen7}{\the\pgfarrowwidth}
    \dimen7 2\pgfarrowwidth\pgfmathdivide{\pgfmathresult}{\the\dimen7}
    \dimen7 \pgfmathresult\pgfarrowlinewidth
    \pgfarrowssettipend{+\dimen7}
    \pgfarrowssetbackend{+-\pgfarrowlength}
    \dimen9 -0.5\pgfarrowlinewidth
    \pgfarrowssetvisualbackend{+\dimen9}
    \pgfarrowssetlineend{+-0.5\pgfarrowlinewidth}
    \pgfarrowshullpoint{+\dimen7}{+0pt}
    \pgfarrowsupperhullpoint{+-\pgfarrowlength}{+0.5\pgfarrowwidth}
    \pgfarrowssavethe\pgfarrowlinewidth
    \pgfarrowssavethe\pgfarrowlength
    \pgfarrowssavethe\pgfarrowwidth
  },
  drawing code = {\pgfsetdash{}{+0pt}
    \ifdim\pgfarrowlinewidth=\pgflinewidth\else\pgfsetlinewidth{+\pgfarrowlinewidth}\fi
    \pgfpathmoveto{\pgfqpoint{-\pgfarrowlength}{0.5\pgfarrowwidth}}
    \pgfpathlineto{\pgfqpoint{0pt}{0pt}}
    \pgfpathlineto{\pgfqpoint{-\pgfarrowlength}{-0.5\pgfarrowwidth}}
    \pgfusepathqstroke
  }
}
\definecolor{blue3}{rgb}{0,0,0.82}
\definecolor{green1}{rgb}{0,0.56,0}
\definecolor{red1}{rgb}{0.56,0,0}
\definecolor{magenta1}{rgb}{0.56,0,0.56}
\definecolor{xfigc32}{rgb}{0.000,0.000,0.000}
\clip(1313,-9720) rectangle (12368,-2816);
\tikzset{inner sep=+0pt, outer sep=+0pt}
\pgfsetlinewidth{+60\XFigu}
\pgfsetstrokecolor{green1}
\pgfsetdash{}{+0pt}
\draw (9690,-6000)--(11370,-6000);
\pgfsetstrokecolor{blue3}
\draw  (9028,-3277) circle [radius=+352];
\pgfsetarrows{[line width=7.5\XFigu, width=90\XFigu, length=180\XFigu]}
\pgfsetarrowsend{xfiga0}
\pgfsetlinewidth{+30\XFigu}
\pgfsetdash{{+120\XFigu}{+120\XFigu}}{++0pt}
\pgfsetstrokecolor{black}
\draw (8400,-3450) arc[start angle=+6.90, end angle=+-92.32, radius=+1742.1];
\draw (3600,-3600) arc[start angle=+-172.06, end angle=+-91.21, radius=+1392.7];
\draw (9900,-6150) arc[start angle=+-54.69, end angle=+-125.31, radius=+2205.8];
\pgfsetlinewidth{+75\XFigu}
\pgfsetdash{}{+0pt}
\pgfsetstrokecolor{red1}
\pgfsetfillpattern{xfigp0}{red1}
\draw[pattern,preaction={fill=black}]  (7274,-6000) circle [radius=+75];
\pgfsetstrokecolor{black}
\pgfsetfillpattern{xfigp0}{black}
\draw[pattern,preaction={fill=black}]  (3117,-6000) circle [radius=+75];
\draw[pattern,preaction={fill=black}]  (3117,-8400) circle [radius=+75];
\draw[pattern,preaction={fill=black}]  (5196,-9600) circle [radius=+75];
\draw[pattern,preaction={fill=black}]  (7274,-8400) circle [radius=+75];
\pgfsetstrokecolor{red1}
\pgfsetfillpattern{xfigp0}{red1}
\draw[pattern,preaction={fill=black}]  (5196,-4800) circle [radius=+75];
\pgfsetlinewidth{+60\XFigu}
\pgfsetdash{}{+0pt}
\pgfsetstrokecolor{green1}
\draw  (4461,-3459) circle [radius=+480];
\pgfsetlinewidth{+75\XFigu}
\pgfsetdash{}{+0pt}
\pgfsetstrokecolor{black}
\pgfsetfillpattern{xfigp0}{black}
\draw[pattern,preaction={fill=black}]  (2325,-3450) circle [radius=+60];
\draw[pattern,preaction={fill=black}]  (3975,-3450) circle [radius=+60];
\pgfsetlinewidth{+45\XFigu}
\pgfsetdash{}{+0pt}
\pgfsetstrokecolor{magenta1}
\draw  (1823,-3455) circle [radius=+480];
\draw  (9177,-5972) circle [radius=+480];
\draw  (7331,-3325) circle [radius=+480];
\pgfsetlinewidth{+75\XFigu}
\pgfsetdash{}{+0pt}
\pgfsetstrokecolor{black}
\draw[pattern,preaction={fill=black}]  (7800,-3300) circle [radius=+60];
\draw[pattern,preaction={fill=black}]  (9675,-6000) circle [radius=+60];
\pgfsetlinewidth{+60\XFigu}
\pgfsetdash{}{+0pt}
\pgfsetstrokecolor{blue3}
\draw  (11850,-6000) circle [radius=+480];
\pgfsetlinewidth{+75\XFigu}
\pgfsetdash{}{+0pt}
\pgfsetstrokecolor{black}
\draw[pattern,preaction={fill=black}]  (11348,-5970) circle [radius=+60];
\pgfsetstrokecolor{red}
\pgfsetarrowsend{}
\draw (5213,-4800)--(5288,-4800);
\pgfsetstrokecolor{black}
\pgfsetfillcolor{cyan!60}
\filldraw (3117,-6000)--(5196,-4800)--(7274,-6000)--(7274,-8400)--(5196,-9600)--(3117,-8400)
  --(3117,-6000);
\pgfsetstrokecolor{red}
\pgfsetarrowsend{xfiga0}
\draw (5196,-4800)--(7274,-6000);
\draw (5196,-4800)--(7274,-6000);
\draw (5250,-4800)--(7350,-6000);
\pgfsetlinewidth{+45\XFigu}
\pgfsetstrokecolor{blue3}
\pgfsetdash{}{+0pt}
\pgfsetarrowsend{}
\draw (7800,-3300)--(8700,-3300);
\pgfsetlinewidth{+60\XFigu}
\pgfsetdash{}{+0pt}
\draw (2340,-3450)--(4020,-3450);
\pgfsetlinewidth{+75\XFigu}
\pgfsetdash{}{+0pt}
\pgfsetstrokecolor{black}
\pgfsetfillpattern{xfigp0}{black}
\draw[pattern,preaction={fill=black}]  (8580,-3300) circle [radius=+60];
\endtikzpicture}%

%% file: example9.tikz
{\pgfkeys{/pgf/fpu/.try=false}%
\ifx\XFigwidth\undefined\dimen1=0pt\else\dimen1\XFigwidth\fi
\divide\dimen1 by 10468
\ifx\XFigheight\undefined\dimen3=0pt\else\dimen3\XFigheight\fi
\divide\dimen3 by 6244
\ifdim\dimen1=0pt\ifdim\dimen3=0pt\dimen1=3946sp\dimen3\dimen1
  \else\dimen1\dimen3\fi\else\ifdim\dimen3=0pt\dimen3\dimen1\fi\fi
\tikzpicture[x=+\dimen1, y=+\dimen3]
{\ifx\XFigu\undefined\catcode`\@11
\def\temp{\alloc@1\dimen\dimendef\insc@unt}\temp\XFigu\catcode`\@12\fi}
\XFigu3946sp
\ifdim\XFigu<0pt\XFigu-\XFigu\fi
\catcode`\@11
\pgfutil@ifundefined{pgf@pattern@name@xfigp0}{
\pgfdeclarepatternformonly{xfigp0}
{\pgfqpoint{-1bp}{-1bp}}{\pgfqpoint{9bp}{5bp}}{\pgfqpoint{8bp}{4bp}}
{	\pgfsetdash{}{0pt}\pgfsetlinewidth{0.45bp}
	\pgfpathqmoveto{-1bp}{4.5bp}\pgfpathqlineto{9bp}{-0.5bp}
	\pgfusepathqstroke
}
}{}
\catcode`\@12
\definecolor{green1}{rgb}{0,0.56,0}
\clip(855,-9556) rectangle (11323,-3312);
\tikzset{inner sep=+0pt, outer sep=+0pt}
\pgfsetfillcolor{black}
\pgftext[base,left,at=\pgfqpointxy{10875}{-5400}] {\fontsize{24}{28.8}\normalfont $z_3$}
\pgfsetlinewidth{+30\XFigu}
\pgfsetstrokecolor{red}
\draw (6630,-4490) arc[start angle=+147.53, end angle=+212.47, radius=+1695];
\draw (6630,-4490) arc[start angle=+120.51, end angle=+239.49, radius=+1056.2];
\pgfsetlinewidth{+7.5\XFigu}
\pgfsetstrokecolor{black}
\pgfsetfillpattern{xfigp0}{black}
\draw[pattern,preaction={fill=black}]  (6630,-6310) circle [radius=+65];
\draw[pattern,preaction={fill=black}]  (6630,-4490) circle [radius=+65];
\draw[pattern,preaction={fill=black}]  (8450,-4490) circle [radius=+65];
\draw[pattern,preaction={fill=black}]  (8450,-6310) circle [radius=+65];
\pgfsetlinewidth{+30\XFigu}
\pgfsetstrokecolor{blue}
\draw (8450,-4490) arc[start angle=+126.87, end angle=+233.13, radius=+650];
\pgfsetlinewidth{+7.5\XFigu}
\pgfsetstrokecolor{black}
\draw[pattern,preaction={fill=black}]  (8450,-5530) circle [radius=+65];
\pgfsetlinewidth{+30\XFigu}
\pgfsetstrokecolor{blue}
\draw (8450,-4490) arc[start angle=+53.13, end angle=+-53.13, radius=+650];
\pgfsetlinewidth{+45\XFigu}
\pgfsetstrokecolor{red}
\draw (2925,-4200) arc[start angle=+90.00, end angle=+270.00, radius=+1200];
\pgfsetstrokecolor{black}
\draw[pattern,preaction={fill=black}]  (3000,-4200) circle [radius=+75];
\draw[pattern,preaction={fill=black}]  (2100,-5400) circle [radius=+75];
\draw[pattern,preaction={fill=black}]  (3900,-5400) circle [radius=+75];
\draw[pattern,preaction={fill=black}]  (3000,-6600) circle [radius=+75];
\pgfsetlinewidth{+7.5\XFigu}
\draw[pattern,preaction={fill=black}]  (10700,-4490) circle [radius=+65];
\draw[pattern,preaction={fill=black}]  (10700,-6310) circle [radius=+65];
\draw[pattern,preaction={fill=black}]  (11250,-5550) circle [radius=+65];
\draw[pattern,preaction={fill=black}]  (10200,-5550) circle [radius=+65];
\pgfsetlinewidth{+45\XFigu}
\pgfsetstrokecolor{red}
\draw (2175,-5475)--(2925,-6525);
\draw (3075,-6525)--(3825,-5475);
\pgfsetstrokecolor{blue}
\draw (2175,-5325)--(2925,-4275);
\draw (3075,-4275)--(3825,-5325);
\pgfsetlinewidth{+7.5\XFigu}
\pgfsetstrokecolor{black}
\pgfsetdash{{+60\XFigu}{+60\XFigu}}{++0pt}
\draw (5100,-4200)--(5100,-6900);
\draw (7650,-4875)--(7650,-6000);
\pgfsetlinewidth{+45\XFigu}
\pgfsetstrokecolor{green1}
\pgfsetdash{}{+0pt}
\draw (2175,-5400)--(3825,-5400);
\pgfsetlinewidth{+30\XFigu}
\pgfsetdash{}{+0pt}
\draw (10275,-5550)--(11250,-5550);
\pgfsetlinewidth{+7.5\XFigu}
\pgfsetstrokecolor{black}
\pgfsetdash{{+60\XFigu}{+60\XFigu}}{++0pt}
\draw (9525,-4875)--(9525,-6000);
\draw (1575,-7725)--(4275,-7725);
\pgfsetfillcolor{black}
\pgftext[base,left,at=\pgfqpointxy{8700}{-4500}] {\fontsize{24}{28.8}\normalfont $u_2$}
\pgftext[base,left,at=\pgfqpointxy{6900}{-4575}] {\fontsize{24}{28.8}\normalfont $u_1$}
\pgftext[base,left,at=\pgfqpointxy{6825}{-6375}] {\fontsize{24}{28.8}\normalfont $v_1$}
\pgftext[base,left,at=\pgfqpointxy{2925}{-4575}] {\fontsize{24}{28.8}\normalfont $u$}
\pgftext[base,left,at=\pgfqpointxy{3225}{-6750}] {\fontsize{24}{28.8}\normalfont $v$}
\pgftext[base,left,at=\pgfqpointxy{4125}{-5550}] {\fontsize{24}{28.8}\normalfont $z$}
\pgftext[base,left,at=\pgfqpointxy{2250}{-5625}] {\fontsize{24}{28.8}\normalfont $w$}
\pgftext[base,left,at=\pgfqpointxy{6150}{-6975}] {\fontsize{24}{28.8}\normalfont $\Gamma^1$}
\pgftext[base,left,at=\pgfqpointxy{10875}{-4575}] {\fontsize{24}{28.8}\normalfont $u_3$}
\pgftext[base,left,at=\pgfqpointxy{8625}{-5625}] {\fontsize{24}{28.8}\normalfont $w_2$}
\pgftext[base,left,at=\pgfqpointxy{8625}{-6375}] {\fontsize{24}{28.8}\normalfont $v_2$}
\pgftext[base,left,at=\pgfqpointxy{10875}{-6375}] {\fontsize{24}{28.8}\normalfont $v_3$}
\pgftext[base,left,at=\pgfqpointxy{1500}{-8325}] {\fontsize{24}{28.8}\normalfont $\lf=(f_1, f_2, f_3)$}
\pgftext[base,left,at=\pgfqpointxy{1500}{-8850}] {\fontsize{24}{28.8}\normalfont $f_2(u_2) = a, f_2(v_2)=b, f_2(w_2) = c$}
\pgftext[base,left,at=\pgfqpointxy{1485}{-9420}] {\fontsize{24}{28.8}\normalfont $f_3(u_3) = a', f_3(v_3)=b', f_3(w_3) = c', f_3(z_3)=d'$}
\pgftext[base,left,at=\pgfqpointxy{870}{-3615}] {\fontsize{24}{28.8}\normalfont $\curve$}
\pgftext[base,left,at=\pgfqpointxy{870}{-4110}] {\fontsize{24}{28.8}\normalfont $({\color{red}\pi_1},{\color{blue}\pi_2},{\color{teal}\pi_3} )$}
\pgftext[base,left,at=\pgfqpointxy{8325}{-6975}] {\fontsize{24}{28.8}\normalfont $\Gamma^2$}
\pgftext[base,left,at=\pgfqpointxy{10560}{-6975}] {\fontsize{24}{28.8}\normalfont $\Gamma^3$}
\pgftext[base,left,at=\pgfqpointxy{9600}{-5625}] {\fontsize{24}{28.8}\normalfont $w_3$}
\pgfsetlinewidth{+30\XFigu}
\pgfsetdash{}{+0pt}
\pgfsetstrokecolor{red}
\draw (6630,-4490) arc[start angle=+31.89, end angle=+-31.89, radius=+1722.5];
\endtikzpicture}%

%% file: example10.tikz
{\pgfkeys{/pgf/fpu/.try=false}%
\ifx\XFigwidth\undefined\dimen1=0pt\else\dimen1\XFigwidth\fi
\divide\dimen1 by 12895
\ifx\XFigheight\undefined\dimen3=0pt\else\dimen3\XFigheight\fi
\divide\dimen3 by 4489
\ifdim\dimen1=0pt\ifdim\dimen3=0pt\dimen1=3946sp\dimen3\dimen1
  \else\dimen1\dimen3\fi\else\ifdim\dimen3=0pt\dimen3\dimen1\fi\fi
\tikzpicture[x=+\dimen1, y=+\dimen3]
{\ifx\XFigu\undefined\catcode`\@11
\def\temp{\alloc@1\dimen\dimendef\insc@unt}\temp\XFigu\catcode`\@12\fi}
\XFigu3946sp
\ifdim\XFigu<0pt\XFigu-\XFigu\fi
\definecolor{blue3}{rgb}{0,0,0.82}
\definecolor{red3}{rgb}{0.82,0,0}
\clip(5385,-10186) rectangle (18280,-5697);
\tikzset{inner sep=+0pt, outer sep=+0pt}
\pgfsetfillcolor{black}
\pgftext[base,left,at=\pgfqpointxy{17850}{-7275}] {\fontsize{24}{28.8}\normalfont $y_2$}
\pgfsetlinewidth{+30\XFigu}
\pgfsetstrokecolor{red3}
\draw (13381,-6870) arc[start angle=+97.6, end angle=+262.4, radius=+484.3];
\draw (13500,-7830) arc[start angle=+-53.13, end angle=+53.13, radius=+600];
\draw (13500,-7830) arc[start angle=+-82.4, end angle=+82.4, radius=+484.3];
\pgfsetlinewidth{+15\XFigu}
\pgfsetdash{}{+0pt}
\pgfsetstrokecolor{black}
\pgfsetfillcolor{cyan}
\filldraw  (13440,-7830) circle [radius=+55];
\filldraw  (13440,-6870) circle [radius=+55];
\pgfsetlinewidth{+30\XFigu}
\pgfsetstrokecolor{red3}
\pgfsetdash{}{+0pt}
\draw (13500,-6870)--(14340,-6870);
\pgfsetlinewidth{+15\XFigu}
\pgfsetdash{}{+0pt}
\pgfsetstrokecolor{black}
\filldraw  (7871,-6887) circle [radius=+68];
\pgfsetlinewidth{+30\XFigu}
\pgfsetdash{}{+0pt}
\pgfsetstrokecolor{red3}
\draw (7186,-7846)--(7871,-6887)--(9515,-6887);
\pgfsetlinewidth{+15\XFigu}
\pgfsetdash{}{+0pt}
\pgfsetstrokecolor{black}
\filldraw  (10816,-6887) circle [radius=+68];
\filldraw  (9515,-6887) circle [radius=+68];
\pgfsetlinewidth{+30\XFigu}
\pgfsetstrokecolor{red3}
\pgfsetdash{}{+0pt}
\draw (9600,-6900)--(10800,-6900);
\pgfsetlinewidth{+15\XFigu}
\pgfsetdash{}{+0pt}
\pgfsetstrokecolor{black}
\filldraw  (8419,-8189) circle [radius=+68];
\pgfsetlinewidth{+30\XFigu}
\pgfsetstrokecolor{red3}
\pgfsetdash{}{+0pt}
\draw (7186,-7846)--(8419,-8119);
\pgfsetlinewidth{+15\XFigu}
\pgfsetdash{}{+0pt}
\pgfsetstrokecolor{black}
\filldraw  (7186,-7846) circle [radius=+68];
\pgfsetlinewidth{+30\XFigu}
\pgfsetstrokecolor{blue3}
\pgfsetdash{}{+0pt}
\draw (7186,-7846)--(9515,-6887);
\draw (7871,-6887)--(8419,-8119);
\pgfsetstrokecolor{red3}
\pgfsetdash{}{+0pt}
\draw (8419,-8119)--(9515,-6887);
\pgfsetlinewidth{+15\XFigu}
\pgfsetdash{}{+0pt}
\pgfsetstrokecolor{black}
\filldraw  (14400,-6900) circle [radius=+68];
\filldraw  (15622,-7608) circle [radius=+49];
\filldraw  (16503,-7853) circle [radius=+49];
\filldraw  (16111,-6923) circle [radius=+49];
\filldraw  (18216,-6923) circle [radius=+49];
\filldraw  (17287,-6923) circle [radius=+49];
\pgfsetlinewidth{+30\XFigu}
\pgfsetstrokecolor{blue3}
\pgfsetdash{}{+0pt}
\draw (16111,-6923)--(16503,-7803);
\draw (15622,-7608)--(17287,-6923);
\pgfsetlinewidth{+7.5\XFigu}
\pgfsetstrokecolor{black}
\pgfsetdash{{+60\XFigu}{+60\XFigu}}{++0pt}
\draw (11700,-6600)--(11700,-8700);
\draw (15000,-6900)--(15000,-7875);
\pgfsetfillcolor{black}
\pgftext[base,left,at=\pgfqpointxy{5700}{-6000}] {\fontsize{24}{28.8}\normalfont $\curve$}
\pgftext[base,left,at=\pgfqpointxy{5700}{-6450}] {\fontsize{24}{28.8}\normalfont $\pi=({\color{red}\pi_1}, {\color{blue}\pi_2})$}
\pgftext[base,left,at=\pgfqpointxy{8550}{-8325}] {\fontsize{24}{28.8}\normalfont $w$}
\pgftext[base,left,at=\pgfqpointxy{9525}{-7200}] {\fontsize{24}{28.8}\normalfont $x$}
\pgftext[base,left,at=\pgfqpointxy{10725}{-7200}] {\fontsize{24}{28.8}\normalfont $y$}
\pgftext[base,left,at=\pgfqpointxy{13275}{-8775}] {\fontsize{24}{28.8}\normalfont $\Gamma^1$}
\pgftext[base,left,at=\pgfqpointxy{16425}{-8775}] {\fontsize{24}{28.8}\normalfont $\Gamma^2$}
\pgftext[base,left,at=\pgfqpointxy{17175}{-7275}] {\fontsize{24}{28.8}\normalfont $x_2$}
\pgftext[base,left,at=\pgfqpointxy{7500}{-6975}] {\fontsize{24}{28.8}\normalfont $u$}
\pgftext[base,left,at=\pgfqpointxy{6900}{-8100}] {\fontsize{24}{28.8}\normalfont $v$}
\pgftext[base,left,at=\pgfqpointxy{14325}{-6675}] {\fontsize{24}{28.8}\normalfont $y_1$}
\pgftext[base,left,at=\pgfqpointxy{16050}{-6750}] {\fontsize{24}{28.8}\normalfont $u_2$}
\pgftext[base,left,at=\pgfqpointxy{5400}{-9525}] {\fontsize{24}{28.8}\normalfont $\lf=(f_1, f_2)$}
\pgftext[base,left,at=\pgfqpointxy{5400}{-10050}] {\fontsize{24}{28.8}\normalfont $f_1(u_1) = a, f_1(v_1)=b, f_1(y_1)=c$}
\pgftext[base,left,at=\pgfqpointxy{13275}{-6675}] {\fontsize{24}{28.8}\normalfont $v_1$}
\pgftext[base,left,at=\pgfqpointxy{13275}{-8175}] {\fontsize{24}{28.8}\normalfont $u_1$}
\pgftext[base,left,at=\pgfqpointxy{15375}{-7875}] {\fontsize{24}{28.8}\normalfont $v_2$}
\pgftext[base,left,at=\pgfqpointxy{16650}{-7950}] {\fontsize{24}{28.8}\normalfont $w_2$}
\pgfsetlinewidth{+30\XFigu}
\pgfsetdash{}{+0pt}
\pgfsetstrokecolor{red3}
\draw (13440,-6870) arc[start angle=+127.06, end angle=+232.94, radius=+601.5];
\endtikzpicture}%

%% file: logmap.tikz
{\pgfkeys{/pgf/fpu/.try=false}%
\ifx\XFigwidth\undefined\dimen1=0pt\else\dimen1\XFigwidth\fi
\divide\dimen1 by 12222
\ifx\XFigheight\undefined\dimen3=0pt\else\dimen3\XFigheight\fi
\divide\dimen3 by 12285
\ifdim\dimen1=0pt\ifdim\dimen3=0pt\dimen1=3946sp\dimen3\dimen1
  \else\dimen1\dimen3\fi\else\ifdim\dimen3=0pt\dimen3\dimen1\fi\fi
\tikzpicture[x=+\dimen1, y=+\dimen3]
{\ifx\XFigu\undefined\catcode`\@11
\def\temp{\alloc@1\dimen\dimendef\insc@unt}\temp\XFigu\catcode`\@12\fi}
\XFigu3946sp
\ifdim\XFigu<0pt\XFigu-\XFigu\fi
\catcode`\@11
\pgfutil@ifundefined{pgf@pattern@name@xfigp0}{
\pgfdeclarepatternformonly{xfigp0}
{\pgfqpoint{-1bp}{-1bp}}{\pgfqpoint{9bp}{5bp}}{\pgfqpoint{8bp}{4bp}}
{	\pgfsetdash{}{0pt}\pgfsetlinewidth{0.45bp}
	\pgfpathqmoveto{-1bp}{4.5bp}\pgfpathqlineto{9bp}{-0.5bp}
	\pgfusepathqstroke
}
}{}
\catcode`\@12
\pgfdeclarearrow{
  name = xfiga0,
  parameters = {
    \the\pgfarrowlinewidth \the\pgfarrowlength \the\pgfarrowwidth},
  defaults = {
	  line width=+7.5\XFigu, length=+120\XFigu, width=+60\XFigu},
  setup code = {
    \dimen7 2.15\pgfarrowlength\pgfmathveclen{\the\dimen7}{\the\pgfarrowwidth}
    \dimen7 2\pgfarrowwidth\pgfmathdivide{\pgfmathresult}{\the\dimen7}
    \dimen7 \pgfmathresult\pgfarrowlinewidth
    \pgfarrowssettipend{+\dimen7}
    \pgfarrowssetbackend{+-\pgfarrowlength}
    \dimen9 -0.5\pgfarrowlinewidth
    \pgfarrowssetvisualbackend{+\dimen9}
    \pgfarrowssetlineend{+-0.5\pgfarrowlinewidth}
    \pgfarrowshullpoint{+\dimen7}{+0pt}
    \pgfarrowsupperhullpoint{+-\pgfarrowlength}{+0.5\pgfarrowwidth}
    \pgfarrowssavethe\pgfarrowlinewidth
    \pgfarrowssavethe\pgfarrowlength
    \pgfarrowssavethe\pgfarrowwidth
  },
  drawing code = {\pgfsetdash{}{+0pt}
    \ifdim\pgfarrowlinewidth=\pgflinewidth\else\pgfsetlinewidth{+\pgfarrowlinewidth}\fi
    \pgfpathmoveto{\pgfqpoint{-\pgfarrowlength}{0.5\pgfarrowwidth}}
    \pgfpathlineto{\pgfqpoint{0pt}{0pt}}
    \pgfpathlineto{\pgfqpoint{-\pgfarrowlength}{-0.5\pgfarrowwidth}}
    \pgfusepathqstroke
  }
}
\clip(1485,-14046) rectangle (13707,-1761);
\tikzset{inner sep=+0pt, outer sep=+0pt}
\pgfsetlinewidth{+7.5\XFigu}
\pgfsetstrokecolor{black}
\pgfsetfillcolor{cyan}
\filldraw  (10950,-3525) circle [radius=+75];
\pgfsetlinewidth{+30\XFigu}
\draw (4050,-2695) arc[start angle=+135.0, end angle=+45.0, radius=+300.5];
\draw (2525,-2948) arc[start angle=+-138.5, end angle=+-41.5, radius=+508.3];
\draw (2693,-3033) arc[start angle=+130.6, end angle=+49.4, radius=+390.7];
\draw (3880,-4560) arc[start angle=+-138.4, end angle=+-41.6, radius=+511.8];
\draw (4050,-4645) arc[start angle=+130.6, end angle=+49.4, radius=+391.8];
\pgfsetbeveljoin
\pgfsetlinewidth{+45\XFigu}
\draw (2761,-4942)--(2743,-4936)--(2725,-4929)--(2708,-4920)--(2690,-4910)--(2673,-4899)
  --(2655,-4886)--(2637,-4872)--(2620,-4857)--(2603,-4839)--(2585,-4821)--(2568,-4800)
  --(2551,-4778)--(2534,-4755)--(2517,-4730)--(2501,-4704)--(2485,-4677)--(2469,-4648)
  --(2453,-4618)--(2438,-4587)--(2424,-4555)--(2410,-4522)--(2397,-4488)--(2384,-4454)
  --(2371,-4419)--(2360,-4383)--(2348,-4347)--(2338,-4311)--(2328,-4274)--(2318,-4236)
  --(2309,-4198)--(2300,-4159)--(2292,-4120)--(2285,-4086)--(2278,-4051)--(2272,-4016)
  --(2266,-3980)--(2260,-3943)--(2255,-3905)--(2249,-3866)--(2244,-3827)--(2240,-3786)
  --(2235,-3745)--(2231,-3704)--(2227,-3661)--(2224,-3618)--(2221,-3575)--(2218,-3531)
  --(2216,-3487)--(2214,-3443)--(2213,-3399)--(2212,-3355)--(2211,-3311)--(2211,-3268)
  --(2212,-3225)--(2212,-3183)--(2214,-3141)--(2216,-3100)--(2218,-3060)--(2221,-3022)
  --(2224,-2984)--(2228,-2947)--(2232,-2912)--(2236,-2877)--(2241,-2844)--(2247,-2812)
  --(2253,-2782)--(2259,-2752)--(2266,-2724)--(2274,-2697)--(2281,-2671)--(2291,-2641)
  --(2302,-2613)--(2314,-2586)--(2326,-2560)--(2339,-2535)--(2353,-2512)--(2368,-2489)
  --(2384,-2468)--(2401,-2448)--(2418,-2429)--(2437,-2412)--(2456,-2396)--(2476,-2381)
  --(2497,-2367)--(2518,-2355)--(2540,-2345)--(2562,-2335)--(2585,-2327)--(2608,-2321)
  --(2632,-2316)--(2655,-2312)--(2679,-2310)--(2703,-2309)--(2726,-2309)--(2750,-2311)
  --(2773,-2314)--(2797,-2318)--(2820,-2322)--(2844,-2328)--(2867,-2335)--(2890,-2343)
  --(2913,-2352)--(2940,-2363)--(2966,-2376)--(2994,-2390)--(3021,-2405)--(3049,-2422)
  --(3076,-2440)--(3105,-2459)--(3133,-2480)--(3161,-2502)--(3189,-2525)--(3217,-2549)
  --(3244,-2574)--(3270,-2600)--(3296,-2627)--(3321,-2654)--(3345,-2682)--(3368,-2710)
  --(3389,-2738)--(3409,-2766)--(3428,-2794)--(3444,-2822)--(3460,-2849)--(3473,-2877)
  --(3485,-2903)--(3495,-2929)--(3503,-2955)--(3510,-2980)--(3515,-3005)--(3518,-3032)
  --(3520,-3058)--(3520,-3085)--(3518,-3111)--(3514,-3138)--(3509,-3165)--(3503,-3193)
  --(3495,-3221)--(3486,-3249)--(3476,-3278)--(3465,-3307)--(3453,-3336)--(3441,-3366)
  --(3428,-3395)--(3415,-3424)--(3402,-3454)--(3389,-3483)--(3377,-3512)--(3365,-3541)
  --(3353,-3569)--(3343,-3597)--(3333,-3625)--(3324,-3653)--(3316,-3681)--(3309,-3709)
  --(3304,-3737)--(3300,-3764)--(3297,-3791)--(3295,-3818)--(3294,-3846)--(3294,-3875)
  --(3294,-3905)--(3295,-3936)--(3296,-3967)--(3298,-3999)--(3301,-4031)--(3304,-4064)
  --(3306,-4097)--(3309,-4130)--(3312,-4163)--(3315,-4197)--(3318,-4230)--(3320,-4263)
  --(3322,-4295)--(3324,-4327)--(3325,-4358)--(3325,-4388)--(3325,-4417)--(3323,-4446)
  --(3321,-4474)--(3318,-4500)--(3314,-4527)--(3309,-4552)--(3302,-4577)--(3295,-4601)
  --(3286,-4625)--(3275,-4648)--(3264,-4671)--(3251,-4694)--(3237,-4716)--(3221,-4738)
  --(3205,-4760)--(3187,-4780)--(3168,-4800)--(3148,-4820)--(3127,-4838)--(3105,-4855)
  --(3082,-4871)--(3059,-4886)--(3036,-4899)--(3012,-4911)--(2988,-4922)--(2964,-4931)
  --(2940,-4939)--(2916,-4944)--(2893,-4949)--(2870,-4951)--(2847,-4953)--(2825,-4952)
  --(2804,-4950)--(2782,-4947)--cycle;
\pgfsetdash{}{+0pt}
\draw (3467,-2608)--(3473,-2595)--(3480,-2583)--(3488,-2570)--(3497,-2557)--(3507,-2544)
  --(3518,-2531)--(3531,-2517)--(3545,-2504)--(3560,-2491)--(3577,-2478)--(3595,-2465)
  --(3614,-2453)--(3634,-2442)--(3655,-2431)--(3677,-2420)--(3701,-2410)--(3725,-2401)
  --(3751,-2392)--(3770,-2386)--(3790,-2380)--(3811,-2374)--(3833,-2368)--(3856,-2363)
  --(3880,-2357)--(3905,-2352)--(3931,-2347)--(3958,-2342)--(3986,-2337)--(4014,-2333)
  --(4044,-2329)--(4073,-2326)--(4103,-2323)--(4133,-2320)--(4162,-2318)--(4192,-2316)
  --(4221,-2315)--(4249,-2315)--(4277,-2315)--(4305,-2315)--(4331,-2316)--(4357,-2317)
  --(4382,-2319)--(4407,-2322)--(4431,-2325)--(4457,-2329)--(4482,-2334)--(4507,-2339)
  --(4532,-2345)--(4558,-2352)--(4583,-2360)--(4608,-2369)--(4634,-2378)--(4659,-2389)
  --(4684,-2400)--(4710,-2412)--(4735,-2425)--(4759,-2438)--(4783,-2452)--(4807,-2466)
  --(4830,-2481)--(4853,-2496)--(4875,-2511)--(4897,-2527)--(4917,-2542)--(4938,-2558)
  --(4958,-2574)--(4977,-2591)--(4996,-2608)--(5014,-2623)--(5032,-2640)--(5050,-2656)
  --(5068,-2674)--(5086,-2692)--(5105,-2711)--(5124,-2731)--(5142,-2752)--(5161,-2774)
  --(5180,-2797)--(5199,-2821)--(5217,-2846)--(5236,-2871)--(5253,-2897)--(5271,-2924)
  --(5288,-2951)--(5304,-2979)--(5320,-3007)--(5335,-3035)--(5349,-3064)--(5362,-3093)
  --(5375,-3122)--(5386,-3152)--(5398,-3182)--(5408,-3212)--(5418,-3243)--(5426,-3271)
  --(5433,-3299)--(5440,-3329)--(5447,-3359)--(5453,-3390)--(5459,-3422)--(5465,-3455)
  --(5469,-3489)--(5474,-3524)--(5478,-3559)--(5481,-3596)--(5484,-3633)--(5487,-3670)
  --(5489,-3708)--(5490,-3746)--(5490,-3784)--(5490,-3822)--(5490,-3860)--(5488,-3897)
  --(5487,-3934)--(5484,-3970)--(5481,-4006)--(5478,-4040)--(5473,-4074)--(5469,-4107)
  --(5464,-4139)--(5458,-4170)--(5452,-4201)--(5445,-4230)--(5438,-4259)--(5429,-4290)
  --(5420,-4320)--(5410,-4349)--(5399,-4378)--(5388,-4408)--(5375,-4437)--(5362,-4465)
  --(5347,-4494)--(5332,-4522)--(5316,-4551)--(5298,-4578)--(5280,-4606)--(5261,-4633)
  --(5242,-4659)--(5221,-4685)--(5200,-4710)--(5179,-4734)--(5157,-4757)--(5135,-4780)
  --(5112,-4801)--(5089,-4822)--(5066,-4841)--(5042,-4860)--(5019,-4878)--(4995,-4895)
  --(4971,-4911)--(4947,-4927)--(4922,-4942)--(4897,-4956)--(4871,-4971)--(4844,-4984)
  --(4817,-4998)--(4789,-5011)--(4759,-5024)--(4730,-5037)--(4699,-5049)--(4667,-5061)
  --(4635,-5073)--(4603,-5084)--(4570,-5095)--(4536,-5105)--(4503,-5115)--(4469,-5124)
  --(4436,-5133)--(4402,-5141)--(4370,-5149)--(4337,-5156)--(4306,-5163)--(4275,-5168)
  --(4245,-5174)--(4216,-5178)--(4187,-5183)--(4160,-5186)--(4133,-5189)--(4107,-5192)
  --(4081,-5194)--(4052,-5196)--(4023,-5198)--(3995,-5199)--(3966,-5199)--(3938,-5199)
  --(3910,-5198)--(3882,-5197)--(3855,-5195)--(3827,-5192)--(3800,-5189)--(3774,-5185)
  --(3748,-5180)--(3722,-5175)--(3698,-5170)--(3674,-5163)--(3651,-5157)--(3629,-5150)
  --(3608,-5143)--(3588,-5135)--(3569,-5127)--(3551,-5119)--(3534,-5110)--(3517,-5101)
  --(3502,-5092)--(3485,-5081)--(3468,-5070)--(3452,-5058)--(3436,-5046)--(3421,-5033)
  --(3406,-5019)--(3392,-5004)--(3378,-4989)--(3365,-4973)--(3352,-4957)--(3341,-4940)
  --(3330,-4923)--(3321,-4905)--(3313,-4888)--(3305,-4870)--(3299,-4853)--(3294,-4835)
  --(3291,-4818)--(3288,-4801)--(3287,-4784)--(3286,-4768)--(3287,-4751)--(3289,-4736)
  --(3292,-4720)--(3295,-4704)--(3300,-4688)--(3306,-4672)--(3313,-4655)--(3321,-4638)
  --(3330,-4620)--(3340,-4603)--(3352,-4585)--(3364,-4567)--(3378,-4550)--(3393,-4532)
  --(3408,-4515)--(3425,-4498)--(3442,-4482)--(3460,-4466)--(3478,-4450)--(3497,-4435)
  --(3517,-4421)--(3537,-4407)--(3557,-4394)--(3579,-4381)--(3600,-4369)--(3621,-4357)
  --(3642,-4346)--(3665,-4336)--(3688,-4325)--(3712,-4314)--(3737,-4303)--(3763,-4291)
  --(3790,-4280)--(3818,-4268)--(3846,-4256)--(3876,-4244)--(3905,-4231)--(3935,-4218)
  --(3965,-4205)--(3996,-4191)--(4026,-4177)--(4056,-4163)--(4085,-4148)--(4114,-4133)
  --(4143,-4118)--(4170,-4102)--(4198,-4087)--(4224,-4070)--(4250,-4053)--(4275,-4036)
  --(4300,-4018)--(4323,-4000)--(4346,-3982)--(4368,-3963)--(4391,-3943)--(4413,-3922)
  --(4436,-3900)--(4458,-3878)--(4480,-3854)--(4502,-3830)--(4524,-3806)--(4545,-3780)
  --(4566,-3755)--(4586,-3729)--(4606,-3703)--(4625,-3677)--(4644,-3651)--(4661,-3626)
  --(4678,-3600)--(4693,-3576)--(4708,-3552)--(4721,-3528)--(4733,-3506)--(4745,-3484)
  --(4755,-3463)--(4764,-3443)--(4772,-3424)--(4779,-3406)--(4786,-3388)--(4792,-3369)
  --(4797,-3350)--(4800,-3332)--(4803,-3314)--(4804,-3297)--(4805,-3281)--(4804,-3265)
  --(4802,-3250)--(4798,-3235)--(4794,-3221)--(4789,-3208)--(4782,-3195)--(4775,-3182)
  --(4766,-3171)--(4757,-3160)--(4748,-3149)--(4737,-3139)--(4726,-3130)--(4715,-3122)
  --(4703,-3114)--(4690,-3106)--(4678,-3099)--(4665,-3092)--(4652,-3085)--(4637,-3078)
  --(4622,-3070)--(4605,-3063)--(4588,-3056)--(4570,-3049)--(4551,-3041)--(4531,-3034)
  --(4510,-3027)--(4489,-3020)--(4466,-3012)--(4443,-3005)--(4419,-2999)--(4394,-2992)
  --(4369,-2985)--(4344,-2979)--(4319,-2973)--(4293,-2968)--(4268,-2962)--(4242,-2957)
  --(4217,-2953)--(4191,-2948)--(4164,-2944)--(4140,-2940)--(4114,-2936)--(4088,-2932)
  --(4062,-2928)--(4034,-2925)--(4006,-2921)--(3977,-2917)--(3948,-2913)--(3918,-2909)
  --(3889,-2905)--(3859,-2901)--(3829,-2896)--(3800,-2892)--(3772,-2887)--(3744,-2882)
  --(3718,-2878)--(3692,-2873)--(3668,-2868)--(3646,-2862)--(3624,-2857)--(3605,-2852)
  --(3586,-2846)--(3569,-2840)--(3554,-2834)--(3538,-2827)--(3524,-2819)--(3511,-2811)
  --(3499,-2803)--(3488,-2794)--(3479,-2784)--(3471,-2774)--(3464,-2764)--(3458,-2753)
  --(3453,-2742)--(3449,-2731)--(3447,-2719)--(3445,-2707)--(3445,-2696)--(3445,-2684)
  --(3447,-2673)--(3449,-2661)--(3451,-2650)--(3455,-2639)--(3459,-2628)--(3463,-2618)--cycle;
\pgfsetlinewidth{+30\XFigu}
\pgfsetdash{}{+0pt}
\draw (4925,-7620) arc[start angle=+-134.22, end angle=+-58.47, radius=+626.9];
\draw (5095,-7705) arc[start angle=+135.0, end angle=+45.0, radius=+300.5];
\draw (2260,-7845) arc[start angle=+-138.5, end angle=+-41.5, radius=+508.3];
\draw (2428,-7930) arc[start angle=+130.6, end angle=+49.4, radius=+390.7];
\pgfsetlinewidth{+45\XFigu}
\pgfsetdash{}{+0pt}
\draw (3150,-7350)--(3147,-7348)--(3140,-7345)--(3129,-7339)--(3111,-7329)--(3088,-7317)
  --(3060,-7303)--(3028,-7286)--(2994,-7269)--(2959,-7251)--(2923,-7234)--(2889,-7218)
  --(2857,-7203)--(2826,-7190)--(2797,-7178)--(2771,-7168)--(2746,-7160)--(2722,-7153)
  --(2700,-7147)--(2679,-7143)--(2658,-7140)--(2638,-7138)--(2620,-7137)--(2603,-7136)
  --(2586,-7137)--(2569,-7138)--(2551,-7141)--(2533,-7144)--(2516,-7149)--(2498,-7154)
  --(2480,-7161)--(2462,-7169)--(2444,-7178)--(2426,-7189)--(2409,-7201)--(2391,-7214)
  --(2374,-7228)--(2358,-7244)--(2342,-7261)--(2327,-7279)--(2312,-7298)--(2298,-7318)
  --(2284,-7339)--(2271,-7362)--(2259,-7385)--(2247,-7410)--(2236,-7435)--(2225,-7463)
  --(2217,-7484)--(2209,-7507)--(2202,-7530)--(2194,-7555)--(2187,-7581)--(2180,-7608)
  --(2173,-7636)--(2167,-7666)--(2161,-7697)--(2154,-7728)--(2149,-7762)--(2143,-7796)
  --(2138,-7831)--(2133,-7867)--(2129,-7904)--(2125,-7942)--(2121,-7981)--(2117,-8021)
  --(2115,-8060)--(2112,-8101)--(2110,-8141)--(2109,-8182)--(2107,-8223)--(2107,-8263)
  --(2107,-8304)--(2107,-8345)--(2108,-8385)--(2109,-8425)--(2110,-8465)--(2112,-8505)
  --(2115,-8544)--(2118,-8584)--(2121,-8623)--(2125,-8663)--(2129,-8700)--(2133,-8737)
  --(2138,-8775)--(2144,-8814)--(2149,-8852)--(2156,-8892)--(2162,-8931)--(2169,-8971)
  --(2177,-9012)--(2185,-9053)--(2193,-9094)--(2202,-9135)--(2211,-9177)--(2221,-9218)
  --(2231,-9260)--(2241,-9301)--(2251,-9342)--(2262,-9383)--(2273,-9423)--(2284,-9463)
  --(2295,-9502)--(2307,-9539)--(2318,-9576)--(2330,-9612)--(2341,-9647)--(2353,-9681)
  --(2364,-9714)--(2375,-9745)--(2387,-9775)--(2398,-9803)--(2409,-9831)--(2420,-9857)
  --(2431,-9882)--(2441,-9906)--(2452,-9928)--(2463,-9950)--(2476,-9976)--(2489,-10001)
  --(2503,-10024)--(2517,-10046)--(2530,-10066)--(2544,-10085)--(2558,-10103)
  --(2572,-10120)--(2587,-10135)--(2601,-10148)--(2616,-10161)--(2630,-10171)
  --(2645,-10181)--(2659,-10189)--(2673,-10196)--(2688,-10201)--(2702,-10205)
  --(2715,-10207)--(2729,-10209)--(2742,-10209)--(2755,-10208)--(2768,-10206)
  --(2780,-10203)--(2792,-10199)--(2804,-10194)--(2815,-10188)--(2826,-10182)
  --(2838,-10175)--(2852,-10165)--(2865,-10154)--(2880,-10141)--(2894,-10126)
  --(2909,-10110)--(2925,-10091)--(2941,-10070)--(2959,-10047)--(2977,-10020)
  --(2997,-9992)--(3017,-9962)--(3038,-9930)--(3059,-9897)--(3079,-9865)--(3098,-9836)
  --(3114,-9809)--(3128,-9787)--(3138,-9770)--(3144,-9759)--(3148,-9753)--(3150,-9750);
\draw (3150,-7950) arc[start angle=+138.89, end angle=+221.11, radius=+912.5];
\pgfsetlinewidth{+30\XFigu}
\pgfsetdash{}{+0pt}
\draw (5056,-12319) arc[start angle=+28.10, end angle=+-28.10, radius=+1043.5];
\draw (8880,-12278) arc[start angle=+150.16, end angle=+209.84, radius=+1024.9];
\draw (11931,-12247) arc[start angle=+29.79, end angle=+-25.63, radius=+1095.3];
\draw (11743,-7385) arc[start angle=+29.87, end angle=+-25.66, radius=+643.3];
\draw (9600,-7350) arc[start angle=+150.26, end angle=+209.74, radius=+604.7];
\pgfsetdash{}{+0pt}
\draw (9450,-3525) arc[start angle=+-133.60, end angle=+-46.40, radius=+1087.5];
\draw (9450,-3525) arc[start angle=+133.60, end angle=+46.40, radius=+1087.5];
\pgfsetdash{}{+0pt}
\draw (4575,-9150) arc[start angle=+28.07, end angle=+-28.07, radius=+637.5];
\pgfsetlinewidth{+45\XFigu}
\pgfsetdash{}{+0pt}
\draw (3150,-9750) arc[start angle=+113.84, end angle=+66.16, radius=+1762.9];
\draw (3150,-9150) arc[start angle=+-113.84, end angle=+-66.16, radius=+1762.9];
\pgfsetlinewidth{+30\XFigu}
\pgfsetdash{}{+0pt}
\draw (3150,-9150) arc[start angle=+150.26, end angle=+209.74, radius=+604.7];
\pgfsetdash{}{+0pt}
\draw (5100,-9505) arc[start angle=+130.6, end angle=+49.4, radius=+391.8];
\draw (4930,-9420) arc[start angle=+-138.4, end angle=+-41.6, radius=+511.8];
\pgfsetlinewidth{+45\XFigu}
\pgfsetdash{}{+0pt}
\draw (9600,-7950) arc[start angle=+138.89, end angle=+221.11, radius=+912.5];
\pgfsetlinewidth{+30\XFigu}
\pgfsetdash{}{+0pt}
\draw (8878,-7930) arc[start angle=+130.6, end angle=+49.4, radius=+390.7];
\draw (8710,-7845) arc[start angle=+-138.5, end angle=+-41.5, radius=+508.3];
\draw (12342,-9521) arc[start angle=+133.2, end angle=+51.1, radius=+389.5];
\draw (12169,-9441) arc[start angle=+-136.3, end angle=+-40.0, radius=+513.9];
\draw (12275,-7723) arc[start angle=+137.6, end angle=+46.4, radius=+296.3];
\draw (12101,-7644) arc[start angle=+-132.30, end angle=+-56.04, radius=+623.5];
\pgfsetfillpattern{xfigp0}{black}
\draw[pattern,preaction={fill=black}]  (2700,-4350) circle [radius=+38];
\draw (3880,-2610) arc[start angle=+-134.22, end angle=+-58.47, radius=+626.9];
\pgfsetlinewidth{+7.5\XFigu}
\pgfsetfillcolor{cyan}
\filldraw  (9450,-3525) circle [radius=+75];
\pgfsetdash{{+60\XFigu}{+60\XFigu}}{++0pt}
\draw (2100,-5850)--(12300,-5850);
\pgfsetlinewidth{+45\XFigu}
\pgfsetdash{}{+0pt}
\draw (4575,-9750)--(4577,-9752)--(4581,-9757)--(4588,-9766)--(4599,-9778)--(4612,-9794)
  --(4628,-9812)--(4646,-9833)--(4665,-9854)--(4685,-9874)--(4704,-9894)--(4723,-9913)
  --(4742,-9930)--(4760,-9946)--(4780,-9961)--(4800,-9975)--(4816,-9985)--(4832,-9996)
  --(4850,-10007)--(4868,-10018)--(4886,-10030)--(4905,-10043)--(4924,-10056)
  --(4944,-10070)--(4963,-10084)--(4984,-10098)--(5004,-10112)--(5024,-10126)
  --(5045,-10139)--(5066,-10152)--(5088,-10163)--(5109,-10174)--(5131,-10183)
  --(5154,-10190)--(5177,-10196)--(5201,-10199)--(5225,-10201)--(5250,-10200)
  --(5272,-10198)--(5295,-10194)--(5318,-10189)--(5343,-10183)--(5368,-10175)
  --(5394,-10167)--(5420,-10157)--(5448,-10147)--(5476,-10135)--(5504,-10124)
  --(5533,-10111)--(5562,-10098)--(5591,-10085)--(5620,-10072)--(5649,-10058)
  --(5678,-10044)--(5706,-10030)--(5733,-10016)--(5760,-10002)--(5787,-9988)
  --(5812,-9974)--(5837,-9960)--(5860,-9945)--(5883,-9931)--(5904,-9916)--(5925,-9900)
  --(5946,-9882)--(5967,-9864)--(5986,-9845)--(6005,-9826)--(6023,-9806)--(6041,-9785)
  --(6058,-9764)--(6074,-9743)--(6091,-9721)--(6106,-9699)--(6122,-9677)--(6137,-9655)
  --(6152,-9633)--(6167,-9610)--(6182,-9587)--(6196,-9565)--(6210,-9542)--(6224,-9518)
  --(6238,-9495)--(6251,-9472)--(6264,-9448)--(6277,-9424)--(6289,-9400)--(6300,-9375)
  --(6310,-9352)--(6319,-9328)--(6328,-9304)--(6336,-9279)--(6344,-9254)--(6352,-9229)
  --(6359,-9203)--(6366,-9177)--(6373,-9151)--(6380,-9125)--(6387,-9098)--(6393,-9071)
  --(6400,-9045)--(6406,-9018)--(6412,-8991)--(6417,-8964)--(6423,-8937)--(6428,-8910)
  --(6432,-8884)--(6437,-8857)--(6440,-8831)--(6444,-8804)--(6446,-8778)--(6448,-8752)
  --(6450,-8726)--(6450,-8700)--(6450,-8674)--(6448,-8648)--(6446,-8622)--(6444,-8596)
  --(6440,-8569)--(6437,-8543)--(6432,-8516)--(6428,-8490)--(6423,-8463)--(6417,-8436)
  --(6412,-8409)--(6406,-8382)--(6400,-8355)--(6393,-8329)--(6387,-8302)--(6380,-8275)
  --(6373,-8249)--(6366,-8223)--(6359,-8197)--(6352,-8171)--(6344,-8146)--(6336,-8121)
  --(6328,-8096)--(6319,-8072)--(6310,-8048)--(6300,-8025)--(6289,-8000)--(6277,-7976)
  --(6264,-7952)--(6250,-7928)--(6236,-7905)--(6222,-7882)--(6207,-7858)--(6191,-7835)
  --(6176,-7813)--(6159,-7790)--(6143,-7767)--(6127,-7745)--(6110,-7723)--(6093,-7701)
  --(6076,-7679)--(6060,-7657)--(6043,-7636)--(6026,-7615)--(6009,-7594)--(5993,-7574)
  --(5976,-7555)--(5959,-7536)--(5942,-7518)--(5925,-7500)--(5906,-7482)--(5887,-7464)
  --(5868,-7448)--(5848,-7432)--(5828,-7416)--(5808,-7401)--(5788,-7387)--(5768,-7373)
  --(5748,-7359)--(5727,-7346)--(5707,-7332)--(5687,-7319)--(5666,-7307)--(5645,-7294)
  --(5625,-7281)--(5604,-7269)--(5583,-7257)--(5562,-7245)--(5540,-7233)--(5519,-7222)
  --(5497,-7211)--(5475,-7200)--(5453,-7190)--(5430,-7180)--(5407,-7171)--(5383,-7162)
  --(5359,-7153)--(5335,-7143)--(5310,-7134)--(5285,-7125)--(5260,-7115)--(5234,-7106)
  --(5209,-7097)--(5184,-7088)--(5158,-7080)--(5133,-7072)--(5108,-7065)--(5084,-7059)
  --(5060,-7054)--(5037,-7050)--(5014,-7048)--(4992,-7047)--(4971,-7048)--(4950,-7050)
  --(4928,-7055)--(4906,-7062)--(4883,-7073)--(4861,-7086)--(4837,-7102)--(4812,-7121)
  --(4786,-7143)--(4759,-7167)--(4732,-7192)--(4705,-7218)--(4679,-7243)--(4655,-7267)
  --(4633,-7290)--(4614,-7309)--(4599,-7324)--(4588,-7336)--(4581,-7344)--(4577,-7348)
  --(4575,-7350);
\draw (4575,-7950)--(4578,-7952)--(4585,-7955)--(4596,-7961)--(4613,-7969)--(4634,-7980)
  --(4658,-7992)--(4684,-8004)--(4710,-8018)--(4736,-8030)--(4760,-8042)--(4783,-8054)
  --(4804,-8064)--(4824,-8075)--(4844,-8084)--(4862,-8094)--(4881,-8103)--(4900,-8113)
  --(4917,-8121)--(4935,-8130)--(4954,-8139)--(4973,-8149)--(4992,-8159)--(5012,-8170)
  --(5033,-8181)--(5053,-8192)--(5074,-8204)--(5095,-8216)--(5115,-8228)--(5135,-8240)
  --(5154,-8253)--(5173,-8265)--(5190,-8277)--(5207,-8290)--(5222,-8302)--(5237,-8314)
  --(5250,-8326)--(5263,-8338)--(5275,-8351)--(5287,-8365)--(5299,-8379)--(5309,-8394)
  --(5319,-8410)--(5327,-8426)--(5335,-8443)--(5343,-8461)--(5349,-8478)--(5354,-8496)
  --(5358,-8515)--(5361,-8533)--(5363,-8551)--(5365,-8568)--(5365,-8586)--(5365,-8603)
  --(5364,-8620)--(5363,-8638)--(5360,-8655)--(5357,-8673)--(5353,-8691)--(5349,-8710)
  --(5343,-8729)--(5336,-8749)--(5328,-8769)--(5318,-8789)--(5308,-8809)--(5296,-8829)
  --(5283,-8848)--(5269,-8867)--(5253,-8885)--(5237,-8903)--(5220,-8919)--(5202,-8934)
  --(5183,-8949)--(5163,-8963)--(5146,-8972)--(5129,-8982)--(5111,-8991)--(5091,-9000)
  --(5069,-9009)--(5045,-9019)--(5020,-9028)--(4991,-9037)--(4961,-9047)--(4927,-9057)
  --(4892,-9068)--(4854,-9078)--(4815,-9089)--(4776,-9099)--(4737,-9110)--(4700,-9119)
  --(4666,-9128)--(4637,-9135)--(4613,-9141)--(4596,-9145)--(4584,-9148)--(4578,-9149)
  --(4575,-9150);
\draw (9600,-7350)--(9597,-7348)--(9590,-7345)--(9579,-7339)--(9561,-7329)--(9538,-7317)
  --(9510,-7303)--(9478,-7286)--(9444,-7269)--(9409,-7251)--(9373,-7234)--(9339,-7218)
  --(9307,-7203)--(9276,-7190)--(9247,-7178)--(9221,-7168)--(9196,-7160)--(9172,-7153)
  --(9150,-7147)--(9129,-7143)--(9108,-7140)--(9088,-7138)--(9070,-7137)--(9053,-7136)
  --(9036,-7137)--(9019,-7138)--(9001,-7141)--(8983,-7144)--(8966,-7149)--(8948,-7154)
  --(8930,-7161)--(8912,-7169)--(8894,-7178)--(8876,-7189)--(8859,-7201)--(8841,-7214)
  --(8824,-7228)--(8808,-7244)--(8792,-7261)--(8777,-7279)--(8762,-7298)--(8748,-7318)
  --(8734,-7339)--(8721,-7362)--(8709,-7385)--(8697,-7410)--(8686,-7435)--(8675,-7463)
  --(8667,-7484)--(8659,-7507)--(8652,-7530)--(8644,-7555)--(8637,-7581)--(8630,-7608)
  --(8623,-7636)--(8617,-7666)--(8611,-7697)--(8604,-7728)--(8599,-7762)--(8593,-7796)
  --(8588,-7831)--(8583,-7867)--(8579,-7904)--(8575,-7942)--(8571,-7981)--(8567,-8021)
  --(8565,-8060)--(8562,-8101)--(8560,-8141)--(8559,-8182)--(8557,-8223)--(8557,-8263)
  --(8557,-8304)--(8557,-8345)--(8558,-8385)--(8559,-8425)--(8560,-8465)--(8562,-8505)
  --(8565,-8544)--(8568,-8584)--(8571,-8623)--(8575,-8663)--(8579,-8700)--(8583,-8737)
  --(8588,-8775)--(8594,-8814)--(8599,-8852)--(8606,-8892)--(8612,-8931)--(8619,-8971)
  --(8627,-9012)--(8635,-9053)--(8643,-9094)--(8652,-9135)--(8661,-9177)--(8671,-9218)
  --(8681,-9260)--(8691,-9301)--(8701,-9342)--(8712,-9383)--(8723,-9423)--(8734,-9463)
  --(8745,-9502)--(8757,-9539)--(8768,-9576)--(8780,-9612)--(8791,-9647)--(8803,-9681)
  --(8814,-9714)--(8825,-9745)--(8837,-9775)--(8848,-9803)--(8859,-9831)--(8870,-9857)
  --(8881,-9882)--(8891,-9906)--(8902,-9928)--(8913,-9950)--(8926,-9976)--(8939,-10001)
  --(8953,-10024)--(8967,-10046)--(8980,-10066)--(8994,-10085)--(9008,-10103)
  --(9022,-10120)--(9037,-10135)--(9051,-10148)--(9066,-10161)--(9080,-10171)
  --(9095,-10181)--(9109,-10189)--(9123,-10196)--(9138,-10201)--(9152,-10205)
  --(9165,-10207)--(9179,-10209)--(9192,-10209)--(9205,-10208)--(9218,-10206)
  --(9230,-10203)--(9242,-10199)--(9254,-10194)--(9265,-10188)--(9276,-10182)
  --(9288,-10175)--(9302,-10165)--(9315,-10154)--(9330,-10141)--(9344,-10126)
  --(9359,-10110)--(9375,-10091)--(9391,-10070)--(9409,-10047)--(9427,-10020)
  --(9447,-9992)--(9467,-9962)--(9488,-9930)--(9509,-9897)--(9529,-9865)--(9548,-9836)
  --(9564,-9809)--(9578,-9787)--(9588,-9770)--(9594,-9759)--(9598,-9753)--(9600,-9750);
\draw (11765,-7984)--(11768,-7985)--(11775,-7989)--(11787,-7994)--(11803,-8002)--(11825,-8012)
  --(11849,-8023)--(11875,-8035)--(11901,-8047)--(11927,-8059)--(11952,-8071)
  --(11975,-8081)--(11996,-8091)--(12017,-8101)--(12036,-8110)--(12056,-8119)
  --(12074,-8128)--(12094,-8136)--(12111,-8144)--(12129,-8153)--(12148,-8161)
  --(12167,-8171)--(12187,-8180)--(12207,-8190)--(12228,-8200)--(12249,-8211)
  --(12270,-8222)--(12291,-8233)--(12312,-8245)--(12333,-8256)--(12352,-8268)
  --(12371,-8279)--(12389,-8291)--(12406,-8303)--(12422,-8314)--(12436,-8325)
  --(12450,-8337)--(12463,-8348)--(12477,-8361)--(12489,-8375)--(12501,-8389)
  --(12512,-8403)--(12522,-8418)--(12531,-8434)--(12540,-8451)--(12548,-8468)
  --(12554,-8485)--(12560,-8503)--(12565,-8521)--(12569,-8539)--(12572,-8557)
  --(12574,-8574)--(12575,-8592)--(12575,-8609)--(12575,-8626)--(12574,-8643)
  --(12573,-8659)--(12571,-8675)--(12568,-8691)--(12565,-8708)--(12560,-8725)
  --(12555,-8743)--(12549,-8761)--(12542,-8780)--(12534,-8798)--(12526,-8817)
  --(12516,-8835)--(12505,-8853)--(12493,-8871)--(12480,-8888)--(12466,-8904)
  --(12452,-8920)--(12437,-8935)--(12421,-8949)--(12404,-8963)--(12386,-8976)
  --(12370,-8986)--(12353,-8996)--(12335,-9006)--(12315,-9016)--(12294,-9026)
  --(12271,-9036)--(12245,-9046)--(12217,-9057)--(12187,-9067)--(12154,-9079)
  --(12118,-9090)--(12081,-9102)--(12042,-9114)--(12003,-9126)--(11965,-9138)
  --(11928,-9149)--(11894,-9158)--(11866,-9167)--(11842,-9173)--(11825,-9178)
  --(11813,-9181)--(11807,-9183)--(11804,-9184);
\draw (11827,-9783)--(11829,-9785)--(11833,-9790)--(11841,-9798)--(11851,-9810)--(11865,-9825)
  --(11882,-9844)--(11901,-9863)--(11920,-9884)--(11940,-9904)--(11960,-9923)
  --(11980,-9942)--(11999,-9958)--(12018,-9973)--(12038,-9987)--(12059,-10001)
  --(12075,-10011)--(12092,-10021)--(12110,-10031)--(12128,-10042)--(12147,-10054)
  --(12166,-10066)--(12186,-10078)--(12206,-10091)--(12226,-10105)--(12247,-10118)
  --(12268,-10131)--(12289,-10144)--(12310,-10157)--(12332,-10169)--(12353,-10180)
  --(12375,-10189)--(12398,-10198)--(12421,-10204)--(12444,-10209)--(12468,-10212)
  --(12492,-10212)--(12517,-10211)--(12539,-10208)--(12561,-10204)--(12585,-10198)
  --(12609,-10191)--(12634,-10182)--(12659,-10173)--(12686,-10162)--(12713,-10151)
  --(12740,-10139)--(12768,-10126)--(12797,-10113)--(12825,-10099)--(12854,-10085)
  --(12882,-10070)--(12911,-10055)--(12939,-10041)--(12967,-10026)--(12994,-10011)
  --(13020,-9996)--(13046,-9981)--(13071,-9966)--(13095,-9951)--(13118,-9936)
  --(13140,-9920)--(13161,-9904)--(13181,-9888)--(13202,-9870)--(13221,-9851)
  --(13240,-9831)--(13258,-9811)--(13276,-9790)--(13293,-9769)--(13309,-9748)
  --(13325,-9726)--(13340,-9704)--(13355,-9681)--(13370,-9659)--(13384,-9636)
  --(13398,-9613)--(13412,-9590)--(13426,-9566)--(13440,-9543)--(13453,-9520)
  --(13466,-9496)--(13479,-9472)--(13492,-9448)--(13504,-9424)--(13515,-9400)
  --(13527,-9375)--(13537,-9350)--(13546,-9326)--(13554,-9302)--(13562,-9278)
  --(13570,-9253)--(13577,-9228)--(13583,-9202)--(13590,-9176)--(13596,-9150)
  --(13602,-9124)--(13608,-9097)--(13614,-9070)--(13619,-9043)--(13625,-9016)
  --(13630,-8989)--(13635,-8962)--(13640,-8935)--(13644,-8908)--(13648,-8881)
  --(13652,-8854)--(13655,-8827)--(13658,-8801)--(13660,-8774)--(13662,-8748)
  --(13663,-8722)--(13663,-8696)--(13663,-8670)--(13662,-8644)--(13660,-8618)
  --(13657,-8592)--(13653,-8566)--(13649,-8540)--(13644,-8513)--(13639,-8487)
  --(13634,-8460)--(13628,-8434)--(13621,-8407)--(13615,-8381)--(13608,-8354)
  --(13601,-8327)--(13594,-8301)--(13586,-8274)--(13579,-8248)--(13571,-8221)
  --(13563,-8196)--(13555,-8170)--(13547,-8144)--(13538,-8119)--(13529,-8095)
  --(13520,-8070)--(13511,-8046)--(13501,-8023)--(13490,-8000)--(13478,-7976)
  --(13465,-7952)--(13451,-7928)--(13437,-7905)--(13422,-7882)--(13407,-7859)
  --(13391,-7836)--(13375,-7814)--(13358,-7792)--(13341,-7769)--(13324,-7747)
  --(13306,-7725)--(13289,-7704)--(13271,-7682)--(13254,-7661)--(13236,-7640)
  --(13219,-7619)--(13201,-7599)--(13184,-7579)--(13166,-7559)--(13149,-7540)
  --(13131,-7522)--(13114,-7504)--(13096,-7487)--(13076,-7469)--(13057,-7453)
  --(13037,-7437)--(13017,-7421)--(12996,-7406)--(12976,-7392)--(12956,-7379)
  --(12935,-7365)--(12914,-7352)--(12894,-7339)--(12873,-7327)--(12852,-7315)
  --(12831,-7303)--(12810,-7291)--(12789,-7279)--(12768,-7267)--(12746,-7256)
  --(12725,-7244)--(12703,-7234)--(12681,-7223)--(12659,-7213)--(12637,-7203)
  --(12614,-7194)--(12591,-7185)--(12568,-7177)--(12544,-7168)--(12520,-7160)
  --(12495,-7152)--(12470,-7144)--(12445,-7135)--(12419,-7127)--(12394,-7119)
  --(12368,-7111)--(12342,-7103)--(12317,-7096)--(12292,-7089)--(12267,-7083)
  --(12242,-7078)--(12218,-7074)--(12195,-7071)--(12172,-7069)--(12150,-7069)
  --(12129,-7071)--(12108,-7074)--(12086,-7080)--(12064,-7088)--(12042,-7099)
  --(12020,-7113)--(11996,-7130)--(11972,-7150)--(11947,-7172)--(11921,-7197)
  --(11895,-7222)--(11869,-7249)--(11843,-7275)--(11820,-7300)--(11799,-7323)
  --(11781,-7343)--(11766,-7359)--(11756,-7371)--(11749,-7379)--(11745,-7383)
  --(11743,-7385);
\pgfsetfillcolor{black}
\pgftext[base,left,at=\pgfqpointxy{2025}{-2100}] {\fontsize{24}{28.8}\normalfont $S_{0}$}
\pgftext[base,left,at=\pgfqpointxy{9450}{-2325}] {\fontsize{24}{28.8}\normalfont $G=(V,E)$}
\pgftext[base,left,at=\pgfqpointxy{10875}{-3900}] {\fontsize{24}{28.8}\normalfont $v$}
\pgftext[base,left,at=\pgfqpointxy{9300}{-3900}] {\fontsize{24}{28.8}\normalfont $u$}
\pgfsetlinewidth{+30\XFigu}
\pgfsetdash{}{+0pt}
\draw (9600,-9150) arc[start angle=+150.26, end angle=+209.74, radius=+604.7];
\draw (11804,-9184) arc[start angle=+30.13, end angle=+-25.73, radius=+639.8];
\pgfsetlinewidth{+45\XFigu}
\pgfsetdash{}{+0pt}
\draw (3150,-7350) arc[start angle=+-113.84, end angle=+-66.16, radius=+1762.9];
\pgfsetlinewidth{+30\XFigu}
\pgfsetdash{}{+0pt}
\draw (3150,-7350) arc[start angle=+150.26, end angle=+209.74, radius=+604.7];
\pgfsetlinewidth{+45\XFigu}
\pgfsetdash{}{+0pt}
\draw (3150,-7950) arc[start angle=+113.84, end angle=+66.16, radius=+1762.9];
\pgfsetlinewidth{+30\XFigu}
\pgfsetdash{}{+0pt}
\draw (4575,-7350) arc[start angle=+28.07, end angle=+-28.07, radius=+637.5];
\pgfsetlinewidth{+45\XFigu}
\pgfsetdash{}{+0pt}
\draw (2721,-12319) arc[start angle=+-113.86, end angle=+-66.14, radius=+2886.4];
\pgfsetlinewidth{+30\XFigu}
\pgfsetdash{}{+0pt}
\draw (2721,-12319) arc[start angle=+150.24, end angle=+209.76, radius=+990.1];
\pgfsetlinewidth{+45\XFigu}
\pgfsetdash{}{+0pt}
\draw (2721,-13302) arc[start angle=+113.86, end angle=+66.14, radius=+2886.4];
\pgftext[base,left,at=\pgfqpointxy{11250}{-8325}] {\fontsize{24}{28.8}\normalfont $\Delta^e_v$}
\pgftext[base,left,at=\pgfqpointxy{13050}{-7425}] {\fontsize{24}{28.8}\normalfont $Y_{v,s}$}
\pgfsetdash{}{+0pt}
\draw (9600,-9150) arc[start angle=+-105.26, end angle=+-74.74, radius=+4275.7];
\draw (9600,-9750) arc[start angle=+105.26, end angle=+74.74, radius=+4275.7];
\pgfsetarrows{[line width=7.5\XFigu]}
\pgfsetarrowsend{xfiga0}
\pgfsetlinewidth{+7.5\XFigu}
\pgfsetdash{{+60\XFigu}{+60\XFigu}}{++0pt}
\draw (3825,-7350) arc[start angle=+132.15, end angle=+44.12, radius=+4967.9];
\draw (9975,-6975) arc[start angle=+26.6, end angle=+-63.4, radius=+335.4];
\draw (3750,-9675) arc[start angle=+-131.75, end angle=+-49.46, radius=+5414.5];
\pgfsetlinewidth{+30\XFigu}
\pgfsetstrokecolor{red}
\pgfsetdash{}{+0pt}
\pgfsetarrowsend{}
\draw (9225,-12750)--(11550,-12750);
\pgftext[base,left,at=\pgfqpointxy{6000}{-7500}] {\fontsize{24}{28.8}\normalfont $Y_{v,t}$}
\pgftext[base,left,at=\pgfqpointxy{10425}{-8025}] {\fontsize{24}{28.8}\normalfont $I_e$}
\pgfsetdash{}{+0pt}
\pgfsetstrokecolor{black}
\draw (11804,-9184) arc[start angle=+154.61, end angle=+209.79, radius=+647.2];
\pgftext[base,left,at=\pgfqpointxy{3300}{-7275}] {\fontsize{24}{28.8}\normalfont $W_{e,t}$}
\draw (8880,-12278) arc[start angle=+27.97, end angle=+-27.97, radius=+1087.5];
\draw (11923,-12337) arc[start angle=+153.51, end angle=+210.66, radius=+1064.9];
\pgftext[base,left,at=\pgfqpointxy{8100}{-7350}] {\fontsize{24}{28.8}\normalfont $Y_{u,s}$}
\pgfsetlinewidth{+15\XFigu}
\pgfsetdash{}{+0pt}
\pgfsetstrokecolor{blue}
\draw (11775,-12375) arc[start angle=+154.83, end angle=+214.70, radius=+905];
\pgfsetlinewidth{+45\XFigu}
\pgfsetdash{}{+0pt}
\pgfsetstrokecolor{black}
\draw (11977,-12278) arc[start angle=+104.1, end angle=+255.9, radius=+525.9];
\pgfsetlinewidth{+15\XFigu}
\pgfsetdash{}{+0pt}
\pgfsetstrokecolor{blue}
\draw (9075,-12375) arc[start angle=+20.77, end angle=+-20.77, radius=+1163.1];
\pgfsetlinewidth{+45\XFigu}
\pgfsetdash{}{+0pt}
\pgfsetstrokecolor{black}
\draw (8880,-12278) arc[start angle=+73.7, end angle=+-73.7, radius=+531.4];
\draw (9600,-7350) arc[start angle=+73.7, end angle=+-73.7, radius=+312.5];
\draw (11775,-7350) arc[start angle=+104.0, end angle=+256.0, radius=+309.2];
\pgfsetlinewidth{+30\XFigu}
\draw (11743,-7385) arc[start angle=+153.60, end angle=+210.60, radius=+628.1];
\pgfsetarrowsend{xfiga0}
\pgfsetlinewidth{+7.5\XFigu}
\pgfsetdash{{+60\XFigu}{+60\XFigu}}{++0pt}
\draw (9975,-6975) arc[start angle=+98.06, end angle=+43.36, radius=+1729.6];
\pgftext[base,left,at=\pgfqpointxy{1500}{-7275}] {\fontsize{24}{28.8}\normalfont $Y_{u,t}$}
\pgftext[base,left,at=\pgfqpointxy{3375}{-9075}] {\fontsize{24}{28.8}\normalfont $W_{e',t}$}
\pgfsetlinewidth{+30\XFigu}
\pgfsetstrokecolor{red}
\pgfsetdash{}{+0pt}
\pgfsetarrowsend{}
\draw (9825,-7650)--(11550,-7650);
\pgfsetlinewidth{+15\XFigu}
\pgfsetdash{}{+0pt}
\pgfsetstrokecolor{blue}
\draw (3000,-12450) arc[start angle=+149.68, end angle=+221.74, radius=+640.7];
\draw (4725,-12450) arc[start angle=+23.50, end angle=+-23.50, radius=+940.5];
\pgfsetarrowsend{xfiga0}
\pgfsetlinewidth{+7.5\XFigu}
\pgfsetdash{{+60\XFigu}{+60\XFigu}}{++0pt}
\pgfsetstrokecolor{black}
\draw (4050,-13350) arc[start angle=+-117.362, end angle=+-59.737, radius=+6149];
\pgfsetarrowsend{}
\pgfsetlinewidth{+30\XFigu}
\pgfsetdash{}{+0pt}
\draw (9600,-7350) arc[start angle=+28.07, end angle=+-28.07, radius=+637.5];
\pgftext[base,left,at=\pgfqpointxy{9525}{-8400}] {\fontsize{24}{28.8}\normalfont $\Delta^e_{u}$}
\pgftext[base,left,at=\pgfqpointxy{10275}{-9225}] {\fontsize{24}{28.8}\normalfont $W_{e',s}$}
\draw (9600,-9150) arc[start angle=+29.74, end angle=+-29.74, radius=+604.7];
\endtikzpicture}%

%% file: canonical_measures_2.tikz
{\pgfkeys{/pgf/fpu/.try=false}%
\ifx\XFigwidth\undefined\dimen1=0pt\else\dimen1\XFigwidth\fi
\divide\dimen1 by 13155
\ifx\XFigheight\undefined\dimen3=0pt\else\dimen3\XFigheight\fi
\divide\dimen3 by 5519
\ifdim\dimen1=0pt\ifdim\dimen3=0pt\dimen1=3946sp\dimen3\dimen1
  \else\dimen1\dimen3\fi\else\ifdim\dimen3=0pt\dimen3\dimen1\fi\fi
\tikzpicture[x=+\dimen1, y=+\dimen3]
{\ifx\XFigu\undefined\catcode`\@11
\def\temp{\alloc@1\dimen\dimendef\insc@unt}\temp\XFigu\catcode`\@12\fi}
\XFigu3946sp
\ifdim\XFigu<0pt\XFigu-\XFigu\fi
\catcode`\@11
\pgfutil@ifundefined{pgf@pattern@name@xfigp0}{
\pgfdeclarepatternformonly{xfigp0}
{\pgfqpoint{-1bp}{-1bp}}{\pgfqpoint{9bp}{5bp}}{\pgfqpoint{8bp}{4bp}}
{	\pgfsetdash{}{0pt}\pgfsetlinewidth{0.45bp}
	\pgfpathqmoveto{-1bp}{4.5bp}\pgfpathqlineto{9bp}{-0.5bp}
	\pgfusepathqstroke
}
}{}
\catcode`\@12
\pgfdeclarearrow{
  name = xfiga0,
  parameters = {
    \the\pgfarrowlinewidth \the\pgfarrowlength \the\pgfarrowwidth},
  defaults = {
	  line width=+7.5\XFigu, length=+120\XFigu, width=+60\XFigu},
  setup code = {
    \dimen7 2.15\pgfarrowlength\pgfmathveclen{\the\dimen7}{\the\pgfarrowwidth}
    \dimen7 2\pgfarrowwidth\pgfmathdivide{\pgfmathresult}{\the\dimen7}
    \dimen7 \pgfmathresult\pgfarrowlinewidth
    \pgfarrowssettipend{+\dimen7}
    \pgfarrowssetbackend{+-\pgfarrowlength}
    \dimen9 -0.5\pgfarrowlinewidth
    \pgfarrowssetvisualbackend{+\dimen9}
    \pgfarrowssetlineend{+-0.5\pgfarrowlinewidth}
    \pgfarrowshullpoint{+\dimen7}{+0pt}
    \pgfarrowsupperhullpoint{+-\pgfarrowlength}{+0.5\pgfarrowwidth}
    \pgfarrowssavethe\pgfarrowlinewidth
    \pgfarrowssavethe\pgfarrowlength
    \pgfarrowssavethe\pgfarrowwidth
  },
  drawing code = {\pgfsetdash{}{+0pt}
    \ifdim\pgfarrowlinewidth=\pgflinewidth\else\pgfsetlinewidth{+\pgfarrowlinewidth}\fi
    \pgfpathmoveto{\pgfqpoint{-\pgfarrowlength}{0.5\pgfarrowwidth}}
    \pgfpathlineto{\pgfqpoint{0pt}{0pt}}
    \pgfpathlineto{\pgfqpoint{-\pgfarrowlength}{-0.5\pgfarrowwidth}}
    \pgfusepathqstroke
  }
}
\definecolor{blue3}{rgb}{0,0,0.82}
\definecolor{red3}{rgb}{0.82,0,0}
\clip(5235,-6097) rectangle (18390,-578);
\tikzset{inner sep=+0pt, outer sep=+0pt}
\pgfsetfillcolor{black}
\pgftext[base,left,at=\pgfqpointxy{5250}{-1200}] {\fontsize{32}{38.4}\normalfont $\curve$}
\pgfsetbeveljoin
\pgfsetlinewidth{+45\XFigu}
\pgfsetdash{}{+0pt}
\pgfsetstrokecolor{black}
\draw (7289,-1327)--(7316,-1342)--(7344,-1358)--(7373,-1375)--(7402,-1393)--(7431,-1412)
  --(7461,-1432)--(7492,-1452)--(7523,-1474)--(7555,-1497)--(7587,-1521)--(7620,-1545)
  --(7653,-1571)--(7686,-1597)--(7719,-1624)--(7752,-1652)--(7785,-1680)--(7817,-1708)
  --(7849,-1737)--(7881,-1767)--(7912,-1796)--(7942,-1825)--(7971,-1855)--(7999,-1884)
  --(8026,-1914)--(8052,-1943)--(8077,-1971)--(8101,-2000)--(8123,-2028)--(8144,-2055)
  --(8164,-2083)--(8183,-2109)--(8200,-2136)--(8217,-2162)--(8232,-2188)--(8247,-2217)
  --(8262,-2246)--(8275,-2275)--(8286,-2304)--(8297,-2334)--(8306,-2363)--(8315,-2393)
  --(8322,-2424)--(8328,-2455)--(8333,-2486)--(8336,-2517)--(8339,-2549)--(8341,-2580)
  --(8342,-2612)--(8342,-2644)--(8341,-2676)--(8339,-2708)--(8337,-2740)--(8333,-2771)
  --(8330,-2802)--(8325,-2833)--(8321,-2864)--(8316,-2894)--(8310,-2924)--(8304,-2953)
  --(8298,-2982)--(8292,-3011)--(8286,-3039)--(8280,-3068)--(8274,-3097)--(8267,-3127)
  --(8261,-3159)--(8254,-3190)--(8248,-3223)--(8241,-3256)--(8235,-3289)--(8229,-3323)
  --(8223,-3358)--(8217,-3393)--(8211,-3429)--(8206,-3465)--(8201,-3501)--(8197,-3537)
  --(8193,-3572)--(8190,-3608)--(8187,-3643)--(8185,-3677)--(8183,-3711)--(8183,-3744)
  --(8183,-3776)--(8183,-3806)--(8185,-3836)--(8187,-3864)--(8190,-3892)--(8194,-3918)
  --(8199,-3943)--(8204,-3967)--(8211,-3991)--(8220,-4019)--(8231,-4047)--(8243,-4074)
  --(8257,-4100)--(8272,-4125)--(8288,-4150)--(8305,-4175)--(8322,-4199)--(8341,-4222)
  --(8359,-4245)--(8378,-4268)--(8397,-4289)--(8415,-4310)--(8432,-4330)--(8449,-4349)
  --(8464,-4368)--(8478,-4386)--(8491,-4403)--(8502,-4419)--(8511,-4435)--(8519,-4451)
  --(8525,-4467)--(8529,-4483)--(8532,-4499)--(8532,-4515)--(8531,-4531)--(8529,-4548)
  --(8525,-4564)--(8519,-4581)--(8511,-4598)--(8503,-4614)--(8492,-4630)--(8481,-4646)
  --(8468,-4661)--(8455,-4675)--(8441,-4688)--(8426,-4701)--(8410,-4712)--(8394,-4722)
  --(8379,-4731)--(8362,-4738)--(8346,-4745)--(8330,-4750)--(8314,-4754)--(8296,-4757)
  --(8278,-4759)--(8260,-4759)--(8241,-4758)--(8222,-4756)--(8203,-4752)--(8183,-4747)
  --(8163,-4741)--(8144,-4733)--(8124,-4725)--(8105,-4715)--(8087,-4704)--(8069,-4693)
  --(8052,-4681)--(8036,-4669)--(8021,-4656)--(8007,-4643)--(7995,-4629)--(7983,-4616)
  --(7972,-4602)--(7964,-4592)--(7957,-4581)--(7950,-4570)--(7943,-4559)--(7936,-4546)
  --(7929,-4533)--(7922,-4518)--(7915,-4502)--(7907,-4484)--(7899,-4465)--(7891,-4445)
  --(7882,-4422)--(7872,-4398)--(7862,-4371)--(7851,-4343)--(7839,-4313)--(7827,-4281)
  --(7814,-4247)--(7800,-4210)--(7785,-4172)--(7770,-4132)--(7753,-4091)--(7736,-4047)
  --(7717,-4001)--(7698,-3952)--(7678,-3902)--(7656,-3849)--(7633,-3793)--(7618,-3758)
  --(7603,-3721)--(7587,-3684)--(7570,-3644)--(7553,-3604)--(7536,-3563)--(7518,-3520)
  --(7499,-3476)--(7480,-3430)--(7460,-3384)--(7440,-3336)--(7419,-3287)--(7398,-3237)
  --(7377,-3187)--(7355,-3135)--(7333,-3082)--(7311,-3029)--(7289,-2975)--(7266,-2921)
  --(7244,-2866)--(7221,-2811)--(7199,-2755)--(7176,-2700)--(7154,-2645)--(7132,-2590)
  --(7111,-2535)--(7090,-2481)--(7069,-2428)--(7049,-2375)--(7029,-2323)--(7010,-2272)
  --(6992,-2221)--(6974,-2172)--(6957,-2124)--(6941,-2078)--(6925,-2032)--(6910,-1988)
  --(6896,-1946)--(6883,-1904)--(6871,-1864)--(6860,-1826)--(6849,-1789)--(6839,-1753)
  --(6831,-1719)--(6823,-1686)--(6816,-1654)--(6809,-1621)--(6803,-1590)--(6799,-1559)
  --(6795,-1531)--(6793,-1503)--(6791,-1477)--(6791,-1452)--(6792,-1429)--(6793,-1406)
  --(6796,-1385)--(6800,-1365)--(6804,-1347)--(6810,-1329)--(6817,-1313)--(6825,-1299)
  --(6833,-1285)--(6843,-1273)--(6854,-1262)--(6865,-1253)--(6878,-1244)--(6891,-1237)
  --(6905,-1232)--(6920,-1227)--(6935,-1224)--(6951,-1222)--(6968,-1221)--(6986,-1221)
  --(7004,-1222)--(7022,-1224)--(7041,-1228)--(7060,-1232)--(7080,-1237)--(7099,-1243)
  --(7120,-1249)--(7140,-1257)--(7161,-1265)--(7181,-1274)--(7202,-1283)--(7224,-1293)
  --(7245,-1304)--(7267,-1315)--cycle;
\pgfsetlinewidth{+30\XFigu}
\pgfsetdash{}{+0pt}
\draw (7688,-2930) arc[start angle=+158.4, end angle=+40.3, radius=+208.6];
\draw (7040,-1870) arc[start angle=+-141.6, end angle=+-38.4, radius=+487.8];
\draw (7216,-1988) arc[start angle=+110.38, end angle=+82.30, radius=+1101.4];
\pgfsetlinewidth{+45\XFigu}
\pgfsetdash{}{+0pt}
\pgfsetstrokecolor{red3}
\draw (8395,-4521) arc[start angle=+-112.65, end angle=+-67.35, radius=+3060.5];
\pgfsetstrokecolor{black}
\draw  (10988,-4285) ellipse [x radius=+1355,y radius=+471];
\pgfsetlinewidth{+30\XFigu}
\pgfsetdash{}{+0pt}
\draw (10340,-4227) arc[start angle=+-109.80, end angle=+-70.20, radius=+2000.5];
\draw (10575,-4285) arc[start angle=+119.77, end angle=+60.23, radius=+890.3];
\pgfsetlinewidth{+45\XFigu}
\pgfsetdash{}{+0pt}
\pgfsetstrokecolor{blue3}
\draw (11930,-2872) arc[start angle=+36.94, end angle=+-36.94, radius=+1175.5];
\pgfsetstrokecolor{black}
\draw (12313,-2940)--(12289,-2956)--(12264,-2971)--(12240,-2984)--(12214,-2997)--(12189,-3009)
  --(12162,-3019)--(12136,-3029)--(12109,-3037)--(12082,-3044)--(12055,-3050)
  --(12028,-3054)--(12000,-3057)--(11973,-3059)--(11947,-3059)--(11921,-3058)
  --(11895,-3056)--(11870,-3052)--(11846,-3047)--(11823,-3041)--(11801,-3033)
  --(11780,-3024)--(11760,-3014)--(11742,-3003)--(11724,-2991)--(11708,-2977)
  --(11693,-2963)--(11679,-2947)--(11666,-2930)--(11653,-2911)--(11641,-2890)
  --(11629,-2868)--(11619,-2844)--(11608,-2820)--(11598,-2793)--(11588,-2766)
  --(11578,-2738)--(11568,-2709)--(11557,-2680)--(11546,-2650)--(11534,-2620)
  --(11522,-2591)--(11508,-2561)--(11494,-2533)--(11478,-2506)--(11462,-2479)
  --(11444,-2454)--(11424,-2431)--(11404,-2410)--(11381,-2390)--(11357,-2372)
  --(11332,-2356)--(11304,-2343)--(11275,-2331)--(11243,-2322)--(11216,-2316)
  --(11188,-2312)--(11159,-2308)--(11128,-2306)--(11095,-2304)--(11061,-2304)
  --(11026,-2304)--(10988,-2306)--(10950,-2307)--(10910,-2310)--(10868,-2313)
  --(10826,-2317)--(10783,-2321)--(10739,-2326)--(10694,-2330)--(10649,-2335)
  --(10604,-2340)--(10559,-2345)--(10514,-2350)--(10469,-2354)--(10425,-2359)
  --(10382,-2363)--(10340,-2366)--(10299,-2369)--(10259,-2371)--(10220,-2373)
  --(10183,-2373)--(10147,-2373)--(10113,-2372)--(10080,-2370)--(10049,-2367)
  --(10020,-2363)--(9992,-2358)--(9966,-2351)--(9946,-2345)--(9928,-2339)--(9911,-2331)
  --(9894,-2323)--(9879,-2314)--(9864,-2305)--(9850,-2294)--(9838,-2283)--(9826,-2271)
  --(9815,-2259)--(9806,-2246)--(9797,-2233)--(9790,-2219)--(9784,-2204)--(9779,-2189)
  --(9775,-2174)--(9773,-2158)--(9772,-2142)--(9773,-2126)--(9774,-2110)--(9778,-2093)
  --(9782,-2076)--(9788,-2060)--(9796,-2043)--(9805,-2027)--(9815,-2011)--(9827,-1995)
  --(9840,-1979)--(9854,-1964)--(9870,-1949)--(9888,-1934)--(9906,-1920)--(9926,-1907)
  --(9948,-1894)--(9970,-1882)--(9994,-1870)--(10019,-1859)--(10045,-1848)--(10073,-1839)
  --(10102,-1830)--(10132,-1822)--(10164,-1814)--(10196,-1808)--(10231,-1802)
  --(10265,-1797)--(10300,-1793)--(10337,-1789)--(10376,-1786)--(10416,-1784)
  --(10457,-1783)--(10499,-1782)--(10543,-1782)--(10589,-1782)--(10635,-1784)
  --(10683,-1786)--(10732,-1788)--(10783,-1792)--(10834,-1796)--(10886,-1800)
  --(10939,-1806)--(10993,-1812)--(11048,-1818)--(11102,-1825)--(11158,-1833)
  --(11213,-1842)--(11269,-1851)--(11325,-1860)--(11380,-1871)--(11435,-1881)
  --(11490,-1892)--(11544,-1904)--(11597,-1916)--(11650,-1928)--(11702,-1941)
  --(11752,-1954)--(11801,-1967)--(11849,-1981)--(11896,-1995)--(11941,-2009)
  --(11985,-2023)--(12027,-2038)--(12068,-2052)--(12107,-2067)--(12145,-2082)
  --(12180,-2097)--(12214,-2113)--(12247,-2128)--(12277,-2143)--(12306,-2159)
  --(12333,-2175)--(12364,-2194)--(12393,-2214)--(12419,-2234)--(12443,-2254)
  --(12465,-2274)--(12485,-2295)--(12503,-2317)--(12519,-2339)--(12533,-2361)
  --(12545,-2383)--(12555,-2406)--(12564,-2429)--(12570,-2452)--(12575,-2476)
  --(12577,-2499)--(12578,-2523)--(12578,-2546)--(12576,-2570)--(12572,-2593)
  --(12566,-2616)--(12560,-2639)--(12551,-2661)--(12542,-2683)--(12531,-2705)
  --(12520,-2726)--(12507,-2746)--(12494,-2766)--(12480,-2786)--(12465,-2804)
  --(12449,-2822)--(12433,-2839)--(12417,-2856)--(12400,-2872)--(12383,-2887)
  --(12366,-2901)--(12348,-2915)--(12331,-2928)--cycle;
\pgfsetlinewidth{+30\XFigu}
\pgfsetdash{}{+0pt}
\draw (11695,-2400) arc[start angle=+-143.6, end angle=+-46.8, radius=+434.9];
\draw (11812,-2518) arc[start angle=+150.3, end angle=+29.7, radius=+237.8];
\draw (10575,-2047) arc[start angle=+-108.45, end angle=+-71.55, radius=+1210.2];
\pgfsetarrows{[line width=7.5\XFigu, width=48\XFigu, length=145\XFigu]}
\pgfsetarrowsend{xfiga0}
\pgfsetlinewidth{+45\XFigu}
\pgfsetdash{}{+0pt}
\pgfsetstrokecolor{red3}
\draw (15528,-2424) arc[start angle=+143.08, end angle=+36.92, radius=+606];
\draw (15528,-2424) arc[start angle=+-143.08, end angle=+-36.92, radius=+606];
\pgfsetlinewidth{+15\XFigu}
\pgfsetdash{}{+0pt}
\pgfsetstrokecolor{black}
\pgfsetfillpattern{xfigp0}{black}
\draw[pattern,preaction={fill=black}]  (15528,-2424) circle [radius=+61];
\draw[pattern,preaction={fill=black}]  (16497,-2424) circle [radius=+61];
\draw[pattern,preaction={fill=black}]  (15528,-4119) circle [radius=+61];
\draw[pattern,preaction={fill=black}]  (16497,-4604) circle [radius=+61];
\pgfsetlinewidth{+45\XFigu}
\pgfsetdash{}{+0pt}
\pgfsetstrokecolor{blue3}
\draw (16497,-3635) arc[start angle=+53.08, end angle=+-53.08, radius=+606];
\pgfsetlinewidth{+15\XFigu}
\pgfsetdash{}{+0pt}
\pgfsetstrokecolor{black}
\draw[pattern,preaction={fill=black}]  (16497,-3635) circle [radius=+61];
\pgfsetdash{{+90\XFigu}{+90\XFigu}}{++0pt}
\pgfsetarrowsend{}
\draw (15286,-3150)--(16739,-3150);
\pgfsetarrows{[width=60\XFigu, length=180\XFigu]}
\pgfsetarrowsend{xfiga0}
\draw (9450,-3225) arc[start angle=+-6.08, end angle=+27.32, radius=+2124.4];
\draw (9450,-3225) arc[start angle=+15.26, end angle=+-52.13, radius=+1282.7];
\pgfsetarrowsend{}
\pgfsetlinewidth{+45\XFigu}
\pgfsetdash{}{+0pt}
\pgfsetstrokecolor{red3}
\draw (7924,-2165) arc[start angle=+111.82, end angle=+71.12, radius=+3304];
\pgfsetlinewidth{+30\XFigu}
\pgfsetdash{}{+0pt}
\pgfsetstrokecolor{black}
\draw (10752,-2047) arc[start angle=+122.6, end angle=+57.4, radius=+382.7];
\pgfsetlinewidth{+15\XFigu}
\pgfsetdash{{+90\XFigu}{+90\XFigu}}{++0pt}
\draw (14100,-600)--(14100,-5100);
\pgfsetmiterjoin
\pgfsetarrowsend{xfiga0}
\draw (12600,-4050)--(12450,-3825)--(12225,-3825);
\pgfsetfillcolor{black}
\pgftext[base,left,at=\pgfqpointxy{12225}{-3525}] {\fontsize{27}{32.4}\normalfont $e_3$}
\pgfsetfillcolor{blue3}
\pgftext[base,left,at=\pgfqpointxy{12525}{-4125}] {\fontsize{27}{32.4}\normalfont $\pi_2$}
\pgfsetfillcolor{black}
\pgftext[base,left,at=\pgfqpointxy{17025}{-2550}] {\fontsize{27}{32.4}\normalfont $\Gamma^1$}
\pgftext[base,left,at=\pgfqpointxy{17100}{-4125}] {\fontsize{27}{32.4}\normalfont $\Gamma^2$}
\pgftext[base,left,at=\pgfqpointxy{12075}{-4950}] {\fontsize{27}{32.4}\normalfont $C_w$}
\pgftext[base,left,at=\pgfqpointxy{8925}{-5100}] {\fontsize{27}{32.4}\normalfont $e_2$}
\pgftext[base,left,at=\pgfqpointxy{8550}{-1875}] {\fontsize{27}{32.4}\normalfont $e_1$}
\pgftext[base,left,at=\pgfqpointxy{18375}{-3075}] {\fontsize{27}{32.4}\normalfont ${\color{red}\ell_1}(e_1)={\color{red}\ell_1}(e_2)=1/2$}
\pgftext[base,left,at=\pgfqpointxy{18375}{-3525}] {\fontsize{27}{32.4}\normalfont ${\color{blue}\ell_2}(e_3)=1$}
\pgfsetfillcolor{red3}
\pgftext[base,left,at=\pgfqpointxy{9525}{-3300}] {\fontsize{27}{32.4}\normalfont $\pi_1$}
\pgfsetfillcolor{black}
\pgftext[base,left,at=\pgfqpointxy{11925}{-1875}] {\fontsize{27}{32.4}\normalfont $C_v$}
\pgftext[base,left,at=\pgfqpointxy{6900}{-1050}] {\fontsize{27}{32.4}\normalfont $C_u$}
\pgftext[base,left,at=\pgfqpointxy{9900}{-5925}] {\fontsize{32}{38.4}\normalfont $\mu_{\hcurve} = \frac16 \Bigl(\mu_{C_v} + \mu_{C_u}+\mu_{C_w} + d\theta_{e_1}+d\theta_{e_2}\Bigr)$}
\pgfsetarrowsend{}
\pgfsetlinewidth{+30\XFigu}
\pgfsetdash{}{+0pt}
\draw (7452,-2872) arc[start angle=+-110.67, end angle=+-60.52, radius=+906.5];
\endtikzpicture}%

%% file: tree-vc.pdf_t
\begin{picture}(0,0)%
\includegraphics{tree-vc.pdf}%
\end{picture}%
\setlength{\unitlength}{3947sp}%
\begingroup\makeatletter\ifx\SetFigFont\undefined%
\gdef\SetFigFont#1#2#3#4#5{%
  \reset@font\fontsize{#1}{#2pt}%
  \fontfamily{#3}\fontseries{#4}\fontshape{#5}%
  \selectfont}%
\fi\endgroup%
\begin{picture}(16505,7435)(-608,-7286)
\put(3826,-2911){\makebox(0,0)[lb]{\smash{{\SetFigFont{29}{34.8}{\familydefault}{\mddefault}{\updefault}{\color[rgb]{0,0,0}$a_e$}%
}}}}
\put(3301,-5011){\makebox(0,0)[lb]{\smash{{\SetFigFont{29}{34.8}{\familydefault}{\mddefault}{\updefault}{\color[rgb]{0,0,0}$a_1,a_2$}%
}}}}
\put(11326,-5086){\makebox(0,0)[lb]{\smash{{\SetFigFont{29}{34.8}{\familydefault}{\mddefault}{\updefault}{\color[rgb]{0,0,0}$e_1$}%
}}}}
\put(12901,-3286){\makebox(0,0)[lb]{\smash{{\SetFigFont{29}{34.8}{\familydefault}{\mddefault}{\updefault}{\color[rgb]{0,0,0}$e$}%
}}}}
\put(4951,-586){\makebox(0,0)[lb]{\smash{{\SetFigFont{29}{34.8}{\familydefault}{\mddefault}{\updefault}{\color[rgb]{0,0,0}$p$}%
}}}}
\put(3751,-286){\makebox(0,0)[lb]{\smash{{\SetFigFont{29}{34.8}{\familydefault}{\mddefault}{\updefault}{\color[rgb]{0,0,0}$\gamma_p$}%
}}}}
\put(12001,-4111){\makebox(0,0)[lb]{\smash{{\SetFigFont{29}{34.8}{\familydefault}{\mddefault}{\updefault}{\color[rgb]{0,0,0}$e_2$}%
}}}}
\put(14551,-3286){\makebox(0,0)[lb]{\smash{{\SetFigFont{29}{34.8}{\familydefault}{\mddefault}{\updefault}{\color[rgb]{0,0,0}$C_{v_4}$}%
}}}}
\put(11851,-2536){\makebox(0,0)[lb]{\smash{{\SetFigFont{29}{34.8}{\familydefault}{\mddefault}{\updefault}{\color[rgb]{0,0,0}$v_1$}%
}}}}
\put(10951,-3811){\makebox(0,0)[lb]{\smash{{\SetFigFont{29}{34.8}{\familydefault}{\mddefault}{\updefault}{\color[rgb]{0,0,0}$v_2$}%
}}}}
\put(11701,-5836){\makebox(0,0)[lb]{\smash{{\SetFigFont{29}{34.8}{\familydefault}{\mddefault}{\updefault}{\color[rgb]{0,0,0}$v_3$}%
}}}}
\put(13651,-3811){\makebox(0,0)[lb]{\smash{{\SetFigFont{29}{34.8}{\familydefault}{\mddefault}{\updefault}{\color[rgb]{0,0,0}$v_4$}%
}}}}
\put(10951,-811){\makebox(0,0)[lb]{\smash{{\SetFigFont{29}{34.8}{\familydefault}{\mddefault}{\updefault}{\color[rgb]{0,0,0}$C_{v_1}$}%
}}}}
\put(9151,-3136){\makebox(0,0)[lb]{\smash{{\SetFigFont{29}{34.8}{\familydefault}{\mddefault}{\updefault}{\color[rgb]{0,0,0}$C_{v_2}$}%
}}}}
\put(10876,-6736){\makebox(0,0)[lb]{\smash{{\SetFigFont{29}{34.8}{\familydefault}{\mddefault}{\updefault}{\color[rgb]{0,0,0}$C_{v_3}$}%
}}}}
\end{picture}%

%% file: ProofPicSphereChain.tikz
{\pgfkeys{/pgf/fpu/.try=false}%
\ifx\XFigwidth\undefined\dimen1=0pt\else\dimen1\XFigwidth\fi
\divide\dimen1 by 6016
\ifx\XFigheight\undefined\dimen3=0pt\else\dimen3\XFigheight\fi
\divide\dimen3 by 6162
\ifdim\dimen1=0pt\ifdim\dimen3=0pt\dimen1=3946sp\dimen3\dimen1
  \else\dimen1\dimen3\fi\else\ifdim\dimen3=0pt\dimen3\dimen1\fi\fi
\tikzpicture[x=+\dimen1, y=+\dimen3]
{\ifx\XFigu\undefined\catcode`\@11
\def\temp{\alloc@1\dimen\dimendef\insc@unt}\temp\XFigu\catcode`\@12\fi}
\XFigu3946sp
\ifdim\XFigu<0pt\XFigu-\XFigu\fi
\catcode`\@11
\pgfutil@ifundefined{pgf@pattern@name@xfigp0}{
\pgfdeclarepatternformonly{xfigp0}
{\pgfqpoint{-1bp}{-1bp}}{\pgfqpoint{9bp}{5bp}}{\pgfqpoint{8bp}{4bp}}
{	\pgfsetdash{}{0pt}\pgfsetlinewidth{0.45bp}
	\pgfpathqmoveto{-1bp}{4.5bp}\pgfpathqlineto{9bp}{-0.5bp}
	\pgfusepathqstroke
}
}{}
\catcode`\@12
\pgfdeclarearrow{
  name = xfiga0,
  parameters = {
    \the\pgfarrowlinewidth \the\pgfarrowlength \the\pgfarrowwidth},
  defaults = {
	  line width=+7.5\XFigu, length=+120\XFigu, width=+60\XFigu},
  setup code = {
    \dimen7 2.15\pgfarrowlength\pgfmathveclen{\the\dimen7}{\the\pgfarrowwidth}
    \dimen7 2\pgfarrowwidth\pgfmathdivide{\pgfmathresult}{\the\dimen7}
    \dimen7 \pgfmathresult\pgfarrowlinewidth
    \pgfarrowssettipend{+\dimen7}
    \pgfarrowssetbackend{+-\pgfarrowlength}
    \dimen9 -0.5\pgfarrowlinewidth
    \pgfarrowssetvisualbackend{+\dimen9}
    \pgfarrowssetlineend{+-0.5\pgfarrowlinewidth}
    \pgfarrowshullpoint{+\dimen7}{+0pt}
    \pgfarrowsupperhullpoint{+-\pgfarrowlength}{+0.5\pgfarrowwidth}
    \pgfarrowssavethe\pgfarrowlinewidth
    \pgfarrowssavethe\pgfarrowlength
    \pgfarrowssavethe\pgfarrowwidth
  },
  drawing code = {\pgfsetdash{}{+0pt}
    \ifdim\pgfarrowlinewidth=\pgflinewidth\else\pgfsetlinewidth{+\pgfarrowlinewidth}\fi
    \pgfpathmoveto{\pgfqpoint{-\pgfarrowlength}{0.5\pgfarrowwidth}}
    \pgfpathlineto{\pgfqpoint{0pt}{0pt}}
    \pgfpathlineto{\pgfqpoint{-\pgfarrowlength}{-0.5\pgfarrowwidth}}
    \pgfusepathqstroke
  }
}
\clip(967,-7773) rectangle (6983,-1611);
\tikzset{inner sep=+0pt, outer sep=+0pt}
\pgfsetfillcolor{black}
\pgftext[base,left,at=\pgfqpointxy{6750}{-6825}] {\fontsize{24}{28.8}\normalfont $l_3$}
\pgftext[base,left,at=\pgfqpointxy{1050}{-7500}] {\fontsize{24}{28.8}\normalfont $u$}
\pgfsetlinewidth{+45\XFigu}
\pgfsetstrokecolor{black}
\draw  (3732,-2653) circle [radius=+546];
\pgfsetbeveljoin
\draw (4395,-4413)--(4381,-4425)--(4367,-4437)--(4352,-4448)--(4336,-4457)--(4318,-4467)
  --(4299,-4475)--(4279,-4483)--(4257,-4490)--(4234,-4496)--(4210,-4501)--(4184,-4506)
  --(4158,-4509)--(4130,-4512)--(4101,-4514)--(4071,-4515)--(4040,-4515)--(4008,-4514)
  --(3975,-4512)--(3942,-4509)--(3909,-4505)--(3874,-4501)--(3840,-4495)--(3805,-4489)
  --(3771,-4482)--(3736,-4475)--(3700,-4467)--(3665,-4458)--(3629,-4448)--(3593,-4438)
  --(3557,-4427)--(3524,-4416)--(3491,-4405)--(3458,-4393)--(3423,-4381)--(3388,-4368)
  --(3352,-4354)--(3316,-4340)--(3279,-4325)--(3241,-4309)--(3203,-4293)--(3165,-4275)
  --(3126,-4258)--(3087,-4239)--(3048,-4220)--(3009,-4201)--(2970,-4181)--(2931,-4161)
  --(2893,-4140)--(2856,-4120)--(2819,-4099)--(2783,-4077)--(2749,-4056)--(2715,-4035)
  --(2682,-4014)--(2651,-3993)--(2621,-3973)--(2593,-3952)--(2566,-3932)--(2540,-3911)
  --(2516,-3891)--(2493,-3872)--(2471,-3852)--(2451,-3833)--(2432,-3814)--(2412,-3792)
  --(2394,-3770)--(2377,-3748)--(2362,-3726)--(2348,-3704)--(2336,-3681)--(2325,-3659)
  --(2315,-3636)--(2307,-3613)--(2301,-3589)--(2296,-3566)--(2292,-3543)--(2291,-3520)
  --(2290,-3496)--(2291,-3473)--(2294,-3451)--(2298,-3428)--(2304,-3406)--(2311,-3385)
  --(2319,-3364)--(2328,-3344)--(2339,-3324)--(2350,-3305)--(2363,-3286)--(2377,-3268)
  --(2391,-3251)--(2407,-3235)--(2423,-3218)--(2440,-3203)--(2458,-3188)--(2478,-3172)
  --(2499,-3157)--(2521,-3142)--(2545,-3127)--(2570,-3113)--(2596,-3099)--(2623,-3085)
  --(2651,-3072)--(2680,-3060)--(2709,-3048)--(2740,-3037)--(2771,-3027)--(2802,-3017)
  --(2834,-3008)--(2866,-3001)--(2898,-2994)--(2929,-2989)--(2960,-2984)--(2990,-2981)
  --(3020,-2978)--(3048,-2977)--(3076,-2977)--(3103,-2978)--(3128,-2981)--(3153,-2984)
  --(3176,-2989)--(3198,-2994)--(3220,-3001)--(3244,-3010)--(3266,-3022)--(3288,-3034)
  --(3310,-3049)--(3330,-3064)--(3351,-3082)--(3370,-3101)--(3390,-3121)--(3410,-3142)
  --(3429,-3165)--(3448,-3188)--(3467,-3212)--(3485,-3236)--(3504,-3261)--(3522,-3285)
  --(3541,-3309)--(3559,-3333)--(3578,-3356)--(3596,-3379)--(3615,-3400)--(3634,-3421)
  --(3654,-3441)--(3674,-3459)--(3695,-3477)--(3715,-3492)--(3737,-3507)--(3759,-3522)
  --(3782,-3536)--(3807,-3549)--(3832,-3563)--(3858,-3576)--(3885,-3589)--(3913,-3602)
  --(3942,-3615)--(3971,-3628)--(4000,-3641)--(4029,-3654)--(4058,-3667)--(4087,-3680)
  --(4115,-3693)--(4143,-3706)--(4170,-3720)--(4196,-3734)--(4221,-3747)--(4244,-3762)
  --(4267,-3776)--(4288,-3791)--(4308,-3806)--(4327,-3822)--(4345,-3839)--(4360,-3855)
  --(4375,-3873)--(4389,-3891)--(4401,-3909)--(4414,-3929)--(4425,-3950)--(4435,-3971)
  --(4444,-3993)--(4453,-4016)--(4460,-4039)--(4466,-4063)--(4471,-4087)--(4475,-4111)
  --(4478,-4135)--(4479,-4159)--(4480,-4183)--(4479,-4207)--(4477,-4230)--(4473,-4252)
  --(4469,-4274)--(4463,-4295)--(4457,-4314)--(4449,-4333)--(4440,-4351)--(4430,-4368)
  --(4419,-4384)--(4407,-4399)--cycle;
\draw (6351,-2831)--(6331,-2813)--(6311,-2796)--(6290,-2780)--(6267,-2763)--(6243,-2747)
  --(6219,-2731)--(6193,-2716)--(6167,-2700)--(6139,-2686)--(6111,-2671)--(6083,-2658)
  --(6054,-2645)--(6025,-2633)--(5996,-2622)--(5967,-2612)--(5939,-2603)--(5911,-2595)
  --(5884,-2589)--(5857,-2583)--(5832,-2578)--(5808,-2575)--(5784,-2572)--(5761,-2571)
  --(5739,-2570)--(5717,-2571)--(5696,-2573)--(5675,-2575)--(5655,-2579)--(5634,-2584)
  --(5615,-2590)--(5595,-2597)--(5576,-2606)--(5558,-2615)--(5540,-2626)--(5522,-2638)
  --(5506,-2650)--(5490,-2664)--(5475,-2679)--(5462,-2694)--(5449,-2710)--(5437,-2727)
  --(5427,-2745)--(5418,-2762)--(5409,-2781)--(5402,-2799)--(5396,-2818)--(5391,-2838)
  --(5386,-2858)--(5383,-2878)--(5380,-2899)--(5378,-2922)--(5377,-2945)--(5377,-2969)
  --(5377,-2994)--(5378,-3020)--(5380,-3048)--(5383,-3076)--(5386,-3104)--(5389,-3134)
  --(5393,-3164)--(5398,-3194)--(5403,-3225)--(5408,-3255)--(5414,-3286)--(5420,-3316)
  --(5426,-3346)--(5431,-3375)--(5437,-3404)--(5443,-3432)--(5449,-3460)--(5455,-3488)
  --(5461,-3515)--(5466,-3542)--(5472,-3570)--(5477,-3597)--(5482,-3625)--(5487,-3654)
  --(5492,-3682)--(5497,-3711)--(5502,-3740)--(5506,-3770)--(5510,-3799)--(5514,-3828)
  --(5517,-3857)--(5520,-3885)--(5522,-3913)--(5524,-3940)--(5525,-3966)--(5525,-3991)
  --(5525,-4014)--(5525,-4037)--(5524,-4058)--(5522,-4079)--(5519,-4098)--(5516,-4116)
  --(5513,-4134)--(5508,-4152)--(5502,-4169)--(5496,-4185)--(5488,-4201)--(5480,-4216)
  --(5471,-4230)--(5461,-4244)--(5450,-4256)--(5438,-4268)--(5426,-4278)--(5413,-4287)
  --(5400,-4296)--(5386,-4303)--(5372,-4308)--(5358,-4313)--(5344,-4316)--(5330,-4318)
  --(5316,-4319)--(5302,-4318)--(5288,-4317)--(5275,-4314)--(5261,-4311)--(5246,-4305)
  --(5231,-4299)--(5216,-4291)--(5200,-4282)--(5184,-4272)--(5167,-4261)--(5150,-4249)
  --(5132,-4237)--(5114,-4224)--(5096,-4210)--(5078,-4196)--(5059,-4183)--(5041,-4169)
  --(5023,-4156)--(5005,-4144)--(4987,-4132)--(4970,-4121)--(4952,-4111)--(4935,-4102)
  --(4918,-4094)--(4900,-4087)--(4882,-4081)--(4864,-4076)--(4845,-4072)--(4825,-4068)
  --(4806,-4066)--(4785,-4065)--(4765,-4065)--(4744,-4065)--(4724,-4067)--(4704,-4070)
  --(4685,-4074)--(4666,-4078)--(4648,-4084)--(4631,-4090)--(4615,-4097)--(4600,-4105)
  --(4585,-4113)--(4572,-4122)--(4560,-4132)--(4549,-4143)--(4538,-4154)--(4528,-4166)
  --(4519,-4179)--(4510,-4193)--(4502,-4208)--(4495,-4224)--(4489,-4241)--(4484,-4259)
  --(4480,-4278)--(4476,-4297)--(4474,-4316)--(4473,-4336)--(4474,-4356)--(4475,-4376)
  --(4477,-4396)--(4480,-4416)--(4484,-4435)--(4490,-4455)--(4496,-4475)--(4502,-4493)
  --(4509,-4511)--(4518,-4530)--(4527,-4549)--(4538,-4569)--(4549,-4589)--(4562,-4610)
  --(4576,-4631)--(4591,-4653)--(4608,-4674)--(4625,-4696)--(4643,-4717)--(4662,-4738)
  --(4682,-4759)--(4703,-4779)--(4724,-4799)--(4745,-4818)--(4767,-4836)--(4790,-4853)
  --(4812,-4870)--(4835,-4886)--(4859,-4902)--(4877,-4913)--(4895,-4924)--(4914,-4934)
  --(4934,-4945)--(4955,-4955)--(4976,-4964)--(4998,-4973)--(5020,-4981)--(5044,-4989)
  --(5068,-4996)--(5093,-5001)--(5119,-5006)--(5146,-5010)--(5173,-5013)--(5200,-5014)
  --(5228,-5014)--(5257,-5012)--(5285,-5009)--(5314,-5005)--(5343,-4999)--(5372,-4991)
  --(5401,-4982)--(5429,-4971)--(5458,-4959)--(5487,-4945)--(5515,-4929)--(5544,-4912)
  --(5572,-4892)--(5601,-4871)--(5629,-4849)--(5652,-4829)--(5675,-4808)--(5699,-4786)
  --(5722,-4763)--(5746,-4738)--(5771,-4711)--(5795,-4684)--(5820,-4655)--(5845,-4625)
  --(5871,-4593)--(5896,-4561)--(5922,-4527)--(5948,-4492)--(5974,-4456)--(6000,-4419)
  --(6026,-4381)--(6052,-4343)--(6078,-4304)--(6103,-4265)--(6128,-4225)--(6152,-4186)
  --(6176,-4146)--(6199,-4106)--(6222,-4067)--(6244,-4027)--(6265,-3989)--(6286,-3950)
  --(6305,-3913)--(6324,-3876)--(6342,-3839)--(6359,-3804)--(6375,-3769)--(6390,-3736)
  --(6405,-3703)--(6418,-3671)--(6430,-3640)--(6442,-3610)--(6453,-3580)--(6465,-3544)
  --(6477,-3509)--(6487,-3475)--(6495,-3441)--(6503,-3408)--(6510,-3376)--(6515,-3344)
  --(6519,-3313)--(6522,-3283)--(6523,-3253)--(6524,-3224)--(6523,-3196)--(6521,-3168)
  --(6518,-3141)--(6514,-3115)--(6509,-3090)--(6502,-3066)--(6495,-3043)--(6487,-3021)
  --(6478,-2999)--(6468,-2979)--(6457,-2959)--(6446,-2941)--(6434,-2923)--(6421,-2906)
  --(6408,-2890)--(6395,-2874)--(6381,-2859)--(6366,-2845)--cycle;
\draw  (4829,-2687) circle [radius=+546];
\pgfsetlinewidth{+30\XFigu}
\draw (4950,-4650) arc[start angle=+135.0, end angle=+45.0, radius=+265.2];
\draw (4800,-4575) arc[start angle=+-134.2, end angle=+-58.5, radius=+553.1];
\draw (5775,-3300) arc[start angle=+126.9, end angle=+53.1, radius=+375];
\draw (5700,-3225) arc[start angle=+-143.1, end angle=+-36.9, radius=+375];
\pgfsetlinewidth{+7.5\XFigu}
\pgfsetfillcolor{cyan}
\filldraw  (1050,-7125) circle [radius=+75];
\pgfsetfillcolor{black}
\pgftext[base,left,at=\pgfqpointxy{1050}{-5925}] {\fontsize{24}{28.8}\normalfont $G=(V,E)$}
\pgftext[base,left,at=\pgfqpointxy{1725}{-6675}] {\fontsize{24}{28.8}\normalfont $l$}
\pgfsetfillcolor{cyan}
\filldraw  (6525,-6750) circle [radius=+75];
\pgfsetlinewidth{+30\XFigu}
\pgfsetdash{}{+0pt}
\draw (5400,-7200)--(5775,-6750);
\pgfsetfillcolor{black}
\pgftext[base,left,at=\pgfqpointxy{5700}{-6600}] {\fontsize{24}{28.8}\normalfont $w_1$}
\pgfsetdash{}{+0pt}
\draw (5400,-7275) arc[start angle=+-133.60, end angle=+-46.40, radius=+1087.5];
\pgftext[base,left,at=\pgfqpointxy{5325}{-6975}] {\fontsize{24}{28.8}\normalfont $l_1$}
\pgfsetlinewidth{+7.5\XFigu}
\pgfsetfillcolor{cyan}
\filldraw  (6900,-7200) circle [radius=+75];
\filldraw  (5400,-7275) circle [radius=+75];
\pgfsetfillcolor{black}
\pgftext[base,left,at=\pgfqpointxy{5250}{-7650}] {\fontsize{24}{28.8}\normalfont $u$}
\pgfsetlinewidth{+30\XFigu}
\draw (1200,-7125) arc[start angle=+-133.60, end angle=+-46.40, radius=+1087.5];
\draw (2625,-3450) arc[start angle=+-138.4, end angle=+-41.6, radius=+451.6];
\draw (2775,-3525) arc[start angle=+130.6, end angle=+49.4, radius=+345.7];
\pgfsetarrows{[line width=7.5\XFigu]}
\pgfsetarrowsend{xfiga0}
\pgfsetlinewidth{+7.5\XFigu}
\pgfsetdash{{+60\XFigu}{+60\XFigu}}{++0pt}
\draw (2475,-6675) arc[start angle=+139.17, end angle=+34.64, radius=+1757];
\pgfsetdash{}{+0pt}
\pgfsetfillcolor{cyan}
\filldraw  (2700,-7125) circle [radius=+75];
\pgfsetlinewidth{+30\XFigu}
\pgfsetfillpattern{xfigp0}{black}
\draw[pattern,preaction={fill=black}]  (3750,-2625) circle [radius=+38];
\draw[pattern,preaction={fill=black}]  (3825,-4200) circle [radius=+38];
\draw[pattern,preaction={fill=black}]  (4800,-2625) circle [radius=+38];
\pgfsetlinewidth{+7.5\XFigu}
\pgfsetfillcolor{cyan}
\filldraw  (5850,-6750) circle [radius=+75];
\pgfsetlinewidth{+30\XFigu}
\pgfsetdash{}{+0pt}
\pgfsetarrowsend{}
\draw (5925,-6750)--(6450,-6750);
\draw (6900,-7125)--(6600,-6750);
\pgfsetfillcolor{black}
\pgftext[base,left,at=\pgfqpointxy{2625}{-7500}] {\fontsize{24}{28.8}\normalfont $v$}
\pgftext[base,left,at=\pgfqpointxy{3675}{-2925}] {\fontsize{18}{21.6}\normalfont $1_{w_1}$}
\pgftext[base,left,at=\pgfqpointxy{4725}{-2925}] {\fontsize{18}{21.6}\normalfont $1_{w_2}$}
\pgftext[base,left,at=\pgfqpointxy{3600}{-1950}] {\fontsize{24}{28.8}\normalfont $\P^1_{w_1}$}
\pgftext[base,left,at=\pgfqpointxy{4725}{-1950}] {\fontsize{24}{28.8}\normalfont $\P^1_{w_2}$}
\pgftext[base,left,at=\pgfqpointxy{1125}{-2250}] {\fontsize{24}{28.8}\normalfont $\widetilde\rsf$}
\pgftext[base,left,at=\pgfqpointxy{6075}{-7050}] {\fontsize{24}{28.8}\normalfont $l_2$}
\pgftext[base,left,at=\pgfqpointxy{6825}{-7575}] {\fontsize{24}{28.8}\normalfont $v$}
\pgftext[base,left,at=\pgfqpointxy{6225}{-6600}] {\fontsize{24}{28.8}\normalfont $w_2$}
\pgftext[base,left,at=\pgfqpointxy{3450}{-4050}] {\fontsize{24}{28.8}\normalfont $p_t$}
\pgftext[base,left,at=\pgfqpointxy{4725}{-5925}] {\fontsize{24}{28.8}\normalfont $\widetilde G=(\widetilde V,\widetilde E)$}
\pgfsetdash{}{+0pt}
\draw (1200,-7125) arc[start angle=+133.60, end angle=+46.40, radius=+1087.5];
\endtikzpicture}%

%% file: ProofPicSubcases.tikz
{\pgfkeys{/pgf/fpu/.try=false}%
\ifx\XFigwidth\undefined\dimen1=0pt\else\dimen1\XFigwidth\fi
\divide\dimen1 by 12826
\ifx\XFigheight\undefined\dimen3=0pt\else\dimen3\XFigheight\fi
\divide\dimen3 by 8031
\ifdim\dimen1=0pt\ifdim\dimen3=0pt\dimen1=3946sp\dimen3\dimen1
  \else\dimen1\dimen3\fi\else\ifdim\dimen3=0pt\dimen3\dimen1\fi\fi
\tikzpicture[x=+\dimen1, y=+\dimen3]
{\ifx\XFigu\undefined\catcode`\@11
\def\temp{\alloc@1\dimen\dimendef\insc@unt}\temp\XFigu\catcode`\@12\fi}
\XFigu3946sp
\ifdim\XFigu<0pt\XFigu-\XFigu\fi
\catcode`\@11
\pgfutil@ifundefined{pgf@pattern@name@xfigp0}{
\pgfdeclarepatternformonly{xfigp0}
{\pgfqpoint{-1bp}{-1bp}}{\pgfqpoint{9bp}{5bp}}{\pgfqpoint{8bp}{4bp}}
{	\pgfsetdash{}{0pt}\pgfsetlinewidth{0.45bp}
	\pgfpathqmoveto{-1bp}{4.5bp}\pgfpathqlineto{9bp}{-0.5bp}
	\pgfusepathqstroke
}
}{}
\catcode`\@12
\definecolor{blue3}{rgb}{0,0,0.82}
\clip(1868,-8694) rectangle (14694,-663);
\tikzset{inner sep=+0pt, outer sep=+0pt}
\pgfsetfillcolor{black}
\pgftext[base,left,at=\pgfqpointxy{12525}{-7050}] {\fontsize{24}{28.8}\normalfont $(f)$}
\pgfsetlinewidth{+30\XFigu}
\pgfsetdash{}{+0pt}
\pgfsetstrokecolor{black}
\draw (6571,-1308) arc[start angle=+28.03, end angle=+-28.03, radius=+1525.5];
\draw (6583,-1305) arc[start angle=+150.30, end angle=+209.70, radius=+1447.1];
\draw (7350,-1575) arc[start angle=+18.92, end angle=+-18.92, radius=+1387.5];
\draw (9375,-1500) arc[start angle=+162.41, end angle=+197.59, radius=+1612.9];
\draw (9992,-1305) arc[start angle=+28.03, end angle=+-28.03, radius=+1525.5];
\draw (10004,-1279) arc[start angle=+150.30, end angle=+209.70, radius=+1447.1];
\pgfsetlinewidth{+45\XFigu}
\pgfsetdash{}{+0pt}
\draw (6583,-1305) arc[start angle=+-113.79, end angle=+-66.21, radius=+4226.2];
\draw (11083,-6039) arc[start angle=+113.79, end angle=+66.21, radius=+4226.2];
\pgfsetlinewidth{+30\XFigu}
\pgfsetdash{}{+0pt}
\draw (11071,-4608) arc[start angle=+28.03, end angle=+-28.03, radius=+1525.5];
\draw (11083,-4605) arc[start angle=+150.30, end angle=+209.70, radius=+1447.1];
\draw (11850,-4875) arc[start angle=+18.92, end angle=+-18.92, radius=+1387.5];
\draw (13875,-4800) arc[start angle=+162.41, end angle=+197.59, radius=+1612.9];
\draw (14492,-4605) arc[start angle=+28.03, end angle=+-28.03, radius=+1525.5];
\draw (14504,-4579) arc[start angle=+150.30, end angle=+209.70, radius=+1447.1];
\pgfsetlinewidth{+45\XFigu}
\pgfsetdash{}{+0pt}
\draw (11083,-4605) arc[start angle=+-113.79, end angle=+-66.21, radius=+4226.2];
\draw (2083,-6039) arc[start angle=+113.79, end angle=+66.21, radius=+4226.2];
\pgfsetlinewidth{+30\XFigu}
\pgfsetdash{}{+0pt}
\draw (2071,-4608) arc[start angle=+28.03, end angle=+-28.03, radius=+1525.5];
\draw (2083,-4605) arc[start angle=+150.30, end angle=+209.70, radius=+1447.1];
\draw (2850,-4875) arc[start angle=+18.92, end angle=+-18.92, radius=+1387.5];
\draw (4875,-4800) arc[start angle=+162.41, end angle=+197.59, radius=+1612.9];
\draw (5492,-4605) arc[start angle=+28.03, end angle=+-28.03, radius=+1525.5];
\draw (5504,-4579) arc[start angle=+150.30, end angle=+209.70, radius=+1447.1];
\pgfsetlinewidth{+45\XFigu}
\pgfsetdash{}{+0pt}
\draw (2083,-4605) arc[start angle=+-113.79, end angle=+-66.21, radius=+4226.2];
\draw (6583,-6039) arc[start angle=+113.79, end angle=+66.21, radius=+4226.2];
\pgfsetlinewidth{+30\XFigu}
\pgfsetdash{}{+0pt}
\draw (6571,-4608) arc[start angle=+28.03, end angle=+-28.03, radius=+1525.5];
\draw (6583,-4605) arc[start angle=+150.30, end angle=+209.70, radius=+1447.1];
\draw (7350,-4875) arc[start angle=+18.92, end angle=+-18.92, radius=+1387.5];
\draw (9375,-4800) arc[start angle=+162.41, end angle=+197.59, radius=+1612.9];
\draw (9992,-4605) arc[start angle=+28.03, end angle=+-28.03, radius=+1525.5];
\draw (10004,-4579) arc[start angle=+150.30, end angle=+209.70, radius=+1447.1];
\pgfsetlinewidth{+45\XFigu}
\pgfsetdash{}{+0pt}
\draw (6583,-4605) arc[start angle=+-113.79, end angle=+-66.21, radius=+4226.2];
\draw (11083,-8664) arc[start angle=+113.79, end angle=+66.21, radius=+4226.2];
\pgfsetlinewidth{+30\XFigu}
\pgfsetdash{}{+0pt}
\draw (11071,-7233) arc[start angle=+28.03, end angle=+-28.03, radius=+1525.5];
\draw (11083,-7230) arc[start angle=+150.30, end angle=+209.70, radius=+1447.1];
\draw (11850,-7500) arc[start angle=+18.92, end angle=+-18.92, radius=+1387.5];
\draw (13875,-7425) arc[start angle=+162.41, end angle=+197.59, radius=+1612.9];
\draw (14492,-7230) arc[start angle=+28.03, end angle=+-28.03, radius=+1525.5];
\draw (14504,-7204) arc[start angle=+150.30, end angle=+209.70, radius=+1447.1];
\pgfsetlinewidth{+45\XFigu}
\pgfsetdash{}{+0pt}
\draw (11083,-7230) arc[start angle=+-113.79, end angle=+-66.21, radius=+4226.2];
\draw (2083,-8664) arc[start angle=+113.79, end angle=+66.21, radius=+4226.2];
\pgfsetlinewidth{+30\XFigu}
\pgfsetdash{}{+0pt}
\draw (2071,-7233) arc[start angle=+28.03, end angle=+-28.03, radius=+1525.5];
\draw (2083,-7230) arc[start angle=+150.30, end angle=+209.70, radius=+1447.1];
\draw (2850,-7500) arc[start angle=+18.92, end angle=+-18.92, radius=+1387.5];
\draw (4875,-7425) arc[start angle=+162.41, end angle=+197.59, radius=+1612.9];
\draw (5492,-7230) arc[start angle=+28.03, end angle=+-28.03, radius=+1525.5];
\draw (5504,-7204) arc[start angle=+150.30, end angle=+209.70, radius=+1447.1];
\pgfsetlinewidth{+45\XFigu}
\pgfsetdash{}{+0pt}
\draw (2083,-7230) arc[start angle=+-113.79, end angle=+-66.21, radius=+4226.2];
\draw (6583,-8664) arc[start angle=+113.79, end angle=+66.21, radius=+4226.2];
\pgfsetlinewidth{+30\XFigu}
\pgfsetdash{}{+0pt}
\draw (6571,-7233) arc[start angle=+28.03, end angle=+-28.03, radius=+1525.5];
\draw (6583,-7230) arc[start angle=+150.30, end angle=+209.70, radius=+1447.1];
\draw (7350,-7500) arc[start angle=+18.92, end angle=+-18.92, radius=+1387.5];
\draw (9375,-7425) arc[start angle=+162.41, end angle=+197.59, radius=+1612.9];
\draw (9992,-7230) arc[start angle=+28.03, end angle=+-28.03, radius=+1525.5];
\draw (10004,-7204) arc[start angle=+150.30, end angle=+209.70, radius=+1447.1];
\pgfsetlinewidth{+45\XFigu}
\pgfsetdash{}{+0pt}
\draw (6583,-7230) arc[start angle=+-113.79, end angle=+-66.21, radius=+4226.2];
\pgfsetdash{{+38\XFigu}{+38\XFigu}}{++0pt}
\pgfsetstrokecolor{blue3}
\draw (3600,-4950) arc[start angle=+23.50, end angle=+-23.50, radius=+940.5];
\draw (4350,-4950) arc[start angle=+149.68, end angle=+221.74, radius=+640.7];
\draw (3000,-7500) arc[start angle=+20.77, end angle=+-20.77, radius=+1163.1];
\draw (2700,-7425) arc[start angle=+16.26, end angle=+-16.26, radius=+1875];
\draw (12075,-7500) arc[start angle=+20.77, end angle=+-20.77, radius=+1163.1];
\draw (8025,-7575) arc[start angle=+23.50, end angle=+-23.50, radius=+940.5];
\draw (8775,-7575) arc[start angle=+149.68, end angle=+221.74, radius=+640.7];
\draw (11700,-7425) arc[start angle=+16.26, end angle=+-16.26, radius=+1875];
\pgfsetlinewidth{+15\XFigu}
\pgfsetdash{{+90\XFigu}{+90\XFigu}}{++0pt}
\pgfsetstrokecolor{black}
\draw (4425,-1350) arc[start angle=+-165.77, end angle=+-74.20, radius=+1570.5];
\pgfsetlinewidth{+30\XFigu}
\pgfsetdash{}{+0pt}
\pgfsetfillpattern{xfigp0}{black}
\draw[pattern,preaction={fill=black}]  (3187,-5325) circle [radius=+38];
\draw[pattern,preaction={fill=black}]  (6975,-5025) circle [radius=+38];
\draw[pattern,preaction={fill=black}]  (7200,-5550) circle [radius=+38];
\draw[pattern,preaction={fill=black}]  (3525,-7950) circle [radius=+38];
\draw[pattern,preaction={fill=black}]  (8250,-7875) circle [radius=+38];
\draw[pattern,preaction={fill=black}]  (11850,-7950) circle [radius=+38];
\draw[pattern,preaction={fill=black}]  (2475,-7725) circle [radius=+38];
\draw[pattern,preaction={fill=black}]  (12075,-7950) circle [radius=+38];
\draw[pattern,preaction={fill=black}]  (11512,-5325) circle [radius=+38];
\draw[pattern,preaction={fill=black}]  (14100,-5250) circle [radius=+38];
\draw[pattern,preaction={fill=black}]  (3975,-5400) circle [radius=+38];
\draw[pattern,preaction={fill=black}]  (8475,-8025) circle [radius=+38];
\pgfsetlinewidth{+7.5\XFigu}
\pgfsetdash{{+60\XFigu}{+60\XFigu}}{++0pt}
\draw (2625,-3300)--(12825,-3300);
\pgfsetfillcolor{black}
\pgftext[base,left,at=\pgfqpointxy{6750}{-4950}] {\fontsize{20}{24}\normalfont $q$}
\pgftext[base,left,at=\pgfqpointxy{14025}{-5025}] {\fontsize{20}{24}\normalfont $y$}
\pgftext[base,left,at=\pgfqpointxy{8025}{-1275}] {\fontsize{24}{28.8}\normalfont $B_{e,t}$}
\pgftext[base,left,at=\pgfqpointxy{3225}{-5250}] {\fontsize{20}{24}\normalfont $q$}
\pgftext[base,left,at=\pgfqpointxy{7200}{-5400}] {\fontsize{20}{24}\normalfont $y$}
\pgftext[base,left,at=\pgfqpointxy{2250}{-7575}] {\fontsize{20}{24}\normalfont $q$}
\pgftext[base,left,at=\pgfqpointxy{3600}{-7875}] {\fontsize{20}{24}\normalfont $y$}
\pgftext[base,left,at=\pgfqpointxy{3600}{-4350}] {\fontsize{24}{28.8}\normalfont $(a)$}
\pgftext[base,left,at=\pgfqpointxy{8475}{-7875}] {\fontsize{20}{24}\normalfont $y$}
\pgftext[base,left,at=\pgfqpointxy{7800}{-7950}] {\fontsize{20}{24}\normalfont $q$}
\pgftext[base,left,at=\pgfqpointxy{11400}{-5100}] {\fontsize{20}{24}\normalfont $q$}
\pgftext[base,left,at=\pgfqpointxy{11475}{-7875}] {\fontsize{20}{24}\normalfont $q$}
\pgftext[base,left,at=\pgfqpointxy{12225}{-7950}] {\fontsize{20}{24}\normalfont $y$}
\pgftext[base,left,at=\pgfqpointxy{8100}{-4350}] {\fontsize{24}{28.8}\normalfont $(b)$}
\pgftext[base,left,at=\pgfqpointxy{12525}{-4350}] {\fontsize{24}{28.8}\normalfont $(c)$}
\pgftext[base,left,at=\pgfqpointxy{3975}{-5250}] {\fontsize{20}{24}\normalfont $y$}
\pgftext[base,left,at=\pgfqpointxy{6750}{-1200}] {\fontsize{24}{28.8}\normalfont $A^e_{u,t}$}
\pgftext[base,left,at=\pgfqpointxy{3975}{-1050}] {\fontsize{28}{33.6}\normalfont $W_{e,t}$}
\pgftext[base,left,at=\pgfqpointxy{9300}{-1200}] {\fontsize{24}{28.8}\normalfont $A^e_{v,t}$}
\pgftext[base,left,at=\pgfqpointxy{3600}{-7050}] {\fontsize{24}{28.8}\normalfont $(d)$}
\pgftext[base,left,at=\pgfqpointxy{8100}{-7050}] {\fontsize{24}{28.8}\normalfont $(e)$}
\pgfsetlinewidth{+45\XFigu}
\pgfsetdash{}{+0pt}
\draw (6583,-2739) arc[start angle=+113.79, end angle=+66.21, radius=+4226.2];
\endtikzpicture}%

%% file: ProofPicAuxiliaryCase.tikz
{\pgfkeys{/pgf/fpu/.try=false}%
\ifx\XFigwidth\undefined\dimen1=0pt\else\dimen1\XFigwidth\fi
\divide\dimen1 by 6773
\ifx\XFigheight\undefined\dimen3=0pt\else\dimen3\XFigheight\fi
\divide\dimen3 by 5712
\ifdim\dimen1=0pt\ifdim\dimen3=0pt\dimen1=3946sp\dimen3\dimen1
  \else\dimen1\dimen3\fi\else\ifdim\dimen3=0pt\dimen3\dimen1\fi\fi
\tikzpicture[x=+\dimen1, y=+\dimen3]
{\ifx\XFigu\undefined\catcode`\@11
\def\temp{\alloc@1\dimen\dimendef\insc@unt}\temp\XFigu\catcode`\@12\fi}
\XFigu3946sp
\ifdim\XFigu<0pt\XFigu-\XFigu\fi
\catcode`\@11
\pgfutil@ifundefined{pgf@pattern@name@xfigp0}{
\pgfdeclarepatternformonly{xfigp0}
{\pgfqpoint{-1bp}{-1bp}}{\pgfqpoint{9bp}{5bp}}{\pgfqpoint{8bp}{4bp}}
{	\pgfsetdash{}{0pt}\pgfsetlinewidth{0.45bp}
	\pgfpathqmoveto{-1bp}{4.5bp}\pgfpathqlineto{9bp}{-0.5bp}
	\pgfusepathqstroke
}
}{}
\catcode`\@12
\pgfdeclarearrow{
  name = xfiga0,
  parameters = {
    \the\pgfarrowlinewidth \the\pgfarrowlength \the\pgfarrowwidth},
  defaults = {
	  line width=+7.5\XFigu, length=+120\XFigu, width=+60\XFigu},
  setup code = {
    \dimen7 2.15\pgfarrowlength\pgfmathveclen{\the\dimen7}{\the\pgfarrowwidth}
    \dimen7 2\pgfarrowwidth\pgfmathdivide{\pgfmathresult}{\the\dimen7}
    \dimen7 \pgfmathresult\pgfarrowlinewidth
    \pgfarrowssettipend{+\dimen7}
    \pgfarrowssetbackend{+-\pgfarrowlength}
    \dimen9 -0.5\pgfarrowlinewidth
    \pgfarrowssetvisualbackend{+\dimen9}
    \pgfarrowssetlineend{+-0.5\pgfarrowlinewidth}
    \pgfarrowshullpoint{+\dimen7}{+0pt}
    \pgfarrowsupperhullpoint{+-\pgfarrowlength}{+0.5\pgfarrowwidth}
    \pgfarrowssavethe\pgfarrowlinewidth
    \pgfarrowssavethe\pgfarrowlength
    \pgfarrowssavethe\pgfarrowwidth
  },
  drawing code = {\pgfsetdash{}{+0pt}
    \ifdim\pgfarrowlinewidth=\pgflinewidth\else\pgfsetlinewidth{+\pgfarrowlinewidth}\fi
    \pgfpathmoveto{\pgfqpoint{-\pgfarrowlength}{0.5\pgfarrowwidth}}
    \pgfpathlineto{\pgfqpoint{0pt}{0pt}}
    \pgfpathlineto{\pgfqpoint{-\pgfarrowlength}{-0.5\pgfarrowwidth}}
    \pgfusepathqstroke
  }
}
\clip(810,-7623) rectangle (7583,-1911);
\tikzset{inner sep=+0pt, outer sep=+0pt}
\pgfsetfillcolor{black}
\pgftext[base,left,at=\pgfqpointxy{1125}{-5925}] {\fontsize{24}{28.8}\normalfont $G=(V,E)$}
\pgftext[base,left,at=\pgfqpointxy{5850}{-7500}] {\fontsize{24}{28.8}\normalfont $u$}
\pgftext[base,left,at=\pgfqpointxy{6600}{-6525}] {\fontsize{24}{28.8}\normalfont $w$}
\pgftext[base,left,at=\pgfqpointxy{6075}{-6750}] {\fontsize{24}{28.8}\normalfont $l_1$}
\pgfsetlinewidth{+30\XFigu}
\pgfsetstrokecolor{black}
\draw (6000,-7050) arc[start angle=+133.60, end angle=+46.40, radius=+1087.5];
\pgfsetlinewidth{+7.5\XFigu}
\pgfsetfillcolor{cyan}
\filldraw  (6750,-6750) circle [radius=+75];
\filldraw  (7500,-7125) circle [radius=+75];
\pgfsetfillcolor{black}
\pgftext[base,left,at=\pgfqpointxy{7425}{-7500}] {\fontsize{24}{28.8}\normalfont $v$}
\pgftext[base,left,at=\pgfqpointxy{7200}{-6825}] {\fontsize{24}{28.8}\normalfont $l_2$}
\pgfsetlinewidth{+30\XFigu}
\draw (6000,-7125) arc[start angle=+-133.60, end angle=+-46.40, radius=+1087.5];
\pgfsetlinewidth{+7.5\XFigu}
\pgfsetfillcolor{cyan}
\filldraw  (3075,-7125) circle [radius=+75];
\pgfsetfillcolor{black}
\pgftext[base,left,at=\pgfqpointxy{3000}{-7500}] {\fontsize{24}{28.8}\normalfont $v$}
\pgfsetlinewidth{+30\XFigu}
\draw (1575,-7125) arc[start angle=+133.60, end angle=+46.40, radius=+1087.5];
\pgftext[base,left,at=\pgfqpointxy{1425}{-7500}] {\fontsize{24}{28.8}\normalfont $u$}
\draw (1575,-7125) arc[start angle=+-133.60, end angle=+-46.40, radius=+1087.5];
\draw (3150,-3300) arc[start angle=+-138.4, end angle=+-41.6, radius=+451.6];
\draw (3300,-3375) arc[start angle=+130.6, end angle=+49.4, radius=+345.7];
\draw (5400,-2850) arc[start angle=+-143.1, end angle=+-36.9, radius=+375];
\draw (5475,-2925) arc[start angle=+126.9, end angle=+53.1, radius=+375];
\draw (4800,-4350) arc[start angle=+-134.2, end angle=+-58.5, radius=+553.1];
\draw (4950,-4425) arc[start angle=+135.0, end angle=+45.0, radius=+265.2];
\pgfsetarrows{[line width=7.5\XFigu]}
\pgfsetarrowsend{xfiga0}
\pgfsetlinewidth{+7.5\XFigu}
\pgfsetdash{{+60\XFigu}{+60\XFigu}}{++0pt}
\draw (2925,-6675) arc[start angle=+139.17, end angle=+34.64, radius=+1757];
\pgfsetlinewidth{+45\XFigu}
\pgfsetdash{}{+0pt}
\draw  (4565,-2928) circle [radius=+546];
\pgfsetlinewidth{+30\XFigu}
\pgfsetfillpattern{xfigp0}{black}
\draw[pattern,preaction={fill=black}]  (3900,-4425) circle [radius=+38];
\draw[pattern,preaction={fill=black}]  (4575,-2925) circle [radius=+38];
\pgfsetlinewidth{+7.5\XFigu}
\pgfsetfillcolor{cyan}
\filldraw  (1575,-7125) circle [radius=+75];
\pgfsetbeveljoin
\pgfsetlinewidth{+45\XFigu}
\draw (4257,-4784)--(4239,-4788)--(4221,-4791)--(4203,-4792)--(4184,-4793)--(4164,-4792)
  --(4143,-4790)--(4122,-4786)--(4100,-4781)--(4077,-4775)--(4053,-4768)--(4029,-4759)
  --(4004,-4748)--(3978,-4737)--(3952,-4724)--(3926,-4709)--(3899,-4694)--(3872,-4677)
  --(3845,-4659)--(3817,-4640)--(3790,-4620)--(3763,-4599)--(3736,-4577)--(3709,-4554)
  --(3682,-4531)--(3656,-4507)--(3629,-4482)--(3603,-4457)--(3577,-4430)--(3551,-4403)
  --(3525,-4376)--(3502,-4350)--(3479,-4324)--(3456,-4297)--(3432,-4269)--(3408,-4241)
  --(3384,-4211)--(3360,-4180)--(3335,-4149)--(3311,-4116)--(3286,-4083)--(3261,-4049)
  --(3236,-4014)--(3212,-3978)--(3187,-3942)--(3163,-3906)--(3140,-3869)--(3116,-3833)
  --(3094,-3796)--(3072,-3759)--(3050,-3723)--(3030,-3686)--(3010,-3651)--(2992,-3616)
  --(2974,-3581)--(2957,-3548)--(2942,-3515)--(2928,-3482)--(2914,-3451)--(2902,-3421)
  --(2891,-3392)--(2881,-3363)--(2872,-3335)--(2864,-3308)--(2858,-3282)--(2851,-3254)
  --(2846,-3226)--(2843,-3198)--(2841,-3171)--(2840,-3145)--(2840,-3120)--(2842,-3094)
  --(2845,-3070)--(2850,-3046)--(2856,-3023)--(2863,-3000)--(2872,-2978)--(2882,-2957)
  --(2893,-2937)--(2905,-2918)--(2919,-2899)--(2934,-2882)--(2949,-2866)--(2966,-2850)
  --(2984,-2836)--(3002,-2824)--(3021,-2812)--(3040,-2801)--(3061,-2791)--(3081,-2783)
  --(3102,-2775)--(3124,-2768)--(3146,-2763)--(3169,-2758)--(3192,-2754)--(3217,-2750)
  --(3243,-2747)--(3270,-2746)--(3298,-2745)--(3326,-2745)--(3356,-2746)--(3386,-2747)
  --(3416,-2750)--(3448,-2754)--(3479,-2759)--(3512,-2764)--(3544,-2771)--(3576,-2778)
  --(3608,-2787)--(3639,-2796)--(3670,-2806)--(3700,-2817)--(3729,-2829)--(3757,-2841)
  --(3783,-2854)--(3809,-2867)--(3833,-2881)--(3855,-2895)--(3876,-2910)--(3896,-2925)
  --(3914,-2941)--(3930,-2957)--(3946,-2974)--(3961,-2994)--(3976,-3015)--(3988,-3037)
  --(3999,-3060)--(4009,-3084)--(4018,-3109)--(4026,-3135)--(4033,-3163)--(4039,-3191)
  --(4045,-3220)--(4049,-3250)--(4054,-3280)--(4058,-3310)--(4062,-3341)--(4066,-3371)
  --(4069,-3401)--(4073,-3431)--(4078,-3460)--(4083,-3489)--(4088,-3517)--(4095,-3544)
  --(4102,-3571)--(4110,-3597)--(4119,-3623)--(4129,-3647)--(4140,-3670)--(4152,-3694)
  --(4165,-3717)--(4180,-3741)--(4195,-3766)--(4211,-3790)--(4228,-3815)--(4246,-3840)
  --(4264,-3866)--(4282,-3891)--(4301,-3917)--(4320,-3943)--(4339,-3969)--(4357,-3994)
  --(4375,-4020)--(4393,-4045)--(4409,-4070)--(4425,-4095)--(4439,-4119)--(4453,-4144)
  --(4465,-4167)--(4476,-4191)--(4486,-4214)--(4494,-4238)--(4501,-4261)--(4506,-4283)
  --(4510,-4305)--(4513,-4327)--(4515,-4350)--(4515,-4373)--(4515,-4397)--(4513,-4420)
  --(4510,-4444)--(4506,-4468)--(4501,-4492)--(4494,-4515)--(4487,-4539)--(4478,-4562)
  --(4468,-4584)--(4457,-4606)--(4446,-4627)--(4433,-4647)--(4420,-4666)--(4406,-4684)
  --(4391,-4700)--(4376,-4715)--(4360,-4729)--(4344,-4742)--(4327,-4753)--(4310,-4763)
  --(4293,-4771)--(4275,-4778)--cycle;
\draw (6033,-2420)--(6010,-2407)--(5987,-2395)--(5962,-2384)--(5936,-2373)--(5910,-2362)
  --(5882,-2352)--(5854,-2343)--(5824,-2334)--(5794,-2326)--(5764,-2319)--(5733,-2312)
  --(5702,-2306)--(5671,-2301)--(5640,-2297)--(5610,-2293)--(5581,-2291)--(5552,-2290)
  --(5524,-2289)--(5497,-2289)--(5471,-2291)--(5447,-2293)--(5423,-2296)--(5400,-2299)
  --(5379,-2304)--(5358,-2309)--(5338,-2316)--(5318,-2323)--(5299,-2331)--(5280,-2340)
  --(5262,-2351)--(5245,-2362)--(5228,-2374)--(5213,-2388)--(5197,-2402)--(5183,-2418)
  --(5170,-2434)--(5158,-2451)--(5147,-2468)--(5137,-2486)--(5129,-2505)--(5121,-2524)
  --(5115,-2543)--(5110,-2562)--(5106,-2582)--(5103,-2602)--(5101,-2622)--(5100,-2642)
  --(5101,-2662)--(5102,-2683)--(5104,-2704)--(5107,-2726)--(5111,-2749)--(5117,-2773)
  --(5123,-2797)--(5130,-2822)--(5137,-2848)--(5146,-2875)--(5156,-2903)--(5166,-2930)
  --(5177,-2959)--(5188,-2987)--(5200,-3016)--(5212,-3044)--(5224,-3073)--(5236,-3101)
  --(5249,-3128)--(5261,-3156)--(5273,-3182)--(5285,-3209)--(5297,-3235)--(5309,-3260)
  --(5321,-3286)--(5332,-3311)--(5344,-3336)--(5355,-3362)--(5367,-3388)--(5378,-3415)
  --(5389,-3441)--(5400,-3468)--(5411,-3496)--(5422,-3523)--(5433,-3551)--(5443,-3579)
  --(5452,-3606)--(5461,-3633)--(5470,-3659)--(5477,-3685)--(5484,-3710)--(5491,-3734)
  --(5496,-3757)--(5500,-3780)--(5504,-3801)--(5507,-3821)--(5509,-3840)--(5510,-3859)
  --(5510,-3877)--(5510,-3895)--(5508,-3913)--(5505,-3931)--(5502,-3948)--(5497,-3964)
  --(5491,-3980)--(5484,-3996)--(5476,-4010)--(5468,-4024)--(5458,-4037)--(5447,-4049)
  --(5436,-4060)--(5424,-4070)--(5412,-4079)--(5399,-4087)--(5386,-4093)--(5373,-4098)
  --(5359,-4102)--(5346,-4105)--(5332,-4107)--(5318,-4107)--(5304,-4107)--(5288,-4105)
  --(5272,-4103)--(5255,-4099)--(5238,-4094)--(5220,-4088)--(5201,-4081)--(5182,-4073)
  --(5162,-4065)--(5141,-4056)--(5120,-4047)--(5100,-4038)--(5079,-4029)--(5058,-4020)
  --(5037,-4011)--(5017,-4003)--(4997,-3996)--(4978,-3989)--(4958,-3984)--(4940,-3979)
  --(4921,-3975)--(4902,-3972)--(4884,-3970)--(4865,-3969)--(4845,-3970)--(4826,-3971)
  --(4806,-3973)--(4786,-3976)--(4766,-3981)--(4746,-3986)--(4727,-3992)--(4708,-4000)
  --(4690,-4008)--(4673,-4016)--(4657,-4026)--(4641,-4036)--(4627,-4046)--(4615,-4057)
  --(4603,-4068)--(4592,-4080)--(4582,-4092)--(4574,-4105)--(4566,-4118)--(4559,-4132)
  --(4553,-4147)--(4547,-4162)--(4543,-4179)--(4540,-4196)--(4538,-4214)--(4537,-4232)
  --(4537,-4251)--(4538,-4271)--(4540,-4290)--(4544,-4309)--(4549,-4329)--(4554,-4348)
  --(4561,-4367)--(4568,-4385)--(4577,-4403)--(4586,-4421)--(4597,-4439)--(4607,-4455)
  --(4618,-4471)--(4631,-4487)--(4644,-4504)--(4659,-4521)--(4675,-4538)--(4692,-4555)
  --(4710,-4573)--(4730,-4590)--(4751,-4608)--(4772,-4625)--(4795,-4641)--(4818,-4658)
  --(4842,-4673)--(4867,-4689)--(4892,-4703)--(4917,-4717)--(4942,-4729)--(4968,-4741)
  --(4994,-4753)--(5020,-4763)--(5047,-4773)--(5067,-4780)--(5087,-4786)--(5108,-4792)
  --(5129,-4798)--(5151,-4803)--(5174,-4807)--(5198,-4811)--(5222,-4814)--(5246,-4816)
  --(5271,-4817)--(5297,-4817)--(5323,-4816)--(5350,-4814)--(5377,-4810)--(5404,-4805)
  --(5431,-4799)--(5458,-4791)--(5486,-4781)--(5513,-4771)--(5539,-4758)--(5566,-4744)
  --(5592,-4729)--(5617,-4712)--(5643,-4693)--(5667,-4673)--(5691,-4651)--(5715,-4628)
  --(5739,-4603)--(5762,-4576)--(5785,-4547)--(5802,-4523)--(5820,-4498)--(5838,-4471)
  --(5856,-4443)--(5873,-4413)--(5891,-4382)--(5909,-4350)--(5927,-4316)--(5945,-4281)
  --(5962,-4245)--(5980,-4207)--(5998,-4168)--(6015,-4128)--(6032,-4088)--(6049,-4046)
  --(6066,-4003)--(6083,-3960)--(6099,-3917)--(6115,-3873)--(6130,-3829)--(6145,-3785)
  --(6159,-3740)--(6173,-3696)--(6186,-3653)--(6199,-3610)--(6211,-3567)--(6222,-3525)
  --(6233,-3484)--(6243,-3444)--(6252,-3405)--(6260,-3366)--(6268,-3329)--(6276,-3293)
  --(6282,-3258)--(6288,-3224)--(6293,-3190)--(6298,-3158)--(6302,-3127)--(6306,-3089)
  --(6309,-3052)--(6311,-3017)--(6312,-2982)--(6312,-2948)--(6311,-2916)--(6309,-2884)
  --(6306,-2853)--(6302,-2822)--(6297,-2793)--(6291,-2765)--(6283,-2737)--(6275,-2711)
  --(6266,-2685)--(6256,-2661)--(6246,-2637)--(6234,-2615)--(6222,-2594)--(6208,-2574)
  --(6195,-2555)--(6180,-2538)--(6166,-2521)--(6150,-2506)--(6135,-2491)--(6119,-2477)
  --(6102,-2464)--(6086,-2452)--(6069,-2441)--(6051,-2430)--cycle;
\pgfsetfillcolor{black}
\pgftext[base,left,at=\pgfqpointxy{3750}{-4275}] {\fontsize{24}{28.8}\normalfont $p_t$}
\pgftext[base,left,at=\pgfqpointxy{4500}{-3225}] {\fontsize{18}{21.6}\normalfont $1_w$}
\pgftext[base,left,at=\pgfqpointxy{4125}{-2250}] {\fontsize{24}{28.8}\normalfont $\P^1_w$}
\pgftext[base,left,at=\pgfqpointxy{825}{-2625}] {\fontsize{24}{28.8}\normalfont $\widehat\rsf_{t_0}$}
\pgftext[base,left,at=\pgfqpointxy{5400}{-5925}] {\fontsize{24}{28.8}\normalfont $\widehat G=(\widehat V,\widehat E)$}
\pgftext[base,left,at=\pgfqpointxy{2250}{-6675}] {\fontsize{24}{28.8}\normalfont $l$}
\pgfsetlinewidth{+7.5\XFigu}
\pgfsetfillcolor{cyan}
\filldraw  (6000,-7125) circle [radius=+75];
\endtikzpicture}%

%% file: ProofPicTwoCylinders.tikz
{\pgfkeys{/pgf/fpu/.try=false}%
\ifx\XFigwidth\undefined\dimen1=0pt\else\dimen1\XFigwidth\fi
\divide\dimen1 by 7823
\ifx\XFigheight\undefined\dimen3=0pt\else\dimen3\XFigheight\fi
\divide\dimen3 by 7975
\ifdim\dimen1=0pt\ifdim\dimen3=0pt\dimen1=3946sp\dimen3\dimen1
  \else\dimen1\dimen3\fi\else\ifdim\dimen3=0pt\dimen3\dimen1\fi\fi
\tikzpicture[x=+\dimen1, y=+\dimen3]
{\ifx\XFigu\undefined\catcode`\@11
\def\temp{\alloc@1\dimen\dimendef\insc@unt}\temp\XFigu\catcode`\@12\fi}
\XFigu3946sp
\ifdim\XFigu<0pt\XFigu-\XFigu\fi
\catcode`\@11
\pgfutil@ifundefined{pgf@pattern@name@xfigp0}{
\pgfdeclarepatternformonly{xfigp0}
{\pgfqpoint{-1bp}{-1bp}}{\pgfqpoint{9bp}{5bp}}{\pgfqpoint{8bp}{4bp}}
{	\pgfsetdash{}{0pt}\pgfsetlinewidth{0.45bp}
	\pgfpathqmoveto{-1bp}{4.5bp}\pgfpathqlineto{9bp}{-0.5bp}
	\pgfusepathqstroke
}
}{}
\catcode`\@12
\clip(1485,-9211) rectangle (9308,-1236);
\tikzset{inner sep=+0pt, outer sep=+0pt}
\pgfsetfillcolor{black}
\pgftext[base,left,at=\pgfqpointxy{7125}{-2100}] {\fontsize{24}{28.8}\normalfont $G=(V,E)$}
\pgfsetlinewidth{+30\XFigu}
\pgfsetstrokecolor{black}
\draw (2250,-2400) arc[start angle=+130.6, end angle=+49.4, radius=+345.7];
\pgfsetbeveljoin
\pgfsetlinewidth{+45\XFigu}
\draw (2541,-4045)--(2523,-4041)--(2506,-4036)--(2488,-4030)--(2470,-4022)--(2453,-4013)
  --(2435,-4002)--(2417,-3990)--(2399,-3976)--(2381,-3961)--(2363,-3944)--(2344,-3926)
  --(2326,-3906)--(2308,-3884)--(2290,-3861)--(2272,-3837)--(2254,-3812)--(2237,-3785)
  --(2220,-3757)--(2203,-3729)--(2187,-3699)--(2171,-3668)--(2156,-3637)--(2141,-3605)
  --(2127,-3573)--(2113,-3540)--(2099,-3506)--(2087,-3472)--(2074,-3437)--(2062,-3402)
  --(2050,-3366)--(2040,-3333)--(2030,-3299)--(2020,-3265)--(2011,-3230)--(2001,-3194)
  --(1992,-3156)--(1983,-3118)--(1974,-3079)--(1965,-3039)--(1957,-2999)--(1949,-2957)
  --(1941,-2915)--(1934,-2873)--(1927,-2830)--(1921,-2787)--(1915,-2744)--(1909,-2700)
  --(1904,-2657)--(1900,-2615)--(1896,-2573)--(1893,-2531)--(1890,-2491)--(1888,-2451)
  --(1886,-2412)--(1886,-2375)--(1886,-2338)--(1886,-2303)--(1887,-2269)--(1889,-2237)
  --(1891,-2205)--(1894,-2175)--(1898,-2146)--(1902,-2119)--(1907,-2092)--(1914,-2064)
  --(1921,-2036)--(1929,-2010)--(1938,-1985)--(1949,-1960)--(1960,-1937)--(1972,-1915)
  --(1985,-1894)--(2000,-1875)--(2015,-1856)--(2031,-1839)--(2048,-1822)--(2066,-1808)
  --(2085,-1794)--(2104,-1782)--(2124,-1771)--(2145,-1761)--(2166,-1753)--(2188,-1747)
  --(2209,-1741)--(2231,-1737)--(2254,-1735)--(2276,-1733)--(2298,-1733)--(2321,-1734)
  --(2343,-1736)--(2365,-1739)--(2388,-1743)--(2410,-1748)--(2433,-1754)--(2457,-1762)
  --(2482,-1770)--(2507,-1780)--(2533,-1791)--(2559,-1803)--(2585,-1816)--(2611,-1831)
  --(2638,-1846)--(2665,-1863)--(2691,-1881)--(2718,-1899)--(2744,-1919)--(2770,-1939)
  --(2796,-1960)--(2820,-1982)--(2844,-2004)--(2866,-2027)--(2888,-2050)--(2908,-2072)
  --(2927,-2095)--(2944,-2118)--(2960,-2141)--(2974,-2163)--(2987,-2186)--(2998,-2208)
  --(3008,-2230)--(3016,-2251)--(3023,-2273)--(3029,-2298)--(3033,-2323)--(3035,-2348)
  --(3035,-2374)--(3034,-2400)--(3032,-2426)--(3028,-2453)--(3022,-2481)--(3016,-2509)
  --(3008,-2538)--(3000,-2567)--(2992,-2596)--(2982,-2625)--(2973,-2655)--(2964,-2684)
  --(2954,-2713)--(2945,-2741)--(2937,-2770)--(2929,-2798)--(2923,-2826)--(2917,-2853)
  --(2912,-2880)--(2908,-2908)--(2906,-2935)--(2905,-2960)--(2905,-2986)--(2906,-3013)
  --(2908,-3040)--(2911,-3068)--(2914,-3096)--(2919,-3125)--(2923,-3155)--(2929,-3185)
  --(2935,-3216)--(2941,-3247)--(2947,-3278)--(2953,-3309)--(2959,-3340)--(2965,-3372)
  --(2971,-3402)--(2976,-3432)--(2980,-3462)--(2984,-3491)--(2987,-3519)--(2989,-3547)
  --(2990,-3574)--(2990,-3600)--(2988,-3625)--(2986,-3650)--(2983,-3674)--(2978,-3696)
  --(2973,-3718)--(2966,-3739)--(2958,-3761)--(2948,-3782)--(2938,-3803)--(2927,-3823)
  --(2914,-3844)--(2900,-3864)--(2885,-3883)--(2869,-3902)--(2852,-3920)--(2835,-3937)
  --(2816,-3953)--(2798,-3968)--(2778,-3982)--(2758,-3995)--(2738,-4007)--(2718,-4017)
  --(2698,-4025)--(2677,-4033)--(2657,-4039)--(2637,-4043)--(2617,-4046)--(2598,-4048)
  --(2578,-4048)--(2559,-4047)--cycle;
\pgfsetlinewidth{+30\XFigu}
\pgfsetfillpattern{xfigp0}{black}
\draw[pattern,preaction={fill=black}]  (2325,-3450) circle [radius=+38];
\pgfsetfillcolor{black}
\pgftext[base,left,at=\pgfqpointxy{2250}{-3300}] {\fontsize{24}{28.8}\normalfont $p_0$}
\draw (3375,-4050) arc[start angle=+135.0, end angle=+45.0, radius=+265.2];
\draw (3225,-3975) arc[start angle=+-134.2, end angle=+-58.5, radius=+553.1];
\draw (3600,-2475) arc[start angle=+126.9, end angle=+53.1, radius=+375];
\draw (3525,-2400) arc[start angle=+-143.1, end angle=+-36.9, radius=+375];
\pgfsetlinewidth{+45\XFigu}
\draw (3913,-1904)--(3888,-1898)--(3862,-1893)--(3835,-1888)--(3807,-1884)--(3779,-1881)
  --(3750,-1879)--(3720,-1877)--(3689,-1876)--(3658,-1876)--(3627,-1876)--(3595,-1878)
  --(3564,-1880)--(3532,-1883)--(3502,-1887)--(3472,-1892)--(3442,-1897)--(3414,-1904)
  --(3387,-1910)--(3361,-1918)--(3337,-1925)--(3313,-1934)--(3291,-1943)--(3270,-1952)
  --(3251,-1962)--(3232,-1973)--(3214,-1984)--(3197,-1996)--(3181,-2009)--(3165,-2023)
  --(3151,-2038)--(3137,-2053)--(3124,-2070)--(3112,-2087)--(3101,-2105)--(3092,-2123)
  --(3083,-2142)--(3076,-2162)--(3070,-2182)--(3065,-2202)--(3061,-2222)--(3059,-2242)
  --(3058,-2262)--(3058,-2282)--(3059,-2302)--(3062,-2322)--(3065,-2342)--(3070,-2361)
  --(3075,-2381)--(3082,-2401)--(3089,-2421)--(3098,-2441)--(3108,-2462)--(3119,-2483)
  --(3131,-2505)--(3145,-2528)--(3159,-2551)--(3174,-2575)--(3190,-2599)--(3207,-2623)
  --(3225,-2648)--(3243,-2672)--(3262,-2697)--(3281,-2721)--(3300,-2745)--(3319,-2769)
  --(3338,-2793)--(3357,-2816)--(3376,-2838)--(3395,-2861)--(3413,-2883)--(3431,-2904)
  --(3449,-2926)--(3466,-2947)--(3484,-2969)--(3502,-2991)--(3520,-3013)--(3537,-3035)
  --(3555,-3058)--(3573,-3082)--(3591,-3105)--(3608,-3129)--(3626,-3153)--(3643,-3177)
  --(3659,-3201)--(3675,-3225)--(3689,-3248)--(3704,-3271)--(3717,-3293)--(3729,-3315)
  --(3740,-3336)--(3750,-3356)--(3759,-3376)--(3767,-3395)--(3774,-3413)--(3780,-3430)
  --(3785,-3448)--(3789,-3466)--(3792,-3484)--(3794,-3501)--(3795,-3519)--(3795,-3536)
  --(3793,-3553)--(3790,-3569)--(3787,-3586)--(3782,-3601)--(3776,-3616)--(3769,-3631)
  --(3761,-3644)--(3752,-3657)--(3742,-3669)--(3732,-3679)--(3721,-3689)--(3709,-3697)
  --(3697,-3705)--(3684,-3711)--(3672,-3716)--(3658,-3720)--(3645,-3724)--(3629,-3726)
  --(3613,-3727)--(3595,-3728)--(3577,-3728)--(3558,-3727)--(3538,-3725)--(3518,-3722)
  --(3496,-3720)--(3474,-3716)--(3452,-3713)--(3429,-3709)--(3407,-3706)--(3384,-3703)
  --(3362,-3700)--(3341,-3697)--(3320,-3695)--(3299,-3694)--(3279,-3693)--(3260,-3694)
  --(3241,-3695)--(3222,-3697)--(3204,-3700)--(3185,-3704)--(3166,-3709)--(3148,-3715)
  --(3129,-3722)--(3111,-3731)--(3093,-3740)--(3075,-3750)--(3058,-3761)--(3042,-3773)
  --(3027,-3786)--(3012,-3798)--(2999,-3812)--(2987,-3825)--(2976,-3839)--(2966,-3853)
  --(2958,-3867)--(2951,-3881)--(2945,-3895)--(2939,-3910)--(2935,-3925)--(2932,-3940)
  --(2930,-3956)--(2929,-3972)--(2929,-3989)--(2930,-4007)--(2933,-4025)--(2936,-4043)
  --(2941,-4061)--(2948,-4079)--(2955,-4097)--(2963,-4115)--(2973,-4133)--(2983,-4150)
  --(2994,-4166)--(3006,-4182)--(3019,-4197)--(3033,-4212)--(3048,-4227)--(3062,-4240)
  --(3077,-4252)--(3093,-4265)--(3110,-4278)--(3129,-4290)--(3149,-4303)--(3170,-4315)
  --(3192,-4327)--(3215,-4339)--(3240,-4350)--(3265,-4361)--(3292,-4371)--(3318,-4381)
  --(3346,-4390)--(3373,-4398)--(3401,-4406)--(3429,-4412)--(3457,-4418)--(3485,-4423)
  --(3513,-4427)--(3541,-4431)--(3569,-4433)--(3590,-4435)--(3611,-4436)--(3633,-4436)
  --(3655,-4436)--(3678,-4435)--(3701,-4433)--(3724,-4431)--(3748,-4428)--(3773,-4423)
  --(3797,-4418)--(3822,-4411)--(3847,-4403)--(3872,-4394)--(3898,-4384)--(3922,-4372)
  --(3947,-4358)--(3971,-4344)--(3995,-4328)--(4018,-4310)--(4041,-4291)--(4063,-4271)
  --(4084,-4250)--(4105,-4226)--(4124,-4202)--(4143,-4176)--(4161,-4149)--(4177,-4120)
  --(4194,-4090)--(4209,-4058)--(4224,-4024)--(4235,-3996)--(4245,-3967)--(4255,-3937)
  --(4265,-3905)--(4275,-3872)--(4284,-3837)--(4293,-3801)--(4301,-3764)--(4309,-3726)
  --(4317,-3686)--(4324,-3645)--(4331,-3603)--(4338,-3560)--(4344,-3516)--(4350,-3471)
  --(4355,-3426)--(4360,-3380)--(4364,-3334)--(4368,-3287)--(4371,-3241)--(4374,-3194)
  --(4377,-3148)--(4379,-3102)--(4380,-3056)--(4381,-3011)--(4382,-2967)--(4382,-2924)
  --(4381,-2881)--(4381,-2840)--(4379,-2800)--(4378,-2760)--(4376,-2722)--(4373,-2685)
  --(4371,-2650)--(4367,-2615)--(4364,-2582)--(4360,-2550)--(4356,-2519)--(4350,-2481)
  --(4343,-2445)--(4336,-2410)--(4328,-2376)--(4319,-2343)--(4310,-2312)--(4300,-2282)
  --(4289,-2252)--(4277,-2224)--(4264,-2197)--(4251,-2171)--(4237,-2146)--(4222,-2123)
  --(4207,-2101)--(4191,-2080)--(4175,-2060)--(4158,-2041)--(4140,-2024)--(4123,-2008)
  --(4104,-1994)--(4086,-1980)--(4067,-1968)--(4049,-1957)--(4030,-1947)--(4011,-1938)
  --(3991,-1930)--(3972,-1922)--(3953,-1916)--(3933,-1910)--cycle;
\pgfsetlinewidth{+30\XFigu}
\draw (7725,-3075) arc[start angle=+133.60, end angle=+46.40, radius=+1087.5];
\draw (7725,-3075) arc[start angle=+-133.60, end angle=+-46.40, radius=+1087.5];
\pgfsetlinewidth{+7.5\XFigu}
\pgfsetfillcolor{cyan}
\filldraw  (9225,-3075) circle [radius=+75];
\filldraw  (7725,-3075) circle [radius=+75];
\pgfsetfillcolor{black}
\pgftext[base,left,at=\pgfqpointxy{8325}{-2625}] {\fontsize{24}{28.8}\normalfont $e$}
\pgfsetlinewidth{+30\XFigu}
\pgfsetfillpattern{xfigp0}{black}
\draw[pattern,preaction={fill=black}]  (2475,-7575) circle [radius=+38];
\pgfsetfillcolor{black}
\pgftext[base,left,at=\pgfqpointxy{2400}{-7425}] {\fontsize{24}{28.8}\normalfont $p_0$}
\draw (2250,-6450) arc[start angle=+-138.4, end angle=+-41.6, radius=+451.6];
\draw (2400,-6525) arc[start angle=+130.6, end angle=+49.4, radius=+345.7];
\draw (4650,-8325) arc[start angle=+135.0, end angle=+45.0, radius=+265.2];
\draw (4500,-8250) arc[start angle=+-134.2, end angle=+-58.5, radius=+553.1];
\draw (4875,-6750) arc[start angle=+126.9, end angle=+53.1, radius=+375];
\draw (4800,-6675) arc[start angle=+-143.1, end angle=+-36.9, radius=+375];
\pgfsetdash{}{+0pt}
\draw (9225,-6828) arc[start angle=+-11.77, end angle=+-93.10, radius=+654.4];
\draw (7725,-6825) arc[start angle=+-168.48, end angle=+-87.27, radius=+657];
\pgfsetdash{}{+0pt}
\draw (7725,-6675) arc[start angle=+133.60, end angle=+46.40, radius=+1087.5];
\pgfsetlinewidth{+7.5\XFigu}
\pgfsetfillcolor{cyan}
\filldraw  (7725,-6750) circle [radius=+75];
\filldraw  (8475,-6375) circle [radius=+75];
\filldraw  (8475,-7350) circle [radius=+75];
\filldraw  (9225,-6750) circle [radius=+75];
\pgfsetlinewidth{+30\XFigu}
\pgfsetfillpattern{xfigp0}{black}
\draw[pattern,preaction={fill=black}]  (3750,-6450) circle [radius=+38];
\pgfsetlinewidth{+45\XFigu}
\draw  (3734,-6410) circle [radius=+546];
\draw  (3629,-8151) circle [radius=+546];
\pgfsetlinewidth{+30\XFigu}
\draw[pattern,preaction={fill=black}]  (3675,-8175) circle [radius=+38];
\pgfsetlinewidth{+45\XFigu}
\draw (2691,-8170)--(2673,-8166)--(2656,-8161)--(2638,-8155)--(2620,-8147)--(2603,-8138)
  --(2585,-8127)--(2567,-8115)--(2549,-8101)--(2531,-8086)--(2513,-8069)--(2494,-8051)
  --(2476,-8031)--(2458,-8009)--(2440,-7986)--(2422,-7962)--(2404,-7937)--(2387,-7910)
  --(2370,-7882)--(2353,-7854)--(2337,-7824)--(2321,-7793)--(2306,-7762)--(2291,-7730)
  --(2277,-7698)--(2263,-7665)--(2249,-7631)--(2237,-7597)--(2224,-7562)--(2212,-7527)
  --(2200,-7491)--(2190,-7458)--(2180,-7424)--(2170,-7390)--(2161,-7355)--(2151,-7319)
  --(2142,-7281)--(2133,-7243)--(2124,-7204)--(2115,-7164)--(2107,-7124)--(2099,-7082)
  --(2091,-7040)--(2084,-6998)--(2077,-6955)--(2071,-6912)--(2065,-6869)--(2059,-6825)
  --(2054,-6782)--(2050,-6740)--(2046,-6698)--(2043,-6656)--(2040,-6616)--(2038,-6576)
  --(2036,-6537)--(2036,-6500)--(2036,-6463)--(2036,-6428)--(2037,-6394)--(2039,-6362)
  --(2041,-6330)--(2044,-6300)--(2048,-6271)--(2052,-6244)--(2057,-6217)--(2064,-6189)
  --(2071,-6161)--(2079,-6135)--(2088,-6110)--(2099,-6085)--(2110,-6062)--(2122,-6040)
  --(2135,-6019)--(2150,-6000)--(2165,-5981)--(2181,-5964)--(2198,-5947)--(2216,-5933)
  --(2235,-5919)--(2254,-5907)--(2274,-5896)--(2295,-5886)--(2316,-5878)--(2338,-5872)
  --(2359,-5866)--(2381,-5862)--(2404,-5860)--(2426,-5858)--(2448,-5858)--(2471,-5859)
  --(2493,-5861)--(2515,-5864)--(2538,-5868)--(2560,-5873)--(2583,-5879)--(2607,-5887)
  --(2632,-5895)--(2657,-5905)--(2683,-5916)--(2709,-5928)--(2735,-5941)--(2761,-5956)
  --(2788,-5971)--(2815,-5988)--(2841,-6006)--(2868,-6024)--(2894,-6044)--(2920,-6064)
  --(2946,-6085)--(2970,-6107)--(2994,-6129)--(3016,-6152)--(3038,-6175)--(3058,-6197)
  --(3077,-6220)--(3094,-6243)--(3110,-6266)--(3124,-6288)--(3137,-6311)--(3148,-6333)
  --(3158,-6355)--(3166,-6376)--(3173,-6398)--(3179,-6423)--(3183,-6448)--(3185,-6473)
  --(3185,-6499)--(3184,-6525)--(3182,-6551)--(3178,-6578)--(3172,-6606)--(3166,-6634)
  --(3158,-6663)--(3150,-6692)--(3142,-6721)--(3132,-6750)--(3123,-6780)--(3114,-6809)
  --(3104,-6838)--(3095,-6866)--(3087,-6895)--(3079,-6923)--(3073,-6951)--(3067,-6978)
  --(3062,-7005)--(3058,-7033)--(3056,-7060)--(3055,-7085)--(3055,-7111)--(3056,-7138)
  --(3058,-7165)--(3061,-7193)--(3064,-7221)--(3069,-7250)--(3073,-7280)--(3079,-7310)
  --(3085,-7341)--(3091,-7372)--(3097,-7403)--(3103,-7434)--(3109,-7465)--(3115,-7497)
  --(3121,-7527)--(3126,-7557)--(3130,-7587)--(3134,-7616)--(3137,-7644)--(3139,-7672)
  --(3140,-7699)--(3140,-7725)--(3138,-7750)--(3136,-7775)--(3133,-7799)--(3128,-7821)
  --(3123,-7843)--(3116,-7864)--(3108,-7886)--(3098,-7907)--(3088,-7928)--(3077,-7948)
  --(3064,-7969)--(3050,-7989)--(3035,-8008)--(3019,-8027)--(3002,-8045)--(2985,-8062)
  --(2966,-8078)--(2948,-8093)--(2928,-8107)--(2908,-8120)--(2888,-8132)--(2868,-8142)
  --(2848,-8150)--(2827,-8158)--(2807,-8164)--(2787,-8168)--(2767,-8171)--(2748,-8173)
  --(2728,-8173)--(2709,-8172)--cycle;
\draw (5188,-6179)--(5163,-6173)--(5137,-6168)--(5110,-6163)--(5082,-6159)--(5054,-6156)
  --(5025,-6154)--(4995,-6152)--(4964,-6151)--(4933,-6151)--(4902,-6151)--(4870,-6153)
  --(4839,-6155)--(4807,-6158)--(4777,-6162)--(4747,-6167)--(4717,-6172)--(4689,-6179)
  --(4662,-6185)--(4636,-6193)--(4612,-6200)--(4588,-6209)--(4566,-6218)--(4545,-6227)
  --(4526,-6237)--(4507,-6248)--(4489,-6259)--(4472,-6271)--(4456,-6284)--(4440,-6298)
  --(4426,-6313)--(4412,-6328)--(4399,-6345)--(4387,-6362)--(4376,-6380)--(4367,-6398)
  --(4358,-6417)--(4351,-6437)--(4345,-6457)--(4340,-6477)--(4336,-6497)--(4334,-6517)
  --(4333,-6537)--(4333,-6557)--(4334,-6577)--(4337,-6597)--(4340,-6617)--(4345,-6636)
  --(4350,-6656)--(4357,-6676)--(4364,-6696)--(4373,-6716)--(4383,-6737)--(4394,-6758)
  --(4406,-6780)--(4420,-6803)--(4434,-6826)--(4449,-6850)--(4465,-6874)--(4482,-6898)
  --(4500,-6923)--(4518,-6947)--(4537,-6972)--(4556,-6996)--(4575,-7020)--(4594,-7044)
  --(4613,-7068)--(4632,-7091)--(4651,-7113)--(4670,-7136)--(4688,-7158)--(4706,-7179)
  --(4724,-7201)--(4741,-7222)--(4759,-7244)--(4777,-7266)--(4795,-7288)--(4812,-7310)
  --(4830,-7333)--(4848,-7357)--(4866,-7380)--(4883,-7404)--(4901,-7428)--(4918,-7452)
  --(4934,-7476)--(4950,-7500)--(4964,-7523)--(4979,-7546)--(4992,-7568)--(5004,-7590)
  --(5015,-7611)--(5025,-7631)--(5034,-7651)--(5042,-7670)--(5049,-7688)--(5055,-7705)
  --(5060,-7723)--(5064,-7741)--(5067,-7759)--(5069,-7776)--(5070,-7794)--(5070,-7811)
  --(5068,-7828)--(5065,-7844)--(5062,-7861)--(5057,-7876)--(5051,-7891)--(5044,-7906)
  --(5036,-7919)--(5027,-7932)--(5017,-7944)--(5007,-7954)--(4996,-7964)--(4984,-7972)
  --(4972,-7980)--(4959,-7986)--(4947,-7991)--(4933,-7995)--(4920,-7999)--(4904,-8001)
  --(4888,-8002)--(4870,-8003)--(4852,-8003)--(4833,-8002)--(4813,-8000)--(4793,-7997)
  --(4771,-7995)--(4749,-7991)--(4727,-7988)--(4704,-7984)--(4682,-7981)--(4659,-7978)
  --(4637,-7975)--(4616,-7972)--(4595,-7970)--(4574,-7969)--(4554,-7968)--(4535,-7969)
  --(4516,-7970)--(4497,-7972)--(4479,-7975)--(4460,-7979)--(4441,-7984)--(4423,-7990)
  --(4404,-7997)--(4386,-8006)--(4368,-8015)--(4350,-8025)--(4333,-8036)--(4317,-8048)
  --(4302,-8061)--(4287,-8073)--(4274,-8087)--(4262,-8100)--(4251,-8114)--(4241,-8128)
  --(4233,-8142)--(4226,-8156)--(4220,-8170)--(4214,-8185)--(4210,-8200)--(4207,-8215)
  --(4205,-8231)--(4204,-8247)--(4204,-8264)--(4205,-8282)--(4208,-8300)--(4211,-8318)
  --(4216,-8336)--(4223,-8354)--(4230,-8372)--(4238,-8390)--(4248,-8408)--(4258,-8425)
  --(4269,-8441)--(4281,-8457)--(4294,-8472)--(4308,-8487)--(4323,-8502)--(4337,-8515)
  --(4352,-8527)--(4368,-8540)--(4385,-8553)--(4404,-8565)--(4424,-8578)--(4445,-8590)
  --(4467,-8602)--(4490,-8614)--(4515,-8625)--(4540,-8636)--(4567,-8646)--(4593,-8656)
  --(4621,-8665)--(4648,-8673)--(4676,-8681)--(4704,-8687)--(4732,-8693)--(4760,-8698)
  --(4788,-8702)--(4816,-8706)--(4844,-8708)--(4865,-8710)--(4886,-8711)--(4908,-8711)
  --(4930,-8711)--(4953,-8710)--(4976,-8708)--(4999,-8706)--(5023,-8703)--(5048,-8698)
  --(5072,-8693)--(5097,-8686)--(5122,-8678)--(5147,-8669)--(5173,-8659)--(5197,-8647)
  --(5222,-8633)--(5246,-8619)--(5270,-8603)--(5293,-8585)--(5316,-8566)--(5338,-8546)
  --(5359,-8525)--(5380,-8501)--(5399,-8477)--(5418,-8451)--(5436,-8424)--(5452,-8395)
  --(5469,-8365)--(5484,-8333)--(5499,-8299)--(5510,-8271)--(5520,-8242)--(5530,-8212)
  --(5540,-8180)--(5550,-8147)--(5559,-8112)--(5568,-8076)--(5576,-8039)--(5584,-8001)
  --(5592,-7961)--(5599,-7920)--(5606,-7878)--(5613,-7835)--(5619,-7791)--(5625,-7746)
  --(5630,-7701)--(5635,-7655)--(5639,-7609)--(5643,-7562)--(5646,-7516)--(5649,-7469)
  --(5652,-7423)--(5654,-7377)--(5655,-7331)--(5656,-7286)--(5657,-7242)--(5657,-7199)
  --(5656,-7156)--(5656,-7115)--(5654,-7075)--(5653,-7035)--(5651,-6997)--(5648,-6960)
  --(5646,-6925)--(5642,-6890)--(5639,-6857)--(5635,-6825)--(5631,-6794)--(5625,-6756)
  --(5618,-6720)--(5611,-6685)--(5603,-6651)--(5594,-6618)--(5585,-6587)--(5575,-6557)
  --(5564,-6527)--(5552,-6499)--(5539,-6472)--(5526,-6446)--(5512,-6421)--(5497,-6398)
  --(5482,-6376)--(5466,-6355)--(5450,-6335)--(5433,-6316)--(5415,-6299)--(5398,-6283)
  --(5379,-6269)--(5361,-6255)--(5342,-6243)--(5324,-6232)--(5305,-6222)--(5286,-6213)
  --(5266,-6205)--(5247,-6197)--(5228,-6191)--(5208,-6185)--cycle;
\pgfsetfillcolor{black}
\pgftext[base,left,at=\pgfqpointxy{2625}{-4425}] {\fontsize{24}{28.8}\normalfont $p_{e'}$}
\pgftext[base,left,at=\pgfqpointxy{8325}{-3750}] {\fontsize{24}{28.8}\normalfont $e'$}
\pgftext[base,left,at=\pgfqpointxy{7350}{-3375}] {\fontsize{24}{28.8}\normalfont $u$}
\pgftext[base,left,at=\pgfqpointxy{8325}{-6225}] {\fontsize{24}{28.8}\normalfont $w$}
\pgftext[base,left,at=\pgfqpointxy{7650}{-6450}] {\fontsize{24}{28.8}\normalfont $e_1$}
\pgftext[base,left,at=\pgfqpointxy{8400}{-7800}] {\fontsize{24}{28.8}\normalfont $w'$}
\pgftext[base,left,at=\pgfqpointxy{9000}{-6450}] {\fontsize{24}{28.8}\normalfont $e_2$}
\pgftext[base,left,at=\pgfqpointxy{7350}{-6900}] {\fontsize{24}{28.8}\normalfont $u$}
\pgftext[base,left,at=\pgfqpointxy{9000}{-7425}] {\fontsize{24}{28.8}\normalfont $e'_2$}
\pgftext[base,left,at=\pgfqpointxy{7650}{-7500}] {\fontsize{24}{28.8}\normalfont $e'_1$}
\pgftext[base,left,at=\pgfqpointxy{1650}{-5550}] {\fontsize{24}{28.8}\normalfont $\widetilde S$}
\pgftext[base,left,at=\pgfqpointxy{3375}{-9075}] {\fontsize{24}{28.8}\normalfont $\P^1_{w'}$}
\pgftext[base,left,at=\pgfqpointxy{1500}{-1575}] {\fontsize{24}{28.8}\normalfont $S_{0}$}
\pgftext[base,left,at=\pgfqpointxy{3675}{-5700}] {\fontsize{24}{28.8}\normalfont $\P^1_w$}
\pgftext[base,left,at=\pgfqpointxy{2850}{-1950}] {\fontsize{24}{28.8}\normalfont $p_e$}
\pgftext[base,left,at=\pgfqpointxy{9075}{-3450}] {\fontsize{24}{28.8}\normalfont $v$}
\pgftext[base,left,at=\pgfqpointxy{8775}{-6900}] {\fontsize{24}{28.8}\normalfont $v$}
\pgftext[base,left,at=\pgfqpointxy{7125}{-5700}] {\fontsize{24}{28.8}\normalfont $\~G=(\~V,\~E)$}
\pgfsetlinewidth{+30\XFigu}
\draw (2100,-2325) arc[start angle=+-138.4, end angle=+-41.6, radius=+451.6];
\endtikzpicture}%

%% file: greenfunctions.bbl
\begin{bibdiv}
\begin{biblist}

\bib{AB15}{article}{
      author={Amini, Omid},
      author={Baker, Matthew},
       title={Linear series on metrized complexes of algebraic curves},
        date={2015},
     journal={Mathematische Annalen},
      volume={362},
      number={1-2},
       pages={55\ndash 106},
}

\bib{ABBF}{article}{
      author={Amini, Omid},
      author={Bloch, Spencer},
      author={{Burgos Gil}, Jos\'e~Ignacio},
      author={Fres{\'a}n, Javier},
       title={Feynman amplitudes and limits of heights},
        date={2016},
     journal={Izvestiya: Mathematics},
      volume={80},
      number={5},
       pages={813\ndash 848, special volume in honor of {J}ean\ndash {P}ierre
  {S}erre},
}

\bib{ACG}{book}{
      author={Arbarello, Enrico},
      author={Cornalba, Maurizio},
      author={Griffiths, Phillip},
       title={Geometry of algebraic curves: volume {II}},
   publisher={Springer Science \& Business Media},
        date={2011},
      volume={268},
}

\bib{ACGH}{book}{
      author={Arbarello, Enrico},
      author={Cornalba, Maurizio},
      author={Griffiths, Phillip~A.},
      author={Harris, Joe},
       title={Geometry of algebraic curves},
   publisher={Springer},
        date={1984},
      volume={267},
}

\bib{ACP}{article}{
      author={Abramovich, Dan},
      author={Caporaso, Lucia},
      author={Payne, Sam},
       title={The tropicalization of the moduli space of curves},
        date={2015},
        ISSN={0012-9593},
     journal={Ann. Sci. \'{E}c. Norm. Sup\'{e}r. (4)},
      volume={48},
      number={4},
       pages={765\ndash 809},
  url={https://doi-org.accesdistant.sorbonne-universite.fr/10.24033/asens.2258},
}

\bib{AI}{article}{
      author={Amini, Omid},
      author={Iriarte, Hern\'an},
       title={Geometry of higher rank valuations},
        date={2022},
     journal={Preprint},
}

\bib{AN}{article}{
      author={Amini, Omid},
      author={Nicolussi, Noema},
       title={Moduli of hybrid curves {I}: variations of canonical measures},
        date={2020},
     journal={Arxiv preprint arXiv:2007.07130},
}

\bib{AN-AG-hybrid}{article}{
      author={Amini, Omid},
      author={Nicolussi, Noema},
       title={Moduli of hybrid curves {III}: Algebraic geometry of hybrid
  curves},
        date={2022},
     journal={Preprint in preparation},
}

\bib{AP}{article}{
      author={Amini, Omid},
      author={Piquerez, Matthieu},
       title={Hodge theory for tropical varieties},
        date={2020},
     journal={arXiv preprint arXiv:2007.07826},
}

\bib{AP-CS}{article}{
      author={Amini, Omid},
      author={Piquerez, Matthieu},
       title={Tropical {C}lemens-{S}chmid sequence and existence of tropical
  cycles with a given cohomology class},
        date={2020},
     journal={arXiv preprint arXiv:2012.13142},
}

\bib{AP-homology}{article}{
      author={Amini, Omid},
      author={Piquerez, Matthieu},
       title={Homology of tropical fans},
        date={2021},
     journal={arXiv preprint arXiv:2105.01504},
}

\bib{Ara74}{article}{
      author={Arakelov, Souren},
       title={Intersection theory of divisors on an arithmetic surface},
        date={1974},
     journal={Mathematics of the USSR-Izvestiya},
      volume={8},
      number={6},
       pages={1167},
}

\bib{BF11}{article}{
      author={Baker, Matthew},
      author={Faber, Xander},
       title={Metric propertis of the tropical {A}bel--{J}acobi map},
        date={2011},
     journal={Journal of Algebraic Combinatorics},
      volume={33},
       pages={349\ndash 381},
}

\bib{BGS94}{article}{
      author={Bost, J-B},
      author={Gillet, Henri},
      author={Soul{\'e}, Christophe},
       title={Heights of projective varieties and positive {G}reen forms},
        date={1994},
     journal={Journal of the American Mathematical Society},
       pages={903\ndash 1027},
}

\bib{BJ17}{article}{
      author={Boucksom, S{\'e}bastien},
      author={Jonsson, Mattias},
       title={Tropical and non-{A}rchimedean limits of degenerating families of
  volume forms},
        date={2017},
     journal={Journal de l'{\'E}cole {P}olytechnique -- Math{\'e}matiques},
      volume={4},
       pages={87\ndash 139},
}

\bib{BerkolaikoKuchment}{book}{
      author={Berkolaiko, Gregory},
      author={Kuchment, Peter},
       title={Introduction to quantum graphs},
   publisher={American Mathematical Soc.},
        date={2013},
}

\bib{bls}{article}{
      author={Berkolaiko, Gregory},
      author={Latushkin, Yuri},
      author={Sukhtaiev, Selim},
       title={Limits of quantum graph operators with shrinking edges},
        date={2019},
     journal={Advances in Mathematics},
      volume={352},
       pages={632\ndash 669},
}

\bib{BM}{book}{
      author={Bondy, John~Adrian},
      author={Murty, U. S.~Ramachandra},
       title={Graph theory},
   publisher={Springer},
        date={2008},
      volume={244},
}

\bib{Bol}{book}{
      author={Bollob{\'a}s, B{\'e}la},
       title={Modern graph theory},
   publisher={Springer},
        date={1998},
      volume={184},
}

\bib{Borisov}{article}{
      author={Borisov, Denis~Ivanovich},
       title={Analyticity of resolvents of elliptic operators on quantum graphs
  with small edges},
        date={2021},
     journal={Advances in Mathematics, to appear; arxiv:2102.12958.},
}

\bib{Bost99}{article}{
      author={Bost, Jean-Beno\^it},
       title={Potential theory and {L}efschetz theorems for arithmetic
  surfaces},
        date={1999},
     journal={Annales scientifiques de l'\'Ecole normale sup{\'e}rieure},
      volume={32},
      number={2},
       pages={241\ndash 312},
}

\bib{BR07}{article}{
      author={Baker, Matthew},
      author={Rumely, Robert},
       title={Harmonic analysis on metrized graphs},
        date={2007},
     journal={Canadian Journal of Mathematics},
      volume={59},
      number={2},
       pages={225\ndash 275},
}

\bib{BR10}{book}{
      author={Baker, Matthew},
      author={Rumely, Robert},
       title={Potential theory and dynamics on the berkovich projective line},
      series={Mathematical Surveys and Monographs},
   publisher={American Mathematical Soc.},
     address={Providence, RI},
        date={2010},
      volume={159},
}

\bib{cacciapuoti}{article}{
      author={Cacciapuoti, Claudio},
       title={Scale invariant effective hamiltonians for a graph with a small
  compact core},
        date={2019},
     journal={Symmetry},
      volume={11},
}

\bib{Cinkir}{article}{
      author={Cinkir, Zubeyir},
       title={{Zhang's conjecture and the effective Bogomolov conjecture over
  function fields}},
        date={2011},
     journal={Inventiones mathematicae},
      volume={183},
      number={3},
       pages={517\ndash 562},
}

\bib{Cin14}{article}{
      author={Cinkir, Zubeyir},
       title={Explicit computation of certain {A}rakelov--{G}reen functions},
        date={2014},
     journal={Kyoto Journal of Mathematics},
      volume={54},
      number={4},
       pages={759\ndash 774},
}

\bib{CKS}{article}{
      author={Cattani, Eduardo},
      author={Kaplan, Aroldo},
      author={Schmid, Wilfried},
       title={Degeneration of {H}odge structures},
        date={1986},
     journal={Annals of Mathematics},
      volume={123},
      number={3},
       pages={457\ndash 535},
}

\bib{CR93}{article}{
      author={Chinburg, Ted},
      author={Rumely, Robert},
       title={The capacity pairing},
        date={1993},
     journal={Journal f{\"u}r die Reine und Angewandte Mathematik},
      volume={1993},
      number={434},
       pages={1\ndash 44},
}

\bib{HGB19}{article}{
      author={D'Hoker, Eric},
      author={Green, Michael~B},
      author={Pioline, Boris},
       title={Higher genus modular graph functions, string invariants, and
  their exact asymptotics},
        date={2019},
     journal={Communications in Mathematical Physics},
      volume={366},
      number={3},
       pages={927\ndash 979},
}

\bib{Diestel}{book}{
      author={Diestel, Reinhard},
       title={Graph theory},
   publisher={Springer},
        date={2000},
      volume={173},
}

\bib{deJong}{article}{
      author={de~Jong, Robin},
       title={Faltings delta-invariant and semistable degeneration},
        date={2019},
     journal={Journal of Differential Geometry},
      volume={111},
      number={2},
       pages={241\ndash 301},
}

\bib{DKY20}{article}{
      author={DeMarco, Laura},
      author={Krieger, Holly},
      author={Ye, Hexi},
       title={{Uniform Manin-Mumford for a family of genus 2 curves}},
        date={2020},
     journal={Annals of Mathematics},
      volume={191},
      number={3},
       pages={949\ndash 1001},
}

\bib{Edm70}{inproceedings}{
      author={Edmonds, Jack},
       title={Submodular functions, matroids, and certain polyhedra},
        date={1970},
   booktitle={Proc. calgary internat. conf., calgary, alta., 1969},
   publisher={Gordon and Breach, New York},
       pages={69\ndash 87},
}

\bib{Faltings21}{article}{
      author={Faltings, Gerd},
       title={Arakelov geometry on degenerating curves},
        date={2021},
     journal={Journal f{\"u}r die reine und angewandte Mathematik},
      volume={771},
       pages={65\ndash 84},
}

\bib{Fal84}{article}{
      author={Faltings, Gerd},
       title={Calculus on arithmetic surfaces},
        date={1984},
     journal={Annals of mathematics},
      volume={119},
       pages={387\ndash 424},
}

\bib{FK92}{book}{
      author={Farkas, Hershel},
      author={Kra, Irwin},
       title={Riemann surfaces},
   publisher={Springer},
        date={1992},
}

\bib{Ful71}{article}{
      author={Fulkerson, Delbert~Ray},
       title={Blocking and anti-blocking pairs of polyhedra},
        date={1971},
     journal={Mathematical programming},
      volume={1},
      number={1},
       pages={168\ndash 194},
}

\bib{GH14}{book}{
      author={Griffiths, Phillip},
      author={Harris, Joseph},
       title={Principles of algebraic geometry},
   publisher={John Wiley \& Sons},
        date={2014},
}

\bib{GS90}{article}{
      author={Gillet, Henri},
      author={Soul{\'e}, Christophe},
       title={Arithmetic intersection theory},
        date={1990},
     journal={Publications Math{\'e}matiques de l'IH{\'E}S},
      volume={72},
       pages={93\ndash 174},
}

\bib{Hain}{article}{
      author={Hain, Richard},
       title={Biextensions and heights associated to curves of odd genus},
        date={1990},
     journal={Duke Mathematical Journal},
      volume={61},
      number={3},
       pages={859\ndash 898},
}

\bib{HK14}{article}{
      author={Hubbard, John~H},
      author={Koch, Sarah},
       title={An analytic construction of the {D}eligne-{M}umford
  compactification of the moduli space of curves},
        date={2014},
     journal={Journal of Differential Geometry},
      volume={98},
      number={2},
       pages={261\ndash 313},
}

\bib{Hof84}{book}{
      author={Hoffman, Jerome~William},
       title={The {H}odge theory of stable curves},
   publisher={American Mathematical Soc.},
        date={1984},
      volume={308},
}

\bib{Iri22}{article}{
      author={Iriarte, Hern\'an},
       title={Higher rank polyhedral and tropical geometry},
        date={2022},
     journal={Preprint},
}

\bib{JK06}{article}{
      author={Jorgenson, Jay},
      author={Kramer, Jurg},
       title={Bounds on canonical green's functions},
        date={2006},
     journal={Compositio Mathematica},
      volume={142},
      number={3},
       pages={679\ndash 700},
}

\bib{JK09}{article}{
      author={Jorgenson, Jay},
      author={Kramer, J{\"u}rg},
       title={Bounds on {F}altings's delta function through covers},
        date={2009},
     journal={Annals of mathematics},
       pages={1\ndash 43},
}

\bib{JW92}{article}{
      author={Ji, Lizhen},
      author={Wentworth, Richard},
       title={Spectral convergence on degenerating surfaces},
        date={1992},
     journal={Duke Mathematical Journal},
      volume={66},
      number={3},
       pages={469\ndash 501},
}

\bib{KNU}{article}{
      author={Kato, Kazuya},
      author={Nakayama, Chikara},
      author={Usui, Sampei},
       title={$sl(2)$-orbit theorem for degeneration of mixed {H}odge
  structure},
        date={2008},
     journal={Journal of Algebraic Geometry},
      volume={17},
       pages={401\ndash 479},
}

\bib{Lang}{book}{
      author={Lang, Serge},
       title={Introduction to {A}rakelov theory},
   publisher={Springer},
        date={2012},
}

\bib{Mas76}{article}{
      author={Masur, Howard},
       title={The extension of the {W}eil-{P}etersson metric to the boundary of
  {T}eichmuller space},
        date={1976},
     journal={Duke Mathematical Journal},
      volume={43},
      number={3},
       pages={623\ndash 635},
}

\bib{OR11}{article}{
      author={Osserman, Brian},
      author={Rabinoff, Joseph},
       title={Lifting non-proper tropical intersections},
        date={2011},
     journal={Tropical and non-Archimedean geometry},
      volume={605},
       pages={15\ndash 44},
}

\bib{Pearl}{article}{
      author={Pearlstein, Gregory},
       title={Degenerations of mixed {H}odge structure},
        date={2001},
     journal={Duke Mathematical Journal},
      volume={110},
      number={2},
       pages={217\ndash 251},
}

\bib{Pearl2}{article}{
      author={Pearlstein, Gregory},
       title={$\textrm{SL}_2$-orbits and degenerations of mixed {H}odge
  structure},
        date={2006},
     journal={Journal of Differential Geometry},
      volume={74},
      number={1},
       pages={1\ndash 67},
}

\bib{Pos09}{article}{
      author={Postnikov, Alexander},
       title={Permutohedra, associahedra, and beyond},
        date={2009},
     journal={International Mathematics Research Notices},
      volume={2009},
      number={6},
       pages={1026\ndash 1106},
}

\bib{PRW08}{article}{
      author={Postnikov, Alex},
      author={Reiner, Victor},
      author={Williams, Lauren},
       title={Faces of {G}eneralized {P}ermutohedra},
        date={2008},
     journal={Documenta Mathematica},
      volume={13},
       pages={207\ndash 273},
}

\bib{Rad52}{article}{
      author={Rado, Richard},
       title={An inequality},
        date={1952},
     journal={Journal of the London Mathematical Society},
      volume={1},
      number={1},
       pages={1\ndash 6},
}

\bib{Smit88}{article}{
      author={Smit, Dirk-Jan},
       title={String theory and algebraic geometry of moduli spaces},
        date={1988},
     journal={Communications in mathematical physics},
      volume={114},
      number={4},
       pages={645\ndash 685},
}

\bib{SW19}{article}{
      author={Shokrieh, Farbod},
      author={Wu, Chenxi},
       title={Canonical measures on metric graphs and a {K}azhdan's theorem},
        date={2019},
     journal={Inventiones mathematicae},
      volume={215},
       pages={819\ndash 862},
}

\bib{Went08}{article}{
      author={Wentworth, Richard~A},
       title={Precise constants in bosonization formulas on {R}iemann surfaces.
  i},
        date={2008},
     journal={Communications in mathematical physics},
      volume={282},
      number={2},
       pages={339\ndash 355},
}

\bib{Went91}{article}{
      author={Wentworth, Richard},
       title={The asymptotics of the {A}rakelov--{G}reen's function and
  {F}altings' delta invariant},
        date={1991},
     journal={Communications in mathematical physics},
      volume={137},
      number={3},
       pages={427\ndash 459},
}

\bib{Wil17}{article}{
      author={Wilms, Robert},
       title={New explicit formulas for {F}altings' delta-invariant},
        date={2017},
     journal={Inventiones mathematicae},
      volume={209},
      number={2},
       pages={481\ndash 539},
}

\bib{Wol90}{article}{
      author={Wolpert, Scott~A},
       title={The hyperbolic metric and the geometry of the universal curve},
        date={1990},
     journal={Journal of Differential Geometry},
      volume={31},
      number={2},
       pages={417\ndash 472},
}

\bib{Zhang2}{article}{
      author={Zhang, Shou-Wu},
       title={{Gross--Schoen cycles and dualising sheaves}},
        date={2010},
     journal={Inventiones mathematicae},
      volume={179},
      number={1},
       pages={1},
}

\bib{Zhang}{article}{
      author={Zhang, Shou-Wu},
       title={Admissible pairing on a curve},
        date={1993},
     journal={Inventiones mathematicae},
      volume={112},
      number={1},
       pages={171\ndash 193},
}

\end{biblist}
\end{bibdiv}
